\documentclass{amsart}
\usepackage{amssymb}
\usepackage{mathrsfs}
\usepackage{stmaryrd}
\usepackage{bbm}
\usepackage{oldgerm}
\usepackage[francais]{babel}
\usepackage[T1]{fontenc}
\usepackage[latin1]{inputenc}
\usepackage[all]{xy}
\usepackage{hyperref}

\makeindex

\newtheorem{teo}[subsection]{Théorème}
\newtheorem{prop}[subsection]{Proposition}
\newtheorem{cor}[subsection]{Corollaire}
\newtheorem{lem}[subsection]{Lemme}

\theoremstyle{definition}

\newtheorem{defi}[subsection]{Définition}
\newtheorem{rema}[subsection]{Remarque}
\newtheorem{remas}[subsection]{Remarques}

\numberwithin{equation}{subsection}

\newcommand{\mC}{{\mathbb C}}
\newcommand{\mD}{{\mathbb D}}
\newcommand{\mH}{{\mathbb H}}
\newcommand{\mQ}{{\mathbb Q}}
\newcommand{\mN}{{\mathbb N}}
\newcommand{\mT}{{\mathbb T}}
\newcommand{\mV}{{\mathbb V}}
\newcommand{\mZ}{{\mathbb Z}}
\newcommand{\mF}{{\mathbb F}}

\newcommand{\mK}{{\mathbb K}}

\newcommand{\bA}{{\bf A}}

\newcommand{\bD}{{\bf D}}
\newcommand{\bL}{{\bf L}}
\newcommand{\bT}{{\bf T}}

\newcommand{\Sch}{{\bf Sch}}
\newcommand{\Mon}{{\bf Mon}}
\newcommand{\Ab}{{\bf Ab}}

\newcommand{\bMod}{{\bf Mod}}
\newcommand{\bRep}{{\bf Rep}}
\newcommand{\bMH}{{\bf MH}}
\newcommand{\bTors}{{\bf Tors}}
\newcommand{\bFPH}{{\bf FPH}}

\newcommand{\bFl}{{\bf Fl}}
\newcommand{\bGr}{{\bf Gr}}
\newcommand{\bOp}{{\bf Op}}

\newcommand{\zar}{{\rm zar}}

\newcommand{\HT}{{\rm HT}}

\newcommand{\Spec}{{\rm Spec}}
\newcommand{\Max}{{\rm Max}}
\newcommand{\Spf}{{\rm Spf}}
\newcommand{\ob}{{\rm Ob}}
\newcommand{\pr}{{\rm pr}}

\newcommand{\lib}{{\rm lib}}

\newcommand{\disc}{{\rm disc}}
\newcommand{\coker}{{\rm coker}}
\newcommand{\Ann}{{\rm Ann}}

\newcommand{\integre}{{\rm int}}
\newcommand{\sat}{{\rm sat}}
\newcommand{\tot}{{\rm tot}}
\newcommand{\qpp}{{\rm qpp}}
\newcommand{\p}{{\rm p}}
\newcommand{\trqpp}{{\rm \text -qpp}}
\newcommand{\trp}{{\rm \text -p}}
\newcommand{\cont}{{\rm cont}}

\newcommand{\gp}{{\rm gp}}
\newcommand{\id}{{\rm id}}
\newcommand{\Tr}{{\rm Tr}}
\newcommand{\rg}{{\rm rg}}

\newcommand{\Ind}{{\rm Ind}}

\newcommand{\Hom}{{\rm Hom}}
\newcommand{\Homgr}{{\rm Homgr}}
\newcommand{\End}{{\rm End}}
\newcommand{\GL}{{\rm GL}}
\newcommand{\Mat}{{\rm Mat}}
\newcommand{\Aut}{{\rm Aut}}

\newcommand{\Gal}{{\rm Gal}}

\newcommand{\rA}{{\rm A}}
\newcommand{\rC}{{\rm C}}
\newcommand{\rK}{{\rm K}}
\newcommand{\rD}{{\rm D}}
\newcommand{\rF}{{\rm F}}
\newcommand{\rH}{{\rm H}}
\newcommand{\rT}{{\rm T}}
\newcommand{\rL}{{\rm L}}

\newcommand{\rP}{{\rm P}}
\newcommand{\rR}{{\rm R}}
\newcommand{\rS}{{\rm S}}
\newcommand{\rV}{{\rm V}}
\newcommand{\rW}{{\rm W}}

\newcommand{\rb}{{\rm b}}

\newcommand{\oF}{{\overline{F}}}
\newcommand{\oK}{{\overline{K}}}
\newcommand{\oL}{{\overline{L}}}

\newcommand{\oR}{{\overline{R}}}
\newcommand{\oS}{{\overline{S}}}

\newcommand{\oX}{{\overline{X}}}

\newcommand{\ok}{{\overline{k}}}
\newcommand{\op}{{\overline{p}}}

\newcommand{\ou}{{\overline{u}}}

\newcommand{\ox}{{\overline{x}}}
\newcommand{\oy}{{\overline{y}}}
\newcommand{\oz}{{\overline{z}}}

\newcommand{\ofc}{{\overline{\fc}}}

\newcommand{\oeta}{{\overline{\eta}}}
\newcommand{\ophi}{{\overline{\phi}}}
\newcommand{\osigma}{{\overline{\sigma}}}
\newcommand{\ogamma}{{\overline{\gamma}}}
\newcommand{\olambda}{{\overline{\lambda}}}

\newcommand{\uA}{{\underline{A}}}
\newcommand{\uG}{{\underline{G}}}
\newcommand{\uY}{{\underline{Y}}}

\newcommand{\um}{{\underline{m}}}
\newcommand{\un}{{\underline{n}}}
\newcommand{\upp}{{\underline{p}}}
\newcommand{\uy}{{\underline{y}}}
\newcommand{\utheta}{{\underline{\theta}}}

\newcommand{\hY}{{\widehat{Y}}}

\newcommand{\hRun}{{\widehat{R_1}}}

\newcommand{\hRi}{{\widehat{R_\infty}}}
\newcommand{\hRpi}{{\widehat{R_{p^\infty}}}}
\newcommand{\hOmega}{\widehat{\Omega}}

\newcommand{\hmZ}{{\widehat{\mZ}}}

\newcommand{\halpha}{{\widehat{\alpha}}}

\newcommand{\cA}{{\mathscr A}}
\newcommand{\cB}{{\mathscr B}}
\newcommand{\cC}{{\mathscr C}}
\newcommand{\cD}{{\mathscr D}}
\newcommand{\cE}{{\mathscr E}}
\newcommand{\cF}{{\mathscr F}}
\newcommand{\cG}{{\mathscr G}}
\newcommand{\cI}{{\mathscr I}}

\newcommand{\cK}{{\mathscr K}}
\newcommand{\cL}{{\mathscr L}}
\newcommand{\cP}{{\mathscr P}}
\newcommand{\co}{{\mathscr O}}
\newcommand{\cR}{{\mathscr R}}
\newcommand{\cS}{{\mathscr S}}

\newcommand{\cH}{{\mathscr H}}
\newcommand{\cM}{{\mathscr M}}
\newcommand{\cN}{{\mathscr N}}
\newcommand{\cQ}{{\mathscr Q}}
\newcommand{\cZ}{{\mathscr Z}}

\newcommand{\cHom}{{\mathscr Hom}}
\newcommand{\cEnd}{{\mathscr End}}

\newcommand{\fD}{{\mathfrak D}}

\newcommand{\fS}{{\mathfrak S}}

\newcommand{\fc}{{\mathfrak c}}
\newcommand{\fm}{{\mathfrak m}}
\newcommand{\fo}{{\mathfrak o}}
\newcommand{\fp}{{\mathfrak p}}
\newcommand{\fq}{{\mathfrak q}}

\newcommand{\tta}{{\tt a}}
\newcommand{\ttt}{{\tt t}}

\newcommand{\hA}{{\widehat{A}}}
\newcommand{\hB}{{\widehat{B}}}
\newcommand{\hM}{{\widehat{M}}}
\newcommand{\hR}{{\widehat{R}}}
\newcommand{\hoR}{{\widehat{\oR}}}

\newcommand{\hfS}{{\widehat{\fS}}}
\newcommand{\hcS}{{\widehat{\cS}}}

\newcommand{\hcC}{{\widehat{\cC}}}

\newcommand{\htta}{{\widehat{\tta}}}

\newcommand{\hotimes}{{\widehat{\otimes}}}

\newcommand{\tU}{{\widetilde{U}}}

\newcommand{\tX}{{\widetilde{X}}}
\newcommand{\tY}{{\widetilde{Y}}}
\newcommand{\tth}{{\widetilde{h}}}
\newcommand{\tlt}{{\widetilde{t}}}

\newcommand{\tq}{{\widetilde{q}}}

\newcommand{\tx}{{\widetilde{x}}}
\newcommand{\ty}{{\widetilde{y}}}
\newcommand{\tz}{{\widetilde{z}}}
\newcommand{\tOmega}{{\widetilde{\Omega}}}
\newcommand{\tbeta}{{\widetilde{\beta}}}

\newcommand{\tgamma}{{\widetilde{\gamma}}}
\newcommand{\tdelta}{{\widetilde{\delta}}}
\newcommand{\tkappa}{{\widetilde{\kappa}}}

\newcommand{\tpi}{{\widetilde{\pi}}}

\newcommand{\tchi}{{\widetilde{\chi}}}
\newcommand{\tlambda}{{\widetilde{\lambda}}}

\newcommand{\tcC}{{\widetilde{\cC}}}

\newcommand{\trK}{{\widetilde{\rK}}}
\newcommand{\tmK}{{\widetilde{\mK}}}
\newcommand{\trT}{{\widetilde{\rT}}}

\newcommand{\coS}{{\check{\oS}}}
\newcommand{\coX}{{\check{\oX}}}

\begin{document}

\title[Sur la correspondance de Simpson $p$-adique. I]
{Sur la correspondance de Simpson $p$-adique.\\
I~: étude locale}
\author{Ahmed Abbes et Michel Gros}
\address{A.A. Laboratoire Alexander Grothendieck, ERL 9216 du CNRS, Institut des Hautes \'Etudes Scientifiques, 35 route de Chartres, 91440 Bures-sur-Yvette, France}
\address{M.G. CNRS UMR 6625, IRMAR, Université de Rennes 1,
Campus de Beaulieu, 35042 Rennes cedex, France}
\email{abbes@ihes.fr}
\email{michel.gros@univ-rennes1.fr}

\maketitle

\setcounter{tocdepth}{1}
\tableofcontents

\section{Introduction}

Le présent article est consacré à la construction et à l'étude de la {\em correspondance de Simpson $p$-adique},
suivant l'approche générale résumée dans \cite{ag0}, pour un schéma logarithmique affine d'un type particulier \eqref{higgs1-dlog1}.
La section~\ref{higgs1-NC} contient les principales notations et conventions générales.  
On trouvera en particulier dans \ref{higgs1-not6} les définitions relatives aux modules de Higgs. 
La section~\ref{higgs1-kun} contient divers sorites sur la cohomologie continue des groupes profinis 
et en particulier, faute de référence adéquate, un traitement approprié de la formule de Künneth.
La section~\ref{higgs1-eph} rappelle et précise quant à elle les relations existant entre torseurs, pour la topologie 
de Zariski et fibrés principaux homogènes, 
ainsi qu'entre les notions équivariantes (sous un groupe abstrait) associées.
Nous rappelons ensuite dans la section~\ref{higgs1-LOG} quelques notions de géométrie logarithmique 
qui joueront un rôle important dans la suite du travail, afin de fixer les notations et de donner 
des repères aux lecteurs non familiers de cette théorie.  
Dans la section~\ref{higgs1-pur}, nous introduisons le cadre logarithmique \eqref{higgs1-dlog1}, 
les anneaux \eqref{higgs1-dlog5} et les groupes de Galois \eqref{higgs1-gal2} utilisés tout le long de l'article 
puis nous établissons quelques-unes de leurs propriétés (\ref{higgs1-dlog4}, \ref{higgs1-gal1} et \ref{higgs1-pur81}). 
Nous rappelons ensuite l'énoncé du théorème de presque pureté de Faltings \eqref{higgs1-pur44} et 
quelques corollaires \eqref{higgs1-pur4}--\eqref{higgs1-cg8}.
La section~\ref{higgs1-ext} est consacrée à l'extension de Faltings~; on y trouve deux variantes 
\eqref{higgs1-ext-log12b} et \eqref{higgs1-log-ext17b}. Pour la commodité du lecteur, 
nous reprenons dans la section~\ref{higgs1-CG}, avec un peu de détails, le  calcul de cohomologie galoisienne, dû à Faltings, 
qui est au coeur de son approche de la théorie de Hodge $p$-adique. 
Nous introduisons dans la section~\ref{higgs1-EIPF} les épaississements infinitésimaux $p$-adiques de Fontaine,  
et nous les munissons de structures logarithmiques suivant Tsuji \eqref{higgs1-eip9}. 
La partie la plus originale de cet article commence à partir de la section~\ref{higgs1-TOR}, 
avec l'introduction du {\em torseur de Higgs-Tate} \eqref{higgs1-tor2}. Le théorème \ref{higgs1-ext246} établit un important lien 
entre ce torseur et l'extension de Faltings. La section~\ref{higgs1-CGII} est consacrée à l'étude des cohomologies 
de de Rham \eqref{higgs1-cohbis8}  et galoisienne \eqref{higgs1-cohbis9}  de diverses algèbres associées à ce torseur.  
Dans la section~\ref{higgs1-dolbeault}, nous définissons les principaux foncteurs \eqref{higgs1-dolb2c} et \eqref{higgs1-dolb3c} reliant  
la catégorie des représentations généralisées à celle des modules de Higgs. 
Nous y  introduisons  aussi les notions de représentation de Dolbeault et de module de Higgs soluble. 
Nous développons en fait deux variantes, une entière (\ref{higgs1-dolb21}  et \ref{higgs1-dolb31}) et une 
rationnelle plus subtile (\ref{higgs1-drt2} et \ref{higgs1-drt3}).
Nous montrons alors que pour chaque variante, ces notions donnent lieu à deux catégories équivalentes (\ref{higgs1-dolb7} et \ref{higgs1-drt16}). 
La section~\ref{higgs1-RP} est consacrée à l'étude des petites représentations et des petits modules de Higgs 
suivant l'approche de Faltings \cite{faltings3}. Nous établissons aussi des liens 
entre ces notions et les notions de représentation de Dolbeault et de module de Higgs soluble 
(\ref{higgs1-higgs131}, \ref{higgs1-drt20} et \ref{higgs1-drt21}). 
La section \ref{higgs1-descente} contient un énoncé de descente pour les petites représentations du à Faltings \eqref{higgs1-desc4}. 
Nous en déduisons de nouveaux liens entre les différentes notions de représentations et de module de Higgs 
introduites précédemment (\ref{higgs1-desc5}, \ref{higgs1-dsct4} et \ref{higgs1-dsct3}).   
La dernière section fait le lien \eqref{higgs1-hyodo10} avec la théorie de Hyodo
pour les représentations de Hodge-Tate.

\section{Notations et conventions}\label{higgs1-NC}

{\em Tous les anneaux considérés dans cet article possèdent un élément unité~;
les homomorphismes d'anneaux sont toujours supposés transformer l'élément unité en l'élément unité.
Nous considérons surtout des anneaux commutatifs, et lorsque nous parlons d'anneau 
sans préciser, il est sous-entendu qu'il s'agit d'un anneau commutatif~; en particulier, 
il est sous-entendu, lorsque nous parlons d'un topos annelé $(X,A)$ sans préciser, que $A$ est commutatif.}

\subsection{}\label{higgs1-not1}\index{10200@$K$, $\oK$, $C$}
\index{10201@$\hA$ ($A$ groupe abélien)}
Dans cet article, $p$ désigne un nombre premier, $K$ un corps de valuation discrète complet de caractéristique
$0$, à corps résiduel parfait $k$ de caractéristique $p$, $\oK$ une clôture algébrique de $K$.
On note $\co_K$ l'anneau de valuation de $K$,  
$\co_\oK$ la clôture intégrale de $\co_K$ dans $\oK$, $\fm_\oK$ l'idéal maximal de $\co_\oK$,
$\ok$ le corps résiduel de $\co_\oK$ et $v$ la valuation de $\oK$ normalisée par $v(p)=1$. 
On désigne par $\co_C$ le séparé complété $p$-adique de $\co_\oK$, par $C$ son corps de fractions 
et par $\fm_C$ son idéal maximal.

On choisit un système compatible $(\beta_n)_{n>0}$ de racines $n$-ièmes de $p$ dans $\co_\oK$. 
Pour tout nombre rationnel $\varepsilon>0$, 
on pose $p^\varepsilon=(\beta_n)^{\varepsilon n}$ où $n$ est un entier $>0$ tel que $\varepsilon n$ soit entier.

On désigne par $G_K=\Gal(\oK/K)$ le groupe de Galois de $\oK$ sur $K$ et par $\hmZ(1)$ et $\mZ_p(1)$ 
les $\mZ[G_K]$-modules 
\begin{eqnarray}
\hmZ(1)&=&\underset{\underset{n\geq 1}{\longleftarrow}}{\lim}\ \mu_{n}(\co_{\oK}),\label{higgs1-not4a}\\
\mZ_p(1)&=&\underset{\underset{n\geq 0}{\longleftarrow}}{\lim}\ \mu_{p^n}(\co_{\oK}),\label{higgs1-not4b}
\end{eqnarray}  
où $\mu_n(\co_{\oK})$ désigne le sous-groupe des racines $n$-ièmes de l'unité dans $\co_\oK$. 
Pour tout $\mZ_p[G_K]$-module $M$ et tout entier $n$, on pose 
$M(n)=M\otimes_{\mZ_p}\mZ_p(1)^{\otimes n}$.

Pour tout groupe abélien $A$, on note $\hA$ son séparé complété $p$-adique.

\subsection{}\label{higgs1-not54}
On munit  $\mZ_p$ de la topologie $p$-adique, ainsi que toutes les $\mZ_p$-algèbres adiques 
({\em i.e.}, les $\mZ_p$-algèbres complètes et séparées pour la topologie $p$-adique).
Soient $A$ une $\mZ_p$-algèbre adique, $i\colon A\rightarrow A[\frac 1 p]$ l'homomorphisme canonique. 
On appelle {\em topologie $p$-adique} sur $A[\frac 1 p]$ l'unique topologie 
compatible avec sa structure de groupe additif  pour laquelle les sous-groupes $i(p^nA)$,
pour $n\in \mN$, forment un système fondamental de voisinages de $0$ (\cite{tg} chap.~III §1.2, prop.~1). 
Elle fait de $A[\frac 1 p]$ un anneau topologique. 
Soient $M$ un $A[\frac 1 p]$-module de type fini, $M^\circ$ un sous-$A$-module de type fini de $M$ 
qui l'engendre sur $A[\frac 1 p]$. On appelle {\em topologie $p$-adique} sur $M$ l'unique topologie 
compatible avec sa structure de groupe additif  pour laquelle les sous-groupes $p^nM^\circ$,
pour $n\in \mN$, forment un système fondamental de voisinages de $0$. Cette topologie ne dépend pas
du choix de $M^\circ$. En effet, si $M'$ est un autre sous-$A$-module de type fini de $M$ qui l'engendre sur $A[\frac 1 p]$,
alors il existe $m\geq 0$ tel que $p^mM^\circ\subset M'$ et $p^mM'\subset M^\circ$. 
Il est clair que $M$ est un $A[\frac 1 p]$-module topologique.

\subsection{}\label{higgs1-not33}
Soient $A$ un anneau, $n$ un entier $\geq1$.  
On désigne par $\rW(A)$ (resp. $\rW_n(A)$) l'anneau des vecteurs de Witt 
(resp. vecteurs de Witt de longueur $n$) à coefficients dans $A$ relatif à $p$. 
On a un homomorphisme d'anneaux
\begin{equation}\label{higgs1-not33a}
\Phi_n\colon 
\begin{array}[t]{clcr}
\rW_n(A)&\rightarrow& A,\\ 
(x_1,\dots,x_n)&\mapsto&x_1^{p^{n-1}}+p x_2^{p^{n-2}}+\dots+p^{n-1}x_n.
\end{array}
\end{equation}
appelé $n$-ième composante fantôme. 
On dispose aussi des morphismes de restriction, de décalage et de Frobenius
\begin{eqnarray}
\rR\colon \rW_{n+1}(A)&\rightarrow& \rW_n(A),\label{higgs1-not33b}\\
\rV\colon \rW_n(A)&\rightarrow& \rW_{n+1}(A),\label{higgs1-not33c}\\
\rF\colon \rW_{n+1}(A)&\rightarrow& \rW_n(A).\label{higgs1-not33d}
\end{eqnarray}
Lorsque $A$ est de caractéristique $p$, $\rF$ induit un endomorphisme de $\rW_n(A)$, encore noté $\rF$. 

\subsection{}\label{higgs1-not60}
Pour toute catégorie abélienne $\bA$, on désigne par 
$\bD(\bA)$ sa catégorie dérivée et par $\bD^-(\bA)$, $\bD^+(\bA)$ et $\bD^\rb(\bA)$ les sous-catégories
pleines de $\bD(\bA)$ formées des complexes à cohomologie bornée
supérieurement, inférieurement et des deux côtés, respectivement. 
Sauf mention expresse du contraire, les complexes de $\bA$ sont à différentielles de degré $+1$,
le degré étant écrit en exposant. 

\subsection{}\label{higgs1-not621}
Soit $(X,A)$ un topos annelé. On note $\bMod(A)$ ou $\bMod(A,X)$ 
la catégorie des $A$-modules de $X$. Si $M$ est un $A$-module, on désigne par $\rS_A(M)$ 
(resp. $\wedge_A(M)$, resp. $\Gamma_A(M)$) l'algèbre symétrique (resp. extérieure, resp. à puissances divisées) 
de $M$ (\cite{illusie1} I 4.2.2.6) et pour tout entier $n\geq 0$, par $\rS_A^n(M)$ (resp. $\wedge_A^n(M)$,
resp. $\Gamma_A^n(M)$) sa partie homogène de degré $n$. 
On omettra l'anneau $A$ des notations lorsqu'il n'y a aucun risque d'ambiguïté. 
Les formations de ces algèbres commutent à la localisation au-dessus d'un objet de $X$.

\subsection{}\label{higgs1-koszul1}\index{10210@$\mK_\bullet^A(u)$, $\mK_\bullet^A(u,C)$, $\mK^\bullet_A(u,C)$} 
Soient $A$ un anneau, $L$ un $A$-module, $u\colon L\rightarrow A$ une forme linéaire. 
Pour tout $x\in \wedge(L)$, on note $d_u(x)$ le produit intérieur de $x$ et  $u$
(\cite{alg1-3} III § 11.7 exemple page 161). D'après ({\em loc. cit.} page 162), on a 
\begin{equation}\label{higgs1-koszul1a}
d_u(x_1\wedge\dots\wedge x_n)=\sum_{i=1}^n(-1)^{i+1}u(x_i)x_1\wedge\dots\wedge x_{i-1}\wedge x_{i+1}
\wedge\dots\wedge x_n
\end{equation}
pour $x_1,\dots,x_n\in L$. L'application $d_u\colon \wedge(L)\rightarrow \wedge(L)$ est une anti-dérivation 
de degré $-1$ et de carré nul (\cite{alg1-3} III § 11.8 exemple page 165). L'algèbre $\wedge(L)$ munie de l'anti-dérivation
$d_u $ s'appelle algèbre (ou complexe) de {\em Koszul} de $u$. On la note
$\mK^A_\bullet(u)$; on a donc $\mK_n^A(u)=\wedge^nL$ et les différentielles
de $\mK^A_\bullet(u)$ sont de degré $-1$ (cf. \cite{alg10} § 9.1). 

Pour tout complexe de $A$-modules $C$, on définit le complexe de chaînes (\cite{alg10} § 5.1)
\begin{equation}\label{higgs1-koszul1c}
\mK_\bullet^A(u,C)=\mK^A_\bullet(u)\otimes_AC
\end{equation}
et le complexe de cochaînes  
\begin{equation}\label{higgs1-koszul1d}
\mK^\bullet_A(u,C)=\Homgr_A(\mK^A_\bullet(u),C).
\end{equation}
D'après (\cite{alg10} § 9.1 cor.~2 de prop.~1), si $\Ann(C)$ est l'annulateur de $C$, alors  $u(L)+\Ann(C)$
annule $\rH^*(\mK^\bullet_A(u,C))$ et $\rH_*(\mK_\bullet^A(u,C))$.

Supposons que $L$ soit somme directe de $L_1,\dots,L_r$ et notons $u_i\colon L_i\rightarrow A$ 
la restriction de $u$ à $L_i$. Alors l'isomorphisme canonique (\cite{alg1-3} III § 7.7 prop.~10)
\begin{equation}\label{higgs1-koszul1b}
{^g\otimes}_{1\leq i\leq r} \wedge( L_i)\stackrel{\sim}{\rightarrow}\wedge(L)
\end{equation}
est un isomorphisme de complexes $\otimes_{1\leq i\leq r} \mK^A_\bullet(u_i)\stackrel{\sim}{\rightarrow}\mK^A_\bullet(u)$,  
où le symbole ${^g\otimes}$ désigne le produit tensoriel gauche (cf. \cite{alg1-3} III § 4.7 remarques page 49).

Comme $d_u$ est une anti-dérivation, le produit dans l'algèbre $\wedge(L)$ induit un morphisme de complexes
\begin{equation}\label{higgs1-koszul1e}
\mK_\bullet^A(u)\otimes_A\mK_\bullet^A(u)\rightarrow \mK_\bullet^A(u).
\end{equation}
Supposant $L$ projectif de rang $n$ et composant avec le morphisme canonique 
$\mK_\bullet^A(u)\rightarrow \wedge^nL[-n]$, on en déduit un morphisme de complexes
\begin{equation}\label{higgs1-koszul1f}
\mK_\bullet^A(u)\rightarrow \Homgr_A(\mK_\bullet^A(u), \wedge^nL[-n]),
\end{equation}
qui est bijectif (\cite{alg1-3} III § 7.8 page 87). Pour tout complexe de $A$-modules $C$, 
on en déduit un isomorphisme de complexes 
(\cite{alg10} § 9.1 page 149)
\begin{equation}\label{higgs1-koszul1g}
\mK_\bullet^A(u,C)\stackrel{\sim}{\rightarrow}\mK_A^\bullet(u,C\otimes_A \wedge^nL[-n]).
\end{equation}
Par passage à l'homologie, on a donc, pour tout entier $i$, un isomorphisme 
\begin{equation}\label{higgs1-koszul1h}
\rH_i(\mK_\bullet^A(u,C))\stackrel{\sim}{\rightarrow}\rH^{n-i}(\mK_A^\bullet(u,C\otimes_A \wedge^nL)).
\end{equation}

\subsection{}\label{higgs1-koszul2}
Soient $A$ un anneau, $L$ un $A$-module, 
$u\colon \rS(L)\otimes_AL\rightarrow \rS(L)$ la forme linéaire telle que $u(s\otimes x)=sx$ 
pour tous $s\in \rS(L)$ et $x\in L$. Par l'isomorphisme canonique (\cite{alg1-3} III § 7.5 prop.~8)
\begin{equation}\label{higgs1-koszul2a}
\wedge(\rS(L)\otimes_AL)\stackrel{\sim}{\rightarrow}\rS(L)\otimes_A\wedge(L),
\end{equation}
la différentielle du complexe de Koszul $\mK^{\rS(L)}_\bullet(u)$ \eqref{higgs1-koszul1} est transportée en l'application 
\begin{equation}
d\colon \rS(L)\otimes_A\wedge(L)\rightarrow \rS(L)\otimes_A\wedge(L)
\end{equation} 
définie pour tous $x_1,\dots,x_n,y_1,\dots,y_m\in L$ par 
\begin{equation}\label{higgs1-koszul2b}
d((x_1\dots x_n)\otimes (y_1\wedge\dots\wedge y_m))
=\sum_{i=1}^m(-1)^{i+1}y_ix_1\dots x_n\otimes
(y_1\wedge\dots\wedge y_{i-1}\wedge y_{i+1}\wedge \dots\wedge y_m).
\end{equation}
Pour tout complexe de $\rS(L)$-modules $C$, on pose
\begin{eqnarray}
\mK_\bullet^{\rS(L)}(C)&=&\mK^{\rS(L)}_\bullet(u,C),\label{higgs1-koszul2dd}\\
\mK^\bullet_{\rS(L)}(C)&=&\mK_{\rS(L)}^\bullet(u,C),\label{higgs1-koszul2d}
\end{eqnarray} 
(le morphisme $u$ étant canonique, peut être omis de la notation). 

Soient $L'$ un $A$-module, $u'\colon \rS(L')\otimes_AL'\rightarrow \rS(L')$ la forme linéaire telle que $u'(s'\otimes x')=s'x'$ 
pour $s'\in \rS(L')$ et $x'\in L'$. L'isomorphisme \eqref{higgs1-koszul1b} induit un isomorphisme 
\begin{equation}\label{higgs1-koszul2c}
(\rS(L)\otimes_A\Lambda(L)){^g\boxtimes}_A(\rS(L')\otimes_A\Lambda(L'))
\stackrel{\sim}{\rightarrow}\rS(L\oplus L')\otimes_A\Lambda(L\oplus L'),
\end{equation}
où le produit tensoriel extérieur gauche est pris relativement au diagramme co-cartésien canonique
\[
\xymatrix{
A\ar[r]\ar[d]&{\rS(L)}\ar[d]\\
{\rS(L')}\ar[r]&{\rS(L\oplus L')}}
\]
Il résulte de \eqref{higgs1-koszul2b} que \eqref{higgs1-koszul2c} est un isomorphisme de complexes 
\begin{equation}\label{higgs1-koszul2cc}
\mK^{\rS(L)}_\bullet(u){^g\boxtimes}_A\mK^{\rS(L')}_\bullet(u')
\stackrel{\sim}{\rightarrow}\mK^{\rS(L\oplus L')}_\bullet(u\oplus u').
\end{equation}
On en déduit, pour tout complexe de $\rS(L\oplus L')$-modules $C$, des isomorphismes de complexes 
\begin{eqnarray}
\mK_\bullet^{\rS(L\oplus L')}(C)&\stackrel{\sim}{\rightarrow}&\mK_\bullet^{\rS(L)}(\mK_\bullet^{\rS(L')}(C)),\label{higgs1-koszul2e}\\
\mK^\bullet_{\rS(L\oplus L')}(C)&\stackrel{\sim}{\rightarrow}&\mK_{\rS(L)}^\bullet(\mK_{\rS(L')}^\bullet(C)).\label{higgs1-koszul2f}
\end{eqnarray}

\begin{defi}\label{higgs1-not6}\index{Module de Higgs}\index{Champ de Higgs}\index{10220@$\bMH(A,E)$ ($A$ un anneau, $E$ un $A$-module)}
Soient $(X,A)$ un topos annelé, $E$ un $A$-module. 

{\rm (i)}\ On appelle {\em $A$-module de Higgs à coefficients dans $E$}
un couple $(M,\theta)$ formé d'un $A$-module $M$ et d'un morphisme $A$-linéaire 
\begin{equation}\label{higgs1-not6a}
\theta\colon M\rightarrow M\otimes_AE
\end{equation}
tel que $\theta\wedge \theta=0$. On dit alors que $\theta$ est un {\em $A$-champ de Higgs} sur $M$ 
à coefficients dans $E$. 

{\rm (ii)}\ Si $(M_1,\theta_1)$ et $(M_2,\theta_2)$ sont deux $A$-modules de Higgs,
un morphisme de $(M_1,\theta_1)$ dans $(M_2,\theta_2)$ est un morphisme $A$-linéaire 
$u\colon M_1\rightarrow M_2$ tel que $(u\otimes\id_E)\circ \theta_1=\theta_2\circ u$. 
\end{defi}

Les $A$-modules de Higgs à coefficients dans $E$ forment une catégorie que l'on note 
$\bMH(A,E)$. On peut compléter la terminologie et faire les remarques suivantes.

\addtocounter{subsubsection}{1}
\addtocounter{equation}{1}

\subsubsection{}\label{higgs1-not61}\index{Complexe de Dolbeault d'un module de Higgs}
\index{10265@$\mK^\bullet(M,\theta)$ ($(M,\theta)$ module de Higgs)}
Soit $(M,\theta)$ un $A$-module de Higgs à coefficients dans $E$. Pour tout $i\geq 1$, on désigne par
\begin{equation}\label{higgs1-not61a}
\theta_i\colon M\otimes_A \wedge^iE \rightarrow M\otimes_A \wedge^{i+1}E
\end{equation}
le morphisme $A$-linéaire défini pour toutes sections locales 
$m$ de $M$ et $\omega$ de $\wedge^iE$ par $\theta_i(m\otimes \omega)=\theta(m)\wedge \omega$.
On a $\theta_{i+1}\circ \theta_i=0$. Suivant Simpson (\cite{simpson4} page 24),
on appelle complexe de {\em Dolbeault} de $(M,\theta)$
et l'on note $\mK^\bullet(M,\theta)$ le complexe de cochaînes de $A$-modules 
\begin{equation}\label{higgs1-not6b}
M\stackrel{\theta}{\longrightarrow}M\otimes_AE\stackrel{\theta_1}{\longrightarrow} M\otimes_A\wedge^2E \dots,
\end{equation}
où $M$ est placé en degré $0$ et les différentielles sont de degré $1$. 

\addtocounter{subsubsection}{2}
\addtocounter{equation}{1}

\subsubsection{} \label{higgs1-not62}\index{Invariants caracteristiques d'un champ de Higgs@Invariants caractéristiques d'un champ de Higgs}
Soit $(M,\theta)$ un $A$-module de Higgs à coefficients dans $E$ tel que $M$ soit un $A$-module
localement libre de type fini. Considérons, pour un entier $i\geq 1$, le morphisme composé 
\begin{equation}\label{higgs1-not6c}
\xymatrix{
{\wedge^iM}\ar[r]^-(0.5){\wedge^i\theta}&{\wedge^i(M\otimes_AE)}\ar[r]&
{\wedge^iM\otimes_A\rS^iE}},
\end{equation}
où la seconde flèche est le morphisme canonique (\cite{illusie1} V 4.5).  
On appelle {\em $i$-ième invariant caractéristique} de $\theta$ et l'on note $\lambda_i(\theta)$
la trace du morphisme \eqref{higgs1-not6c} vue comme section de 
$\Gamma(X,\rS^iE)$.

\addtocounter{subsubsection}{1}

\subsubsection{} \label{higgs1-not63}\index{Produit tensoriel!1@de modules de Higgs}
\addtocounter{equation}{1}
Soient $(M_1,\theta_1),(M_2,\theta_2)$ deux $A$-modules de Higgs à coefficients dans $E$. 
On appelle champ de Higgs {\em total} sur $M_1\otimes_AM_2$  
le morphisme  $A$-linéaire 
\begin{equation}\label{higgs1-not6ab}
\theta_\tot\colon M_1\otimes_AM_2\rightarrow M_1\otimes_AM_2\otimes_AE
\end{equation} 
défini par 
\begin{equation}
\theta_\tot=\theta_1\otimes \id_{M_2}+\id_{M_1}\otimes \theta_2.
\end{equation}
On dit que $(M_1\otimes_AM_2,\theta_\tot)$ est le {\em produit tensoriel} de $(M_1,\theta_1)$ et $(M_2,\theta_2)$.

\addtocounter{subsubsection}{2}
\addtocounter{equation}{1}

\subsubsection{}\label{higgs1-not65}
Supposons $E$ localement libre de type fini sur $A$ et posons $F=\cHom_A(E,A)$. 
Pour tout $A$-module $M$, le morphisme canonique 
\begin{equation}
\cEnd_A(M)\otimes_AE\rightarrow \cHom_A(M,M\otimes_AE)
\end{equation}
étant un isomorphisme, la donnée d'un $A$-champ de Higgs $\theta$ sur $M$ 
est équivalente à la donnée d'une structure de $\rS(F)$-module sur $M$
compatible avec sa structure de $A$-module. 
D'autre part, en vertu de (\cite{alg1-3} § 11.5 prop.~7), le $A$-module $\cHom_A(\wedge(F),A)$ 
s'identifie à l'algèbre duale graduée de $\wedge(F)$, et on a un isomorphisme canonique d'algèbres graduées
\begin{equation}
\wedge(E)\rightarrow \cHom_A(\wedge(F),A). 
\end{equation}
On vérifie que celui-ci induit un isomorphisme de complexes de $A$-modules \eqref{higgs1-koszul2d}
\begin{equation}
\mK^\bullet(M,\theta)\stackrel{\sim}{\rightarrow} \mK_{\rS(F)}^\bullet(M).
\end{equation}

\subsection{}\label{higgs1-not149}
Soient $f\colon (X',A')\rightarrow (X,A)$ un morphisme de topos annelés,
$E$ un $A$-module, $E'$ un $A'$-module, $\gamma\colon f^*(E)\rightarrow E'$ un morphisme $A'$-linéaire,
$(M,\theta)$ un $A$-module de Higgs à coefficients dans $E$. Le morphisme composé 
\begin{equation}
\theta'\colon \xymatrix{
{f^*(M)}\ar[r]^-(0.4){f^*(\theta)}&{f^*(M)\otimes_{A'}f^*(E)}\ar[r]^-(0.4){\id\otimes \gamma}&{f^*(M)\otimes_{A'}E'}}
\end{equation}
est alors un $A'$-champ de Higgs à coefficients dans $E'$. On dit que le $A'$-module
de Higgs $(f^*(M),\theta')$ est {\em l'image inverse} de $(M,\theta)$ par $(f,\gamma)$.

\subsection{}\label{higgs1-not15}\index{connexion integrable@$\lambda$-connexion (intégrable)}
Soient $(X,A)$ un topos annelé, $B$ une $A$-algèbre, 
$\Omega^1_{B/A}$ le $B$-module des différentielles de Kähler de $B$ sur $A$ (\cite{illusie1} II 1.1.2), 
\begin{equation}
\Omega_{B/A}=\oplus_{n\in \mN}\Omega^n_{B/A}=\wedge_B(\Omega^1_{B/A})
\end{equation} 
l'algèbre extérieure de $\Omega^1_{B/A}$. Il existe alors une et une unique $A$-anti-dérivation 
$d\colon \Omega_{B/A}\rightarrow \Omega_{B/A}$ de degré $1$, de carré nul, qui prolonge
la $A$-dérivation universelle $d\colon B\rightarrow \Omega^1_{B/A}$. Cela résulte par exemple de 
(\cite{alg10} §2.10 prop.~13) en remarquant que $\Omega^1_{B/A}$ est le faisceau associé au préfaisceau 
$U\mapsto \Omega^1_{B(U)/A(U)}$ $(U\in \ob(X))$.

Soient $M$ un $B$-module, $\lambda\in \Gamma(X,A)$. 
Une {\em $\lambda$-connexion} sur $M$ relativement à l'extension $B/A$
est la donnée d'un morphisme $A$-linéaire 
\begin{equation}
\nabla\colon M\rightarrow \Omega^1_{B/A}\otimes_BM
\end{equation}
tel que pour toutes sections locales $x$ de $B$ et $s$ de $M$, on ait 
\begin{equation}\label{higgs1-not15a}
\nabla(xs)=\lambda d(x)\otimes s+x\nabla(s).
\end{equation} 
On dit aussi que $(M,\nabla)$ est un $B$-module à $\lambda$-connexion relativement à l'extension $B/A$.
On omettra l'extension $B/A$ de la terminologie lorsqu'il n'y a aucun risque de confusion. 
Le morphisme $\nabla$ se prolonge en un unique morphisme $A$-linéaire gradué de degré $1$ que l'on note aussi 
\begin{equation}\label{higgs1-not15b}
\nabla\colon \Omega_{B/A}\otimes_BM\rightarrow   \Omega_{B/A}\otimes_BM,
\end{equation} 
tel que pour toutes sections locales $\omega$ de  $\Omega^i_{B/A}$ et $s$ de $\Omega^j_{B/A}\otimes_BM$  
($i,j\in \mN$), on ait 
\begin{equation}\label{higgs1-not15c}
\nabla(\omega\wedge s)=\lambda  d(\omega)\wedge s+(-1)^i\omega \wedge \nabla(s).
\end{equation} 
Itérons cette formule~: 
\begin{equation}\label{higgs1-not15d}
\nabla \circ \nabla(\omega\wedge s)= \omega \wedge \nabla\circ \nabla(s).
\end{equation} 
On dit que $\nabla$ est {\em intégrable} si $\nabla \circ \nabla=0$.

Soient $(M,\nabla)$, $(M',\nabla')$ deux modules à $\lambda$-connexions. 
Un morphisme de $(M,\nabla)$ dans $(M',\nabla')$
est la donnée d'un morphisme $B$-linéaire $u\colon M\rightarrow M'$ 
tel que  $(\id \otimes u)\circ \nabla=\nabla'\circ u$. 

Les $1$-connexions sont classiquement appelées {\em connexions}. 
Les $0$-connexions intégrables sont les $B$-champs de Higgs à coefficients dans $\Omega^1_{B/A}$ \eqref{higgs1-not6}.

\subsection{}\label{higgs1-not18}
Soient $f\colon (X',A')\rightarrow (X,A)$ un morphisme de topos annelés, $B$ une $A$-algèbre, 
$B'$ une $A'$-algèbre, $\alpha\colon f^*(B)\rightarrow B'$ un homomorphisme de $A'$-algèbres,
$\lambda\in \Gamma(X,A)$, $(M,\nabla)$ un module à $\lambda$-connexion relativement à l'extension $B/A$. 
Notons $\lambda'$ l'image canonique de $\lambda$ dans $\Gamma(X',A')$, 
$d'\colon B'\rightarrow \Omega^1_{B'/A'}$ la $A'$-dérivation universelle de $B'$ et 
\begin{equation}
\gamma\colon f^*(\Omega^1_{B/A}) \rightarrow \Omega^1_{B'/A'}
\end{equation}
le morphisme $\alpha$-linéaire canonique. On voit aussitôt que 
$f^*(\nabla)$ est une $\lambda'$-connexion sur $f^*(M)$ relativement à l'extension $f^*(B)/A'$,
qui est intégrable si $\nabla$ l'est.
Par ailleurs, il existe un unique morphisme $A'$-linéaire 
\begin{equation}
\nabla'\colon B'\otimes_{f^*(B)}f^*(M)\rightarrow  \Omega^1_{B'/A'}\otimes_{f^*(B)}f^*(M)
\end{equation}
tel que pour toutes sections locales $x'$ de $B'$ et $t$ de $f^*(M)$, on ait 
\begin{equation}
\nabla'(x'\otimes t)=\lambda' d'(x')\otimes t+x'(\gamma\otimes\id)(f^*(\nabla)(t)).
\end{equation}
C'est une $\lambda'$-connexion sur $B'\otimes_{f^*(B)}f^*(M)$ relativement à l'extension 
$B'/A'$, qui est intégrable si $\nabla$ l'est. 

\subsection{}\label{higgs1-not16}
Soient $(X,A)$ un topos annelé, $B$ une $A$-algèbre, $\lambda\in \Gamma(X,A)$, $(M,\nabla)$ un
module à $\lambda$-connexion relativement à l'extension $B/A$. 
Supposons qu'il existe un $A$-module $E$
et un isomorphisme $B$-linéaire $\gamma\colon E\otimes_AB\stackrel{\sim}{\rightarrow}\Omega^1_{B/A}$ 
tels que pour toute section locale $\omega$ de $E$, on ait $d(\gamma(\omega\otimes 1))=0$.
Notons $\vartheta\colon M\rightarrow E\otimes_AM$  le morphisme induit par $\nabla$ et $\gamma$.
Alors pour que la $\lambda$-connexion $\nabla$ soit intégrable, il faut et il suffit que 
$\vartheta$ soit un $A$-champ de Higgs sur $M$ à coefficients dans $E$. En effet, le diagramme 
\begin{equation}
\xymatrix{
{E\otimes_AM}\ar[r]^-(0.5){-\id\wedge \vartheta}\ar[d]&{(\wedge^2E)\otimes_AM}\ar[d]\\
{\Omega^1_{B/A}\otimes_BM}\ar[r]^-(0.4)\nabla&{\Omega^2_{B/A}\otimes_BM}}
\end{equation}
où les flèches verticales sont les isomorphismes induits par $\gamma$, est clairement commutatif.

\subsection{}\label{higgs1-not17}\index{Produit tensoriel!2@d'une $\lambda$-connexion et d'un module de Higgs}
Soient $(X,A)$ un topos annelé, $B$ une $A$-algèbre, $\lambda\in \Gamma(X,A)$, $(M,\nabla)$ un module 
à $\lambda$-connexion intégrable relativement à l'extension $B/A$. Supposons qu'il existe un $A$-module $E$
et un $B$-isomorphisme $\gamma\colon E\otimes_AB\stackrel{\sim}{\rightarrow}\Omega^1_{B/A}$ 
tels que pour toute section locale $\omega$ de $E$, on ait $d(\gamma(\omega\otimes 1))=0$ (cf. \ref{higgs1-not16}).
Soit $(N,\theta)$ un $A$-module de Higgs à coefficients dans $E$.  
Il existe un unique morphisme $A$-linéaire 
\begin{equation}\label{higgs1-not17a}
\nabla'\colon M\otimes_AN\rightarrow  \Omega^1_{B/A}\otimes_BM\otimes_AN
\end{equation}
tel que pour toutes sections locales $x$ de $M$ et $y$ de $N$, on ait 
\begin{equation}\label{higgs1-not17b}
\nabla'(x\otimes y)=\nabla(x)\otimes_A y+ (\gamma\otimes_B\id_{M\otimes_AN})(x\otimes_A \theta(y)).
\end{equation}
C'est une $\lambda$-connexion intégrable sur $M\otimes_AN$ relativement à l'extension $B/A$.

\subsection{}\label{higgs1-not20}
Soient $A$ un anneau adique, $I$ un idéal de définition de $A$, $B$ une $A$-algèbre adique, 
{\em i.e.}, une $A$-algèbre $B$ complète et séparée pour la topologie $(IB)$-adique. On rappelle que la topologie 
canonique du $B$-module  $\Omega^1_{B/A}$ est déduite de celle de $B$ (\cite{ega4} 0.20.4.5). 
On désigne par $\hOmega^1_{B/A}$ son séparé complété et on note
\begin{equation}\label{higgs1-not20a}
d\colon B\rightarrow \hOmega^1_{B/A}
\end{equation}
la $A$-dérivation continue universelle de $B$. Soient  $M$ un $B$-module  
complet et séparé pour la topologie $(IB)$-adique, $\lambda\in A$. Une {\em $\lambda$-connexion adique} 
(ou {\em $I$-adique}) sur $M$ relativement à l'extension $B/A$ est la donnée d'un morphisme $A$-linéaire 
\begin{equation}\label{higgs1-not20b}
\nabla\colon M\rightarrow \hOmega^1_{B/A}\hotimes_BM
\end{equation}
tel que pour tous $x\in B$ et $t\in M$, on ait 
\begin{equation}\label{higgs1-not20c}
\nabla(xt)=\lambda d(x)\hotimes t+x\nabla(t).
\end{equation} 
On dit aussi que $(M,\nabla)$ est un $B$-module à $\lambda$-connexion adique (ou $I$-adique) relativement à l'extension $B/A$.
Pour tout entier $n\geq 1$, on pose $A_n=A/I^n$, $B_n=B\otimes_AA_n$ et $M_n=M\otimes_AA_n$ et 
on note $\olambda_n$ la classe de $\lambda$ dans $A_n$. 
Alors, $\nabla$ induit une $\olambda_n$-connexion (usuelle) relativement à l'extension $B_n/A_n$,
\begin{equation}\label{higgs1-not20d}
\nabla_n\colon M_n\rightarrow \Omega^1_{B_n/A_n}\otimes_{B_n}M_n.
\end{equation}
De plus, $\nabla$ s'identifie à la limite projective des morphismes $\nabla_n$. 
On dit que $\nabla$ est {\em intégrable} si $\nabla_n$ est intégrable pour tout $n\geq 1$. 

Soient $(M,\nabla)$, $(M',\nabla')$ deux modules à $\lambda$-connexions adiques. 
Un morphisme de $(M,\nabla)$ dans $(M',\nabla')$
est la donnée d'un morphisme $B$-linéaire $u\colon M\rightarrow M'$ 
tel que  $(\id \hotimes u)\circ \nabla=\nabla'\circ u$.

\subsection{}\label{higgs1-not21}
Soient $A$ un anneau adique, $\lambda\in A$, $B$ une $A$-algèbre adique,
$B'$ une $B$-algèbre adique, $(M,\nabla)$ un $B$-module à $\lambda$-connexion adique 
relativement à l'extension $B/A$. On désigne par 
$d'\colon B'\rightarrow \hOmega^1_{B'/A}$ la $A$-dérivation continue universelle de $B'$ et par
\begin{equation}\label{higgs1-not21a}
\gamma\colon \hOmega^1_{B/A}\rightarrow \hOmega^1_{B'/A}
\end{equation}
le morphisme canonique. Il existe un unique morphisme $A$-linéaire 
\begin{equation}\label{higgs1-not21b}
\nabla'\colon B'\hotimes_{B}M\rightarrow  \hOmega^1_{B'/A}\hotimes_{B}M
\end{equation}
tel que pour tous $x'\in B'$ et $t\in M$, on ait 
\begin{equation}\label{higgs1-not21c}
\nabla'(x'\hotimes t)=\lambda d'(x')\hotimes t+x'(\gamma\hotimes\id_M)(\nabla(t)).
\end{equation}
C'est une $\lambda$-connexion adique sur $B'\hotimes_{B}M$ relativement à l'extension 
$B'/A$, qui est intégrable si $\nabla$ l'est. 

\subsection{}\label{higgs1-not22}
Soient $A$ un anneau adique, $I$ un idéal de définition de $A$, $\lambda\in A$, $B$ une $A$-algèbre adique, 
$(M,\nabla)$ un $B$-module à $\lambda$-connexion adique relativement à l'extension $B/A$. 
Pour tout entier $n\geq 1$, posons $A_n=A/I^n$ et $B_n=B\otimes_AA_n$. 
Supposons les conditions suivantes remplies~:
\begin{itemize}
\item[(i)] $I$ est un idéal de type fini de $A$, et  $\Omega^1_{B_1/A_1}$ est un $B_1$-module de type fini.
\item[(ii)] Il existe un $A$-module {\em libre de type fini} $E$
et un isomorphisme $B$-linéaire $\gamma\colon E\otimes_AB\stackrel{\sim}{\rightarrow}\hOmega^1_{B/A}$ 
tels que $\gamma(E)\subset d(B)$. 
\end{itemize}
On observera que la topologie de $\hOmega^1_{B/A}$ est la topologie $(IB)$-adique et qu'on a 
$\hOmega^1_{B/A}\otimes_AA_n=\Omega^1_{B_n/A_n}$ pour tout $n\geq 1$, 
d'après (\cite{ac} chap.~III § 2.12 cor.~1 de prop.~14). 
Par ailleurs, on a $\hOmega^1_{B/A}\hotimes_BM=\hOmega^1_{B/A}\otimes_BM=E\otimes_AM$.
On note $\vartheta\colon M\rightarrow E\otimes_AM$  
le morphisme induit par $\nabla$ et $\gamma$. 
Alors, pour que la $\lambda$-connexion adique $\nabla$ soit intégrable, 
il faut et il suffit que $\vartheta$ soit un $A$-champ de Higgs sur $M$ à coefficients dans $E$. 
En effet, pour que $\vartheta$ soit un $A$-champ de Higgs sur $M$ à coefficients dans $E$, il faut et il suffit que
pour tout $n\geq 1$, $\vartheta\otimes_AA_n$ soit un $A_n$-champ de Higgs sur $M\otimes_AA_n$ à coefficients dans 
$E\otimes_AA_n$. L'assertion recherchée résulte alors de \ref{higgs1-not16}. 

\subsection{}\label{higgs1-not23}
Soient $A$ un anneau adique, $\lambda\in A$, $B$ une $A$-algèbre adique, 
$(M,\nabla)$ un $B$-module à $\lambda$-connexion adique intégrable relativement à l'extension $B/A$. 
Supposons les conditions (i) et (ii) de \ref{higgs1-not22} remplies. 
Soit $(N,\theta)$ un $A$-module de Higgs à coefficients dans $E$; on munit 
$N$ de la topologie induite par celle de $A$. Par passage à la limite de \ref{higgs1-not17}, 
il existe un unique morphisme $A$-linéaire 
\begin{equation}\label{higgs1-not23a}
\nabla'\colon M\hotimes_AN\rightarrow  \hOmega^1_{B/A}\hotimes_BM\hotimes_AN
\end{equation}
tel que pour tous $x\in M$ et $y\in N$, on ait 
\begin{equation}\label{higgs1-not23b}
\nabla'(x\hotimes y)=\nabla(x)\hotimes_A y+ (\gamma\otimes_B\id_{M\hotimes_A N})(x\hotimes_A \theta(y)).
\end{equation}
C'est une $\lambda$-connexion adique intégrable sur $M\hotimes_AN$ relativement à l'extension $B/A$.

\subsection{}
Soit $X$ une variété projective lisse complexe. 
Un fibré harmonique sur $X$ est la donnée d'un triplet $(M,D,\langle\ ,\ \rangle)$,
où $M$ est un fibré vectoriel complexe de classe $\cC^\infty$ sur $X$, $D$ est une connexion intégrable sur $M$
et $\langle\ ,\ \rangle$ est une métrique hermitienne sur $M$, satisfaisant à une condition 
décrite plus bas (\cite{lepotier} §1).
On peut écrire de manière unique $D=\nabla+\alpha$, où $\nabla$ est une connexion 
hermitienne et $\alpha$ est une forme différentielle de  degré $1$, à valeurs dans $\End_{\co_X}(M)$,
auto-adjointe par rapport à $\langle\ ,\ \rangle$. On décompose $\nabla$ et $\alpha$
selon leurs types 
\begin{equation}
\nabla=\partial+\overline{\partial}, \ \ \alpha=\theta+\theta^*,
\end{equation}
où $\partial$ et $\overline{\partial}$ sont de type $(1,0)$ et $(0,1)$, respectivement, et 
$\theta$ est une forme différentielle de type $(1,0)$ à valeurs dans $\End_{\co_X}(M)$. 
La condition requise est que l'opérateur $D''=\overline{\partial}+\theta$ est intégrable,
{\em i.e.}, $D''^2=0$. Cette condition équivaut à dire que $\overline{\partial}^2=0$, 
$\overline{\partial}\theta=0$ et $\theta\wedge \theta=0$. Ainsi l'opérateur 
$\overline{\partial}$ définit sur $M$ une structure de fibré vectoriel holomorphe,
et $\theta$ est alors un champ de Higgs sur $M$ à coefficients dans $\Omega^1_{X/\mC}$.
Le complexe de Dolbeault de $(M,\theta)$
\begin{equation}
0\rightarrow \rA^0(M)\stackrel{D''}{\longrightarrow} \rA^1(M)\stackrel{D''}{\longrightarrow} \dots
\end{equation}
s'obtient en prolongeant $D''$ aux formes différentielles de classe $\cC^\infty$.
Il induit par restriction aux formes différentielles holomorphes le complexe $\mK(M,\theta)$;
d'où la terminologie.

\section{Rappels sur la cohomologie continue des groupes profinis}\label{higgs1-kun}

\subsection{}\label{higgs1-not55}\index{10300@$\bRep_A(G)$ , $\bRep_A^\cont(G)$, $\bRep_A^\disc(G)$}
Soient $G$ un groupe profini, $A$ un anneau topologique muni d'une action continue de $G$ 
par des homomorphismes d'anneaux. Une {\em $A$-représentation} de $G$ est la donnée d'un $A$-module 
$M$ et d'une action $A$-semi-linéaire de $G$ sur $M$, {\em i.e.},  telle que 
pour tous $g\in G$, $a\in A$ et $m\in M$, on ait $g(am)=g(a)g(m)$.
On dit que la $A$-représentation est {\em continue} si $M$ est un $A$-module topologique 
et si l'action de $G$ sur $M$ est continue. Soient $M$, $N$ deux $A$-représentations 
(resp. deux $A$-représentations continues) de $G$. 
Un morphisme de $M$ dans $N$ est la donnée d'un morphisme $A$-linéaire et $G$-équivariant 
(resp. $A$-linéaire, continu et $G$-équivariant) de $M$ dans $N$. 
On note $\bRep_A(G)$ (resp. $\bRep_A^\cont(G)$)
la catégorie des $A$-représentations (resp. $A$-représentations continues) de $G$.
Si $M$ et $N$ sont deux $A$-représentations de $G$, 
les $A$-modules $M\otimes_AN$ et $\Hom_A(M,N)$ sont naturellement des $A$-représentations de $G$. 

Supposons l'action de $G$ sur $A$ triviale. 
Les objets de $\bRep_A^\cont(G)$ sont alors aussi appelés des {\em $A$-$G$-modules topologiques}. 
Un $A$-$G$-module topologique dont la topologie est discrète est appelé {\em $A$-$G$-module discret}. 
On désigne par $\bRep_A^\disc(G)$ la sous-catégorie pleine de $\bRep_A^\cont(G)$ 
formée des $A$-$G$-modules discrets. 

Si $R$ est un anneau sans topologie, il sera sous-entendu que les 
$R$-$G$-modules topologiques sont définis relativement à la topologie discrète de $R$ (et l'action triviale de $G$ sur $R$).

\subsection{}\label{higgs1-not59}
Soient $A$ un anneau, $G$ un groupe profini, $M$ un $A$-module muni de la topologie discrète. 
Le {\em $A$-$G$-module induit} de $M$,
noté $\Ind_{A,G}(M)$, est le $A$-module des applications continues de $G$ dans $M$, 
muni de l'action de $G$ définie pour tous $f\in \Ind_{A,G}(M)$ et $g\in G$, par 
\begin{equation}\label{higgs1-not59a}
(g\cdot f)(x)=f(x\cdot g). 
\end{equation}
C'est un $A$-$G$-module discret. On définit ainsi un foncteur exact
\begin{equation}\label{higgs1-not59b}
\Ind_{A,G}\colon \bMod(A)\rightarrow \bRep_A^\disc(G), \ \ \ M\mapsto \Ind_{A,G}(M),
\end{equation}
qui est un adjoint à droite du foncteur d'oubli de l'action de $G$ (\cite{tsuji2} 11.1). 
Il transforme donc les $A$-modules injectifs en $A$-$G$-modules injectifs.
La catégorie $\bRep_A^\disc(G)$ a suffisamment d'injectifs. 
Pour qu'un objet de $\bRep_A^\disc(G)$ soit injectif, il faut et il suffit qu'il soit un facteur
direct d'un objet de la forme $\Ind_{A,G}(I)$, où $I$ est un $A$-module injectif (\cite{tsuji2} 11.2).
Un $A$-$G$-module discret est dit {\em induit} s'il est isomorphe au $A$-$G$-module induit d'un $A$-module. 

On désigne par 
$\Gamma (G,-)$ le foncteur exact à gauche 
\begin{equation}\label{higgs1-not592b}
\Gamma (G,-)\colon  \bRep_A^\disc(G)\rightarrow \bMod(A), \ \ \ M\mapsto M^G, 
\end{equation} 
et par
\begin{eqnarray}
\rR\Gamma (G,-)\colon  \bD^+(\bRep_A^\disc(G))&\rightarrow& \bD^+(\bMod(A)),\label{higgs1-not592a}\\
\rH^q (G,-)\colon  \bRep_A^\disc(G)&\rightarrow& \bMod(A), \ \ \ (q\geq 0), \label{higgs1-not592ab}
\end{eqnarray}
ses foncteurs dérivés droits (cf. \ref{higgs1-not60}).

\subsection{}\label{higgs1-not593}
Soient $A$ un anneau, $G$ un groupe profini, $H$ un sous-groupe fermé et distingué de $G$. 
Les groupes $H$ et $G/H$ sont alors profinis. On désigne encore par $\Gamma(H,-)$ le foncteur exact à gauche
\begin{equation}\label{higgs1-not593a}
\Gamma(H,-)\colon \bRep_A^\disc(G) \rightarrow \bRep_A^\disc(G/H), \ \ \ M\mapsto M^H,
\end{equation}
et par $\rR\Gamma(H,-)$ et $\rH^q(H,-)$ $(q\geq 0)$ ses foncteurs dérivés droits. 
Cet abus de notation est justifié par le fait que le diagramme 
\begin{equation}\label{higgs1-not593b}
\xymatrix{
{\bD^+(\bRep_A^\disc(G))}\ar[rr]^-(0.5){\rR\Gamma(H,-)}\ar[d]&&{\bD^+(\bRep_A^\disc(G/H))}\ar[d]\\
{\bD^+(\bRep_A^\disc(H))}\ar[rr]^-(0.5){\rR\Gamma(H,-)}&&{\bD^+(\bMod(A))}}
\end{equation}
où les flèches verticales sont induites par les foncteurs d'oubli, est commutatif à isomorphisme canonique près
(\cite{tsuji2} 11.5). 

\begin{prop}[\cite{tsuji2} 11.7]\label{higgs1-not594}
Soient $A$ un anneau, $G$ un groupe profini, $H$ un sous-groupe fermé et distingué de $G$, 
$M$ un $A$-$G$-module discret. On a alors un isomorphisme canonique fonctoriel de $\bD^+(\bMod(A))$
\begin{equation}\label{higgs1-not594a}
\rR\Gamma(G/H,\rR\Gamma(H,M))\stackrel{\sim}{\rightarrow}\rR\Gamma(G,M).
\end{equation}
\end{prop}

\begin{rema}\label{higgs1-not595}
Soient $A$ un anneau, $G$ un groupe profini, $H$ un sous-groupe fermé et distingué de $G$, 
$M$ un $A$-$G$-module discret. Alors pour tout entier $n\geq 0$, le morphisme de restriction 
\begin{equation}\label{higgs1-not595b}
\rho_n\colon \rH^n(G,M) \rightarrow \rH^n(H,M)
\end{equation}
coïncide avec le composé 
\begin{equation}\label{higgs1-not595a}
u_n\colon \rH^n(G,M)\stackrel{\sim}{\rightarrow} \rH^n(G/H,\tau_{\leq n}\rR\Gamma(H,M))\rightarrow 
\Gamma(G/H,\rH^n(H,M))\rightarrow \rH^n(H,M),
\end{equation}
où $\tau_{\leq n} \rR\Gamma(H,M)$ est la filtration canonique de $\rR\Gamma(H,M)$ (\cite{deligne} 1.4.6), 
la première flèche est induite par l'isomorphisme \eqref{higgs1-not594a} et les autres flèches sont les 
morphismes canoniques. 
En effet, $(\rho_n)$ et $(u_n)$ sont deux morphismes de $\partial$-foncteurs universels qui coïncident en degré $0$. 
\end{rema}

\subsection{}\label{higgs1-not580}
On désigne par $\fD$ la catégorie dont les objets sont les ensembles ordonnés 
$[n]=\{0,\dots,n\}$ (pour $n\in \mN$) et les morphismes sont les applications croissantes.
On réserve la notation $\Delta$ pour un groupe de Galois qui apparaîtra dans \eqref{higgs1-gal2}.  
Pour tous $n\in \mN$ et $i\in [n]$, on note $d^i_n\colon [n-1]\rightarrow [n]$ l'injection croissante qui oublie $i$,
et $s^i_n\colon [n+1]\rightarrow [n]$ la surjection croissante qui répète $i$. On omet 
les indices $n$ dans $d^i_n$ et $s^i_n$ lorsqu'il n'y a aucun risque de confusion.

Soient $\cA$ une catégorie additive, $q$ un entier $\geq 1$, $X$ un objet $q$-cosimplicial de $\cA$ (\cite{illusie1} I 1.1).
On appelle {\em sous-objet diagonal} de $X$, et l'on note $\Delta X$, 
l'objet cosimplicial de $\cA$ défini par $[n]\mapsto X^{n,\dots,n}$. 
On appelle {\em complexe de cochaînes} de $X$ le complexe de cochaînes $q$-uple $\tX$ défini 
pour tout $(n_1,\dots,n_q)\in \mN^q$, par $\tX^{n_1,\dots,n_q}=X^{n_1,\dots,n_q}$, et pour $1\leq i\leq q$,
par la différentielle 
\begin{equation}\label{higgs1-not580a}
d^i=\sum_j(-1)^jX(\id,\dots,d^j,\dots,\id),
\end{equation}
où $d^j$ est placé à la $i$-ième place (\cite{illusie1} I 1.2.2). 
Pour tout complexe de cochaînes $q$-uple $M$ de $\cA$,
on appelle {\em complexe simple associé}, le complexe de cochaînes $\int M$ défini pour tout $n\in \mN$, par 
\begin{equation}\label{higgs1-not580b}
(\int M)^n=\oplus_{\sum_{i=1}^q n_i=n}M^{n_1,\dots,n_q},
\end{equation}
la différentielle $d$ étant donnée, en restriction à $M^{n_1,\dots,n_q}$, par 
$\sum_j(-1)^{\sum_{i<j}n_i}d^j$ (\cite{illusie1} I 1.2.1).

D'après le théorème d'Eilenberg-Zilber-Cartier (\cite{illusie1} I 1.2.2), pour tout objet $2$-cosimplicial $X$ de $\cA$,
il existe deux homomorphismes fonctoriels, la flèche ``d'Alexander-Whitney''
\begin{equation}\label{higgs1-not580c}
\int \tX\rightarrow (\Delta X)^\sim
\end{equation}
et le ``shuffle map''
\begin{equation}\label{higgs1-not580d}
(\Delta X)^\sim\rightarrow \int \tX,
\end{equation}
induisant l'identité en degré $0$, et inverses l'un de l'autre à homotopie fonctorielle près. On prendra garde
que le sens des flèches est inversé par rapport à {\em loc. cit.} qui traite des objets simpliciaux.
Plus précisément, 
la flèche \eqref{higgs1-not580c} envoie $X^{m,n}$ dans $X^{m+n,m+n}$ par $X(d^{m+n}\dots d^{m+1},d^0\dots d^0)$,
où $d^0$ est composé $m$ fois, et la flèche  \eqref{higgs1-not580d} envoie $X^{m+n,m+n}$ dans $X^{m,n}$ par
\begin{equation}\label{higgs1-not580e}
\sum_{(\mu,\nu)}\varepsilon(\mu,\nu)X(s^{\nu_1}\dots s^{\nu_n},s^{\mu_1}\dots s^{\mu_m}),
\end{equation}
la somme étant prise sur tous les $(m,n)$-shuffles de l'ensemble $[m+n-1]$, 
où $\varepsilon(\mu,\nu)$ est le signe du shuffle.

\subsection{}\label{higgs1-not584}
Soient $A$ un anneau, $X$, $Y$ deux $A$-modules cosimpliciaux. On désigne par $X\otimes_AY$ le $A$-module 
$2$-cosimplicial défini par $(m,n)\mapsto X^m\otimes_AY^n$, et par $X\otimes^\Delta_AY$ le sous-objet 
diagonal $\Delta(X\otimes_AY)$. On notera que le complexe simple associé au complexe de cochaînes de $X\otimes_AY$
n'est autre que le produit tensoriel $\tX\otimes_A\tY$ des complexes de $A$-modules $\tX$ et $\tY$ (cf. \cite{alg10} § 4.1).
On a donc la flèche ``d'Alexander-Whitney'' 
\begin{equation}\label{higgs1-not584a}
\tX\otimes_A\tY \rightarrow (X\otimes^\Delta_AY)^\sim,
\end{equation} 
et le ``shuffle map''   
\begin{equation}\label{higgs1-not584b}
(X\otimes^\Delta_AY)^\sim\rightarrow \tX\otimes_A\tY.
\end{equation}

\subsection{}\label{higgs1-not78}\index{10380@$\rK^\bullet(G,M)$, $\trK^\bullet(G,M)$}
Soient $A$ un anneau topologique, $G$ un groupe profini, $M$ un $A$-$G$-module topologique. 
On associe à $M$ un $A[G]$-module cosimplicial $\rK^\bullet(G,M)$ \eqref{higgs1-not580},
bifonctoriel contravariant en $G$ et covariant en $M$, de la façon suivante~: pour tout entier $n\in \mN$, 
$\rK^n(G,M)$ est le $A$-module des applications continues de $G^{[n]}$ dans $M$, 
équipé de l'action de $G$ définie pour tous $g\in G$ et $f\in \rK^n(G,M)$ par 
\begin{equation}\label{higgs1-not78a}
(g\cdot f)(g_0,\dots,g_n)=g\cdot f(g^{-1}g_0,\dots,g^{-1}g_n).
\end{equation}
Tout morphisme $\delta\colon [n]\rightarrow [m]$
de $\fD$ \eqref{higgs1-not580} induit un morphisme de $A[G]$-modules
\begin{equation}\label{higgs1-not78b}
\rK^n(G,M)\rightarrow \rK^m(G,M),
\end{equation}
défini par composition avec l'application 
\begin{equation}\label{higgs1-not78c}
G^{[m]}\rightarrow G^{[n]}, \ \ \ (g_0,\dots,g_m)\mapsto (g_{\delta(0)},\dots,g_{\delta(n)}).
\end{equation}
Le complexe de cochaînes $\trK^\bullet(G,M)$ associé à $\rK^\bullet(G,M)$ \eqref{higgs1-not580}
est appelé {\em complexe de cochaînes homogènes continues de $G$ à valeurs dans $M$}. 

On appelle complexe de {\em cochaînes non homogènes continues} de $G$ à valeurs dans $M$,
et l'on note $\rC^\bullet_\cont(G,M)$, 
le complexe de $A$-modules défini comme suit~: pour tout $n\in \mN$, 
$\rC^n_\cont(G,M)$ est l'ensemble de toutes les applications continues de $G^n$ dans $M$, et la différentielle 
$d\colon \rC^n_\cont(G,M)\rightarrow \rC^{n+1}_\cont(G,M)$ est définie par la formule 
\begin{eqnarray}\label{higgs1-not780a}
d(f)(g_1,\dots,g_{n+1})&=&g_1\cdot f(g_2,\dots,g_{n+1})\\
&&+\sum_{i=1}^{n}(-1)^if(g_1,\dots,g_ig_{i+1},\dots,g_{n+1})\nonumber\\
&&+(-1)^{n+1}f(g_1,\dots,g_n).\nonumber
\end{eqnarray}
Les groupes de cohomologie de $\rC^\bullet_\cont(G,M)$ sont notés $\rH^\bullet_\cont(G,M)$ et 
appelés les groupes de {\em cohomologie continue} de $G$ à valeurs dans $M$ (cf. \cite{tate} § 2).
Les applications continues (pour tout $n\in \mN$)
\begin{equation}\label{higgs1-not78e}
G^n\rightarrow G^{[n]}, \ \ \ (g_1,\dots,g_n)\mapsto (1,g_1, g_1g_2,\dots,g_1\dots g_n)
\end{equation} 
induisent un isomorphisme de complexes de $A$-modules
\begin{equation}\label{higgs1-not78f}
\trK^\bullet(G,M)^G\stackrel{\sim}{\rightarrow}\rC^\bullet_\cont(G,M).
\end{equation}

Supposons $M$ discret. Pour tout $n\in \mN$, 
$\rK^n(G,M)$ est alors un $A$-$G$-module discret et induit (\cite{tsuji2} 11.9). 
Notons $\varepsilon\colon  M\rightarrow \rK^0(G,M)$ 
le morphisme de $A$-$G$-modules discrets défini pour tout $m\in M$ par $\varepsilon(m)(g_0)=m$. 
On a $d^{0}\circ \varepsilon =0$, où $d^0$ est la différentielle de degré $0$ du complexe $\trK^\bullet(G,M)$.
Le complexe augmenté de $A$-modules $M\stackrel{\varepsilon}{\rightarrow}\trK^{\bullet}(G,M)$ 
est homotope à $0$ (\cite{tsuji2} 11.8). 
On en déduit un isomorphisme canonique fonctoriel de $\bD^+(\bMod(A))$
\begin{equation}\label{higgs1-not78g}
\rR\Gamma(G,M)\stackrel{\sim}{\rightarrow}\rC^\bullet_\cont(G,M).
\end{equation}
On peut donc omettre l'indice ``cont'' des notations $\rC^\bullet_\cont$ et $\rH^\bullet_\cont$
sans induire d'ambiguïté.

\subsection{}\label{higgs1-limproj1}\index{10390@$\cA^\mN$ ($\cA$ catégorie abélienne)}
Pour toute catégorie abélienne $\cA$, on note $\cA^\mN$ la catégorie des systèmes projectifs de $\cA$ indexés 
par l'ensemble ordonné $\mN$. Alors $\cA^\mN$ est une catégorie abélienne, dont les noyaux et les conoyaux 
se calculent composante par composante. Si $\cA$ a suffisamment d'injectifs, il en est de même de $\cA^\mN$ 
(\cite{jannsen} 1.1); un objet $(A_n,d_n)_{n\in \mN}$ de $\cA^\mN$ est injectif si et seulement si pour tout $n\in \mN$, 
$A_n$ est injectif et $d_n\colon A_{n+1}\rightarrow A_n$ est inversible à gauche.   

Soit $h\colon \cA\rightarrow \cB$  un foncteur exact à gauche entre catégories abéliennes. 
On note $h^\mN\colon \cA^\mN\rightarrow \cB^\mN$ son prolongement naturel. 
Supposons que $\cA$ ait suffisamment d'injectifs.
Alors $\rR^ih^\mN=(\rR^i h)^\mN$ pour tout $i\geq 0$ (\cite{jannsen} 1.2).  
Supposons de plus que les limites projectives indexées par $\mN$ soient représentables dans $\cB$. 
On désigne alors par  
\begin{equation}\label{higgs1-limproj1a}
\underset{\underset{n}{\longleftarrow}}{\lim}\ h\colon \cA^\mN\rightarrow \cB
\end{equation}
le foncteur qui associe à $(A_n,d_n)_{n\in \mN}$ la limite projective 
\begin{equation}\label{higgs1-limproj1b}
\underset{\underset{n}{\longleftarrow}}{\lim}\ (h(A_n),h(d_n)),
\end{equation}
et par  $\rR^+(\underset{\underset{n}{\longleftarrow}}{\lim}\ h)$
son foncteur dérivé droit. Si $h$ transforme les injectifs en injectifs, 
il en est de même de $h^\mN$. Si de plus le foncteur dérivé droit $\rR^+(\underset{\underset{n}{\longleftarrow}}{\lim})$ 
du foncteur 
\begin{equation}
\underset{\underset{n}{\longleftarrow}}{\lim}\colon \cB^\mN\rightarrow \cB
\end{equation} 
existe, on a alors un isomorphisme canonique 
\begin{equation}\label{higgs1-limproj1d}
\rR^+(\underset{\underset{n}{\longleftarrow}}{\lim} \ h)\stackrel{\sim}{\rightarrow}  
\rR^+(\underset{\underset{n}{\longleftarrow}}{\lim})\circ (\rR^+ h)^\mN.
\end{equation}

\subsection{}\label{higgs1-limproj2}
Soient $A$ un anneau, $G$ un groupe profini.    
Les limites projectives dans $\bMod(A)$ sont représentables~; le foncteur 
\begin{equation}\label{higgs1-limproj2a}
\underset{\underset{n}{\longleftarrow}}{\lim}\colon \bMod(A)^\mN\rightarrow \bMod(A)
\end{equation}
admet un foncteur dérivé droit~; et on a $\rR^i(\underset{\underset{n}{\longleftarrow}}{\lim})=0$ pour tout $i\geq 2$
(cf. \cite{jannsen} 1.4 et \cite{roos} 2.1). Si $(M_n)_{n\in \mN}$ est un système projectif de $A$-modules
vérifiant la condition de Mittag-Leffler, on a (\cite{jannsen} 1.15 et \cite{roos} 3.1)
\begin{equation}\label{higgs1-limproj2e}
\rR^1\underset{\underset{n}{\longleftarrow}}{\lim} \ M_n=0.
\end{equation}

On désigne par 
\begin{equation}\label{higgs1-limproj2b}
\Gamma(G,-)\colon (\bRep_A^\disc(G))^\mN\rightarrow \bMod(A)
\end{equation}
le foncteur $\underset{\underset{n}{\longleftarrow}}{\lim} \ \Gamma(G,-)$ \eqref{higgs1-limproj1a} et par 
$\rR^+\Gamma(G,-)$ son foncteur dérivé droit. 
Pour tout système projectif $(M_n)_{n\in \mN}$ de $\bRep_A^\disc(G)$
et pour tout entier $i\geq 0$, on a, d'après \eqref{higgs1-limproj1d}, une suite exacte
\begin{equation}\label{higgs1-limproj2c}
0\rightarrow \rR^1\underset{\underset{n}{\longleftarrow}}{\lim} \ \rH^{i-1}(G,M_n)\rightarrow 
\rH^i(G,(M_n)_{n\in \mN})\rightarrow 
\underset{\underset{n}{\longleftarrow}}{\lim} \ \rH^{i}(G,M_n)\rightarrow 0.
\end{equation}
Supposons que $(M_n)_{n\in \mN}$ vérifie la condition de Mittag-Leffler et notons $M$ sa 
limite projective en tant que $A$-$G$-module topologique.  
D'après (\cite{jannsen} 2.2), on a un isomorphisme canonique dans $\bD^+(\bMod(A))$
\begin{equation}\label{higgs1-limproj2d}
\rC^\bullet_\cont(G,M)\stackrel{\sim}{\rightarrow}
\rR^+\Gamma(G,(M_n)_{n\in \mN}).
\end{equation}

\subsection{}\label{higgs1-not5951}
Soient $A$ un anneau, $G$ et $H$ deux groupes. Si $M$ est un $A[G]$-module 
et $N$ est un $A[H]$-module, on désigne par $M\boxtimes N$ le $A[G\times H]$-module 
de $A$-module sous-jacent $M\otimes_AN$ tel que pour tous $(g,h)\in G\times H$ et $(x,y)\in M\times N$, on ait 
\begin{equation}\label{higgs1-not5951a}
(g,h)\cdot(x\otimes y)=(g\cdot x)\otimes (h\cdot y).
\end{equation}
Si $G$ et $H$ sont profinis, si $M$ est un $A$-$G$-module discret et si $N$ est un $A$-$H$-module discret, 
$M\boxtimes N$ est alors un $A$-$(G\times H)$-module discret.

Si $C$ est un complexe de $A[G]$-modules et $D$ est un complexe de $A[H]$-modules, 
on désigne par $C\boxtimes D$ le complexe simple associé 
au bicomplexe de $A[G\times H]$-modules défini par $(m,n)\mapsto C^m\boxtimes D^n$.

Si $X$ est un $A[G]$-module cosimplicial et $Y$ est un $A[H]$-module cosimplicial, 
on désigne par
$X\boxtimes Y$ le $A[G\times H]$-module $2$-cosimplicial défini par $(m,n)\mapsto X^m\boxtimes Y^n$,
et par $X\boxtimes^\Delta Y$ le sous-objet diagonal $\Delta(X\boxtimes Y)$ \eqref{higgs1-not580}.

\subsection{}\label{higgs1-not582}
Soient $A$ un anneau topologique, $G$ et $H$ deux groupes profinis, $M$ un $A$-$G$-module topologique, 
$N$ un $A$-$H$-module topologique. Avec les notations de \ref{higgs1-not78} et \ref{higgs1-not5951}, pour tout $n\in \mN$, 
on a un morphisme de $A[G\times H]$-modules 
\begin{equation}\label{higgs1-not582c}
\rK^n(G,M)\boxtimes \rK^n(H,N)\rightarrow \rK^n(G\times H,M\boxtimes N)
\end{equation}
défini pour $\varphi \in \rK^n(G,M)$ et $\psi \in \rK^n(H,N)$, par
\begin{equation}\label{higgs1-not582d}
\varphi\otimes \psi \mapsto (((g_0,h_0),\dots,(g_n,h_n))\mapsto \varphi(g_0,\dots,g_n)\otimes \psi(h_0,\dots,h_n)).
\end{equation} 
On obtient ainsi un morphisme de $A[G\times H]$-modules cosimpliciaux \eqref{higgs1-not580}
\begin{equation}\label{higgs1-not582e}
\rK^\bullet(G,M)\boxtimes^\Delta \rK^\bullet(H,N)\rightarrow \rK^\bullet(G\times H,M\boxtimes N).
\end{equation}
On en déduit un morphisme de $A$-modules cosimpliciaux 
\begin{equation}\label{higgs1-not582f}
\rK^\bullet(G,M)^G\otimes^\Delta_A \rK^\bullet(H,N)^H\rightarrow \rK^\bullet(G\times H,M\boxtimes N)^{G\times H}.
\end{equation}
Prenant les complexes de cochaînes associés et composant avec la flèche \eqref{higgs1-not584a}, on obtient 
un morphisme de complexes de $A$-modules \eqref{higgs1-not78f}
\begin{equation}\label{higgs1-not582g}
\rC^\bullet_\cont(G,M)\otimes_A \rC^\bullet_\cont(H,N)\rightarrow \rC^\bullet_\cont(G\times H,M\boxtimes N).
\end{equation}
Celui-ci induit pour tous entiers $m,n\geq 0$, le {\em cross-produit}
\begin{equation}\label{higgs1-not582h}
\rH^m_\cont(G,M)\otimes_A\rH^n_\cont(H,N)\rightarrow \rH^{m+n}_\cont(G\times H,M\boxtimes N), 
\ \ \ x\otimes y\mapsto x\times y.
\end{equation}
Les morphismes \eqref{higgs1-not582e}, \eqref{higgs1-not582g} et \eqref{higgs1-not582h} sont fonctoriels contravariants en $G$ et $H$ et covariants en $M$ et $N$.

\subsection{}\label{higgs1-not583}
Conservons les hypothèses de \eqref{higgs1-not582}, supposons de plus $G=H$. 
Composant le morphisme \eqref{higgs1-not582g} avec le morphisme 
\begin{equation}\label{higgs1-not583a}
\rC^\bullet_\cont(G\times G,M\boxtimes N)\rightarrow
\rC^\bullet_\cont(G,M\otimes N)
\end{equation}
induit par l'homomorphisme diagonal $G\rightarrow G\times G$, on obtient un morphisme 
de complexes de $A$-modules 
\begin{equation}\label{higgs1-not583b}
\rC^\bullet_\cont(G,M)\otimes_A \rC^\bullet_\cont(G,N)\rightarrow \rC^\bullet_\cont(G,M\otimes N).
\end{equation}
Celui-ci induit pour tous entiers $m,n\geq 0$, le {\em cup-produit}
\begin{equation}\label{higgs1-not583c}
\rH^m_\cont(G,M)\otimes_A\rH^n_\cont(G,N)\rightarrow \rH^{m+n}_\cont(G,M\otimes_AN), \ \ \ x\otimes y\mapsto x\cup y.
\end{equation}
Les morphismes \eqref{higgs1-not583b} et \eqref{higgs1-not583c} sont fonctoriels 
contravariants en $G$ et covariants en $M$ et $N$.

\begin{rema}\label{higgs1-not5831}
Sous les hypothèses de \eqref{higgs1-not582}, pour tous $m,n\geq 0$, le diagramme 
\begin{equation}\label{higgs1-not5831a}
\xymatrix{
{\rH^m_\cont(G,M)\otimes_A\rH^n_\cont(H,N)}\ar[d]_{\pi_1^*\times \pi_2^*}\ar[rd]^\times &\\
{\rH^m_\cont(G\times H,M)\otimes_A\rH^n_\cont(G\times H,N)}\ar[r]^-(0.4)\cup&{\rH^{m+n}_\cont(G\times H,M\boxtimes N)}}
\end{equation}
où $\pi_1\colon G\times H\rightarrow G$ et $\pi_2\colon G\times H\rightarrow H$ sont les projections canoniques, est commutatif. Ceci résulte du caractère fonctoriel du cross-produit et du fait que le composé
\begin{equation}
\xymatrix{
{G\times H}\ar[r]^-(0.4)\delta&{(G\times H)\times (G\times H)}\ar[r]^-(0.4){\pi_1\times \pi_2}&{G\times H}},
\end{equation}
où $\delta$ est l'homomorphisme diagonal, est l'identité. 
\end{rema}

\begin{lem}\label{higgs1-cg2}
Soient $A$ un anneau, $G$ un groupe profini, $M$ un $A$-$G$-module discret,
$N$ un $A$-module plat. On a alors un isomorphisme canonique bifonctoriel de $\bD^+(\bMod(A))$
\begin{equation}\label{higgs1-cg2a}
\rR\Gamma(G,M)\otimes^\rL_AN\stackrel{\sim}{\rightarrow} \rR\Gamma(G,M\otimes_AN).
\end{equation}
\end{lem}

Montrons d'abord que pour tout entier $q\geq 0$, on a un isomorphisme 
\begin{equation}\label{higgs1-cg2b}
\rH^q(G,M)\otimes_AN\stackrel{\sim}{\rightarrow} \rH^q(G,M\otimes_AN).
\end{equation}
On a un isomorphisme canonique 
\begin{equation}\label{higgs1-cg2c}
\rH^q(G,M)\stackrel{\sim}{\rightarrow} \underset{\underset{H}{\longrightarrow}}{\lim}\ \rH^q(G/H,M^H),
\end{equation}
où $H$ parcourt l'ensemble des sous-groupes ouverts distingués de $G$; et de même pour $M\otimes_AN$.
D'autre part, $N$ étant $A$-plat, pour tout sous-groupe $H$ de $G$, le morphisme canonique
\begin{equation}\label{higgs1-cg2d}
M^H\otimes_AN\rightarrow (M\otimes_AN)^H
\end{equation}
est un isomorphisme. On peut donc se réduire au cas où $G$ est fini. Soit 
\begin{equation}\label{higgs1-cg2e}
\dots \rightarrow P_i\rightarrow P_{i-1}\rightarrow \dots \rightarrow P_0\rightarrow \mZ\rightarrow 0
\end{equation}
une résolution du $\mZ[G]$-module $\mZ$ par des $\mZ[G]$-modules libres de type fini. 
Pour tout $i\geq 0$, le morphisme canonique 
\begin{equation}\label{higgs1-cg2f}
\Hom_{\mZ[G]}(P_i,M)\otimes_AN\rightarrow \Hom_{\mZ[G]}(P_i,M\otimes_AN)
\end{equation}
est alors un isomorphisme~; d'où l'assertion recherchée puisque $\rH^*(G,M)$ est
la cohomologie du complexe de cochaînes $\Hom_{\mZ[G]}(P_\bullet,M)$. 

Avec les notations de \ref{higgs1-not78}, le complexe $\trK(G,M)\otimes_AN$ est une résolution de $M\otimes_AN$
par des $A$-$G$-modules discret et $G$-acycliques en vertu de \eqref{higgs1-cg2b}. Par ailleurs, le morphisme canonique 
\begin{equation}
\rC^\bullet(G,M)\otimes_A N\rightarrow (\trK(G,M)\otimes_AN)^G
\end{equation}
est un isomorphisme d'après \eqref{higgs1-cg2d}~; d'où la proposition.

\begin{lem}\label{higgs1-cg21}
Soit $A$ un anneau de Dedekind, $G$ un groupe profini, $M$ un $A$-$G$-module discret, $G$-acyclique
et plat en tant que $A$-module, $N$ un $A$-module. Alors $M\otimes_AN$ est $G$-acyclique et le morphisme canonique
\begin{equation}\label{higgs1-cg21a}
M^G\otimes_AN\rightarrow (M\otimes_AN)^G
\end{equation}
est un isomorphisme.
\end{lem}
En effet, on peut se borner au cas où $N$ est de type fini (\cite{sga4} VI 5.3). Il existe  
alors deux $A$-modules projectifs de type fini $N_1$ et $N_2$ et une suite exacte 
$0\rightarrow N_2\rightarrow N_1\rightarrow N\rightarrow 0$. 
On se réduit ainsi au cas où $N$ est projectif de type fini, et même au cas où $N$ est libre de type fini,
auquel cas la proposition est immédiate.  

\begin{lem}\label{higgs1-kun5}
Soient $A$ un anneau de Dedekind, $C$ et $D$ deux complexes de $A$-modules
tels que $D$ soit $A$-plat (i.e., toutes ses composantes sont $A$-plates). 
Si $C$ ou $D$ est acyclique, il en est de même de $C\otimes_AD$.    
\end{lem}

Cela résulte de la formule de Künneth (\cite{alg10} §4.7 théo.~3). 

\begin{lem}\label{higgs1-kun6}
Soient $A$ un anneau de Dedekind, $C$, $C'$ et $D$ des complexes de $A$-modules, 
$u\colon C\rightarrow C'$ un quasi-isomorphisme.  
Si $D$ est $A$-plat ou si $C$ et $C'$ sont $A$-plats, alors
$u\otimes \id\colon C\otimes_AD\rightarrow C'\otimes_AD$ est un quasi-isomorphisme. 
\end{lem}

Cela résulte de \ref{higgs1-kun5} et (\cite{alg10} §4.3 lem.~2).

\subsection{}\label{higgs1-kun7}
Soient $A$ un anneau de Dedekind, $G$ et $H$ deux groupes profinis, $M$ un $A$-$G$-module discret
tel que le $A$-module sous-jacent soit plat, $N$ un $A$-$H$-module discret. Comme le complexe $\trK^\bullet(G,M)$
est $A$-plat \eqref{higgs1-not78}, $\trK^\bullet(G,M)\boxtimes \trK^\bullet(H,N)$ \eqref{higgs1-not5951}
est une résolution de $M\boxtimes N$ par des $A$-$(G\times H)$-modules discrets et $(G\times H)$-acycliques, 
en vertu de \ref{higgs1-not594}, \ref{higgs1-cg2}, \ref{higgs1-cg21} et \ref{higgs1-kun6}. 
C'est aussi une résolution acyclique de $M\boxtimes N$ pour le foncteur \eqref{higgs1-not593a}
\begin{equation}
\Gamma(H,-)\colon \bRep_A^\disc(G\times H) \rightarrow \bRep_A^\disc(G), \ \ \ L\mapsto L^H,
\end{equation}
d'après \ref{higgs1-cg2} et \eqref{higgs1-not593b}. On a un isomorphisme canonique 
\begin{equation}\label{higgs1-kun7a}
\trK^\bullet(G,M)\otimes_A \rC^\bullet(H,N) \stackrel{\sim}{\rightarrow} (\trK^\bullet(G,M)\otimes_A \trK^\bullet(H,N))^H.
\end{equation}
Par ailleurs, le morphisme canonique 
\begin{equation}\label{higgs1-kun7b}
M\otimes_A \rC^\bullet(H,N)\rightarrow  \trK^\bullet(G,M)\otimes_A \rC^\bullet(H,N)
\end{equation}
est un quasi-isomorphisme \eqref{higgs1-kun6}~; $\trK^\bullet(G,M)\otimes_A \rC^\bullet(H,N)$ est un complexe de 
$A$-$G$-modules discrets et $G$-acycliques, et le morphisme canonique 
\begin{equation}\label{higgs1-kun7c}
\rC^\bullet(G,M)\otimes_A \rC^\bullet(H,N)\rightarrow 
(\trK^\bullet(G,M)\otimes_A \rC^\bullet(H,N))^G
\end{equation}
est un isomorphisme d'après \ref{higgs1-cg21}. Par suite, le morphisme canonique \eqref{higgs1-not582g} 
\begin{equation}\label{higgs1-kun7d}
\rC^\bullet(G,M)\otimes_A \rC^\bullet(H,N)\rightarrow \rC^\bullet(G\times H,M\boxtimes N)
\end{equation}
est un quasi-isomorphisme. 
Le composé de \eqref{higgs1-kun7d} et de l'inverse de \eqref{higgs1-kun7c}, induit alors un isomorphisme de $\bD^+(\bMod(A))$
\begin{equation}\label{higgs1-kun7e}
\rR\Gamma(G,M\otimes^\rL_A \rR\Gamma(H,N))\stackrel{\sim}{\rightarrow}\rR\Gamma(G\times H,M\boxtimes N).
\end{equation}
On en déduit, compte tenu de \eqref{higgs1-cg2a}, un isomorphisme de $\bD^+(\bMod(A))$
\begin{equation}\label{higgs1-kun7f}
\rR\Gamma(G,\rR\Gamma(H,M\boxtimes N))\stackrel{\sim}{\rightarrow}\rR\Gamma(G\times H,M\boxtimes N).
\end{equation}
Celui-ci est égal à l'isomorphisme \eqref{higgs1-not594a}.

\begin{lem}\label{higgs1-tsj1}
Soient $G$ un groupe profini isomorphe à $\mZ_p$, $\gamma$ un générateur topologique de $G$, 
$A$ un anneau, $M$ un $A$-$G$-module discret  
dont tous les éléments sont de torsion $p$-primaire. 
On a alors une suite exacte de $A$-$G$-modules discrets 
\begin{equation}\label{higgs1-tsj1a}
0\longrightarrow M\stackrel{\varepsilon}{\longrightarrow} \Ind_{A,G}(M) \stackrel{d_\gamma}{\longrightarrow} \Ind_{A,G}(M)
\longrightarrow 0,
\end{equation}
où pour tous $x\in M$,  $f\in \Ind_{A,G}(M)$ et $g\in G$, on a 
$\varepsilon(x)(g)=g\cdot x$ et $d_\gamma (f)(g)=\gamma  f(\gamma  ^{-1}g)-f(g)$. 
\end{lem}
Il n'y a que la surjectivité de $d_\gamma$ qui nécessite une preuve. Soit $f\colon G\rightarrow M$ une application continue. 
Comme $G$ est compact, il existe un entier $n\geq0$ tel que $f$ se factorise à travers 
$G/\gamma^{p^n\mZ_p}$, que $p^n f(G)=0$ et que
pour tout $x\in f(G)$, on ait $\gamma^{p^n}\cdot x=x$. On désigne par $h\colon \mZ_{>0}\rightarrow M$ 
l'application définie pour tout $a\in  \mZ_{>0}$, par 
\begin{equation}
h(a)=-\sum_{i=1}^a\gamma^{a-i}\cdot f(\gamma^i).
\end{equation} 
Celle-ci se factorise à travers $\mZ/p^{2n}\mZ$. Composant avec l'homomorphisme surjectif $G\rightarrow \mZ/p^{2n}\mZ$
qui envoie $\gamma$ sur la classe de $1$, 
on obtient une application continue $\tth\colon G\rightarrow M$ telle que $d_\gamma(\tth)=f$.

\begin{lem}\label{higgs1-tsj11}
Conservons les hypothèses de \eqref{higgs1-tsj1}, notons de plus 
\begin{equation}\label{higgs1-tsj11a}
\partial_M\colon M\rightarrow \rH^1(G,M)
\end{equation}
le morphisme bord de la suite exacte \eqref{higgs1-tsj1a}. Alors~:
\begin{itemize}
\item[{\rm (i)}] Pour tout $x\in M$, il existe un et un unique homomorphisme 
croisé continu $\nu_x\colon G\rightarrow M$ tel que $\nu_x(\gamma)=x$. 
\item[{\rm (ii)}] Pour tout $x\in M$, la classe de $\nu_x$ dans $\rH^1(G,M)$ est égale à $-\partial_M(x)$. 
\end{itemize}
\end{lem}

(i) En effet, supposons qu'un tel homomorphisme croisé continu $\nu_x$ existe. 
Pour tout entier $a\geq 1$, on a alors
\begin{equation}\label{higgs1-tsj11b}
\nu_x(\gamma^a)=(\gamma^{a-1}+\dots+1)\cdot x.
\end{equation}
Soit $n$ un entier $\geq 0$ tel que $p^n x=0$ et $\gamma^{p^n}\cdot x=x$. Pour tout $g\in G$, on a 
\begin{equation}\label{higgs1-tsj11c}
\nu_x(\gamma^{p^{2n}}g)=g\cdot \nu_x(\gamma^{p^{2n}})+\nu_x(g)=\nu_x(g).
\end{equation}
Par continuité, $\nu_x$ se factorise donc à travers l'homomorphisme canonique $G\rightarrow G/\gamma^{p^{2n}\mZ_p}$. 
Par ailleurs, l'application 
\begin{equation}\label{higgs1-tsj11d}
\mZ_{>0}\rightarrow M,\ \ \ a\mapsto (\gamma^{a-1}+\dots+1)\cdot x,
\end{equation}
se factorise à travers $\mZ/p^{2n}\mZ$. Composant avec l'homomorphisme surjectif 
$G\rightarrow \mZ/p^{2n}\mZ$ qui envoie $\gamma$ sur la classe de $1$, on obtient un homomorphisme 
croisé continu de $G$ à valeur dans $M$; d'où l'assertion. 

(ii) En effet, pour tout $g\in G$, on a \eqref{higgs1-tsj1a}
\begin{equation}\label{higgs1-tsj11e}
d_\gamma(\nu_x)(g)=\gamma\cdot \nu_x(\gamma^{-1}g)-\nu_x(g)=-\nu_x(\gamma)=-x.
\end{equation}

\begin{lem}\label{higgs1-tsj12}
Soient $G$ un groupe profini isomorphe à $\mZ_p$, $\gamma$ un générateur topologique de $G$, 
$A$ un anneau de torsion $p$-primaire, $M$ un $A$-$G$-module discret,  $x\in M^G$. 
On a alors 
\begin{equation}\label{higgs1-tsj12a}
\partial_M(x)=\partial_A(1) \cup x,
\end{equation}
où $\partial_M$ est le morphisme bord de la suite exacte \eqref{higgs1-tsj1a} et $\cup$ est le cup-produit. 
\end{lem}

Cela résulte de \ref{higgs1-tsj11}. En effet, notons $\nu_x\colon G\rightarrow M$ et $\mu_1\colon G\rightarrow A$ 
les homomorphismes croisés continus tels que $\nu_x(\gamma)=x$ et $\mu_1(\gamma)=1$. 
Pour tout entier $a\geq 1$, on a $\nu_x(\gamma^a)=ax=\mu_1(\gamma^a) x$ d'après \eqref{higgs1-tsj11b}. 
On en déduit que $\nu_x(g)=\mu_1(g) x$ pour tout $g\in G$; d'où la proposition.

\begin{prop}\label{higgs1-tsj2}
Soit $n$ un entier $\geq 1$, et pour tout entier $1\leq i\leq n$, soient $G_i$ un groupe profini isomorphe à $\mZ_p$,
$\gamma_i$ un générateur topologique de $G_i$. Posons $G=\prod_{i=1}^nG_i$. Soient $A$ un anneau, 
$M$ un $A$-$G$-module discret dont tous les éléments sont de torsion $p$-primaire.
On définit par récurrence, pour tout entier $0\leq i\leq n$, un complexe de $A$-$G$-modules discrets $K_i^\bullet$   
en posant $K_0^\bullet=M[0]$ et pour tout $1\leq i\leq n$, $K_i^\bullet$ est la fibre du morphisme $\gamma_i-1\colon 
K_{i-1}^\bullet\rightarrow K_{i-1}^\bullet$. Il existe alors un isomorphisme canonique dans $\bD^+(\bMod(A))$
\begin{equation}\label{higgs1-tsj2a}
\rR\Gamma(G,M)\stackrel{\sim}{\rightarrow}K^\bullet_n.
\end{equation}
\end{prop}

Pour tout entier $0\leq j\leq n$, on pose $G_{\leq j}=\prod_{1\leq i\leq j}G_i$ et $G_{>j}=\prod_{j<i\leq n}^nG_i$. Montrons 
par récurrence que pour tout entier $0\leq j\leq n$, on a un isomorphisme canonique dans $\bD^+(\bRep_{A}^\disc(G_{>j}))$
\begin{equation}\label{higgs1-tsj2b}
\rR\Gamma(G_{\leq j},M)\stackrel{\sim}{\rightarrow}K^\bullet_j.
\end{equation}
L'assertion est immédiate pour $j=0$. Supposons l'assertion établie pour un entier $0\leq j\leq n-1$. 
Notons $\Ind_{A,G_{j+1}}(K^\bullet_j)$ l'image du complexe $K^\bullet_j$ par le foncteur $\Ind_{A,G_{j+1}}$ \eqref{higgs1-not59b}. 
On le munit de l'action de $G_{>j}$ définie pour tous $q\in \mZ$, $f\in \Ind_{A,G_{j+1}}(K^q_j)$,
$g=(g_0,g_1)\in G_{>j}=G_{j+1}\times G_{>j+1}$ et $x\in G_{j+1}$, par 
\begin{equation}\label{higgs1-tsj2c}
(g\cdot f)(x)=g_1\cdot f(x\cdot g_0)\in K^q_j. 
\end{equation}
Il résulte aussitôt de \ref{higgs1-tsj1} qu'on a une suite exacte de complexes de $\bRep_{A}^\disc(G_{>j})$
\begin{equation}\label{higgs1-tsj2f}
\xymatrix{
0\ar[r]&{K^\bullet_j}\ar[r]^-(0.5){\varepsilon_j}&{\Ind_{A,G_{j+1}}(K^\bullet_j)}
\ar[r]^{d_{\gamma_{j+1}}}&{\Ind_{A,G_{j+1}}(K^\bullet_j)}\ar[r]&0},
\end{equation}
où pour tous $q\in \mZ$, $c\in K^q_j$, $f\in \Ind_{A,G_{j+1}}(K^q_j)$ et $x\in G_{j+1}$, on a 
\begin{eqnarray}
\varepsilon_j^q(c)(x)&=&x\cdot c,\\
d^q_{\gamma_{j+1}}(f)(x)&=&\gamma_{j+1}  f(\gamma_{j+1}^{-1}x)-f(x).
\end{eqnarray} 
D'après l'hypothèse de récurrence, on a donc un triangle distingué dans $\bD^+(\bRep_{A}^\disc(G_{>j}))$
\begin{equation}\label{higgs1-tsj2d}
\xymatrix{
{\rR\Gamma(G_{\leq j},M)}\ar[r]&{\Ind_{A,G_{j+1}}(K^\bullet_j)}\ar[r]^{d_{\gamma_{j+1}}}&{\Ind_{A,G_{j+1}}(K^\bullet_j)}
\ar[r]^-(0.4){+1}&}.
\end{equation}
Compte tenu de \ref{higgs1-not59} et \eqref{higgs1-not593b}, pour tout entier $q$, le morphisme 
\begin{equation}
K^q_j\rightarrow \Gamma(G_{j+1},\Ind_{A,G_{j+1}}(K^q_j)), \ \ \ u\mapsto (x\mapsto u)
\end{equation}
induit un isomorphisme de $\rD^+(\bRep_{A}^\disc(G_{>j+1}))$
\begin{equation}
K^q_j[0]\stackrel{\sim}{\rightarrow}\rR\Gamma(G_{j+1},\Ind_{A,G_{j+1}}(K^q_j)).
\end{equation}
En vertu de \ref{higgs1-not594}, le triangle \eqref{higgs1-tsj2d} induit donc
un triangle distingué dans $\bD^+(\bRep_{A}^\disc(G_{>j+1}))$
\begin{equation}
\xymatrix{
{\rR\Gamma(G_{\leq j+1},M)}\ar[r]&
{K^\bullet_j}\ar[rr]^{\gamma_{j+1}-1}&&{K^\bullet_j}\ar[r]^-(0.4){+1}&};
\end{equation}
d'où la proposition. 

\begin{cor}\label{higgs1-cg05}
Pour tout entier $d\geq 1$, la $p$-dimension cohomologique du groupe profini $\mZ_p^d$ est égale à $d$. 
\end{cor}

En effet, pour tout $\mZ$-$\mZ_p^d$-module discret de torsion $p$-primaire $M$, 
le complexe $K^\bullet_d$ qui lui est associé dans \ref{higgs1-tsj2} est concentré 
en degrés $[0,d]$. La $p$-dimension cohomologique du groupe profini $\mZ_p^d$ 
est donc inférieure ou égale à $d$ en vertu de \ref{higgs1-tsj2} et (\cite{serre1}, I prop.~11). 
Pour le $\mZ$-$\mZ_p^d$-module discret trivial $\mF_p$, 
le complexe associé $K^\bullet_d$ est à différentielles nulles et $\rH^d(\mZ_p^d, \mF_p)$ 
est donc isomorphe à $\mF_p$ \eqref{higgs1-tsj2a}~; d'où la proposition.

\begin{cor}\label{higgs1-cg01}
Soient $n$ un entier $\geq 1$, $G$ un groupe profini isomorphe à $\mZ_p^n$, 
$e_1,\dots,e_n$ une $\mZ_p$-base de $G$, $A$ un anneau topologique, 
$M$ un $A$-$G$-module topologique. 
Notons $\varphi\colon G\rightarrow \Aut_A(M)$ la représentation de $G$
sur $M$, $\rS_A(G)$ l'algèbre symétrique du $A$-module $G\otimes_{\mZ}A$
et $M^\triangleright$ le $\rS_A(G)$-module dont le 
$A$-module sous-jacent est $M$ et tel que pour tout $1\leq i\leq n$, 
l'action de $e_i$ sur $M$ soit donnée par $\varphi(e_i)-\id_M$. 
Supposons que l'une des conditions suivantes soit remplie~:
\begin{itemize}
\item[{\rm (i)}] $M$ est un $A$-$G$-module discret de torsion $p$-primaire.
\item[{\rm (ii)}] $M$ est muni de la topologie $p$-adique et est complet et séparé pour cette topologie. 
\end{itemize} 
On a alors un isomorphisme canonique fonctoriel en $M$ de $\bD^+(\bMod(A))$
\begin{equation}\label{higgs1-cg01a}
\rC^\bullet_\cont(G,M)\stackrel{\sim}{\rightarrow} \mK^\bullet_{\rS_A(G)}(M^\triangleright),
\end{equation}
où le complexe de gauche est défini dans \eqref{higgs1-not78} et celui de droite dans \eqref{higgs1-koszul2d}.
\end{cor}

Le cas (i) résulte de \ref{higgs1-tsj2} et \eqref{higgs1-koszul2f}. 
Considérons ensuite le cas (ii) et posons  pour tout $r\geq 0$, $M_r=M/p^rM$. En vertu de \eqref{higgs1-limproj2d}, 
on a un isomorphisme canonique de $\bD^+(\bMod(A))$
\begin{equation}\label{higgs1-cg01g}
\rC_\cont^\bullet(G,M)\stackrel{\sim}{\rightarrow}\rR^+\Gamma(G,(M_r)_{r\in \mN}).
\end{equation}
D'autre part, d'après le cas (i), on a un système compatible d'isomorphismes de $\bD^+(\bMod(A))$
\begin{equation}\label{higgs1-cg01h}
\rR\Gamma(G,M_r)\stackrel{\sim}{\rightarrow} 
\mK^\bullet_{\rS_A(G)}(M^\triangleright/p^rM^\triangleright).
\end{equation}
Comme pour tout entier $n$, le système projectif
$(\mK^n_{\rS_A(G)}(M^\triangleright/p^rM^\triangleright))_{r\geq 0}$ vérifie la condition de Mittag-Leffler, 
on obtient, compte tenu de \eqref{higgs1-limproj1d} et \eqref{higgs1-limproj2e}, un isomorphisme 
\begin{equation}\label{higgs1-cg01i}
\rC_\cont^\bullet(G,M)\stackrel{\sim}{\rightarrow}\underset{\underset{r\geq 0}{\longleftarrow}}{\lim} \ 
\mK^\bullet_{\rS_A(G)}(M^\triangleright/p^rM^\triangleright).
\end{equation}
La proposition s'ensuit car 
$\mK^n_{\rS_A(G)}(M^\triangleright)$ est complet et séparé pour la topologie $p$-adique pour tout entier $n$.

\begin{rema}\label{higgs1-cg02}
Conservons les hypothèses de \ref{higgs1-cg01}. D'après  \ref{higgs1-not65}, le $\rS_A(G)$-module $M^\triangleright$ 
correspond à un $A$-champ de Higgs $\theta$ sur $M$ à coefficients dans $\Hom_{\mZ}(G,A)$,
et le complexe $\mK^\bullet_{\rS_A(G)}(M^\triangleright)$ s'identifie au complexe de Dolbeault de $(M,\theta)$. 
\end{rema}

\subsection{}\label{higgs1-kun1}
Soient $G$ un groupe profini isomorphe à $\mZ_p$, $\gamma$ un générateur topologique de $G$,
$H$ un groupe profini, $A$ un anneau de torsion $p$-primaire, $M$ un $A$-$H$-module discret,
que l'on considère aussi comme un $A$-$(G\times H)$-module discret via la projection canonique $G\times H\rightarrow H$. 
On désigne par $\Ind_{A,G}(\rC^\bullet(H,M))$ l'image du complexe de cochaînes continues
$\rC^\bullet(H,M)$ \eqref{higgs1-not78}  par le foncteur $\Ind_{A,G}$ \eqref{higgs1-not59b}. 
On munit $\rC^\bullet(H,M)$ de l'action triviale de $G$. 
D'après \ref{higgs1-tsj1}, on a une suite exacte de complexes de $\bRep_{A}^\disc(G)$
\begin{equation}\label{higgs1-kun1a}
0\longrightarrow \rC^\bullet(H,M)\stackrel{\varepsilon^\bullet}{\longrightarrow}
\Ind_{A,G}(\rC^\bullet(H,M))\stackrel{d_\gamma^\bullet}{\longrightarrow}
\Ind_{A,G}(\rC^\bullet(H,M))\longrightarrow 0,
\end{equation}
où pour tous $q\in \mZ$, $c\in \rC^q(H,M)$, $f\in \Ind_{A,G}(\rC^q(H,M))$ et $x\in G$, on a
\begin{eqnarray}
\varepsilon^q(c)(x)&=&x,\label{higgs1-kun1b}\\
d_\gamma^q(f)(x)&=&\gamma f(\gamma^{-1}x)-f(x).\label{higgs1-kun1c}
\end{eqnarray} 
Pour tout entier $q\geq 0$, $\varepsilon^q$ induit un isomorphisme de $\bD^+(\bMod(A))$
\begin{equation}\label{higgs1-kun1d}
\rC^q(H,M)[0]\stackrel{\sim}{\rightarrow}\rR\Gamma(G,\Ind_{A,G}(\rC^q(H,M))).
\end{equation}
En vertu de \ref{higgs1-not594}, \eqref{higgs1-kun1a} induit donc
un triangle distingué dans $\bD^+(\bMod(A))$
\begin{equation}\label{higgs1-kun1e}
\xymatrix{
{\rR\Gamma(G\times H,M)}\ar[r]&
{\rC^\bullet(H,M)}\ar[r]^{0}&{\rC^\bullet(H,M)}\ar[r]^-(0.4){+1}&},
\end{equation}
où $0$ est le morphisme nul. On en déduit pour tout entier $n\geq 0$, une suite exacte
\begin{equation}\label{higgs1-kun1f}
0\longrightarrow \rH^{n-1}(H,M)\stackrel{\alpha_n}{\longrightarrow} 
\rH^n(G\times H,M)\stackrel{\beta_n}{\longrightarrow} \rH^n(H,M)\longrightarrow 0.
\end{equation}
Par ailleurs, toujours d'après \ref{higgs1-tsj1}, on a une suite exacte canonique \eqref{higgs1-tsj1a}
\begin{equation}
0\rightarrow A\rightarrow \Ind_{A,G}(A) \rightarrow \Ind_{A,G}(A)\rightarrow 0.
\end{equation}
Elle induit un isomorphisme 
\begin{equation}\label{higgs1-kun1g}
\partial_A\colon A\stackrel{\sim}{\rightarrow} \rH^1(G,A).
\end{equation}

\begin{prop}\label{higgs1-kun2}
Les hypothèses étant celles de \eqref{higgs1-kun1}, soit de plus $n$ un entier $\geq 0$. Alors~:
\begin{itemize}
\item[{\rm (i)}] $\beta_n$ est le morphisme de restriction relativement à l'injection canonique $H\rightarrow G\times H$. 
\item[{\rm (ii)}] Pour tout $x\in \rH^{n-1}(H,M)$, on a 
\begin{equation}\label{higgs1-kun2aa}
\alpha_n(x)=\partial_A(1)\times x,
\end{equation}
où le cross-produit est défini dans \eqref{higgs1-not582h}.
\end{itemize}
\end{prop}

Notons $\tau_{\leq n}\rC^\bullet(H,M)$ la filtration canonique de $\rC^\bullet(H,M)$ (\cite{deligne} 1.4.6). 
D'après \ref{higgs1-tsj1}, on a un diagramme commutatif de complexes de $\bRep_{A}^\disc(G)$
\begin{equation}\label{higgs1-kun2a}
\xymatrix{
0\ar[r]&{\tau_{\leq n}\rC^\bullet(H,M)}\ar[r]^-(0.5){\psi_n}\ar[d]_{u_n}&
{\Ind_{A,G}(\tau_{\leq n}\rC^\bullet(H,M))}\ar[r]^{\phi_n}\ar[d]_{\Ind_{A,G}(u_n)}&
{\Ind_{A,G}(\tau_{\leq n}\rC^\bullet(H,M))}\ar[r]\ar[d]_{\Ind_{A,G}(u_n)}& 0\\
0\ar[r]&{\rH^n(H,M)[-n]}\ar[r]&{\Ind_{A,G}(\rH^n(H,M))[-n]}\ar[r]&
{\Ind_{A,G}(\rH^n(H,M))[-n]}\ar[r]& 0}
\end{equation}
où $u_n$ est le morphisme canonique et les lignes horizontales 
sont les suites exactes définies comme dans \eqref{higgs1-kun1a}.
En fait, le foncteur $\Ind_{A,G}$ étant exact, 
la ligne horizontale supérieure se déduit de la suite exacte \eqref{higgs1-kun1a} en appliquant le foncteur $\tau_{\leq n}$. 

(i)  On a clairement $\beta_n=\rH^n(G,\Ind_{A,G}(u_n)\circ \psi_n)$. 
Compte tenu de \eqref{higgs1-kun2a}, on en déduit que $\beta_n=\rH^n(G,u_n)$. 
La proposition résulte donc de \ref{higgs1-not595}. 

(ii)  Il résulte par récurrence de \ref{higgs1-cg05} que pour tout entier $i\geq 2$, on a 
\begin{equation}\label{higgs1-kun2b}
\rH^n(G,\tau_{\leq n-i} \rC^\bullet(H,M))=0.
\end{equation}
Par suite, le morphisme canonique
\begin{equation}\label{higgs1-kun2c}
\rH^n(G,\tau_{\leq n-1} \rC^\bullet(H,M))\rightarrow\rH^1(G,\rH^{n-1}(H,M))
\end{equation}
est un isomorphisme. On en déduit, compte tenu de \eqref{higgs1-not594a}, un morphisme canonique 
\begin{equation}\label{higgs1-kun2d}
v_n \colon \rH^1(G,\rH^{n-1}(H,M))\rightarrow \rH^n(G\times H,M).
\end{equation}
Considérons par ailleurs le diagramme commutatif
\begin{equation}\label{higgs1-kun2e}
\xymatrix{
{\Ind_{A,G}(\tau_{\leq n-1}\rC^\bullet(H,M))}\ar[r]^-(0.5){\delta_{\leq n-1}}\ar[d]_{\Ind_{A,G}(u_{n-1})}
&{\tau_{\leq n-1}\rC^\bullet(H,M)[1]}\ar[d]^{u_{n-1}[1]}\\
{\Ind_{A,G}(\rH^{n-1}(H,M))[-n+1]}\ar[r]^-(0.5){\delta_{n-1}}&{\rH^{n-1}(H,M)[-n+2]}}
\end{equation}
induit par le diagramme \eqref{higgs1-kun2a} (pour $n-1$ au lieu de $n$). Posons 
\begin{equation}
\partial_{n-1}=\rH^{n-1}(G,\delta_{n-1})\colon \rH^{n-1}(H,M)\rightarrow \rH^1(G,\rH^{n-1}(H,M)),
\end{equation}\label{higgs1-kun2f}
qui est en fait un isomorphisme. Comme \eqref{higgs1-kun2c} est un isomorphisme, on a 
\begin{equation}
\alpha_n=v_n\circ \rH^{n-1}(G,u_{n-1}[1]\circ \delta_{\leq n-1})=v_n\circ \partial_{n-1}.
\end{equation}
D'après \ref{higgs1-tsj12}, pour tout $x\in \rH^{n-1}(H,M)$, on a 
\begin{equation}\label{higgs1-kun2g}
\partial_{n-1}(x)=\partial_A(1)\cup x.
\end{equation}
Il suffit donc de montrer que le diagramme 
\begin{equation}\label{higgs1-kun2h}
\xymatrix{
{\rH^1(G,A)\otimes_A\rH^{n-1}(H,M)}\ar[d]_{\cup}\ar[rd]^{\times}&\\
{\rH^1(G,\rH^{n-1}(H,M))}\ar[r]^{v_n}&{\rH^n(G\times H,M)}}
\end{equation}
où le cup-produit de gauche est défini dans \eqref{higgs1-not583c} et le cross-produit de droite est défini dans \eqref{higgs1-not582h}
est commutatif. D'après \ref{higgs1-cg01}, on a un isomorphisme canonique 
\begin{equation}\label{higgs1-kun2i}
\partial_{\mZ_p}\colon \mZ_p\stackrel{\sim}{\rightarrow} \rH^1(G,\mZ_p),
\end{equation}
compatible avec l'isomorphisme $\partial_A$ \eqref{higgs1-kun1g} via l'homomorphisme canonique 
$\mZ_p\rightarrow A$. Il suffit encore de montrer que le diagramme
\begin{equation}\label{higgs1-kun2j}
\xymatrix{
{\rH^1(G,\mZ_p)\otimes_{\mZ_p}\rH^{n-1}(H,M)}\ar[d]_{\cup}\ar[rd]^{\times}&\\
{\rH^1(G,\rH^{n-1}(H,M))}\ar[r]^{v_n}&{\rH^n(G\times H,M)}}
\end{equation}
est commutatif. D'après \ref{higgs1-kun7}, on a un isomorphisme canonique de $\bD^+(\bMod(A))$
\begin{equation}\label{higgs1-kun2k}
\rR\Gamma(G,\tau_{\leq n-1}\rC^\bullet(H,M))\stackrel{\sim}{\rightarrow} \rC^\bullet(G,\mZ_p)\otimes_{\mZ_p}\tau_{\leq n-1}\rC^\bullet(H,M),
\end{equation}
et compte tenu de \eqref{higgs1-kun2c}, le morphisme $v_n$ s'obtient en appliquant le foncteur $\rH^n$ au morphisme 
\begin{equation}\label{higgs1-kun2l}
\rC^\bullet(G,\mZ_p)\otimes_{\mZ_p}\tau_{\leq n-1}\rC^\bullet(H,M)\rightarrow 
\rC^\bullet(G\times H,M)
\end{equation}
déduit de \eqref{higgs1-not582g}~; d'où la commutativité de \eqref{higgs1-kun2j}.

\begin{rema}
La proposition \ref{higgs1-kun2} fournit une ``formule de Künneth'' 
pour la cohomologie de produits de groupes profinis à valeurs dans certains 
modules discrets que nous n'avons pu trouver dans la littérature.  
Nous n'avons traité qu'un cas simple, le seul nécessaire pour la suite de ce travail, où l'un des groupes est $\mZ_p$. 
Pour le cas général, on trouvera dans \cite{jannsen2} un énoncé plus faible dans lequel 
la compatibilité avec le cross-produit n'est pas abordée.
\end{rema}

\begin{prop}\label{higgs1-kun20}
Soient $A$ une $\mZ_p$-algèbre complète et séparée pour la topologie $p$-adique, $d$ un entier $\geq 1$. 
Il existe alors un et un unique isomorphisme de $A$-algèbres graduées
\begin{equation}
\wedge (\rH^1_\cont(\mZ_p^d,A))\stackrel{\sim}{\rightarrow} \oplus_{n\geq 0}\rH^n_\cont(\mZ_p^d,A),
\end{equation}
le but étant muni du cup-produit, dont la composante de degré $1$ est l'identité. 
\end{prop}
Considérons d'abord le cas où $A$ est de torsion $p$-primaire. 
Procédons par récurrence sur $d$. Le cas $d=1$ est immédiat compte tenu de \ref{higgs1-cg05}
et de l'isomorphisme canonique $A\stackrel{\sim}{\rightarrow} \rH^1_\cont(\mZ_p,A)$. 
Supposons $d\geq 2$ et l'énoncé établi pour $d-1$. 
Posons $G=\mZ_p$ et $H=\mZ_p^{d-1}$. Il résulte de \ref{higgs1-kun2} et \ref{higgs1-not5831} que
le cross-produit \eqref{higgs1-not582h} induit un isomorphisme de $A$-algèbres graduées
\begin{equation}
(A \oplus \rH^1(G,A)){^g \otimes}_A (\oplus_{n\geq 0}\rH^n(H,A)) \stackrel{\sim}{\rightarrow} 
\oplus_{n\geq 0}\rH^n(G\times H,A),
\end{equation}
où le symbole ${^g\otimes}_A$ désigne le produit tensoriel gauche (cf. \cite{alg1-3} III § 4.7 remarques page 49). 
L'assertion recherchée s'ensuit compte tenu de l'hypothèse de récurrence. 
On notera par ailleurs qu'on a un isomorphisme canonique 
$A^d\stackrel{\sim}{\rightarrow} \rH^1_\cont(\mZ_p^d,A)$. 

Considérons ensuite le cas général. Il résulte de ce qui précède que pour tout entier $n\geq 0$, 
le système projectif $(\rH^n(\mZ_p^d,A/p^rA))_{r\in\mN}$ vérifie la condition de Mittag-Leffler.  
En vertu de \eqref{higgs1-limproj2e}, \eqref{higgs1-limproj2c} et \eqref{higgs1-limproj2d}, le morphisme canonique 
\begin{equation}
\rH^n_\cont(\mZ_p^d,A) \rightarrow \underset{\underset{r\geq 0}{\longleftarrow}}{\lim}\ 
\rH^n(\mZ_p^d,A/p^rA)
\end{equation}
est donc un isomorphisme. La proposition se déduit alors du cas déjà démontré.

\subsection{}\label{higgs1-cnab1}
Soit $G$ un groupe profini. Un {\em $G$-ensemble} est un espace topologique discret 
muni d'une action continue de $G$. Les $G$-ensembles forment naturellement une catégorie. 
Un $G$-groupe $M$ est un groupe de cette catégorie. On lui associe son sous-groupe des $G$-invariants
$\rH^0(G,M)=M^G$ et son premier ensemble de cohomologie
$\rH^1(G,M)$; on renvoie à (\cite{serre1} I §5) pour la définition et les principales propriétés de cet ensemble pointé. 

\subsection{}\label{higgs1-cnab2}
Soient $G$ un groupe profini, $A$ un anneau muni de la topologie discrète et d'une action continue de $G$, 
$M$ un $A$-module libre de rang $r\geq 1$, muni de la topologie discrète, $(e_1,\dots,e_r)$ une base de $M$ sur $A$. 
On note $\Mat_r(A)$ la $A$-algèbre des matrices 
carrées de taille $r$ à coefficients dans $A$ et $\GL_r(A)$ le groupe des éléments inversibles de $\Mat_r(A)$. 
On observera que $\GL_r(A)$ est naturellement un $G$-groupe.  
La donnée d'une $A$-représentation continue $\rho$ de $G$ sur $M$ est équivalente à la 
donnée pour tout $g\in G$ d'un élément $U_g$ de $\GL_r(A)$ tel que 
l'application $g\mapsto U_g$ soit continue et que pour tout $g,h\in G$, on ait
\begin{equation}\label{higgs1-cnab2a}
U_{gh}=U_g\cdot {^gU_h}.
\end{equation}
La matrice $U_g$ exprime alors les coordonnées des vecteurs $e_1,\dots,e_r$ dans la base $g(e_1),\dots,g(e_r)$. 
Un changement de la base $(e_i)_{1\leq i\leq r}$ change le cocycle $g\mapsto U_g$ en un cocycle cohomologue.
L'application qui a $\rho$ associe la classe $[\rho]$ du cocycle $g\mapsto U_g$ dans $\rH^1(G,\GL_r(A))$  
est une bijection de l'ensemble des classes d'isomorphismes de $A$-représentations continues de $G$ sur $M$,
dans l'ensemble $\rH^1(G,\GL_r(A))$. La $A$-représentation de $G$ sur $M$ qui fixe les $(e_i)_{1\leq i\leq r}$ 
correspond à l'élément distingué de $\rH^1(G,\GL_r(A))$.

Soient $a\in A^G$ tel que $a$ soit nilpotent dans $A$, $\rho$ une $A$-représentation de $G$ sur $M$
telle que $\rho(g)(e_i)-e_i\in aM$ pour tous $g\in G$ et $1\leq i\leq r$. Alors le cocycle $g\mapsto U_g$
défini ci-dessus prend ses valeurs dans le sous-groupe $\id_r+a\Mat_r(A)$ de $\GL_r(A)$. 
Si on change la base $(e_i)_{1\leq i\leq r}$ en une base $(e'_i)_{1\leq i\leq r}$ telle que $e'_i-e_i\in aM$
pour tout $1\leq i\leq r$, le cocycle $g\mapsto U_g$ se transforme en un cocycle cohomologue dans $\id_r+a\Mat_r(A)$. 
L'application qui à $\rho$ associe la classe $[\rho]$ du cocycle $g\mapsto U_g$ dans $\rH^1(G,\id_r+a\Mat_r(A))$  
est une bijection de l'ensemble des classes d'isomorphismes de $A$-représentations continues de $G$ sur $M$
qui fixent $e_i$ modulo $aM$ pour tout $1\leq i\leq r$, modulo les isomorphismes qui fixent 
$e_i$ modulo $aM$ pour tout $1\leq i\leq r$, dans l'ensemble $\rH^1(G,\id_r+a\Mat_r(A))$.

\subsection{}\label{higgs1-cnab3}
Soient $G$ un groupe profini, $A$ un anneau muni d'une action de $G$, $a\in A^G$, 
$m,n,q,r$ des entiers $\geq 1$ tels que $q\geq n\geq m$  et $m+n\geq q$. Supposons que l'action de $G$ sur $A$
soit continue pour la topologie $a$-adique et que la multiplication par $a^n$ dans $A$ induise un isomorphisme  
\begin{equation}\label{higgs1-cnab3aa}
A/a^{q-n}A\stackrel{\sim}{\rightarrow} a^nA/a^qA.
\end{equation}
La seconde condition est remplie par exemple si $a$ n'est pas un diviseur dans $A$. 
Considérons la suite exacte canonique de $G$-groupes 
\begin{equation}\label{higgs1-cnab3a}
1\rightarrow \id_r+a^n\Mat_r(A/a^qA)\rightarrow \id_r+a^m\Mat_r(A/a^qA) \rightarrow 
\id_r+a^m\Mat_r(A/a^{n}A)\rightarrow 1.
\end{equation}
Alors $\id_r+a^n\Mat_r(A/a^qA)$ est contenu dans le centre de $\id_r+a^m\Mat_r(A/a^qA)$.
D'autre part, d'après \eqref{higgs1-cnab3aa}, on a un isomorphisme canonique de $G$-groupes abéliens
\begin{equation}\label{higgs1-cnab3b}
\Mat_r(A/a^{q-n}A)\stackrel{\sim}{\rightarrow}\id_r+a^n\Mat_r(A/a^qA).
\end{equation}
D'après (\cite{serre1} I prop.~43), on a une suite exacte canonique d'ensembles pointés
\begin{eqnarray}\label{higgs1-cnab3d}
\lefteqn{\rH^1(G,\Mat_r(A/a^{q-n}A))\longrightarrow \rH^1(G, \id_r+a^m\Mat_r(A/a^qA)) \longrightarrow}\\
&&\rH^1(G,  \id_r+a^m\Mat_r(A/a^{n}A))\stackrel{\partial}{\longrightarrow} \rH^2(G,\Mat_r(A/a^{q-n}A)).\nonumber
\end{eqnarray}

\subsection{}\label{higgs1-cnab4} 
Conservons les hypothèses de \eqref{higgs1-cnab3}, soient de plus $N$ un $(A/a^qA)$-module libre de rang $r$,
$f_1,\dots,f_r$ une base de $N$, $\rho,\rho'$ deux représentations continues de $G$ sur $N$ (muni de la topologie
discrète) telles que $\rho(g)(f_i)-f_i\in a^mN$ et $\rho'(g)(f_i)-f_i\in a^mN$ pour tous $g\in G$ et $1\leq i\leq r$. 
Notons $g\mapsto V_g$ et $g\mapsto V'_g$ les cocycles de $G$ 
à valeurs dans $\id_r+a^m\Mat_r(A/a^qA)$ associés à $N$ et $N'$, respectivement, 
par le choix de la base $(f_i)_{1\leq i\leq d}$. 
En vertu de \eqref{higgs1-cnab3d} et (\cite{serre1} I prop.~42),  les conditions suivantes sont équivalentes~:
\begin{itemize}
\item[{\rm (i)}] Il existe un automorphisme $A$-linéaire  
\begin{equation}\label{higgs1-cnab4a}
u\colon N/a^nN\stackrel{\sim}{\rightarrow}N/a^nN
\end{equation}
tel que $u\circ \rho(g)=\rho'(g)\circ u$ pour tout $g\in G$ et que $u(f_i)-f_i\in a^mN$ pour tout $1\leq i\leq r$. 
\item[{\rm (ii)}] Il existe un cocycle $g\mapsto W_g$ de $G$ à valeurs dans $\Mat_r(A/a^{q-n}A)$
et une matrice $U\in \id_r+a^m\Mat_r(A/a^qA)$ tels que pour tout $g\in G$, on ait  \eqref{higgs1-cnab3b}
\begin{equation}\label{higgs1-cnab4b}
V'_g=U^{-1} (\id_r+a^nW_g)V_g {^gU}.
\end{equation} 
\end{itemize}
Nous dirons alors que $\rho'$ est déduit de $\rho$ par torsion par le cocycle $g\mapsto W_g$. 
Deux cocycles cohomologues définissent des représentations isomorphes
par un isomorphisme compatible à \eqref{higgs1-cnab4a}. Nous dirons alors aussi que 
$\rho'$ est déduit de $\rho$ par torsion par la classe $\fc\in \rH^1(G,\Mat_r(A/a^{q-n}A))$ 
du cocycle $g\mapsto W_g$. 

Soient $n',m', q'$ des entiers $\geq 1$ tels que $n\geq n'\geq m'$, $m\geq m'\geq q-n$ et $q'=q-n+n'$.
Supposons que la multiplication par $a^{n'}$ dans $A$ induise un isomorphisme  
\begin{equation}\label{higgs1-cnab4c}
A/a^{q'-n'}A\stackrel{\sim}{\rightarrow} a^{n'}A/a^{q'}A.
\end{equation}
On a alors un diagramme commutatif 
\begin{equation}\label{higgs1-cnab4d}
\xymatrix{
{\Mat_r(A/a^{q-n}A)}\ar[r]\ar[d]_{\cdot a^{n-n'}}&{\id_r+a^m\Mat_r(A/a^qA)}\ar[r]\ar[d]^\alpha&
{\id_r+a^m\Mat_r(A/a^nA)}\ar[d]^\beta\\
{\Mat_r(A/a^{q'-n'}A)}\ar[r]&{\id_r+a^{m'}\Mat_r(A/a^{q'}A)}\ar[r]&{\id_r+a^{m'}\Mat_r(A/a^{n'}A)}}
\end{equation}
où les lignes correspondent aux suites exactes \eqref{higgs1-cnab3a} et $\alpha$ et $\beta$ sont induits par
les homomorphismes de réduction $A/a^qA\rightarrow A/a^{q'}A$ et $A/a^nA\rightarrow A/a^{n'}A$. 
Par suite, si $\rho'$ se déduit de $\rho$ par torsion par une classe $\fc\in \rH^1(G,\Mat_r(A/a^{q-n}A))$, 
la représentation induite par $\rho'$ sur $N/a^{q'} N$ se déduit de 
celle induite par $\rho$ sur $N/a^{q'} N$ par torsion par la classe 
\begin{equation}\label{higgs1-cnab4e}
a^{n-n'}\cdot\fc\in \rH^1(G,\Mat_r(A/a^{q'-n'}A)).
\end{equation}

\subsection{}\label{higgs1-cnab5}
Conservons les hypothèses de \eqref{higgs1-cnab3}, soient de plus 
$M$ un $(A/a^nA)$-module libre de rang $r$, $e_1,\dots,e_r$ une base de $M$, 
$\rho$ une représentation continue de $G$ sur $M$ (muni de la topologie discrète)
telle que $\rho(g)(e_i)-e_i\in a^mM$ pour tous $g\in G$ et tous $1\leq i\leq r$.  On note $g\mapsto U_g$ le cocycle de $G$ 
à valeurs dans $\id_r+a^m\Mat_r(A/a^nA)$ associé à $M$ par le choix de la base $(e_i)_{1\leq i\leq d}$,
et $[M]$ sa classe dans $\rH^1(G,  \id_r+a^m\Mat_r(A/a^{n}A))$.
En vertu de \eqref{higgs1-cnab3d} et (\cite{serre1} I prop.~41), la classe 
\begin{equation}\label{higgs1-cnab5b}
\partial([M])\in \rH^2(G,\Mat_r(A/a^{q-n}A))
\end{equation}
est l'obstruction à relever $M$ en une $(A/a^qA)$-représentation continue de $G$ dont le 
$(A/a^qA)$-module sous-jacent est libre de type fini. 

Soient $n'$ un entier tel que $n\geq n'\geq m$, $q'=q-n+n'$ et que la multiplication par $a^{n'}$ dans $A$ induise un isomorphisme  
\begin{equation}\label{higgs1-cnab5bc}
A/a^{q'-n'}A\stackrel{\sim}{\rightarrow} a^{n'}A/a^{q'}A.
\end{equation}
On note 
\begin{equation}\label{higgs1-cnab5c}
\partial'\colon \rH^1(G,  \id_r+a^m\Mat_r(A/a^{n'}A))\rightarrow \rH^2(G,\End_A(M/a^{q'-n'}M))
\end{equation}
l'application bord définie comme dans \eqref{higgs1-cnab3d}. Par fonctorialité, on a 
\begin{equation}\label{higgs1-cnab5d}
\partial'([M/a^{n'}M])=a^{n-n'}\partial([M]). 
\end{equation}

\section{Sorites sur les objets à groupe d'opérateurs}\label{higgs1-eph}

\subsection{}\label{higgs1-eph1}\index{10410@$\Sch$}
On note $\Sch$ la catégorie des schémas. On choisit un clivage normalisé de la catégorie des flèches de $\Sch$
 (\cite{sga1} VI §11), autrement dit, pour tout morphisme de schémas $f\colon X'\rightarrow X$, on choisit 
un foncteur de changement de base
\begin{equation}\label{higgs1-eph1b}
f^\bullet \colon \Sch_{/X}\rightarrow \Sch_{/X'}, \ \ \ Y\mapsto Y\times_XX'
\end{equation}
tel que pour tout schéma $X$, $f=\id_X$ implique $f^\bullet=\id_{\Sch_{/X}}$. 
Pour tout schéma $X$, on désigne par $X_\zar$ le topos de Zariski de $X$ et par
\begin{equation}\label{higgs1-eph1a}
\varphi_X\colon \Sch_{/X}\rightarrow X_\zar, \ \ \ Y\mapsto \Hom_X(-,Y)
\end{equation}
le foncteur canonique.

\begin{rema}\label{higgs1-eph11}
Si $f\colon X'\rightarrow X$ est un morphisme de schémas, le diagramme 
\begin{equation}\label{higgs1-eph11a}
\xymatrix{
{\Sch_{/X}}\ar[r]^{\varphi_X}\ar[d]_{f^\bullet}&{X_\zar}\ar[d]^{f^*}\\
{\Sch_{/X'}}\ar[r]^{\varphi_{X'}}&{X'_\zar}}
\end{equation}
où les flèches horizontales sont les foncteurs canoniques \eqref{higgs1-eph1a}, n'est pas commutatif en général. 
Toutefois, on a un morphisme canonique de foncteurs
\begin{equation}\label{higgs1-eph11b}
f^*\circ \varphi_X\rightarrow \varphi_{X'}\circ f^\bullet.
\end{equation}
Si, de plus, $f$ est une immersion ouverte, le morphisme \eqref{higgs1-eph11b}
est un isomorphisme. En effet, tout ouvert de $X'$ est un ouvert de $X$, et $f^*$ est le foncteur de restriction.  
Dans ce cas, le diagramme \eqref{higgs1-eph11a} est donc commutatif, à isomorphisme canonique près.
\end{rema}

\subsection{}\label{higgs1-eph2}\index{Fibre principal@$G$-fibré principal homogène}
\index{10430@$\bFPH(G/X)$}\index{10431@$\bTors(G,X_\zar)$}
Soient $X$ un schéma, $G$ un $X$-schéma en groupes. Dans cet article, par 
{\em $G$-fibré principal homogène sur $X$}, on sous-entend un $G$-pseudo-torseur (à droite) de 
$\Sch_{/X}$, localement trivial pour la topologie de Zariski sur $X$ (\cite{giraud2} III 1.1.5). 
On désigne par $\bFPH(G/X)$ la catégorie des $G$-fibrés principaux homogènes sur $X$
et par $\bTors(\varphi_X(G),X_\zar)$ la catégorie des $\varphi_X(G)$-torseurs (à droite) de $X_\zar$ 
(\cite{giraud2} III 1.4.1). Le foncteur $\varphi_X$  \eqref{higgs1-eph1a} induit un foncteur que l'on note aussi 
\begin{equation}\label{higgs1-eph2a}
\varphi_X\colon \bFPH(G/X)\rightarrow \bTors(\varphi_X(G),X_\zar), \ \ \ Y\mapsto \Hom_X(-,Y).
\end{equation}
 
\begin{prop}\label{higgs1-eph3}
Soient $X$ un schéma, $G$ un $X$-schéma en groupes. Alors~:
\begin{itemize}
\item[{\rm (i)}] Le foncteur \eqref{higgs1-eph2a} est pleinement fidèle. 
\item[{\rm (ii)}] Si, de plus, $X$ est cohérent (i.e., quasi-compact et quasi-séparé)
et que $G$ est affine sur $X$, alors le foncteur \eqref{higgs1-eph2a} est une équivalence de catégories. 
\end{itemize}
\end{prop}
(i) Soient $Y,Z$ deux objets de $\bFPH(G/X)$. On note $\Hom_G(Y,Z)$ l'ensemble des morphismes de 
$Y$ dans $Z$ dans $\bFPH(G/X)$ et $\Hom_{\varphi_X(G)}(\varphi_X(Y),\varphi_X(Z))$ 
l'ensemble des morphismes de $\varphi_X(Y)$ dans 
$\varphi_X(Z)$ dans $\bTors(\varphi_X(G),X_\zar)$. Montrons que l'application 
\begin{equation}\label{higgs1-eph3a}
\Hom_G(Y,Z)\rightarrow \Hom_{\varphi_X(G)}(\varphi_X(Y),\varphi_X(Z))
\end{equation}
induite par le foncteur $\varphi_X$ \eqref{higgs1-eph2a} est bijective. 
Soit $(U_i)_{i\in I}$ un recouvrement ouvert de Zariski de $X$.
Pour tout $i\in I$, on pose  $G_i=G\times_XU_i$, $Y_i=Y\times_XU_i$ et $Z_i=Z\times_XU_i$.
Pour tout $(i,j)\in I^2$, on pose $U_{ij}=U_i\cap U_j$,
$G_{ij}=G\times_XU_{ij}$, $Y_{ij}=Y\times_XU_{ij}$ et $Z_{ij}=Z\times_XU_{ij}$.
On a alors un diagramme commutatif d'applications d'ensembles
\begin{equation}\label{higgs1-eph3b}
\xymatrix{
{\Hom_G(Y,Z)}\ar[d]\ar[r]&{\Hom_{\varphi_X(G)}(\varphi_X(Y),\varphi_X(Z))}\ar[d]\\
{\prod_{i\in I}\Hom_{G_i}(Y_i,Z_i)}\ar@<0.5ex>[d]
\ar@<-0.5ex>[d]\ar[r]&{\prod_{i\in I}\Hom_{\varphi_{X_i}(G_i)}(\varphi_{X_i}(Y_i),\varphi_{X_i}(Z_i))}\ar@<0.5ex>[d]
\ar@<-0.5ex>[d]\\
{\prod_{(i,j)\in I^2}\Hom_{G_{ij}}(Y_{ij},Z_{ij})}\ar[r]&
{\prod_{(i,j)\in I^2}\Hom_{\varphi_{X_{ij}}(G_{ij})}(\varphi_{X_{ij}}(Y_{ij}),\varphi_{X_{ij}}(Z_{ij}))}}
\end{equation}
où les flèches horizontales sont induites par le foncteur \eqref{higgs1-eph2a} et les flèches verticales sont définies
par restriction \eqref{higgs1-eph11}. Les colonnes verticales étant clairement exactes,   
on peut se borner au cas où $Y=G$. La source et le but de 
l'application \eqref{higgs1-eph3a} s'identifient alors à l'ensemble $\Hom_X(X,Z)$, d'où l'assertion. 

(ii) Il suffit de montrer que le foncteur en question est essentiellement surjectif. 
Soient $F$ un objet de $\bTors(\varphi_X(G),X_\zar)$, $(U_i)_{i\in I}$ un recouvrement de $X$
par des ouverts affines tel que pour tout $i\in I$, $F|U_i$ soit trivial. 
Posons $X'=\sqcup_{i\in I}U_i$, $X''=X'\times_XX'$, $G'=G\times_XX'$
et $G''=G\times_XX''$ et  notons $f\colon X'\rightarrow X$ le morphisme canonique
et $\pr_1,\pr_2\colon X''\rightarrow X'$ les projections canoniques. 
Comme $X$ est quasi-compact, on peut supposer $I$ fini. 
Comme $X$ est quasi-séparé, $f$ est alors quasi-compact. 
On a un isomorphisme $\varphi_{X'}(G')\stackrel{\sim}{\rightarrow}f^*(F)$ de $\bTors(\varphi_{X'}(G'),X'_\zar)$. 
La donnée de descente canonique sur $f^*(F)$ induit une donnée de descente sur $\varphi_{X'}(G')$ relativement à $f$
(en tant que $\varphi_{X'}(G')$-torseur de $X'_\zar$), c'est-à-dire un isomorphisme de $\bTors(\varphi_{X''}(G''),X''_\zar)$
\begin{equation}
\psi\colon \pr_1^*(\varphi_{X'}(G'))\stackrel{\sim}{\rightarrow}\pr_2^*(\varphi_{X'}(G'))
\end{equation}
vérifiant une relation de cocycles. On notera que $\pr_1^*(\varphi_{X'}(G'))$ et $\pr_2^*(\varphi_{X'}(G'))$
sont canoniquement isomorphes à $\varphi_{X''}(G'')$ \eqref{higgs1-eph11} et que $\psi$ est en général différent
de la donnée de descente triviale (induite par $\varphi_{X}(G)$). Compte tenu de (i), $\psi$ induit une donnée de 
descente sur $G'$ relativement à $f$ (en tant que $G'$-fibré principal homogène sur $X'$), c'est-à-dire
un isomorphisme de $\bFPH(G''/X'')$
\begin{equation}
\varphi\colon\pr_1^\bullet(G')\stackrel{\sim}{\rightarrow}\pr_2^\bullet(G'),
\end{equation}
vérifiant une relation de cocycles. En vertu de (\cite{sga1} VIII 2.1), 
il existe un $G$-fibré principal homogène $Y$ sur $X$ qui correspond à $\varphi$. 
Comme le $G$-torseur $\varphi_{X}(Y)$ de $X_\zar$ correspond à la donnée de descente $\psi$,
il existe alors un isomorphisme  
$\varphi_{X}(Y)\stackrel{\sim}{\rightarrow} F$ de $\bTors(\varphi_{X}(G),X_\zar)$, d'où l'assertion.

\subsection{}\label{higgs1-eph150}
Soient $f\colon X'\rightarrow X$ un morphisme de schémas, $\rT$ un $\co_X$-module,
$\cL$ un $\rT$-torseur de $X_\zar$.  Dans la suite, nous utilisons pour les $\co_X$-modules la notation 
$f^{-1}$ pour désigner l'image inverse au sens des faisceaux abéliens et nous réservons la notation 
$f^*$ pour l'image inverse au sens des modules. On appelle {\em image inverse affine} de $\cL$ par $f$,
et l'on note $f^+(\cL)$, le $f^*(\rT)$-torseur de $X'_\zar$ déduit du $f^{-1}(\rT)$-torseur $f^*(\cL)$
par extension de son groupe structural par l'homomorphisme canonique $f^{-1}(\rT)\rightarrow f^*(\rT)$~:
\begin{equation}\label{higgs1-eph150a}
f^+(\cL)=f^*(\cL)\wedge^{f^{-1}(\rT)}f^*(\rT),
\end{equation}
autrement dit, le quotient de $f^*(\cL)\times f^*(\rT)$ par l'action diagonale de $f^{-1}(\rT)$ (\cite{giraud2} III 1.4.6). 
Soient $\rT'$ un $\co_X$-module, $\cL'$ un $\rT'$-torseur de $X_\zar$, 
$u\colon \rT\rightarrow \rT'$ un morphisme $\co_X$-linéaire, $v\colon \cL\rightarrow \cL'$ un morphisme
$u$-équivariant de $X_\zar$. D'après (\cite{giraud2} III 1.3.6), il existe un et un unique 
morphisme $f^*(u)$-équivariant  
\begin{equation}\label{higgs1-eph150c}
f^+(v)\colon f^+(\cL)\rightarrow f^+(\cL')
\end{equation}
qui s'insère dans le diagramme commutatif 
\begin{equation}\label{higgs1-eph150b}
\xymatrix{
{f^*(\cL)}\ar[r]^-(0.5){f^*(v)}\ar[d]&{f^*(\cL')}\ar[d]\\
{f^+(\cL)}\ar[r]^-(0.5){f^+(v)}&{f^+(\cL')}}
\end{equation}
où les flèches verticales sont les morphismes canoniques. 
La correspondance ainsi définie $(\rT,\cL) \mapsto (f^*(\rT),f^+(\cL))$ est un foncteur de la catégorie des 
torseurs de $X_\zar$ sous un $\co_X$-module dans la catégorie des torseurs de $X'_\zar$ sous un $\co_{X'}$-module. 

Soient $g\colon X''\rightarrow X'$ un morphisme de schémas, $h=f\circ g\colon X''\rightarrow X$. On a un isomorphisme canonique de foncteurs 
\begin{equation}\label{higgs1-eph150d}
h^*\stackrel{\sim}{\rightarrow}g^*\circ f^*.
\end{equation}
Comme $g^*$ commute aux limites inductives, on a un isomorphisme canonique 
\begin{equation}\label{higgs1-eph150e}
g^*(f^+(\cL))\stackrel{\sim}{\rightarrow}g^*(f^*(\cL))\wedge^{h^{-1}(\rT)}g^{-1}(f^*(\rT)).
\end{equation}
Compte tenu de (\cite{giraud2} III 1.3.5), celui-ci induit un isomorphisme $h^*(\rT)$-équivariant canonique 
\begin{equation}\label{higgs1-eph150f}
g^+(f^+(\cL))\stackrel{\sim}{\rightarrow}h^+(\cL). 
\end{equation}
On vérifie aussitôt que cet isomorphisme est fonctoriel et qu'il vérifie une relation de cocycles du type 
(\cite{sga1} VI 7.4 B)) pour la composition de trois morphismes de schémas. 

\subsection{}\label{higgs1-eph151}
Soient $f\colon X'\rightarrow X$ un morphisme de schémas, $\rT$ un $\co_X$-module localement libre de type fini,
$\Omega=\cHom_{\co_X}(\rT,\co_X)$ son dual, $\bT=\Spec(\rS_{\co_X}(\Omega))$ le $X$-fibré correspondant,
$\bL$ un $\bT$-fibré principal homogène sur $X$. On a un diagramme commutatif de morphismes de groupes de $X'_\zar$
\begin{equation}\label{higgs1-eph151a}
\xymatrix{
{f^{-1}(\rT)}\ar[r]^-(0.5)\sim\ar[d]&{f^*(\varphi_X(\bT))}\ar[d]\\
{f^*(\rT)}\ar[r]^-(0.5)\sim&{\varphi_{X'}(\bT\times_XX')}}
\end{equation}
où la flèche verticale de droite est le morphisme \eqref{higgs1-eph11b} et 
les autres flèches sont les morphismes canoniques. 
Par ailleurs, on a un morphisme canonique $f^{-1}(\rT)$-équivariant \eqref{higgs1-eph11b}  
\begin{equation}\label{higgs1-eph151b}
f^*(\varphi_X(\bL))\rightarrow\varphi_{X'}(\bL\times_XX').
\end{equation}
D'après (\cite{giraud2} III 1.3.6), celui-ci induit un isomorphisme de $f^*(\rT)$-torseurs 
\begin{equation}\label{higgs1-eph151c}
f^+(\varphi_X(\bL))\stackrel{\sim}{\rightarrow}\varphi_{X'}(\bL\times_XX').
\end{equation}

\subsection{}\label{higgs1-eph10}\index{Fonction affine@Fonction affine sur un torseur} 
Soient $X$ un schéma, $\rT$ un $\co_X$-module localement libre de type fini, 
$\Omega=\cHom_{\co_X}(\rT,\co_X)$ son dual, $\cL$ un $\rT$-torseur  de $X_\zar$. 
On appelle {\em fonction affine sur $\cL$} la donnée d'un morphisme $f\colon \cL\rightarrow \co_X$ de $X_\zar$
vérifiant les conditions équivalentes suivantes~:
\begin{itemize}
\item[(i)] pour tout ouvert $U$ de $X$ et tout $s\in \cL(U)$, l'application 
\begin{equation}\label{higgs1-eph10b}
\rT(U)\rightarrow \co_X(U), \ \ \ t\mapsto f(s+t)-f(s)
\end{equation}
est $\co_X(U)$-linéaire. 
\item[(ii)]  il existe une section $\omega\in \Omega(X)$, dite {\em terme linéaire} de $f$,  
telle que pour tout ouvert $U$ de $X$ et tous $s\in \cL(U)$ et $t\in \rT(U)$, on ait
\begin{equation}\label{higgs1-eph10c}
f(s+t)=f(s)+\omega(t).
\end{equation}
\end{itemize}
On dit alors aussi que le morphisme $f$ est {\em affine}.
En effet, (ii) implique clairement (i). Inversement, supposons la condition (i) remplie. 
Soit $(U_i)_{i\in I}$ un recouvrement ouvert 
de $X$ tel que pour tout $i\in I$, il existe $s_i\in \cL(U_i)$. Il existe alors $\omega_i\in \Omega(U_i)$ tel que 
pour tout ouvert $U$ de $U_i$ et tout $t\in \rT(U)$, on ait 
\begin{equation}\label{higgs1-eph10d}
f(s_i+t)=f(s_i)+\omega_i(t)\in \co_X(U).
\end{equation}
Pour tout $(i,j)\in I^2$, soit $t_{ij}\in \rT(U_i\cap U_j)$ tel que $s_i=s_j+t_{ij}\in \cL(U_i\cap U_j)$. 
Pour tout ouvert $U$ de $U_i\cap U_j$ et tout $t\in \rT(U)$, on a
\begin{eqnarray}
f(s_i)+\omega_i(t)&=&f(s_i+t)=f(s_j+t_{ij}+t)\nonumber\\
&=&f(s_j)+\omega_j(t_{ij})+\omega_j(t)=f(s_i)+\omega_j(t)\in \co_X(U).\label{higgs1-eph10e}
\end{eqnarray}
On a donc $\omega_i|U_i\cap U_j=\omega_j|U_i\cap U_j$. Par suite, les sections $(\omega_i)_{i\in I}$ se recollent 
et définissent une section $\omega\in \Omega(X)$. Soient $U$ un ouvert de $X$, $s\in \cL(U)$, $t\in \rT(U)$.
Montrons que $f(s+t)=f(s)+\omega(t)$.  
On peut supposer qu'il existe $i\in I$ tel que $U\subset U_i$. Soit $t'\in \rT(U)$ tel que $s=s_i+t'$. On a alors
\begin{equation}\label{higgs1-eph10f}
f(s+t)=f(s_i+t'+t)=f(s_i)+\omega(t')+\omega(t)=f(s)+\omega(t).
\end{equation}

\begin{rema}\label{higgs1-eph101}
La condition \ref{higgs1-eph10}(ii) revient à dire qu'il existe $\omega\in \Omega(X)$ tel que 
$f$ soit $\rT$-équivariant lorsque l'on munit $\co_X$ de la structure de $\rT$-objet définie par $\omega$, 
plus précisément, par le morphisme 
\begin{equation}
\rT\times \co_X\rightarrow \co_X,\ \ \ (t,x)\mapsto \omega(t)+x.
\end{equation}
\end{rema}

\subsection{}\label{higgs1-eph13}
Soient $X$ un schéma, $\rT$ un $\co_X$-module localement libre de type fini, 
$\Omega=\cHom_{\co_X}(\rT,\co_X)$ son dual, $\cL$ un $\rT$-torseur  de $X_\zar$. 
La condition \ref{higgs1-eph10}(i) étant clairement locale pour la topologie de Zariski sur $X$, 
on désigne par $\cF$ le sous-faisceau de $\cHom_{X_\zar}(\cL,\co_X)$ formé des fonctions affines sur $\cL$, 
autrement dit, pour tout ouvert $U$ de $X$, $\cF(U)$ est l'ensemble des fonctions affines sur $\cL|U$.
On appelle $\cF$ le faisceau des {\em fonctions affines sur $\cL$}.
Il est naturellement muni d'une structure de $\co_X$-module. 
On a un morphisme $\co_X$-linéaire canonique $c\colon \co_X\rightarrow \cF$
dont l'image est formée des fonctions constantes. Le ``terme linéaire'' définit un morphisme $\co_X$-linéaire
$\nu \colon \cF\rightarrow \Omega$. On vérifie que la suite 
\begin{equation}\label{higgs1-eph13b}
0\longrightarrow \co_X\stackrel{c}{\longrightarrow}\cF\stackrel{\nu}{\longrightarrow} \Omega\longrightarrow 0
\end{equation}
est exacte~; pour ce faire, on peut supposer $\cL$ trivial. 
D'après (\cite{illusie1} I 4.3.1.7), cette suite induit pour tout entier $n\geq 1$, 
une suite exacte  \eqref{higgs1-not621}
\begin{equation}\label{higgs1-eph13d}
0\rightarrow \rS^{n-1}_{\co_X}(\cF)\rightarrow \rS^{n}_{\co_X}(\cF)\rightarrow \rS^n_{\co_X}(\Omega)\rightarrow 0.
\end{equation}
Les $\co_X$-modules $(\rS^{n}_{\co_X}(\cF))_{n\in \mN}$ forment donc un système inductif filtrant, 
dont la limite inductive 
\begin{equation}\label{higgs1-eph13c}
\cC=\underset{\underset{n\geq 0}{\longrightarrow}}\lim\ \rS^n_{\co_X}(\cF)
\end{equation}
est naturellement munie d'une structure de $\co_X$-algèbre. 
Pour tout entier $n\geq 0$, le morphisme canonique $\rS^{n}_{\co_X}(\cF)\rightarrow \cC$ est injectif. 
On notera que pour toute $\co_X$-algèbre $\cB$, si $u\colon \cF\rightarrow \cB$ est un morphisme $\co_X$-linéaire 
tel que $u\circ c$ soit l'homomorphisme structural,
il existe un et un unique homomorphisme de $\co_X$-algèbres $\cC\rightarrow \cB$ qui prolonge $u$. 

Il existe un et un unique homomorphisme de $\co_X$-algèbres 
\begin{equation}\label{higgs1-eph13e}
\mu\colon \cC \rightarrow \rS_{\co_X}(\Omega)\otimes_{\co_X}\cC 
\end{equation}
tel que pour toute section locale $x$ de $\cF$, on ait
\begin{equation}\label{higgs1-eph13f}
\mu(x)=\nu(x)\otimes 1+1\otimes x.
\end{equation}

Notons $\bT=\Spec(\rS_{\co_X}(\Omega))$ le $X$-fibré vectoriel associé à $\Omega$ et posons 
\begin{equation}\label{higgs1-eph13g}
\bL=\Spec(\cC).
\end{equation} 
On a un isomorphisme canonique de groupes de $X_\zar$
\begin{equation}\label{higgs1-eph13j}
\rT\stackrel{\sim}{\rightarrow}\varphi_X(\bT).
\end{equation}
Pour tout $s\in \cL(X)$, le morphisme $\rho_s\colon \cF\rightarrow \co_X$ qui à toute section locale $f$ de $\cF$
associe $f(s)$, est un scindage de la suite exacte \eqref{higgs1-eph13b}. Celui-ci se prolonge
en un et un unique homomorphisme de $\co_X$-algèbres
$\varrho_s\colon\cC\rightarrow \co_X$, qui induit une section $\sigma_s\in \bL(X)$. 
La correspondance $s\mapsto \sigma_s$ définit un morphisme de $X_\zar$ 
\begin{equation}\label{higgs1-eph13i}
\iota\colon \cL\rightarrow \varphi_X(\bL).
\end{equation}

\begin{prop}\label{higgs1-eph12}
Sous les hypothèses de \eqref{higgs1-eph13}, le morphisme $\bT\times_X\bL\rightarrow \bL$ induit par $\mu$ \eqref{higgs1-eph13e}
fait de $\bL$ un $\bT$-fibré principal homogène sur $X$, et le morphisme canonique 
$\iota\colon\cL\rightarrow \varphi_X(\bL)$ \eqref{higgs1-eph13i} est un isomorphisme de $\rT$-torseurs. 
\end{prop}

Les questions étant locales, on peut se borner au cas où $\cL$ est trivial. Soient $s\in \cL(X)$, 
$\rho_s\colon \cF\rightarrow \co_X$ le scindage associé de la suite exacte \eqref{higgs1-eph13b}. 
Le morphisme $\lambda_s\colon \Omega\rightarrow \cF$ déduit de $\id_\cF-c\circ \rho_s$
se prolonge en un homomorphisme de $\co_X$-algèbres
\begin{equation}
\psi\colon \rS_{\co_X}(\Omega)\rightarrow \cC
\end{equation}
compatible avec les filtrations $(\oplus_{0\leq i\leq n}\rS^i_{\co_X}(\Omega))_{n\in \mN}$ et 
$(\rS^n_{\co_X}(\cF))_{n\in \mN}$. Il résulte de \eqref{higgs1-eph13d} que $\psi$ est un isomorphisme. 
Notons 
\begin{equation}
\delta\colon \rS_{\co_X}(\Omega)\rightarrow \rS_{\co_X}(\Omega)\otimes_{\co_X}\rS_{\co_X}(\Omega)
\end{equation}
l'homomorphisme de $\co_X$-algèbres tel que pour toute section locale $\omega$ de $\Omega$, on ait
$\delta(\omega)=\omega\otimes 1+1\otimes \omega$. On voit aussitôt que le diagramme
\begin{equation}
\xymatrix{
{\rS_{\co_X}(\Omega)}\ar[r]^\psi\ar[d]_\delta&{\cC}\ar[d]^\mu\\
{\rS_{\co_X}(\Omega)\otimes_{\co_X}\rS_{\co_X}(\Omega)}\ar[r]^-(0.5){\id\otimes \psi}
&{\rS_{\co_X}(\Omega)\otimes_{\co_X}\cC}}
\end{equation}
est commutatif. Par suite, le morphisme $\bT\times_X\bL\rightarrow \bL$ induit par $\mu$ 
fait de $\bL$ un $\bT$-fibré principal homogène sur $X$, et le morphisme 
\begin{equation}
\Psi\colon \bL\rightarrow \bT
\end{equation}
induit par $\psi$ est un isomorphisme de $\bT$-torseurs. Comme $\rho_s\circ \lambda_s=0$, 
$\Psi(\iota(s))=0$ dans $\bT(X)$.
Pour tout $t\in \rT(X)$, on a $\rho_{s+t}=\rho_s+t\circ \nu$ et donc $\rho_{s+t}\circ\lambda_s=t \circ \nu
\circ\lambda_s= t$. On en déduit que $\Psi(\iota(s+t))=t$ dans $\bT(X)$, et par suite que 
\begin{equation}
\iota(s+t)= \iota(s)+t.
\end{equation}
Il s'ensuit que $\iota$ est un morphisme de $\rT$-torseurs et donc un isomorphisme. 

\begin{defi}\label{higgs1-eph121}
Sous les hypothèses de \eqref{higgs1-eph13}, on dit que $\bL$ est le $\bT$-fibré principal homogène canonique 
sur $X$ qui représente $\cL$. 
\end{defi}

D'après \ref{higgs1-eph3}(i), il existe au plus un $\bT$-fibré principal homogène sur $X$ qui représente
$\cL$, à isomorphisme canonique près. La construction du \ref{higgs1-eph13} en fournit un canonique.

\subsection{}\label{higgs1-eph122}
Soient $X$ un schéma, $\rT$ et $\rT'$ deux $\co_X$-modules localement libre de type fini,
$\cL$ un $\rT$-torseur de $X_\zar$, 
$\cL'$ un $\rT'$-torseur de $X_\zar$. 
On pose $\Omega=\cHom_{\co_X}(\rT,\co_X)$, 
$\Omega'=\cHom_{\co_X}(\rT',\co_X)$, $\bT=\Spec(\rS_{\co_X}(\Omega))$ et $\bT'=\Spec(\rS_{\co_X}(\Omega'))$. 
On désigne par $\cF$ le faisceau des fonctions affines sur $\cL$ \eqref{higgs1-eph13}, par 
$\cF'$ le faisceau des fonctions affines sur $\cL'$, par 
$\bL$ le $\bT$-fibré principal homogène canonique sur $X$ qui représente $\cL$ \eqref{higgs1-eph121}
et par $\bL'$ le $\bT'$-fibré principal homogène canonique sur $X$ qui représente $\cL'$. 
Soient $u\colon \rT\rightarrow \rT'$ un morphisme $\co_X$-linéaire, $u^\vee\colon \Omega'\rightarrow \Omega$
le morphisme dual de $u$, $v\colon \cL\rightarrow \cL'$ un morphisme $u$-équivariant de $X_\zar$. 
Si $h\colon \cL'\rightarrow \co_X$ est une fonction affine de terme linéaire $\omega' \in \Omega'(X)$, 
le morphisme composé $h'=h\circ v\colon \cL\rightarrow \co_X$ est affine, de terme linéaire 
$u^\vee(\omega')$. La correspondance $h\mapsto h'$ ainsi définie induit un morphisme $\co_X$-linéaire 
\begin{equation}\label{higgs1-eph122a}
w\colon \cF'\rightarrow \cF
\end{equation}
qui s'insère dans un diagramme commutatif 
\begin{equation}\label{higgs1-eph122b}
\xymatrix{
0\ar[r]&{\co_X}\ar[r]\ar@{=}[d]&{\cF'}\ar[r]\ar[d]^{w}&{\Omega'}\ar[r]\ar[d]^{u^\vee}&0\\
0\ar[r]&{\co_X}\ar[r]&{\cF}\ar[r]&{\Omega}\ar[r]&0}
\end{equation}
où les lignes sont les suites exactes canoniques \eqref{higgs1-eph13b}.

Le morphisme $u^\vee$ induit un morphisme de $X$-schémas en groupes $\alpha\colon \bT\rightarrow \bT'$. 
Le morphisme $w$ induit un $X$-morphisme $\alpha$-équivariant 
\begin{equation}\label{higgs1-eph122c}
\beta\colon \bL\rightarrow \bL'.
\end{equation}
Le diagramme
\begin{equation}\label{higgs1-eph122d}
\xymatrix{
{\cL}\ar[r]^-(0.45){v}\ar[d]_{\iota}&{\cL'}\ar[d]^{\iota'}\\
{\varphi_X(\bL)}\ar[r]^-(0.5){\varphi_X(\beta)}&{\varphi_X(\bL')}}
\end{equation}
où $\iota$ et $\iota'$ sont les isomorphismes canoniques \eqref{higgs1-eph13i}, est clairement commutatif. 

La correspondance qui à un $\rT$-torseur de $X_\zar$ associe 
le $\bT$-fibré principal homogène canonique sur $X$ qui le représente définit donc un foncteur 
\begin{equation}\label{higgs1-eph122e}
\bTors(\rT,X_\zar)\rightarrow \bFPH(\bT/X), \ \ \ \cL\mapsto \bL.
\end{equation}
C'est un quasi-inverse du foncteur \eqref{higgs1-eph2a}
\begin{equation}\label{higgs1-eph122f}
\varphi_X\colon \bFPH(\bT/X)\rightarrow \bTors(\rT,X_\zar), \ \ \ \bL\mapsto \varphi_X(\bL),
\end{equation}
en vertu de \ref{higgs1-eph12}, \ref{higgs1-eph3}(i) et \eqref{higgs1-eph122d}.

\subsection{}\label{higgs1-eph15}
Soient $f\colon X'\rightarrow X$ un morphisme de schémas, $\rT$ un $\co_X$-module localement libre de type fini,
$\cL$ un $\rT$-torseur de $X_\zar$. On pose $\Omega=\cHom_{\co_X}(\rT,\co_X)$ et $\bT=\Spec(\rS_{\co_X}(\Omega))$.
On désigne par $\cF$ le faisceau des fonctions affines sur $\cL$ \eqref{higgs1-eph10},
par $\cF^+$ le faisceau des fonctions affines sur $f^+(\cL)$ \eqref{higgs1-eph150},
par $\bL$ le $\bT$-fibré principal homogène canonique sur $X$ qui représente $\cL$  \eqref{higgs1-eph121}
et par $\bL^+$ le $(\bT\times_XX')$-fibré principal homogène canonique sur $X'$ qui représente $f^+(\cL)$. 
Soient $\ell\colon \cL\rightarrow \co_X$ un morphisme affine, $\omega\in \Omega(X)$ son terme linéaire,
$\omega'=f^*(\omega)\in f^*(\Omega)(X')$. Munissant $\co_{X'}$ de la structure de 
$f^*(\rT)$-objet définie par $\omega'$ \eqref{higgs1-eph101}, il existe un et un unique morphisme $f^*(\rT)$-équivariant 
$\ell'\colon f^+(\cL)\rightarrow \co_{X'}$ qui s'insère dans le diagramme commutatif 
\begin{equation}\label{higgs1-eph15b}
\xymatrix{
{f^*(\cL)}\ar[r]^-(0.5){f^*(\ell)}\ar[d]&{f^{-1}(\co_X)}\ar[d]\\
{f^+(\cL)}\ar[r]^-(0.5){\ell'}&{\co_{X'}}}
\end{equation}
où les flèches verticales sont les morphismes canoniques (\cite{giraud2} III 1.3.6). Le morphisme $h'$ est donc affine,
de terme linéaire $\omega'$. La correspondance $\ell\mapsto \ell'$ ainsi définie induit un morphisme $\co_X$-linéaire
\begin{equation}\label{higgs1-eph15c}
\lambda_\sharp\colon \cF\rightarrow f_*(\cF^+)
\end{equation}
qui s'insère dans un diagramme commutatif 
\begin{equation}\label{higgs1-eph15d}
\xymatrix{
\co_X\ar[r]\ar[d]&\cF\ar[r]\ar[d]^{\lambda_\sharp}&{\Omega}\ar[d]\\
{f_*(\co_{X'})}\ar[r]&{f_*(\cF^+)}\ar[r]&{f_*(f^*(\Omega))}}
\end{equation}
où les autres flèches sont les morphismes canoniques. Le morphisme adjoint 
\begin{equation}\label{higgs1-eph15e}
\lambda\colon f^*(\cF)\rightarrow  \cF^+
\end{equation}
s'insère donc dans un diagramme commutatif 
\begin{equation}\label{higgs1-eph15f}
\xymatrix{
0\ar[r]&\co_{X'}\ar@{=}[d]\ar[r]&f^*(\cF)\ar[r]\ar[d]^{\lambda}&{f^*(\Omega)}\ar@{=}[d]\ar[r]&0\\
0\ar[r]&{\co_{X'}}\ar[r]&{\cF^+}\ar[r]&{f^*(\Omega)}\ar[r]&0}
\end{equation}
où les lignes sont les suites exactes canoniques \eqref{higgs1-eph13b}. Par suite, $\lambda$ est un isomorphisme.
On en déduit un isomorphisme de  $(\bT\times_XX')$-fibrés principaux homogènes 
\begin{equation}\label{higgs1-eph15g}
\bL^+\stackrel{\sim}{\rightarrow} \bL\times_XX'.
\end{equation}

Le diagramme
\begin{equation}\label{higgs1-eph15h}
\xymatrix{
{f^*(\cL)}\ar[r]^{a}\ar[d]_{f^*(\iota)}&{f^+(\cL)}\ar[d]^{\iota^+}\\
{f^*(\varphi_X(\bL))}\ar[r]^{b}&{\varphi_{X'}(\bL^+)}}
\end{equation}
où $\iota$ et $\iota^+$ sont les isomorphismes canoniques \eqref{higgs1-eph13i},
$a$ est le morphisme canonique 
et $b$ est le morphisme induit par \eqref{higgs1-eph11b} et \eqref{higgs1-eph15g}, est commutatif. En effet, il suffit de 
montrer que le diagramme
\begin{equation}\label{higgs1-eph15i}
\xymatrix{
{\cL}\ar[r]^-(0.5){a_\sharp}\ar[d]_{\iota}&{f_*(\cL^+)}\ar[d]^{f_*(\iota^+)}\\
{\varphi_X(\bL)}\ar[r]^(0.4){b_\sharp}&{f_*(\varphi_{X'}(\bL^+))}}
\end{equation}
où $a_\sharp$ et $b_\sharp$ sont les morphismes adjoints de $a$ et $b$, est commutatif, 
ou encore que le diagramme déduit de ce dernier en évaluant les faisceaux en $X$ est commutatif. 
Soient $s\in \cL(X)$, $\rho_s\colon \cF\rightarrow \co_X$ le morphisme qui à toute section locale $\ell$ de $\cF$
associe $\ell(s)$.  Posons $s'=a_\sharp(s)\in \cL^+(X')$ et notons $\rho_{s'}\colon \cF^+\rightarrow \co_{X'}$ 
le morphisme qui à toute section locale $\ell'$ de $\cF^+$ associe $\ell'(s')$. Il résulte aussitôt de la définition 
du morphisme $\lambda_\sharp$ \eqref{higgs1-eph15c} que le diagramme 
\begin{equation}\label{higgs1-eph15aa}
\xymatrix{
{\cF}\ar[r]^-(0.5){\lambda_\sharp}\ar[d]_{\rho_s}&{f_*(\cF^+)}\ar[d]^{f_*(\rho_{s'})}\\
{\co_X}\ar[r]&{f_*(\co_{X'})}}
\end{equation}
est commutatif. On en déduit la relation $\rho_{s'}\circ \lambda=f^*(\rho_s)$,
qui implique aussitôt la commutativité du diagramme \eqref{higgs1-eph15i}.

Soient $g\colon X''\rightarrow X'$ un morphisme de schémas, $h=f\circ g \colon X''\rightarrow X$.
On désigne par $\cF^\dagger$ le faisceau des fonctions affines sur $h^+(\cL)$.
On a alors un isomorphisme canonique 
\begin{equation}\label{higgs1-eph15ab}
\theta \colon h^*(\cF)\stackrel{\sim}{\rightarrow}\cF^\dagger.
\end{equation}
Compte tenu de \eqref{higgs1-eph150f}, on a aussi un isomorphisme canonique 
\begin{equation}\label{higgs1-eph15ac}
\lambda^+ \colon g^*(\cF^+)\stackrel{\sim}{\rightarrow}\cF^\dagger.
\end{equation}
On vérifie aussitôt qu'on a
\begin{equation}\label{higgs1-eph15ad}
\theta= \lambda^+ \circ g^*(\lambda).
\end{equation}

\subsection{}\label{higgs1-eph14}
Soient $f\colon X'\rightarrow X$ un morphisme de schémas, $\rT$ un $\co_X$-module localement libre de type fini,
$\rT'$ un $\co_{X'}$-module localement libre de type fini, $\cL$ un $\rT$-torseur de $X_\zar$, 
$\cL'$ un $\rT'$-torseur de $X'_\zar$. On pose $\Omega=\cHom_{\co_X}(\rT,\co_X)$, 
$\Omega'=\cHom_{\co_{X'}}(\rT',\co_{X'})$, $\bT=\Spec(\rS_{\co_X}(\Omega))$ et $\bT'=\Spec(\rS_{\co_{X'}}(\Omega'))$. 
On désigne par $\cF$ le faisceau des fonctions affines sur $\cL$ \eqref{higgs1-eph13}, par 
$\cF'$ le faisceau des fonctions affines sur $\cL'$, par 
$\bL$ le $\bT$-fibré principal homogène canonique sur $X$ qui représente $\cL$ \eqref{higgs1-eph121} et par
$\bL'$ le $\bT'$-fibré principal homogène canonique sur $X'$ qui représente $\cL'$. 
Le faisceau $f_*(\cL')$ est naturellement un $f_*(\rT')$-objet de $X_\zar$.  
Soient $u\colon \rT\rightarrow f_*(\rT')$ un morphisme $\co_X$-linéaire, 
$v\colon \cL\rightarrow f_*(\cL')$ un morphisme $u$-équivariant. On désigne 
par $\gamma\colon f^{-1}(\rT)\rightarrow f^*(\rT)$ 
le morphisme canonique, par $u^\sharp\colon f^*(\rT)\rightarrow \rT'$ le morphisme adjoint de $u$,
par $u^\vee\colon \Omega'\rightarrow f^*(\Omega)$ le morphisme dual de $u^\sharp$ et 
par $v^\sharp\colon f^*(\cL)\rightarrow \cL'$ le morphisme adjoint de $v$. 
Comme $v^\sharp$ est $(u^\sharp\circ \gamma)$-équivariant, il se factorise de manière unique à travers un morphisme 
$u^\sharp$-équivariant 
\begin{equation}\label{higgs1-eph14c}
v^+\colon f^+(\cL)\rightarrow \cL'.
\end{equation}
D'après \eqref{higgs1-eph122a} et \eqref{higgs1-eph15e}, celui-ci induit un morphisme $\co_{X'}$-linéaire 
\begin{equation}\label{higgs1-eph14d}
w\colon \cF'\rightarrow f^*(\cF)
\end{equation}
qui s'insère dans un diagramme commutatif 
\begin{equation}\label{higgs1-eph14e}
\xymatrix{
0\ar[r]&{\co_{X'}}\ar[r]\ar@{=}[d]&{\cF'}\ar[r]\ar[d]^w&{\Omega'}\ar[r]\ar[d]^{u^\vee}&0\\
0\ar[r]&{\co_{X'}}\ar[r]&{f^*(\cF)}\ar[r]&{f^*(\Omega)}\ar[r]&0}
\end{equation}
où les lignes sont les suites exactes canoniques \eqref{higgs1-eph13b}.

Le morphisme $u^\vee$ induit un morphisme de $X'$-schémas en groupes 
\begin{equation}\label{higgs1-eph14f}
\alpha\colon \bT\times_XX'\rightarrow \bT'.
\end{equation}
Le morphisme $w$ induit un $X'$-morphisme $\alpha$-équivariant 
\begin{equation}\label{higgs1-eph14g}
\beta\colon \bL\times_XX'\rightarrow \bL'.
\end{equation}
Il résulte de \eqref{higgs1-eph122d} et \eqref{higgs1-eph15h} que le diagramme
\begin{equation}\label{higgs1-eph14h}
\xymatrix{
{f^*(\cL)}\ar[r]^-(0.45){v^\sharp}\ar[d]_{f^*(\iota)}&{\cL'}\ar[d]^{\iota'}\\
{f^*(\varphi_X(\bL))}\ar[r]^-(0.5){\delta}&{\varphi_{X'}(\bL')}}
\end{equation}
où $\iota$ et $\iota'$ sont les isomorphismes canoniques \eqref{higgs1-eph13i}
et $\delta$ est le morphisme induit par \eqref{higgs1-eph11b} et $\beta$, est commutatif.

\subsection{}\label{higgs1-eph141}
Conservons les hypothèses de \eqref{higgs1-eph14}, soient, de plus, $g\colon X''\rightarrow X'$ un morphisme de schémas, 
$\rT''$ un $\co_{X''}$-module localement libre de type fini,
$\cL''$ un $\rT''$-torseur de $X''_\zar$, $\cF''$ le faisceau des fonctions affines sur $\cL''$.
Soient $u'\colon \rT'\rightarrow g_*(\rT'')$ un morphisme $\co_{X'}$-linéaire, 
$v'\colon \cL'\rightarrow g_*(\cL'')$ un morphisme $u'$-équivariant. D'après \ref{higgs1-eph14},
le couple $(u',v')$ induit un morphisme $\co_{X''}$-linéaire 
\begin{equation}\label{higgs1-eph141a}
w'\colon \cF''\rightarrow g^*(\cF').
\end{equation}
De même, le couple $(f_*(u')\circ u,f_*(v')\circ v)$ induit un morphisme $\co_{X''}$-linéaire 
\begin{equation}\label{higgs1-eph141b}
t\colon \cF''\rightarrow g^*(f^*(\cF)).
\end{equation}
On a alors 
\begin{equation}\label{higgs1-eph141c}
t=g^*(w)\circ w'.
\end{equation}
Cela résulte de \eqref{higgs1-eph15ad} et d'une chasse au diagramme commutatif
\begin{equation}\label{higgs1-eph141d}
\xymatrix{
{g^*(f^*(\cL))}\ar[rr]^-(0.5){g^*(v^\sharp)}\ar[rd]&&{g^*(\cL')}\ar[rr]^-(0.5){v'^\sharp}\ar[rd]&&{\cL''}\\
&{g^*(f^+(\cL))}\ar[ru]_{g^*(v^+)}\ar[rd]&&{g^+(\cL')}\ar[ru]_{v'^+}&\\
&&{g^+(f^+(\cL))}\ar[ru]_{g^+(v^+)}&&}
\end{equation}
où $v'^\sharp$ est le morphisme l'adjoint de $v'$, $v'^+$ le morphisme induit par $v'^\sharp$
et les flèches non libellées sont les morphismes canoniques.

\subsection{}\label{higgs1-eph6}
Soient $\cC$, $\cF$ deux catégories, 
\begin{equation}\label{higgs1-eph6a}
\pi\colon \cF\rightarrow \cC
\end{equation}
un foncteur fibrant (\cite{sga1} VI 6.1), $\Delta$ un groupe (abstrait).
Pour tout $X\in \ob(\cC)$, on désigne par $\cF_X$ la catégorie fibre de $\pi$ au-dessus de $X$ (\cite{sga1} VI §4). 
On note encore $\Delta$ le groupoïde associé à $\Delta$, {\em i.e.},
la catégorie ayant un seul objet, de classe de morphismes $\Delta$.

Soient $X$ un objet de $\cC$, $\varphi$ une action à gauche de $\Delta$ sur $X$, autrement dit, 
$\varphi\colon \Delta\rightarrow \cC$ est un foncteur qui fait correspondre $X$ à l'unique objet de $\Delta$.
Pour $\sigma\in \Delta$, on note (abusivement) $\sigma$ l'automorphisme $\varphi(\sigma)$ de $X$. 
Le changement de base de $\pi$ par $\varphi$ (\cite{sga1} VI §3) 
\begin{equation}\label{higgs1-eph6b}
\pi_\Delta\colon \cF_\Delta\rightarrow \Delta
\end{equation}
est un foncteur fibrant. 
La catégorie $\cF_\Delta$ a mêmes objets que la catégorie fibre $\cF_X$, mais elle a plus de morphismes. 
Les sections cartésiennes de $\pi_\Delta$ sont appelées les 
{\em $X$-objets $\Delta$-équivariants} de $\cF$ (ou {\em objets $\Delta$-équivariants} de $\cF_X$). 
La donnée d'un $X$-objet $\Delta$-équivariant de $\cF$ est donc équivalente à la donnée 
d'un objet $Y$ de $\cF_X$ et d'une action à gauche de $\Delta$ sur $Y$ en tant qu'objet de $\cF$, 
compatible avec son action sur $X$ par le foncteur $\pi$. Les 
$X$-objets $\Delta$-équivariants de $\cF$ forment naturellement une catégorie, à savoir, 
la catégorie des sections cartésiennes de  $\pi_\Delta$ (\cite{sga4} VI 6.10). 

Choisissons un clivage normalisé de $\cF$ sur $\cC$ (\cite{sga1} VI §7; cf. aussi \cite{egr1} 1.1.2),
autrement dit, choisissons pour tout morphisme $f\colon Z\rightarrow Y$ de $\cC$ un foncteur image inverse 
\begin{equation}\label{higgs1-eph6c}
f^*\colon \cF_Y\rightarrow \cF_Z,
\end{equation}
tel que pour tout $Y\in \ob(\cC)$, $f=\id_Y$ implique que $f^*=\id_{\cF_Y}$. 
Pour tout couple de morphismes composables $(f,g)$ de $\cC$, on a un isomorphisme canonique de foncteurs
\begin{equation}\label{higgs1-eph6cd}
c_{g,f}\colon g^*f^*\stackrel{\sim}{\rightarrow}(fg)^*,
\end{equation}
vérifiant des relations de compatibilité (\cite{sga1} VI 7.4). 
D'après (\cite{sga1} VI §12), la donnée d'un objet $\Delta$-équivariant de $\cF_X$ est équivalente à 
la donnée d'un objet $Y$ de $\cF_X$ et pour tout $\sigma\in \Delta$, d'un isomorphisme 
\begin{equation}\label{higgs1-eph6d}
\tau_\sigma^Y\colon Y\stackrel{\sim}{\rightarrow}\sigma^*(Y)
\end{equation}
tels que $\tau_\id^Y=\id_Y$ et que pour tout $(\sigma,\sigma')\in \Delta^2$, on ait 
\begin{equation}\label{higgs1-eph6e}
\tau^Y_{\sigma\sigma'}=c_{\sigma',\sigma}\circ (\sigma'^**\tau^Y_\sigma)\circ \tau^Y_{\sigma'}.
\end{equation}
On laissera le soin au lecteur de décrire explicitement les morphismes de $X$-objets $\Delta$-équivariants de $\cF$.

\subsection{}\label{higgs1-eph7}
Soient $\cC$ une catégorie dans laquelle les produits fibrés sont représentables, 
$X$ un objet de $\cC$ muni d'une action à gauche d'un groupe (abstrait) $\Delta$. 
On note encore $\Delta$ le groupoïde associé à $\Delta$,
et $\varphi\colon \Delta\rightarrow \cC$ l'action de $\Delta$ sur $X$. 
On désigne par $\bFl(\cC)$  la catégorie des flèches de $\cC$ et par
\begin{equation}\label{higgs1-eph7a}
\bFl(\cC)\rightarrow \cC
\end{equation}
le foncteur-but, qui est un foncteur fibrant (\cite{sga1} VI §11 a)). On en déduit par changement de base
par $\varphi$ (\cite{sga1} VI §3) un foncteur fibrant 
\begin{equation}\label{higgs1-eph7b}
\bFl(\cC)_\Delta\rightarrow \Delta,
\end{equation}
dont les sections cartésiennes sont appelées les {\em $X$-objets $\Delta$-équivariants} de $\cC$
(ou {\em objets $\Delta$-équiva\-riants} de $\cC_{/X}$). La donnée d'un tel objet est équivalente à la donnée 
d'un objet $Y$ de $\cC_{/X}$ et d'une action à gauche de $\Delta$ sur $Y$ en tant qu'objet de $\cC$, 
compatible avec son action sur $X$.

On désigne par $\bGr(\cC)$ la catégorie suivante. Les objets de $\bGr(\cC)$ sont les couples 
formés d'un morphisme $G\rightarrow Y$ de $\cC$ et d'une structure sur $G$ de groupe de $\cC_{/Y}$ dans 
le sens de (\cite{giraud2} III 1.1.1). On omettra dans la suite la structure de groupe des notations. 
Soient $G\rightarrow Y$, $G'\rightarrow Y'$ deux objets de  $\bGr(\cC)$. Un morphisme de $G'\rightarrow Y'$ 
dans $G\rightarrow Y$ est la donnée d'un diagramme commutatif de~$\cC$ 
\begin{equation}\label{higgs1-eph7c}
\xymatrix{
G'\ar[r]\ar[d]&G\ar[d]\\
Y'\ar[r]&Y}
\end{equation}
tel que le morphisme induit $G'\rightarrow G\times_YY'$ soit un morphisme de groupes de $\cC_{/Y'}$. 
Le foncteur but
\begin{equation}\label{higgs1-eph7d}
\bGr(\cC)\rightarrow \cC, \ \ \ (G\rightarrow Y)\mapsto Y
\end{equation}
est fibrant~; sa catégorie fibre au-dessus d'un objet $Y$ de $\cC$ est la catégorie
des groupes de $\cC_{/Y}$. On en déduit par  changement de base par $\varphi$ un foncteur fibrant
\begin{equation}\label{higgs1-eph7e}
\bGr(\cC)_\Delta\rightarrow \Delta,
\end{equation}
dont les sections cartésiennes sont appelées les {\em $X$-groupes $\Delta$-équivariants} de $\cC$
(ou {\em groupes $\Delta$-équivariants} de $\cC_{/X}$). La donnée d'un tel objet est équivalente à la donnée 
d'un groupe $G$ de $\cC_{/X}$ et d'une action à gauche de $\Delta$ sur $G$ en tant qu'objet de $\cC$,
compatible avec son action sur $X$ et avec la structure de groupe de $G$ dans un sens que nous n'explicitons pas.

On désigne par $\bOp(\cC)$ la catégorie suivante (cf. \cite{giraud2} III 1.1.6). 
Les objets de $\bOp(\cC)$ sont les triplets 
formés d'un objet $G\rightarrow Y$ de $\bGr(\cC)$, d'un morphisme $Z\rightarrow Y$ de $\cC$
et d'une opération (à droite) $m$ de $G$ sur $Z$ au-dessus de $Y$, {\em i.e.}, d'un $Y$-morphisme
$m\colon Z\times_YG\rightarrow Z$ assujetti aux conditions algébriques habituelles (cf. \cite{giraud2} III 1.1.2). 
Soient $(G\rightarrow Y,Z\rightarrow Y)$, $(G'\rightarrow Y',Z'\rightarrow Y')$ deux objets de $\bOp(\cC)$. 
Un morphisme de $(G\rightarrow Y,Z\rightarrow Y,m)$ dans $(G'\rightarrow Y',Z'\rightarrow Y',m')$ 
est la donnée de deux diagrammes commutatifs de $\cC$ 
\begin{equation}
\xymatrix{
G'\ar[r]\ar[d]&G\ar[d]\\
Y'\ar[r]&Y} \ \ \ \xymatrix{
Z'\ar[r]\ar[d]&Z\ar[d]\\
Y'\ar[r]&Y}
\end{equation}
tels que le diagramme 
\begin{equation}
\xymatrix{
{Z'\times_{Y'}G'}\ar[r]\ar[d]_{m'}&{Z\times_YG\times_YY'}\ar[d]^{m\times_YY'}\\
Z'\ar[r]&{Z\times_YY'}}
\end{equation}
soit commutatif. Le foncteur
\begin{equation}
\bOp(\cC)\rightarrow \cC, \ \ \ (G\rightarrow Y,Z\rightarrow Y,m)\mapsto Y
\end{equation}
est fibrant. Sa catégorie fibre au-dessus d'un objet $Y$ de $\cC$ est la catégorie
des objets à groupe d'opérateurs (à droite) de $\cC_{/Y}$. 
On en déduit par  changement de base par $\varphi$ un foncteur fibrant
\begin{equation}
\bOp(\cC)_\Delta\rightarrow \Delta,
\end{equation}
dont les sections cartésiennes sont appelées les {\em $X$-objets à groupe d'opérateurs   
$\Delta$-équivariants} de $\cC$ (ou {\em objets à groupe d'opérateurs $\Delta$-équivariants} de $\cC_{/X}$). 
La donnée d'un tel objet est équivalente à la donnée d'un groupe $\Delta$-équivariant $G$ de $\cC_{/X}$, 
d'un objet $\Delta$-équivariant $Y$ de $\cC_{/X}$ et d'une opération $m$ de $G$ sur $Y$, compatible 
avec les structures $\Delta$-équivariantes dans un sens que nous n'explicitons pas. 
On dit aussi que $(Y,m)$ est un {\em $G$-objet $\Delta$-équivariant de $\cC_{/X}$}.

\subsection{}\label{higgs1-eph5}
Soit $X$ un schéma muni d'une action à gauche d'un groupe (abstrait) $\Delta$. 
On définit les {\em objets $\Delta$-équivariants de $X_\zar$} en prenant dans \eqref{higgs1-eph6} 
pour $\cF$ la catégorie fibrée, clivée et normalisée
\begin{equation}
\cZ\rightarrow \Sch
\end{equation}
obtenue en associant à tout schéma $Y$ le topos $Y_\zar$ et à tout morphisme de schémas 
$f\colon Z\rightarrow Y$ le foncteur image inverse $f^*\colon  Y_\zar\rightarrow Z_\zar$. 
On définit les {\em groupes $\Delta$-équivariants de $X_\zar$} en prenant dans \eqref{higgs1-eph6} 
pour $\cF$ la catégorie fibrée, clivée et normalisée
\begin{equation}
\cG\rightarrow \Sch
\end{equation}
obtenue en associant à tout schéma $Y$ la catégorie des groupes de $Y_\zar$ et à tout morphisme de schémas 
$f\colon Z\rightarrow Y$ le foncteur image inverse $f^*$. 
On définit les {\em $\co_X$-modules $\Delta$-équivariants de $X_\zar$} en prenant dans \eqref{higgs1-eph6} 
pour $\cF$ la catégorie fibrée, clivée et normalisée
\begin{equation}
\cM\rightarrow \Sch
\end{equation}
obtenue en associant à tout schéma $Y$ la catégorie des $\co_Y$-modules de $Y_\zar$ 
et à tout morphisme de schémas $f\colon Z\rightarrow Y$ le foncteur image inverse au sens des modules. 
On notera que les images inverses d'un $\co_X$-module par un automorphisme de $X$ 
aux sens des modules et des groupes abéliens
étant égaux, tout $\co_X$-module $\Delta$-équivariant de $X_\zar$ est 
aussi un groupe $\Delta$-équivariant. 

On désigne par $\cP$ la catégorie suivante. Les objets de $\cP$ sont les triplets $(Y,G,P)$
où $Y$ est un schéma, $G$ est un groupe de $Y_\zar$, $P$ est un $G$-objet (à droite) de $Y_\zar$ (\cite{giraud2} III 1.1.2). 
Soient $(Y,G,P)$, $(Y',G',P')$ deux objets de $\cP$. Un morphisme de $(Y',G',P')$ dans $(Y,G,P)$
est la donnée d'un morphisme de schémas $f\colon Y'\rightarrow Y$, d'un morphisme de groupes
$\gamma\colon G'\rightarrow f^*(G)$ de $Y'_\zar$ et d'un morphisme de $\gamma$-équivariant 
$\delta\colon P'\rightarrow f^*(P)$ de $Y'_\zar$. Le foncteur 
\begin{equation}
\cP\rightarrow \Sch, \ \ \ (Y,G,P)\mapsto Y
\end{equation}
est fibrant. Prenant dans \eqref{higgs1-eph6} pour $\cF$ la catégorie fibrée ci-dessus, on obtient la notion
d'{\em objets à groupe d'opérateurs $\Delta$-équivariants de $X_\zar$}. La donnée d'un tel objet
est équivalente à la donnée d'un groupe $\Delta$-équivariant $G$ de $X_\zar$, d'un objet 
$\Delta$-équivariant $P$ de $X_\zar$ et d'une opération $m$ de $G$ sur $P$, 
compatible avec les structures $\Delta$-équivariantes dans un sens que nous n'explicitons pas. 
On dit aussi que $P$ est un {\em $G$-objet $\Delta$-équivariant de $X_\zar$}. 
Si, de plus, $P$ est un $G$-torseur de $X_\zar$, on dit encore que c'est un 
{\em $G$-torseur $\Delta$-équivariant de $X_\zar$}. 

\begin{remas}\label{higgs1-eph55}
Soient $X$ et $X'$ deux schémas munis d'actions à gauche d'un groupe (abstrait) $\Delta$,
$f\colon X'\rightarrow X$ un morphisme $\Delta$-équivariant. 
\begin{itemize}
\item[(i)] Pour tout objet à groupe d'opérateurs $\Delta$-équivariant
$(G,P)$ de $X_\zar$, $(f^*(G),f^*(P))$ est naturellement 
un objet à groupe d'opérateurs $\Delta$-équivariant de $X'_\zar$.
\item[(ii)] Pour tout objet à groupe d'opérateurs $\Delta$-équivariant $(G',P')$ de $X'_\zar$, 
$(f_*(G'),f_*(P'))$ est naturellement un objet à groupe d'opérateurs $\Delta$-équivariant de $X_\zar$.
En effet, pour tout $\sigma\in \Delta$, 
le morphisme de changement de base pour les topos de Zariski $\sigma^*f_*\rightarrow f_*\sigma^*$,
déduit de la relation $f\sigma=\sigma f$ (\cite{egr1} 1.2.2),
est un isomorphisme. L'assertion s'ensuit compte tenu de (\cite{egr1} 1.2.4(i)).
\end{itemize}
\end{remas}

\subsection{}\label{higgs1-eph4}\index{Fibre principal@$G$-fibré principal homogène!1@$\Delta$-équivariant}
Soit $X$ un schéma muni d'une action à gauche d'un groupe (abstrait) $\Delta$, 
$G$ un $X$-schéma en groupes $\Delta$-équivariant. 
Un {\em $G$-fibré principal homogène $\Delta$-équivariant sur $X$} est un $G$-objet $\Delta$-équivariant
$(Y,m)$ de $\Sch_{/X}$ \eqref{higgs1-eph7} tel que $Y$ soit un $G$-fibré principal homogène sur $X$ \eqref{higgs1-eph2}. 

Compte tenu de \ref{higgs1-eph11}, le foncteur \eqref{higgs1-eph1a}
\begin{equation}\label{higgs1-eph4a}
\varphi_X\colon \Sch_{/X}\rightarrow X_\zar, \ \ \ Y\mapsto \Hom_X(-,Y)
\end{equation}
transforme les $X$-schémas (resp. $X$-schémas en groupes) $\Delta$-équivariants \eqref{higgs1-eph7}
en objets (resp. groupes) $\Delta$-équivariants de $X_\zar$ \eqref{higgs1-eph5}.
De même, posons $\uG=\varphi_X(G)$, le foncteur \eqref{higgs1-eph2a} 
\begin{equation}\label{higgs1-eph4b}
\varphi_X\colon \bFPH(G/X)\rightarrow \bTors(\uG,X_\zar), \ \ \ Y\mapsto \varphi_X(Y)
\end{equation}
transforme les $G$-fibrés principaux homogènes $\Delta$-équivariants sur $X$ en
$\uG$-torseurs $\Delta$-équivariants de $X_\zar$.
Inversement, soient $Y$ un $G$-fibré principal homogène sur $X$,  
$\uY=\varphi_X(Y)$. La donnée sur le $\uG$-torseur $\uY$ 
d'une structure $\Delta$-équivariante
détermine sur $Y$ une et une unique structure de $G$-fibré principal homogène 
$\Delta$-équivariant au-dessus de $X$. En effet, pour tout $\sigma\in \Delta$, notons 
\begin{eqnarray}
\tau_\sigma^{G}\colon G&\stackrel{\sim}{\rightarrow} &\sigma^\bullet(G),\label{higgs1-eph4c}\\
\tau_\sigma^{\uY}\colon \uY&\stackrel{\sim}{\rightarrow}& \sigma^*(\uY)\label{higgs1-eph4d}
\end{eqnarray}
les isomorphismes induits par les structures $\Delta$-équivariantes de $G$ et $\uY$ \eqref{higgs1-eph6d}.   
Par définition, $\tau_\sigma^{G}$ est un isomorphisme de $X$-schémas en groupes. 
D'après \ref{higgs1-eph11}, il induit un isomorphisme de groupes de $X_\zar$
\begin{equation}\label{higgs1-eph4e}
\tau_\sigma^{\uG}\colon \uG\stackrel{\sim}{\rightarrow} \sigma^*(\uG).
\end{equation}
Les isomorphismes $(\tau_\sigma^{\uG})_{\sigma\in \Delta}$ font de $\uG$ un groupe $\Delta$-équivariant de $X_\zar$. 
Pour tout $\sigma\in \Delta$, $\tau_\sigma^\uY$ est un isomorphisme de $\uG$-torseurs, où $\sigma^*(\uY)$
est considéré comme $\uG$-torseur via $\tau_\sigma^\uG$. 
D'après \ref{higgs1-eph3}(i), $\tau_\sigma^\uY$ est l'image par le foncteur \eqref{higgs1-eph2a} d'un isomorphisme de $\bFPH(G/X)$
\begin{equation}\label{higgs1-eph4g}
\tau^Y_\sigma\colon Y\stackrel{\sim}{\rightarrow} \sigma^\bullet(Y),
\end{equation}
où $\sigma^\bullet(Y)$ est considéré comme un $G$-fibré principal homogène via $\tau_\sigma^G$. 
Les isomorphismes $(\tau_\sigma^Y)_{\sigma\in \Delta}$ satisfont les relations de compatibilité \eqref{higgs1-eph6e}.
Ils font de $Y$ un $G$-fibré principal homogène $\Delta$-équivariant sur $X$ (cf. \ref{higgs1-eph6}).

\subsection{}\label{higgs1-eph17}
Soient $X$ un schéma muni d'une action à gauche d'un groupe (abstrait) $\Delta$, 
$\rT$ un $\co_X$-module localement libre de type fini et $\Delta$-équivariant \eqref{higgs1-eph5}, 
$\cL$ un $\rT$-torseur $\Delta$-équivariant de $X_\zar$. On pose $\Omega=\cHom_{\co_X}(\rT,\co_X)$
et $\bT=\Spec(\rS_{\co_X}(\Omega))$, et on désigne par $\cF$ le faisceau des fonctions affines sur
$\cL$ \eqref{higgs1-eph13} et par $\bL$ le $\bT$-fibré principal homogène canonique sur $X$ 
qui représente $\cL$ \eqref{higgs1-eph121}. Pour tout $\sigma\in \Delta$, on note 
\begin{equation}
\tau_\sigma^\rT\colon \rT\stackrel{\sim}{\rightarrow}\sigma^*(\rT) \ \ \ {\rm et}\ \ \ 
\tau_\sigma^\cL\colon \cL\stackrel{\sim}{\rightarrow}\sigma^*(\cL)
\end{equation}
les isomorphismes qui définissent les structures $\Delta$-équivariantes sur $\rT$ et sur $\cL$ (cf. \ref{higgs1-eph6}). 
D'après \ref{higgs1-eph14}, les inverses de $\tau_\sigma^\rT$ et $\tau_\sigma^\cL$ induisent
un morphisme $\co_X$-linéaire 
\begin{equation}
\tau_\sigma^\cF\colon \cF\rightarrow \sigma^*(\cF)
\end{equation}
qui s'insère dans un diagramme commutatif 
\begin{equation}
\xymatrix{
0\ar[r]&{\co_{X}}\ar[r]\ar@{=}[d]&{\cF}\ar[r]\ar[d]^{\tau_\sigma^\cF}&{\Omega}\ar[r]\ar[d]^{\tau_\sigma^\Omega}&0\\
0\ar[r]&{\co_{X}}\ar[r]&{\sigma^*(\cF)}\ar[r]&{\sigma^*(\Omega)}\ar[r]&0}
\end{equation}
où les lignes sont les suites exactes canoniques \eqref{higgs1-eph13b} et $\tau_\sigma^\Omega$ est l'inverse de l'isomorphisme
induit par $\tau_\sigma^\rT$. Il s'ensuit que $\tau_\sigma^\cF$ est un isomorphisme. 
D'après \eqref{higgs1-eph141c}, les $(\tau_\sigma^\cF)_{\sigma\in \Delta}$ 
satisfont les relations de compatibilité \eqref{higgs1-eph6e}. Ils font donc de $\cF$ un $\co_X$-module $\Delta$-équivariant. 
De même, les $(\tau_\sigma^\Omega)_{\sigma\in \Delta}$ font de $\Omega$ un 
$\co_X$-module $\Delta$-équivariant. On en déduit sur $\bT$ une structure de $X$-schéma en groupes 
$\Delta$-équivariant et sur $\bL$ une structure de $\bT$-fibré principal homogène $\Delta$-équivariant sur $X$ \eqref{higgs1-eph4}. 
D'après \ref{higgs1-eph12} et \eqref{higgs1-eph14h}, on a un isomorphisme d'objets à groupe d'opérateurs 
$\Delta$-équivariants de $X_\zar$
\begin{equation}
(\rT,\cL)\stackrel{\sim}{\rightarrow}(\varphi_X(\bT),\varphi_X(\bL)).
\end{equation} 

\subsection{}\label{higgs1-eph19}
Soient $X$ et $X'$ deux schémas munis d'actions à gauche d'un groupe (abstrait) $\Delta$,
$f\colon X'\rightarrow X$ un morphisme $\Delta$-équivariant,
$\rT$ un $\co_X$-module localement libre de type fini et $\Delta$-équivariant,
$\rT'$ un $\co_{X'}$-module localement libre de type fini et $\Delta$-équivariant, $\cL$ un $\rT$-torseur 
$\Delta$-équivariant de $X_\zar$, $\cL'$ un $\rT'$-torseur $\Delta$-équivariant de $X'_\zar$. 
On désigne par $\cF$ le faisceau des fonctions affines sur $\cL$ \eqref{higgs1-eph13} et par 
$\cF'$ le faisceau des fonctions affines sur $\cL'$. 
Le couple $(f_*(\rT'),f_*(\cL'))$ est naturellement un objet à groupe d'opérateurs 
$\Delta$-équivariant de $X_\zar$ \eqref{higgs1-eph55}.  
Soient $u\colon \rT\rightarrow f_*(\rT')$ un morphisme $\co_X$-linéaire et $\Delta$-équivariant, 
$v\colon \cL\rightarrow f_*(\cL')$ un morphisme $u$-équivariant et $\Delta$-équivariant.
Le couple $(u,v)$ induit un morphisme $\co_{X'}$-linéaire \eqref{higgs1-eph14d}
\begin{equation}
w\colon \cF'\rightarrow f^*(\cF).
\end{equation}
Il résulte aussitôt de \eqref{higgs1-eph141c} que $w$ est $\Delta$-équivariant
lorsque l'on munit $\cF$ et $\cF'$ des structures $\Delta$-équivariantes canoniques \eqref{higgs1-eph17}. 

Notons $u^\sharp\colon f^*(\rT)\rightarrow \rT'$ le morphisme adjoint de $u$,
$v^\sharp\colon f^*(\cL)\rightarrow \cL'$ le morphisme adjoint de $v$ et  $v^+\colon f^+(\cL)\rightarrow \cL'$ 
le morphisme induit par $v^\sharp$ \eqref{higgs1-eph14c}. 
Si $u^\sharp$ est un isomorphisme, $v^+$ est un isomorphisme de $f^*(\rT)$-torseurs, 
et $w$ est un isomorphisme en vertu de \ref{higgs1-eph15}.

\section{Lexique de géométrie logarithmique}\label{higgs1-LOG}

Nous rappelons quelques notions de géométrie logarithmique qui joueront un rôle important dans 
la suite de cet article, dans le but de fixer les notations et de donner des repères aux lecteurs non familiers 
avec cette théorie.  
Nous renvoyons à \cite{kato1,kato2,gr1, ogus} pour les développements systématiques de la théorie. 

\subsection{}\label{higgs1-log1}
On sous-entend par {\em monoïde}  un monoïde commutatif et unitaire. 
Les homomorphismes de monoïdes sont toujours supposés transformer l'élément unité en l'élément unité.
Si $M$ est un monoïde, on désigne par $M^\gp$ le groupe associé, par $M^\times$ le groupe des unités de $M$,
par $M^\sharp$ l'ensemble des orbites $M/M^\times$ (qui est aussi le quotient de $M$ par $M^\times$ dans  la 
catégorie des monoïdes) et par $i_M\colon M\rightarrow M^\gp$ l'homomorphisme canonique. On pose
$M^\integre=i_M(M)$ et 
\begin{equation}\label{higgs1-log1b}
M^\sat=\{x\in M^\gp | x^n\in i_M(M)\ {\rm pour\ un \ entier}\ n\geq 1\}.
\end{equation}

On dit qu'un monoïde $M$ est {\em intègre} si l'homomorphisme canonique $i_M\colon M\rightarrow M^\gp$
est injectif, que $M$ est {\em fin} s'il est intègre et de type fini, que $M$ est {\em saturé} s'il est intègre et est isomorphe à 
$M^\sat$ et que $M$ est {\em torique} s'il est fin et saturé et si $M^\gp$ est libre sur $\mZ$.

Si $M$ est intègre, $M^\sharp$ est intègre, et pour que $M$ soit saturé, il faut et il suffit que $M^\sharp$ soit saturé. 

On dit qu'un morphisme de monoïdes $u\colon M\rightarrow N$
est {\em strict} si le morphisme induit $u^\sharp\colon M^\sharp\rightarrow N^\sharp$ est un isomorphisme. 

\subsection{}\label{higgs1-log111}
Soit $u\colon M\rightarrow N$ un morphisme de monoïdes intègres. 
On dit que $u$ est {\em exact} si le diagramme 
\begin{equation}\label{higgs1-log1a}
\xymatrix{
M\ar[r]^u\ar[d]&N\ar[d]\\
{M^\gp}\ar[r]^{u^\gp}&{N^\gp}}
\end{equation}
est cartésien. On dit que $u$ est {\em intègre} si pour tout monoïde 
intègre $M'$ et tout homomorphisme $v\colon M\rightarrow M'$, la somme amalgamée 
$M'\oplus_MN$ est intègre. On dit que $u$ est {\em saturé} s'il est intègre et si pour tout monoïde 
saturé $M'$ et tout homomorphisme $v\colon M\rightarrow M'$, la somme amalgamée 
$M'\oplus_MN$ est saturée.

\subsection{}\label{higgs1-log2}
Soit $T$ un topos. On désigne par $\Mon_T$ la catégorie des monoïdes (commutatifs et unitaires) de $T$
et par $\Ab_T$ la catégorie des groupes abéliens de $T$. 
Le foncteur d'injection canonique de $\Ab_T$ dans $\Mon_T$ admet un adjoint à droite 
\begin{equation}\label{higgs1-log2a}
\Mon_T\rightarrow \Ab_T, \ \ \ \cM\mapsto \cM^\times.
\end{equation}
Il est immédiat de voir que pour tout $U\in \ob(T)$, le morphisme d'adjonction $\cM^\times(U)\rightarrow \cM(U)$
induit un isomorphisme $\cM^\times(U)\simeq \cM(U)^\times$. On dit qu'un monoïde $\cM$ 
de $T$ est {\em affûté} ({\em sharp} en anglais) si $\cM^\times=1_T$. 
Le foncteur d'injection canonique de la sous-catégorie pleine des monoïdes affûtés de $T$ dans $\Mon_T$ 
admet un adjoint à gauche 
\begin{equation}\label{higgs1-log2b}
\cM\mapsto \cM^\sharp=\cM/\cM^\times.
\end{equation}

Le foncteur d'injection canonique de $\Ab_T$ dans $\Mon_T$ admet un adjoint à gauche 
\begin{equation}\label{higgs1-log2c}
\Mon_T\rightarrow \Ab_T, \ \ \ \cM\mapsto \cM^\gp.
\end{equation}
On dit qu'un monoïde $\cM$ de $T$ est {\em intègre} si le morphisme d'adjonction $\cM\rightarrow \cM^\gp$
est un monomorphisme. On désigne par $\Mon_{T,\integre}$ la sous-catégorie pleine de $\Mon_T$ formée
des monoïdes intègres de $T$. Le foncteur d'injection canonique de $\Mon_{T,\integre}$ dans $\Mon_T$ admet un adjoint
à gauche 
\begin{equation}\label{higgs1-log2d}
\Mon_T\rightarrow \Mon_{T,\integre},\ \ \ \cM\rightarrow \cM^\integre.
\end{equation}

\subsection{}\label{higgs1-log3}
Soient $\cC$ un site, $\tcC$ le topos des faisceaux d'ensembles sur $\cC$ (relativement à un univers fixé).  
Pour tout préfaisceau de monoïdes $\cP$ sur $\cC$, on désigne 
par $\cP^\gp$ (resp. $\cP^\integre$) le préfaisceau de monoïdes sur $\cC$ qui à 
$U\in \ob(\cC)$ associe le monoïde $\cP(U)^\gp$ (resp.  $\cP(U)^\integre$)
et par $\cP^a$ le faisceau de monoïdes associé à $\cP$. On a alors un isomorphisme canonique fonctoriel
\begin{equation}
(\cP^\gp)^a\stackrel{\sim}{\rightarrow}(\cP^a)^\gp\label{higgs1-log3a}.
\end{equation}
Comme le foncteur $\cP\mapsto \cP^a$ est exact, il transforme les préfaisceaux de monoïdes intègres sur $\cC$
en des monoïdes intègres de $\tcC$. On a un isomorphisme canonique fonctoriel
\begin{equation}\label{higgs1-log3b}
(\cP^\integre)^a\stackrel{\sim}{\rightarrow}(\cP^a)^\integre.
\end{equation}
Par suite, pour qu'un monoïde 
$\cM$ de $\tcC$ soit intègre, il faut et il suffit que pour tout $U\in \ob(\cC)$, le monoïde $\cM(U)$ soit intègre. 

\subsection{}\label{higgs1-log4}
Soit $T$ un topos. On dit qu'un morphisme de monoïdes intègres $u\colon \cM\rightarrow \cN$ de $T$ est {\em exact} si 
le diagramme canonique 
\begin{equation}\label{higgs1-log4a}
\xymatrix{
\cM\ar[r]^u\ar[d]&\cN\ar[d]\\
{\cM^\gp}\ar[r]^{u^\gp}&{\cN^\gp}}
\end{equation}
est cartésien. Il revient au même de demander que pour tout $U\in \ob(T)$, l'homomorphisme
$u(U)\colon \cM(U)\rightarrow \cN(U)$ soit exact. En effet, si le diagramme 
\eqref{higgs1-log4a} est cartésien, il en est de même des diagrammes obtenus en l'évaluant en tout $U\in \ob(T)$.
Par suite, le diagramme 
\begin{equation}\label{higgs1-log4aa}
\xymatrix{
{\cM(U)}\ar[r]^u\ar[d]&{\cN(U)}\ar[d]\\
{\cM(U)^\gp}\ar[r]^{u^\gp}&{\cN(U)^\gp}}
\end{equation}
est cartésien car le morphisme canonique $\cM(U)^\gp\rightarrow \cM^\gp(U)$ est injectif (et il en est de même 
pour $\cN$). L'implication réciproque résulte de \eqref{higgs1-log3a} et du fait que le foncteur 
$\cP\mapsto \cP^a$ est exact. 

Soient $\cM$ un monoïde intègre de $T$, $n$ un entier $\geq 1$. On dit que $\cM$ est {\em $n$-saturé} si 
l'endomorphisme d'élévation à la puissance $n$ dans $\cM$, défini par $x\mapsto x^n$, 
est exact. On dit que $\cM$ est {\em saturé} s'il est $n$-saturé pour tout entier $n\geq 1$.  
Il revient au même de demander que pour tout $U\in \ob(T)$, le monoïde $\cM(U)$ est saturé. 
On désigne par $\Mon_{T,\sat}$ la sous-catégorie pleine de $\Mon_T$
formée des monoïdes saturés de $T$. Le foncteur d'injection canonique de $\Mon_{T,\sat}$ dans $\Mon_T$ admet 
un adjoint à gauche
\begin{equation}
\Mon_T\rightarrow \Mon_{T,\sat},\ \ \ \cM\mapsto \cM^\sat.
\end{equation}

\subsection{}\label{higgs1-log45}
Soient $\cC$ un site, $\tcC$ le topos des faisceaux d'ensembles sur $\cC$ (relativement à un univers fixé).  
Pour tout préfaisceau de monoïdes $\cP$ sur $\cC$, on désigne par $\cP^\sat$ le préfaisceau de monoïdes 
saturés défini, pour $U\in \ob(\cC)$, par $U\mapsto \cP(U)^\sat$.
Alors le foncteur ``faisceau de monoïdes associé'',
$\cP\mapsto \cP^a$, transforme les préfaisceaux de monoïdes saturés sur $\cC$ 
en des monoïdes saturés de $\tcC$, et on a un isomorphisme canonique fonctoriel
\begin{equation}\label{higgs1-log45a}
(\cP^\sat)^a\stackrel{\sim}{\rightarrow}(\cP^a)^\sat.
\end{equation}
Par suite, pour qu'un monoïde $\cM$ de $\tcC$ soit saturé, il faut et il suffit que  pour tout $U\in \ob(\cC)$,
le monoïde $\cM(U)$ soit saturé. 

\subsection{}
Soient $f\colon T'\rightarrow T$ un morphisme de topos, $\cM$ un monoïde de $T$. 

(i)\ Si $\cM$ est intègre (resp. saturé), il en est de même de $f^*(\cM)$. 

(ii)\ On a des isomorphismes canoniques fonctoriels 
\begin{eqnarray}
f^*(\cM)^\gp&\stackrel{\sim}{\rightarrow}&f^*(\cM^\gp),\\
f^*(\cM)^\integre&\stackrel{\sim}{\rightarrow}&f^*(\cM^\integre).
\end{eqnarray}
Si de plus $\cM$ est intègre, on a un isomorphisme canonique 
\begin{equation}
f^*(\cM)^\sat\stackrel{\sim}{\rightarrow}f^*(\cM^\sat).
\end{equation}

(iii)\ Si $u\colon \cM\rightarrow \cN$ est un morphisme exact de monoïdes intègres de $T$, 
alors $f^*(u)\colon f^*(\cM)\rightarrow f^*(\cN)$ est exact. 

\subsection{}
Soit $T$ un topos. 

(i)\ Si $M$ est un monoïde intègre (resp. saturé), 
le faisceau de monoïdes constant $M_T$ de valeur $M$ sur $T$ est intègre (resp. saturé). 

(ii) Supposons que $T$ ait suffisamment de points. 
Pour qu'un monoïde $\cM$ de $T$ soit intègre (resp. saturé), il faut et il suffit que pour tout point $p$ de $T$, 
le monoïde $\cM_p$ soit intègre (resp. saturé). Pour qu'un morphisme de monoïdes intègres 
$u\colon \cM\rightarrow \cN$ de $T$ soit exact, il faut et il suffit que 
pour tout point $p$ de $T$, l'homomorphisme $u_p\colon \cM_p\rightarrow \cN_p$ soit exact. 

(iii) Soit $\cM$ un monoïde intègre de $T$. Alors $\cM^\sharp$ est intègre, 
et pour que $\cM$ soit saturé, il faut et il suffit que $\cM^\sharp$ soit saturé. 

\subsection{}\label{higgs1-log5}
Une structure {\em pré-logarithmique} sur un schéma $X$ est une paire $(\cP,\beta)$, où $\cP$ est un faisceau
de monoïdes abéliens sur le site étale de $X$ et $\beta$ est un homomorphisme de $\cP$ dans le monoïde 
multiplicatif $\co_X$. Une structure pré-logarithmique $(\cP,\beta)$ est dite {\em logarithmique} si $\beta$ induit un
isomorphisme $\beta^{-1}(\co_X^\times)\stackrel{\sim}{\rightarrow}\co_X^\times$. Les structures pré-logarithmiques
sur $X$ forment naturellement une catégorie qui contient la sous-catégorie pleine des structures logarithmiques sur $X$. 
L'injection canonique de la catégorie des structures logarithmiques sur $X$ dans celle des structures pré-logarithmiques  sur
$X$ admet un adjoint à gauche. Il associe à une structure pré-logarithmique $(\cP,\beta)$ la structure logarithmique 
$(\cM,\alpha)$, où $\cM$ est défini par le diagramme co-cartésien 
\begin{equation}\label{higgs1-log5a}
\xymatrix{
{\beta^{-1}(\co_X^\times)}\ar[d]\ar[r]&{\cP}\ar[d]\\
{\co_X^\times}\ar[r]&{\cM}}
\end{equation}
On dit que $(\cM,\alpha)$ est la structure logarithmique {\em associée} à $(\cP,\beta)$. 

\subsection{}\label{higgs1-log50}
Soient $f\colon X\rightarrow Y$ un morphisme de schémas.
Nous utilisons pour les faisceaux de monoïdes la notation $f^{-1}$ pour désigner l'image inverse 
au sens des faisceaux de monoïdes et nous réservons la notation $f^*$ 
pour l'image inverse au sens des structures logarithmiques définie comme suit.
L'{\em image inverse} par $f$ d'une structure logarithmique $(\cM,\alpha)$ sur $Y$ 
est la structure logarithmique $(f^*(\cM),\beta)$ sur $X$ associée à la structure pré-logarithmique définie par
l'homomorphisme composé $f^{-1}(\cM)\rightarrow f^{-1}(\co_Y)\rightarrow \co_X$. Il résulte aussitôt de la définition que
l'homomorphisme canonique 
\begin{equation}
f^{-1}(\cM^\sharp)\rightarrow (f^*(\cM))^\sharp
\end{equation}
est un isomorphisme.

\subsection{}\label{higgs1-log6}
Un {\em schéma pré-logarithmique} (resp. {\em logarithmique}) est un triplet $(X,\cM_X,\alpha_X)$ 
formé d'un schéma $X$ et d'une structure pré-logarithmique (resp. logarithmique) $(\cM_X,\alpha_X)$ sur $X$. 
Lorsqu'il n'y a aucun risque d'ambiguïté, on se permettra d'omettre $\alpha_X$ et même $\cM_X$ des notations. 
Un morphisme de schémas pré-logarithmiques (resp. logarithmiques) $(X,\cM_X,\alpha_X)\rightarrow (Y,\cM_Y,\alpha_Y)$
est la donnée d'une paire $(f,f^\flat)$ formée d'un morphisme de schémas $f\colon X\rightarrow Y$ et d'un homomorphisme 
$f^\flat\colon f^{-1}(\cM_Y)\rightarrow  \cM_X$ tels que le diagramme 
\begin{equation}
\xymatrix{
{f^{-1}(\cM_Y)}\ar[rr]^{f^{-1}(\alpha_Y)}\ar[d]_{f^\flat}&&{f^{-1}(\co_Y)}\ar[d]\\
{\cM_X}\ar[rr]^{\alpha_X}&&{\co_X}}
\end{equation}
soit commutatif.

On dit qu'un schéma logarithmique $(X,\cM_X,\alpha_X)$ est {\em intègre} (resp. {\em saturé}) 
si $\cM_X$ est intègre (resp. saturé).  
On dit qu'un morphisme de schémas logarithmiques $f\colon (X,\cM_X,\alpha_X)\rightarrow (Y,\cM_Y,\alpha_Y)$
est {\em strict} si $(\cM_X,\alpha_X)$ est l'image inverse de  $(\cM_Y,\alpha_Y)$ par $f$, 
ou de façon équivalente, si l'homomorphisme canonique $f^{-1}(\cM_Y^\sharp)\rightarrow \cM_X^\sharp$ est un isomorphisme. 

\subsection{} \label{higgs1-not-mon1}\index{10512@$\bA_M$ ($M$ monoïde)}\index{10513@$e^x$ ($x$ élément d'un monoïde)}
Soit $M$ un monoïde. Pour tout entier $n\geq 1$, on désigne (abusivement) par $\varpi_n\colon M\rightarrow M$ 
l'homomorphisme de Frobenius d'ordre $n$ de $M$ ({\em i.e.}, l'élévation à la puissance $n$ dans $M$ en notation multiplicative). 

Pour tout anneau $R$, on désigne par $R[M]$ la $R$-algèbre de $M$ et par 
$e\colon M\rightarrow R[M]$ l'homomorphisme canonique, où $R[M]$ est considéré comme un
monoïde multiplicatif. Pour tout $x\in M$, on notera $e^x$ au lieu de $e(x)$. 

On désigne par $\bA_M$ le schéma $\Spec(\mZ[M])$ muni de la structure logarithmique associée à 
la structure pré-logarithmique définie par $e\colon M\rightarrow \mZ[M]$. 

Pour tout homomorphisme de monoïdes $\vartheta\colon M\rightarrow N$, on note 
$\bA_\vartheta\colon \bA_N\rightarrow \bA_M$ le morphisme de schémas logarithmiques associé. 

\subsection{}\label{higgs1-log7}
Soient $(X,\cM_X)$ un schéma logarithmique,  $M$ un monoïde,
$M_X$ le faisceau (étale) constant de monoïdes de valeur $M$ sur $X$. 
Les données suivantes sont équivalentes~:
\begin{itemize}
\item[(i)] un homomorphisme $\gamma\colon M\rightarrow \Gamma(X,\cM_X)$; 
\item[(ii)] un homomorphisme $\tgamma\colon M_X\rightarrow \cM_X$; 
\item[(iii)] un morphisme de schémas logarithmiques $\gamma^a\colon (X,\cM_X)\rightarrow \bA_M$. 
\end{itemize}
De plus, les conditions suivantes sont équivalentes~:
\begin{itemize}
\item[(a)] $\cM_X$ est associé à la structure pré-logarithmique qu'il induit sur $M_X$; 
\item[(b)] Le morphisme $\gamma^a\colon (X,\cM_X)\rightarrow \bA_M$ est strict. 
\end{itemize}
On dit alors que $(M,\gamma)$ est une {\em carte} pour $(X,\cM_X)$. 
On dit que la carte $(M,\gamma)$ est {\em cohérente}
(resp. {\em intègre}, resp. {\em fine}, resp. {\em saturée}, resp. {\em torique}) si le monoïde $M$ est de type fini (resp. intègre,
resp. fin, resp. saturé, resp. torique). 

\subsection{}\label{higgs1-log9}
Soit $f\colon (X,\cM_X)\rightarrow (Y,\cM_Y)$ un morphisme de schémas logarithmiques. Une {\em carte} pour $f$ est 
un triplet $((M,\gamma),(N,\delta),\vartheta\colon N\rightarrow M)$
formé d'une carte $(M,\gamma)$ pour $(X,\cM_X)$, d'une carte $(N,\delta)$ pour $(Y,\cM_Y)$
et d'un homomorphisme $\vartheta\colon N\rightarrow M$ tels que le diagramme 
\begin{equation}\label{higgs1-log9a}
\xymatrix{
N\ar[r]^-(0.5){\delta}\ar[d]_{\vartheta}&{\Gamma(Y,\cM_Y)}\ar[d]^{f^\flat}\\
M\ar[r]^-(0.5){\gamma}&{\Gamma(X,\cM_X)}}
\end{equation}
soit commutatif, ou ce qui revient au même que le diagramme associé 
de morphismes de schémas logarithmiques  \eqref{higgs1-not-mon1}
\begin{equation}\label{higgs1-log9b}
\xymatrix{
{(X,\cM_X)}\ar[r]^-(0.5){\gamma^a}\ar[d]_f&{\bA_M}\ar[d]^{\bA_\vartheta}\\
{(Y,\cM_Y)}\ar[r]^-(0.5){\delta^a}&{\bA_N}}
\end{equation}
soit commutatif. On dit que la carte $((M,\gamma),(N,\delta),\vartheta\colon N\rightarrow M)$ est {\em cohérente}
(resp. {\em fine}) si $M$ et $N$ sont de type fini (resp. fins).

\subsection{}\label{higgs1-log8}
Soit $(X,\cM_X)$ un schéma logarithmique. 
On dit que $(X,\cM_X)$ est {\em cohérent} si chaque point géométrique $\ox$ 
de $X$ admet un voisinage étale $U$ dans $X$ tel que $(U,\cM_X|U)$ admette une carte cohérente. 
On dit que $(X,\cM_X)$ est {\em fin} s'il est cohérent et intègre. 

Pour que $(X,\cM_X)$ soit fin (resp. fin et saturé), il faut et il suffit que tout point géométrique $\ox$ de $X$ 
admette un voisinage étale $U$ dans $X$ tel que $(U,\cM_X|U)$ admette une carte fine (resp. fine et saturée).

On dit que $(X,\cM_X)$ est {\em torique} si chaque point géométrique $\ox$ 
de $X$ admet un voisinage étale $U$ dans $X$ tel que $(U,\cM_X|U)$ admette une carte torique.

\begin{lem}[\cite{tsuji1} 1.3.2]\label{higgs1-log14}
Pour tout schéma logarithmique fin et saturé $(X,\cM_X)$, les conditions suivantes sont équivalentes~:
\begin{itemize}
\item[{\rm (i)}] Il existe une carte fine et saturée $\gamma\colon P\rightarrow \Gamma(X,\cM_X)$ 
pour $(X,\cM_X)$ \eqref{higgs1-log7} telle que l'homomorphisme composé 
\begin{equation}\label{higgs1-log14a}
P\rightarrow \Gamma(X,\cM_X) \rightarrow \Gamma(X,\cM_X)/\Gamma(X,\co_X^\times)
\end{equation} 
soit un isomorphisme. 

\item[{\rm (ii)}] Il existe une carte cohérente $\gamma\colon P\rightarrow \Gamma(X,\cM_X)$ 
pour $(X,\cM_X)$  telle que l'homomorphisme composé \eqref{higgs1-log14a} soit surjectif. 

\item[{\rm (iii)}] Le monoïde $\Gamma(X,\cM_X)/\Gamma(X,\co_X^\times)$ est de type fini et 
l'identité de $\Gamma(X,\cM_X)$ est une carte pour $(X,\cM_X)$. 
\end{itemize}
\end{lem}

L'implication (i)$\Rightarrow$(ii) est évidente. Montrons (ii)$\Rightarrow$(iii). Posons $Q=\Gamma(X,\cM_X)$
et notons $\cP$ et $\cQ$ les structures logarithmiques sur $X$ associées aux structures pré-logarithmiques 
définies par les homomorphismes $\gamma$ et l'identité de $\Gamma(X,\cM_X)$, respectivement,
et $\theta\colon \cP\rightarrow \cQ$ l'homomorphisme induit par $\gamma$. 
Comme l'homomorphisme composé $\cP\stackrel{\theta}{\rightarrow}\cQ\rightarrow \cM$ est un isomorphisme, 
il suffit de montrer que $\theta$ est surjectif. 
Soit $\ox$ un point géométrique de $X$. On désigne par $\beta\colon Q\rightarrow \cM_\ox$ l'homomorphisme canonique
et par $\alpha\colon P\rightarrow \cM_\ox$ l'homomorphisme composé $\beta\circ \gamma$. 
On a alors un digramme commutatif d'homomorphismes de monoïdes
\begin{equation}
\xymatrix{
{P/\alpha^{-1}(\co_{X,\ox}^\times)}\ar[r]^{\ogamma}\ar[d]&{Q/\beta^{-1}(\co_{X,\ox}^\times)}\ar[d]\\
{\cP_\ox/\co_{X,\ox}^\times}\ar[r]^{\theta^\sharp_\ox}&{\cQ_\ox/\co_{X,\ox}^\times}}
\end{equation}
où $\ogamma$ (resp. $\theta^\sharp_\ox$) est induit par $\gamma$ (resp. $\theta_\ox$)
et les flèches verticales sont les isomorphismes canoniques. 
Comme $\Gamma(X,\co_X^\times)\subset \beta^{-1}(\co_{X,\ox}^\times)$, $\ogamma$ est surjectif. 
Par suite, $\theta^\sharp_\ox$ est surjectif. Donc $\theta_\ox$ est surjectif~; d'où l'assertion. 
Montrons enfin (iii)$\Rightarrow$(i). Posons $P=\Gamma(X,\cM_X)/\Gamma(X,\co_X^\times)$. 
Comme $\Gamma(X,\cM_X)$ est saturé, $P$ est fin et saturé. 
Le sous-groupe de torsion de $P^\gp$ est donc contenu dans $P$. Comme $P$ est affuté, 
$P^\gp$ est un groupe abélien libre de type fini. Soit 
$\delta\colon P^\gp\rightarrow \Gamma(X,\cM_X)^\gp$ une section de l'homomorphisme canonique 
$\Gamma(X,\cM_X)^\gp\rightarrow P^\gp$.
La restriction de $\delta$ à $P$ induit un homomorphisme $\gamma\colon P\rightarrow \Gamma(X,\cM_X)$  
qui répond à la question en vertu de (\cite{tsuji1} 1.3.1).

\begin{lem}[\cite{tsuji1} 1.3.3]\label{higgs1-log15}
Soient $(X,\cM_X)$ un schéma logarithmique fine et saturé dont le schéma sous-jacent $X$ 
est noethérien, $x\in X$, $\gamma\colon P\rightarrow \Gamma(X,\cM_X)$ 
une carte fine et saturée pour $(X,\cM_X)$. Il existe alors un voisinage ouvert (de Zariski) $U$ de $x$ dans $X$ 
tel que pour tout voisinage ouvert (de Zariski) $V$ de $x$ dans $U$, 
le schéma logarithmique $(V,\cM_X|V)$ vérifie les conditions équivalentes de \eqref{higgs1-log14}.
\end{lem}

On désigne par $\Spec(P)$ l'ensemble des idéaux premiers de $P$ (cf. \cite{kato2} 5.1),  
par $P_X$ le faisceau (étale) constant de valeur $P$ sur $X$ et par $\tgamma \colon P_X\rightarrow \cM_X$
l'homomorphisme induit par $\gamma$. Pour tout point géométrique $\oy$ de $X$, on pose 
\begin{equation}
\fp_\oy=P-\tgamma_\oy^{-1}(\co^\times_{X,\oy})\subset (P_X)_\oy=P,
\end{equation}
qui est un idéal premier de $P$. Pour tout $t\in P$, on désigne par $Y_t$ le co-support  
de l'image de $t$ dans $\Gamma(X,\cM_X/\co^\times_X)$, ou ce qui revient au même puisque $\cM_X$ est intègre, 
le co-support  de l'image de $t$ dans $\Gamma(X,\cM^\gp_X/\co^\times_X)$ (\cite{sga4} IV 8.5.2). C'est un ouvert du
topos étale de $X$, que l'on identifie à un ouvert au sens habituel de $X$ (\cite{sga4} VIII 6.1). On pose $X_t=X-Y_t$. 
On vérifie aussitôt que la fibre de $X_t$ au-dessus d'un point 
géométrique $\oy$ de $X$ est non-vide si et seulement si $t\in \fp_\oy$. 
Pour tout $\fp\in \Spec(P)$, on pose $X_\fp=\cap_{t\in \fp}X_t$. 
Pour que la fibre de $X_\fp$ au-dessus d'un point 
géométrique $\oy$ de $X$ soit non-vide, il faut et il suffit que $\fp\subset \fp_\oy$.
Comme $X$ est noethérien 
et que $\Spec(P)$ est fini (cf. \cite{kato2} 5.5), il existe un voisinage ouvert (de Zariski) $U$ de $x$ dans $X$
vérifiant la propriété suivante. Pour tout $\fp\in \Spec(P)$ tel que $U_\fp=U\cap X_\fp\not=\emptyset$, 
$x$ appartient à toutes les composantes irréductibles de~$U_\fp$. Comme tout 
voisinage ouvert (de Zariski) $V$ de $x$ dans $U$ vérifie la même propriété, il suffit de montrer que $(U,\cM_X|U)$
vérifie les conditions équivalentes de \ref{higgs1-log14}.

Soit $\ox$ un point géométrique de $X$ au-dessus de $x$. Comme le composé des homomorphismes canoniques
\begin{equation}
P\rightarrow \Gamma(U,\cM_X)/\Gamma(U,\co_X^\times)\stackrel{u}{\rightarrow} \Gamma(U,\cM_X/\co_X^\times)  
\stackrel{v}{\rightarrow} \cM_{X,\ox}/\co_{X,\ox}^\times
\end{equation}
est surjectif et que $u$ est injectif, il suffit de montrer que $v$ est injectif. Soient $a,b\in \Gamma(U,\cM_X/\co_X^\times)$
telles que $v(a)=v(b)$, $\oy$ un point géométrique de $U$. 
Il existe alors un point géométrique $\oz$ de $U$ 
et une flèche de spécialisation $\varphi\colon \oz\rightsquigarrow \oy$ 
tels que l'image canonique $z$ de $\oz$ dans $U$ soit un point générique de $U_{\fp_\oy}$. 
On a $\tgamma_\oz=\varphi^*\circ \tgamma_\oy$, où $\varphi^*\colon \cM_{X,\oy}\rightarrow \cM_{X,\oz}$ est 
l'homomorphisme de spécialisation associé à $\varphi$ (\cite{sga4} VIII 7.7). 
On en déduit que $\fp_\oz\subset \fp_\oy$. Comme $z\in U_{\fp_\oy}$, on a $\fp_\oy\subset \fp_\oz$
et donc $\fp_\oy= \fp_\oz$. Le diagramme commutatif  
\begin{equation}
\xymatrix{
{P/\tgamma_\oy^{-1}(\co_{X,\oy}^\times)}\ar[r]^-(0.5){\tgamma_\oy}\ar@{=}[d]&{\cM_{X,\oy}/\co_{X,\oy}^\times}\ar[d]^{\varphi^*}\\
{P/\tgamma_\oz^{-1}(\co_{X,\oz}^\times)}\ar[r]^-(0.5){\tgamma_\oz}&{\cM_{X,\oz}/\co_{X,\oz}^\times}}
\end{equation}
où les flèches horizontales sont des isomorphismes montre alors que $\varphi^*\colon \cM_{X,\oy}/\co_{X,\oy}^\times
\rightarrow \cM_{X,\oz}/\co_{X,\oz}^\times$ est un isomorphisme. Comme $\ox$ est une spécialisation de $\oz$ dans $U$, 
les images canoniques $a_\oz$ et $b_\oz$ de $a$ et $b$ dans $\cM_{X,\oz}/\co_{X,\oz}^\times$ sont égales. 
Donc les images canoniques  $a_\oy$ et $b_\oy$ de $a$ et $b$ dans $\cM_{X,\oy}/\co_{X,\oy}^\times$ sont aussi égales. 
On en déduit que $a=b$; d'où l'assertion. 

\subsection{}\label{higgs1-log11}
Soit $f\colon (X,\cM_X)\rightarrow (Y,\cM_Y)$ un morphisme de schémas logarithmiques intègres.
On dit que $f$ est {\em intègre} (resp. {\em saturé}) si pour tout point géométrique $\ox$ de $X$, 
l'homomorphisme $\cM_{Y,f(\ox)}\rightarrow \cM_{X,\ox}$ est intègre (resp. saturé), ou ce qui revient au même,
si l'homomorphisme $\cM^\sharp_{Y,f(\ox)}\rightarrow \cM^\sharp_{X,\ox}$ est intègre (resp. saturé). 

\subsection{}\label{higgs1-log12}
Soit $f\colon X\rightarrow Y$ un morphisme de schémas logarithmiques intègres (resp. fins).
Pour que $f$ soit intègre, il faut et il suffit que pour tout schéma logarithmique intègre (resp. fin) $Z$ 
et tout morphisme $Z\rightarrow Y$,
le produit fibré dans la catégorie des schémas logarithmiques $Z\times_YX$ soit intègre (resp. fin) (\cite{kato1} 4.3.1).

\subsection{}\label{higgs1-log13}
Soit $f\colon X\rightarrow Y$ un morphisme intègre de schémas logarithmiques fins et saturés.
Pour que $f$ soit saturé, il faut et il suffit que pour tout schéma logarithmique fin et saturé $Z$ 
et tout morphisme $Z\rightarrow Y$,
le produit fibré dans la catégorie des schémas logarithmiques $Z\times_YX$ soit fin et saturé (\cite{tsuji4} II 2.13 page 24).

\subsection{}\label{higgs1-log85}\index{10521@$\Omega^1_{(X,\cM_X)/(Y,\cM_Y)}$}
Soit $f\colon (X,\cM_X)\rightarrow (Y,\cM_Y)$ un morphisme de schémas pré-logarithmiques. 
On pose 
\begin{equation}\label{higgs1-log85a}
\Omega^1_{(X,\cM_X)/(Y,\cM_Y)}=\frac{\Omega^1_{X/Y}\oplus(\co_X\otimes_\mZ\cM_X^\gp)}{\cF},
\end{equation}
où $\Omega^1_{X/Y}$ est le $\co_X$-module des $1$-différentielles relatives de $X$ sur $Y$
et $\cF$ est le sous-$\co_X$-module engendré localement par les sections de la forme 
\begin{itemize}
\item[(i)] $(d(\alpha_X(a)),0)-(0,\alpha_X(a)\otimes a)$ pour toute section locale $a$ de $\cM_X$;
\item[(ii)] $(0,1\otimes a)$ pour toute section locale $a$ de l'image de l'homomorphisme 
$f^\flat\colon f^{-1}(\cM_Y)\rightarrow \cM_X$. 
\end{itemize}
On note aussi 
\begin{equation}
d\colon \co_X\rightarrow \Omega^1_{(X,\cM_X)/(Y,\cM_Y)}
\end{equation}
le morphisme induit par la dérivation universelle $d\colon \co_X\rightarrow \Omega^1_{X/Y}$, et on désigne par 
\begin{equation}
d\log\colon \cM_X\rightarrow \Omega^1_{(X,\cM_X)/(Y,\cM_Y)}
\end{equation}
l'homomorphisme défini pour une section locale $a$ de $\cM_X$ par  
\begin{equation}
d\log(a)=1\otimes a.
\end{equation}
Alors le triplet $(\Omega^1_{(X,\cM_X)/(Y,\cM_Y)},d,d\log)$ est universel pour les dérivations logarithmiques. 
On appelle $\Omega^1_{(X,\cM_X)/(Y,\cM_Y)}$ le $\co_X$-module des {\em $1$-différentielles 
logarithmiques} de $(X,\cM_X)$ sur $(Y,\cM_Y)$ (ou de $f$). Il satisfait aux mêmes propriétés de fonctorialité
que le module $\Omega^1_{X/Y}$. 

Si $\cM_X^a$ (resp. $\cM^a_Y$) désigne la structure 
logarithmique sur $X$ (resp. $Y$) associée à $\cM_X$ (resp. $\cM_Y$), alors $f$ induit un morphisme 
$f^a\colon (X,\cM^a_X)\rightarrow (Y,\cM^a_Y)$, et on a un isomorphisme canonique 
\begin{equation}
\Omega^1_{(X,\cM_X)/(Y,\cM_Y)}\stackrel{\sim}{\rightarrow}\Omega^1_{(X,\cM^a_X)/(Y,\cM^a_Y)}.
\end{equation}

Si $f$ est un morphisme de schémas logarithmiques cohérents, alors $\Omega^1_{(X,\cM_X)/(Y,\cM_Y)}$
est un $\co_X$-module quasi-cohérent. Si, de plus, le morphisme de schémas sous-jacent à $f$ est localement 
de présentation finie, alors $\Omega^1_{(X,\cM_X)/(Y,\cM_Y)}$ est un $\co_X$-module de présentation finie. 

Si $f$ est un morphisme strict de schémas logarithmiques, 
le morphisme canonique 
\begin{equation}
\Omega^1_{X/Y}\rightarrow \Omega^1_{(X,\cM_X)/(Y,\cM_Y)}
\end{equation}
est un isomorphisme.

\subsection{}
Un morphisme de schémas logarithmiques $f\colon (X,\cM_X)\rightarrow (Y,\cM_Y)$ est une {\em immersion fermée}
(resp. {\em immersion fermée exacte})
si le morphisme de schémas sous-jacents $X\rightarrow Y$ est une immersion fermée et si le morphisme 
$f^*(\cM_Y)\rightarrow \cM_X$ est un épimorphisme (resp. un isomorphisme). 

\subsection{}\label{higgs1-log10}
Considérons un diagramme commutatif de morphismes de schémas logarithmiques 
\begin{equation}
\xymatrix{
{(Z',\cM_{Z'})}\ar[r]^{u'}\ar[d]_j&{(X,\cM_X)}\ar[d]^f\\
{(Z,\cM_Z)}\ar[r]^g\ar@{.>}[ru]^u&{(Y,\cM_Y)}}
\end{equation}
où $Z'$ est un sous-schéma fermé de $Z$ défini par un idéal 
$\cI$ de $\co_Z$ de carré nul et $j$ est une immersion fermée exacte. 
On désigne par $\rP_f(j,u')$ l'ensemble des 
$(Y,\cM_Y)$-morphismes $u\colon (Z,\cM_Z)\rightarrow (X,\cM_X)$ tels que $u'=u\circ j$. 
Alors $\rP_f(j,u')$ est un pseudo-torseur sous 
\begin{equation}
\Hom_{\co_{Z'}}(u'^*\Omega^1_{(X,\cM_X)/(Y,\cM_Y)},\cI).
\end{equation}
Plus précisément  (\cite{kato1} 3.9), 
si $\rP_f(j,u')$ est non vide, tout élément $u\in \rP_f(j,u')$ détermine uniquement un isomorphisme 
\begin{equation}
\varphi_u\colon \rP_f(j,u') \stackrel{\sim}{\rightarrow} \Hom_{\co_{Z'}}(u'^*\Omega^1_{(X,\cM_X)/(Y,\cM_Y)},\cI) 
\end{equation}
tel que pour tout $v\in \rP_f(j,u')$, toute section locale $a$  de $\co_X$ et toute section locale $b$ de $\cM_X$, on ait
\begin{eqnarray}
\varphi_u(v)(u'^*(da))&=&v^\sharp(v^{-1}(a))-u^\sharp(u^{-1}(a)),\\
\varphi_u(v)(u'^*(d\log b))&=& \beta-1,
\end{eqnarray}
où $\beta$ est l'unique section locale de $1+\cI\subset \co_Z^\times\subset \cM_Z$ telle que  
$v^\flat(v^{-1}(b))=\beta \cdot u^\flat(u^{-1}(b))$.

\subsection{}
Soit $f\colon X\rightarrow Y$ un morphisme de schémas logarithmiques. On dit que $f$  est {\em formellement lisse} 
(resp. {\em formellement étale}) si pour tout schéma logarithmique $Y$ dont le schéma sous-jacent est 
affine, toute immersion fermée exacte d'idéal nilpotent $j\colon Y'_0\rightarrow Y'$, 
et tout morphisme  $Y'\rightarrow Y$, l'application 
\begin{equation}
\Hom_Y(Y',X)\rightarrow \Hom(Y'_0,X)
\end{equation}
déduite de $j$, est surjective (resp. bijective). On dit que $f$  est {\em lisse} 
(resp. {\em étale}) s'il est formellement lisse (resp. formellement étale), 
si les schémas logarithmiques $X$ et $Y$ sont cohérents et si le morphisme de schémas
sous-jacent à $f$ est localement de présentation finie. 

\subsection{}\label{higgs1-log105}
Soient $f\colon (X,\cM_X)\rightarrow (Y,\cM_Y)$ un morphisme de schémas logarithmiques fins,
$(N,\delta)$ une carte fine pour $(Y,\cM_Y)$ \eqref{higgs1-log7}. D'après (\cite{kato1} 3.5), pour que $f$ soit lisse (resp. étale), 
il faut et il suffit que localement pour la topologie étale sur $X$, $f$
admette une carte fine $((M,\gamma),(N,\delta),\vartheta\colon N\rightarrow M)$ \eqref{higgs1-log9} satisfaisant aux conditions suivantes~:
\begin{itemize}
\item[(i)] le noyau et le sous-groupe de torsion du conoyau (resp. le noyau et le conoyau) de l'homomorphisme
$\vartheta^\gp\colon N^\gp\rightarrow M^\gp$ sont finis d'ordres inversibles dans $X$;
\item[(ii)] le morphisme induit $X\rightarrow Y\times_{\bA_N}\bA_M$ \eqref{higgs1-log9b} est étale dans le sens classique.
\end{itemize}

\section{Le théorème de presque pureté de Faltings}\label{higgs1-pur}

\subsection{}\label{higgs1-dlog0}\index{10601@$(S,\cM_S)$}
On pose $S=\Spec(\co_K)$ et on note $\eta$ (resp. $s$) son point générique 
(resp. fermé) et $\oeta$ le point géométrique générique correspondant à $\oK$ \eqref{higgs1-not1}. 
On munit $S$ de la structure logarithmique $\cM_S$ définie par son point fermé, 
autrement dit, $\cM_S=j_*(\co_\eta^\times)\cap \co_S$, où $j\colon \eta\rightarrow S$ est l'injection canonique.
On fixe une uniformisante $\pi$ de $\co_K$ et on désigne par $\iota\colon \mN\rightarrow \Gamma(S,\cM_S)$
l'homomorphisme défini par $\iota(1)=\pi$, qui est une carte pour $(S,\cM_S)$.

\subsection{}\label{higgs1-dlog1}\index{10602@$(X,\cM_X)$}\index{10603@$R$}
\index{10604@$\vartheta\colon \mN\rightarrow P$, $\lambda=\vartheta(1)$}\index{10605@$L$, $L_\lambda$}\index{10605@$X^\circ$}
\index{10606@$\alpha\colon P\rightarrow R$}
Dans toute la suite de cet article, on considère un morphisme de schémas logarithmiques
\begin{equation}
f\colon (X,\cM_X)\rightarrow (S,\cM_S),
\end{equation} 
une carte {\em torique} $(P,\gamma)$ pour $(X,\cM_X)$ \eqref{higgs1-log7} et 
un homomorphisme de monoïdes $\vartheta\colon \mN\rightarrow P$ tels que les conditions suivantes soient remplies~:
\begin{itemize}
\item[(C$_1$)] Le schéma $X=\Spec(R)$ est affine et connexe. 
\item[(C$_2$)] Le schéma $X_s=X\times_Ss$ est non-vide. 
\item[(C$_3$)] Le triplet $((P,\gamma),(\mN,\iota),\vartheta)$ est une carte pour $f$ \eqref{higgs1-log9}, autrement dit,
le diagramme d'homomorphismes de monoïdes 
\begin{equation}\label{higgs1-dlog1aa}
\xymatrix{
P\ar[r]^-(0.5)\gamma&{\Gamma(X,\cM_X)}\\
\mN\ar[r]^-(0.5){\iota}\ar[u]^\vartheta&{\Gamma(S,\cM_S)}\ar[u]_{f^\flat}}
\end{equation}
est commutatif, ou ce qui revient au même le diagramme associé de morphismes de schémas logarithmiques 
\begin{equation}\label{higgs1-dlog1a}
\xymatrix{
{(X,\cM_X)}\ar[r]^-(0.5){\gamma^a}\ar[d]_f&{\bA_P}\ar[d]^{\bA_\vartheta}\\
{(S,\cM_S)}\ar[r]^-(0.5){\iota^a}&{\bA_\mN}}
\end{equation}
est commutatif.
\item[(C$_4$)] L'homomorphisme $\vartheta$ est saturé \eqref{higgs1-log111}. 
\item[(C$_5$)] L'homomorphisme $\vartheta^\gp\colon \mZ\rightarrow P^\gp$ est injectif, 
le sous-groupe de torsion de $\coker(\vartheta^\gp)$ est d'ordre premier à $p$ et le morphisme de schémas usuels
\begin{equation}\label{higgs1-dlog1b}
X\rightarrow S\times_{\bA_\mN}\bA_P
\end{equation}
déduit de \eqref{higgs1-dlog1a} est étale.  
\item[(C$_6$)] Posons $\lambda=\vartheta(1)\in P$, 
\begin{eqnarray}
L&=&\Hom_{\mZ}(P^\gp,\mZ),\label{higgs1-dlog1f}\\
\rH(P)&=&\Hom(P,\mN).\label{higgs1-dlog1fg}
\end{eqnarray} 
On notera que $\rH(P)$ est un monoïde fin, saturé et affûté et que l'homomorphisme canonique 
$\rH(P)^\gp\rightarrow \Hom((P^\sharp)^\gp,\mZ)$ est un isomorphisme (\cite{ogus} I 2.2.3). 
On suppose qu'il existe $h_1,\dots,h_r\in \rH(P)$, qui sont $\mZ$-linéairement indépendants dans $L$, tels que  
\begin{equation}
\ker(\lambda)\cap \rH(P)=\{\sum_{i=1}^ra_ih_i | \ (a_1,\dots,a_r)\in \mN^r\},
\end{equation}
où l'on considère $\lambda$ comme un homomorphisme $L\rightarrow \mZ$. 
\end{itemize}

\vspace{2mm}

On rappelle que $P^\gp$ est un $\mZ$-module libre de type fini. On notera qu'on a  
\begin{equation} \label{higgs1-dlog1g}
S\times_{\bA_\mN}\bA_P=\Spec(\co_K[P]/(\pi-e^\lambda)),
\end{equation} 
où $\lambda=\vartheta(1)$ (cf. \ref{higgs1-not-mon1}). On pose 
\begin{equation}\label{higgs1-dlog1ff}
L_\lambda=\Hom_{\mZ}(P^\gp/\lambda\mZ,\mZ),
\end{equation}
que l'on identifie au noyau de l'homomorphisme $L\rightarrow \mZ$ donné par 
$y\mapsto \langle y,\lambda\rangle$  \eqref{higgs1-dlog1f}. On note $d$ le rang de $L_\lambda$. 
On désigne par $X^\circ$ le sous-schéma ouvert maximal de $X$ où la structure logarithmique $\cM_X$
est triviale. On a $X^\circ=X\times_{\bA_P}\bA_{P^\gp}$, qui est un sous-schéma ouvert affine 
de $X_\eta$. 

On note $\alpha\colon P\rightarrow R$ l'homomorphisme induit par la carte $(P,\gamma)$. 
On observera que pour tout $t\in P$, $\alpha(t)$ est inversible sur $X^\circ$; en particulier, $\alpha(t)\not=0$.

\begin{prop}\label{higgs1-dlog19}
{\rm (i)}\ Le morphisme $f$ est lisse et saturé. 

{\rm (ii)}\ Le schéma $X$ est intègre, normal, Cohen-Macaulay et plat sur $S$.

{\rm (iii)}\  Le schéma $X\otimes_{\co_K}\co_{\oK}$ est normal.

{\rm (iv)}\ Le schéma $X\otimes_{\co_K}\ok$ est réduit.

{\rm (v)}\ Le schéma (usuel) $X\times_S\eta$ est lisse sur $\eta$, $X^\circ\times_S\eta$ est l'ouvert 
complémentaire dans $X\times_S\eta$ d'un diviseur à croisements normaux stricts $D$ et  $\cM_X|(X\times_S\eta)$ 
est la structure logarithmique sur $X\times_S\eta$ définie par $D$. 
\end{prop}

(i) Cela résulte aussitôt de \ref{higgs1-log105} et (\cite{tsuji4} chap.~II 3.5). 

(ii) Les trois dernières propriétés résultent de (i), (\cite{kato1} 4.5) et (\cite{kato2} 8.2 et 4.1); cf. aussi 
(\cite{tsuji1} 1.5.1). 
Comme de plus $X$ est noethérien et connexe (C$_1$), il est alors intègre. 

(iii) Pour toute extension finie $K'$ de $K$,
d'anneau de valuation $\co_{K'}$, munissons $S'=\Spec(\co_{K'})$ de la structure logarithmique $\cM_{S'}$ 
définie par son point fermé et notons $f'\colon (X',\cM_{X'})\rightarrow (S',\cM_{S'})$ le morphisme déduit de 
$f$ par changement de base dans la catégorie des schémas logarithmiques 
par le morphisme canonique $(S',\cM_{S'})\rightarrow (S,\cM_S)$. 
Alors $f'$ est lisse et saturé \eqref{higgs1-log13}, et par suite $X'$ est normal en vertu de (ii). 
Comme $X'=X\times_SS'$, l'assertion s'ensuit
par passage à la limite inductive sur les extensions finies de $K$ contenues dans $\oK$ (\cite{ega1n} 0.6.5.12(ii)).

(iv) Cela résulte de (i)  et (\cite{tsuji4} chap.~II 4.2). 

(v) Soit $F$ la face de $P$ engendrée par $\lambda$, c'est-à-dire l'ensemble des éléments 
$x\in P$ tels qu'il existe $y\in P$ et $n\in \mN$ tels que $x+y=n\lambda$ (\cite{ogus} I 1.4.2).
On note $F^{-1}P$ la localisation de $P$ par $F$ (\cite{ogus} I 1.4.4). 
Il résulte aussitôt des propriétés universelles des localisations de monoïdes et d'anneaux que 
l'homomorphisme canonique $\mZ[P]\rightarrow \mZ[F^{-1}P]$ induit un isomorphisme
\begin{equation}\label{higgs1-dlog19a}
\mZ[P]_{\lambda}\stackrel{\sim}{\rightarrow} \mZ[F^{-1}P].
\end{equation}
Soient $P/F$ (resp. $\Lambda$) le conoyau dans la catégorie
des monoïdes de l'injection canonique $F\rightarrow P$ 
(resp. de l'homomorphisme $\vartheta\colon \mN\rightarrow P$) (cf. \cite{ogus} I 1.1.6). 
On a des isomorphismes canoniques 
\begin{equation}\label{higgs1-dlog19c}
\Lambda^\sharp \stackrel{\sim}{\rightarrow}  P/F\stackrel{\sim}{\rightarrow}(F^{-1}P)^\sharp.
\end{equation}
L'homomorphisme canonique 
\begin{equation}\label{higgs1-dlog19d}
\Hom(P/F,\mN)\rightarrow \ker(\lambda)\cap \rH(P)
\end{equation}
est un isomorphisme. 
En tant que somme amalgamée de l'homomorphisme saturé $\vartheta$ et de $\mN\rightarrow 0$, 
$\Lambda$ est saturé \eqref{higgs1-log111}. Par suite, $P/F$ est saturé \eqref{higgs1-dlog19c}. Donc en vertu de (\cite{ogus} I 2.2.3), 
on a un isomorphisme canonique 
\begin{equation}\label{higgs1-dlog19f} 
P/F\stackrel{\sim}{\rightarrow} \Hom(\Hom(P/F,\mN),\mN).
\end{equation}
La condition \ref{higgs1-dlog1}(C$_6$) implique alors que $P/F$ est un monoïde libre de type fini. 
Par suite, il existe un homomorphisme $P/F\rightarrow F^{-1}P$ qui relève 
l'isomorphisme canonique $P/F\stackrel{\sim}{\rightarrow}(F^{-1}P)^\sharp$ \eqref{higgs1-dlog19c},
de sorte que l'homomorphisme induit $P/F\rightarrow \mZ[F^{-1}P]$ est une carte pour $\bA_{F^{-1}P}$ (\cite{tsuji1} 1.3.1). 
On en déduit par \eqref{higgs1-dlog19a} une carte 
\begin{equation}
(X\times_S\eta,\cM_X|(X\times_S\eta))\rightarrow \bA_{P/F}. 
\end{equation}
D'autre part, le schéma logarithmique $(X,\cM_X)$ est régulier en vertu de (i) et 
(\cite{kato2} 8.2); cf. aussi (\cite{niziol} 2.3) et la preuve de (\cite{tsuji1} 1.5.1). 
Il résulte alors de (\cite{ega4} 0.16.3.7 et 0.17.1.7) et des définitions (\cite{kato2} 2.1 et \cite{niziol} 2.2) 
que le schéma $X\times_S\eta$ est régulier et donc lisse sur $\eta$, que 
$X^\circ\times_S\eta$ est l'ouvert complémentaire dans $X\times_S\eta$
d'un diviseur à croisements normaux stricts $D$ et  que $\cM_X|(X\times_S\eta)$ 
est la structure logarithmique sur $X\times_S\eta$ définie par $D$.

\subsection{}\label{higgs1-dlog2}\index{10609@$K_n$, $\pi_n$}
Pour tout entier $n\geq 1$, on pose
\begin{eqnarray}\label{higgs1-dlog2a}
\co_{K_n}=\co_K[\zeta]/(\zeta^{n}-\pi),
\end{eqnarray}
qui est un anneau de valuation discrète. On note $K_n$ le corps des fractions de $\co_{K_n}$
et $\pi_n$ la classe de $\zeta$ dans $\co_{K_n}$, qui est une uniformisante de $\co_{K_n}$.  
On pose $S_n=\Spec(\co_{K_n})$
que l'on munit de la structure logarithmique $\cM_{S_n}$ définie par son point fermé.
On désigne par $\tau_n\colon (S_n,\cM_{S_n})\rightarrow (S,\cM_S)$ le morphisme canonique et
par $\iota_n\colon \mN\rightarrow \Gamma(S_n,\cM_{S_n})$
l'homomorphisme défini par $\iota_n(1)=\pi_n$.
On notera que $\iota_n$ est une carte pour $(S_n,\cM_{S_n})$ et que le diagramme 
\begin{equation}\label{higgs1-dlog2b}
\xymatrix{
{(S_n,\cM_{S_n})}\ar[r]^-(0.5){\iota^a_n}\ar[d]_{\tau_n}&{\bA_\mN}\ar[d]^{\bA_{\varpi_n}}\\
{(S,\cM_S)}\ar[r]^-(0.5){\iota^a}&{\bA_\mN}}
\end{equation}
est cartésien (cf. \ref{higgs1-not-mon1}). 

Pour tous entiers $m,n\geq 1$, on a $\varpi_{mn}=\varpi_m\circ \varpi_n$. 
On en déduit un morphisme canonique 
\begin{equation}\label{higgs1-dlog2c}
\tau_{m,n}\colon (S_{mn},\cM_{S_{mn}})\rightarrow (S_n,\cM_{S_n})
\end{equation}
tel que $\tau_{mn}=\tau_n\circ \tau_{m,n}$. Pour tous entiers $r,m,n\geq 1$, on a 
$\tau_{rm,n}=\tau_{m,n}\circ\tau_{r,mn}$. 
Donc les schémas logarithmiques $(S_n,\cM_{S_n})$ pour $n\geq 1$ forment un système projectif cofiltrant
indexé par l'ensemble $\mZ_{\geq 1}$ ordonné par la relation de divisibilité. 

\subsection{}\label{higgs1-dlog3}\index{10610@$(X_n,\cM_{X_n})$, $A_n$}
Pour tout entier $n\geq 1$, on pose 
\begin{equation}\label{higgs1-dlog3a}
(X_n,\cM_{X_n})=(X,\cM_{X})\times_{\bA_P, \bA_{\varpi_n}}\bA_P
\end{equation}
et on note $\rho_n\colon (X_n,\cM_{X_n})\rightarrow (X,\cM_{X})$ la projection canonique (cf. \ref{higgs1-not-mon1}). 
On note aussi (abusivement)  $\rho_n\colon X_n\rightarrow X$ le morphisme de schémas 
sous-jacent à $\rho_n$. Alors $\rho_n$ est fini, le schéma $X_n$ est affine d'anneau 
\begin{equation}\label{higgs1-dlog3b}
A_n=R\otimes_{\mZ[P],\varpi_n}\mZ[P],
\end{equation}
et la projection canonique $(X_n,\cM_{X_n})\rightarrow \bA_P$ est stricte. Comme le diagramme \eqref{higgs1-dlog2b}
est cartésien, il existe un unique morphisme
\begin{equation}\label{higgs1-dlog3cc}
f_n\colon (X_n,\cM_{X_n})\rightarrow (S_n,\cM_{S_n}),
\end{equation} 
qui s'insère dans le diagramme commutatif
\begin{equation}\label{higgs1-dlog3d}
\xymatrix{
{(X_n,\cM_{X_n})}\ar[rrr]\ar[ddd]_{\rho_n}\ar[rd]^{f_n}&&&{\bA_P}\ar[ld]_{\bA_{\vartheta}}\ar[ddd]^{\bA_{\varpi_n}}\\
&{(S_n,\cM_{S_n})}\ar[r]\ar[d]&{\bA_\mN}\ar[d]^{\bA_{\varpi_n}}&\\
&{(S,\cM_{S})}\ar[r]&{\bA_\mN}&\\
{(X,\cM_{X})}\ar[rrr]\ar[ru]^f&&&{\bA_P}\ar[lu]_{\bA_\vartheta}}
\end{equation}

On désigne par $X^\circ_n$ le sous-schéma ouvert maximal de $X_n$ où la structure logarithmique $\cM_{X_n}$ 
est triviale. On a $X^\circ_n=X_n\times_{\bA_P}\bA_{P^\gp}$ et $f_n(X_n^\circ)=\Spec(K_n)$. 

Pour tous entiers $m,n\geq 1$, on a $\varpi_{mn}=\varpi_m\circ \varpi_n$.
On en déduit un morphisme canonique 
\begin{equation}\label{higgs1-dlog3c}
\rho_{m,n}\colon (X_{mn},\cM_{X_{mn}})\rightarrow (X_n,\cM_{X_n})
\end{equation}
tel que $\rho_{mn}=\rho_n\circ \rho_{m,n}$. Pour tous entiers $r, m, n\geq 1$, on a 
$\rho_{rm,n}=\rho_{m,n}\circ\rho_{r,mn}$.
Donc les schémas logarithmiques $(X_n,\cM_{X_n})$ forment un système projectif cofiltrant
indexé par l'ensemble $\mZ_{\geq 1}$ ordonné par la relation de divisibilité. 

\begin{prop}\label{higgs1-dlog4}
Soit $n$ un entier $\geq 1$. 
\begin{itemize}
\item[{\rm (i)}] Le morphisme de schémas usuels $X_n\rightarrow S_n\times_{\bA_\mN}\bA_P$ 
déduit de \eqref{higgs1-dlog3d} est étale, 
et le morphisme de schémas logarithmiques $f_n$ est lisse et saturé. 
\item[{\rm (ii)}] Le schéma $X_n$ est normal, Cohen-Macaulay et plat sur $S_n$. 
\item[{\rm (iii)}] Le schéma $X_n\otimes_{\co_{K_n}}\co_{\oK}$ est normal et n'a qu'un nombre fini de points génériques~;
en particulier, $X_n\otimes_{\co_{K_n}}\co_{\oK}$ est une somme finie de schémas intègres et normaux.
\item[{\rm (iv)}] Le morphisme $X_n\otimes_{\co_{K_n}}K_n\rightarrow X\otimes_{\co_K}K$ déduit de $\rho_n$ est plat. 
Si, de plus, $X$ est régulier, $\rho_n$ est plat. 
\item[{\rm (v)}] Si $n$ est une puissance de $p$, $X_n$ est intègre, et 
l'image inverse de toute composante connexe de $X\otimes_{\co_K}\co_\oK$ 
par le morphisme canonique $X_n\otimes_{\co_{K_n}}\co_\oK\rightarrow X\otimes_{\co_K}\co_\oK$
est intègre.
\item[{\rm (vi)}] Les carrés du diagramme commutatif canonique
\begin{equation}\label{higgs1-dlog4a}
\xymatrix{
X_n\ar[d]_{\rho_n}&{X_n^\circ}\ar[l]\ar[d]\ar[r]&{\bA_{P^\gp}}\ar[d]^{\bA_{\varpi_n}}\\
X&{X^\circ}\ar[l]\ar[r]&{\bA_{P^\gp}}}
\end{equation}
sont cartésiens. En particulier, $X_n^\circ$ est un espace principal homogène pour la topologie étale, 
au-dessus de $X^\circ$, sous le groupe $\Hom_\mZ(P^\gp,\mu_{n}(\oK))$.
\end{itemize}
\end{prop}

(i) En effet, les carrés du diagramme commutatif de morphismes de schémas usuels
\begin{equation}\label{higgs1-dlog4b}
\xymatrix{
{X_n}\ar[r]\ar[d]_{\rho_n}&{S_n\times_{\bA_\mN}\bA_P}\ar[r]\ar[d]&{\bA_P}\ar[d]^{\bA_{\varpi_n}}\\
X\ar[r]&{S\times_{\bA_\mN}\bA_P}\ar[r]&{\bA_P}}
\end{equation}
déduit de \eqref{higgs1-dlog3d} sont cartésiens.

(ii) Cela résulte de (i) par la même preuve que \ref{higgs1-dlog19}(ii). 

(iii) Il résulte de (i) par la même preuve que \ref{higgs1-dlog19}(iii) que $X_n\otimes_{\co_{K_n}}\co_{\oK}$ est normal. 
D'autre part, comme $X_n\otimes_{\co_{K_n}}\co_{\oK}$ est plat sur $\co_\oK$, ses points génériques sont 
les points génériques du schéma $X_n\otimes_{\co_{K_n}}\oK$, qui est noethérien. Par suite, 
l'ensemble des points génériques de $X_n\otimes_{\co_{K_n}}\co_{\oK}$ est fini.
La dernière assertion s'ensuit compte tenu de (\cite{ega1n} 0.2.1.6). En effet, 
$X_n\otimes_{\co_{K_n}}\co_\oK$ étant normal, deux composantes irréductibles distinctes 
de $X_n\otimes_{\co_{K_n}}\co_\oK$ ne se rencontrent pas. 

(iv) Cela résulte de (\cite{ega4} 0.17.3.5) car $X\otimes_{\co_K}K$ est régulier en vertu de \ref{higgs1-dlog19}(v), 
$X_n$ est Cohen-Macaulay d'après (ii) et $\rho_n$ est fini.

(v) D'après \eqref{higgs1-dlog4b},  on a un diagramme cartésien de $S$-morphismes
\begin{equation}\label{higgs1-dlog4ee}
\xymatrix{
{X_n}\ar[r]\ar[d]_{\rho_n}&{\Spec(\co_{K_n}[P]/(\pi_n-e^\lambda))}\ar[d]^{\gamma_n}\\
{X}\ar[r]&{\Spec(\co_K[P]/(\pi-e^\lambda))}}
\end{equation}
où $\gamma_n$ est induit par l'homomorphisme $\varpi_n$. D'après \ref{higgs1-dlog1}(C$_2$), il existe $x\in X_s$.
Comme $n$ est une puissance de $p$, le morphisme $\rho_n\otimes_{\co_K}k$ est un homéomorphisme universel.  
Donc $\rho_n^{-1}(x)$ contient un seul point que l'on note $x_n$. D'autre part, $X_n$ étant plat sur $S_n$, 
tout point générique de $X_n$ est un point générique de $X_n\otimes_{\co_{K_n}}K_n$;
donc son image par $\rho_n$ est le point générique de $X$ en vertu de (iv). 
Comme $\rho_n$ est fermé, on en déduit que $x_n$ est une spécialisation de tous 
les points génériques de $X_n$. Comme $X_n$ est noethérien et normal d'après (ii),  il est intègre.

La preuve de la seconde assertion est similaire à celle de la première. 
En effet, on déduit facilement de \eqref{higgs1-dlog4ee}  un diagramme cartésien de $\co_{\oK}$-morphismes
\begin{equation}\label{higgs1-dlog4e}
\xymatrix{
{X_n\otimes_{\co_{K_n}}\co_{\oK}}\ar[r]\ar[d]_{\beta_n}&{\Spec(\co_{\oK}[P]/(\pi_n-e^\lambda))}\ar[d]^{\alpha_n}\\
{X\otimes_{\co_K}\co_{\oK}}\ar[r]&{\Spec(\co_{\oK}[P]/(\pi-e^\lambda))}}
\end{equation}
D'après (iii), $X\otimes_{\co_K}\co_{\oK}$ est une somme finie de schémas intègres et normaux~;
en particulier, il n'a qu'un nombre fini de composantes connexes, qui sont alors ouvertes.
Comme $X$ est intègre et plat sur $S$ d'après \ref{higgs1-dlog19}(ii), 
$G_K$ agit transitivement sur l'ensemble des composantes connexes de $X\otimes_{\co_K}\co_{\oK}$. 
Soient $Y$ une composante connexe de $X\otimes_{\co_K}\co_{\oK}$, $Y_n=\beta_n^{-1}(Y)$. 
Compte tenu de \ref{higgs1-dlog1}(C$_2$), $Y\otimes_{\co_\oK}\ok$ est non-vide. 
Soit $y\in Y\otimes_{\co_\oK}\ok$. 
Comme $n$ est une puissance de $p$, le morphisme $\alpha_n\otimes_{\co_\oK}\ok$ est un homéomorphisme universel.  
Donc $\beta_n^{-1}(y)$ contient un seul point que l'on note $y_n$. D'autre part, tout point
générique de $Y_n$ est au-dessus de l'unique point générique de $X_n$ et est donc au-dessus du point générique de $X$. 
Par suite, l'image de tout point générique de $Y_n$ par $\beta_n$ est le point générique de $Y$. 
Comme $\beta_n$ est fermé, on en déduit que $y_n$ est une spécialisation de tous 
les points génériques de $Y_n$. D'autre part, $Y_n$ est une réunion disjointe finie de sous-schémas ouverts et fermés 
qui sont intègres et normaux d'après (iii). Il est donc intègre. 

(vi) Montrons que le diagramme canonique 
\begin{equation}\label{higgs1-dlog4c}
\xymatrix{
P\ar[r]^{\varpi_n}\ar[d]&P\ar[d]\\
{P^\gp}\ar[r]^{\varpi_n}&{P^\gp}}
\end{equation}
est co-cartésien. En effet, la somme amalgamée $P^\gp\oplus_{P,\varpi_n}P$
est le quotient de $P^\gp\oplus P$ par la relation de congruence $E$ définie par l'ensemble  
des paires $((y,x),(y',x'))$ d'éléments de $P^\gp\oplus P$ telles qu'il existe $z,z'\in P$
tels que $y+z=y'+z'\in P^\gp$ et $x+nz=x'+nz'\in P$ (\cite{ogus} I 1.1.5). 
Il suffit donc de montrer que $E$ est la relation de congruence définie par l'homomorphisme
\begin{equation}
P^\gp\oplus P\rightarrow  P^\gp,\ \ \ (y,x)\rightarrow x-ny.
\end{equation}
Si $((y,x),(y',x'))\in E$ alors $x-ny=x'-ny'$. Inversement, supposons $x=x'+n(y-y')\in P^\gp$.
Comme $P$ est intègre, il existe $z,z'\in P$ tels que $y+z=y'+z'\in P^\gp$; on a donc 
$x+nz=x'+nz'\in P$, ce qui prouve l'assertion. 
Par suite, le diagramme 
\begin{equation}\label{higgs1-dlog4d}
\xymatrix{
{\bA_{P^\gp}}\ar[r]\ar[d]_{\bA_{\varpi_n}}&{\bA_P}\ar[d]^{\bA_{\varpi_n}}\\
{\bA_{P^\gp}}\ar[r]&{\bA_P}}
\end{equation}
induit par \eqref{higgs1-dlog4c} est cartésien. Donc les carrés du diagramme \eqref{higgs1-dlog4a} sont cartésiens.   
La seconde assertion résulte de la première et du fait que le noyau de l'isogénie étale 
$\bA_{\varpi_n}\otimes_\mZ K$ de $\bA_{P^\gp}\otimes_\mZ K$ 
correspond au $\mZ[G_K]$-module $\Hom_\mZ(P^\gp,\mu_{n}(\oK))$.

\subsection{}\label{higgs1-dlog5}\index{10614@$R_n$, $R_\infty$, $R_{p^\infty}$, $\oR$}\index{10615@$F$, $F_n$, $F_\infty$, $F_{p^\infty}$, $\oF$}
\index{10616@$B_n$, $B_\infty$, $H_n$, $H_\infty$}
On note $\kappa$ le point générique de $X$ et $F$ le corps résiduel de $X$ en $\kappa$
({\em i.e.}, le corps des fractions de $R$). 
Dans la suite de cet article, on se donne un point géométrique 
$\tkappa$ de $X\otimes_{\co_K}\co_{\oK}$ au-dessus de $\kappa$, autrement dit, le spectre d'une  
extension séparablement close $F^a$ de $F$ contenant $\oK$. 
On désigne par $\oF$ l'union des extensions finies $N$ de $F$, contenues dans $F^a$, telles que  la clôture
intégrale de $R$ dans $N$ soit étale au-dessus de $X^\circ$, et par $\oR$ la clôture intégrale de $R$ dans $\oF$.

Pour tous entiers $m,n\geq 1$, le morphisme $\rho_{m,n}\colon X_{mn}\rightarrow X_n$ 
est fini et surjectif. Par suite, en vertu de (\cite{ega4} 8.3.8(i)), il existe un $X$-morphisme
\begin{equation}\label{higgs1-dlog5b}
\tkappa\rightarrow \underset{\underset{n\geq 1}{\longleftarrow}}{\lim}\ X_{n},
\end{equation} 
où la limite projective est indexée par l'ensemble $\mZ_{\geq 1}$ ordonné par la relation de divisibilité.
{\em On se donne un tel morphisme qu'on suppose fixé dans toute la suite de cet article}. 
L'ensemble des entiers $n!$, pour $n\geq 0$, étant cofinal dans $\mZ_{\geq 1}$ pour la relation de divisibilité, 
il revient au même de se donner un morphisme 
\begin{equation}\label{higgs1-dlog5bb}
\tkappa\rightarrow \underset{\underset{n\geq 0}{\longleftarrow}}{\lim}\ X_{n!},
\end{equation} 
où la limite projective est indexée par l'ensemble $\mN$ ordonné par la relation d'ordre habituelle. 

D'après \ref{higgs1-dlog4}(vi), le morphisme \eqref{higgs1-dlog5b} se factorise à travers un $X$-morphisme 
\begin{equation}\label{higgs1-dlog5c}
\Spec(\oR)\rightarrow \underset{\underset{n\geq 1}{\longleftarrow}}{\lim}\ X_{n}. 
\end{equation} 
On en déduit un système inductif de $R$-homomorphismes 
$u_{n}\colon A_{n}\rightarrow \oR$, indexé par l'ensemble $\mZ_{\geq 1}$ ordonné par la relation de divisibilité. 
On note $B_{n}$ l'image de $u_{n}$ et on pose 
\begin{equation}\label{higgs1-dlog5dd}
B_\infty=\underset{\underset{n\geq 1}{\longrightarrow}}{\lim}\ B_{n},
\end{equation}
que l'on identifie à une sous-$R$-algèbre de $\oR$. On note $H_n$ le corps des fractions de $B_n$
et $H_\infty$ le corps des fractions de $B_\infty$.

D'autre part, le morphisme \eqref{higgs1-dlog5c} induit un morphisme
\begin{equation}\label{higgs1-dlog5cd}
\Spec(\co_{\oK})\rightarrow \underset{\underset{n\geq 1}{\longleftarrow}}{\lim}\ S_{n}.
\end{equation}
On peut donc étendre les $u_n$ en un système inductif de $(R\otimes_{\co_K}\co_\oK)$-homomorphismes 
\begin{equation}
v_{n}\colon A_{n}\otimes_{\co_{K_{n}}}\co_{\oK}\rightarrow \oR,
\end{equation} 
indexé par l'ensemble $\mZ_{\geq 1}$ ordonné par la relation de divisibilité. 
On note $R_{n}$ l'image de $v_{n}$ et on pose 
\begin{equation}\label{higgs1-dlog5d}
R_\infty=\underset{\underset{n\geq 1}{\longrightarrow}}{\lim}\ R_{n},
\end{equation}
que l'on identifie à une sous-$(R\otimes_{\co_K}\co_\oK)$-algèbre de $\oR$. 
On note $F_n$ le corps des fractions de $R_n$ et $F_\infty$ le corps des fractions de $R_\infty$.
On pose 
\begin{equation}\label{higgs1-dlog5f}
R_{p^\infty}=\underset{\underset{n\geq 0}{\longrightarrow}}{\lim}\ R_{p^n},
\end{equation}
où la limite inductive est indexée par l'ensemble $\mN$ ordonné par la relation d'ordre habituelle. 
On identifie $R_{p^\infty}$ à une sous-$(R\otimes_{\co_K}\co_\oK)$-algèbre de $R_\infty$, et 
on note $F_{p^\infty}$ le corps des fractions de $R_{p^\infty}$. 

\begin{prop}\label{higgs1-gal1}
{\rm (i)}\ Pour tout $n\geq 1$, $\Spec(B_{n})$ est une composante connexe ouverte de $X_{n}$ et 
$\Spec(R_{n})$ est une composante connexe ouverte de $X_{n}\otimes_{\co_{K_{n}}}\co_{\oK}$.

{\rm (ii)}\ Les anneaux $B_n$, $R_n$ $(n\geq 1)$, $B_\infty$, $R_\infty$ et $R_{p^\infty}$ sont intègres et normaux. 

{\rm (iii)}\ Pour tout $n\geq 0$, on a 
\begin{eqnarray}
\Spec(B_{p^n})&=&X_{p^n},\\
\Spec(R_{p^n})&=&(X_{p^n}\otimes_{\co_{K_n}}\co_\oK)\times_{(X\otimes_{\co_K}\co_\oK)}\Spec(R_1).
\end{eqnarray}

{\rm (iv)}\ Les extensions  $F_n$ $(n\geq 1)$, $F_\infty$ et $F_{p^\infty}$ de $F$ sont galoisiennes et on a des homomorphismes injectifs canoniques \eqref{higgs1-dlog1ff}
\begin{eqnarray}
\Gal(F_\infty/F_1)&\rightarrow& L_\lambda\otimes \hmZ(1),\label{higgs1-gal1a}\\
\Gal(F_{p^\infty}/F_1)&\stackrel{\sim}{\rightarrow}& L_\lambda\otimes \mZ_p(1), \label{higgs1-gal1b}
\end{eqnarray}
le second étant un isomorphisme. De plus, le diagramme 
\begin{equation}\label{higgs1-gal1c}
\xymatrix{
{\Gal(F_\infty/F_1)}\ar[r]\ar[d]&{L_\lambda\otimes \hmZ(1)}\ar[d]\\
{\Gal(F_{p^\infty}/F_1)}\ar[r]&{L_\lambda\otimes \mZ_p(1)}}
\end{equation}
où les flèches verticales sont les morphismes canoniques est commutatif.
\end{prop}

D'après \ref{higgs1-dlog4}(vi), pour tout $n\geq 1$, le morphisme \eqref{higgs1-dlog5b} induit un morphisme
\begin{equation}
\tkappa\rightarrow X_n^\circ\otimes_{K_n}\oK.
\end{equation}
On désigne par $\tkappa_n$ son image, qui est un point générique de $X_n^\circ\otimes_{K_n}\oK$,
et par $\kappa_n$ l'image de $\tkappa_n$ dans $X_n^\circ$, qui est un point générique de $X_n^\circ$.  
On notera que $H_n$ est le corps résiduel de 
$X_n$ en $\kappa_n$ et que $F_n$ est le corps résiduel de $X_n\otimes_{\co_{K_n}}\co_\oK$ en $\tkappa_n$.

(i)\ Comme $X_n$ est noethérien et normal d'après \ref{higgs1-dlog4}(ii), 
$\Spec(B_n)$ est la composante connexe de $X_{n}$ contenant $\kappa_n$, qui est évidemment ouverte dans $X_n$. 
Compte tenu de \ref{higgs1-dlog4}(iii),
$\Spec(R_{n})$ est la composante connexe de $X_{n}\otimes_{\co_{K_{n}}}\co_{\oK}$ contenant  $\tkappa_n$,
qui est ouverte dans $X_{n}\otimes_{\co_{K_{n}}}\co_{\oK}$.

(ii) Ces anneaux sont clairement intègres.
Pour tout $n\geq 1$, $B_n$ et $R_n$ sont normaux en vertu de (i) et \ref{higgs1-dlog4}(ii)-(iii).  
Il en est donc de même de $B_\infty$, $R_\infty$ et $R_{p^\infty}$. 

(iii) Cela résulte de (i) et \ref{higgs1-dlog4}(v). 

(iv) Il résulte de \ref{higgs1-dlog4}(vi) que pour tout entier $n\geq 1$, $H_n$ est une extension galoisienne
de $F$, de groupe de Galois canoniquement isomorphe à un sous-groupe de $L\otimes_\mZ\mu_n(\co_{\oK})$;
plus précisément, $\Gal(H_n/F)$ est le sous-groupe de décomposition de $\kappa_n$. 
D'après \ref{higgs1-dlog4}(v), si $n$ est une puissance de $p$, on a 
\begin{equation}\label{higgs1-gal1e}
\Gal(H_n/F)\simeq L\otimes_\mZ\mu_n(\co_{\oK}).
\end{equation}
Notons $M_n$ l'image de l'homomorphisme canonique $F\otimes_KK_n\rightarrow H_n$, de sorte que $M_n$ est 
une extension galoisienne de $F$ de groupe de Galois un sous-groupe de $\Gal(K_n/K)$. On a alors un diagramme
commutatif d'extensions de corps
\begin{equation}
\xymatrix{
&H_n\ar[r]&F_n\\
F\ar[r]&M_n\ar[r]\ar[u]&F_1\ar[u]\\
K\ar[r]\ar[u]&K_n\ar[u]\ar[r]&\oK\ar[u]}
\end{equation}
où les homomorphismes $F\otimes_KK_n\rightarrow M_n$, $M_n\otimes_{K_n}\oK\rightarrow F_1$
et $H_n\otimes_{M_n}F_1\rightarrow F_n$ sont surjectifs. On en déduit que $F_1$ est une extension galoisienne de $F$,
de groupe de Galois canoniquement isomorphe à un sous-groupe de $G_K$.
Par suite, $F_n$ est une extension galoisienne de $F$, 
en tant que composée des extensions galoisiennes $H_n$ et $F_1$ de $F$. 
En particulier, $F_n$ est une extension galoisienne de $F_1$,
de groupe de Galois canoniquement isomorphe à un sous-groupe de $\Gal(H_n/M_n)$.  
Il résulte de \ref{higgs1-dlog4}(v) que si $n$ est une puissance de $p$, on a 
\begin{equation}\label{higgs1-gal1ee}
\Gal(F_n/F_1)\simeq \Gal(H_n/M_n).
\end{equation}

D'autre part, on a un diagramme commutatif
\begin{equation}\label{higgs1-gal1g}
\xymatrix{
{\Gal(H_n/F)}\ar[rr]\ar@{->>}[d]&&{L\otimes_\mZ\mu_n(\co_{\oK})}\ar[d]^{\lambda_n}\\
{\Gal(M_n/F)}\ar@{^(->}[r]&{\Gal(K_n/K)}\ar[r]^-(0.4)\sim&{\mu_n(\co_{\oK})}}
\end{equation}
où $\lambda_n$ est le morphisme défini par $\lambda\in P$ \eqref{higgs1-dlog1f} et 
les flèches non libellées sont les morphismes canoniques. 
On en déduit que $\Gal(F_n/F_1)$ 
est canoniquement isomorphe à un sous-groupe de $\ker(\lambda_n)$.
Si $n$ est une puissance de $p$, les isomorphismes \eqref{higgs1-gal1e} et \eqref{higgs1-gal1ee} et une chasse au diagramme \eqref{higgs1-gal1g}
montrent que 
\begin{equation}\label{higgs1-gal1h}
\Gal(F_n/F_1)\simeq \ker(\lambda_n).
\end{equation}
La proposition s'ensuit par passage à la limite projective sur $n$.

\begin{rema}
On notera que la condition \ref{higgs1-dlog1}(C$_6$) n'est pas utilisée dans les preuves de \ref{higgs1-dlog4} et \ref{higgs1-gal1}. 
Pour l'énoncé \ref{higgs1-dlog19}, elle ne sert que dans la preuve de (v).
\end{rema}

\subsection{}\label{higgs1-gal2}\index{10620@$\Delta$, $\Delta_\infty$, $\Delta_{p^\infty}$}
\index{10621@$\Gamma$, $\Gamma_\infty$, $\Gamma_{p^\infty}$}\index{10622@$\Sigma$, $\Sigma_0$}
On pose $\Gamma=\Gal(\oF/F)$, $\Gamma_\infty=\Gal(F_\infty/F)$, $\Gamma_{p^\infty}=\Gal(F_{p^\infty}/F)$,
$\Delta=\Gal(\oF/F_1)$, $\Delta_\infty=\Gal(F_\infty/F_1)$, $\Delta_{p^\infty}=\Gal(F_{p^\infty}/F_1)$, 
$\Sigma=\Gal(\oF/F_\infty)$ et $\Sigma_0=\Gal(F_{\infty}/F_{p^\infty})$.

\begin{equation}\label{higgs1-gal2a}
\xymatrix{
R\ar[r]\ar@/^3pc/[rrrrrrr]|{\Gamma}&R_1\ar[rr]^{\Delta_{p^\infty}}\ar@/^2pc/[rrrr]|{\Delta_\infty}\ar@/_2pc/[rrrrrr]|{\Delta}&&
{R_{p^\infty}}\ar[rr]^{\Sigma_0}&&{R_\infty}\ar[rr]^\Sigma&&\oR\\
{\co_K}\ar[r]^{G_K}\ar[u]&{\co_\oK}\ar[u]}
\end{equation}
D'après \ref{higgs1-gal1}(iv), $\Delta_\infty$ est canoniquement isomorphe à un 
sous-groupe de $L_\lambda\otimes_{\mZ}\hmZ(1)$, 
$\Delta_{p^\infty}$ est canoniquement isomorphe à $L_\lambda\otimes_{\mZ}\mZ_p(1)$,
et $\Sigma_0$ est un groupe profini d'ordre premier à $p$.

On désigne par $K^+$ l'extension de $K$ contenue dans $\oK$ telle que 
$\Gal(\oK/K^+)$ soit l'image de l'homomorphisme canonique $\Gal(F_1/F)\rightarrow G_K$.
Pour que $K^+=K$, il faut et il suffit que $X\times_S\oeta$ soit intègre, ou ce qui revient au même, connexe.   

\subsection{}\label{higgs1-gal3}
Soit $M$ un $\mZ$-$\Delta_\infty$-module discret de torsion $p$-primaire. 
Comme la $p$-dimension cohomologique de $\Sigma_0$ est nulle
(\cite{serre1} I cor.~2 de prop.~14), pour tout $q\geq 0$, le morphisme canonique 
\begin{equation}\label{higgs1-gal3a}
\rH^q(\Delta_{p^\infty},M^{\Sigma_0})\rightarrow  \rH^q(\Delta_\infty,M)
\end{equation}
est un isomorphisme. Par suite, la $p$-dimension cohomologique de $\Delta_\infty$ est égale à celle 
de $\Delta_{p^\infty}$, c'est-à-dire au rang $d$ de $L_\lambda$ \eqref{higgs1-cg05}. 

\begin{remas}\label{higgs1-higgs551}
Soient $A$ une $\mZ_p$-algèbre complète et séparée pour la topologie $p$-adique,  
$M$ un $A$-module complet et séparé pour la topologie $p$-adique. Alors~:
\begin{itemize}
\item[(i)] Les homomorphismes canoniques
\begin{equation}\label{higgs1-higgs551a}
\Hom_{\mZ_p}(\Delta_{p^\infty},M)\rightarrow 
\Hom_{\mZ}(\Delta_{p^\infty},M)\rightarrow \Hom_{\mZ}(\Delta_\infty,M)
\end{equation}
sont bijectifs. En effet, comme la multiplication par $p$ dans $\Sigma_0$ est un isomorphisme, 
pour tout homomorphisme $\psi\colon \Sigma_0\rightarrow M$,  on a $\psi(\Sigma_0)\subset \cap_{n\geq 0} p^nM=0$.  
\item[(ii)] Le morphisme canonique
\begin{equation}\label{higgs1-higgs551b}
\Hom_{\mZ}(\Delta_\infty,A)\otimes_{A}M\rightarrow \Hom_{\mZ}(\Delta_\infty,M)
\end{equation}
est bijectif. Cela résulte de (i)  car $\Delta_{p^\infty}$ est un $\mZ_p$-module libre de type fini.
\item[(iii)] Tout homomorphisme de $\Delta_\infty$ (resp. $\Delta_{p^\infty}$) dans $M$ est continu lorsque 
l'on munit $\Delta_\infty$ (resp. $\Delta_{p^\infty}$) de la topologie profinie et $M$ de la topologie $p$-adique. 
En effet, tout homomorphisme  $\psi\colon \Delta_\infty \rightarrow M$
se factorise à travers un homomorphisme $\varphi\colon  \Delta_{p^\infty} \rightarrow M$ d'après (i), 
et $\varphi$ est clairement continu pour les topologies $p$-adiques. 
En particulier, munissant $M$ de l'action triviale de $\Delta_{p^\infty}$, les morphismes canoniques 
\begin{eqnarray}
\Hom_{\mZ}(\Delta_{p^\infty},M)&\rightarrow& \rH^1_\cont(\Delta_{p^\infty},M),\\
\Hom_{\mZ}(\Delta_{\infty},M)&\rightarrow& \rH^1_\cont(\Delta_{\infty},M),
\end{eqnarray}
sont des isomorphismes. 
\end{itemize}
\end{remas}

\begin{lem}\label{higgs1-cg07}
Pour tout $a\in \co_\oK$, l'homomorphisme canonique  \eqref{higgs1-gal2}
\begin{equation}\label{higgs1-cg07a}
R_{p^\infty}/aR_{p^\infty}\rightarrow (R_\infty/aR_\infty)^{\Sigma_0}
\end{equation}
est un isomorphisme.
\end{lem}
Soient $N$ une extension galoisienne finie de $F_{p^\infty}$ contenue dans $F_\infty$, 
$A$ la clôture intégrale de $R_{p^\infty}$ dans $N$, $G=\Gal(N/F_{p^\infty})$. Comme on a $R_{p^\infty}=A\cap F_{p^\infty}$
et $A=R_\infty\cap N$, les homomorphismes canoniques $R_{p^\infty}/aR_{p^\infty}\rightarrow A/aA\rightarrow 
R_\infty/aR_\infty$ sont injectifs. D'autre part, l'ordre de $G$ étant premier à $p$ (donc inversible dans $A$), 
l'homomorphisme 
\begin{equation}
R_{p^\infty}=A^G\rightarrow (A/aA)^G
\end{equation}
est surjectif. L'assertion s'ensuit.

\begin{lem}\label{higgs1-pur8}
Les anneaux $\hRun$, $\hRpi$, $\hRi$ et $\hoR$ sont $\co_C$-plats et les homomorphismes 
canoniques $\hRun\rightarrow \hRpi$, $\hRpi\rightarrow \hRi$ et $\hRi\rightarrow \hoR$ sont injectifs. 
\end{lem}

Comme $R_1$, $R_{p^\infty}$ et $R_\infty$ sont normaux d'après \ref{higgs1-gal1}(ii), on a, pour tout $n\geq 0$,
\[
p^nR_1=(p^nR_{p^\infty})\cap R_1, \ \ \ 
p^nR_{p^\infty}=(p^nR_\infty)\cap R_{p^\infty} \ \ \ {\rm et} \ \ \ p^nR_\infty=(p^n\oR)\cap R_\infty.
\]
On en déduit que les homomorphismes 
canoniques $\hRun\rightarrow \hRpi$, $\hRpi\rightarrow \hRi$ et $\hRi\rightarrow \hoR$ sont injectifs.
D'autre part, d'après (\cite{ac} chap.~III §2.11 prop.~14 et cor.~1), pour tout $n\geq 0$,  on a
\begin{equation}\label{higgs1-pur8a}
\hRun/p^n\hRun\simeq R_1/p^nR_1.
\end{equation} 
Soient $x\in \hRun$ tel que $px=0$, $\ox$ la classe de $x$ dans $\hRun/p^n\hRun$ $(n\geq 1)$. 
Comme $R_1$ est $\co_\oK$-plat, il résulte de \eqref{higgs1-pur8a} que $\ox\in p^{n-1}\hRun/p^n\hRun$.
On en déduit que $x\in \cap_{n\geq 0}p^n\hRun=\{0\}$. 
Par suite, $p$ n'est pas diviseur de zéro dans $\hRun$, et donc $\hRun$ est $\co_C$-plat (\cite{ac} Chap.~VI §3.6 lem.~1). 
Le même argument montre que $\hRi$, $\hRpi$ et $\hoR$ sont plats sur $\co_C$.

\begin{prop}\label{higgs1-pur81}
L'anneau $\hRun$ est normal. 
\end{prop}

On note d'abord que $\hRun$ est une $\co_C$-algèbre topologiquement de présentation finie (\cite{egr1} 1.10.4),
et par suite que $\hRun[\frac 1 p]$ est une algèbre affinoïde sur $C$. 

Montrons que $\hRun[\frac 1 p]$ est normal. 
On identifie $\hRun$ au séparé complété $p$-adique de $B=R_1\otimes_{\co_\oK}\co_{C}$ 
et on note $\varphi\colon B\rightarrow \hRun$ l'homomorphisme canonique. 
Soient $\fq$ un idéal maximal de $\hRun[\frac 1 p]$, $\fp=\varphi^{-1}(\fq)$.  
En vertu de (\cite{egr1} 1.12.18), l'homomorphisme canonique $B_\fp\rightarrow (\hRun)_{\fq}$ 
induit un isomorphisme entre les séparés complétés de ces anneaux locaux pour les topologies 
définies par leurs idéaux maximaux respectifs. Le schéma $\Spec(B[\frac 1 p])$ muni de la structure logarithmique 
image inverse de $\cM_X$ est lisse sur $\Spec(C)$ muni de la structure logarithmique triviale. 
Par suite, $B[\frac 1 p]$ est normal en vertu de (\cite{kato2} 4.1 et 8.2).  
Comme $B[\frac 1 p]$ est un anneau excellent, 
on déduit de ce qui précède que les séparés complétés de $\hRun[\frac 1 p]$ 
en chacun de ses idéaux maximaux sont normaux (\cite{ega4} 7.8.3(v)).
Comme $\hRun[\frac 1 p]$ est un anneau excellent (\cite{bkkn} 3.3.3), 
ses localisés en chacun de ses idéaux maximaux sont normaux (\cite{ega4} 7.8.3(v)). 
Par suite, $\hRun[\frac 1 p]$ est normal (\cite{ega4} 7.8.3(iv)). 

Pour tout $h\in \hRun[\frac 1 p]$, posons 
\begin{equation}
|h|_{\sup}=\sup_{x\in \Max(\hRun[\frac 1 p])}|h(x)|,
\end{equation}
où $\Max(\hRun[\frac 1 p])$ est le spectre maximal de $\hRun[\frac 1 p]$, ou ce qui revient au même, 
l'ensemble des points rigides de $\Spf(\hRun)$ (\cite{egr1} 3.3.2). 
C'est une semi-norme multiplicative sur $\hRun[\frac 1 p]$ (\cite{bgr} 6.2.1/1). 
Posons 
\begin{equation}
B_{\sup}=\{h\in\hRun[\frac 1 p]\ | \ |h|_{\sup}\leq 1\}.
\end{equation}
Comme $\hRun$ est $\co_C$-plat \eqref{higgs1-pur8} et que $\hRun[\frac 1 p]$ est réduit, on a 
$\hRun\subset B_{\sup}$ et $B_{\sup}$ est la fermeture intégrale de $\hRun$ dans $\hRun[\frac 1 p]$ 
(\cite{bgr} 6.3.4/1 et 6.2.2/3). Comme $R_1\otimes_{\co_\oK}\ok$ est réduit d'après \ref{higgs1-dlog19}(iv) et \ref{higgs1-gal1}(i), 
$\hRun=B_{\sup}$ en vertu de (\cite{bgr} 6.4.3/4; cf. aussi \cite{blr4} 1.1).
Donc $\hRun$ est intégralement clos dans $\hRun[\frac 1 p]$ et est par suite normal.

\begin{teo}[Faltings, \cite{faltings2} § 2b]\label{higgs1-pur44}\index{Theoreme de presque-purete de Faltings@Théorème de presque pureté de Faltings}
Pour toute extension finie $N$ de $H_\infty$ contenue dans $\oF$, la clôture intégrale
de $B_\infty$ dans $N$ est presque-étale sur $B_\infty$ \eqref{higgs1-dlog5dd}. 
\end{teo}

Signalons ici que Scholze dispose, dans le cadre de sa théorie des perfectoïdes,  d'une généralisation de ce  résultat  
(\cite{scholze} 1.10 et 7.9).

\begin{cor}\label{higgs1-pur4}
L'extension $\oF$ de $F_\infty$ est la réunion d'un système inductif filtrant de 
sous-extensions galoisiennes finies $E$ de $F_\infty$ telles que la clôture intégrale
de $R_\infty$ dans $E$ soit presque-étale sur $R_\infty$ \eqref{higgs1-dlog5d}. 
\end{cor}

Soient $N$ une extension finie de $H_\infty$ contenue dans $\oF$, 
$E$ l'image de l'homomorphisme canonique $N\otimes_{H_\infty}F_\infty\rightarrow \oF$, 
$\cN$ (resp. $\cE$) la clôture intégrale de $R$ dans $N$ (resp. $E$). On sait \eqref{higgs1-pur44} que 
$\cN$ est presque-étale sur $B_\infty$. Donc $\cN\otimes_{B_\infty}R_\infty$ est presque-étale sur $R_\infty$
d'après (\cite{tsuji2} 7.4).  Par suite, $\cE$ est presque-étale sur $R_\infty$ en vertu de (\cite{tsuji2} 7.11 et 7.4).
Si $N/H_\infty$ est galoisienne, il en est de même de $E/F_\infty$, d'où la proposition.

\begin{cor}\label{higgs1-pur7}
Pour tout sous-anneau $A$ de $R_\infty$, 
le morphisme canonique 
\begin{equation}
\Omega^1_{R_\infty/A}\otimes_{R_\infty}\oR\rightarrow \Omega^1_{\oR/A}
\end{equation}
est un presque-isomorphisme.
\end{cor} 

Cela résulte de \ref{higgs1-pur4} et (\cite{faltings1} I 2.4(i)). 

\begin{cor}\label{higgs1-pur5}
Soit $M$ un $\oR$-module muni d'une action $\oR$-semi-linéaire continue de $\Sigma$
pour la topologie discrète de $M$. Alors $\rH^i(\Sigma,M)$ est presque nul pour tout $i\geq 1$, 
et le morphisme canonique $M^{\Sigma}\otimes_{R_\infty}\oR\rightarrow M$ est un presque-isomorphisme.
\end{cor}

Soient $N$ une extension galoisienne finie de $F_\infty$ contenue dans $\oF$, $D$ la clôture intégrale
de $R_\infty$ dans $N$, $G=\Gal(N/F_\infty)$, $\Sigma_N=\Gal(\oF/N)$. Supposons que $D$ 
soit presque-étale sur $R_\infty$. Alors $D$ est un presque $G$-torseur sur $R_\infty$ (\cite{tsuji2}  12.9). 
On en déduit par (\cite{tsuji2} 12.5 et 12.8) que, 
pour tout $i\geq 1$, $\rH^i(G,M^{\Sigma_N})$ est presque nul, et le morphisme canonique 
$M^{\Sigma}\otimes_{R_\infty}D\rightarrow M^{\Sigma_N}$ est un presque-isomorphisme.
La proposition s'ensuit par passage à la limite inductive en vertu de \ref{higgs1-pur4}.

\begin{cor}\label{higgs1-pur6}
Soit $M$ un $\oR$-module muni d'une action $\oR$-semi-linéaire continue de $\Delta$
pour la topologie discrète de $M$. Alors le morphisme canonique 
$\rH^i(\Delta_\infty,M^\Sigma)\rightarrow \rH^i(\Delta,M)$
est un presque-isomorphisme pour tout $i\geq 0$.
\end{cor}

Cela résulte de \ref{higgs1-pur5} et de la suite spectrale 
\begin{equation}
E_1^{ij}=\rH^i(\Delta_\infty,\rH^j(\Sigma,M))\Rightarrow \rH^{i+j}(\Delta,M).
\end{equation}

\begin{cor}\label{higgs1-pur55}
Soient $(M_n)_{n\in \mN}$ 
un système projectif de $\hoR$-représentations de $\Sigma$ \eqref{higgs1-not55}, $M$ sa limite projective. 
On suppose que, pour tout $n\geq 0$, $M_n$ est annulé par une puissance de $p$, que l'action de $\Sigma$ sur $M_n$ est 
continue pour la topologie discrète et que le morphisme $M_{n+1}\rightarrow M_n$ est surjectif. 
Alors $\rH^i_{\cont}(\Sigma,M)$ est presque nul pour tout entier $i\geq 1$.
\end{cor}

En effet, d'après \eqref{higgs1-limproj2c} et \eqref{higgs1-limproj2d}, on a  
\begin{equation}\label{higgs1-pur55a}
0\rightarrow \rR^1\underset{\underset{n}{\longleftarrow}}{\lim}\ \rH^{i-1}(\Sigma,M_n)\rightarrow 
\rH^i_\cont(\Sigma,M)\rightarrow \underset{\underset{n}{\longleftarrow}}{\lim}\ \rH^{i}(\Sigma,M_n)\rightarrow 0.
\end{equation}
Pour tout $q\geq 1$, $\underset{\underset{n}{\longleftarrow}}{\lim}\ \rH^{q}(\Sigma,M_n)$ et 
$\rR^1\underset{\underset{n}{\longleftarrow}}{\lim}\ \rH^{q}(\Sigma,M_n)$ 
sont presque nuls en vertu de \ref{higgs1-pur5} et (\cite{gr} 2.4.2(ii)). 
Pour tout $n\geq 0$, notons $C_n$ le noyau du morphisme surjectif $M_{n+1}\rightarrow M_n$, de sorte que 
l'on a une suite exacte 
\begin{equation}\label{higgs1-pur55c}
\xymatrix{
{M_{n+1}^\Sigma}\ar[r]^{\psi_n}&{M_n^\Sigma}\ar[r]&{\rH^1(\Sigma,C_n)}}.
\end{equation}
Alors $\coker(\psi_n)$ est presque nul d'après \ref{higgs1-pur5}. On en déduit que 
$\rR^1\underset{\underset{n}{\longleftarrow}}{\lim}\ M_n^\Sigma$ est presque nul 
en vertu de (\cite{gr} 2.4.2(iii) et 2.4.3), d'où la proposition.

\begin{cor}\label{higgs1-pur12}
Pour tout $a\in \co_\oK$, l'homomorphisme canonique
\begin{equation}
R_\infty/a R_\infty \rightarrow (\oR/a \oR)^{\Sigma}
\end{equation}
est un presque-isomorphisme.
\end{cor}

Soient $N$ une extension galoisienne finie de $F_\infty$ contenue dans $\oF$, $D$ la clôture intégrale
de $R_\infty$ dans $N$, $G=\Gal(N/F_\infty)$, 
$\Tr_G$ l'endomorphisme $R_\infty$-linéaire de $D$ (ou de $D/aD$) induit par 
$\sum_{\sigma\in G}\sigma$. Comme on a $D=\oR\cap N$ et $R_\infty=D\cap F_\infty$ d'après \ref{higgs1-gal1}(ii), 
les homomorphismes $R_\infty/aR_\infty\rightarrow D/aD\rightarrow \oR/a\oR$ sont injectifs. 
Supposons que $D$ soit presque-étale sur $R_\infty$.  
Alors $D$ est un presque $G$-torseur sur $R_\infty$ en vertu de (\cite{tsuji2} 12.9). Par suite, le quotient
\[
\frac{(D/aD)^G}{\Tr_G(D/aD)}
\]
est presque nul en vertu de (\cite{tsuji2} 12.8). Comme $\Tr_G(D)\subset R_\infty$, l'homomorphisme 
$R_\infty/aR_\infty\rightarrow (D/aD)^G$ est un presque-isomorphisme. 
La proposition s'en déduit par passage à la limite inductive d'après \ref{higgs1-pur4}. 

\begin{cor}\label{higgs1-pur122}
L'homomorphisme canonique $\hRi\rightarrow \hoR^\Sigma$ est un presque-isomorphisme.
\end{cor}
En effet, l'homomorphisme en question est la limite projective des homomorphismes ($r\geq 0$)
\begin{equation}
R_\infty/p^r R_\infty \rightarrow (\oR/p^r \oR)^{\Sigma}.
\end{equation}
Cela se vérifie aisément ou se déduit de \eqref{higgs1-limproj2d}. 
La proposition résulte donc de \ref{higgs1-pur12} et (\cite{gr} 2.4.2(ii)).

\begin{cor}\label{higgs1-pur121}
Pour tout élément non nul $a$ de $\co_\oK$ et tout entier $i\geq 0$, les morphismes canoniques 
\begin{eqnarray}
\rH^i(\Delta_{p^\infty},R_{p^\infty}/a R_{p^\infty})&\rightarrow& \rH^i(\Delta,\oR/a \oR),\label{higgs1-pur121a}\\
\rH^i(\Delta_\infty,R_\infty/a R_\infty)&\rightarrow& \rH^i(\Delta,\oR/a \oR),\label{higgs1-pur121b}
\end{eqnarray}
sont des presque-isomorphismes. 
\end{cor}

En effet, \eqref{higgs1-pur121b} est un presque-isomorphisme en vertu de \ref{higgs1-pur6} et \ref{higgs1-pur12}. 
D'autre part, le morphisme canonique 
\begin{equation}
\rH^i(\Delta_{p^\infty},R_{p^\infty}/aR_{p^\infty})\rightarrow \rH^i(\Delta_\infty,R_{\infty}/aR_\infty)
\end{equation}
est un isomorphisme d'après \ref{higgs1-cg07} et \eqref{higgs1-gal3a}.

\begin{cor}\label{higgs1-cg8}
Pour tout entier $i\geq 0$, les morphismes canoniques 
\begin{eqnarray}
\rH^i_\cont(\Delta_{p^\infty},\hRpi)&\rightarrow& \rH^i_\cont(\Delta,\hoR),\label{higgs1-cg8b}\\
\rH^i_\cont(\Delta_\infty,\hRi)&\rightarrow& \rH^i_\cont(\Delta,\hoR),\label{higgs1-cg8a}
\end{eqnarray}
sont des presque-isomorphismes. 
\end{cor}
Pour tout entier $r\geq 0$, notons 
\begin{equation}
\psi_r\colon \rH^i(\Delta_\infty,R_\infty/p^rR_\infty)\rightarrow \rH^i(\Delta,\oR/p^r\oR)
\end{equation}
l'homomorphisme canonique et $A_r$ (resp. $C_r$) son noyau (resp. conoyau).   
On sait que les $\co_\oK$-modules $A_r$ et $C_r$ sont presque nuls d'après \ref{higgs1-pur121}. Donc les $\co_\oK$-modules
\[
\underset{\underset{r\geq 0}{\longleftarrow}}{\lim}\ A_r, \ \ \ \underset{\underset{r\geq 0}{\longleftarrow}}{\lim}\ C_r, 
\ \ \ \rR^1\underset{\underset{r\geq 0}{\longleftarrow}}{\lim}\ A_r, \ \ \ 
\rR^1\underset{\underset{r\geq 0}{\longleftarrow}}{\lim}\ C_r
\]
sont presque nuls en vertu de (\cite{gr} 2.4.2(ii)). On en déduit que les morphismes 
\[
\underset{\underset{r\geq 0}{\longleftarrow}}{\lim}\ \psi_r
\ \ \ {\rm et} \ \ \ \rR^1\underset{\underset{r\geq 0}{\longleftarrow}}{\lim}\ \psi_r
\]
sont des presque-isomorphismes, et il en est alors de même de \eqref{higgs1-cg8a} 
compte tenu de \eqref{higgs1-limproj2c} et \eqref{higgs1-limproj2d}. On démontre de même que \eqref{higgs1-cg8b}
est un presque-isomorphisme.

\section{L'extension de Faltings}\label{higgs1-ext}

\subsection{}\label{higgs1-ext1} 
Soient $K_0$ le corps des fractions de $\rW(k)$ \eqref{higgs1-not33}, 
$\cD_{K/K_0}$ la différente de l'extension $K/K_0$.
D'après (\cite{fontaine2} théo.~1'), il existe un et un unique morphisme  $\co_{\oK}$-linéaire et $G_K$-équivariant
\begin{equation}\label{higgs1-ext1a} 
\phi\colon \oK\otimes_{\mZ_p}\mZ_p(1)\rightarrow \Omega^1_{\co_\oK/\co_K},
\end{equation}
tel que pour tous $\zeta\in \mZ_p(1)$, $a\in \co_{\oK}$ et $r\in \mN$, 
si $\zeta_r\in \mu_{p^r}(\co_\oK)$ est l'image canonique de $\zeta$, on ait
\begin{equation}
\phi(p^{-r}a\otimes \zeta)=a\cdot  d\log(\zeta_r).
\end{equation}
Elle est surjective de noyau $\rho^{-1}\co_\oK(1)$, 
où $\rho$ est un élément de $\co_{\oK}$ de valuation $\frac{1}{p-1}+v(\cD_{K/K_0})$.

\subsection{}\label{higgs1-log-ext6}\index{10633@$d\log(\pi_n)$}
Le morphisme \eqref{higgs1-dlog5cd} induit un $\co_K$-homomorphisme 
\begin{equation}\label{higgs1-log-ext6a}
\underset{\underset{n\geq 1}{\longrightarrow}}{\lim}\ \co_{K_n}\rightarrow \co_{\oK}, 
\end{equation}
où la limite inductive est indexée par l'ensemble 
$\mZ_{\geq 1}$ ordonné par la relation de divisibilité.
Pour tout $n\geq 1$, on identifie dans la suite $\co_{K_n}$ à une sous-$\co_K$-algèbre de $\co_{\oK}$; en particulier, 
on considère $\pi_{n}$ comme un élément de $\co_{\oK}$  \eqref{higgs1-dlog2}. 
On a $\pi_1=\pi$ et $\pi_{mn}^m=\pi_n$ pour tous $m,n\geq 1$.

Pour tout entier $n\geq 1$, il existe un élément de $\Omega^1_{\co_{\oK}/\co_K}$, que l'on note $d\log(\pi_n)$, 
tel que pour tout entier $m\geq 1$, divisible par $p$, on ait
\begin{equation}\label{higgs1-log-ext6b}
d\log(\pi_n)=\frac{m}{\pi_{mn}}d\pi_{mn} \in \Omega^1_{\co_{\oK}/\co_K}.
\end{equation}
En effet, $m\in \pi_{mn}\co_{K_{mn}}$, de sorte que $m/\pi_{mn}\in \co_{K_{mn}}$, et  
pour tout $m'\geq 1$, on a 
\begin{equation}\label{higgs1-log-ext6c}
\frac{m}{\pi_{mn}}d\pi_{mn}=
\frac{m}{\pi_{mn}}m'\pi_{mm'n}^{m'-1}d\pi_{mm'n}=\frac {mm'}{\pi_{mm'n}}d\pi_{mm'n} 
\in \Omega^1_{\co_{\oK}/\co_K}.
\end{equation}  
On prendra garde que l'élément $d\log(\pi_n)$ ne dépend pas seulement de $\pi_n$, 
mais aussi de l'homomorphisme \eqref{higgs1-log-ext6a}. 

Pour tous entiers $m,n\geq 1$, on a 
\begin{eqnarray}
\pi_nd\log(\pi_n)&=&d\pi_n,\label{higgs1-log-ext6d}\\
d\log(\pi_n)&=&md\log(\pi_{mn}).\label{higgs1-log-ext6e}
\end{eqnarray}

Comme les morphismes canoniques $\Omega^1_{\co_{K_n}/\co_K}\rightarrow \Omega^1_{\co_{\oK}/\co_K}$
sont injectifs (\cite{fontaine2} 2.4 lem.~4), pour tout entier $n\geq 1$, 
l'annulateur de $d\log(\pi_n)$ dans  $\Omega^1_{\co_{\oK}/\co_K}$ est $n\pi\co_{\oK}$ (\cite{fontaine2} 2.1 lem.~1).

\subsection{}\label{higgs1-log-ext3}\index{10640@$\mN^{(n)}$, $\mN_\infty$}\index{10641@$a_n$, $a_\infty$}
Considérons le système inductif de monoïdes $(\mN^{(n)})_{n\geq 1}$, 
indexé par l'ensemble $\mZ_{\geq 1}$ ordonné par la relation de divisibilité, 
défini par $\mN^{(n)}=\mN$ pour tout $n\geq 1$ et dont l'homomorphisme de transition
$\lambda_{n,mn}\colon \mN^{(n)}\rightarrow \mN^{(mn)}$ (pour $m, n\geq 1$) est l'homomorphisme $\varpi_m$  
de Frobenius d'ordre $m$ de $\mN$ \eqref{higgs1-not-mon1}. 
On note $\mN_\infty$ sa limite inductive~:
\begin{equation}\label{higgs1-log-ext3b}
\mN_\infty=\underset{\underset{n\geq 1}{\longrightarrow}}{\lim}\ \mN^{(n)}.
\end{equation} 
Pour tout $n\geq 1$, on note $a_n\colon \mN^{(n)}\rightarrow \co_{\oK}$ l'homomorphisme défini 
par $a_n(1)=\pi_n$. Pour tous entiers $m,n\geq 1$, le diagramme 
\begin{equation}\label{higgs1-log-ext3a}
\xymatrix{
{\mN^{(n)}}\ar[r]^{a_n}\ar[d]_{\lambda_{n,mn}}&{\co_{\oK}}\ar@{=}[d]\\
{\mN^{(mn)}}\ar[r]^-(0.5){a_{mn}}&{\co_{\oK}}}
\end{equation}
est commutatif. Les $a_{n}$ définissent donc par passage à la limite inductive un homomorphisme 
\begin{equation}\label{higgs1-log-ext3c}
a_\infty\colon \mN_\infty\rightarrow \co_{\oK}.
\end{equation}
On notera $\mN^{(1)}$ (resp. $a_1$) simplement $\mN$ (resp. $a$).

\begin{lem}\label{higgs1-log-ext7}
Soit $n$ un entier $\geq 1$. Alors~: 

{\rm (i)}\ On a des isomorphismes canoniques
\begin{eqnarray}
\Omega^1_{(\co_{\oK},\mN^{(n)})/(\co_{K},\mN)}&\simeq& \frac{\Omega^1_{\co_{\oK}/\co_K}\oplus \co_{\oK}/n\co_\oK}
{(d\pi_n-\pi_n)\co_{\oK}},\label{higgs1-log-ext7a}\\
\Omega^1_{(\co_{\oK},\mN^{(n)})/(\co_{\oK},\mN)}&\simeq& \co_{\oK}/(n\co_{\oK}+\pi_n\co_{\oK}).\label{higgs1-log-ext7b}
\end{eqnarray}

{\rm (ii)}\ Si $p$ divise $n$, le noyau du morphisme canonique 
\begin{equation}\label{higgs1-log-ext7c}
\Omega^1_{\co_{\oK}/\co_K}\rightarrow \Omega^1_{(\co_{\oK},\mN^{(n)})/(\co_{K},\mN)}
\end{equation}
est engendré par $d\log(\pi)$ \eqref{higgs1-log-ext6b}. 
\end{lem}

(i) On a un diagramme commutatif à carré cartésien
\begin{equation}\label{higgs1-log-ext7d}
\xymatrix{
{\Spec(\co_{\oK})}\ar[r]^-(0.5){j_n}\ar[rd]&{\Spec(\co_{\oK}[\xi]/(\xi^n-\pi)}\ar[r]\ar[d]\ar@{}[rd]|\Box&
{\Spec(\co_{\oK}[\mN^{(n)}])}\ar[d]\\
&{\Spec(\co_{K})}\ar[r]&{\Spec(\co_{K}[\mN])}}
\end{equation}
où $j_n$ est l'immersion fermée définie par l'équation $\xi-\pi_n$. D'autre part, on a un isomorphisme canonique 
\begin{equation}
\Omega^1_{(\co_{\oK}[\mN^{(n)}],\mN^{(n)})/(\co_{K}[\mN],\mN)}\simeq 
\Omega^1_{\co_{\oK}/\co_K}\otimes_{\co_\oK}\co_{\oK}[\mN^{(n)}]\oplus \co_{\oK}[\mN^{(n)}]/n\co_{\oK}[\mN^{(n)}].
\end{equation}
L'isomorphisme \eqref{higgs1-log-ext7a} s'en déduit aussitôt. 
L'isomorphisme \eqref{higgs1-log-ext7b} se démontre par un diagramme analogue à \eqref{higgs1-log-ext7d}. 
 
(ii)  Soit $\omega\in \Omega^1_{\co_{\oK}/\co_K}$ tel que son image dans $\Omega^1_{(\co_{\oK},\mN^{(n)})/(\co_{K},\mN)}$
soit nulle. D'après \eqref{higgs1-log-ext7a}, il existe $x\in \co_{\oK}$ tel que $\omega=x(d\pi_n-\pi_n)\in 
\Omega^1_{\co_{\oK}/\co_K}\oplus \co_{\oK}/n\co_K$. Par suite, $x\pi_n\in n\co_{\oK}$ et $\omega=xd\pi_n\in \co_{\oK}
d\log(\pi)$. Inversement, on a 
\begin{equation}
d\log(\pi)=(n/\pi_n)(d\pi_n-\pi_n)\in \Omega^1_{\co_{\oK}/\co_K}\oplus \co_{\oK}/n\co_K.
\end{equation} 
Donc l'image de $d\log(\pi)$ par le morphisme \eqref{higgs1-log-ext7c} est nulle en vertu de \eqref{higgs1-log-ext7a}.

\begin{prop}\label{higgs1-log-ext8}
Le morphisme canonique 
\begin{equation}\label{higgs1-log-ext8a}
\Omega^1_{\co_{\oK}/\co_K}\rightarrow \Omega^1_{(\co_{\oK},\mN_\infty)/(\co_{K},\mN)}
\end{equation}
est surjectif et son noyau est engendré par $d\log(\pi)$. En particulier, le morphisme $\phi$ 
\eqref{higgs1-ext1a} induit un morphisme $\co_{\oK}$-linéaire surjectif  
\begin{equation}\label{higgs1-log-ext8b}
\oK\otimes_{\mZ_p}\mZ_p(1)\rightarrow \Omega^1_{(\co_{\oK},\mN_\infty)/(\co_{K},\mN)},
\end{equation}
de noyau $(\pi\rho)^{-1}\co_{\oK}(1)$.
\end{prop}

On notera d'abord que la seconde proposition est une conséquence immédiate de la première et de \ref{higgs1-ext1}. 
Montrons la première proposition.
Par la propriété universelle des modules de $1$-différentielles logarithmiques, le morphisme canonique 
\begin{equation}\label{higgs1-log-ext8c}
\underset{\underset{n\geq 1}{\longrightarrow}}{\lim}\ \Omega^1_{(\co_{\oK},\mN^{(n)})/(\co_{K},\mN)}
\rightarrow \Omega^1_{(\co_{\oK},\mN_\infty)/(\co_{K},\mN)}
\end{equation}
où la limite inductive est indexée par l'ensemble 
$\mZ_{\geq 1}$ ordonné par la relation de divisibilité, est un isomorphisme. 
Il résulte donc de \ref{higgs1-log-ext7}(ii) que le noyau du morphisme \eqref{higgs1-log-ext8a} est engendré par $d\log(\pi)$. 
D'autre part, on a un isomorphisme canonique 
\begin{equation}\label{higgs1-log-ext8d}
\underset{\underset{n\geq 1}{\longrightarrow}}{\lim}\ \Omega^1_{(\co_{\oK},\mN^{(n)})/(\co_{\oK},\mN)}
\stackrel{\sim}{\rightarrow} \Omega^1_{(\co_{\oK},\mN_\infty)/(\co_{\oK},\mN)}.
\end{equation}
D'après \eqref{higgs1-log-ext7b}, pour tout entier $n\geq 1$, on a
\begin{equation}\label{higgs1-log-ext8e}
\Omega^1_{(\co_{\oK},\mN^{(n)})/(\co_{\oK},\mN)}=\left\{
\begin{array}{clcr}
0&{\rm si} \ (n,p)=1,\\
k_n\otimes_{\co_{K_n}}\co_{\oK}&{\rm si}\ p|n,
\end{array}
\right.
\end{equation}
où $k_n$ est le corps résiduel de $\co_{K_n}$. On notera que pour tous entiers $m,n\geq 1$, 
le morphisme canonique 
\begin{equation}\label{higgs1-log-ext8f}
\Omega^1_{(\co_{\oK},\mN^{(n)})/(\co_{\oK},\mN)}\rightarrow \Omega^1_{(\co_{\oK},\mN^{(nm)})/(\co_{\oK},\mN)}
\end{equation}
s'identifie à $m$ fois le morphisme canonique 
$k_n\otimes_{\co_{K_n}}\co_{\oK} \rightarrow k_{mn}\otimes_{\co_{K_{mn}}}\co_{\oK}$. 
On en déduit que $\Omega^1_{(\co_{\oK},\mN_\infty)/(\co_{\oK},\mN)}=0$ et 
par suite que le morphisme \eqref{higgs1-log-ext8a} est surjectif.

\begin{rema}\label{higgs1-log-ext85}
Pour tout entier $n\geq 1$, l'image canonique 
de l'élément $d\log(\pi_n)\in \Omega^1_{\co_\oK/\co_K}$ \eqref{higgs1-log-ext6b}
dans  $\Omega^1_{(\co_\oK,\mN_\infty)/(\co_K,\mN)}$ 
est égale à l'élément $d\log(1^{(n)})\in \Omega^1_{(\co_\oK,\mN^{(n)})/(\co_K,\mN)}$, où $1^{(n)}$
désigne l'image de $1$ par l'isomorphisme canonique $\mN\stackrel{\sim}{\rightarrow}\mN^{(n)}$. 
En effet, on a dans $\Omega^1_{(\co_\oK,\mN^{(pn)})/(\co_K,\mN)}$
\begin{equation}
d\log(\pi_{n})-d\log(1^{(n)})
=\frac{p}{\pi_{pn}}(d\pi_{pn}-\pi_{pn}d\log(1^{(pn)}))=0.
\end{equation}
On notera que l'égalité $d\log(\pi_n)=d\log(1^{(n)})$ vaut dans  
$\Omega^1_{(\co_\oK,\mN^{(mn)})/(\co_K,\mN)}$ pour tout entier $m\geq 1$ divisible par $p$
(mais en général elle ne vaut pas lorsque $m=1$).
\end{rema}

\begin{lem}\label{higgs1-log-ext84}
Pour toute $\co_\oK$-algèbre plate $A$, les $\co_{\oK}$-modules  
\begin{equation}
\Omega^1_{\co_\oK/\co_K}\otimes_{\co_{\oK}}A \ \ \ et \ \ \  
\Omega^1_{(\co_\oK,\mN_\infty)/(\co_K,\mN)}\otimes_{\co_{\oK}}A
\end{equation} 
n'ont pas de $\fm_\oK$-torsion non nulle. 
\end{lem}
Compte tenu de \eqref{higgs1-ext1a} et \eqref{higgs1-log-ext8b},
il suffit de montrer que le $\co_\oK$-module $A[\frac 1 p]/A$ n'a pas de $\fm_\oK$-torsion non nulle. 
La $\fm_\oK$-torsion étant contenue dans la $p$-torsion, 
il suffit encore de montrer que le $\co_\oK$-module $A/pA$, n'a pas de $\fm_\oK$-torsion non nulle. 
Un calcul de valuations montre que $\co_\oK/p\co_\oK$ n'a pas de $\fm_\oK$-torsion non nulle, autrement dit, le morphisme 
\begin{equation}
\co_\oK/p\co_\oK\rightarrow \oplus_{n\geq 1}\co_\oK/p\co_\oK, \ \ \ x\mapsto (p^{1/n}x)_{n\geq 1}
\end{equation}
est injectif. Il en est alors de même du morphisme obtenu par extension des scalaires de $\co_\oK$ à $A$,
d'où l'assertion.

\subsection{}\label{higgs1-log-ext2}\index{10645@$P^{(n)}$, $P_\infty$, $t^{(n)}\in P^{(n)}$}
\index{10646@$\alpha_{n}$, $\alpha_\infty$}
Considérons le système inductif de monoïdes $(P^{(n)})_{n\geq 1}$, 
indexé par l'ensemble $\mZ_{\geq 1}$ ordonné par la relation de divisibilité, 
défini par $P^{(n)}=P$ pour tout $n\geq 1$ et dont l'homomorphisme de transition
$i_{n,mn}\colon P^{(n)}\rightarrow P^{(mn)}$ (pour $m, n\geq 1$) est l'homomorphisme $\varpi_m$  
de Frobenius d'ordre $m$ de $P$ \eqref{higgs1-not-mon1}. 
On note $P_\infty$ sa limite inductive,
\begin{equation}\label{higgs1-log-ext2b}
P_\infty=\underset{\underset{n\geq 1}{\longrightarrow}}{\lim}\ P^{(n)}.
\end{equation} 
Pour tout $n\geq 1$, on note
\begin{equation}\label{higgs1-ext-log2d}
P\stackrel{\sim}{\rightarrow}P^{(n)}, \ \ \ t\mapsto t^{(n)},
\end{equation} 
l'isomorphisme canonique. Pour tous $m, n\geq 1$ et $t\in P$, on a  
$i_{n,mn}(t^{(n)})=(t^{(mn)})^m$ et par suite
\begin{equation}\label{higgs1-ext-log2e}
t^{(n)}=(t^{(mn)})^m \in P_\infty. 
\end{equation}

Pour tout $n\geq 1$, on désigne par $\alpha_{n}\colon P^{(n)}\rightarrow R_{n}$ l'homomorphisme induit par 
le morphisme strict canonique $(X_{n},\cM_{X_{n}})\rightarrow \bA_P$ \eqref{higgs1-dlog3a}. 
Le lecteur se souviendra ici que $R_n$ dépend du choix du morphisme \eqref{higgs1-dlog5b}.
Pour tous entiers $m,n\geq 1$, le diagramme 
\begin{equation}\label{higgs1-log-ext2a}
\xymatrix{
{P^{(n)}}\ar[r]^{\alpha_{n}}\ar[d]_{i_{n,mn}}&R_{n}\ar[d]\\
{P^{(mn)}}\ar[r]^-(0.5){\alpha_{mn}}&R_{mn}}
\end{equation}
est commutatif. Les $\alpha_{n}$ définissent donc par passage à la limite inductive un homomorphisme 
\begin{equation}\label{higgs1-log-ext2c}
\alpha_\infty\colon P_\infty\rightarrow R_\infty.
\end{equation}
On note encore $\alpha_\infty\colon P_\infty\rightarrow \oR$ le composé de $\alpha_\infty$ et de l'injection canonique 
$R_\infty\rightarrow \oR$. On notera $P^{(1)}$ simplement $P$. On observera que $\alpha_1$ se factorise à travers $\alpha$ \eqref{higgs1-dlog1}.

\begin{prop}\label{higgs1-log-ext4}
{\rm (i)}\ La suite canonique 
\begin{equation}\label{higgs1-log-ext4a}
0\rightarrow \Omega^1_{(R,P)/(\co_K,\mN)}\otimes_{R}R_\infty\rightarrow 
\Omega^1_{(R_\infty,P_\infty)/(\co_{\oK},\mN)}\rightarrow 
\Omega^1_{(R_\infty,P_\infty)/(R_1,P)}\rightarrow 0
\end{equation}
est exacte.

{\rm (ii)}\ Pour tout entier $m\geq 0$, le morphisme 
\begin{equation}\label{higgs1-log-ext4b}
\Hom(p^{-m}\mZ/\mZ,\Omega^1_{(R_\infty,P_\infty)/(R_1,P)})\rightarrow
\Omega^1_{(R,P)/(\co_K,\mN)}\otimes_{R}(R_\infty/p^m R_\infty)
\end{equation}
déduit de \eqref{higgs1-log-ext4a} par le diagramme du serpent est un isomorphisme. 

{\rm (iii)}\ Il existe un isomorphisme $R_\infty$-linéaire canonique
\begin{equation}\label{higgs1-log-ext4bb}
(P^\gp/\lambda\mZ)\otimes_{\mZ} (R_\infty[\frac 1 p]/R_\infty)
\stackrel{\sim}{\rightarrow} \Omega^1_{(R_\infty,P_\infty)/(R_1,P)}. 
\end{equation}
\end{prop}

(i) D'après \ref{higgs1-gal1}(i), 
on a un isomorphisme canonique 
\begin{equation}\label{higgs1-log-ext4cc}
\Omega^1_{(R_1,P)/(\co_\oK,\mN)} \stackrel{\sim}{\rightarrow} \Omega^1_{(R,P)/(\co_K,\mN)}\otimes_RR_1.
\end{equation}
D'autre part, on a un isomorphisme canonique
\begin{equation}\label{higgs1-log-ext4c}
\underset{\underset{n\geq 1}{\longrightarrow}}{\lim}\ \Omega^1_{(R_{n},P^{(n)})/(\co_{\oK},\mN)}
\stackrel{\sim}{\rightarrow} \Omega^1_{(R_\infty,P_\infty)/(\co_{\oK},\mN)},
\end{equation}
où la limite inductive est indexée par l'ensemble 
$\mZ_{\geq 1}$ ordonné par la relation de divisibilité. 
Il suffit donc de montrer que pour tout $n\geq 1$, le morphisme canonique 
\begin{equation}\label{higgs1-log-ext4d}
\Omega^1_{(R_1,P)/(\co_{\oK},\mN)}\otimes_{R_1}R_{n}\rightarrow \Omega^1_{(R_{n},P^{(n)})/(\co_{\oK},\mN)}
\end{equation}
est injectif. Rappelons que 
$\Spec(R_{n})$ est une composante connexe ouverte de $X_{n}\otimes_{\co_{K_{n}}}\co_{\oK}$ \eqref{higgs1-gal1},
que le morphisme canonique $X\rightarrow \Spec(\co_K[P]/(\pi-e^\lambda))$ est étale \eqref{higgs1-dlog1b}
et qu'on a un diagramme cartésien de $\co_{\oK}$-morphismes \eqref{higgs1-dlog4e}
\begin{equation}\label{higgs1-log-ext4ee}
\xymatrix{
{X_n\otimes_{\co_{K_n}}\co_{\oK}}\ar[r]\ar[d]&{\Spec(\co_{\oK}[P^{(n)}]/(\pi_n-e^{\lambda^{(n)}}))}\ar[d]\\
{X\otimes_{\co_K}\co_{\oK}}\ar[r]&{\Spec(\co_{\oK}[P]/(\pi-e^\lambda))}}
\end{equation}
où $\lambda^{(n)}$ est l'image de $\lambda$ dans $P^{(n)}$ par l'isomorphisme \eqref{higgs1-ext-log2d}.
Considérons le diagramme commutatif à carrés cartésiens
\begin{equation}\label{higgs1-log-ext4e}
\xymatrix{
{\Spec(\co_{\oK}[P^{(n)}]/(\pi_{n}-e^{\lambda^{(n)}}))}\ar[r]\ar[rd]&
{\Spec(\co_{\oK}[P^{(n)}]/(\pi-e^{n\lambda^{(n)}}))}\ar[r]\ar[d]\ar@{}[rd]|{\Box}&{\Spec(\co_{\oK}[P^{(n)}])}\ar[d]^{i_n^a}\\
&{\Spec(\co_{\oK}[P]/(\pi-e^\lambda))}\ar[r]\ar[d]\ar@{}[rd]|{\Box}&{\Spec(\co_{\oK}[P])}\ar[d]\\
&{\Spec(\co_{\oK})}\ar[r]&{\Spec(\co_{\oK}[\mN])}}
\end{equation}
où $i_n^a$ est le morphisme induit par l'homomorphisme canonique $i_n\colon P\rightarrow P^{(n)}$. 
On a un isomorphisme canonique 
\begin{equation}\label{higgs1-log-ext4g}
\Omega^1_{(\co_{\oK}[P],P)/(\co_{\oK}[\mN],\mN)}\stackrel{\sim}{\rightarrow} (P^\gp/\lambda\mZ)\otimes_{\mZ}\co_{\oK}[P],
\end{equation}
tel que pour tout $x\in P$, l'image de $d\log(x)$ soit la classe de $x$ dans $P^\gp/\lambda\mZ$. Par suite, le morphisme 
\[
\Omega^1_{(\co_{\oK}[P],P)/(\co_{\oK}[\mN],\mN)}\otimes_{\co_{\oK}[P]}\co_{\oK}[P^{(n)}]\rightarrow 
\Omega^1_{(\co_{\oK}[P^{(n)}],P^{(n)})/(\co_{\oK}[\mN],\mN)}
\] 
induit par $i_n^a$ s'identifie au morphisme $\co_{\oK}[P^{(n)}]$-linéaire 
\begin{equation}\label{higgs1-log-ext4h}
(P^\gp/\lambda\mZ)\otimes_{\mZ}\co_{\oK}[P^{(n)}]\rightarrow (P^\gp/n\lambda\mZ)\otimes_{\mZ}\co_{\oK}[P^{(n)}]
\end{equation}
déduit de la multiplication par $n$ sur $P^\gp$. 

Il résulte du diagramme \eqref{higgs1-log-ext4e} et de ce qui le précède qu'on a un isomorphisme canonique 
\begin{equation}\label{higgs1-log-ext4k}
\frac{(P^\gp/n\lambda\mZ)\otimes_{\mZ}R_{n}}{\olambda\otimes \pi_{n}  R_{n}}\stackrel{\sim}{\rightarrow} 
\Omega^1_{(R_{n},P^{(n)})/(\co_{\oK},\mN)},
\end{equation}
où $\olambda$ désigne la classe de $\lambda$ dans $P^\gp/n\lambda\mZ$.
On en déduit un morphisme surjectif canonique 
\begin{equation}\label{higgs1-log-ext4m}
\Omega^1_{(R_{n},P^{(n)})/(\co_{\oK},\mN)}\rightarrow (P^\gp/\lambda\mZ)\otimes_{\mZ}R_{n}.
\end{equation}
D'après \eqref{higgs1-log-ext4h}, le composé des morphismes \eqref{higgs1-log-ext4d} et \eqref{higgs1-log-ext4m}  
s'identifie au morphisme de multiplication par $n$ sur $(P^\gp/\lambda\mZ)\otimes_{\mZ}R_{n}$, 
qui est injectif puisque le sous-groupe de torsion de $P^\gp/\lambda\mZ$ est d'ordre premier à $p$
et que $R_{n}$ est plat sur $\mZ_p$. Par suite, le morphisme \eqref{higgs1-log-ext4d}
est injectif.

(ii) Il suffit de montrer que la multiplication par $p^m$ dans $\Omega^1_{(R_\infty,P_\infty)/(\co_{\oK},\mN)}$
est un isomorphisme. Pour tous entiers $n,n'\geq 1$, on a un diagramme commutatif 
\begin{equation}\label{higgs1-log-ext4n}
\xymatrix{
{(P^\gp/n\lambda\mZ)\otimes_{\mZ}R_{n}}\ar[r]\ar@{->>}[d]&{(P^\gp/nn'\lambda\mZ)\otimes_{\mZ}R_{nn'}}\ar@{->>}[d]\\
{\Omega^1_{(R_{n},P^{(n)})/(\co_{\oK},\mN)}}\ar[r]^-(0.5){u_{n,nn'}}&{\Omega^1_{(R_{nn'},P^{(nn')})/(\co_{\oK},\mN)}}}
\end{equation}
où les flèches verticales sont les morphismes surjectifs déduits de \eqref{higgs1-log-ext4k}, $u_{n,nn'}$
est le morphisme canonique et la flèche horizontale supérieure est induite par la multiplication par $n'$ dans $P^\gp$
et l'homomorphisme canonique $R_{n}\rightarrow R_{nn'}$. On en déduit que 
\begin{equation}\label{higgs1-log-ext4p}
u_{n,nn'}(\Omega^1_{(R_{n},P^{(n)})/(\co_{\oK},\mN)})\subset 
n'\cdot \Omega^1_{(R_{nn'},P^{(nn')})/(\co_{\oK},\mN)}.
\end{equation}
Par suite, la multiplication par $p^m$ dans $\Omega^1_{(R_\infty,P_\infty)/(\co_{\oK},\mN)}$
est surjective. 

Soit $\omega\in \Omega^1_{(R_\infty,P_\infty)/(\co_{\oK},\mN)}$ tel que $p^m\omega=0$. Il existe alors $n\geq 1$
tel que $\omega\in \Omega^1_{(R_{n},P^{(n)})/(\co_{\oK},\mN)}$ et que $p^m\omega=0$ dans 
$\Omega^1_{(R_{n},P^{(n)})/(\co_{\oK},\mN)}$. Considérons la suite exacte canonique 
\begin{equation}
\Omega^1_{(\co_{\oK},\mN^{(n)})/(\co_{\oK},\mN)}\otimes_{\co_{\oK}}R_{n}\stackrel{h_n}{\longrightarrow} 
\Omega^1_{(R_{n},P^{(n)})/(\co_{\oK},\mN)}\longrightarrow \Omega^1_{(R_{n},P^{(n)})/(\co_{\oK},\mN^{(n)})}
\longrightarrow 0.
\end{equation}
Le $R_{n}$-module $\Omega^1_{(R_{n},P^{(n)})/(\co_{\oK},\mN^{(n)})}$ est libre de type fini en vertu de \ref{higgs1-dlog4}(i) et 
\ref{higgs1-gal1}(i). Il est donc $\mZ_p$-plat. Par suite, $\omega$ est dans l'image de $h_n$.
Il résulte alors de \eqref{higgs1-log-ext8e} et \eqref{higgs1-log-ext8f} qu'il existe $n'\geq 1$
tel que l'image de $\omega$ dans $\Omega^1_{(R_{nn'},P^{(nn')})/(\co_{\oK},\mN)}$ soit nulle. Donc $\omega$
est nul dans $\Omega^1_{(R_\infty,P_\infty)/(\co_{\oK},\mN)}$.

(iii) On a un isomorphisme canonique 
\begin{equation}\label{higgs1-log-ext4aa}
\underset{\underset{n\geq 1}{\longrightarrow}}{\lim}\ \Omega^1_{(R_{n},P^{(n)})/(R_1,P)}
\stackrel{\sim}{\rightarrow} \Omega^1_{(R_\infty,P_\infty)/(R_1,P)},
\end{equation}
où la limite inductive est indexée par l'ensemble $\mZ_{\geq 1}$ ordonné par la relation de divisibilité. 
Pour tout entier $n\geq 1$, il résulte du diagramme \eqref{higgs1-log-ext4e} et de ce qui le précède
qu'on a un isomorphisme canonique 
\begin{equation}\label{higgs1-log-ext4ab}
M_n=\frac{(P^\gp/nP^\gp)\otimes_{\mZ}R_{n}}{\tlambda\otimes \pi_{n}  R_{n}}\stackrel{\sim}{\rightarrow} 
\Omega^1_{(R_{n},P^{(n)})/(R_1,P)},
\end{equation}
où $\tlambda$ est la classe de $\lambda$ dans $P^\gp/nP^\gp$. 
On note $d\log(\lambda^{(n)})$ l'image de $\tlambda\otimes 1$ dans $M_n$, ce qui se justifie par l'isomorphisme 
\eqref{higgs1-log-ext4ab}. Pour tout entier $m\geq 1$, la multiplication par $m$ dans $P^\gp$ et l'homomorphisme
canonique $R_n\rightarrow R_{mn}$ induisent un morphisme $M_n\rightarrow M_{mn}$. Les $(M_n)_{n\geq 1}$ 
forment un système inductif pour la relation de divisibilité, et on a un isomorphisme canonique
\begin{equation}\label{higgs1-log-ext4ba}
\underset{\underset{n\geq 1}{\longrightarrow}}{\lim}\ M_n
\stackrel{\sim}{\rightarrow} \Omega^1_{(R_\infty,P_\infty)/(R_1,P)}.
\end{equation}
 
Soit $n$ un entier $\geq 1$. On a dans $M_{pn}$ 
\begin{equation}\label{higgs1-log-ext4ac}
d\log(\lambda^{(n)})=pd\log(\lambda^{(pn)})=\frac{p}{\pi_{pn}}\pi_{pn}d\log(\lambda^{(pn)})=0.
\end{equation}
Par suite, pour tout entier $m\geq 1$ divisible par $p$, le composé des morphismes canoniques
\begin{equation}\label{higgs1-log-ext4ad}
(P^\gp/nP^\gp)\otimes_{\mZ}R_{n} \rightarrow M_n\rightarrow M_{mn}
\end{equation}
se factorise à travers un morphisme $R_n$-linéaire
\begin{equation}\label{higgs1-log-ext4ae}
N_n=((P^\gp/\lambda\mZ)/n (P^\gp/\lambda\mZ))\otimes_{\mZ}R_{n} 
\rightarrow M_{mn}.
\end{equation}
D'autre part, on a un morphisme surjectif canonique
\begin{equation}\label{higgs1-log-ext4af}
M_{mn}\rightarrow 
((P^\gp/\lambda\mZ)/mn (P^\gp/\lambda\mZ))\otimes_{\mZ}R_{mn}.
\end{equation}
Le composé de \eqref{higgs1-log-ext4ae} et \eqref{higgs1-log-ext4af} est induit par la multiplication par $m$ dans 
$P^\gp/\lambda\mZ$ et par l'homomorphisme canonique $R_n\rightarrow R_{mn}$. 
Comme le sous-groupe de torsion de $P^\gp/\lambda\mZ$ est d'ordre premier à $p$
et que $R_{n}$ est $\mZ_p$-plat, la multiplication par $m$ induit un morphisme 
injectif 
\[
((P^\gp/\lambda\mZ)/n (P^\gp/\lambda\mZ))\otimes_{\mZ}R_{n} \rightarrow 
((P^\gp/\lambda\mZ)/mn (P^\gp/\lambda\mZ))\otimes_{\mZ}R_{n}. 
\]
Comme $R_n$ est intègre et normal d'après \ref{higgs1-gal1}(ii), 
l'homomorphisme canonique 
\[
R_n/{mn}R_n\rightarrow R_{mn}/mnR_{mn}
\]
est injectif. On en déduit que le morphisme \eqref{higgs1-log-ext4ae} est injectif. 

Pour tous entiers $m,n\geq 1$, la multiplication par $m$ dans $P^\gp/\lambda\mZ$ et l'homomorphisme
canonique $R_n\rightarrow R_{mn}$ induisent un morphisme $N_n\rightarrow N_{mn}$. Les $(N_n)_{n\geq 1}$ 
forment un système inductif pour la relation de divisibilité. 
Par passage à la limite inductive sur les morphismes \eqref{higgs1-log-ext4ae}, d'abord relativement aux 
entiers $m\geq 1$ multiples de $p$, puis aux entiers $n\geq 1$,  on obtient un isomorphisme
\begin{equation}
\underset{\underset{n\geq 1}{\longrightarrow}}{\lim}\ N_n 
\stackrel{\sim}{\rightarrow} \Omega^1_{(R_\infty,P_\infty)/(R_1,P)}.
\end{equation}
Il est clair que le morphisme canonique 
\begin{equation}
\underset{\underset{s\geq 0}{\longrightarrow}}{\lim}\ N_{p^s} 
\rightarrow \underset{\underset{n\geq 1}{\longrightarrow}}{\lim}\ N_n, 
\end{equation}
où la première limite inductive est indexée par l'ensemble $\mN$ ordonné par la relation d'ordre habituelle,
est un isomorphisme. La proposition s'ensuit.

\begin{cor}\label{higgs1-log-ext44}
Les $\co_\oK$-modules $\Omega^1_{(R_\infty,P_\infty)/(R_1,P)}$ et 
$\Omega^1_{(R_\infty,P_\infty)/(R_1,P)}\otimes_{R_\infty}\oR$ n'ont pas de $\fm_\oK$-torsion non nulle.
\end{cor}
Il suffit de calquer la preuve de \ref{higgs1-log-ext84} en tenant compte de \eqref{higgs1-log-ext4bb} 
et du fait que $R_\infty$ et $\oR$ sont $\co_\oK$-plats (\cite{ac} chap.~VI §3.6 lem.~1).

\begin{cor}\label{higgs1-log-ext43}
Tout élément de $\Omega^1_{(R_\infty,P_\infty)/(R,P)}$ (resp. $\Omega^1_{(\oR,P_\infty)/(R,P)}$) 
est annulé par une puissance de $p$. 
\end{cor}

Considérons la suite exacte
\begin{equation}\label{higgs1-log-ext43b}
\Omega^1_{R_1/R}\otimes_{R_1}R_\infty\rightarrow 
\Omega^1_{(R_\infty,P_\infty)/(R,P)}\rightarrow \Omega^1_{(R_\infty,P_\infty)/(R_1,P)}\rightarrow 0.
\end{equation}
D'après \ref{higgs1-gal1}(i), on a un isomorphisme canonique
\begin{equation}
\Omega^1_{R_1/R}\stackrel{\sim}{\rightarrow} \Omega^1_{\co_{\oK}/\co_K}\otimes_{\co_{\oK}}R_1.
\end{equation}
D'autre part, tout élément de $\Omega^1_{(R_\infty,P_\infty)/(R_1,P)}$ est annulé par une puissance de $p$ 
d'après \eqref{higgs1-log-ext4bb}, et il en est de même de $\Omega^1_{\co_\oK/\co_K}$  \eqref{higgs1-ext1}.
Par suite, tout élément de $\Omega^1_{(R_\infty,P_\infty)/(R,P)}$  
est annulé par une puissance de $p$. 
Comme le morphisme canonique 
\begin{equation}\label{higgs1-log-ext43a}
\Omega^1_{(R_\infty,P_\infty)/(R,P)}\otimes_{R_\infty}\oR\rightarrow \Omega^1_{(\oR,P_\infty)/(R,P)}
\end{equation}
est un presque-isomorphisme en vertu de \ref{higgs1-pur7} et (\cite{kato1} 1.7), on en déduit que
tout élément de $\Omega^1_{(\oR,P_\infty)/(R,P)}$ est annulé par une puissance de $p$.

\begin{rema}\label{higgs1-log-ext45}
L'isomorphisme \eqref{higgs1-log-ext4b} peut se déduire aussi de l'isomorphisme \eqref{higgs1-log-ext4bb} comme suit. 
On notera d'abord qu'on a un isomorphisme canonique 
\begin{equation}\label{higgs1-log-ext45b}
(P^\gp/\lambda\mZ)\otimes_{\mZ}R \stackrel{\sim}{\rightarrow} \Omega^1_{(R,P)/(\co_K,\mN)}.
\end{equation}
Pour tout entier $m\geq 0$, la $p^m$-torsion de $R_\infty[\frac 1 p]/R_\infty$ est canoniquement isomorphe à 
$R_\infty/p^mR_\infty$. Comme  le sous-groupe de torsion de $P^\gp/\lambda\mZ$ est d'ordre premier à $p$, 
\eqref{higgs1-log-ext4bb} induit un isomorphisme 
\begin{equation}\label{higgs1-log-ext45a}
(P^\gp/\lambda\mZ)\otimes_{\mZ}(R_\infty/p^m R_\infty)\stackrel{\sim}{\rightarrow}
\Hom(p^{-m}\mZ/\mZ,\Omega^1_{(R_\infty,P_\infty)/(R_1,P)}).
\end{equation}
Il ressort de la preuve de \ref{higgs1-log-ext4} que celui-ci est l'inverse de l'isomorphisme \eqref{higgs1-log-ext4b}.
\end{rema}

\begin{prop}\label{higgs1-log-ext5}
{\rm (i)}\ Le noyau du morphisme canonique 
\begin{equation}\label{higgs1-log-ext5a}
\Omega^1_{\co_{\oK}/\co_K}\otimes_{\co_{\oK}}R_\infty\rightarrow \Omega^1_{(R_\infty,P_\infty)/(R,P)}
\end{equation}
est engendré par $d\log(\pi)$ \eqref{higgs1-log-ext6b}. 

{\rm (ii)}\ La suite de morphismes canoniques 
\begin{equation}\label{higgs1-log-ext5aa}
0\rightarrow \Omega^1_{(\co_{\oK},\mN_\infty)/(\co_K,\mN)}\otimes_{\co_{\oK}}R_\infty\rightarrow 
\Omega^1_{(R_\infty,P_\infty)/(R,P)}\rightarrow \Omega^1_{(R_\infty,P_\infty)/(R_1,P)}\rightarrow 0
\end{equation}
est exacte et scindée. 
\end{prop}

(i) On a un isomorphisme canonique 
\begin{equation}\label{higgs1-log-ext5b}
\underset{\underset{n\geq 1}{\longrightarrow}}{\lim}\ \Omega^1_{(R_{n},P^{(n)})/(R,P)}
\stackrel{\sim}{\rightarrow} \Omega^1_{(R_\infty,P_\infty)/(R,P)},
\end{equation}
où la limite inductive est indexée par l'ensemble 
$\mZ_{\geq 1}$ ordonné par la relation de divisibilité.  
Il suffit donc de montrer que pour tout entier $n\geq 1$ tel que $p$ divise $n$, le noyau du morphisme canonique 
\begin{equation}\label{higgs1-log-ext5c}
\Omega^1_{\co_{\oK}/\co_K}\otimes_{\co_{\oK}}R_{n}\rightarrow \Omega^1_{(R_{n},P^{(n)})/(R,P)}
\end{equation}
est engendré par $d\log(\pi)$. Par fonctorialité, ce morphisme se factorise à travers le 
morphisme canonique 
\begin{equation}\label{higgs1-log-ext5cc}
\Omega^1_{(\co_{\oK},\mN^{(n)})/(\co_K,\mN)}\otimes_{\co_{\oK}}R_{n}\rightarrow \Omega^1_{(R_{n},P^{(n)})/(R,P)}.
\end{equation}
Il résulte donc de \ref{higgs1-log-ext7}(ii) que $d\log(\pi)$ appartient au noyau du morphisme \eqref{higgs1-log-ext5c}.

Montrons inversement que le noyau du morphisme \eqref{higgs1-log-ext5c} est contenu dans $R_nd\log(\pi)$. 
Rappelons que $\Spec(R_{n})$ est une composante connexe ouverte de $X_{n}\otimes_{\co_{K_{n}}}\co_{\oK}$ \eqref{higgs1-gal1}
et qu'on a un diagramme cartésien de $\co_K$-morphismes \eqref{higgs1-dlog4e}
\begin{equation}\label{higgs1-log-ext5dd}
\xymatrix{
{X_n\otimes_{\co_{K_n}}\co_{\oK}}\ar[r]\ar[d]&{\Spec(\co_{\oK}[P^{(n)}]/(\pi_n-e^{\lambda^{(n)}}))}\ar[d]\\
{X}\ar[r]&{\Spec(\co_K[P]/(\pi-e^\lambda))}}
\end{equation}
où $\lambda^{(n)}$ est l'image de $\lambda$ dans $P^{(n)}$ par l'isomorphisme \eqref{higgs1-ext-log2d}.
Considérons le diagramme commutatif à carré cartésien
\begin{equation}\label{higgs1-log-ext5d}
\xymatrix{
{\Spec(\co_{\oK}[P^{(n)}]/(\pi_{n}-e^{\lambda^{(n)}}))}\ar[r]\ar[rd]&
{\Spec(\co_{\oK}[P^{(n)}]/(\pi-e^{n\lambda^{(n)}}))}\ar[r]\ar[d]\ar@{}[rd]|{\Box}&{\Spec(\co_{\oK}[P^{(n)}])}\ar[d]\\
&{\Spec(\co_{K}[P]/(\pi-e^\lambda)}\ar[r]&{\Spec(\co_{K}[P])}}
\end{equation}
On a un isomorphisme canonique
\begin{equation}\label{higgs1-log-ext5e}
\Omega^1_{(\co_{\oK}[P^{(n)}],P^{(n)})/(\co_{K}[P],P)}\stackrel{\sim}{\rightarrow} 
(P^\gp/nP^\gp)\otimes_{\mZ}\co_{\oK}[P^{(n)}] \oplus 
\Omega^1_{\co_{\oK}/\co_K}\otimes_{\co_{\oK}}\co_{\oK}[P^{(n)}]. 
\end{equation}
On en déduit un isomorphisme 
\begin{equation}\label{higgs1-log-ext5f}
\Omega^1_{(R_{n},P^{(n)})/(R,P)}\stackrel{\sim}{\rightarrow} \frac{(P^\gp/nP^\gp)\otimes_{\mZ}R_{n}\oplus 
\Omega^1_{\co_{\oK}/\co_K}\otimes_{\co_{\oK}}R_{n}}{(\olambda\otimes \pi_{n} -d\pi_{n})R_{n}},
\end{equation}
où $\olambda$ est la classe de $\lambda$ dans $P^\gp/nP^\gp$. On notera pour la preuve de (ii) 
que pour tout entier $m\geq 1$, le morphisme canonique 
\begin{equation}
\Omega^1_{(R_{n},P^{(n)})/(R,P)}\rightarrow \Omega^1_{(R_{mn},P^{(mn)})/(R,P)}
\end{equation}
est induit par le morphisme de multiplication par $m$ dans $P^\gp$, par l'identité de $\Omega^1_{\co_{\oK}/\co_K}$
et par l'homomorphisme canonique $R_n\rightarrow R_{mn}$.

Soit $\omega\in \Omega^1_{\co_{\oK}/\co_K}\otimes_{\co_{\oK}}R_n$ tel que son image par le morphisme \eqref{higgs1-log-ext5c}
soit nulle. D'après \eqref{higgs1-log-ext5f}, il existe $x\in R_{n}$ tel que 
\begin{equation}\label{higgs1-log-ext5g}
\omega=x(\olambda\otimes \pi_{n} -d\pi_{n})\in (P^\gp/nP^\gp)\otimes_{\mZ}R_{n}\oplus 
\Omega^1_{\co_{\oK}/\co_K}\otimes_{\co_{\oK}}R_{n}.
\end{equation}
Par suite, $\olambda\otimes \pi_{n} x=0$ dans $(P^\gp/nP^\gp)\otimes_{\mZ}R_{n}$.  
Comme le sous-groupe de torsion de $P^\gp/\lambda\mZ$ est d'ordre premier à $p$
et que $R_{n}$ est plat sur $\mZ_p$, l'homomorphisme
\begin{equation}\label{higgs1-log-ext5h}
(\lambda\mZ/n\lambda\mZ)\otimes_{\mZ}R_{n}\rightarrow (P^\gp/nP^\gp)\otimes_{\mZ}R_{n}
\end{equation}
est injectif. On en déduit que $\pi_{n} x\in nR_{n}$. Donc $x\in (n/\pi_{n})R_{n}$ car $p$ divise $n$ et
$R_{n}$ est plat sur $\co_{\oK}$ en vertu de \ref{higgs1-dlog4}(ii) et \ref{higgs1-gal1}(i). 
Il résulte alors de \eqref{higgs1-log-ext5g} que $\omega=-xd\pi_n\in R_n d\log(\pi)$.

(ii)  D'après \ref{higgs1-gal1}(i), on a un isomorphisme canonique
\begin{equation}
\Omega^1_{R_1/R}\stackrel{\sim}{\rightarrow} \Omega^1_{\co_{\oK}/\co_K}\otimes_{\co_{\oK}}R_1.
\end{equation}
Donc la suite \eqref{higgs1-log-ext5aa} est exacte en vertu de (i), \ref{higgs1-log-ext8} et de la suite exacte canonique 
\begin{equation}\label{higgs1-log-ext5ab}
\Omega^1_{R_1/R}\otimes_{R_1}R_\infty\rightarrow \Omega^1_{(R_\infty,P_\infty)/(R,P)}\rightarrow 
\Omega^1_{(R_\infty,P_\infty)/(R_1,P)}\rightarrow 0.
\end{equation}
Il reste à construire un scindage de \eqref{higgs1-log-ext5aa}. 
On a un isomorphisme canonique 
\begin{equation}\label{higgs1-log-ext5ac}
\underset{\underset{n\geq 1}{\longrightarrow}}{\lim}\ \Omega^1_{(R_{n},P^{(n)})/(R_1,P)}
\stackrel{\sim}{\rightarrow} \Omega^1_{(R_\infty,P_\infty)/(R_1,P)},
\end{equation}
où la limite inductive est indexée par l'ensemble 
$\mZ_{\geq 1}$ ordonné par la relation de divisibilité.
Donc compte tenu de l'isomorphisme \eqref{higgs1-log-ext5b}, il suffit de construire pour tout entier $n\geq 1$
tel que $p$ divise $n$, un inverse à droite du morphisme canonique
\begin{equation}\label{higgs1-log-ext5ad}
\Omega^1_{(R_{n},P^{(n)})/(R,P)}\rightarrow \Omega^1_{(R_{n},P^{(n)})/(R_1,P)}, 
\end{equation}
tel que la famille de ces inverses à droite soit un morphisme de systèmes inductifs.

Il résulte du diagramme \eqref{higgs1-log-ext4e} et de ce qui le précède qu'on a un isomorphisme 
\begin{equation}\label{higgs1-log-ext5ae}
\Omega^1_{(R_{n},P^{(n)})/(R_1,P)}\stackrel{\sim}{\rightarrow} 
\frac{(P^\gp/nP^\gp)\otimes_{\mZ}R_{n}}{\olambda\otimes \pi_{n} R_{n}}, 
\end{equation}
où $\olambda$ est la classe de $\lambda$ dans $P^\gp/nP^\gp$. 
Pour tout entier $m\geq 1$, le morphisme canonique 
\begin{equation}
\Omega^1_{(R_{n},P^{(n)})/(R,P)}\rightarrow \Omega^1_{(R_{mn},P^{(mn)})/(R,P)}
\end{equation}
est induit par le morphisme de multiplication par $m$ dans $P^\gp$ et par l'homomorphisme canonique 
$R_n\rightarrow R_{mn}$.
Compte tenu de \eqref{higgs1-log-ext5f} et \eqref{higgs1-log-ext5ae}, le morphisme \eqref{higgs1-log-ext5ad} est induit par la projection canonique 
\begin{equation}\label{higgs1-log-ext5af}
(P^\gp/nP^\gp)\otimes_{\mZ}R_{n}\oplus 
\Omega^1_{\co_{\oK}/\co_K}\otimes_{\co_{\oK}}R_{n}\rightarrow (P^\gp/nP^\gp)\otimes_{\mZ}R_{n}.
\end{equation}
Comme le sous-groupe de torsion de $P^\gp/\lambda\mZ$ est d'ordre premier à $p$
et que $R_{n}$ est plat sur $\mZ_p$, le morphisme $R_n$-linéaire
\begin{equation}\label{higgs1-log-ext5ag}
R_n\rightarrow P^\gp\otimes_{\mZ}R_n
\end{equation}
défini par $\lambda$, admet un inverse à gauche $R_n$-linéaire $u\colon P^\gp\otimes_{\mZ}R_n\rightarrow R_n$. 
Considérons le morphisme 
\begin{equation}
v_n\colon P^\gp\otimes_{\mZ}R_{n}\rightarrow (P^\gp/nP^\gp)\otimes_{\mZ}R_{n}\oplus 
\Omega^1_{\co_{\oK}/\co_K}\otimes_{\co_{\oK}}R_{n}
\end{equation}
défini pour $x\in P^\gp\otimes_{\mZ}R_{n}$, de classe $\ox$ dans $(P^\gp/nP^\gp)\otimes_{\mZ}R_{n}$, par 
\begin{equation}
v_n(x)=\ox-d\log(\pi_n)\otimes u(x),
\end{equation}
où $d\log(\pi_n)$ est l'élément de $\Omega^1_{\co_\oK/\co_K}$ défini dans \eqref{higgs1-log-ext6b}. 
Compte tenu de \eqref{higgs1-log-ext5f}, $v_n$ induit un morphisme $R_n$-linéaire 
\begin{equation}
w_n\colon P^\gp\otimes_{\mZ}R_{n}\rightarrow  \Omega^1_{(R_{n},P^{(n)})/(R,P)}.
\end{equation}
Pour tout $t\in P^\gp$, on a $w_n(nt)=-nu(t)d\log(\pi_n)=-u(t)d\log(\pi)$ \eqref{higgs1-log-ext6e}. 
Comme $p$ divise $n$, l'image de $d\log(\pi)$ 
dans $\Omega^1_{(R_{n},P^{(n)})/(R,P)}$ est nulle d'après la preuve de (i); donc $w_n(nt)=0$. 
D'autre part, en vertu de \eqref{higgs1-log-ext6d}, on a 
\begin{equation}
w_n(\lambda\otimes \pi_n)= \olambda\otimes\pi_n-\pi_nd\log(\pi_n)=\olambda\otimes\pi_n-d\pi_n=0.
\end{equation} 
Par suite, compte tenu de \eqref{higgs1-log-ext5ae}, $w_n$ induit un morphisme $R_n$-linéaire 
\begin{equation}
\omega_n\colon \Omega^1_{(R_{n},P^{(n)})/(R_1,P)} \rightarrow  \Omega^1_{(R_{n},P^{(n)})/(R,P)},
\end{equation}
inverse à droite du morphisme canonique \eqref{higgs1-log-ext5ad}. Comme $md\log(\pi_{mn})=d\log(\pi_n)$ pour tout $m\geq 1$,
les $\omega_n$ forment un morphisme de systèmes inductifs pour la relation de divisibilité, ce qui achève la preuve.

\subsection{}\label{higgs1-ext-log10}\index{10660@$\langle \ ,\ \rangle\colon  \Gamma_\infty\times P_\infty \rightarrow \mu_\infty(\co_{\oK})$}
Il existe une unique application 
\begin{equation}\label{higgs1-ext-log10a}
\langle \ ,\ \rangle\colon  \Gamma_\infty\times P_\infty \rightarrow 
\mu_\infty(\co_{\oK})=\underset{\underset{n\geq 1}{\longrightarrow}}{\lim}\ \mu_n(\co_{\oK}), 
\end{equation}
où la limite inductive est indexée par l'ensemble $\mZ_{\geq 1}$ ordonné par la relation de divisibilité,
telle que pour tout $g\in \Gamma_\infty$ et tout $x\in P^{(n)}$ $(n\geq 1)$, 
on ait $\langle g,x\rangle\in \mu_n(\co_{\oK})$ et 
\begin{equation}\label{higgs1-ext-log10b}
g(\alpha_\infty(x))=\langle g,x\rangle\cdot  \alpha_\infty(x),
\end{equation}
où $\alpha_\infty$ est l'homomorphisme \eqref{higgs1-log-ext2c}. 
En effet, $R_\infty$ est intègre, et on a $\alpha_\infty(x)^n\in \alpha(P)\subset R$; 
donc $\alpha_\infty(x)\not=0$ \eqref{higgs1-log-ext2} et $\alpha_\infty(x)^n$ est invariant par $\Gamma_\infty$. 

Pour tout $g\in \Gamma_\infty$, l'application $x\mapsto \langle g,x \rangle$ est un morphisme de monoïdes 
de $P_\infty$ dans $\mu_\infty(\co_{\oK})$. Pour tout $x\in P_\infty$, l'application $g\mapsto \langle g,x \rangle$ est
un $1$-cocycle, autrement dit, pour tous $g,g'\in \Gamma_\infty$, on a 
\begin{equation}\label{higgs1-ext-log10c}
\langle gg',x\rangle=g(\langle g',x\rangle)\cdot \langle g,x\rangle.
\end{equation}

Soit $n$ un entier $\geq 1$. Rappelons qu'on a un isomorphisme canonique $P^{(n)}\stackrel{\sim}{\rightarrow} P$ 
\eqref{higgs1-ext-log2d} et que $\Delta_\infty$ est canoniquement isomorphe 
à un sous-groupe de $L_\lambda\otimes_\mZ\hmZ(1)$ \eqref{higgs1-gal1a}. 
On a donc un homomorphisme canonique $\Delta_\infty\rightarrow L\otimes_\mZ\mu_n(\co_{\oK})$.
D'après \ref{higgs1-dlog4}(vi), le diagramme 
\begin{equation}\label{higgs1-ext-log10d}
\xymatrix{
{\Delta_\infty\times P^{(n)}}\ar[rr]^-(0.5){\langle \ ,\ \rangle}\ar[d]&&{\mu_n(\co_{\oK})}\\
{(L\otimes_\mZ\mu_n(\co_{\oK}))\times P}\ar[rru]&&}
\end{equation}
où la flèche oblique est induite par l'accouplement canonique $L\otimes_\mZ P^\gp\rightarrow \mZ$, 
est commutatif.

\subsection{}\label{higgs1-ext-log101}
Il existe une unique application 
\begin{equation}\label{higgs1-ext-log101a}
\langle \ ,\ \rangle\colon  G_K\times \mN_\infty \rightarrow \mu_\infty(\co_{\oK}), 
\end{equation}
telle que pour tout $g\in G_K$ et tout $x\in \mN^{(n)}$ $(n\geq 1)$, 
on ait $\langle g,x\rangle\in \mu_n(\co_{\oK})$ et 
\begin{equation}\label{higgs1-ext-log101b}
g(a_\infty(x))=\langle g,x\rangle\cdot  a_\infty(x),
\end{equation}
où $a_\infty$ est l'homomorphisme \eqref{higgs1-log-ext3c}. Identifiant $\mN^{(n)}$ et $\mN$ par l'isomorphisme 
canonique $\mN^{(n)}\stackrel{\sim}{\rightarrow}\mN$, 
la relation \eqref{higgs1-ext-log101b} devient $g(\pi_n^x)=\langle g,x\rangle\cdot  \pi_n^x$. 
L'application \eqref{higgs1-ext-log101a} vérifie des propriétés analogues à celles de \eqref{higgs1-ext-log10a}.
D'ailleurs, le diagramme
\begin{equation}
\xymatrix{
{\Gamma_\infty\times \mN_\infty}\ar[r]\ar[d]&{\Gamma_\infty\times P_\infty}\ar[d]^{\langle \ ,\ \rangle}\\
{G_K\times\mN_\infty}\ar[r]^{\langle \ ,\ \rangle}&{\mu_\infty(\co_{\oK})}}
\end{equation}
où les flèches non-libellées sont les homomorphismes canoniques est commutatif, ce qui 
justifie la même notation.

\subsection{}\label{higgs1-ext-log11}\index{10664@$\cM_{\infty}$, $\cL$, $\cN_{\infty}$, $\cN$}
Soit $\cM_{\infty}$ la structure logarithmique sur $\Spec(R_\infty)$ associée à la structure pré-logarithmique
définie par $(P_\infty,\alpha_\infty)$ \eqref{higgs1-log-ext2}. Pour tout $g\in \Gamma_\infty$, notons $\tau_g$ l'automorphisme de 
$\Spec(R_\infty)$ induit par $g$. La structure logarithmique 
$\tau_g^*(\cM_\infty)$ sur $\Spec(R_\infty)$ est associée à la structure pré-logarithmique définie par 
$(P_\infty,g\circ \alpha_\infty)$. L'homomorphisme  
\begin{equation}\label{higgs1-ext-log11a}
P_\infty\rightarrow \Gamma(\Spec(R_\infty),\cM_\infty), \ \ \ x\mapsto \langle g,x\rangle \cdot x
\end{equation}
induit un morphisme de structures logarithmiques sur $\Spec(R_\infty)$
\begin{equation}\label{higgs1-ext-log11b}
a_g\colon \tau_g^*(\cM_\infty)\rightarrow \cM_\infty.
\end{equation}
De même, en vertu de \eqref{higgs1-ext-log10c}, l'homomorphisme  
\begin{equation}\label{higgs1-ext-log11c}
P_\infty\rightarrow \Gamma(\Spec(R_\infty),\tau_g^*(\cM_\infty)), \ \ \ x\mapsto g(\langle g^{-1},x\rangle) \cdot x
\end{equation}
induit un morphisme de structures logarithmiques sur $\Spec(R_\infty)$
\begin{equation}\label{higgs1-ext-log11d}
b_g\colon \cM_\infty\rightarrow \tau_g^*(\cM_\infty).
\end{equation}
On voit aussitôt que $a_g$ et $b_g$ sont des isomorphismes inverses l'un de l'autre \eqref{higgs1-ext-log10c}, et que 
l'application $g\mapsto  (\tau_{g^{-1}},a_{g^{-1}})$ est une action à gauche de $\Gamma_\infty$ sur 
$(\Spec(R_\infty),\cM_\infty)$. 

On note $\cL$ la structure logarithmique sur $\Spec(\oR)$ image inverse de $\cM_\infty$.
Alors l'action précédente se relève en une action à gauche de $\Gamma$ sur $(\Spec(\oR),\cL)$. 

Soit $\cN_\infty$ la structure logarithmique sur $\Spec(\co_{\oK})$ associée à 
la structure pré-logarithmique définie par $(\mN_\infty,a_\infty)$ \eqref{higgs1-log-ext3c}. 
De même, l'application \eqref{higgs1-ext-log101a} définit 
une action à gauche de $G_K$ sur le schéma logarithmique $(\Spec(\co_{\oK}),\cN_\infty)$. 

Soit $\cM$ la structure logarithmique sur $\Spec(R_1)$ associée à la structure pré-logarithmique
définie par $(P,\alpha)$ \eqref{higgs1-log-ext2}. Pour tout $g\in \Gal(F_1/F)$, notons $u_g$ l'automorphisme de 
$\Spec(R_1)$ induit par $g$. 
La structure logarithmique $u_g^*(\cM)$ sur $\Spec(R_1)$ est associée à la structure pré-logarithmique définie par 
$(P,g\circ \alpha)$. Comme $g\circ \alpha=\alpha$, on en déduit un isomorphisme canonique 
\begin{equation}\label{higgs1-ext-log11e}
c_g\colon u_g^*(\cM)\stackrel{\sim}{\rightarrow} \cM.
\end{equation}
L'application $g\mapsto  (u_{g^{-1}},c_{g^{-1}})$ est une action à gauche de $\Gal(F_1/F)$ sur 
le schéma logarithmique $(\Spec(R_1),\cM)$.  

Soit $\cN$ la structure logarithmique sur $\Spec(\co_{\oK})$ associée à la structure pré-logarithmique
définie par $(\mN,a)$ \eqref{higgs1-log-ext3}. On définit de même une action à gauche de $G_K$ sur le schéma logarithmique 
$(\Spec(\co_{\oK}),\cN)$. 

On notera que tous les morphismes du diagramme commutatif
\begin{equation}\label{higgs1-ext-log11f}
\xymatrix{
{(\Spec(\oR),\cL)}\ar[r]&{(\Spec(R_\infty),\cM_\infty)}\ar[r]\ar[d]&{(\Spec(\co_{\oK}),\cN_\infty)}\ar[d]\\
&{(\Spec(R_1),\cM)}\ar[r]&{(\Spec(\co_{\oK}),\cN)}}
\end{equation}
sont $\Gamma$-équivariants.

\subsection{}\label{higgs1-ext-log12}\index{10668@$\cE_\infty$}\index{10669@$\tOmega^1_{R/\co_K}$}\index{10670@$\delta$}
On note $\cE_\infty$ la $\hRi$-représentation de $\Gamma_\infty$ définie par 
\begin{equation}\label{higgs1-ext-log12a}
\cE_\infty=\Hom(\mQ_p/\mZ_p,\Omega^1_{(R_\infty,P_\infty)/(R,P)})\otimes_{\mZ_p}\mZ_p(-1),
\end{equation}
où l'action de $\Gamma_\infty$ provient de son action sur  
$\Omega^1_{(R_\infty,P_\infty)/(R,P)}$ \eqref{higgs1-ext-log11}. 
Compte tenu de \ref{higgs1-log-ext8} et \ref{higgs1-log-ext4}(ii),
appliquant le foncteur $\Hom(\mQ_p/\mZ_p,-)\otimes_{\mZ_p}\mZ_p(-1)$ à la suite exacte scindée 
\eqref{higgs1-log-ext5aa}, on obtient une suite exacte 
de $\hRi$-représentations de $\Gamma_\infty$
\begin{equation}\label{higgs1-ext-log12b}
0\rightarrow (\pi\rho)^{-1}\hRi\rightarrow \cE_\infty \rightarrow \Omega^1_{(R,P)/(\co_K,\mN)}\otimes_{R}\hRi(-1)\rightarrow 0,
\end{equation}
où l'on a écrit $(\pi\rho)^{-1}\hRi$ au lieu de $(\pi\rho)^{-1}\co_C\otimes_{\co_C}\hRi$, 
ce qui est justifié puisque $\hRi$ est $\co_C$-plat \eqref{higgs1-pur8}. 
On notera que $\Omega^1_{(R,P)/(\co_K,\mN)}$ est un $R$-module libre de type fini. 
Pour alléger les notations, on pose 
\begin{eqnarray}
\tOmega^1_{R/\co_K}&=&\Omega^1_{(R,P)/(\co_K,\mN)}, \label{higgs1-ext-log12c}\\
\tOmega^i_{R/\co_K}&=&\wedge^i\tOmega^1_{R/\co_K}, \ \ \ \ (i\geq 1).
\end{eqnarray}
La suite \eqref{higgs1-ext-log12b}, comme la suite \eqref{higgs1-log-ext5aa}, 
sont scindées en tant que suites de $\hRi$-modules (sans actions de $\Gamma_\infty$).

On désigne par 
\begin{equation}\label{higgs1-ext-log12d}
\delta\colon \tOmega^1_{R/\co_K}\otimes_R \hRun\rightarrow \rH^1_\cont(\Delta_\infty, (\pi\rho)^{-1}  \hRi)(1)
\end{equation}
le bord de la suite exacte longue de cohomologie obtenue 
en appliquant  le foncteur $\Gamma(\Delta_\infty,-)(1)$ à la suite exacte \eqref{higgs1-ext-log12b}.
On montrera ultérieurement que $(\hRi)^{\Delta_\infty}=\hRun$ \eqref{higgs1-cg31}.

\subsection{}\label{higgs1-ext-log13}
Pour tout $\zeta\in \mZ_p(1)$, on désigne par $d\log(\zeta)$ l'élément de $\cE_\infty(1)$ défini par 
\begin{equation}\label{higgs1-ext-log13a}
d\log(\zeta)(p^{-n})=d\log(\zeta_n),
\end{equation}
où $\zeta_n$ est l'image canonique de $\zeta$ dans $\mu_{p^n}(\co_\oK)$. 
Il est clair que $d\log(\zeta)$ est l'image de $1\otimes \zeta$ par l'injection 
$(\pi\rho)^{-1}\hRi(1)\rightarrow \cE_\infty(1)$ \eqref{higgs1-ext-log12b}. 

Pour tout $t\in P$, on désigne par $d\log(\tlt)$ l'élément de $\cE_\infty(1)$ défini par 
\begin{equation}\label{higgs1-ext-log13c}
d\log(\tlt)(p^{-n})=d\log(t^{(p^n)}), 
\end{equation}
où $t^{(p^n)}$ est l'image de $t$ dans $P^{(p^n)}$ par l'isomorphisme \eqref{higgs1-ext-log2d}. 
Cet élément est bien défini en vertu de \eqref{higgs1-ext-log2e}. 
Il est clair que l'image de $d\log(\tlt)$ dans $\tOmega^1_{R/\co_K}\otimes_R \hRi$ est $d\log(t)$. 

D'après \eqref{higgs1-log-ext85}, pour tout $n\geq 0$, $d\log(\tlambda)(p^{-n})$
est l'image canonique de l'élément $d\log(\pi_{p^n})\in \Omega^1_{\co_\oK/\co_K}$ \eqref{higgs1-log-ext6b} dans 
$\tOmega^1_{(R_\infty,P_\infty)/(R,P)}$. En particulier, $d\log(\tlambda)\in (\pi\rho)^{-1}\hRi(1)\subset \cE_\infty(1)$
\eqref{higgs1-ext-log12b}. 

L'application $P\rightarrow\cE_\infty(1)$ définie par $t\mapsto d\log(\tlt)$ est un homomorphisme~; 
elle induit donc un homomorphisme que l'on note aussi
\begin{equation}\label{higgs1-ext-log13cd}
P^\gp\rightarrow \cE_\infty(1), \ \ \ t\mapsto d\log(\tlt).
\end{equation} 
Celui-ci s'insère dans un digramme commutatif 
\begin{equation}\label{higgs1-ext-log13d}
\xymatrix{
0\ar[r]&{\mZ \lambda}\ar[r]\ar[d]&{P^\gp}\ar[r]\ar[d]&{P^\gp/\mZ\lambda}\ar[d]\ar[r]&0\\
0\ar[r]&{(\pi\rho)^{-1}\hRi(1)}\ar[r]&{\cE_\infty(1)}\ar[r]&{\tOmega^1_{R/\co_K}\otimes_R\hRi}\ar[r]&0}
\end{equation}
où la flèche verticale de droite provient de l'isomorphisme canonique \eqref{higgs1-log-ext45b}. 
Il est utile de noter aussi 
\begin{equation}\label{higgs1-ext-log13e}
d\log\colon P^\gp\rightarrow \tOmega^1_{R/\co_K}
\end{equation}
l'homomorphisme induit par la dérivation logarithmique $d\log \colon P\rightarrow \tOmega^1_{R/\co_K}$.

\subsection{}\label{higgs1-ext-log14}\index{10671@$\tdelta$}
Soit $t\in P$. On désigne par  
\begin{equation}\label{higgs1-ext-log14a}
\tchi_t\colon \Gamma_\infty\rightarrow \mZ_p(1) 
\end{equation}
l'application qui à tout $g\in \Gamma_\infty$ associe l'élément 
\begin{equation}\label{higgs1-ext-log14b}
\tchi_t(g)=\underset{\underset{n\geq 0}{\longleftarrow}}{\lim}\ 
\langle g,t^{(p^n)}\rangle,
\end{equation}
où $t^{(p^n)}$ est l'image de $t$ dans $P^{(p^n)}$ par l'isomorphisme \eqref{higgs1-ext-log2d}
et $\langle g,t^{(p^n)}\rangle\in \mu_{p^n}(\co_\oK)$ est défini dans \eqref{higgs1-ext-log10a}. 
D'après \eqref{higgs1-ext-log10c}, pour tous $g,g'\in \Gamma_\infty$, on a 
\begin{equation}\label{higgs1-ext-log14c}
\tchi_t(gg')=g(\tchi_t(g'))\tchi_t(g).
\end{equation}
Donc la restriction de $\tchi_t$ à $\Delta_\infty$ est un caractère
à valeurs dans $\mZ_p(1)$; on le note $\chi_t\colon \Delta_\infty\rightarrow \mZ_p(1)$.  

On a clairement $\tchi_0=1$, et pour tous $t,t'\in P$,  
\begin{equation}\label{higgs1-ext-log14d}
\tchi_{tt'}=\tchi_t\cdot\tchi_{t'}.
\end{equation}
Par suite, l'application $P\rightarrow \Hom(\Delta_\infty,\mZ_p(1))$ définie par  $t\mapsto \chi_t$
est un homomorphisme. Elle induit donc un homomorphisme que l'on note encore 
\begin{equation}\label{higgs1-ext-log14e}
P^\gp\rightarrow \Hom(\Delta_\infty,\mZ_p(1)), \ \ \ t\mapsto \chi_t.
\end{equation}
Comme $\chi_\lambda=1$, on en déduit un homomorphisme 
\begin{equation}\label{higgs1-ext-log14f}
P^\gp/\mZ\lambda\rightarrow \Hom(\Delta_\infty,\mZ_p(1)).
\end{equation}
D'après \eqref{higgs1-ext-log10d}, ce dernier est égal au composé
\begin{equation}\label{higgs1-ext-log14g}
P^\gp/\mZ\lambda\rightarrow \Hom(L_\lambda\otimes_\mZ\mZ_p(1),\mZ_p(1))\rightarrow 
\Hom(\Delta_\infty,\mZ_p(1)),
\end{equation}
où la première flèche est induite par le morphisme canonique (de bidualité) \eqref{higgs1-dlog1ff} et la seconde flèche par
le morphisme canonique $\nu\colon \Delta_\infty \rightarrow L_\lambda\otimes_\mZ\mZ_p(1)$ \eqref{higgs1-gal1a}. 
Rappelons que $\nu$ induit un isomorphisme 
$\Delta_{p^\infty}\stackrel{\sim}{\rightarrow} L_\lambda\otimes_\mZ\mZ_p(1)$  et qu'on a un isomorphisme canonique 
\eqref{higgs1-higgs551a}
\[
\Hom_{\mZ_p}(\Delta_{p^\infty},\mZ_p(1))\stackrel{\sim}{\rightarrow} \Hom_\mZ(\Delta_\infty,\mZ_p(1)).
\]
Comme le sous-groupe de torsion de $P^\gp/\mZ\lambda$ est d'ordre premier à $p$, 
l'homomorphisme \eqref{higgs1-ext-log14f} induit donc un isomorphisme 
\begin{equation}\label{higgs1-ext-log14gh}
(P^\gp/\mZ \lambda)\otimes_\mZ\mZ_p\stackrel{\sim}{\rightarrow} \Hom(\Delta_\infty,\mZ_p(1)).
\end{equation}
On en déduit, compte tenu de \eqref{higgs1-log-ext45b} et \eqref{higgs1-higgs551b}, un isomorphisme $\hRun$-linéaire 
\begin{equation}\label{higgs1-ext-log14h}
\tdelta\colon \tOmega^1_{R/\co_K}\otimes_R\hRun\stackrel{\sim}{\rightarrow} \Hom(\Delta_\infty,\hRun(1)).
\end{equation}

\subsection{}\label{higgs1-ext-log15}
Il résulte aussitôt des définitions que, pour tous $t\in P$ et $g\in \Gamma_\infty$, on a 
\begin{equation}\label{higgs1-ext-log15a}
g(d\log(\tlt))=d\log(\tlt)+d\log(\chi_t(g)).
\end{equation}
Comme les deux membres de l'équation sont des homomorphismes de $P$ dans $\cE(1)$, 
l'égalité vaut pour tout $t\in P^\gp$. Par suite, le diagramme 
\begin{equation}\label{higgs1-ext-log15b} 
\xymatrix{
{\tOmega^1_{R/\co_K}\otimes_{R}\hRun}\ar[r]^-(0.5){\tdelta}\ar[d]_{\delta}
&{\Hom(\Delta_\infty,\hRun(1))}\ar@{=}[d]\\
{\rH^1_\cont(\Delta_\infty, (\pi\rho)^{-1}  \hRi(1))}&{\rH^1_\cont(\Delta_\infty,\hRun(1))}\ar[l]}
\end{equation}
où $\delta$ est le morphisme \eqref{higgs1-ext-log12d}, $\tdelta$ est le morphisme \eqref{higgs1-ext-log14h} et 
la flèche horizontale inférieure est induite par l'injection canonique 
$\hRun\rightarrow (\pi\rho)^{-1}  \hRi$, est commutatif. En effet, comme $\tOmega^1_{R/\co_K}$ est engendré 
sur $R$ par les éléments de la forme $d\log(t)$ pour $t\in P$, il suffit de montrer la commutativité de ce diagramme
pour ces éléments, ce qui résulte de \eqref{higgs1-ext-log15a}.

\subsection{}\label{higgs1-ext-log16}
Considérons le diagramme commutatif de morphismes canoniques
\begin{equation}\label{higgs1-ext-log16a}
\xymatrix{
{\Omega^1_{(\co_{\oK},\mN_\infty)/(\co_K,\mN)}\otimes_{\co_{\oK}}\oR}\ar@{=}[d]\ar[r]^-(0.5)u&
{\Omega^1_{(R_\infty,P_\infty)/(R,P)}\otimes_{R_\infty}\oR}\ar[d]_a\ar@{->>}[r]&
{\Omega^1_{(R_\infty,P_\infty)/(R_1,P)}\otimes_{R_\infty}\oR}\ar[d]_b\\
{\Omega^1_{(\co_{\oK},\mN_\infty)/(\co_K,\mN)}\otimes_{\co_{\oK}}\oR}\ar[r]^-(0.5){u'}&
{\Omega^1_{(\oR,P_\infty)/(R,P)}}\ar@{->>}[r]&{\Omega^1_{(\oR,P_\infty)/(R_1,P)}}}
\end{equation}
Comme la suite \eqref{higgs1-log-ext5aa} est exacte et scindée, $u$ est injectif. 
D'autre part, le noyau de $a$ est annulé par $\fm_\oK$ en vertu de  \ref{higgs1-pur7} et \eqref{higgs1-log85a}.  
Par suite, le noyau de $u'$ est annulé par $\fm_\oK$. Comme 
$\Omega^1_{(\co_{\oK},\mN_\infty)/(\co_K,\mN)}\otimes_{\co_{\oK}}\oR$ n'a pas de $\fm_\oK$-torsion 
non nulle d'après \ref{higgs1-log-ext84}, $u'$ est injectif. 

Appliquant le foncteur ``module de Tate'' $\rT_p(-)=\Hom(\mQ_p/\mZ_p,-)$ au diagramme ci-dessus, 
on obtient un diagramme commutatif 
\begin{equation}\label{higgs1-ext-log16b}
\xymatrix{
{\rT_p(\Omega^1_{(\co_{\oK},\mN_\infty)/(\co_K,\mN)}\otimes_{\co_{\oK}}\oR)}\ar@{=}[d]\ar@{^(->}[r]&
{\rT_p(\Omega^1_{(R_\infty,P_\infty)/(R,P)}\otimes_{R_\infty}\oR)}\ar[d]_{\rT_p(a)}\ar[r]^-(0.5)v&
{\rT_p(\Omega^1_{(R_\infty,P_\infty)/(R_1,P)}\otimes_{R_\infty}\oR)}\ar[d]_{\rT_p(b)}\\
{\rT_p(\Omega^1_{(\co_{\oK},\mN_\infty)/(\co_K,\mN)}\otimes_{\co_{\oK}}\oR)}\ar@{^(->}[r]&
{\rT_p(\Omega^1_{(\oR,P_\infty)/(R,P)})}\ar[r]^-(0.5){v'}&{\rT_p(\Omega^1_{(\oR,P_\infty)/(R_1,P)})}}
\end{equation}
Comme la suite \eqref{higgs1-log-ext5aa} est exacte et scindée, $v$ est surjectif. D'autre part, 
le noyau et le conoyau de $b$ sont annulés par $\fm_\oK$, et $\Omega^1_{(R_\infty,P_\infty)/(R_1,P)}\otimes_{R_\infty}\oR$
n'a pas de $\fm_\oK$-torsion non nulle en vertu de \ref{higgs1-log-ext44}. 
Par suite, $b$ est injectif, et donc $\rT_p(b)$ est un isomorphisme. Il s'ensuit que $v'$ est surjectif et que $\rT_p(a)$
est un isomorphisme. 

Il résulte de \eqref{higgs1-log-ext8b} que le morphisme canonique 
\begin{equation} 
\rT_p(\Omega^1_{(\co_\oK,\mN_\infty)/(\co_K,\mN)})\otimes_{\co_C}\hoR\rightarrow 
\rT_p(\Omega^1_{(\co_\oK,\mN_\infty)/(\co_K,\mN)}\otimes_{\co_\oK}\oR)
\end{equation}
est un isomorphisme. De même, il résulte de \eqref{higgs1-log-ext4bb} que le morphisme canonique 
\begin{equation} 
\rT_p(\Omega^1_{(R_\infty,P_\infty)/(R_1,P)})\otimes_{\hRi}\hoR\rightarrow 
\rT_p(\Omega^1_{(R_\infty,P_\infty)/(R_1,P)}\otimes_{R_\infty}\oR)
\end{equation}
est un isomorphisme. Par suite, le morphisme canonique 
\begin{equation} 
\rT_p(\Omega^1_{(R_\infty,P_\infty)/(R,P)})\otimes_{\hRi}\hoR\rightarrow 
\rT_p(\Omega^1_{(R_\infty,P_\infty)/(R,P)}\otimes_{R_\infty}\oR)
\end{equation}
est un isomorphisme.

\subsection{} \label{higgs1-log-ext17}\index{10680@$\cE$ (extension de Faltings)} \index{Extension de Faltings}
On note $\cE$ la $\hoR$-représentation de $\Gamma$ définie par 
\begin{equation}\label{higgs1-ext-log17a}
\cE=\Hom(\mQ_p/\mZ_p,\Omega^1_{(\oR,P_\infty)/(R,P)})\otimes_{\mZ_p}\mZ_p(-1),
\end{equation}
où l'action de $\Gamma$ provient de son action sur  
$\Omega^1_{(\oR,P_\infty)/(R,P)}$ \eqref{higgs1-ext-log11}. Il résulte de \ref{higgs1-ext-log16} et \eqref{higgs1-ext-log12b} 
qu'on a une suite exacte canonique de $\hoR$-représentations de $\Gamma$ 
\begin{equation}\label{higgs1-log-ext17b}
0\rightarrow (\pi\rho)^{-1} \hoR\rightarrow \cE \rightarrow
\tOmega^1_{R/\co_K}\otimes_R \hoR(-1)\rightarrow 0,
\end{equation}
dite {\em extension de Faltings}. On a un isomorphisme de $\hoR$-représentations de $\Gamma$ 
\begin{equation}\label{higgs1-log-ext17c}
\cE_\infty\otimes_{\hRi}\hoR\stackrel{\sim}{\rightarrow} \cE,
\end{equation}
induisant un isomorphisme des extensions \eqref{higgs1-ext-log12b} et \eqref{higgs1-log-ext17b}. 
En particulier, la suite \eqref{higgs1-log-ext17b} est scindée en tant que suite de $\hoR$-modules (sans actions de $\Gamma$). 

\section{Cohomologie galoisienne}\label{higgs1-CG}

\begin{prop}\label{higgs1-cg1}
Soient $n$ un entier $\geq 0$, $\nu\colon \Delta_{p^\infty}\rightarrow \mu_{p^n}(\co_\oK)$ 
un homomorphisme surjectif, $\zeta$ un générateur du groupe $\mu_{p^n}(\co_\oK)$, $a\in \co_\oK$, 
$\fq$ l'idéal de $\co_\oK$ engendré par $a$ et $\zeta-1$. 
Soit $A$ une $\co_\oK$-algèbre complète et séparée pour la topologie $p$-adique et $\co_{\oK}$-plate. 
On note  $A(\nu)$ le $A$-$\Delta_{p^\infty}$-module topologique $A$, muni de la topologie $p$-adique
et de l'action de $\Delta_{p^\infty}$ définie par la multiplication par $\nu$ \eqref{higgs1-not55}. 
\begin{itemize}
\item[{\rm (i)}] Si $\nu=1$ (i.e., $n=0$), alors on a un isomorphisme canonique de $A$-algèbres graduées 
\begin{equation}\label{higgs1-cg1aa} 
\wedge(\Hom_{\mZ_p}(\Delta_{p^\infty},A/aA))
\stackrel{\sim}{\rightarrow}  \rH^*_\cont(\Delta_{p^\infty},A(\nu)/aA(\nu)).
\end{equation}
\item[{\rm (ii)}] Si $\nu\not=1$ (i.e., $n\not=0$), alors $\rH^i_\cont(\Delta_{p^\infty},A(\nu)/aA(\nu))$
est un $A/\fq A$-module libre de type fini pour tout $i\geq 0$ et est nul pour tout $i\geq \rg(L)=d+1$ \eqref{higgs1-dlog1f}.
\item[{\rm (iii)}] Le système projectif $(\rH^*(\Delta_{p^\infty},A(\nu)/p^rA(\nu)))_{r\geq 0}$ 
vérifie la condition de Mittag-Leffler uniformément en $\nu$, autrement dit, si on note, pour tous entiers $r'\geq r\geq 0$,
\begin{equation}
h_{r,r'}^\nu\colon \rH^*(\Delta_{p^\infty},A(\nu)/p^{r'}A(\nu))\rightarrow \rH^*(\Delta_{p^\infty},A(\nu)/p^rA(\nu))
\end{equation}
le morphisme canonique, alors pour tout entier $r\geq 1$,
il existe un entier $r'\geq r$, dépendant de $d$ mais pas de $\nu$, tel que pour tout entier $r''\geq r'$, 
les images de $h_{r,r'}^\nu$ et $h_{r,r''}^\nu$ soient égales. 
\end{itemize}
\end{prop}

Fixons une $\mZ_p$-base $e_1,\dots,e_d$ de $\Delta_{p^\infty}$ 
et notons $\rS_A(\Delta_{p^\infty})$ l'algèbre symétrique du $A$-module $\Delta_{p^\infty}\otimes_{\mZ_p}A$.
L'anneau $A$ est muni d'une structure de $\rS_A(\Delta_{p^\infty})$-algèbre définie par 
l'homomorphisme de $A$-algèbres $\rS_A(\Delta_{p^\infty})\rightarrow A$
qui envoie $e_i$ sur $\nu(e_i)-1$ pour $1\leq i\leq d$; on la note $A(\nu)^\triangleright$.
On observera que le $\rS_A(\Delta_{p^\infty})$-module 
sous-jacent à $A(\nu)^\triangleright$  n'est autre que le $\rS_A(\Delta_{p^\infty})$-module  
associé à $A(\nu)$ et noté avec le même symbole dans \ref{higgs1-cg01}. Considérons la forme linéaire 
\begin{equation}\label{higgs1-cg03b}
u\colon \Delta_{p^\infty}\otimes_{\mZ_p}A\rightarrow A
\end{equation}
qui envoie $e_i\otimes 1$ sur $\nu(e_i)-1$ pour $1\leq i\leq d$. Il résulte aussitôt des définitions \eqref{higgs1-koszul1d}
et \eqref{higgs1-koszul2dd} qu'on a un isomorphisme canonique 
\begin{equation}\label{higgs1-cg03d}
\mK_{\rS_A(\Delta_{p^\infty})}^\bullet((A(\nu)/aA(\nu))^\triangleright)\stackrel{\sim}{\rightarrow}\mK_{A}^\bullet(u,A/aA).
\end{equation}
En vertu de \ref{higgs1-cg01}, on en déduit un isomorphisme canonique dans $\bD^+(\bMod(A))$
\begin{equation}\label{higgs1-cg03e}
\rC^\bullet_\cont(\Delta_{p^\infty},A(\nu)/aA(\nu))\stackrel{\sim}{\rightarrow} \mK^\bullet_{A}(u,A/aA),
\end{equation}
où la source désigne le complexe de cochaînes continues de $G$ à valeurs dans $A(\nu)/aA(\nu)$ \eqref{higgs1-not78}.

(i) Cela résulte de \ref{higgs1-kun20} et \ref{higgs1-higgs551}(iii). 

On notera que la forme $u$ \eqref{higgs1-cg03b} étant nulle, 
\eqref{higgs1-cg03e} fournit un isomorphisme de $A$-modules gradués
\begin{equation}
\rH^*_\cont(\Delta_{p^\infty},A(\nu)/aA(\nu)) \stackrel{\sim}{\rightarrow} 
\wedge(\Hom_{\mZ_p}(\Delta_{p^\infty},A/aA)).
\end{equation}
Mais il n'est pas clair a priori que c'est un isomorphisme de $A$-algèbres graduées. 
On peut toutefois le déduire de \ref{higgs1-kun2} (ce qui correspond essentiellement à la preuve de \ref{higgs1-kun20}).

(ii) Notons aussi $u$ la forme linéaire $\Delta_{p^\infty}\otimes_{\mZ_p}(A/aA)\rightarrow (A/aA)$ 
déduite de $u$ \eqref{higgs1-cg03b}. 
En vertu de \eqref{higgs1-cg03e} et \eqref{higgs1-koszul1h}, il suffit de montrer que $\rH_i(\mK_\bullet^{A/aA}(u))$ 
est un $(A/\fq A)$-module libre de type fini pour tout $i\geq 0$ et est nul pour tout $i\geq d+1$.
La seconde proposition est évidente.
On sait aussi que $\rH_i(\mK_\bullet^{A/aA}(u))$ est annulé par $\fq$ \eqref{higgs1-koszul1}.
Pour $1\leq j\leq d$, posons $\zeta_j=\nu(e_j)$. On peut supposer que $\zeta_1=\zeta\not=1$ et que \eqref{higgs1-not1}
\[
v(\zeta_1-1)\leq v(\zeta_2-1)\leq \dots \leq v(\zeta_d-1).
\]
Procédons par récurrence sur $d$. La proposition pour $d=1$ est une conséquence immédiate de la platitude de
$A$ sur $\co_\oK$. Supposons $d\geq 2$ et l'assertion établie pour $d-1$. 
Notons $G$ le sous-$\mZ_p$-module de $\Delta_{p^\infty}$ engendré par $e_1,\dots,e_{d-1}$
et $u'\colon G\otimes_{\mZ_p}(A/aA)\rightarrow A/aA$ 
la restriction de $u$ à $G\otimes_{\mZ_p}(A/aA)$. D'après \eqref{higgs1-koszul1b}, on a un isomorphisme canonique
\begin{equation}\label{higgs1-cg1d}
\mK_\bullet^{A/aA}(\zeta_d-1)\otimes\mK_\bullet^{A/aA}(u')\stackrel{\sim}{\rightarrow}\mK_\bullet^{A/aA}(u),
\end{equation} 
où $\mK_\bullet^{A/aA}(\zeta_d-1)$ est le complexe de Koszul défini par la forme linéaire $\zeta_d-1$ de $A/aA$. 
En vertu de (\cite{ega3} 1.1.4.1), pour tout entier $i$, on a une suite exacte 
\begin{eqnarray}\label{higgs1-cg1e}
\lefteqn{0\rightarrow \rH_0(\mK_\bullet^{A/aA}(\zeta_d-1)\otimes\rH_i(\mK_\bullet^{A/aA}(u')))}\\
&&\rightarrow 
\rH_i(\mK_\bullet^{A/aA}(u))\rightarrow \rH_1(\mK_\bullet^{A/aA}(\zeta_d-1)
\otimes\rH_{i-1}(\mK_\bullet^{A/aA}(u')))\rightarrow 0.\nonumber
\end{eqnarray}
Par hypothèse de récurrence, $\rH_{i}(\mK_\bullet^{A/aA}(u'))$ est un 
$(A/\fq A)$-module libre de type fini pour tout $i\geq 0$. Comme $(\zeta_d-1)\in \fq$, 
on en déduit une suite exacte 
\begin{equation} \label{higgs1-cg1f}
0\rightarrow \rH_i(\mK_\bullet^{A/aA}(u'))\rightarrow \rH_i(\mK_\bullet^{A/aA}(u))\rightarrow 
\rH_{i-1}(\mK_\bullet^{A/aA}(u'))\rightarrow 0.
\end{equation}
Par suite, $\rH_{i}(\mK_\bullet^{A/aA}(u))$ étant annulé par $\fq$, 
il est libre de type fini sur $A/\fq A$.

(iii) Il résulte de (i) que le système projectif $(\rH^*(\Delta_{p^\infty},A(\nu)/p^rA(\nu)))_{r\geq 0}$ 
vérifie la condition de Mittag-Leffler lorsque $\nu=1$. On peut donc se borner aux caractères $\nu\not=1$. 
Posons $A_r=A/p^rA$ et notons aussi $u$
la forme linéaire $u\otimes_{A} \id_{A_r}\colon \Delta_{p^\infty}\otimes_{\mZ_p}A_r\rightarrow A_r$ \eqref{higgs1-cg03b}.
D'après \eqref{higgs1-cg03e} et \eqref{higgs1-koszul1h}, il suffit de montrer la proposition analogue pour 
le système projectif $(\rH_i(\mK_\bullet^{A_r}(u)))_{r\geq 1}$. Procédons par récurrence sur $d$.
Supposons d'abord $d=1$. L'homomorphisme canonique 
\begin{equation}\label{higgs1-cg1g}
\rH_1(\mK_\bullet^{A_{r+1}}(u))\rightarrow \rH_1(\mK_\bullet^{A_r}(u))
\end{equation}
est clairement bijectif, et comme $v(\zeta-1)\leq 1$, l'homomorphisme canonique 
\begin{equation}\label{higgs1-cg1h}
\rH_0(\mK_\bullet^{A_{r+1}}(u))\rightarrow \rH_0(\mK_\bullet^{A_r}(u))
\end{equation}
est nul. L'assertion est donc satisfaite avec $r'=r+1$. Supposons $d\geq 2$ et l'assertion établie pour $d-1$. 
L'assertion pour $d$ résulte alors facilement de la suite exacte \eqref{higgs1-cg1f} en prenant  $a=p^r$ pour $r\geq 0$ 
(cf. la preuve de \cite{ega3} 0.13.2.1).

\begin{cor}
Sous les hypothèses de \eqref{higgs1-cg1}, l'homomorphisme canonique 
\begin{equation}
\rH^*_\cont(\Delta_{p^\infty},A(\nu))\rightarrow \underset{\underset{r\geq 0}{\longleftarrow}}{\lim}\ 
\rH^*(\Delta_{p^\infty},A(\nu)/p^rA(\nu))
\end{equation}
est un isomorphisme. 
\end{cor}
En effet, d'après \eqref{higgs1-limproj2c} et \eqref{higgs1-limproj2d}, pour tout 
$i\geq 0$, on a une suite exacte 
\begin{eqnarray}
\lefteqn{0\rightarrow \rR^1 \underset{\underset{r\geq 0}{\longleftarrow}}{\lim}\ 
\rH^{i-1}(\Delta_{p^\infty},A(\nu)/p^rA(\nu))\rightarrow}\\
&& \rH^{i}_\cont(\Delta_{p^\infty},A(\nu))
\rightarrow \underset{\underset{r\geq 0}{\longleftarrow}}{\lim}\ 
\rH^{i}(\Delta_{p^\infty},A(\nu)/p^rA(\nu))\rightarrow 0,\nonumber
\end{eqnarray}
dont le terme de gauche est nul en vertu de \ref{higgs1-cg1}(iii) et \eqref{higgs1-limproj2e}. 

\begin{rema}\label{higgs1-cg04}
Sous les hypothèses de \eqref{higgs1-cg1}, si $\nu\not=1$, l'homomorphisme canonique 
\begin{equation}
\rH^*_\cont(\Delta_{p^\infty},A(\nu))\otimes_{A}A/p^rA\rightarrow
\rH^*(\Delta_{p^\infty},A(\nu)/p^rA(\nu))
\end{equation}
n'est pas en général un isomorphisme. 
\end{rema}

\begin{prop}\label{higgs1-cg14}
Soient $n$ un entier $\geq 1$, $\nu\colon \Delta_{p^\infty}\rightarrow \mu_{p^n}(\co_\oK)$ un homomorphisme 
surjectif, $\zeta$ un générateur du groupe $\mu_{p^n}(\co_\oK)$,
$a$ un élément non nul de $\co_\oK$, $b=a(\zeta-1)^{-1}$, $\alpha$ un nombre rationnel. 
Soient $A$ une $\co_\oK$-algèbre complète et séparée pour la topologie $p$-adique et $\co_{\oK}$-plate, 
$N$ un $(A/aA)$-module, $M$ un $(A/aA)$-$\Delta_{p^\infty}$-module discret. 
On note $N(\nu)$ le $A$-$\Delta_{p^\infty}$-module discret $N$, 
muni de l'action de $\Delta_{p^\infty}$ définie par la multiplication par $\nu$. 
Supposons que les conditions suivantes soient remplies~:
\begin{itemize}
\item[{\rm (i)}] $\inf(v(a),\alpha)> v(\zeta-1)$; 
\item[{\rm (ii)}] $N$ est plat sur $\co_\oK/a\co_\oK$;
\item[{\rm (iii)}] $M$ est projectif de type fini sur $A/aA$, et est engendré par un nombre fini d'éléments 
$\Delta_{p^\infty}$-invariants modulo $p^{\alpha} M$.   
\end{itemize}
Alors pour tout $i\geq 0$, on a 
\[
(\zeta-1)\cdot \rH^i(\Delta_{p^\infty},(M/bM)\otimes_{A}N(\nu))=0.
\]
\end{prop}

Comme $M$ est un facteur direct d'un $(A/aA)$-module libre de type fini, on peut se borner au cas où 
il est libre de type fini sur $A/aA$. Il admet donc une $(A/aA)$-base formée d'éléments 
$\Delta_{p^\infty}$-invariants modulo $p^{\alpha} M$. Posons $T=M\otimes_AN(\nu)$. 
Fixons une $\mZ_p$-base $e_1,\dots,e_d$ de $\Delta_{p^\infty}$ 
et notons $\rS_A(\Delta_{p^\infty})$ l'algèbre symétrique du $A$-module $\Delta_{p^\infty}\otimes_{\mZ_p}A$.
En vertu de \ref{higgs1-cg01}, on a un isomorphisme canonique dans $\bD^+(\bMod(A))$
\begin{equation}\label{higgs1-cg14a}
\rC^\bullet(\Delta_{p^\infty},T/bT)\stackrel{\sim}{\rightarrow} 
\mK^\bullet_{\rS_A(\Delta_{p^\infty})}((T/bT)^\triangleright).
\end{equation}
Il existe un entier $i$ tel que $1\leq i\leq d$ et que $\nu(e_i)$ 
soit un générateur de $\mu_{p^n}(\co_\oK)$; on peut supposer $\nu(e_i)=\zeta$.
Notons $\varphi\colon \Delta_{p^\infty}\rightarrow \Aut_{A}(M)$ la représentation de $\Delta_{p^\infty}$
sur $M$. Il existe alors un endomorphisme $A$-linéaire $U$ de $M$ 
tel que $\varphi(e_i)=\id_M+p^{\alpha}U$. Posons $c=p^\alpha(\zeta-1)^{-1}\in \fm_\oK$ et 
$V=\id_{T}+c U\otimes (\zeta \cdot \id_N)$, qui est un automorphisme $A$-linéaire de $T$, de sorte qu'on a 
\begin{equation}\label{higgs1-cg14b}
\varphi(e_i)\otimes (\nu(e_i) \cdot \id_N)-\id_{T}=(\zeta-1)V\in \End_{A}(T).
\end{equation}
Pour tout entier $1\leq j\leq d$, on a 
\begin{eqnarray}\label{higgs1-cg14c}
\lefteqn{(\varphi(e_j)\otimes (\nu(e_j) \cdot \id_N)-\id_{T})\circ V^{-1}}\\
&=&V^{-1}\circ (\varphi(e_j)\otimes (\nu(e_j) \cdot \id_N)-\id_{T}) \in 
\End_{A}(T/bT).\nonumber
\end{eqnarray}
En effet, $\Delta_{p^\infty}$ étant abélien, les produits de ces endomorphismes par $\zeta-1$ sont égaux 
dans $\End_{A}(T)$, ce qui implique la relation \eqref{higgs1-cg14c} car $T$ est plat sur $\co_{\oK}/a\co_\oK$. 
Donc $V^{-1}$ est un automorphisme du $\rS_A(\Delta_{p^\infty})$-module $(T/bT)^\triangleright$,
et il induit, pour tout $q\geq 0$, un automorphisme $\rS_A(\Delta_{p^\infty})$-linéaire de 
$\rH^q(\mK^\bullet_{\rS_A(\Delta_{p^\infty})}((T/bT)^\triangleright))$. Or ce module est annulé par $e_i$   
\eqref{higgs1-koszul1}. Par suite, pour tout $x\in \rH^q(\mK^\bullet_{\rS_A(\Delta_{p^\infty})}((T/bT)^\triangleright))$, 
on a, en vertu de \eqref{higgs1-cg14b},
\begin{equation}\label{higgs1-cg14d}
e_i\cdot V^{-1}(x)=(\zeta-1)x=0.
\end{equation}

\subsection{}\label{higgs1-cg4}
On désigne par $\Lambda$ le conoyau dans la catégorie des monoïdes 
de l'homomorphisme $\vartheta\colon \mN\rightarrow P$ et par $q\colon P\rightarrow \Lambda$
l'homomorphisme canonique. D'après (\cite{ogus} I 1.1.5),
$\Lambda$ est le quotient de $P$ par la relation de congruence $E$ formée des éléments $(x,y)\in P\times P$ 
pour lesquels il existe $a,b\in \mN$ tels que $x+a\lambda=y+b\lambda$. Dire que $E$ est une relation 
de congruence signifie qu'elle est une relation d'équivalence et que $E$ est un sous-monoïde de $P\times P$.  
Le groupe associé à $\Lambda$ s'identifie canoniquement à $P^\gp/\mZ\lambda$. 
Comme $P$ est intègre, $\Lambda$ est intègre~; on peut donc l'identifier à l'image 
de $P$ dans $P^\gp/\mZ\lambda$.

\begin{lem}\label{higgs1-cg5}
Conservons les notations de \eqref{higgs1-cg4}. Alors~:

{\rm (i)}\ Le monoïde $\Lambda$ est saturé.

{\rm (ii)}\ Pour tout $x\in \Lambda$, l'ensemble $q^{-1}(x)$ admet un unique élément minimal $\tx$ 
pour la relation de pré-ordre de $P$ définie par sa structure de monoïde.

{\rm (iii)}\ Pour tout $x\in \Lambda$ et tout entier $n\geq 0$, on a $\widetilde{nx}=n\tx$. 
\end{lem}
(i) En effet, $\Lambda$ est la somme amalgamée 
de l'homomorphisme saturé $\vartheta$ et de $\mN\rightarrow 0$. Il est donc saturé \eqref{higgs1-log111}. 

(ii) On notera d'abord que deux éléments quelconques de $q^{-1}(x)$ sont nécessairement comparables \eqref{higgs1-cg4}. 
Montrons que l'ensemble $q^{-1}(x)$ admet un élément minimal $\tx$. 
Il revient au même de dire que pour tout $t\in P$, il existe $n\in \mN$ tel que l'élément $t-n\lambda$ de $P^\gp$
n'appartient pas à $P$. En effet, si ce n'était pas le cas, pour tout $n\geq 0$, l'élément $\alpha(t)/\pi^n$ de $R_K$ 
appartient à $R$, où $\alpha\colon P\rightarrow R$ est l'homomorphisme défini par la carte $(P,\gamma)$ \eqref{higgs1-dlog1}. 
Comme $\alpha(t)\not=0$ et que $R$ est noethérien et intègre, on en déduit que $\pi$ est inversible dans $R$, 
ce qui contredit l'hypothèse \ref{higgs1-dlog1}(C$_2$). 
Il s'ensuit aussi que $-\lambda$ n'appartient pas à $P$. 
Donc $\tx$ est nécessairement unique car $P$ est intègre et $P^\gp$ est sans torsion.

(iii) Comme $-\lambda$ n'appartient pas à $P$, on a $\widetilde{0}=0$. 
On peut donc se borner au cas où $n$ est un nombre premier. On notera d'abord que $q(\widetilde{nx})=q(n\tx)$.
Si $\widetilde{nx}\not=n\tx$, alors il existe $m\geq 1$ tel que $n\tx\geq m \lambda$. 
Comme l'homomorphisme $\vartheta$ est saturé, il existe alors $m'\in \mN$ tel que $\tx\geq m'\lambda$ et $nm'\geq m$
en vertu de (\cite{tsuji4} 4.1 page 11; cf. aussi \cite{ogus} I 4.8.13). Comme $m'\geq 1$, la relation  $\tx\geq m'\lambda$ 
contredit le caractère minimal de $\tx$. 

\subsection{}\label{higgs1-cg51}\index{10820@$\Xi_{p^n}$, $\Xi_{p^\infty}$}
Conservons les notations de \eqref{higgs1-cg4}. 
Pour tout $n\geq 0$, on note $\Lambda^{(p^n)}$ le monoïde au-dessus de $\Lambda$ défini par le couple
$(\Lambda,\varpi_{p^n})$, autrement dit, $\Lambda^{(p^n)}$ est le monoïde $\Lambda$
et l'homomorphisme structural $\Lambda \rightarrow \Lambda^{(p^n)}$ est le Frobenius d'ordre $p^n$ 
de $\Lambda$ \eqref{higgs1-not-mon1}. On identifie naturellement $\Lambda^{(p^n)}$ au conoyau de l'homomorphisme 
$\vartheta\colon \mN^{(p^n)}\rightarrow P^{(p^n)}$ et on note aussi $q\colon P^{(p^n)}\rightarrow \Lambda^{(p^n)}$
l'homomorphisme canonique. Posons 
\[
P_{p^\infty}=\underset{\underset{\mN}{\longrightarrow}}{\lim}\ P^{(p^n)} \ \ \ {\rm et} \ \ \ 
\Lambda_{p^\infty}=\underset{\underset{\mN}{\longrightarrow}}{\lim}\ \Lambda^{(p^n)}.
\]
On identifie $P_{p^\infty}$ à un sous-monoïde de $P_\infty$ \eqref{higgs1-log-ext2b}. 
On note aussi $q\colon P_{p^\infty}\rightarrow \Lambda_{p^\infty}$ la limite inductive des homomorphismes canoniques 
$q\colon P^{(p^n)}\rightarrow \Lambda^{(p^n)}$. 

Pour tout $t\in P_{p^\infty}$, l'application 
\begin{equation}\label{higgs1-cg3g}
\nu_t\colon \Delta_{\infty}\rightarrow \mu_{p^\infty}(\co_\oK)=
\underset{\underset{n\geq 0}{\longrightarrow}}{\lim}\ \mu_{p^n}(\co_{\oK}), \ \ \ g\mapsto \langle g,t\rangle,
\end{equation} 
où $\langle g,t\rangle$ est défini dans \eqref{higgs1-ext-log10a}, induit, compte tenu de \eqref{higgs1-ext-log10c}, un homomorphisme
\begin{equation}\label{higgs1-cg3gg}
\nu_t\colon \Delta_{p^\infty}\rightarrow \mu_{p^\infty}(\co_\oK).
\end{equation} 
Il est clair que l'application 
\begin{equation}\label{higgs1-cg3hh}
P_{p^\infty}\rightarrow \Hom(\Delta_{p^\infty},\mu_{p^\infty}(\co_\oK)),\ \ \ t\mapsto \nu_t
\end{equation}
est un homomorphisme. Pour tout $n\geq 0$, notons 
$\lambda^{(p^n)}$ l'image de $\lambda$ dans $P^{(p^n)}$ par l'isomorphisme \eqref{higgs1-ext-log2d}. 
Comme $\nu_{\lambda^{(p^n)}}=1$, \eqref{higgs1-cg3hh} induit un homomorphisme que l'on note aussi
\begin{equation}\label{higgs1-cg3k}
\Lambda_{p^\infty}\rightarrow \Hom(\Delta_{p^\infty},\mu_{p^\infty}(\co_\oK)),\ \ \ x\mapsto \nu_x.
\end{equation}
En fait, celui-ci induit, pour tout $n\geq 0$, un homomorphisme 
\begin{equation}\label{higgs1-cg3kk}
\Lambda^{(p^n)}\rightarrow \Hom(\Delta_{p^\infty},\mu_{p^n}(\co_\oK)).
\end{equation}

Posons 
\begin{eqnarray}
\Xi_{p^\infty}&=&\Hom(\Delta_{p^\infty},\mu_{p^\infty}(\co_\oK)),\label{higgs1-cg3kl}\\
\Xi_{p^n}&=&\Hom(\Delta_{p^\infty},\mu_{p^n}(\co_\oK)).\label{higgs1-cg3km}
\end{eqnarray}
On identifie $\Xi_{p^n}$ à un sous-groupe de $\Xi_{p^\infty}$. 

\begin{lem}\label{higgs1-cg52}
Les hypothèses étant celles de \eqref{higgs1-cg51}, soient de plus $x\in \Lambda_{p^\infty}$, $n$ un entier $\geq 0$. 
Alors pour que $\nu_x(\Delta_{p^\infty})\subset \mu_{p^n}(\co_\oK)$ il faut et il suffit que 
$x\in \Lambda^{(p^n)}\subset \Lambda_{p^\infty}$; en particulier, pour que $ \nu_x=1$, il faut et
il suffit que $x\in \Lambda^{(1)}\subset \Lambda_{p^\infty}$.
\end{lem}

Supposons que $x\in \Lambda^{(p^{n+m})}$ pour un entier $m\geq 0$. 
En vertu de  \eqref{higgs1-ext-log10d}, on a un diagramme commutatif
\begin{equation}\label{higgs1-cg52a}
\xymatrix{
{\Lambda^{(p^n)}}\ar[r]^{\varpi_{p^m}}\ar[d]&{\Lambda^{(p^{n+m})}}\ar[d]\ar[rdd]^-(0.4)a&\\
{P^\gp/\mZ\lambda}\ar[r]^{\cdot p^m}\ar[d]&{P^\gp/\mZ\lambda}\ar[d]&\\
{\Hom(\Delta_{p^\infty},\mZ_p(1))}\ar[r]^{\cdot p^m}&{\Hom(\Delta_{p^\infty},\mZ_p(1))}\ar@{->>}[r]^-(0.4)b&
{\Hom(\Delta_{p^\infty},\mu_{p^{n+m}}(\co_\oK))}}
\end{equation}
où les flèches verticales inférieures proviennent de l'identification de 
$\Delta_{p^\infty}$ et $L_\lambda\otimes_\mZ\mZ_p(1)$ \eqref{higgs1-gal2}, $a$ est l'homomorphisme \eqref{higgs1-cg3kk}
et $b$ est le morphisme canonique. Le carré inférieur est cartésien car le sous-groupe de torsion
de $P^\gp/\mZ\lambda$ est d'ordre premier à $p$, et le carré supérieur est cartésien car $\Lambda$ est saturé
en vertu de \ref{higgs1-cg5}(i). Pour que $\nu_x(\Delta_{p^\infty})\subset \mu_{p^n}(\co_\oK)$, il faut et il suffit que 
$a(x)$ soit dans l'image de $p^mb$. La proposition s'obtient alors par chasse au diagramme \eqref{higgs1-cg52a}. 

\begin{lem}\label{higgs1-cg16}\index{10825@$R_{p^\infty}^{(\nu)}$}
Il existe une décomposition canonique de $R_{p^\infty}$ en somme directe de $R_1$-modules de présentation finie,
stables sous l'action de $\Delta_{p^\infty}$,
\begin{equation}\label{higgs1-cg16a}
R_{p^\infty}=\bigoplus_{\nu\in \Xi_{p^\infty}}R_{p^\infty}^{(\nu)},
\end{equation}
telle que l'action de $\Delta_{p^\infty}$ sur le facteur $R_{p^\infty}^{(\nu)}$ soit donnée par le caractère $\nu$. 
De plus, pour tout $n\geq 0$, on a 
\begin{equation}\label{higgs1-cg16b}
R_{p^n}=\bigoplus_{\nu\in \Xi_{p^n}}R_{p^\infty}^{(\nu)}. 
\end{equation}
\end{lem}

Pour tout $n\geq 1$, posons $C_n=A_n\otimes_{\co_{K_n}}\co_\oK$ \eqref{higgs1-dlog3b} et notons 
\begin{equation}\label{higgs1-cg16aa}
\lambda_n\colon L\otimes_\mZ\mu_{n}(\oK)\rightarrow \mu_{n}(\oK)
\end{equation}
l'homomorphisme défini par $\lambda$. On a un isomorphisme canonique  \eqref{higgs1-dlog4e}
\begin{equation}\label{higgs1-cg16d}
C_n\stackrel{\sim}{\rightarrow}\co_{\oK}[P^{(n)}]/(\pi_{n}-e^{\lambda^{(n)}})\otimes_{\co_{\oK}[P]/(\pi-e^\lambda)}C_1,
\end{equation}
où $\lambda^{(n)}$ désigne l'image de $\lambda$ dans $P^{(n)}$ par l'isomorphisme canonique \eqref{higgs1-ext-log2d}. 
Le groupe $L\otimes_\mZ\mu_{n}(\oK)$ agit naturellement
sur $\co_{\oK}[P^{(n)}]$ par des homomorphismes de $(\co_{\oK}[P])$-algèbres.
On en déduit, compte tenu de \eqref{higgs1-cg16d}, une action du groupe $\ker(\lambda_n)$ sur $C_n$ 
par des homomorphismes de $C_1$-algèbres. Posons
\begin{equation}\label{higgs1-cg16ac}
C_{p^\infty}=\underset{\underset{n\geq 0}{\longrightarrow}}{\lim}\ C_{p^n},
\end{equation}
où la limite inductive est indexée par l'ensemble $\mN$ ordonné par la relation d'ordre habituelle.
On a des isomorphismes canoniques \eqref{higgs1-gal1b}
\begin{equation}\label{higgs1-cg16ad}
\Delta_{p^\infty}\stackrel{\sim}{\rightarrow} L_\lambda\otimes\mZ_p(1)\stackrel{\sim}{\rightarrow}
\underset{\underset{n\geq 0}{\longleftarrow}}{\lim}\ \ker(\lambda_{p^n}).
\end{equation}
On obtient par passage à la limite une action de $\Delta_{p^\infty}$ sur $C_{p^\infty}$ 
par des homomorphismes de $C_1$-algèbres. 

En vertu de \ref{higgs1-gal1}(iii), pour tout entier $n\geq 0$, on a un isomorphisme 
$\ker(\lambda_{p^n})$-équivariant canonique 
\begin{equation}\label{higgs1-cg16ab}
R_{p^n}\stackrel{\sim}{\rightarrow}C_{p^n}\otimes_{C_1}R_1.
\end{equation}
On en déduit un isomorphisme $\Delta_{p^\infty}$-équivariant
\begin{equation}\label{higgs1-cg16ae}
R_{p^\infty}\stackrel{\sim}{\rightarrow}C_{p^\infty}\otimes_{C_1}R_1.
\end{equation}
Il suffit donc de montrer qu'il existe une décomposition canonique de $C_{p^\infty}$ en somme directe de 
$C_1$-modules de présentation finie, stables sous l'action de $\Delta_{p^\infty}$,
\begin{equation}\label{higgs1-cg16af}
C_{p^\infty}=\bigoplus_{\nu\in \Xi_{p^\infty}}C_{p^\infty}^{(\nu)},
\end{equation}
telle que l'action de $\Delta_{p^\infty}$ sur le facteur $C_{p^\infty}^{(\nu)}$ soit donnée par le caractère $\nu$; 
de plus, pour tout $n\geq 0$, on a 
\begin{equation}\label{higgs1-cg16ag}
C_{p^n}=\bigoplus_{\nu\in \Xi_{p^n}}C_{p^\infty}^{(\nu)}. 
\end{equation}

Compte tenu de \eqref{higgs1-cg16d}, on peut se réduire au cas où $R=\co_{K}[P]/(\pi-e^{\lambda})$.
On notera que $X$ n'est plus nécessairement connexe après cette réduction.
Soit $\Lambda$ le conoyau dans la catégorie des monoïdes 
de l'homomorphisme $\vartheta\colon \mN\rightarrow P$ et soit $q\colon P\rightarrow \Lambda$
l'homomorphisme canonique (cf. \ref{higgs1-cg4}). On a alors 
\begin{equation}\label{higgs1-cg16e}
R=\bigoplus_{x\in \Lambda}\co_K\cdot \alpha(\tx),
\end{equation}  
où $\tx$ est le relèvement minimal de $x$ dans $P$ défini dans \ref{higgs1-cg5}(ii) et 
$\alpha\colon P\rightarrow R$ est l'homomorphisme induit par la carte $(P,\gamma)$ \eqref{higgs1-dlog1}. 
Pour tout $n\geq 0$, on a d'après \eqref{higgs1-cg16d}
\begin{equation}\label{higgs1-cg16f}
C_{p^n}=\co_{\oK}[P^{(p^n)}]/(\pi_{p^n}-e^{\lambda^{(p^n)}}).
\end{equation}
Notons $\alpha_{p^n}\colon P^{(p^n)}\rightarrow C_{p^n}$ l'homomorphisme 
induit par le morphisme strict canonique \eqref{higgs1-dlog3a}
\[
(X_{n},\cM_{X_{n}})\rightarrow \bA_P.
\] 
Cette notation est compatible avec celle introduite dans \eqref{higgs1-log-ext2} et n'induit aucune confusion. 
On obtient par passage à la limite inductive un homomorphisme 
\begin{equation}
\alpha_{p^\infty}\colon P_{p^\infty}\rightarrow C_{p^\infty}.
\end{equation}
Reprenons les notations de \eqref{higgs1-cg51} et notons, pour tout $x\in \Lambda^{(p^n)}$, $\tx\in P^{(p^n)}$ 
l'unique élément minimal de $q^{-1}(x)$ pour la relation de pré-ordre de $P^{(p^n)}$ 
définie par sa structure de monoïde (cf. \ref{higgs1-cg5}(ii)). On a alors
\begin{equation}\label{higgs1-cg16g}
C_{p^n}=\bigoplus_{x\in \Lambda^{(p^n)}}\co_\oK\cdot \alpha_{p^n}(\tx).
\end{equation}

En vertu de \ref{higgs1-cg5}(iii), les applications $\Lambda^{(p^n)}\rightarrow P^{(p^n)}, x\mapsto \tx$ sont compatibles. 
Elles définissent donc par passage à la limite inductive une application que l'on note aussi 
\begin{equation}\label{higgs1-cg16gh}
\Lambda_{p^\infty}\rightarrow P_{p^\infty}, \ \ \ x\mapsto \tx.
\end{equation}
On a clairement $q(\tx)=x$. 
Comme $\Lambda$ est intègre et que le sous-groupe de torsion de $P^\gp/\mZ\lambda$ est d'ordre premier à $p$, 
les homomorphismes $\Lambda^{(p^n)}\rightarrow \Lambda^{(p^{n+1})}$ sont injectifs. 
Par passage à la limite inductive de la décomposition \eqref{higgs1-cg16g}, on obtient 
\begin{equation}\label{higgs1-cg16h}
C_{p^\infty}=\bigoplus_{x\in \Lambda_{p^\infty}}\co_\oK\cdot \alpha_{p^\infty}(\tx).
\end{equation}  
Chaque élément $\alpha_{p^\infty}(\tx)\in C_{p^\infty}$ est un vecteur propre pour l'action de 
$\Delta_{p^\infty}$. Il résulte aussitôt de \eqref{higgs1-ext-log10d}  que l'action de $\Delta_{p^\infty}$ 
sur $\alpha_{p^\infty}(\tx)$ est donnée par le caractère $\nu_x$. 
Pour tout $\nu\in \Xi_{p^\infty}$, on pose alors 
\begin{equation}\label{higgs1-cg16i}
C_{p^\infty}^{(\nu)}=\bigoplus_{x\in \Lambda_{p^\infty}| \nu_x=\nu}\co_\oK\cdot \alpha_{p^\infty}(\tx),
\end{equation}
de sorte que l'on a 
\begin{equation}
C_{p^\infty}=\bigoplus_{\nu\in \Xi_{p^\infty}}C_{p^\infty}^{(\nu)}.
\end{equation}
Comme l'application $\Lambda_{p^\infty}\rightarrow \Xi_{p^\infty},x\mapsto\nu_x$ 
est un homomorphisme, on voit que $C_{p^\infty}^{(\nu)}$ est un sous-$C_1$-module de $C_{p^\infty}$.
Il résulte de \ref{higgs1-cg52} et \eqref{higgs1-cg16g} que pour tout $n\geq 0$, on a 
\begin{equation}
C_{p^n}=\bigoplus_{\nu\in \Xi_{p^n}}C_{p^\infty}^{(\nu)}.
\end{equation}
Comme $C_{p^n}$ est de présentation finie sur $C_1$ pour tout $n\geq 0$, $C_{p^\infty}^{(\nu)}$
est de présentation finie sur $C_1$ pour tout $\nu\in \Xi_{p^\infty}$.

\begin{teo}\label{higgs1-cg3}
Soient $a$ un élément non nul de $\co_\oK$, $\zeta$ une racine primitive $p$-ième de $1$ dans $\co_\oK$. 
Alors~: 

{\rm (i)}\ Il existe un et un unique homomorphisme de $R_1$-algèbres graduées
\begin{equation}\label{higgs1-cg3a}
\wedge(\Hom_{\mZ}(\Delta_{p^\infty},R_1/aR_1))\rightarrow 
\rH^*(\Delta_{p^\infty},R_{p^\infty}/aR_{p^\infty})
\end{equation}
dont la composante en degré un est le composé des morphismes canoniques
\begin{equation}
\Hom_{\mZ}(\Delta_{p^\infty},R_1/aR_1)\stackrel{\sim}{\rightarrow}
\rH^1(\Delta_{p^\infty},R_1/aR_1)\rightarrow 
\rH^1(\Delta_{p^\infty},R_{p^\infty}/aR_{p^\infty}).
\end{equation}
Celui-ci admet, en tant que morphisme de $R_1$-modules gradués,  un inverse à gauche canonique
\begin{equation}\label{higgs1-cg3b}
\rH^*(\Delta_{p^\infty},R_{p^\infty}/aR_{p^\infty})\rightarrow 
\wedge(\Hom_{\mZ}(\Delta_{p^\infty},R_1/aR_1)),
\end{equation}
dont le noyau est annulé par $\zeta-1$.  

{\rm (ii)}\ Le $R_1$-module $\rH^i(\Delta_{p^\infty},R_{p^\infty}/a R_{p^\infty})$ est presque de présentation finie 
pour tout $i\geq 0$, et est nul pour tout $i\geq \rg(L)=d+1$ \eqref{higgs1-dlog1f}. 

{\rm (iii)}\ Le système projectif $(\rH^*(\Delta_{p^\infty},R_{p^\infty}/p^rR_{p^\infty}))_{r\geq 0}$ vérifie la condition 
de Mittag-Leffler~; plus précisément, si on note, pour tous entiers $r'\geq r\geq 0$,
\begin{equation}
h_{r,r'}\colon \rH^*(\Delta_{p^\infty},R_{p^\infty}/p^{r'} R_{p^\infty})\rightarrow \rH^*(\Delta_{p^\infty},R_{p^\infty}/p^r R_{p^\infty})
\end{equation}
le morphisme canonique, alors pour tout entier $r\geq 1$,
il existe un entier $r'\geq r$, dépendant seulement de $d$ mais pas des autres données dans \eqref{higgs1-dlog1}, 
tel que pour tout entier $r''\geq r'$, les images de $h_{r,r'}$ et $h_{r,r''}$ soient égales.  
\end{teo}

En effet, d'après \ref{higgs1-cg16}, on a une décomposition canonique de $R_{p^\infty}$ en 
somme directe de $R_1[\Delta_{p^\infty}]$-modules
\begin{equation}\label{higgs1-cg3c}
R_{p^\infty}=\bigoplus_{\nu\in \Xi_{p^\infty}}R_{p^\infty}^{(\nu)}\otimes_{\co_{\oK}}\co_{\oK}(\nu),
\end{equation}
où $\Delta_{p^\infty}$ agit trivialement sur $R_{p^\infty}^{(\nu)}$ et agit sur $\co_\oK(\nu)=\co_\oK$ par le caractère $\nu$. 
Comme les $R_{p^\infty}^{(\nu)}$ sont $\co_{\oK}$-plats, on a,  en vertu de \ref{higgs1-cg2}, une décomposition canonique 
en somme directe de $R_1$-modules
\begin{equation}\label{higgs1-cg3d}
\rH^*(\Delta_{p^\infty},R_{p^\infty}/aR_{p^\infty})=\bigoplus_{\nu\in \Xi_{p^\infty}}
\rH^*(\Delta_{p^\infty},\co_{\oK}(\nu)/a\co_\oK(\nu))\otimes_{\co_{\oK}}R_{p^\infty}^{(\nu)}. 
\end{equation}

(i) On a $R_{p^\infty}^{(1)}=R_1$ \eqref{higgs1-cg16b}, de sorte que la composante en $\nu=1$ de la décomposition \eqref{higgs1-cg3d}
est l'image de l'homomorphisme canonique de $R_1$-algèbres graduées 
\begin{equation}
\rH^*(\Delta_{p^\infty},\co_{\oK}/a\co_\oK)\otimes_{\co_{\oK}}R_1 \rightarrow 
\rH^*(\Delta_{p^\infty},R_{p^\infty}/aR_{p^\infty}).
\end{equation}
Par ailleurs, l'homomorphisme canonique de $R_1$-algèbres graduées
\begin{equation}
\wedge(\Hom_{\mZ}(\Delta_{p^\infty},\co_\oK/a\co_\oK))\otimes_{\co_\oK} R_1
\rightarrow \wedge(\Hom_{\mZ}(\Delta_{p^\infty},R_1/aR_1))
\end{equation}
est un isomorphisme. La proposition résulte alors de \eqref{higgs1-cg3d} et \ref{higgs1-cg1} (appliqué avec $A=\co_C$).

(ii) Soient $i,n$ deux entiers $\geq 0$, $\zeta_n$ une racine primitive $p^n$-ième de l'unité. 
Il résulte de \eqref{higgs1-cg3d}, \ref{higgs1-cg1} et \ref{higgs1-cg16} que $\rH^i(\Delta_{p^\infty},R_{p^\infty}/aR_{p^\infty})$ 
est la somme directe d'un $R_1$-module de présentation
finie et d'un $R_1$-module annulé par $\zeta_n-1$. Il est donc presque de présentation finie sur $R_1$.
La seconde assertion est claire puisque 
la $p$-dimension cohomologique de $\Delta_{p^\infty}$ est égale à $d$.

(iii) Cela résulte de \eqref{higgs1-cg3d} et \ref{higgs1-cg1}(iii). 

\begin{cor}\label{higgs1-cg17}
Pour tout élément non nul $a$ de $\co_\oK$ et tout entier $i\geq 0$, le noyau et le conoyau du morphisme canonique 
\begin{equation}
\rH^i(\Delta_{p^\infty},R_1/a R_1)\rightarrow \rH^i(\Delta,\oR/a \oR)
\end{equation}
sont annulés par $\fm_\oK$ et $p^{\frac{1}{p-1}}\fm_\oK$ respectivement. 
\end{cor}

En effet, le morphisme canonique 
\begin{equation}
\rH^i(\Delta_{p^\infty},R_{p^\infty}/a R_{p^\infty})\rightarrow \rH^i(\Delta,\oR/a \oR)
\end{equation}
est un presque-isomorphisme d'après \ref{higgs1-pur121}. D'autre part, il résulte de \eqref{higgs1-cg3d}, 
\ref{higgs1-cg1} (appliqué avec $A=\co_C$), \ref{higgs1-cg2} et du fait que $R_{p^\infty}^{(1)}=R_1$ \eqref{higgs1-cg16b} que le morphisme canonique 
\begin{equation}
\rH^i(\Delta_{p^\infty},R_1/a R_1)\rightarrow \rH^i(\Delta_{p^\infty},R_{p^\infty}/aR_{p^\infty})
\end{equation}
est injectif de conoyau annulé par $\zeta-1$, où $\zeta$ est une
racine primitive $p$-ième de l'unité.
La proposition s'ensuit  car $v(\zeta-1)=\frac{1}{p-1}$. 

\begin{cor}\label{higgs1-cg7}
L'homomorphisme canonique 
\begin{equation}\label{higgs1-cg7a}
\rH^*_\cont(\Delta_{p^\infty},\hRpi)\rightarrow \underset{\underset{r\geq 0}{\longleftarrow}}{\lim}\ 
\rH^*(\Delta_{p^\infty},R_{p^\infty}/p^rR_{p^\infty})
\end{equation}
est un isomorphisme.  
\end{cor}
En effet, d'après \eqref{higgs1-limproj2c} et \eqref{higgs1-limproj2d}, pour tout 
$i\geq 0$, on a une suite exacte 
\[
0\rightarrow \rR^1 \underset{\underset{r\geq 0}{\longleftarrow}}{\lim}\ 
\rH^{i-1}(\Delta_{p^\infty},R_{p^\infty}/p^rR_{p^\infty})\rightarrow \rH^{i}_\cont(\Delta_{p^\infty},\hRpi)
\rightarrow \underset{\underset{r\geq 0}{\longleftarrow}}{\lim}\ 
\rH^{i}(\Delta_{p^\infty},R_{p^\infty}/p^rR_{p^\infty})\rightarrow 0,
\]
dont le terme de gauche est nul en vertu de \ref{higgs1-cg3}(iii) et \eqref{higgs1-limproj2e}. 

\subsection{}\label{higgs1-cg35}
Soit $a$ un élément non nul de $\co_\oK$. 
Comme le sous-groupe de torsion de $P^\gp/\mZ\lambda$ est d'ordre premier à $p$,  
on a un isomorphisme canonique \eqref{higgs1-gal2}
\begin{equation}\label{higgs1-cg35b}
(P^\gp/\mZ\lambda)\otimes_\mZ\mZ_p(-1)\stackrel{\sim}{\rightarrow}
\Hom_{\mZ_p}(\Delta_{p^\infty}, \mZ_p).
\end{equation}
On en déduit, compte tenu de \eqref{higgs1-log-ext45b}, des isomorphismes $\hRun$-linéaires
\begin{eqnarray}
\tOmega^1_{R/\co_K}\otimes_R(R_1/aR_1)(-1)&\stackrel{\sim}{\rightarrow} &
\Hom_{\mZ}(\Delta_{p^\infty},R_1/aR_1),\label{higgs1-cg35c}\\
\tOmega^1_{R/\co_K}\otimes_R\hRun(-1)&\stackrel{\sim}{\rightarrow} &
\Hom_{\mZ}(\Delta_{p^\infty},\hRun).\label{higgs1-cg35d}
\end{eqnarray}
On interprétera dans la suite les buts de ces morphismes comme des groupes de cohomologie \eqref{higgs1-higgs551}. 
Par ailleurs, le composé de \eqref{higgs1-cg35d} avec l'isomorphisme canonique \eqref{higgs1-higgs551a}
\begin{equation}\label{higgs1-cg35e}
\Hom_{\mZ}(\Delta_{p^\infty},\hRun) \stackrel{\sim}{\rightarrow} \Hom_{\mZ}(\Delta_\infty,\hRun)
\end{equation}
n'est autre que le morphisme $\tdelta(-1)$ \eqref{higgs1-ext-log14h}, d'après \ref{higgs1-ext-log14}.

\begin{prop}\label{higgs1-cg6}
Il existe un et un unique homomorphisme de $\hRun$-algèbres graduées
\begin{equation}\label{higgs1-cg6a}
\wedge(\tOmega^1_{R/\co_K}\otimes_R\hRun(-1))\rightarrow \rH^*_\cont(\Delta_{p^\infty},\hRpi)
\end{equation}
dont la composante en degré un est induite par \eqref{higgs1-cg35d}.
Celui-ci admet, en tant que morphisme de $\hRun$-modules gradués, un inverse à gauche canonique 
\begin{equation}\label{higgs1-cg6b}
\rH^*_\cont(\Delta_{p^\infty},\hRpi)\rightarrow \wedge(\tOmega^1_{R/\co_K}\otimes_R\hRun(-1)),
\end{equation}
dont le noyau est annulé par $p^{\frac{1}{p-1}}$.  
\end{prop}

Cela résulte de \ref{higgs1-cg3}(i) et \ref{higgs1-cg7}.

\begin{cor}\label{higgs1-cg61}
On a $(\hRpi)^{\Delta_{p^\infty}}=\hRun$. 
\end{cor}
En effet, d'après \ref{higgs1-cg6}, l'homomorphisme canonique $\hRun\rightarrow (\hRpi)^{\Delta_{p^\infty}}$ admet un inverse à gauche 
$(\hRpi)^{\Delta_{p^\infty}}\rightarrow \hRun$ dont le noyau
est annulé par $p^{\frac{1}{p-1}}$. Comme $\hRpi$ est $\co_C$-plat \eqref{higgs1-pur8}, l'inverse à gauche 
est injectif, d'où la proposition. 

\begin{rema}\label{higgs1-cg31}
D'après \ref{higgs1-cg07}, \eqref{higgs1-gal3a} et \eqref{higgs1-higgs551a}, le théorème \ref{higgs1-cg3} et ses corollaires \ref{higgs1-cg17} et \ref{higgs1-cg7}
valent encore si on remplace $\Delta_{p^\infty}$ par $\Delta_\infty$
et $R_{p^\infty}$ par $R_\infty$. 
Il en est alors de même de la proposition \ref{higgs1-cg6} et de son corollaire \ref{higgs1-cg61}. On a donc $(\hRi)^{\Delta_{\infty}}=\hRun$. 
\end{rema}

\begin{prop}\label{higgs1-cg18}
Soit $a$ un élément non nul de $\co_\oK$. 
\begin{itemize}
\item[{\rm (i)}] Il existe un et un unique homomorphisme de $R_1$-algèbres graduées
\begin{equation}
\wedge(\tOmega^1_{R/\co_K}\otimes_R(R_1/aR_1)(-1))\rightarrow \rH^*(\Delta,\oR/a\oR)
\end{equation}
dont la composante en degré un est induite par \eqref{higgs1-cg35c}. Celui-ci 
est presque-injectif et son conoyau est annulé par $p^{\frac{1}{p-1}}\fm_\oK$.  
\item[{\rm (ii)}] Le $R_1$-module $\rH^i(\Delta,\oR/a \oR)$ est presque de présentation finie 
pour tout $i\geq 0$, et est presque-nul pour tout $i\geq d+1$. 

\item[{\rm (iii)}] Pour tous entiers $r'\geq r\geq 0$, notons
\begin{equation}
\hbar_{r,r'}\colon \rH^*(\Delta,\oR/p^{r'} \oR)\rightarrow \rH^*(\Delta,\oR/p^r \oR)
\end{equation}
le morphisme canonique. Alors pour tout entier $r\geq 1$,
il existe un entier $r'\geq r$, dépendant seulement de $d$ mais pas des autres données dans \eqref{higgs1-dlog1}, 
tel que pour tout entier $r''\geq r'$, les images de $\hbar_{r,r'}$ et $\hbar_{r,r''}$ soient presque-isomorphes.  
\end{itemize}
\end{prop}
Cela résulte de \ref{higgs1-pur121} et \ref{higgs1-cg3}.

\subsection{}\label{higgs1-cg10}
Pour tout $\co_\oK$-module $M$, on pose 
\begin{equation}\label{higgs1-cg10a}
M^\flat=\Hom_{\co_\oK}(\fm_\oK,M).
\end{equation}
Le morphisme canonique $M\rightarrow M^\flat$ est un presque-isomorphisme. 
Pour qu'un morphisme de $\co_\oK$-modules $u\colon M\rightarrow N$ soit un presque-isomorphisme,
il faut et il suffit que le morphisme associé $u^\flat\colon M^\flat\rightarrow N^\flat$ soit un isomorphisme (\cite{tsuji2} 2.5). 

\begin{lem}\label{higgs1-cg11}
Le morphisme canonique 
\begin{equation}
\jmath\colon \hRun\rightarrow (\hRun)^\flat
\end{equation}
est un isomorphisme. 
\end{lem}
Comme $\hRun$ est plat sur $\co_C$ \eqref{higgs1-pur8}, $\jmath$ est injectif. Montrons que $\jmath$ est surjectif.
Soit $u\in (\hRun)^\flat$. Pour tout $\alpha \in \mQ_{>0}$, posons $x_\alpha=u(p^\alpha)$. 
Rappelons que $\hRun[\frac 1 p]$ est une algèbre affinoïde sur $C$ et que 
$\hRun$ est la boule unité pour la norme $|\ |_{\sup}$ sur $\hRun[\frac 1 p]$ (cf. la preuve de \ref{higgs1-pur81}). 
Les relations $x_{\alpha+\beta}=p^\alpha x_{\beta}=p^\beta x_\alpha$ ($\alpha,\beta\in \mQ_{>0}$)
impliquent alors que 
\begin{equation}
|p^{-\alpha}x_\alpha|_{\sup} \leq 1.
\end{equation}
Par suite, $x=p^{-\alpha}x_\alpha\in \hRun$ et est indépendant de $\alpha$. Il est clair que $\jmath(x)=u$. 

\begin{rema}
Le morphisme canonique $\co_C\rightarrow (\co_C)^\flat$ est un isomorphisme. 
La preuve est une variante très simple de celle de \ref{higgs1-cg11}.
\end{rema}

\begin{prop}\label{higgs1-cg9}
Il existe un et un unique homomorphisme de $\hRun$-algèbres graduées
\begin{equation}\label{higgs1-cg9a}
\wedge(\tOmega^1_{R/\co_K}\otimes_R\hRun(-1))\rightarrow \rH^*_\cont(\Delta,\hoR)
\end{equation}
dont la composante en degré un est induite par \eqref{higgs1-cg35d}.
Celui-ci admet, en tant que morphisme de $\hRun$-modules gradués, un inverse à gauche canonique 
\begin{equation}\label{higgs1-cg9b}
\rH^*_\cont(\Delta,\hoR)\rightarrow \wedge(\tOmega^1_{R/\co_K}\otimes_R\hRun(-1)),
\end{equation}
dont le noyau est annulé par $p^{\frac{1}{p-1}}$.  
\end{prop}

En effet, on a un diagramme commutatif d'homomorphismes de $\hRun$-algèbres graduées
\begin{equation}\label{higgs1-cg9c}
\xymatrix{
{\wedge(\tOmega^1_{R/\co_K}\otimes_R\hRun(-1))}\ar[r]^-(0.5){v}\ar[d]_{\iota}&
{\rH^*_\cont(\Delta_{p^\infty},\hRpi)}\ar[r]^-(0.5)u\ar[d]\ar@/_1pc/[l]_-(0.5){w}&
{\rH^*_\cont(\Delta,\hoR)}\ar[d]\\
{(\wedge(\tOmega^1_{R/\co_K}\otimes_R\hRun(-1)))^\flat}\ar[r]^-(0.5){v^\flat}&{(\rH^*_\cont(\Delta_{p^\infty},\hRpi))^\flat}
\ar[r]^-(0.5){u^\flat}\ar@/^1pc/[l]^-(0.5){w^\flat}& {(\rH^*_\cont(\Delta,\hoR))^\flat}}
\end{equation}
où le foncteur $(\ )^\flat$ est défini dans \eqref{higgs1-cg10a}, les flèches verticales et $u$ sont les morphismes canoniques,
$v$ est l'homomorphisme \eqref{higgs1-cg6a} et $w$ est la section \eqref{higgs1-cg6b} de $v$. 
Alors $\iota$ et $u^\flat$ sont des isomorphismes (\ref{higgs1-cg11} et \ref{higgs1-cg8}). 
D'autre part, il résulte de \ref{higgs1-cg6} que le noyau de $w^\flat$ est annulé par $p^{\frac{1}{p-1}}$. 
La proposition s'ensuit par une chasse au diagramme \eqref{higgs1-cg9c} puisque $u\circ v$ est l'homomorphisme \eqref{higgs1-cg9a}.

\begin{cor}\label{higgs1-cg12}
{\rm (i)}\ On a $(\hoR)^{\Delta}=\hRun$.

{\rm (ii)}\ Le sous-module de torsion $p$-primaire $M^1$ de $\rH^1_\cont(\Delta,\hoR)$ est égal au noyau du morphisme 
\begin{equation}\label{higgs1-cg12a}
\iota\colon \rH^1_\cont(\Delta,\hoR)\rightarrow \rH^1_\cont(\Delta,(\pi\rho)^{-1}\hoR)
\end{equation}
déduit de l'injection canonique $\hoR\subset (\pi\rho)^{-1}\hoR$.

{\rm (iii)}\  Le morphisme composé 
\begin{equation}\label{higgs1-cg12b}
\tOmega^1_{R/\co_K}\otimes_R\hRun(-1)\longrightarrow
\rH^1_\cont(\Delta,\hoR)\stackrel{\iota}{\longrightarrow} 
\rH^1_\cont(\Delta,(\pi\rho)^{-1}\hoR),
\end{equation}
où la première flèche est induite par \eqref{higgs1-cg35d}, 
est le bord de la suite exacte longue de cohomologie déduite de la suite exacte courte \eqref{higgs1-log-ext17b}.

{\rm (iv)}\ Pour tout entier $i\geq 0$, notons $M^i$ le sous-module de torsion $p$-primaire de $\rH^i_\cont(\Delta,\hoR)$. 
Il existe alors un et un unique isomorphisme de $\hRun$-algèbres graduées
\begin{equation}\label{higgs1-cg12c}
\wedge(\tOmega^1_{R/\co_K}\otimes_R\hRun(-1))\stackrel{\sim}{\rightarrow} \oplus_{i\geq 0} 
(\rH^i_\cont(\Delta,\hoR)/M^i)
\end{equation}
tel que le composé 
\begin{equation}
\tOmega^1_{R/\co_K}\otimes_R\hRun(-1)\rightarrow
\rH^1_\cont(\Delta,\hoR)/M^1\rightarrow  \rH^1_\cont(\Delta,(\pi\rho)^{-1}\hoR),
\end{equation}
où la première flèche est la composante en degré un de \eqref{higgs1-cg12c} et la seconde flèche est induite par $\iota$, 
soit le bord de la suite exacte longue de cohomologie déduite de la suite exacte courte \eqref{higgs1-log-ext17b}.
\end{cor}

(i) D'après \ref{higgs1-cg9}, l'homomorphisme canonique $\hRun\rightarrow (\hoR)^{\Delta}$ admet un inverse à gauche 
$(\hoR)^{\Delta}\rightarrow \hRun$ dont le noyau
est annulé par $p^{\frac{1}{p-1}}$. Comme $\hoR$ est $\co_C$-plat \eqref{higgs1-pur8}, l'inverse à gauche 
est injectif, d'où la proposition. 

(ii) Il résulte de \ref{higgs1-cg9} que $M^1$ est annulé par $p^{\frac{1}{p-1}}$. D'autre part, la multiplication 
par $\pi\rho$ dans $\hoR$ est le composé de l'injection canonique $\hoR\rightarrow (\pi\rho)^{-1}\hoR$
et de l'isomorphisme 
$(\pi\rho)^{-1}\hoR\stackrel{\sim}{\rightarrow}\hoR$ qui, pour tout $x\in \hoR$, associe $x$ à $(\pi\rho)^{-1}x$. 
Comme $v(\rho)\geq \frac{1}{p-1}$, on en déduit que $\ker(\iota)=M^1$.  

(iii) Cela résulte de la définition de \eqref{higgs1-cg35d}, \eqref{higgs1-ext-log15b} et \eqref{higgs1-log-ext17c}.

(iv) Cela résulte de (ii), (iii) et \ref{higgs1-cg9}.

\begin{prop}\label{higgs1-cg15}
Soient $M$ un $R_1$-$\Delta_{p^\infty}$-module discret, 
$a$ un élément non nul de $\co_\oK$, $\alpha$ un nombre rationnel. Supposons que $\inf(v(a),\alpha)> \frac{1}{p-1}$ 
et que $M$ soit un $(R_1/aR_1)$-module projectif de type fini engendré par un nombre fini d'éléments 
$\Delta_{p^\infty}$-invariants modulo $p^{\alpha} M$.   
Posons $b=ap^{-\frac{1}{p-1}}$. Alors pour tout $i\geq 0$, le noyau et le conoyau du morphisme canonique 
\begin{equation}
\rH^i(\Delta_{p^\infty},M/bM)\rightarrow \rH^i(\Delta,(M/bM)\otimes_{R_1}\oR)
\end{equation}
sont annulés par $\fm_\oK$ et $p^{\frac{1}{p-1}}\fm_\oK$ respectivement. 
\end{prop}
Comme $M$ est un facteur direct d'un $(R_1/aR_1)$-module libre de type fini, on a, d'après \ref{higgs1-cg07}, 
\begin{equation}
((M/bM)\otimes_{R_1}R_{\infty})^{\Sigma_0}=(M/bM)\otimes_{R_1}R_{p^\infty}.
\end{equation}
Donc en vertu de \eqref{higgs1-gal3a}, le morphisme canonique 
\begin{equation}
\rH^i(\Delta_{p^\infty},(M/bM)\otimes_{R_1}R_{p^\infty})\rightarrow \rH^i(\Delta_\infty,(M/bM)\otimes_{R_1}R_{\infty})
\end{equation}
est un isomorphisme. D'autre part, d'après \ref{higgs1-pur12}, le morphisme canonique 
\begin{equation}
(M/bM)\otimes_{R_1}R_{\infty} \rightarrow ((M/bM)\otimes_{R_1}\oR)^{\Sigma}
\end{equation}
est un presque-isomorphisme. Par suite, en vertu de \ref{higgs1-pur6}, le morphisme canonique 
\begin{equation}
\rH^i(\Delta_{\infty},(M/bM)\otimes_{R_1}R_{\infty})\rightarrow \rH^i(\Delta,(M/bM)\otimes_{R_1}\oR)
\end{equation}
est un presque-isomorphisme. Il suffit donc de montrer que le morphisme canonique 
\begin{equation}
\rH^i(\Delta_{p^\infty},M/bM)\rightarrow \rH^i(\Delta_{p^\infty},(M/bM)\otimes_{R_1}R_{p^\infty})
\end{equation}
est injectif de conoyau annulé par $p^{\frac{1}{p-1}}$. 

D'après \ref{higgs1-cg16}, on a une décomposition canonique de $(M/bM)\otimes_{R_1}R_{p^\infty}$ en 
somme directe de $R_1[\Delta_{p^\infty}]$-modules
\begin{equation}
(M/bM)\otimes_{R_1}R_{p^\infty}=\bigoplus_{\nu\in \Xi_{p^\infty}}(M/bM)\otimes_{R_1}
R_{p^\infty}^{(\nu)}\otimes_{R_1}R_1(\nu),
\end{equation}
où $\Delta_{p^\infty}$ agit trivialement sur $R_{p^\infty}^{(\nu)}$ et agit sur $R_1(\nu)=R_1$ par le caractère $\nu$.
Comme $R_{p^\infty}^{(1)}=R_1$ \eqref{higgs1-cg16b}, la proposition résulte de \ref{higgs1-cg14} 
(appliqué à $A=R_1/aR_1$ et $N=R_{p^\infty}^{(\nu)}/aR_{p^\infty}^{(\nu)}$ pour $\nu\not=1$).

\section{\texorpdfstring{\'Epaississements infinitésimaux $p$-adiques de Fontaine}
{\'Epaississements infinitésimaux p-adiques de Fontaine}}\label{higgs1-EIPF}

\subsection{}\label{higgs1-eip1}
Commençons par rappeler la construction suivante due à Grothendieck (\cite{grot1} IV 3.3).
Soient $A$ une $\mZ_{(p)}$-algèbre commutative, $n$ un entier $\geq 1$. 
L'homomorphisme d'anneaux \eqref{higgs1-not33a}
\begin{equation}\label{higgs1-eip1a}
\Phi_{n+1}\colon 
\begin{array}[t]{clcr}
\rW_{n+1}(A/p^nA)&\rightarrow& A/p^nA\\
(x_1,\dots,x_{n+1})&\mapsto& x_1^{p^n}+p x_2^{p^{n-1}}+\dots+p^{n}x_{n+1}
\end{array}
\end{equation}
s'annule sur $\rV^n(A/p^nA)$ et induit donc par passage au quotient un homomorphisme d'anneaux
\begin{equation}\label{higgs1-eip1b}
\Phi'_{n+1}\colon 
\begin{array}[t]{clcr}
\rW_{n}(A/p^nA)&\rightarrow& A/p^nA,\\
(x_1,\dots,x_n)&\mapsto&x_1^{p^n}+p x_2^{p^{n-1}}+\dots+p^{n-1}x_n^p.
\end{array}
\end{equation}
Ce dernier s'annule sur 
\begin{equation}\label{higgs1-eip1c}
\rW_n(pA/p^nA)=\ker(\rW_n(A/p^nA)\rightarrow \rW_n(A/pA))
\end{equation}
et se factorise à son tour en un homomorphisme d'anneaux
\begin{equation}\label{higgs1-eip1d}
\theta_n\colon \rW_{n}(A/pA)\rightarrow A/p^nA.
\end{equation}
Il résulte aussitôt de la définition que le diagramme 
\begin{equation}\label{higgs1-eip1e}
\xymatrix{
{\rW_{n+1}(A/pA)}\ar[r]^-(0.4){\theta_{n+1}}\ar[d]_{\rR\rF}&{A/p^{n+1}A}\ar[d]\\
{\rW_{n}(A/pA)}\ar[r]^-(0.4){\theta_n}&{A/p^nA}}
\end{equation}
où $\rR$ est le morphisme de restriction \eqref{higgs1-not33b}, 
$\rF$ est le Frobenius \eqref{higgs1-not33d} et la flèche non libellée est l'homomorphisme canonique, est commutatif. 

Pour tout homomorphisme de $\mZ_{(p)}$-algèbres commutatives $\varphi\colon A\rightarrow B$, 
le diagramme  
\begin{equation}\label{higgs1-eip1f}
\xymatrix{
{\rW_n(A/pA)}\ar[r]\ar[d]_{\theta_n}&{\rW_n(B/pB)}\ar[d]^{\theta_n}\\
{A/p^nA}\ar[r]&{B/p^nB}}
\end{equation}
où les flèches horizontales sont les morphismes induits par $\varphi$, est commutatif.

\begin{prop}\label{higgs1-eip2}
Soit $A$ une $\mZ_{(p)}$-algèbre commutative vérifiant les conditions suivantes~:
\begin{itemize}
\item[{\rm (i)}] $A$ est $\mZ_{(p)}$-plat.
\item[{\rm (ii)}] $A$ est intégralement clos dans $A[\frac 1 p]$.
\item[{\rm (iii)}] Le Frobenius absolu  de $A/pA$ est surjectif. 
\item[{\rm (iv)}] Il existe un entier $N\geq 1$ et une suite $(p_n)_{0\leq n\leq N}$ d'éléments de $A$
tels que $p_0=p$ et $p_{n+1}^p=p_n$ pour tout $0\leq n\leq N-1$.  
\end{itemize}
Pour tout entier $0\leq n\leq N$, on pose  
\begin{equation}\label{higgs1-eip2a}
\xi_n=[\op_n]-p \in \rW_n(A/pA),
\end{equation}
où $\op_n$ est la classe de $p_n$ dans $A/pA$ et $[\ ]$ désigne le représentant multiplicatif.
Alors pour tous entiers $n\geq 1$ et $i\geq 0$ tels que $n+i\leq N$, la suite
\begin{equation}\label{higgs1-eip2b}
\xymatrix{
{\rW_n(A/pA)}\ar[rr]^-(0.5){\cdot \rR^i(\xi_{n+i})}&&{\rW_n(A/pA)}
\ar[rr]^-(0.5){\theta_n\circ \rF^i}&&{A/p^nA}\ar[r]& 0}
\end{equation}
est exacte.
\end{prop}
En effet, on a $\rF \rR(\xi_{n+1})=\xi_n$ et $\theta_n (\rF^i(\rR^i(\xi_{n+i})))=\theta_n(\xi_n)=0$. 
Pour établir l'exactitude de \eqref{higgs1-eip2b}, procédons par récurrence sur $n$. 
Soient $i$ un entier tel que $0\leq i\leq N-1$, $\alpha,\beta\in A$ tels que 
$\alpha^{p^{i+1}}=p\beta$. Il résulte des hypothèses (i) et (ii) que $\alpha\in p_{i+1}A$.  
Par conséquent, la suite
\begin{equation}\label{higgs1-eip2c}
\xymatrix{
{A/pA}\ar[r]^-(0.5){\cdot \op_{i+1}}&{A/pA}\ar[r]^{\rF^{i+1}}&{A/pA}\ar[r]&0}
\end{equation}
est exacte, et la proposition est vraie pour $n=1$.  
Supposons ensuite la proposition vraie pour $1\leq n\leq N-1$ (et tout $0\leq i\leq N-n$) et montrons-là pour $n+1$.
Soit $i$ un entier tel que $0\leq i\leq N-n-1$. On a un diagramme commutatif à lignes exactes \eqref{higgs1-eip1e}
\begin{equation}\label{higgs1-eip2d}
\xymatrix{
0\ar[r]&{A/pA}\ar[r]^-(0.5){\rV^n}\ar[d]^{\cdot \rR^i(\xi_{i+1})}&{\rW_{n+1}(A/pA)}\ar[r]^{\rR}
\ar[d]^{\cdot \rR^i(\xi_{n+i+1)}}&
{\rW_n(A/pA)}\ar[d]^{\cdot \rR^{i+1}(\xi_{n+i+1})}\ar[r]&0\\
0\ar[r]&{A/pA}\ar[r]^-(0.5){\rV^n}\ar[d]^{\rF^{i+1}}&{\rW_{n+1}(A/pA)}\ar[r]^{\rR}\ar[dd]^{\theta_{n+1}\circ \rF^i}&
{\rW_n(A/pA)}\ar[r]\ar[dd]^{\theta_n\circ \rF^{i+1}}&0\\
&{A/pA}\ar[d]^{\cdot p^n}&&&\\
0\ar[r]&{p^nA/p^{n+1}A}\ar[r]&{A/p^{n+1}A}\ar[r]&{A/p^nA}\ar[r]&0}
\end{equation}
L'hypothèse de récurrence et \eqref{higgs1-eip2c} impliquent alors
la proposition pour $n+1$ par le lemme du serpent.

\subsection{}\label{higgs1-eipo3}\index{10934@$\cR_A$}\index{10935@$\theta\colon \rW(\cR_A)\rightarrow \hA$}\index{10936@$\cA_2(A)$}
Soit $A$ une $\mZ_{(p)}$-algèbre commutative. 
On désigne par $\cR_A$ la limite projective du système projectif $(A/pA)_{\mN}$ 
dont les morphismes de transition sont les itérés de l'endomorphisme de Frobenius de $A/pA$.
\begin{equation}\label{higgs1-eipo3a}
\cR_A=\underset{\underset{x\mapsto x^p}{\longleftarrow}}{\lim}A/pA.
\end{equation} 
C'est un anneau parfait de caractéristique $p$. 
Pour tout entier $n\geq 1$, la projection canonique $\cR_A\rightarrow A/pA$ sur la $(n+1)$-ième composante
du système projectif $(A/pA)_{\mN}$ ({\em i.e.}, la composante d'indice $n$) induit un homomorphisme
\begin{equation}\label{higgs1-eipo3b}
\nu_n\colon \rW(\cR_A)\rightarrow \rW_n(A/pA).
\end{equation}
Comme $\nu_n=\rF\circ\rR\circ \nu_{n+1}$, on obtient par passage à la limite projective un homomorphisme 
\begin{equation}\label{higgs1-eipo3c}
\nu\colon \rW(\cR_A)\rightarrow \underset{\underset{n\geq 0}{\longleftarrow}}{\lim}\ \rW_n(A/pA),
\end{equation}
où les morphismes de transition de la limite projective sont les morphismes $\rF \rR$. On vérifie aussitôt  
qu'il est bijectif. Compte tenu de \eqref{higgs1-eip1e},
les homomorphismes $\theta_n$ induisent par passage à la limite projective un 
homomorphisme 
\begin{equation}\label{higgs1-eipo3d}
\theta\colon \rW(\cR_A)\rightarrow \hA,
\end{equation}
où $\hA$ est le séparé complété $p$-adique de $A$. On retrouve l'homomorphisme défini par Fontaine 
(\cite{fontaine1} 2.2). On pose
\begin{equation}\label{higgs1-eipo3e}
\cA_2(A)=\rW(\cR_A)/\ker(\theta)^2,
\end{equation}
et on note encore $\theta\colon \cA_2(A)\rightarrow \hA$ l'homomorphisme induit par $\theta$ (cf. \cite{fontaine3} 1.2.2). 

Pour tout homomorphisme de $\mZ_{(p)}$-algèbres commutatives $\varphi\colon A\rightarrow B$, 
le diagramme  
\begin{equation}\label{higgs1-eip3f}
\xymatrix{
{\rW(\cR_A)}\ar[r]\ar[d]_{\theta}&{\rW(\cR_B)}\ar[d]^{\theta}\\
{\hA}\ar[r]&{\hB}}
\end{equation}
où les flèches horizontales sont les morphismes induits par $\varphi$, est commutatif \eqref{higgs1-eip1f}.  
La correspondance $A\mapsto \cA_2(A)$ est donc fonctorielle. 

\begin{rema}
La projection canonique $\cR_{\rW(k)}\rightarrow k$ sur la première composante ({\em i.e.}, d'indice 0)
est un isomorphisme. Elle induit donc un isomorphisme $\rW(\cR_{\rW(k)})\stackrel{\sim}{\rightarrow}\rW(k)$,
que nous utilisons pour identifier ces deux anneaux. L'homomorphisme $\theta$ s'identifie alors à
l'endomorphisme de Frobenius de $\rW(k)$. 
\end{rema}

\begin{prop}[\cite{tsuji1} A.1.1 et A.2.2]\label{higgs1-eip4}\index{10939@$\xi\in \rW(\cR_A)$}
Soit $A$ une $\mZ_{(p)}$-algèbre commutative vérifiant les conditions suivantes~:
\begin{itemize}
\item[{\rm (i)}] $A$ est $\mZ_{(p)}$-plat.
\item[{\rm (ii)}] $A$ est intégralement clos dans $A[\frac 1 p]$.
\item[{\rm (iii)}] Le Frobenius absolu  de $A/pA$ est surjectif. 
\item[{\rm (iv)}] Il existe une suite $(p_n)_{n\geq 0}$ d'éléments de $A$
tels que $p_0=p$ et $p_{n+1}^p=p_n$ pour tout $n\geq 0$.  
\end{itemize}
On désigne par $\upp$ l'élément de $\cR_A$ induit par la suite $(p_n)_{n\geq 0}$ et on pose  
\begin{equation}\label{higgs1-eip4a}
\xi=[\upp]-p \in \rW(\cR_A),
\end{equation}
où $[\ ]$ est le représentant multiplicatif.
Alors la suite
\begin{equation}\label{higgs1-eip4b}
0\longrightarrow \rW(\cR_A)\stackrel{\cdot \xi}{\longrightarrow} \rW(\cR_A)\stackrel{\theta}{\longrightarrow} 
\hA \longrightarrow 0
\end{equation}
est exacte.
\end{prop}

En effet, pour tout $n\geq 1$, posant 
\begin{equation}\label{higgs1-eip4c}
\xi_n=[\op_n]-p \in \rW_n(A/pA),
\end{equation}
où $\op_n$ est la classe de $p_n$ dans $A/pA$, la suite 
\begin{equation}\label{higgs1-eip4d}
\rW_n(A/pA)\stackrel{\cdot \xi_n}{\longrightarrow}\rW_n(A/pA)
\stackrel{\theta_n}{\longrightarrow}A/p^nA\longrightarrow 0
\end{equation}
est exacte en vertu de \ref{higgs1-eip2}. Comme l'homomorphisme 
$\rR\rF\colon \rW_{n+1}(A/pA)\rightarrow \rW_n(A/pA)$ est surjectif pour tout $n\geq 0$, 
la suite \eqref{higgs1-eip4b} est exacte au centre et à droite (\cite{ega3} 0.13.2.1(i) et 0.13.2.2). 

Si $a=(a_0,a_1,a_2,\dots)\in \rW(\cR_A)$ est tel que $\xi a=0$, alors  
\begin{equation}\label{higgs1-eip4e}
(\upp a_0,\upp^pa_1,\upp^{p^2}a_2,\dots)=(0,a_0^p,a_1^p,\dots). 
\end{equation}
Pour montrer que $\xi$ n'est pas un diviseur de zéro dans $\rW(\cR_A)$, il suffit donc de montrer que $\upp$ 
n'est pas un diviseur de zéro dans $\cR_A$.  Soit $y=(y_n)_{n\in \mN}\in \cR_A$ tel que $\upp y=0$. 
Pour tout $n\geq 0$, soit $\ty_n$ un relèvement de $y_n$ dans $A$. On a $p_n\ty_n\in pA$. 
Par suite, $\ty_n\in p_n^{p^n-1}A$ car $p$ n'est pas un diviseur de zéro dans $A$. On en déduit que 
\begin{equation}\label{higgs1-eip4f}
y_n=y_{n+1}^p=(\ty_{n+1}^p \mod pA)=0
\end{equation} 
car $p^{n+2}-p\geq p^{n+1}$.

\subsection{}\label{higgs1-eip5}
Soient $Y=(Y,\cM_Y)$ un $\mZ_{(p)}$-schéma logarithmique affine d'anneau $A$, 
$M$ un monoïde, $u\colon M\rightarrow \Gamma(Y,\cM_Y)$ un homomorphisme.  
Considérons le système projectif de monoïdes multiplicatifs $(A)_{n\in \mN}$, 
où les morphismes de transition sont tous égaux à l'élévation à la puissance $p$-ième.
On désigne par $Q$ le produit fibré du diagramme d'homomorphismes de monoïdes 
\begin{equation}\label{higgs1-eip5a}
\xymatrix{
&{M}\ar[d]\\
{\underset{\underset{x\mapsto x^p}{\longleftarrow}}{\lim}\ A}\ar[r]&A}
\end{equation}
où  la flèche horizontale est la projection sur la première composante ({\em i.e.}, d'indice $0$)
et la flèche verticale est le composé de $u$ et de l'homomorphisme canonique $\Gamma(Y,\cM_Y)\rightarrow A$, 
et par $q$ l'homomorphisme composé
\begin{equation}\label{higgs1-eip5b}
Q\longrightarrow
\underset{\underset{x\mapsto x^p}{\longleftarrow}}{\lim}\ A \longrightarrow \cR_A \stackrel{[\ ]}{\longrightarrow} 
\rW(\cR_A),
\end{equation} 
où les première et deuxième flèches sont les homomorphismes canoniques \eqref{higgs1-eipo3a}
et $[\ ]$ est le représentant multiplicatif. Il résulte aussitôt des définitions que le diagramme 
\begin{equation}\label{higgs1-eip5c}
\xymatrix{
Q\ar[r]\ar[d]_{q}&{M}\ar[d]\\
{\rW(\cR_A)}\ar[r]^-(0.4)\theta&{\hA}}
\end{equation}
où les flèches non libellées sont les morphismes canoniques, est commutatif. 

On munit $\hY=\Spec(\hA)$ de la structure logarithmique $\cM_{\hY}$ image inverse de $\cM_Y$ et 
$\Spec(\rW(\cR_A))$ de la structure logarithmique $\cQ$ associée à la structure pré-logarithmique définie par $q$ 
\eqref{higgs1-eip5b}. D'après \eqref{higgs1-eip5c}, $\theta$ induit un morphisme 
\begin{equation}\label{higgs1-eip5d}
(\hY,\cM_\hY)\rightarrow (\Spec(\rW(\cR_A)),\cQ).
\end{equation}

L'énoncé suivant est inspiré de (\cite{tsuji1} 1.4.2).

\begin{prop}\label{higgs1-eip6}
Conservons les hypothèses de \eqref{higgs1-eip5}, notons de plus $Y^\circ$ l'ouvert maximal de $Y$ où la structure logarithmique 
$\cM_Y$ est triviale et supposons les conditions suivantes remplies~:
\begin{itemize}
\item[{\rm (a)}] $A$ est intègre et normal. 
\item[{\rm (b)}] $Y^\circ$ est un $\mQ$-schéma non-vide et simplement connexe. 
\item[{\rm (c)}] $M$ est intègre et il existe un monoïde fin et saturé $M'$ et un homomorphisme $v\colon M'\rightarrow M$
tels que l'homomorphisme induit  $M'\rightarrow M/M^\times$ soit un isomorphisme.  
\end{itemize}
Alors~:
\begin{itemize}
\item[{\rm (i)}] Le monoïde $Q$ est intègre et le groupe $M'^\gp$ est libre. 
\item[{\rm (ii)}] On peut compléter le diagramme \eqref{higgs1-eip5a} en un diagramme commutatif 
\begin{equation}\label{higgs1-eip6a}
\xymatrix{
{M'}\ar[r]^v\ar[d]_-(0.4)w&{M}\ar[d]\\
{\underset{\underset{x\mapsto x^p}{\longleftarrow}}{\lim}\ A}\ar[r]&A}
\end{equation}
Notons $\beta\colon M'\rightarrow Q$ l'homomorphisme induit. 
\item[{\rm (iii)}] La structure logarithmique $\cQ$ sur $\Spec(\rW(\cR_A))$ est associée à la structure pré-logarithmique 
définie par l'homomorphisme composé
\begin{equation}\label{higgs1-eip6b}
M'\stackrel{\beta}{\rightarrow} Q \stackrel{q}{\rightarrow} \rW(\cR_A).
\end{equation}
En particulier, le schéma logarithmique $(\Spec(\rW(\cR_A)),\cQ)$ est fin et saturé.
\item[{\rm (iv)}] Si de plus, l'homomorphisme composé $u\circ v\colon M'\rightarrow \Gamma(Y,\cM_Y)$ est 
une carte pour $Y$ \eqref{higgs1-log7}, alors  le morphisme \eqref{higgs1-eip5d} est strict.
\end{itemize}
\end{prop}

(i) Comme $Y^\circ$ n'est pas vide, l'image canonique de $\Gamma(Y,\cM_Y)$ dans $A$ ne contient par $0$. 
On en déduit aussitôt que $Q$ est intègre. 
D'autre part, comme $M'$ est saturé, le sous-groupe de torsion de $M'^\gp$ est contenu dans $M'$. 
Mais $M'$ est affûté~; donc $M'^\gp$ est sans torsion.  

(ii) Soient $L$ le corps des fractions de $A$, $\oL$ une clôture algébrique de $L$. 
Considérons le système projectif de monoïdes multiplicatifs $(\oL)_{n\in \mN}$, 
où les morphismes de transition sont l'élévation à la puissance $p$-ième. 
D'après (i) et sa preuve, il existe un homomorphisme 
\begin{equation}\label{higgs1-eip6c}
M'^\gp \rightarrow \underset{\underset{x\mapsto x^p}{\longleftarrow}}{\lim}\ \oL
\end{equation}
qui relève l'homomorphisme $M'^\gp\rightarrow \oL$ induit par le composé 
$M'\stackrel{v}{\rightarrow} M\rightarrow A$. 

Soient $t\in \Gamma(Y,\cM_Y)$, $x$ son image canonique dans $A$, $y\in \oL$ tel que $y^p=x$.  
L'extension $A[z]/z^p-x$ de $A$ est étale au-dessus de $Y^\circ$ car $p$ est inversible dans $Y^\circ$ et que
$x$ ne s'annule en aucun point de $Y^\circ$.
Comme $Y^\circ$ est simplement connexe, on en déduit que $y\in L$ et par suite que $y\in A$ puisque $A$ est normal.
On en déduit que la restriction de \eqref{higgs1-eip6c} à $M'$ induit un homomorphisme 
\begin{equation}\label{higgs1-eip6d}
w\colon M'\rightarrow \underset{\underset{x\mapsto x^p}{\longleftarrow}}{\lim}\ A
\end{equation}
qui répond à la question.

(iii) Soit $G$ l'image inverse de $M^\times$ par  l'homomorphisme canonique $Q\rightarrow M$. 
Il résulte aussitôt de la définition \eqref{higgs1-eip5a} que $G$ est un sous-groupe de $Q$. 
D'autre part, l'homomorphisme composé $M'\rightarrow Q/G\rightarrow M/M^\times$,
où la première flèche est déduite de $\beta$, est un isomorphisme. 
Par suite, $M'\rightarrow Q/G$ est un isomorphisme. 
La proposition résulte  alors de (\cite{tsuji1} 1.3.1).  

(iv) Cela résulte aussitôt de (iii).

\subsection{}\label{higgs1-eip7}\index{10940@$i_\oS\colon (\coS,\cM_\coS)\rightarrow (\cA_2(\oS),\cM_{\cA_2(\oS)})$}
Dans la suite de cet article, on fixe une suite 
$(p_n)_{n\geq 0}$ d'éléments de $\co_\oK$ telle que $p_0=p$ et $p_{n+1}^p=p_n$ pour tout $n\geq 0$.  
On désigne par $\upp$ l'élément de $\cR_{\co_\oK}$ induit par la suite $(p_n)_{n\geq 0}$ et on pose  
\begin{equation}\label{higgs1-eip7a}
\xi=[\upp]-p \in \rW(\cR_{\co_\oK}),
\end{equation}
où $[\ ]$ est le représentant multiplicatif. D'après \ref{higgs1-eip4}, la suite 
\begin{equation}\label{higgs1-eip7b}
0\longrightarrow \rW(\cR_{\co_\oK})\stackrel{\cdot \xi}{\longrightarrow} \rW(\cR_{\co_\oK})
\stackrel{\theta}{\longrightarrow} \co_C \longrightarrow 0
\end{equation}
est exacte. Elle induit une suite exacte 
\begin{equation}\label{higgs1-eip7c}
0\longrightarrow \co_C\stackrel{\cdot \xi}{\longrightarrow} \cA_2(\co_\oK)
\stackrel{\theta}{\longrightarrow} \co_C \longrightarrow 0,
\end{equation}
où on a encore noté $\cdot \xi$ le morphisme induit par la multiplication par $\xi$ dans $\cA_2(\co_\oK)$ \eqref{higgs1-eipo3e}. 
L'idéal $\ker(\theta)$ de $\cA_2(\co_\oK)$ est de carré nul. 
C'est un $\co_C$-module libre de base $\xi$. Il sera noté $\xi\co_C$. 
On observera que contrairement à $\xi$, ce module ne dépend pas du choix de la suite $(p_n)_{n\geq 0}$. 

Le groupe de Galois $G_K$ agit naturellement sur $\rW(\cR_{\co_\oK})$ par des automorphismes d'anneaux,
et l'homomorphisme $\theta$ \eqref{higgs1-eip7b} est $G_K$-équivariant. On en déduit une action de $G_K$ 
sur $\cA_2(\co_{\oK})$ par des automorphismes d'anneaux tel que l'homomorphisme $\theta$ \eqref{higgs1-eip7c} soit 
$G_K$-équivariant.

On pose 
\begin{equation}\label{higgs1-eip7e}
\oS=\Spec(\co_\oK)\ \ \ {\rm et}\ \ \ \coS=\Spec(\co_C)
\end{equation} 
que l'on munit des structures logarithmiques images inverses de $\cM_S$ \eqref{higgs1-dlog0}, 
notées respectivement $\cM_\oS$ et $\cM_\coS$. 
Les actions de $G_K$ sur $\co_\oK$ et $\co_C$ s'étendent naturellement en des actions à gauche sur 
les schémas logarithmiques $(\oS,\cM_\oS)$ et $(\coS,\cM_\coS)$, respectivement. 

On pose 
\begin{equation}\label{higgs1-eip7f}
\cA_2(\oS)=\Spec(\cA_2(\co_\oK))
\end{equation} 
que l'on munit de la structure logarithmique $\cM_{\cA_2(\oS)}$
définie comme suit. Soient $Q_\oS$ le monoïde et $q_\oS\colon Q_\oS\rightarrow \rW(\cR_{\co_\oK})$
l'homomorphisme 
définis dans \ref{higgs1-eip5} (notés $Q$ et $q$) en prenant pour $(Y,\cM_Y)$ le schéma logarithmique $(\oS,\cM_\oS)$ 
et pour $u$ l'homomorphisme canonique  
$\Gamma(S,\cM_S)\rightarrow \Gamma(\oS,\cM_\oS)$. On désigne par $\cM_{\cA_2(\oS)}$ la structure logarithmique 
sur $\cA_2(\oS)$ associée à la structure pré-logarithmique définie par l'homomorphisme 
$Q_\oS\rightarrow \cA_2(\co_\oK)$ induit par $q_\oS$. D'après \ref{higgs1-eip6}, 
le schéma logarithmique $(\cA_2(\oS),\cM_{\cA_2(\oS)})$ est fin et saturé, et $\theta$ induit une immersion fermée exacte 
\begin{equation}\label{higgs1-eip7d}
i_\oS\colon (\coS,\cM_\coS)\rightarrow (\cA_2(\oS),\cM_{\cA_2(\oS)}).
\end{equation}

Le groupe de Galois $G_K$ agit naturellement sur le monoïde $Q_\oS$,
et l'homomorphisme $q_\oS\colon Q_\oS\rightarrow \rW(\cR_{\co_\oK})$ est $G_K$-équivariant. 
On en déduit une action à gauche de $G_K$ sur le schéma logarithmique $(\cA_2(\oS),\cM_{\cA_2(\oS)})$. 
Le morphisme $i_\oS$ est $G_K$-équivariant.

\begin{rema}\label{higgs1-eip75}\index{10945@$\xi^{-1} M$ ($M$ un $\co_C$-module)}
On note $\xi^{-1}\co_C$ le $\co_C$-module dual de $\xi\co_C$. 
Pour tout $\co_C$-module $M$, on désigne les $\co_C$-modules $M\otimes_{\co_C}(\xi \co_C)$ 
et $M\otimes_{\co_C}(\xi^{-1} \co_C)$ simplement par $\xi M$ et $\xi^{-1} M$, respectivement. 
On observera que contrairement à $\xi$, ces modules ne dépendent pas du choix de la suite $(p_n)_{n\geq 0}$. 
Il est donc important de ne pas les identifier à $M$. 

Pour des $\hRun$-algèbres $A$, nous considérerons dans la suite de cet article des $A$-modules 
de Higgs à coefficients dans $\xi^{-1}\tOmega^1_{R/\co_K}\otimes_RA$ (cf. \ref{higgs1-not6}). 
Nous dirons abusivement qu'ils sont à coefficients dans $\xi^{-1}\tOmega^1_{R/\co_K}$. 
La catégorie de ces modules sera notée $\bMH(A,\xi^{-1}\tOmega^1_{R/\co_K})$. 
\end{rema}

\begin{prop}\label{higgs1-eip8}
Les homomorphismes de Frobenius absolus de  $R_\infty/pR_\infty$ et $\oR/p\oR$ sont surjectifs.
\end{prop}

Montrons d'abord que l'homomorphisme de Frobenius absolu de $R_\infty/pR_\infty$ est surjectif. 
Rappelons que pour tout entier $n\geq 1$, $\Spec(R_{n})$ est une composante connexe de 
$X_{n}\otimes_{\co_{K_{n}}}\co_{\oK}$ \eqref{higgs1-gal1} et qu'on a 
un diagramme cartésien de $\co_{\oK}$-morphismes \eqref{higgs1-dlog4e}
\begin{equation}\label{higgs1-eip8a}
\xymatrix{
{X_n\otimes_{\co_{K_n}}\co_{\oK}}\ar[r]\ar[d]&{\Spec(\co_{\oK}[P^{(n)}]/(\pi_n-e^{\lambda^{(n)}}))}\ar[d]\\
{X\otimes_{\co_K}\co_{\oK}}\ar[r]&{\Spec(\co_{\oK}[P]/(\pi-e^\lambda))}}
\end{equation}
où pour tout $t\in P$, on a noté $t^{(n)}$ son image dans $P^{(n)}$ par l'isomorphisme canonique \eqref{higgs1-ext-log2d}. 
Comme le morphisme canonique $X\rightarrow S\times_{\bA_\mN}\bA_P$ est étale \eqref{higgs1-dlog1b},
il suffit de montrer que l'homomorphisme de Frobenius absolu de la limite inductive
de $\mF_p$-algèbres (relativement à la relation de divisibilité)
\begin{equation}\label{higgs1-eip8b}
\underset{\underset{n\geq 1}{\longrightarrow}}{\lim}\ (\co_{\oK}/p\co_\oK)[P^{(n)}]/(\pi_{n}-e^{\lambda^{(n)}})
\end{equation}
est surjectif. Ceci résulte du fait que l'homomorphisme de Frobenius absolu de $\co_\oK/p\co_\oK$
est surjectif et que pour tout $t\in P$, on a $e^{t^{(n)}}=e^{pt^{(pn)}}$ dans \eqref{higgs1-eip8b}. 

Montrons ensuite que l'homomorphisme de Frobenius absolu de $\oR/p\oR$ est surjectif.
Soient $N$ une extension finie de $F_\infty$ contenue dans $\oF$, $D$ la clôture intégrale de $R$ dans $N$. 
Supposons que $D$ soit presque-étale sur $R_\infty$. 
Alors $D/pD$ est presque-étale sur $R_\infty/pR_\infty$ (\cite{tsuji2} 7.4(3)). 
Comme l'homomorphisme de Frobenius absolu de $R_\infty/pR_\infty$ est surjectif, 
l'homomorphisme de Frobenius absolu de $D/pD$ est presque surjectif en vertu de (\cite{tsuji2} 7.9). 
Pour tout $x\in D$, il existe $x',y\in D$ tel que $p^{1/2}x=x'^{p}+py$. Alors 
$x''=p^{-1/(2p)}x'\in D$ et on a $x=x''^p+p^{1/2}y$. Par le même argument, il existe $y',z\in D$
tel que $y=y'^p+p^{1/2}z$. Par suite, on a $x\equiv (x''+p^{1/(2p)}y')^p\mod pD$. 
On en déduit par passage à la limite inductive, en vertu de \ref{higgs1-pur4},
que l'homomorphisme de Frobenius absolu de $\oR/p\oR$ est surjectif.

\subsection{}\label{higgs1-eip9}\index{10950@$(Y,\cM_Y)$, $(\hY,\cM_\hY)$}\index{10951@$\cA_2(Y)$}\index{10952@$q_Y\colon Q_Y\rightarrow \rW(\cR_{\oR})$}
D'après \ref{higgs1-eip4} et \ref{higgs1-eip8}, la suite 
\begin{equation}\label{higgs1-eip9a}
0\longrightarrow \rW(\cR_\oR)\stackrel{\cdot \xi}{\longrightarrow} \rW(\cR_\oR)
\stackrel{\theta}{\longrightarrow} \hoR \longrightarrow 0,
\end{equation}
où $\xi$ est l'élément \eqref{higgs1-eip7a}, est exacte. Elle induit une suite exacte 
\begin{equation}\label{higgs1-eip9b}
0\longrightarrow \hoR\stackrel{\cdot \xi}{\longrightarrow} \cA_2(\oR)
\stackrel{\theta}{\longrightarrow} \hoR \longrightarrow 0,
\end{equation}
où on a encore noté $\cdot \xi$ le morphisme induit par la multiplication par $\xi$ dans $\cA_2(\oR)$. 
L'idéal $\ker(\theta)$ de $\cA_2(\oR)$ est de carré nul. C'est un $\hoR$-module libre de base $\xi$, 
canoniquement isomorphe à $\xi\hoR$ (cf. \ref{higgs1-eip75}). 
Le groupe de Galois $\Gamma$ \eqref{higgs1-gal2} agit naturellement sur $\rW(\cR_\oR)$ par des automorphismes d'anneaux,
et l'homomorphisme $\theta$ \eqref{higgs1-eip9a} est $\Gamma$-équivariant. On en déduit une action de $\Gamma$ 
sur $\cA_2(\oR)$ par des automorphismes d'anneaux telle que l'homomorphisme $\theta$ \eqref{higgs1-eip9b} soit 
$\Gamma$-équivariant.

On pose 
\begin{equation}\label{higgs1-eip9e}
Y=\Spec(\oR)\ \ \ {\rm et} \ \ \ \hY=\Spec(\hoR)
\end{equation}
que l'on munit des structures logarithmiques images inverses de $\cM_X$ \eqref{higgs1-dlog1}, 
notées respectivement $\cM_Y$ et $\cM_\hY$. 
Les actions de $\Gamma$ sur $\oR$ et $\hoR$ induisent des actions à gauche sur 
les schémas logarithmiques $(Y,\cM_Y)$ et $(\hY,\cM_\hY)$, respectivement. 

On pose 
\begin{equation}\label{higgs1-eip9c}
\cA_2(Y)=\Spec(\cA_2(\oR))
\end{equation} 
que l'on munit de la structure logarithmique $\cM'_{\cA_2(Y)}$
définie comme suit. Soient $Q_Y$ le monoïde et $q_Y\colon Q_Y\rightarrow \rW(\cR_{\oR})$
l'homomorphisme définis dans \ref{higgs1-eip5} (notés $Q$ et $q$)
en prenant pour $u$ l'homomorphisme canonique $\Gamma(X,\cM_X)\rightarrow \Gamma(Y,\cM_Y)$. 
On désigne par $\cM'_{\cA_2(Y)}$ la structure logarithmique 
sur $\cA_2(Y)$ associée à la structure pré-logarithmique définie par l'homomorphisme 
$Q_Y\rightarrow \cA_2(\oR)$ induit par $q_Y$. L'homomorphisme $\theta$ induit alors un morphisme \eqref{higgs1-eip5d}
\begin{equation}\label{higgs1-eip9d}
i'_Y\colon (\hY,\cM_\hY)\rightarrow (\cA_2(Y),\cM'_{\cA_2(Y)}).
\end{equation}

Le groupe de Galois $\Gamma$ agit naturellement sur le monoïde $Q_Y$,
et l'homomorphisme $q_Y\colon Q_Y\rightarrow \rW(\cR_\oR)$ est $\Gamma$-équivariant. 
On en déduit une action à gauche de $\Gamma$ sur $(\cA_2(Y),\cM'_{\cA_2(Y)})$. Le morphisme $i'_Y$
est $\Gamma$-équivariant. 

Sans les hypothèses de \ref{higgs1-eip6}, on ne saurait pas si 
le schéma logarithmique $(\cA_2(Y),\cM'_{\cA_2(Y)})$ est fin et saturé et 
si $i'_Y$ est une immersion fermée exacte. C'est pourquoi nous équiperons $\cA_2(Y)$ 
d'une autre structure logarithmique $\cM_{\cA_2(Y)}$ dans \ref{higgs1-eplog5}. 

\subsection{}\label{higgs1-eplog5}\index{10955@$i_Y\colon (\hY,\cM_{\hY})\rightarrow (\cA_2(Y),\cM_{\cA_2(Y)})$}
Pour tout $t\in P$, on désigne par $\tlt$ l'élément  de  $Q_Y$ défini par ses projections \eqref{higgs1-eip5a} 
\begin{equation}\label{higgs1-eplog5a}
(\alpha_\infty(t^{(p^n)}))_{n\in \mN}\in \underset{\underset{\mN}{\longleftarrow}}{\lim}\ \oR \ \ \ {\rm et}\ \ \ 
\gamma(t) \in \Gamma(X,\cM_X),
\end{equation}
où $t^{(p^n)}$ est l'image de $t$ dans $P^{(p^n)}$ par l'isomorphisme \eqref{higgs1-ext-log2d}
et $\alpha_\infty$ est l'homomorphisme \eqref{higgs1-log-ext2c}. On verra que cette notation est compatible
avec celle introduite dans \eqref{higgs1-ext-log13c} et n'induit aucune confusion.  
Le lecteur observera ici que $\tlt$ dépend du choix du morphisme \eqref{higgs1-dlog5b}.
L'application 
\begin{equation}\label{higgs1-eplog5b}
P\rightarrow Q_Y, \ \ \ t\mapsto \tlt
\end{equation}
ainsi définie est un morphisme de monoïdes. On désigne par 
\begin{equation}\label{higgs1-eplog5bb}
\tq_Y\colon P\rightarrow  \cA_2(\oR)
\end{equation}
l'homomorphisme composé
\begin{equation}\label{higgs1-eplog5bc}
P\longrightarrow Q_Y\stackrel{q_Y}{\longrightarrow} \rW(\cR_\oR) \longrightarrow  \cA_2(\oR).
\end{equation}
On munit $\cA_2(Y)$ de la structure logarithmique $\cM_{\cA_2(Y)}$ 
associée à la structure pré-logarithmique définie par $\tq_Y$. 
Le morphisme \eqref{higgs1-eplog5b} induit alors un morphisme de structures logarithmiques sur $\cA_2(Y)$ \eqref{higgs1-eip9}
\begin{equation}\label{higgs1-eplog5c}
\cM_{\cA_2(Y)}\rightarrow \cM'_{\cA_2(Y)}.
\end{equation}
Il est clair que l'homomorphisme composé $\theta\circ \tq_Y\colon P\rightarrow\hoR$ 
est induit par $\alpha$ (cf. \ref{higgs1-log-ext2}). Par suite, $\theta$ induit une immersion fermée exacte
\begin{equation}\label{higgs1-eplog5d}
i_Y\colon (\hY,\cM_{\hY})\rightarrow (\cA_2(Y),\cM_{\cA_2(Y)}),
\end{equation}
qui se factorise à travers $i'_Y$ \eqref{higgs1-eip9d}. 

\begin{prop}\label{higgs1-eplog4}
Supposons qu'il existe une carte fine et saturée $M'\rightarrow \Gamma(X,\cM_X)$ 
pour $(X,\cM_X)$ \eqref{higgs1-log7} telle que l'homomorphisme induit
\begin{equation}\label{higgs1-eplog4a}
M'\rightarrow \Gamma(X,\cM_X)/\Gamma(X,\co_X^\times)
\end{equation}
soit un isomorphisme.  Alors le morphisme $\cM_{\cA_2(Y)}\rightarrow \cM'_{\cA_2(Y)}$
\eqref{higgs1-eplog5c} est un isomorphisme. 
\end{prop}

En effet,  d'après \ref{higgs1-eip6}, le schéma logarithmique $(\cA_2(Y),\cM'_{\cA_2(Y)})$ est fin et saturé, et 
le morphisme $i'_Y$ \eqref{higgs1-eip9d} est une immersion fermée exacte. Par suite, 
pour tout point géométrique $\oy$ de $\hY$, notant encore  $\oy$ le point géométrique $i_Y(\oy)=i'_Y(\oy)$
de $\cA_2(Y)$, l'homomorphisme 
\begin{equation}
\cM_{\cA_2(Y),\oy}/\co^\times_{\cA_2(Y),\oy} \rightarrow \cM'_{\cA_2(Y),\oy}/\co^\times_{\cA_2(Y),\oy}
\end{equation}
induit par \eqref{higgs1-eplog5c} est un isomorphisme. Comme $\cM'_{\cA_2(Y),\oy}$ est intègre, on en déduit que 
le morphisme $\cM_{\cA_2(Y),\oy}\rightarrow \cM'_{\cA_2(Y),\oy}$ \eqref{higgs1-eplog5c} est un isomorphisme; d'où la proposition.

\begin{rema} \label{higgs1-eplog8}
On notera que contrairement à $\cM'_{\cA_2(Y)}$, la structure logarithmique $\cM_{\cA_2(Y)}$
dépend de la carte $(P,\gamma)$ pour $(X,\cM_X)$ \eqref{higgs1-dlog1}.  
Toutefois,  d'après \ref{higgs1-log15}, quitte à remplacer $X$ par un recouvrement ouvert affine, 
on peut supposer la condition de \ref{higgs1-eplog4} remplie, auquel cas $\cM_{\cA_2(Y)}$
ne dépend plus de la carte $(P,\gamma)$.
\end{rema}

\begin{rema} \label{higgs1-eplog9}
On désigne par $\tpi$ l'élément de $Q_\oS$ défini par ses projections \eqref{higgs1-eip5a}
\begin{equation}\label{higgs1-eplog9a}
(\pi_{p^n})_{n\in \mN}\in \underset{\underset{\mN}{\longleftarrow}}{\lim}\ \co_\oK \ \ \ {\rm et}\ \ \ 
\pi \in \Gamma(S,\cM_S),
\end{equation}
(cf. \ref{higgs1-dlog2} et \ref{higgs1-eip7}) et par 
\begin{equation}\label{higgs1-eplog9b}
\tq_\oS\colon \mN \rightarrow \cA_2(\co_\oK)
\end{equation} 
l'homomorphisme composé 
\begin{equation}\label{higgs1-eplog9c}
\mN\longrightarrow Q_\oS \stackrel{q_\oS}{\longrightarrow}\rW(\cR_{\co_\oK}) \rightarrow \cA_2(\co_\oK),
\end{equation}
où la première flèche envoie $1$ sur $\tpi$. 
La structure logarithmique sur $\cA_2(\oS)$ associée à la structure 
pré-logarithmique définie par $\tq_\oS$  est canoniquement isomorphe à $\cM_{\cA_2(\oS)}$ \eqref{higgs1-eip7}. 
Cela résulte de \ref{higgs1-eip6} comme dans la preuve de \ref{higgs1-eplog4}. 
\end{rema} 

\subsection{}\label{higgs1-eplog25}
On a un homomorphisme canonique $Q_\oS\rightarrow Q_Y$ qui s'insère dans un diagramme commutatif 
\begin{equation}\label{higgs1-eplog25a}
\xymatrix{
{\mN}\ar[r]^{\vartheta}\ar[d]&P\ar[d]\\
Q_\oS\ar[r]\ar[d]_{q_\oS}&Q_Y\ar[d]^{q_Y}\\
{\rW(\cR_{\co_\oK})}\ar[r]&{\rW(\cR_\oR)}}
\end{equation}
où $\vartheta$ est l'homomorphisme donné dans \eqref{higgs1-dlog1}, 
la flèche verticale supérieure de gauche envoie $1$ sur $\tpi$ \eqref{higgs1-eplog9} et celle de droite est \eqref{higgs1-eplog5b}. 
On en déduit un diagramme commutatif 
\begin{equation}\label{higgs1-eplog25b}
\xymatrix{
{(\hY,\cM_\hY)}\ar[r]^-(0.5){i_Y}\ar[d]&{(\cA_2(Y),\cM_{\cA_2(Y)})}\ar[d]\\
{(\coS,\cM_\coS)}\ar[r]^-(0.5){i_S}&{(\cA_2(\oS),\cM_{\cA_2(\oS)})}}
\end{equation}

\subsection{}\label{higgs1-ext110}\index{10960@$\log([\ ])\colon \mZ_p(1)\rightarrow \cA_2(\oR)$}
On a un homomorphisme canonique 
\begin{equation}\label{higgs1-ext110a}
\mZ_p(1)\rightarrow \cR^\times_\oR.
\end{equation} 
Pour tout $\zeta\in \mZ_p(1)$, on note encore $\zeta$ son image dans $\cR^\times_\oR$.  
Comme $\theta([\zeta]-1)=0$, on obtient un homomorphisme de groupes
\begin{equation}\label{higgs1-ext11b}
\mZ_p(1)\rightarrow \cA_2(\oR),\ \ \ 
\zeta\mapsto\log([\zeta])=[\zeta]-1,
\end{equation}
dont l'image est contenue dans $\ker(\theta)=\xi\hoR$. 

\begin{lem}\label{higgs1-ext111}
L'homomorphisme $\mZ_p(1)\rightarrow \cA_2(\oR)$ \eqref{higgs1-ext11b} est $\mZ_p$-linéaire 
et $\Gamma$-équivariant~; son image engendre l'idéal $p^{\frac{1}{p-1}}\xi\hoR$ de $\cA_2(\oR)$
et le morphisme $\hoR$-linéaire 
\begin{equation}\label{higgs1-ext11c}
\hoR(1)\rightarrow p^{\frac{1}{p-1}}\xi \hoR
\end{equation}
qui à $x\otimes \zeta$, où $x\in \hoR$ et $\zeta\in \mZ_p(1)$, associe $x\cdot \log([\zeta])$ est un isomorphisme. 
\end{lem}

La première assertion est immédiate puisque $\mZ_p(1)$ et $\cA_2(\oR)$ sont complets et séparés pour les topologies
$p$-adiques. Soit $\zeta=(\zeta_n)_{n\geq 0}\in \mZ_p(1)$ tel que $\zeta_1\not=1$ (on a $\zeta_0=1$). 
Notons $\zeta'$ l'image canonique de $(\zeta_{n+1})_{n\geq 0}$ dans $\cR_\oR$
et posons $\omega=\sum_{i=0}^{p-1}[\zeta']^i\in \rW(\cR_\oR)$. Alors $\omega$ est un générateur 
de $\ker(\theta)$ (\cite{tsuji1} A.2.6), et on a dans $\cA_2(\oR)$ 
\begin{equation}
[\zeta]-1=\omega ([\zeta']-1)=\omega \theta([\zeta']-1)=\omega(\zeta_1-1).
\end{equation}
Comme $v(\zeta_1-1)=\frac{1}{p-1}$, le reste de l'assertion s'ensuit.

\subsection{}\label{higgs1-eip10}
L'homomorphisme canonique $\mZ_p(1)\rightarrow (\cR_\oR,\times)$ \eqref{higgs1-ext110a} et l'homomorphisme
trivial $\mZ_p(1)\rightarrow \Gamma(X,\cM_X)$  (de valeur $1$) induisent un homomorphisme
\begin{equation}\label{higgs1-eip10aa}
\mZ_p(1)\rightarrow Q_Y.
\end{equation} 
Pour tous $g\in \Gamma$ et $t\in P$, on a dans $Q_Y$
\begin{equation}\label{higgs1-eip10ab}
g(\tlt)=\tchi_t(g)\cdot \tlt,
\end{equation}
où l'on a (abusivement) noté $\tchi_t\colon \Gamma\rightarrow \mZ_p(1)$ l'application déduite de \eqref{higgs1-ext-log14a}.
On en déduit la relation suivante dans $\cA_2(\oR)$
\begin{equation}\label{higgs1-eip10a}
g(\tq_Y(t))=[\tchi_t(g)]\cdot \tq_Y(t),
\end{equation}
où $[\tchi_t(g)]$ désigne l'image de $\tchi_t(g)$ par l'application composée 
\[
\mZ_p(1)\longrightarrow \cR_\oR\stackrel{[\ ]}{\longrightarrow}\rW(\cR_\oR)\longrightarrow \cA_2(\oR).
\]

Notons $\tau_g$ l'automorphisme de $\cA_2(Y)$ induit par l'action de $g$ sur $\cA_2(\oR)$. 
La structure logarithmique $\tau_g^*(\cM_{\cA_2(Y)})$ sur $\cA_2(Y)$ 
est associée à la structure pré-logarithmique définie par l'homomorphisme
composé $g\circ \tq_Y\colon P\rightarrow \cA_2(\oR)$ \eqref{higgs1-eplog5bb}. L'application
\begin{equation}\label{higgs1-eip10b}
P\rightarrow \Gamma(\cA_2(Y),\cM_{\cA_2(Y)}), \ \ \ t\mapsto [\tchi_t(g)]\cdot t
\end{equation}
est un morphisme de monoïdes \eqref{higgs1-ext-log14d}. 
Elle induit donc un morphisme de structures logarithmiques sur $\cA_2(Y)$
\begin{equation}\label{higgs1-eip10c}
a_g\colon \tau_g^*(\cM_{\cA_2(Y)})\rightarrow \cM_{\cA_2(Y)}.
\end{equation}
De même, en vertu de \eqref{higgs1-ext-log14c}, le morphisme de monoïdes  
\begin{equation}\label{higgs1-eip10d}
P\rightarrow \Gamma(\cA_2(Y),\tau_g^*(\cM_{\cA_2(Y)})), \ \ \ t\mapsto [g(\tchi_t(g^{-1}))] \cdot t
\end{equation}
induit un morphisme de structures logarithmiques sur $\cA_2(Y)$
\begin{equation}\label{higgs1-eip10e}
b_g\colon \cM_{\cA_2(Y)}\rightarrow \tau_g^*(\cM_{\cA_2(Y)}).
\end{equation}
On voit aussitôt que $a_g$ et $b_g$ sont des isomorphismes inverses l'un de l'autre \eqref{higgs1-ext-log14c}, et que 
l'application $g\mapsto  (\tau_{g^{-1}},a_{g^{-1}})$ est une action à gauche de $\Gamma$ sur le schéma logarithmique
$(\cA_2(Y),\cM_{\cA_2(Y)})$. 

On vérifie aussitôt que le morphisme $i_Y$ \eqref{higgs1-eplog5d} et 
le morphisme canonique \eqref{higgs1-eplog25b}
\begin{equation}\label{higgs1-eip10ac}
(\cA_2(Y),\cM_{\cA_2(Y)})\rightarrow (\cA_2(\oS),\cM_{\cA_2(\oS)})
\end{equation}
sont $\Gamma$-équivariants. Par ailleurs, pour tout $g\in \Gamma$, le diagramme 
\begin{equation}\label{higgs1-eip10ad}
\xymatrix{
{\tau_g^*(\cM_{\cA_2(Y)})}\ar[r]\ar[d]_{a_g}&{\tau_g^*(\cM'_{\cA_2(Y)})}\ar[d]^{a'_g}\\
{\cM_{\cA_2(Y)}}\ar[r]&{\cM'_{\cA_2(Y)}}}
\end{equation}
où les flèches horizontales sont induites par l'homomorphisme \eqref{higgs1-eplog5c} et 
$a'_g$ est l'automorphisme de structures logarithmiques sur $\cA_2(Y)$ induit par l'action de $g$ sur $Q_Y$, 
est commutatif. 

\section{Torseurs et algèbres de Higgs-Tate}\label{higgs1-TOR}

Les notations introduites dans \ref{higgs1-eip7}, \ref{higgs1-eip9} et \ref{higgs1-eplog5} sont en vigueur dans la suite de cet article. 

\subsection{}\label{higgs1-deflog1}\index{101001@$(\oX,\cM_{\oX})$, $(\coX,\cM_{\coX})$}\index{101002@$(\tX,\cM_\tX)$}
\index{deformation@$(\cA_2(\oS),\cM_{\cA_2(\oS)})$-déformation}
On pose $\oX=X\times_S\oS$ et $\coX=X\times_S\coS$ que l'on munit des structures logarithmiques 
images inverses de $\cM_X$ \eqref{higgs1-dlog1}, notées respectivement $\cM_\oX$ et $\cM_\coX$. 
On a alors des isomorphismes canoniques
\begin{eqnarray}
(\oX,\cM_{\oX})&\stackrel{\sim}{\rightarrow}&(X,\cM_X)\times_{(S,\cM_S)}(\oS,\cM_\oS),\label{higgs1-deflog1a}\\
(\coX,\cM_{\coX})&\stackrel{\sim}{\rightarrow}&(X,\cM_X)\times_{(S,\cM_S)}(\coS,\cM_\coS),\label{higgs1-deflog1b}
\end{eqnarray} 
le produit étant indifféremment pris dans la catégorie des schémas logarithmiques ou 
dans celle des schémas logarithmiques fins. On en déduit en particulier une action à gauche 
de $G_K$ sur les schémas logarithmiques $(\oX,\cM_{\oX})$ et $(\coX,\cM_{\coX})$. On a un morphisme 
canonique $\Gamma$-équivariant \eqref{higgs1-eip9e}
\begin{equation}\label{higgs1-deflog1c}
(\hY,\cM_{\hY})\rightarrow (\coX,\cM_{\coX}).
\end{equation}

Une {\em $(\cA_2(\oS),\cM_{\cA_2(\oS)})$-déformation lisse} de $(\coX,\cM_{\coX})$
est la donnée d'un morphisme lisse de schémas logarithmiques fins $(\tX,\cM_\tX)\rightarrow (\cA_2(\oS), \cM_{\cA_2(\oS)})$
et d'un $(\coS,\cM_{\coS})$-isomorphisme 
\begin{equation}\label{higgs1-deflog1d}
(\coX,\cM_{\coX})\stackrel{\sim}{\rightarrow}
(\tX,\cM_\tX)\times_{(\cA_2(\oS), \cM_{\cA_2(\oS)})}(\coS,\cM_{\coS}).
\end{equation}
Comme $(\coX,\cM_{\coX})\rightarrow (\coS,\cM_\coS)$ est lisse et que $\coX$ est affine,
une telle déformation existe et est unique à isomorphisme près en vertu de (\cite{kato1}, 3.14). 
Son groupe d'automorphismes est isomorphe à
\begin{equation}\label{higgs1-deflog1e}
\Hom_{\co_{\coX}}(\tOmega^1_{R/\co_K}\otimes_{R}\co_{\coX},\xi \co_{\coX}).
\end{equation}
On rappelle qu'on a posé $\tOmega^1_{R/\co_K}=\Omega^1_{(R,P)/(\co_K,\mN)}$ \eqref{higgs1-ext-log12c}. 

Dans la suite de cette section, nous fixons une telle déformation $(\tX,\cM_\tX)$.

\begin{rema}\label{higgs1-eplog79}
D'après \ref{higgs1-gal1}(i), le schéma $\Spec(R_1)$ est une composante connexe ouverte de $\oX$. 
Par suite, $Z=\Spec(R_1\otimes_{\co_\oK}\co_C)$ est un ouvert et fermé de $\coX$, par lequel se factorise 
le morphisme \eqref{higgs1-deflog1c}. Une $(\cA_2(\oS),\cM_{\cA_2(\oS)})$-déformation de $(Z,\cM_\coX|Z)$
suffit pour les besoins de cet article.   
\end{rema}

\subsection{}\label{higgs1-tor2}\index{101005@$\rT$, $\trT$, $\bT$} \index{101006@$\cS=\oplus_{n\geq 0}\cS^{n}$}
\index{101008@$\cL$ (torseur de Higgs-Tate)}\index{101009@$\cF$ (module des fonctions affines sur $\cL$)}\index{101011@$\cC$ (algèbre de Higgs-Tate)}
\index{101012@$\bL$}\index{Higgs-Tate!1@Torseur de --- ($\cL$)}
On pose 
\begin{equation}\label{higgs1-tor2ab}
\rT=\Hom_{\hoR}(\tOmega^1_{R/\co_K}\otimes_R\hoR,\xi\hoR).
\end{equation} 
On identifie le $\hoR$-module dual à $\xi^{-1}\tOmega^1_{R/\co_K}\otimes_R\hoR$ \eqref{higgs1-eip75}
et on note $\cS$ la $\hoR$-algèbre symétrique associée \eqref{higgs1-not621}
\begin{equation}\label{higgs1-tor2d}
\cS=\oplus_{n\geq 0}\cS^{n}=\rS_{\hoR}(\xi^{-1}\tOmega^1_{R/\co_K}\otimes_R\hoR).
\end{equation}
On désigne par $\hY_\zar$ le topos de Zariski de $\hY=\Spec(\hoR)$ \eqref{higgs1-eip9e}, par $\trT$ le $\co_\hY$-module associé à $\rT$
et par $\bT$ le $\hY$-fibré vectoriel associé à son dual, autrement dit,  
\begin{equation}\label{higgs1-tor2c}
\bT=\Spec(\cS).
\end{equation}
 
Soient $U$ un ouvert de Zariski de $\hY$, $\tU$ l'ouvert correspondant de $\cA_2(Y)$ (cf. \ref{higgs1-eip9}). 
On désigne par $\cL(U)$ 
l'ensemble des morphismes représentés par des flèches pointillées qui complètent  le diagramme 
\begin{equation}\label{higgs1-tor2a}
\xymatrix{
{(U,\cM_\hY|U)}\ar[r]^-(0.5){i_Y|U}\ar[d]&{(\tU,\cM_{\cA_2(Y)}|\tU)}\ar@{.>}[d]\ar@/^2pc/[dd]\\
{(\coX,\cM_{\coX})}\ar[r]\ar[d]&{(\tX,\cM_\tX)}\ar[d]\\
{(\coS,\cM_\coS)}\ar[r]^-(0.5){i_\oS}&{(\cA_2(\oS),\cM_{\cA_2(\oS)})}}
\end{equation}
de façon à le laisser commutatif. D'après \ref{higgs1-log10},
le foncteur $U\mapsto \cL(U)$ est un $\trT$-torseur de $\hY_\zar$.  
On l'appelle {\em torseur de Higgs-Tate} associé à $(\tX,\cM_\tX)$.  
On désigne par $\cF$ le $\hoR$-module des fonctions affines sur $\cL$ (cf. \ref{higgs1-eph13}). 
Celui-ci s'insère dans une suite exacte canonique \eqref{higgs1-eph13b}
\begin{equation}\label{higgs1-tor2e}
0\rightarrow \hoR\rightarrow \cF\rightarrow \xi^{-1}\tOmega^1_{R/\co_K} \otimes_R \hoR\rightarrow 0.
\end{equation} 
D'après (\cite{illusie1} I 4.3.1.7), cette suite induit pour tout entier $n\geq 1$, une suite exacte  \eqref{higgs1-not621}
\begin{equation}\label{higgs1-tor2f}
0\rightarrow \rS^{n-1}_{\hoR}(\cF)\rightarrow \rS^{n}_{\hoR}(\cF)\rightarrow \rS^n_{\hoR}(\xi^{-1}\tOmega^1_{R/\co_K}
\otimes_R\hoR)\rightarrow 0.
\end{equation}
Les $\hoR$-modules $(\rS^{n}_{\hoR}(\cF))_{n\in \mN}$ forment donc un système inductif filtrant, 
dont la limite inductive 
\begin{equation}\label{higgs1-tor2g}
\cC=\underset{\underset{n\geq 0}{\longrightarrow}}\lim\ \rS^n_{\hoR}(\cF)
\end{equation}
est naturellement munie d'une structure de $\hoR$-algèbre. D'après \ref{higgs1-eph12}, le $\hY$-schéma 
\begin{equation}\label{higgs1-tor2b}
\bL=\Spec(\cC)
\end{equation}
est naturellement un $\bT$-fibré principal homogène sur $\hY$ qui représente canoniquement $\cL$. 
On prendra garde que $\cL$, $\cF$, $\cC$ et $\bL$ dépendent de $(\tX,\cM_\tX)$.

\subsection{}\label{higgs1-tor200}\index{101016@${^g\psi}$}
On munit $\hY$ de l'action naturelle à gauche de $\Delta$~; pour tout $g\in \Delta$, 
l'automorphisme de $\hY$ défini par $g$, que l'on note aussi $g$, est induit par l'automorphisme $g^{-1}$ de $\hoR$. 
On considère $\trT$ comme un $\co_\hY$-module $\Delta$-équivariant au moyen 
de la donnée de descente correspondant au $\hRun$-module 
$\Hom_{\hRun}(\tOmega^1_{R/\co_K}\otimes_R\hRun,\xi\hRun)$ (cf. \ref{higgs1-eph5}). Pour tout $g\in \Delta$, on a donc 
un isomorphisme canonique de $\co_\hY$-modules
\begin{equation}\label{higgs1-tor200aa}
\tau_g^\trT\colon \trT\stackrel{\sim}{\rightarrow} g^*(\trT).
\end{equation}
Celui-ci induit un isomorphisme de $\hY$-schémas en groupes
\begin{equation}\label{higgs1-tor200a}
\tau_g^\bT\colon \bT\stackrel{\sim}{\rightarrow} g^\bullet(\bT),
\end{equation}
où $g^\bullet$ désigne le foncteur de changement de base par l'automorphisme $g$ de $\hY$ \eqref{higgs1-eph1b}.
On obtient ainsi une structure $\Delta$-équivariante sur le $\hY$-schéma en groupes $\bT$ (cf. \ref{higgs1-eph7})
et par suite une action à gauche de $\Delta$ sur $\bT$ compatible avec son action sur $\hY$; 
l'automorphisme de $\bT$ défini par un élément $g$ de $\Delta$ est le composé de $\tau^\bT_g$ 
et de la projection canonique $g^\bullet(\bT)\rightarrow \bT$. 
On en déduit une action de $\Delta$ sur $\cS$ par des automorphismes 
d'anneaux, compatible avec son action sur $\hoR$, que l'on appelle {\em action canonique}.
Cette dernière est concrètement  induite par l'action triviale sur 
$\rS_{\hRun}(\xi^{-1}\tOmega^1_{R/\co_K}\otimes_R\hRun)$.

L'action à gauche de $\Delta$ sur le schéma logarithmique $(\cA_2(Y),\cM_{\cA_2(Y)})$ définie dans \ref{higgs1-eip10}
induit sur le $\trT$-torseur $\cL$ une structure $\Delta$-équivariante (cf. \ref{higgs1-eph5}), 
autrement dit, elle induit pour tout $g\in \Delta$, un isomorphisme de $\tau_g^\trT$-équivariant
\begin{equation}\label{higgs1-tor200b}
\tau^{\cL}_g\colon \cL\stackrel{\sim}{\rightarrow} g^*(\cL);
\end{equation}
ces isomorphismes étant soumis à des relations de compatibilité \eqref{higgs1-eph6e}. 
En effet, pour tout ouvert de Zariski $U$ de $\hY$, on prend pour 
\begin{equation}\label{higgs1-tor200c}
\tau^{\cL}_g(U)\colon \cL(U)\stackrel{\sim}{\rightarrow} \cL(g(U))
\end{equation}
l'isomorphisme défini de la façon suivante. 
Soient $\tU$ l'ouvert de $\cA_2(Y)$ correspondant à $U$, $\psi\in \cL(U)$ que l'on considère comme un morphisme
\begin{equation}\label{higgs1-tor200d}
\psi\colon (\tU,\cM_{\cA_2(Y)}|\tU)\rightarrow (\tX,\cM_\tX).
\end{equation}
Comme les morphismes $i_Y$ \eqref{higgs1-eplog5d} et 
$(\hY,\cM_{\hY})\rightarrow (\coX,\cM_{\coX})$ \eqref{higgs1-deflog1c} sont $\Delta$-équivariants, le morphisme composé
\begin{equation}\label{higgs1-tor200e}
(g(\tU),\cM_{\cA_2(Y)}|g(\tU))\stackrel{g^{-1}}{\longrightarrow} (\tU,\cM_{\cA_2(Y)}|\tU)
\stackrel{\psi}{\longrightarrow} (\tX,\cM_\tX)
\end{equation}
prolonge le morphisme canonique $(g(U),\cM_\hY|g(U))\rightarrow (\tX,\cM_\tX)$. 
Il correspond à l'image de $\psi$ par $\tau^{\cL}_g(U)$. On vérifie aussitôt que 
le morphisme $\tau^{\cL}_g$ ainsi défini est un isomorphisme $\tau_g^\trT$-équivariant et que 
ces isomorphismes vérifient les relations de compatibilité requises \eqref{higgs1-eph6e}.

D'après \ref{higgs1-eph17}, les structures $\Delta$-équivariantes sur $\trT$ et $\cL$ induisent une structure
$\Delta$-équivariante sur le $\co_\hY$-module associé à $\cF$, ou, ce qui revient au même,
une action $\hoR$-semi-linéaire de $\Delta$ sur $\cF$, telle que les morphismes de la suite \eqref{higgs1-tor2e} soient 
$\Delta$-équivariants. 
On en déduit sur $\bL$ une structure de $\bT$-fibré principal homogène $\Delta$-équivariant 
sur $\hY$ (cf. \ref{higgs1-eph4}). Pour tout $g\in \Delta$, on a donc un isomorphisme $\tau_g^\bT$-équivariant
\begin{equation}\label{higgs1-tor200f}
\tau_g^{\bL}\colon \bL\stackrel{\sim}{\rightarrow} g^\bullet(\bL).
\end{equation}
Cette structure détermine une action à gauche de $\Delta$ sur $\bL$ compatible avec son action sur $\hY$; 
l'automorphisme de $\bL$ défini par un élément $g$ de $\Delta$ est le composé de $\tau^{\bL}_g$ 
et de la projection canonique $g^\bullet(\bL)\rightarrow \bL$. On obtient ainsi une action de $\Delta$ sur 
$\cC$ par des automorphismes d'anneaux, compatible avec son action sur $\hoR$, 
que l'on appelle {\em action canonique}. Cette dernière est concrètement induite par l'action de $\Delta$ sur $\cF$.

Pour tout $g\in \Delta$, on désigne par 
\begin{equation}\label{higgs1-tor200g}
\bL(\hY)\stackrel{\sim}{\rightarrow}\bL(\hY), \ \ \ \psi\mapsto {^g\psi} 
\end{equation}
le composé des isomorphismes 
\begin{eqnarray}
\tau_g^{\bL}\colon \bL(\hY)&\stackrel{\sim}{\rightarrow}& g^\bullet(\bL)(\hY),\\
g^\bullet(\bL)(\hY)&\stackrel{\sim}{\rightarrow}&\bL(\hY), \ \ \ \psi\mapsto \pr\circ \psi\circ g^{-1}, 
\end{eqnarray}
où $g^{-1}$ agit sur $\hY$ et $\pr\colon g^\bullet(\bL)\rightarrow \bL$ est la projection canonique,
de sorte que le diagramme 
\begin{equation}\label{higgs1-tor200h}
\xymatrix{
{\bL}\ar[r]^g&{\bL}\\
\hY\ar[r]^g\ar[u]^{\psi}&\hY\ar[u]_{{^g\psi}}}
\end{equation}
est commutatif. 

\begin{defi}\label{higgs1-tor201}\index{Higgs-Tate!3@Algèbre de --- ($\cC$)} \index{Higgs-Tate!2@Extension de --- ($\cF$)}
La $\hoR$-algèbre $\cC$ \eqref{higgs1-tor2g}, munie de l'action canonique de $\Delta$ 
\eqref{higgs1-tor200}, est appelée {\em l'algèbre de Higgs-Tate} associée à $(\tX,\cM_\tX)$. 
La $\hoR$-représentation $\cF$ \eqref{higgs1-tor2e} de $\Delta$ est appelée 
l'{\em extension de Higgs-Tate} associée à $(\tX,\cM_\tX)$.
\end{defi}
On notera que sous les hypothèses de \ref{higgs1-eplog4}, ces deux $\hoR$-représentations de $\Delta$ 
ne dépend pas du choix de la carte $(P,\gamma)$ \eqref{higgs1-dlog1}, d'après  \ref{higgs1-eplog4} et \eqref{higgs1-eip10ad}.

\subsection{}\label{higgs1-tor202}\index{101020@$\varphi_\psi$}
Pour tout $u\in \rT=\bT(\hY)$, on note 
\begin{equation}\label{higgs1-tor202a}
\ttt_u\colon \bT\rightarrow \bT 
\end{equation}
la translation par $u$. Pour tout $\psi\in \cL(\hY)$, on désigne par $\ttt_{\psi}$
l'isomorphisme de $\bT$-fibrés principaux homogènes sur $\hY$
\begin{equation}\label{higgs1-tor202b}
\ttt_{\psi}\colon \bT\stackrel{\sim}{\rightarrow} \bL,\ \ \ v\mapsto v+\psi.
\end{equation}
La structure de $\bT$-fibré principal homogène $\Delta$-équivariant sur $\bL$ se transporte 
par $\ttt_\psi$ en une structure de $\bT$-fibré principal homogène $\Delta$-équivariant sur $\bT$.
Pour tout $g\in \Delta$, on a donc un isomorphisme $\tau_g^\bT$-équivariant
\begin{equation}\label{higgs1-tor202c}
\tau_{g,\psi}^{\bT}\colon \bT\stackrel{\sim}{\rightarrow} g^\bullet(\bT).
\end{equation}
Cette structure détermine une action à gauche de $\Delta$ sur $\bT$ compatible avec son action sur $\hY$; 
l'automorphisme de $\bT$ défini par un élément $g$ de $\Delta$ est le composé de $\tau^\bT_{g,\psi}$ 
et de la projection canonique $g^\bullet(\bT)\rightarrow \bT$. On en déduit une action 
\begin{equation}\label{higgs1-tor202g}
\varphi_\psi\colon \Delta\rightarrow \Aut_\hRun(\cS)
\end{equation}
de $\Delta$ sur $\cS$ \eqref{higgs1-tor2d} par des automorphismes d'anneaux, 
compatible avec son action sur $\hoR$~;  pour tout $g\in \Delta$, 
$\varphi_\psi(g)$ est induit par l'automorphisme de $\bT$ défini par $g^{-1}$. 

L'isomorphisme $\ttt_\psi^*\colon \cC\stackrel{\sim}{\rightarrow} \cS$ induit par $\ttt_\psi$ est $\Delta$-équivariant 
si l'on munit $\cC$ de l'action canonique de $\Delta$ et $\cS$ de l'action $\varphi_\psi$.

On vérifie aussitôt que pour tout $g\in \Delta$, le diagramme 
\begin{equation}\label{higgs1-tor202d}
\xymatrix{
\bT\ar[r]^-(0.5){\ttt_\psi}\ar[d]_{\tau_g^\bT}&{\bL}\ar[d]^{\tau_g^{\bL}}\\
{g^\bullet(\bT)}\ar[r]^-(0.5){g^\bullet(\ttt_{{^g\psi}})}&{g^\bullet(\bL)}}
\end{equation}
est commutatif \eqref{higgs1-tor200g}. Par suite, le diagramme 
\begin{equation}\label{higgs1-tor202e}
\xymatrix{
\bT\ar[r]^-(0.5){\tau_{g}^\bT}\ar[dr]_-(0.5){\tau_{g,\psi}^\bT}&{g^\bullet(\bT)}\ar[d]^{g^\bullet(\ttt_{({^g\psi}-\psi)})}\\
&{g^\bullet(\bT)}}
\end{equation}
est aussi commutatif. On en déduit que
\begin{equation}\label{higgs1-tor202f}
\varphi_\psi(g)=\ttt^*_{(\psi-{^g \psi})}\circ g,
\end{equation}
où $\ttt^*_{(\psi-{^g \psi})}$ est l'automorphisme de $\cS$ induit par $\ttt_{(\psi-{^g \psi})}$
et $g$ agit sur $\cS$ par l'action canonique.

\subsection{}\label{higgs1-tor21}\index{101021@$D_u$, $\exp(D_u)$}
L'accouplement
$\rT\otimes_{\hoR} (\xi^{-1}\tOmega^1_{R/\co_K}\otimes_R\hoR)\rightarrow \hoR$ s'étend en un 
accouplement $\rT\otimes_{\hoR} \cS\rightarrow \cS$, où les éléments de $\rT$ agissent 
comme des $\hoR$-dérivations de $\cS$. On définit ainsi un morphisme 
\begin{equation}\label{higgs1-tor21a}
\rT\rightarrow \Gamma(\bT,\mT_{\bT/\hY}), \ \ \ u\mapsto D_u,
\end{equation}
où $\mT_{\bT/\hY}$ est le fibré tangent de $\bT$ sur $\hY$. Celui-ci identifie $\rT$
au module des champs de vecteurs de $\bT$ sur $\hY$ invariants par translation. 
Il induit aussi un isomorphisme 
\begin{equation}\label{higgs1-tor21b}
\rT\otimes_{\hoR}\co_{\bT}\stackrel{\sim}{\rightarrow}\mT_{\bT/\hY}.
\end{equation}
Notons $\Pi=\oplus_{n\geq 0}\Pi_{n}$ l'algèbre à puissances divisées de $\rT$ \eqref{higgs1-not621}.  
On a un accouplement canonique 
\begin{equation}\label{higgs1-tor21c}
\Pi_n\otimes_{\hoR}\cS^{n+m}\rightarrow \cS^{m},
\end{equation}
qui est parfait si $m=0$ (\cite{bo} A.10). Pour tout $u\in \rT=\bT(\hY)$, 
$D_u$ appartient à l'idéal à puissances divisées de $\Pi$, et on peut définir $\exp(D_u)$ 
comme un opérateur différentiel d'ordre infini de $\cS$. Pour tout 
$x\in \cS$, on a la formule de Taylor 
\begin{equation}\label{higgs1-tor21d}
\ttt_u^*(x)=\exp(D_u)(x),
\end{equation}
où $\ttt_u^*$ est l'automorphisme de $\cS$ induit par $\ttt_u$ \eqref{higgs1-tor202a}.

\subsection{}\label{higgs1-tor22}
On a un isomorphisme $\cS$-linéaire canonique 
\begin{equation}\label{higgs1-tor22a}
\Omega^1_{\cS/\hoR}\stackrel{\sim}{\rightarrow} \xi^{-1}\tOmega^1_{R/\co_K} \otimes_R\cS.
\end{equation}
On note
\begin{equation}\label{higgs1-tor22b}
d_{\cS}\colon \cS\rightarrow \xi^{-1}\tOmega^1_{R/\co_K} \otimes_R\cS
\end{equation}
la $\hoR$-dérivation universelle de $\cS$.
Le morphisme \eqref{higgs1-tor21a} est défini explicitement de la façon suivante~: pour tout $u\in \rT=\Hom_{\hoR}(\tOmega^1_{R/\co_K}\otimes_R\hoR,\xi\hoR)$, $D_u$ est la $\hoR$-dérivation composée
\begin{equation}\label{higgs1-tor22c}
\xymatrix{
\cS \ar[r]^-(0.5){d_{\cS}}& {\xi^{-1}\tOmega^1_{R/\co_K}\otimes_R \cS} \ar[rr]^-(0.5){(\xi^{-1}\cdot u)\otimes \id}&& \cS}.
\end{equation}
On peut alors redémontrer directement que l'opérateur différentiel 
\begin{equation}\label{higgs1-tor22d}
\exp(D_u)\colon \cS\rightarrow \cS
\end{equation}
est bien défini et que c'est un isomorphisme de $\hoR$-algèbres \eqref{higgs1-tor21d}. En effet, pour tout $n\geq 0$,
$D_u$ envoie $\cS^n$ sur $\cS^{n-1}$ \eqref{higgs1-tor2d}; il est donc nilpotent sur $\oplus_{0\leq i\leq n}\cS^i$.
Par suite, $\exp(D_u)$ est bien défini comme automorphisme de la $\hoR[\frac 1 p]$-algèbre $\cS[\frac 1 p]$.
D'autre part, pour tout $x\in \cS^1\subset \cS$, on a  
\begin{equation}\label{higgs1-tor22e}
\exp(D_u)(x)=x+(\xi^{-1}u)(x).
\end{equation}
On en déduit que $\exp(D_u)(\cS)\subset \cS$ et par suite que $\exp(D_u)(\cS)=\cS$ 
(puisque $D_{-u}=-D_u$).

\subsection{}\label{higgs1-tor23}
L'isomorphisme \eqref{higgs1-tor22a} induit un isomorphisme 
\begin{equation}\label{higgs1-tor23a}
\Omega^1_{\cC/\hoR}\stackrel{\sim}{\rightarrow} \xi^{-1}\tOmega^1_{R/\co_K} \otimes_R\cC.
\end{equation}
On note
\begin{equation}\label{higgs1-tor23b}
d_{\cC}\colon \cC\rightarrow \xi^{-1}\tOmega^1_{R/\co_K} \otimes_R\cC,
\end{equation}
la $\hoR$-dérivation universelle de $\cC$. Pour tout $x\in \cF$, $d_{\cC}(x)$ est l'image canonique 
de $x$ dans $\xi^{-1}\tOmega^1_{R/\co_K} \otimes_R\hoR$ \eqref{higgs1-tor2e}. 
Par suite, pour tout $g\in \Delta$, le diagramme 
\begin{equation}\label{higgs1-tor23d}
\xymatrix{
{\cC}\ar[rr]^{g}\ar[d]_{d_{\cC}}&&{\cC}\ar[d]^{d_{\cC}}\\
{\xi^{-1}\tOmega^1_{R/\co_K}\otimes_R\cC}\ar[rr]^{\id \otimes g}&&
{\xi^{-1}\tOmega^1_{R/\co_K}\otimes_R\cC}}
\end{equation}
est commutatif. Pour tout $\psi\in \cL(\hY)$, le diagramme 
\begin{equation}\label{higgs1-tor23c}
\xymatrix{
{\cC}\ar[d]_{d_{\cC}}\ar[rr]^-(0.5){\ttt^*_{\psi}}&&{\cS}\ar[d]^{d_{\cS}}\\
{\xi^{-1}\tOmega^1_{R/\co_K}\otimes_{R}\cC}
\ar[rr]^-(0.5){\id\otimes \ttt_{\psi}^*  }&&{\xi^{-1}\tOmega^1_{R/\co_K}\otimes_{R}\cS}}
\end{equation}
où $\ttt_\psi^*$ est l'isomorphisme induit par $\ttt_\psi$ \eqref{higgs1-tor202b}, est clairement commutatif.
On en déduit que le diagramme 
\begin{equation}\label{higgs1-tor13a}
\xymatrix{
{\cS}\ar[rr]^{\varphi_{\psi}(g)}\ar[d]_{d_{\cS}}&&{\cS}\ar[d]^{d_{\cS}}\\
{\xi^{-1}\tOmega^1_{R/\co_K}\otimes_R\cS}\ar[rr]^{\id \otimes\varphi_{\psi}(g)}&&
{\xi^{-1}\tOmega^1_{R/\co_K}\otimes_R\cS}}
\end{equation}
est commutatif, où $\varphi_\psi$ est l'action de $\Delta$ sur $\cS$ définie dans \eqref{higgs1-tor202g}.

\begin{rema}\label{higgs1-tor25}
Soient $(\tX',\cM_{\tX'})$ une autre $(\cA_2(\oS),\cM_{\cA_2(\oS)})$-déformation lisse de 
$(\coX,\cM_{\coX})$, $\cL'$ le torseur de Higgs-Tate, 
$\cC'$ la $\hoR$-algèbre de Higgs-Tate et $\cF'$ la $\hoR$-extension de Higgs-Tate associés. 
On a alors un isomorphisme de $(\cA_2(\oS),\cM_{\cA_2(\oS)})$-déformations
\begin{equation}
u\colon (\tX,\cM_{\tX})\stackrel{\sim}{\rightarrow} (\tX',\cM_{\tX'}).
\end{equation}
L'isomorphisme de $\trT$-torseurs $\cL\stackrel{\sim}{\rightarrow} \cL'$, 
$\psi\mapsto u\circ \psi$ \eqref{higgs1-tor2a} induit un isomorphisme $\hoR$-linéaire et $\Delta$-équivariant 
\begin{equation}
\cF'\stackrel{\sim}{\rightarrow}\cF,
\end{equation}
qui s'insère dans un diagramme commutatif \eqref{higgs1-tor2e}
\begin{equation}
\xymatrix{
0\ar[r]&{\hoR}\ar[r]\ar@{=}[d]&{\cF'}\ar[r]\ar[d]&{\xi^{-1}\tOmega^1_{R/\co_K} \otimes_R \hoR}\ar[r]\ar@{=}[d]&0\\
0\ar[r]&{\hoR}\ar[r]&{\cF}\ar[r]&{\xi^{-1}\tOmega^1_{R/\co_K} \otimes_R \hoR}\ar[r] & 0}
\end{equation} 
On en déduit un $\hoR$-isomorphisme $\Delta$-équivariant
\begin{equation}
\cC'\stackrel{\sim}{\rightarrow} \cC.
\end{equation}
\end{rema}

\subsection{}\label{higgs1-tor8}\index{101030@$\sigma_{\psi,\psi'}$}
Toute section $\psi\in \cL(\hY)$ définit un scindage $v_{\psi}\colon \cF\rightarrow \hoR$
de la suite exacte \eqref{higgs1-tor2e}, à savoir, le morphisme qui à toute fonction affine $\ell$ sur $\cL$ associe $\ell(\psi)$. 
Le morphisme $\id_{\cF}-v_{\psi}$ induit un morphisme $\hoR$-linéaire 
\begin{equation}\label{higgs1-tor8dd}
u_{\psi}\colon \xi^{-1}\tOmega^1_{R/\co_K}\otimes_R\hoR\rightarrow \cF,
\end{equation}
qui n'est autre que la restriction de l'isomorphisme 
$(\ttt_{\psi}^*)^{-1}\colon \cS\rightarrow \cC$ \eqref{higgs1-tor202b} à $\cS^1$. 

Pour tous $\psi,\psi'\in \cL(\hY)$, la différence $\psi- \psi'\in \bT(\hY)=\rT$ \eqref{higgs1-tor2ab} 
détermine un morphisme $\hoR$-linéaire 
\begin{equation}\label{higgs1-tor8a}
\sigma_{\psi,\psi'}\colon \xi^{-1} \tOmega^1_{R/\co_K}\otimes_R\hoR\rightarrow \hoR.
\end{equation}
On a 
\begin{equation}\label{higgs1-tor8ee}
\sigma_{\psi,\psi'}= u_{\psi'}-u_{\psi}.
\end{equation}

D'après \eqref{higgs1-tor200h}, pour tous $g\in \Delta$ et $\psi\in \cL(\hY)$, on a 
\begin{equation}\label{higgs1-tor85b}
u_{^g\psi}=g\circ u_{\psi} \circ g^{-1}.
\end{equation}
Par suite, pour tout $x\in \xi^{-1} \tOmega^1_{R/\co_K}\otimes_R\hRun$, 
l'application $g\mapsto \sigma_{\psi,^g\psi}(x)$ est un cocycle de $\Delta$ à valeurs dans $\hoR$. 
L'application induite 
\begin{equation}\label{higgs1-tor85bb}
\xi^{-1} \tOmega^1_{R/\co_K}\otimes_R\hRun \rightarrow \rH^1(\Delta,\hoR)
\end{equation}
est clairement le bord de la suite exacte longue de cohomologie déduite de la suite exacte \eqref{higgs1-tor2e}. 

D'après \eqref{higgs1-tor202f} et \eqref{higgs1-tor21d}, notant $D_{\psi-{^g\psi}}$ l'image de  $\psi-{^g\psi}$ par le morphisme 
\eqref{higgs1-tor21a},  on a 
\begin{equation}\label{higgs1-tor8b}
\varphi_{\psi}(g)=\exp(D_{\psi-{^g\psi}})\circ g,
\end{equation}
où $g$ désigne l'action canonique de $g$ sur $\cS$ \eqref{higgs1-tor200}. En particulier, d'après \eqref{higgs1-tor22e}, 
pour tout $x\in \xi^{-1}\tOmega^1_{R/\co_K}\subset \cS^1$, on a 
\begin{equation}
\varphi_\psi(g)(x)=x+\sigma_{\psi,{^g\psi}}(x).
\end{equation}

\begin{rema}\label{higgs1-ext21}
On dit qu'un élément $\psi$ de $\cL(\hY)$ est {\em optimal} 
si pour tout $g\in \Delta$, le morphisme $\psi-{^g\psi}\in \bT(\hY)=\rT$ se factorise en 
\begin{equation}\label{higgs1-ext21a}
\tOmega^1_{R/\co_K}\otimes_R\hoR\rightarrow \hoR(1)\stackrel{\sim}{\rightarrow} p^{\frac{1}{p-1}}\xi\hoR\rightarrow \xi\hoR,
\end{equation}
où la deuxième flèche est l'isomorphisme \eqref{higgs1-ext11c} la dernière flèche est l'injection canonique. 

Notons $\cS'$ la $\hoR$-algèbre graduée \eqref{higgs1-not621}
\begin{equation}\label{higgs1-tor10a}
\cS'=\rS_{\hoR}(\tOmega^1_{R/\co_K}\otimes_R\hoR(-1)).
\end{equation} 
L'isomorphisme canonique \eqref{higgs1-ext11c}
\begin{equation}\label{higgs1-tor10b}
\xi^{-1}\tOmega^1_{R/\co_K}\otimes_R\hoR\stackrel{\sim}{\rightarrow} 
p^{\frac{1}{p-1}}\tOmega^1_{R/\co_K}\otimes_R\hoR(-1)
\end{equation}
induit un $\hoR$-homomorphisme 
\begin{equation}\label{higgs1-tor10c}
\jmath\colon \cS\rightarrow \cS'
\end{equation}
tel que $\jmath\otimes_{\mZ_p}\mQ_p$ soit un isomorphisme. 
Pour tout $\psi\in \cL(\hY)$, 
on note encore $\varphi_{\psi}$ l'action de $\Delta$ sur $\cS'[\frac 1 p]$ déduite 
de son action sur $\cS$ définie dans \eqref{higgs1-tor202g}. 
Pour que cette action préserve $\cS'$, il faut et il suffit que $\psi$ soit optimal.
En effet, notons 
\begin{equation}
\varsigma_{\psi,{^g\psi}}\colon \tOmega^1_{R/\co_K}\otimes_R\hoR(-1)\rightarrow p^{-\frac{1}{p-1}}\hoR
\end{equation}
le morphisme $\hoR$-linéaire déduit de $\psi-{^g\psi}\in \bT(\hY)=\rT$ et de l'isomorphisme \eqref{higgs1-ext11c}.
Pour tous $g\in \Delta$ et $x\in \tOmega^1_{R/\co_K}(-1)\subset \cS'$, on a alors
\begin{equation}
\varphi_{\psi}(x)=x+\varsigma_{\psi,{^g\psi}}(x).
\end{equation}
Par suite, pour que $\varphi_{\psi}$ préserve $\cS'$, il faut et il suffit que $\varsigma_{\psi,{^g \psi}}$ 
se factorise à travers $\hoR \subset p^{-\frac{1}{p-1}}\hoR$ pour tout $g\in \Delta$, ce qui est équivalent au fait
que $\psi$ est optimal.
\end{rema}

\subsection{}\label{higgs1-ext24}\index{101037@$(\tX_0,\cM_{\tX_0})$, $\cL_0$, $\cF_0$, $\cC_0$, $\psi_0\in \cL_0(\hY)$}
D'après \ref{higgs1-dlog1}(C$_5$), il existe essentiellement un unique morphisme étale 
\begin{equation}\label{higgs1-ext24ac}
(\tX_0,\cM_{\tX_0})\rightarrow (\cA_2(\oS),\cM_{\cA_2(\oS)}) \times_{\bA_\mN}\bA_P
\end{equation}
qui s'insère dans un diagramme commutatif à carrés cartésiens
\begin{equation}\label{higgs1-ext24a}
\xymatrix{
{(\coX,\cM_\coX)}
\ar[r]\ar[d]&{(\coS,\cM_{\coS})\times_{\bA_\mN}\bA_P}\ar[d]\ar[r]&{(\coS,\cM_{\coS})}\ar[d]^{i_\oS}\\
{(\tX_0,\cM_{\tX_0})}\ar[r]&{(\cA_2(\oS),\cM_{\cA_2(\oS)})
\times_{\bA_\mN}\bA_P}\ar[r]\ar[d]&{(\cA_2(\oS),\cM_{\cA_2(\oS)})}\ar[d]^a\\
&{\bA_P}\ar[r]&{\bA_\mN}}
\end{equation}
où le morphisme $a$ est défini par la carte 
$\mN\rightarrow \Gamma(\cA_2(\oS),\cM_{\cA_2(\oS)}), 1\mapsto \tpi$ \eqref{higgs1-eplog9}. 
On dit que $(\tX_0,\cM_{\tX_0})$ est la $(\cA_2(\oS),\cM_{\cA_2(\oS)})$-déformation lisse de 
$(\coX,\cM_\coX)$ définie par la carte $(P,\gamma)$. On désigne par $\cL_0$ le torseur de Higgs-Tate
\eqref{higgs1-tor2}, par $\cC_0$ la $\hoR$-algèbre de Higgs-Tate 
et par $\cF_0$ la $\hoR$-extension de Higgs-Tate associés à $(\tX_0,\cM_{\tX_0})$ \eqref{higgs1-tor201}.

Le diagramme 
\begin{equation}\label{higgs1-ext24ab}
\xymatrix{
{(\hY,\cM_\hY)}\ar[r]^-(0.4){i_Y}\ar[d]&{(\cA_2(Y),\cM_{\cA_2(Y)})}\ar[r]^-(0.4)b\ar[d]&{\bA_P}\ar[d]\\
{(\coS,\cM_\coS)}\ar[r]^-(0.4){i_\oS}&{(\cA_2(\oS),\cM_{\cA_2(\oS)})}\ar[r]^-(0.4)a&{\bA_\mN}}
\end{equation}
où le morphisme $b$ est défini par la carte canonique $P\rightarrow \Gamma(\cA_2(Y),\cM_{\cA_2(Y)})$
\eqref{higgs1-eplog5} est commutatif d'après \eqref{higgs1-eplog25a}. Il est clair que le diagramme (sans la flèche pointillée) 
\begin{equation}\label{higgs1-ext24b}
\xymatrix{
{(\hY,\cM_{\hY})}\ar[r]\ar[d]_{i_Y}&{(\coX,\cM_\coX)}\ar[r]&{(\tX_0,\cM_{\tX_0})}\ar[d]\\
{(\cA_2(Y),\cM_{\cA_2(Y)})}\ar[rr]^{\phi_0}\ar@{.>}[rru]^{\psi_0}&&{(\cA_2(\oS),\cM_{\cA_2(\oS)})
\times_{\bA_\mN}\bA_P}}
\end{equation}
où $\phi_0$ est défini par \eqref{higgs1-ext24ab} est commutatif. On peut le compléter par une unique flèche 
pointillée  
$\psi_0\in \cL_0(\hY)$ de façon à le laisser commutatif. 
On dit que $\psi_0$ est la section de $\cL_0(\hY)$ définie par la carte $(P,\gamma)$. 
Le lecteur observera ici que $(\tX_0,\cM_{\tX_0})$ et $\psi_0$ dépendent du choix du morphisme \eqref{higgs1-dlog5b},
puisque les cartes $a$ et $b$ en dépendent.

\begin{prop}\label{higgs1-ext245}
Conservons les notations de \eqref{higgs1-ext24}. Pour tous $t\in P^\gp$ et $g\in \Delta$, on a 
\begin{equation}\label{higgs1-ext245a}
(\psi_0-{^g\psi_0})(d\log(t))=-\log([\chi_t(g)]),
\end{equation}  
où  l'on considère $\psi_0-{^g\psi_0}$ comme un élément de $\bT(\hY)=\rT$ \eqref{higgs1-tor2ab},
$d\log(t)$ est l'image canonique de $t$ dans $\tOmega^1_{R/\co_K}$ \eqref{higgs1-ext-log13e} et 
$\log([\chi_t])$ désigne l'homomorphisme composé 
\begin{equation}\label{higgs1-ext245b}
\xymatrix{
\Delta\ar[r]&{\Delta_\infty}\ar[r]^-(0.5){\chi_t}& 
{\mZ_p(1)}\ar[r]^-(0.5){\log([\ ])}&{p^{\frac{1}{p-1}}\xi\hoR}\ar[r]&{\xi\hoR}},
\end{equation}
où la première et la dernière flèches sont les morphismes canoniques, $\chi_t$ est l'homomorphisme \eqref{higgs1-ext-log14e} 
et la troisième flèche est induite par l'isomorphisme \eqref{higgs1-ext11c}. En particulier, $\psi_0$ est optimal \eqref{higgs1-ext21}. 
\end{prop}
  
Comme les deux membres de l'équation \eqref{higgs1-ext245a} sont des homomorphismes de $P^\gp$ dans $\xi\hoR$,
on peut se borner au cas où $t\in P$.   
Les morphismes $\phi_0$ et $\phi_0\circ g^{-1}$, où $\phi_0$ est le morphisme défini dans \eqref{higgs1-ext24b} 
et $g^{-1}$ agit sur $(\cA_2(Y),\cM_{\cA_2(Y)})$, prolongent le même morphisme 
\[
(\hY,\cM_\hY)\rightarrow (\cA_2(\oS),\cM_{\cA_2(\oS)})\times_{\bA_\mN}\bA_P.
\]
D'après les définitions et la condition \ref{higgs1-dlog1}(C$_5$), la différence $\phi_0-\phi_0\circ g^{-1}$ 
correspond au morphisme $\psi_0-{^g\psi_0}\in \rT$. 
D'autre part, on a $g(\tlt)=[\chi_t(g)]\cdot \tlt$ dans $\Gamma(\cA_2(Y),\cM_{\cA_2(Y)})$ \eqref{higgs1-eip10b}.
La première proposition s'ensuit compte tenu de \ref{higgs1-log10} et \eqref{higgs1-ext11b}. 
La seconde proposition est une conséquence immédiate de la première. 

\begin{lem}\label{higgs1-ext2451}
On a un diagramme commutatif 
\begin{equation}
\xymatrix{
{\xi^{-1}\tOmega^1_{R/\co_K}\otimes_R\hRun}\ar[d]_u\ar[r]^-(0.5){\partial}&{\rH^1(\Delta,\hoR)}\\
{\rH^1(\Delta,\xi^{-1}\hoR(1))}\ar[r]_\sim^{-v}&{\rH^1(\Delta,p^{\frac{1}{p-1}}\hoR)}\ar[u]_w}
\end{equation}
où $\partial$ est le bord de la suite exacte longue de cohomologie déduite de la suite exacte courte \eqref{higgs1-tor2e},
$u$ est induit par le morphisme \eqref{higgs1-cg35d}, $v$ est induit par l'isomorphisme \eqref{higgs1-ext11c} et $w$
est induit par l'injection canonique $p^{\frac{1}{p-1}}\hoR \rightarrow \hoR$. 
\end{lem}

En effet, il suffit de montrer que le diagramme
\begin{equation}
\xymatrix{
{\xi^{-1}\tOmega^1_{R/\co_K}\otimes_R\hRun}\ar[d]_{\xi^{-1} \tdelta}\ar[r]^{\partial}&{\rH^1(\Delta,\hoR)}\\
{\Hom_{\mZ}(\Delta_\infty,\xi^{-1}\hRun(1))}\ar[r]_\sim^{-v'}&{\Hom_{\mZ}(\Delta_{\infty},p^{\frac{1}{p-1}}\hRun)}\ar[u]_{w'}}
\end{equation}
où $\tdelta$ est l'isomorphisme \eqref{higgs1-ext-log14h}, 
$v'$ provient de l'isomorphisme $\co_C(1)\stackrel{\sim}{\rightarrow}p^{\frac{1}{p-1}}\xi\co_C$
\eqref{higgs1-ext11c} et $w'$ est induit par l'injection canonique $p^{\frac{1}{p-1}}\hRun \rightarrow \hoR$, est commutatif. 
D'après \ref{higgs1-tor25}, $\partial$ ne dépend pas de la déformation $(\tX,\cM_\tX)$. 
On peut donc se borner au cas où $(\tX,\cM_\tX)$
est la déformation définie par la carte $(P,\gamma)$ \eqref{higgs1-ext24}. Soit $\psi_0$ 
la section de $\cL_0(\hY)$ définie par la carte $(P,\gamma)$. 
D'après \ref{higgs1-tor8}, pour tout $x\in \xi^{-1} \tOmega^1_{R/\co_K}\otimes_R\hRun$, $\partial(x)$ est la classe du cocycle 
$g\mapsto \sigma_{\psi_0,^g\psi_0}(x)$ de $\Delta$ à valeurs dans $\hoR$. 
L'assertion résulte donc de \ref{higgs1-ext245}.

\begin{rema}\label{higgs1-ext2452}
Il résulte aussitôt de \ref{higgs1-ext2451} que pour tout élément non-nul $a$ de $\co_\oK$, on a un diagramme commutatif 
\begin{equation}
\xymatrix{
{\xi^{-1}\tOmega^1_{R/\co_K}\otimes_R(R_1/aR_1)}\ar[d]_{u_a}\ar[r]^-(0.5){\partial_a}&{\rH^1(\Delta,\oR/a\oR)}\\
{\rH^1(\Delta,\xi^{-1}\oR(1)/a\xi^{-1}\oR(1))}\ar[r]_\sim^{-v_a}&
{\rH^1(\Delta,p^{\frac{1}{p-1}}\oR/ap^{\frac{1}{p-1}}\oR)}\ar[u]_{w_a}}
\end{equation}
où $\partial_a$ est induit par le bord de la suite exacte longue de cohomologie déduite de la suite exacte courte 
obtenue en réduisant \eqref{higgs1-tor2e} modulo $a$,
$u_a$ est induit par le morphisme \eqref{higgs1-cg35c}, $v_a$ est induit par l'isomorphisme \eqref{higgs1-ext11c} et $w_a$
est induit par l'injection canonique $p^{\frac{1}{p-1}}\oR \rightarrow \oR$.  
\end{rema}

\subsection{}\label{higgs1-ext26}\index{101040@$\varphi_0$}\index{101042@$\fS$}
On désigne par $\varphi_0=\varphi_{\psi_0}$ \eqref{higgs1-tor202g} l'action de $\Delta$ sur $\cS$ induite par la section 
$\psi_0\in \cL_0(\hY)$ définie par la carte $(P,\gamma)$ \eqref{higgs1-ext24}. 
Soient $t_1,\dots,t_d\in P^\gp$ tels que leurs images dans $(P^\gp/\mZ\lambda)\otimes_{\mZ}\mZ_p$ 
forment une $\mZ_p$-base, de sorte que $(d\log(t_i))_{1\leq i\leq d}$ est une $R$-base
de $\tOmega^1_{R/\co_K}$ \eqref{higgs1-ext-log13e}. 
Pour tout $1\leq i\leq d$, posons $y_i=\xi^{-1}d\log(t_i)\in \xi^{-1}\tOmega^1_{R/\co_K}\subset \cS$ 
et notons $\chi_i$ l'homomorphisme composé 
\begin{equation}\label{higgs1-ext26a}
\xymatrix{
\Delta\ar[r]&{\Delta_\infty}\ar[r]^-(0.5){\chi_{t_i}}& 
{\mZ_p(1)}\ar[r]^-(0.5){\log([\ ])}&{p^{\frac{1}{p-1}}\xi\co_C}},
\end{equation}
où la première flèche est le morphisme canonique, $\chi_{t_i}$ est l'homomorphisme \eqref{higgs1-ext-log14e} 
et la troisième flèche est induite par l'isomorphisme \eqref{higgs1-ext11c}.   
Il résulte alors de \eqref{higgs1-tor8b} et \eqref{higgs1-ext245a} que pour tout $g\in \Delta$, on a 
\begin{equation}\label{higgs1-ext26b}
\varphi_0(g)=\exp(-\sum_{i=1}^d\xi^{-1}\frac{\partial}{\partial y_i}\otimes \chi_i(g)) \circ g.
\end{equation}
Par suite, $\varphi_0$ préserve la sous-$\hRun$-algèbre 
\begin{equation}\label{higgs1-ext26c}
\fS=\rS_{\hRun}(\xi^{-1}\tOmega^1_{R/\co_K}\otimes_R\hRun)
\end{equation} 
de $\cS$, et l'action induite de $\Delta$ sur $\fS$ se factorise à travers $\Delta_{p^\infty}$.

\begin{teo}\label{higgs1-ext246}
Conservons les notations de \eqref{higgs1-ext24} et 
notons $(P^\gp/\mZ\lambda)_{\lib}$ le quotient de $P^\gp/\mZ\lambda$ par son sous-module de torsion. 
Alors la donnée d'un inverse à droite 
\begin{equation}\label{higgs1-ext246a}
w\colon (P^\gp/\mZ\lambda)_{\lib}\rightarrow P^\gp
\end{equation}
du morphisme canonique $P^\gp\rightarrow (P^\gp/\mZ\lambda)_\lib$ 
détermine uniquement un morphisme $\hoR$-linéaire et $\Delta$-équivariant 
\begin{equation}\label{higgs1-ext246b}
\beta_w\colon p^{-\frac{1}{p-1}}\cF_0\rightarrow \cE
\end{equation}
qui s'insère dans un diagramme commutatif 
\begin{equation}\label{higgs1-ext246c}
\xymatrix{
0\ar[r]&{p^{-\frac{1}{p-1}}\hoR}\ar[r]\ar@{^(->}[d]&{p^{-\frac{1}{p-1}}\cF_0}\ar[r]\ar[d]^{\beta_w}&
{p^{-\frac{1}{p-1}}\xi^{-1}\tOmega^1_{R/\co_K} \otimes_R \hoR}\ar[r]\ar[d]^{-c}&0\\
0\ar[r]&{(\pi\rho)^{-1}\hoR}\ar[r]&{\cE}\ar[r]&{\tOmega^1_{R/\co_K} \otimes_R \hoR(-1)}\ar[r] & 0}
\end{equation}
où les lignes proviennent des suites exactes \eqref{higgs1-log-ext17b} et \eqref{higgs1-tor2e} 
et $c$ est induit par l'isomorphisme \eqref{higgs1-ext11c}.
\end{teo}

Comme le sous-groupe de torsion de $P^\gp/\mZ\lambda$
est d'ordre premier à $p$, l'isomorphisme \eqref{higgs1-log-ext45b} induit un isomorphisme 
\[
(P^\gp/\mZ\lambda)_\lib\otimes_\mZ R\stackrel{\sim}{\rightarrow}\tOmega^1_{R/\co_K}.
\]
Par suite, le composé de $w$ avec l'homomorphisme $P^\gp\rightarrow \cE(1), t\mapsto d\log(\tlt)$
\eqref{higgs1-ext-log13cd} induit un morphisme $\hoR$-linéaire que l'on note
\begin{equation}\label{higgs1-ext246e}
\sigma_{\psi_0}\colon \tOmega^1_{R/\co_K}\otimes_R\hoR\rightarrow \cE(1).
\end{equation}
D'après \eqref{higgs1-ext-log13d}, $\sigma_{\psi_0}$ est un scindage de la suite exacte \eqref{higgs1-log-ext17b} tordue par $\mZ_p(1)$. 
Pour tout $x\in \cE(1)$, posons 
\begin{equation}\label{higgs1-ext246f}
\langle \psi_0,x\rangle =x-\sigma_{\psi_0}(\nu(x))\in (\pi\rho)^{-1}\hoR(1),
\end{equation}
où $\nu\colon \cE(1)\rightarrow \tOmega^1_{R/\co_K}\otimes_R\hoR$ est le morphisme canonique. 
On définit une application 
\begin{equation}\label{higgs1-ext246g}
\langle \ ,\ \rangle \colon 
\cL_0(\hY)\times \cE(1) \rightarrow (\pi\rho)^{-1} \hoR(1)
\end{equation}
de la façon suivante~: pour tout $\phi\in \rT$ \eqref{higgs1-tor2ab}, si on note $\ophi$ le morphisme composé
\begin{equation}\label{higgs1-ext246gh}
\tOmega^1_{R/\co_K}\otimes_R\hoR\stackrel{\phi}{\longrightarrow}\xi\hoR\stackrel{\sim}{\longrightarrow}
p^{-\frac{1}{p-1}}\hoR(1)\longrightarrow (\pi\rho)^{-1}\hoR(1),
\end{equation}
où la deuxième flèche est induite par l'isomorphisme \eqref{higgs1-ext11c} et la troisième flèche est l'injection canonique, 
on a 
\begin{equation}\label{higgs1-ext246h}
\langle\psi_0+\phi,x\rangle=\langle\psi_0,x\rangle- \ophi(\nu(x)).
\end{equation}

Pour tout $\psi\in  \cL_0(\hY)$, le morphisme $\cE(1)\rightarrow (\pi\rho)^{-1} \hoR(1),
x\mapsto \langle \psi,x\rangle$ est un scindage de la suite exacte \eqref{higgs1-log-ext17b} tordue par $\mZ_p(1)$.

Montrons que pour tous $\psi\in  \cL_0(\hY)$, $x\in \cE(1)$ et $g\in \Delta$, on a 
\begin{equation}\label{higgs1-ext246i}
g(\langle \psi,x\rangle)=\langle {^g\psi},g(x)\rangle.
\end{equation}
Considérons d'abord le cas $\psi=\psi_0$. Il suffit alors de montrer que pour tout 
$x\in \tOmega^1_{R/\co_K}\otimes_R\hoR$, on a 
\begin{equation}\label{higgs1-ext246j}
\osigma_{\psi_0,{^g\psi_0}}(x)=\sigma_{\psi_0}(x)-g(\sigma_{\psi_0}(g^{-1}x)),
\end{equation}
où $\osigma_{\psi_0,{^g\psi_0}}$ est le morphisme défini dans \eqref{higgs1-ext246gh} 
à partir de $\psi_0-{^g\psi_0}\in \rT$ \eqref{higgs1-tor8}. 
Soient $t\in P$, $y$ son image canonique dans $(P^\gp/\mZ\lambda)_\lib$, $z=w(y)\in P^\gp$. On a alors 
\[
\sigma_{\psi_0}(d\log(t))=d\log(\tz).
\]
D'autre part, en vertu de \ref{higgs1-ext245}, on a 
\[
\osigma_{\psi_0,{^g\psi_0}}(d\log(t))=(\chi_t(g))^{-1}, 
\]
où  $(\chi_t(g))^{-1}\in \mZ_p(1)\subset \hoR(1)\subset (\pi\rho)^{-1}\hoR(1)$. On observera que le 
caractère $\log([\ ])$ disparait de la formule à cause de la définition \eqref{higgs1-ext246gh}. 
Comme l'image de $\chi_t(g)$ par l'injection canonique $(\pi\rho)^{-1}\hoR(1)\rightarrow \cE(1)$ 
est $d\log(\chi_t(g))$ \eqref{higgs1-ext-log13}, il suffit encore de montrer la relation suivante dans $\cE(1)$
\begin{equation}\label{higgs1-ext246m}
g(d\log(\tz))=d\log(\tz)+\log(\chi_t(g)).
\end{equation}
On peut évidemment remplacer $t$ par une puissance, et donc supposer $t-z\in \mZ\lambda$. 
Comme $\Delta$ fixe $d\log(\tlambda)$, la relation \eqref{higgs1-ext246m} résulte de \eqref{higgs1-ext-log15a}. 

Dans le cas général, pour tous $\phi\in \rT$, $x\in \cE(1)$ et $g\in \Delta$, on a 
\begin{eqnarray*}
g(\langle \psi_0+\phi,x\rangle)&=&g(\langle \psi_0,x\rangle)-g\circ \ophi(\nu(x))\\
&=&\langle {^g\psi_0},g(x)\rangle-g\circ \ophi\circ g^{-1}(\nu(g(x)))\\
&=&\langle {^g\psi_0}+g\circ \phi\circ g^{-1},g(x)\rangle \\
&=& \langle {^g(\psi_0+\phi)},g(x)\rangle,
\end{eqnarray*}
ce qui achève la preuve de \eqref{higgs1-ext246i}.

Pour tout $x\in \cE(1)$, l'application $\psi\mapsto \langle \psi,x\rangle$ est une fonction affine sur $\cL_0$ 
à valeurs dans $(\pi\rho)^{-1}\hoR(1)$ (cf. \ref{higgs1-eph10}). Elle est donc naturellement un élément de $(\pi\rho)^{-1}\cF_0(1)$. 
L'application 
\begin{equation}
\cE(1)\rightarrow(\pi\rho)^{-1}\cF_0(1),\ \ \ x\mapsto(\psi\mapsto\langle \psi,x\rangle)
\end{equation}
est un morphisme $\hoR$-linéaire et $\Delta$-équivariant en vertu de \eqref{higgs1-ext246i}. 
Le morphisme induit 
\[
\cE\rightarrow (\pi\rho)^{-1}\cF_0
\] 
s'insère dans un digramme commutatif
\[
\xymatrix{
0\ar[r]&{(\pi\rho)^{-1}\hoR}\ar[r]\ar@{=}[dd]&{\cE}\ar[r]\ar[dd]&
{\tOmega^1_{R/\co_K}\otimes_R\hoR(-1)}\ar[r]\ar[d]^{i_1}&0\\
&&&
{(\pi\rho)^{-1}p^{\frac{1}{p-1}}\tOmega^1_{R/\co_K}\otimes_R\hoR(-1)}\ar[r]\ar[d]^{-c^{-1}}&0\\
0\ar[r]&{(\pi\rho)^{-1}\hoR}\ar[r]&{(\pi\rho)^{-1}\cF_0}\ar[r]&{(\pi\rho)^{-1}\xi^{-1}\tOmega^1_{R/\co_K}\otimes_R\hoR}\ar[r]&0}
\]
où $i_1$ est l'injection canonique. Considérons le diagramme commutatif 
\[
\xymatrix{
0\ar[r]&{p^{-\frac{1}{p-1}}\hoR}\ar[r]\ar[d]_{i_3}&{p^{-\frac{1}{p-1}}\cF_0}\ar[r]\ar[d]&{p^{-\frac{1}{p-1}}\xi^{-1}\tOmega^1_{R/\co_K}\otimes_R\hoR}
\ar@{=}[d]\ar[r]&0\\
0\ar[r]&{(\pi\rho)^{-1}\hoR}\ar[r]\ar@{=}[d]&{\cH}\ar[r]\ar[d]&{p^{-\frac{1}{p-1}}\xi^{-1}\tOmega^1_{R/\co_K}\otimes_R\hoR}
\ar[d]_{i_2}\ar[r]&0\\
0\ar[r]&{(\pi\rho)^{-1}\hoR}\ar[r]&{(\pi\rho)^{-1}\cF_0}\ar[r]&{(\pi\rho)^{-1}\xi^{-1}\tOmega^1_{R/\co_K}\otimes_R\hoR}\ar[r]&0}
\]
où $i_2$ et $i_3$ sont les injections canoniques et $\cH$ est l'image inverse de $(\pi\rho)^{-1}\cF_0$ par $i_2$. 
Comme $i_2\circ c^{-1}=c^{-1} \circ i_1$, on en déduit un isomorphisme $\alpha\colon \cH\stackrel{\sim}{\rightarrow}\cE$
qui s'insère dans un diagramme commutatif 
\[
\xymatrix{
0\ar[r]&{(\pi\rho)^{-1}\hoR}\ar[r]\ar@{=}[d]&{\cH}\ar[r]\ar[d]^\alpha&
{p^{-\frac{1}{p-1}}\xi^{-1}\tOmega^1_{R/\co_K} \otimes_R \hoR}\ar[r]\ar[d]^{-c}&0\\
0\ar[r]&{(\pi\rho)^{-1}\hoR}\ar[r]&{\cE}\ar[r]&{\tOmega^1_{R/\co_K} \otimes_R \hoR(-1)}\ar[r] & 0}
\]
Le morphisme $p^{-\frac{1}{p-1}}\cF_0\rightarrow \cE$ composé de l'injection 
canonique $p^{-\frac{1}{p-1}}\cF_0\rightarrow\cH$ et de $\alpha$ répond alors à la question.

\begin{cor}\label{higgs1-ext247}
Il existe un morphisme $\hoR$-linéaire et $\Delta$-équivariant 
\begin{equation}\label{higgs1-ext247b}
p^{-\frac{1}{p-1}}\cF\rightarrow \cE
\end{equation}
qui s'insère dans un diagramme commutatif 
\begin{equation}\label{higgs1-ext247c}
\xymatrix{
0\ar[r]&{p^{-\frac{1}{p-1}}\hoR}\ar[r]\ar@{^(->}[d]&{p^{-\frac{1}{p-1}}\cF}\ar[r]\ar[d]&
{p^{-\frac{1}{p-1}}\xi^{-1}\tOmega^1_{R/\co_K} \otimes_R \hoR}\ar[r]\ar[d]^{-c}&0\\
0\ar[r]&{(\pi\rho)^{-1}\hoR}\ar[r]&{\cE}\ar[r]&{\tOmega^1_{R/\co_K} \otimes_R \hoR(-1)}\ar[r] & 0}
\end{equation}
où les lignes proviennent des suites exactes \eqref{higgs1-log-ext17b} et \eqref{higgs1-tor2e} 
et $c$ est induit par l'isomorphisme \eqref{higgs1-ext11c}.
\end{cor}

Cela résulte de \ref{higgs1-tor25} et \ref{higgs1-ext246}.
 
\section{Cohomologie galoisienne II}\label{higgs1-CGII}

Dans cette section, nous reprenons en les généralisant légèrement des résultats de T. Tsuji \cite{tsuji3}.  

\subsection{}\label{higgs1-chb1}\index{101100@$\cS^{(r)}$, $\cS^\dagger$}\index{101105@$d_{\cS^{(r)}}$, $d_{\hcS^{(r)}}$, $d_{\cS^\dagger}$}
\index{101106@$\tmK^\bullet(\hcS^{(r)},p^rd_{\hcS^{(r)}})$}
On note $\cS$ la $\hoR$-algèbre symétrique \eqref{higgs1-tor2d}
\begin{equation}
\cS=\rS_{\hoR}(\xi^{-1}\tOmega^1_{R/\co_K}\otimes_R\hoR)
\end{equation}
et $\hcS$ son séparé complété $p$-adique. 
Pour tout nombre rationnel $r\geq 0$, on désigne par $\cS^{(r)}$ la sous-$\hoR$-algèbre de $\cS$ définie par \eqref{higgs1-not621}
\begin{equation}\label{higgs1-chb1a}
\cS^{(r)}=\rS_{\hoR}(p^r\xi^{-1}\tOmega^1_{R/\co_K}\otimes_R\hoR)
\end{equation}
et par $\hcS^{(r)}$ son séparé complété $p$-adique que l'on suppose toujours muni de la topologie $p$-adique.
On munit $\hcS^{(r)}\otimes_{\mZ_p}\mQ_p$ de la topologie $p$-adique \eqref{higgs1-not54}. 
Compte tenu de \ref{higgs1-pur8} et de sa preuve, $\cS^{(r)}$ et $\hcS^{(r)}$ sont $\co_C$-plats. 
Pour tous nombres rationnels $r'\geq r\geq 0$, on a un homomorphisme injectif canonique 
$\alpha^{r,r'}\colon \cS^{(r')}\rightarrow \cS^{(r)}$. On vérifie aussitôt que l'homomorphisme 
induit $\halpha^{r,r'}\colon\hcS^{(r')}\rightarrow \hcS^{(r)}$ est injectif. On pose 
\begin{equation}\label{higgs1-chb1b}
\cS^\dagger=\underset{\underset{r\in \mQ_{>0}}{\longrightarrow}}{\lim} \hcS^{(r)},
\end{equation}
que l'on identifie à une sous-$\hoR$-algèbre de $\hcS$ par la limite inductive des homomorphismes 
$\halpha^{0,r}$.

On a un $\cS^{(r)}$-isomorphisme canonique 
\begin{equation}\label{higgs1-chb1c}
\Omega^1_{\cS^{(r)}/\hoR}\stackrel{\sim}{\rightarrow}\xi^{-1}\tOmega^1_{R/\co_K}\otimes_R\cS^{(r)}.
\end{equation} 
On désigne par
\begin{equation}\label{higgs1-chb1d}
d_{\cS^{(r)}}\colon \cS^{(r)}\rightarrow \xi^{-1}\tOmega^1_{R/\co_K}\otimes_R\cS^{(r)}
\end{equation}
la $\hoR$-dérivation universelle de $\cS^{(r)}$ et par 
\begin{equation}\label{higgs1-chb1e}
d_{\hcS^{(r)}}\colon \hcS^{(r)}\rightarrow \xi^{-1}\tOmega^1_{R/\co_K}\otimes_R\hcS^{(r)}
\end{equation}
son prolongement aux complétés (on notera que le $R$-module $\tOmega^1_{R/\co_K}$ est libre de type fini). 
Comme $\xi^{-1}\tOmega^1_{R/\co_K}\otimes_R \hoR\subset d_{\cS^{(r)}}(\cS^{(r)})$, 
$d_{\cS^{(r)}}$ et $d_{\hcS^{(r)}}$ sont également des $\hoR$-champs de Higgs à coefficients dans 
$\xi^{-1}\tOmega^1_{R/\co_K}$ (cf. \ref{higgs1-not16}, \ref{higgs1-not22} et \ref{higgs1-eip75}). 
On désigne par $\mK^\bullet(\hcS^{(r)},p^rd_{\hcS^{(r)}})$ le complexe de Dolbeault du module de Higgs 
$(\hcS^{(r)},p^rd_{\hcS^{(r)}})$ \eqref{higgs1-not61} et 
par $\tmK^\bullet(\hcS^{(r)},p^rd_{\hcS^{(r)}})$ le complexe de Dolbeault augmenté~:
\begin{equation}\label{higgs1-chb1f}
\hoR\rightarrow \mK^0(\hcS^{(r)},p^rd_{\hcS^{(r)}})\rightarrow \mK^1(\hcS^{(r)},p^rd_{\hcS^{(r)}})\rightarrow \dots
\rightarrow \mK^n(\hcS^{(r)},p^rd_{\hcS^{(r)}})\rightarrow \dots,
\end{equation}
où $\hoR$ est placé en degré $-1$ et la différentielle $\hoR\rightarrow\hcS^{(r)}$ est l'homomorphisme canonique.

Pour tous nombres rationnels $r'\geq r\geq 0$, on a 
\begin{equation}\label{higgs1-chb1g}
p^{r'}(\id \times \alpha^{r,r'}) \circ d_{\cS^{(r')}}=p^rd_{\cS^{(r)}}\circ \alpha^{r,r'}.
\end{equation}
L'homomorphisme $\halpha^{r,r'}$ induit donc un morphisme de complexes 
\begin{equation}\label{higgs1-chb1i}
\nu^{r,r'}\colon \tmK^\bullet(\hcS^{(r')},p^{r'}d_{\hcS^{(r')}})\rightarrow \tmK^\bullet(\hcS^{(r)},p^rd_{\hcS^{(r)}}).
\end{equation}

D'après \eqref{higgs1-chb1g}, les dérivations $p^rd_{\hcS^{(r)}}$ induisent une $\hoR$-dérivation 
\begin{equation}\label{higgs1-chb1j}
d_{\cS^\dagger}\colon \cS^\dagger\rightarrow \xi^{-1}\tOmega^1_{R/\co_K}\otimes_R\cS^\dagger,
\end{equation}
qui n'est autre que la restriction de $d_{\hcS}$ à $\cS^\dagger$. 
On désigne par $\mK^\bullet(\cS^\dagger,d_{\cS^\dagger})$ 
le complexe de Dolbeault  du module de Higgs $(\cS^\dagger,d_{\cS^\dagger})$.
Comme $\hoR$ est $\co_C$-plat \eqref{higgs1-pur8}, pour tout nombre rationnel $r\geq 0$, on a
\begin{equation}\label{higgs1-chb1k}
\ker(d_{\cS^\dagger})=\ker(d_{\hcS^{(r)}})=\hoR.
\end{equation}

\begin{prop}\label{higgs1-cohbis70}
Pour tous nombres rationnels $r'>r>0$, il existe un nombre rationnel $\alpha\geq 0$ dépendant de $r$ et $r'$, 
mais pas des données \eqref{higgs1-dlog1}, tel que 
\begin{equation}\label{higgs1-cohbis70a}
p^\alpha \nu^{r,r'}\colon \tmK^\bullet(\hcS^{(r')},p^{r'}d_{\hcS^{(r')}})\rightarrow 
\tmK^\bullet(\hcS^{(r)},p^rd_{\hcS^{(r)}})
\end{equation}
soit homotope à $0$ par une homotopie $\hoR$-linéaire. 
\end{prop}

Soient $t_1,\dots,t_d\in P^\gp$ tels que leurs images dans $(P^\gp/\mZ\lambda)\otimes_{\mZ}\mZ_p$ 
forment une $\mZ_p$-base, de sorte que $(d\log(t_i))_{1\leq i\leq d}$ forment une $R$-base 
de $\tOmega^1_{R/\co_K}$ \eqref{higgs1-ext-log13e}. 
Pour tout $1\leq i\leq d$, posons $y_i=\xi^{-1}d\log(t_i)\in \xi^{-1}\tOmega^1_{R/\co_K}\subset \cS$. On désigne par
\begin{equation}
h^{-1}\colon \hcS^{(r')}\otimes_{\mZ_p}\mQ_p\rightarrow \hoR \otimes_{\mZ_p}\mQ_p
\end{equation}
le morphisme $\hoR$-linéaire défini par 
\begin{equation}
h^{-1}(\sum_{\un=(n_1,\dots,n_d)\in \mN^d}a_\un\prod_{1\leq i\leq d}y_i^{n_i})=a_0.
\end{equation}
Pour tout entier $m\geq 0$, il existe un et un unique morphisme $\hoR$-linéaire
\begin{equation}
h^{m}\colon \xi^{-m-1}\tOmega^{m+1}_{R/\co_K}\otimes_R\hcS^{(r')}\otimes_{\mZ_p}\mQ_p\rightarrow 
\xi^{-m}\tOmega^{m}_{R/\co_K}\otimes_R\hcS^{(r)}\otimes_{\mZ_p}\mQ_p
\end{equation}
tel que pour tout $1\leq i_1<\dots<i_{m+1}\leq d$, on ait 
\begin{eqnarray}
\ \ \ \ \ \lefteqn{h^{m}(\sum_{\un=(n_1,\dots,n_d)\in \mN^d}a_\un\prod_{1\leq i\leq d}y_i^{n_i}\otimes 
\xi^{-1}d\log(t_{i_1})\wedge \dots\wedge \xi^{-1}d\log(t_{i_{m+1}}))}\\
&=&\sum_{\un=(n_1,\dots,n_d)\in J_{i_1-1}}
\frac{a_\un}{n_{i_1}+1}\prod_{1\leq i\leq d}y_i^{n_i+\delta_{ii_1}}\otimes 
\xi^{-1}d\log(t_{i_2})\wedge \dots\wedge \xi^{-1}d\log(t_{i_{m+1}}),\nonumber
\end{eqnarray}
où $J_{i_1-1}$ est le sous-ensemble de $\mN^d$ formé des éléments $\un=(n_1,\dots,n_d)$
tels que $n_1=\dots=n_{i_1-1}=0$. 

Soit $\alpha$ un nombre rationnel tel que 
\begin{equation}
\alpha\geq \sup_{x\in \mQ_{\geq 0}}(\log_p(x+1)+(x+1)r-xr'),
\end{equation}
où $\log_p$ est la fonction logarithme de base $p$. Pour tout entier $m\geq 0$, on a clairement
\begin{equation}
p^\alpha h^m(\xi^{-m-1}\tOmega^{m+1}_{R/\co_K}\otimes_R\hcS^{(r')})\subset 
\xi^{-m}\tOmega^{m}_{R/\co_K}\otimes_R\hcS^{(r)}. 
\end{equation}
On vérifie aussitôt que les morphismes $(p^\alpha h^m)_{m\geq -1}$ définissent une homotopie reliant $0$ au morphisme 
$p^\alpha \nu^{r,r'}$.

\begin{cor}\label{higgs1-cohbis7}
Pour tous nombres rationnels $r'>r>0$, le morphisme canonique \eqref{higgs1-chb1i}
\begin{equation}\label{higgs1-cohbis7a}
\tmK^\bullet(\hcS^{(r')},p^{r'}d_{\hcS^{(r')}})\otimes_{\mZ_p}\mQ_p\rightarrow 
\tmK^\bullet(\hcS^{(r)},p^rd_{\hcS^{(r)}})\otimes_{\mZ_p}\mQ_p
\end{equation}
est homotope à $0$ par une homotopie continue. 
\end{cor}

\begin{cor}\label{higgs1-cohbis8}
Le complexe $\mK^\bullet(\cS^\dagger,d_{\cS^\dagger})\otimes_{\mZ_p}\mQ_p$ est une résolution 
de $\hoR[\frac 1 p]$.
\end{cor}

En effet, on a un isomorphisme canonique de complexes 
\begin{equation}
\underset{\underset{r\in \mQ_{>0}}{\longrightarrow}}{\lim}\ 
\mK^\bullet(\hcS^{(r)},p^rd_{\hcS^{(r)}})\otimes_{\mZ_p}\mQ_p
\stackrel{\sim}{\rightarrow} \mK^\bullet(\cS^\dagger,d_{\cS^\dagger})\otimes_{\mZ_p}\mQ_p.
\end{equation}
La proposition résulte donc de \ref{higgs1-cohbis7}.

\subsection{}\label{higgs1-cohbis5}\index{101040@$\varphi_0$}
On désigne par $(\tX_0,\cM_{\tX_0})$ la $(\cA_2(\oS),\cM_{\cA_2(\oS)})$-déformation lisse de 
$(\coX,\cM_\coX)$ définie par la carte $(P,\gamma)$ \eqref{higgs1-ext24ac}, 
par $\cL_0$ le torseur de Higgs-Tate associé \eqref{higgs1-tor2}, 
par $\psi_0\in \cL_0(\hY)$ la section définie par la même carte \eqref{higgs1-ext24b}
et par $\varphi_0=\varphi_{\psi_0}$ l'action de $\Delta$ sur $\cS$ induite par $\psi_0$ \eqref{higgs1-tor202g}. 
D'après \ref{higgs1-cohbis3} ci-dessous, pour tout nombre rationnel $r\geq 0$, la sous-$\hoR$-algèbre $\cS^{(r)}$ de $\cS$
est stable par $\varphi_0$. Sauf mention expresse du contraire, 
on munit $\cS^{(r)}$, $\hcS^{(r)}$ et $\cS^\dagger$ des actions de $\Delta$ induites par $\varphi_0$ . 
Les dérivations $d_{\cS}$ et $d_{\hcS}$ sont $\Delta$-équivariantes
d'après \eqref{higgs1-tor13a}. Il en est donc de même des dérivations $d_{\cS^{(r)}}$, $d_{\hcS^{(r)}}$ et $d_{\cS^\dagger}$ 
compte tenu de \eqref{higgs1-chb1g}.

\begin{prop}\label{higgs1-cohbis3}
Pour tout nombre rationnel $r\geq 0$, la sous-$\hoR$-algèbre $\cS^{(r)}$ de $\cS$
est stable par l'action $\varphi_0$ de $\Delta$ sur $\cS$, et les actions induites  
de $\Delta$ sur $\cS^{(r)}$ et $\hcS^{(r)}$ sont continues pour les topologies $p$-adiques. 
\end{prop}

Reprenons les hypothèses et notations de \ref{higgs1-ext26}. D'après \eqref{higgs1-ext26b}, 
pour tout $g\in \Delta$ et tout $1\leq i\leq d$, on a 
\begin{equation}\label{higgs1-cohbis3a}
\varphi_{0}(g)(y_i)=y_i-\xi^{-1}\chi_i(g).
\end{equation} 
Comme $\xi^{-1}\chi_i(g)\in p^{\frac{1}{p-1}}\co_C$ \eqref{higgs1-ext111}, $\cS^{(r)}$ est stable par $\varphi_0(g)$.
Soit $\zeta$ un générateur de $\mZ_p(1)$. Il existe $a_g\in \mZ_p$ tel que $\chi_i(g)=[\zeta^{a_g}]-1\in \cA_2(\co_\oK)$. 
Par linéarité, on a $\log([\zeta^{a_g}])\in p^{v_p(a_g)}\xi \co_C$, et par suite $\varphi_{0}(g)(y_i)- y_i \in p^{v_p(a_g)}\cS$. 
Pour tout entier $n\geq 0$, l'ensemble des $g\in \Delta$ tels que $v_p(a_g)\geq n$ 
étant un sous-groupe ouvert de $\Delta$, on en déduit que le stabilisateur de 
la classe de $p^ry_i$ dans $\cS^{(r)}/p^n\cS^{(r)}$ est ouvert dans $\Delta$.
La seconde assertion s'ensuit car l'action de $\Delta$ sur $\oR/p^n\oR$ est continue pour la topologie discrète.

\begin{teo}[\cite{tsuji3} 5.3.4] \label{higgs1-cohbis9}
Soit $r$ un nombre rationnel $>0$. Alors~:
\begin{itemize}
\item[{\rm (i)}] Le morphisme canonique 
\begin{equation}
\hRun\otimes_{\mZ_p}\mQ_p\rightarrow (\hcS^{(r)}\otimes_{\mZ_p}\mQ_p)^\Delta
\end{equation}
est un isomorphisme. 
\item[{\rm (ii)}] Pour tout entier $i\geq 1$, 
\begin{equation}
\underset{\underset{r\in \mQ_{>0}}{\longrightarrow}}{\lim}\ 
\rH^i_\cont(\Delta,\hcS^{(r)}\otimes_{\mZ_p}\mQ_p)=0.
\end{equation}
\end{itemize}
\end{teo}

La preuve de ce théorème sera donnée dans \ref{higgs1-cohbis15}. Nous l'avons repoussée vers la fin de cette 
section car elle nécessite l'introduction de notations assez lourdes. On notera que cet énoncé est légèrement 
plus général que celui de Tsuji (\cite{tsuji3} 5.3.4). 

\begin{cor}\label{higgs1-cohbis99}
Pour tout nombre rationnel $r>0$, on a $(\cS^\dagger)^\Delta=(\hcS^{(r)})^\Delta=\hRun$. 
\end{cor}
En effet, $\hcS^{(r)}$ est $\co_C$-plat 
et l'homomorphisme canonique $\hoR\rightarrow \hcS^{(r)}$ est injectif. Il résulte de \ref{higgs1-cohbis9}
qu'on a $(\hcS^{(r)})^\Delta=\hcS^{(r)} \cap \hRun[\frac 1 p]$. Par ailleurs, on a 
$\hcS^{(r)}\cap  \hoR[\frac 1 p]=\hoR$ et $(\hoR)^\Delta=\hRun$ d'après \ref{higgs1-cg12}(i). On en déduit 
que $(\hcS^{(r)})^\Delta=\hRun$ et par suite que $(\cS^\dagger)^\Delta=\hRun$.

\subsection{}\label{higgs1-cohbis10}\index{101042@$\fS$}
On désigne par $\fS$ la sous-$\hRun$-algèbre de $\cS$ définie par \eqref{higgs1-not621}
\begin{equation}\label{higgs1-cohbis10d}
\fS=\rS_{\hRun}(\xi^{-1}\tOmega^1_{R/\co_K}\otimes_R\hRun). 
\end{equation} 
D'après \ref{higgs1-ext26}, l'action $\varphi_0$ de $\Delta$ sur $\cS$ préserve $\fS$, 
et l'action induite sur $\fS$ se factorise à travers $\Delta_{p^\infty}$. On note encore $\varphi_0$ l'action 
de $\Delta_{p^\infty}$ sur $\fS$ ainsi définie. 
Calquant la preuve de \ref{higgs1-cohbis3}, on montre que cette dernière est continue pour la topologie $p$-adique sur $\fS$.
On peut aussi déduire cette propriété directement de \ref{higgs1-cohbis3} en observant que 
pour tout entier $n\geq 0$, l'homomorphisme canonique $\fS/p^n\fS\rightarrow \cS/p^n\cS$ est injectif (cf. la preuve de \ref{higgs1-pur8}).  
On pose
\begin{eqnarray}
\cS_\infty&=&\fS\otimes_{\hRun}\hRi,\label{higgs1-cohbis10e}\\
\cS_{p^\infty}&=&\fS\otimes_{\hRun}\hRpi.\label{higgs1-cohbis10f}
\end{eqnarray}
L'action $\varphi_0$ de $\Delta_{p^\infty}$ sur $\fS$ induit des actions de $\Delta_{p^\infty}$ sur $\cS_{p^\infty}$ et
de  $\Delta_{\infty}$ sur $\cS_{\infty}$.

Pour tout nombre rationnel $r\geq 0$, on désigne par $\cS^{(r)}_\infty$ 
la sous-$\hRi$-algèbre de $\cS_\infty$ définie par 
\begin{equation}\label{higgs1-cohbis10b}
\cS^{(r)}_\infty=\rS_{\hRi}(p^r\xi^{-1}\tOmega^1_{R/\co_K}\otimes_R\hRi),
\end{equation}
et par $\hcS^{(r)}_\infty$ son séparé complété $p$-adique.
Compte tenu de \ref{higgs1-pur8} et de sa preuve, $\cS^{(r)}_\infty$ et $\hcS^{(r)}_\infty$ sont $\co_C$-plats.
Pour tout nombre rationnel $r'\geq r$, on a un homomorphisme injectif canonique 
$\cS^{(r')}_\infty\rightarrow \cS^{(r)}_\infty$. On vérifie aussitôt que l'homomorphisme induit 
$\hcS^{(r')}_\infty\rightarrow \hcS^{(r)}_\infty$ est injectif. La preuve de \ref{higgs1-cohbis3} montre que 
$\cS^{(r)}_\infty$ est stable par l'action de $\Delta_\infty$ sur $\cS_\infty=\cS^{(0)}_\infty$, 
et les actions induites de $\Delta_\infty$ sur $\cS^{(r)}_\infty$ et $\hcS^{(r)}_\infty$ 
sont continues pour les topologies $p$-adiques. Sauf mention expresse du contraire, 
on munit $\cS^{(r)}_\infty$ et $\hcS^{(r)}_\infty$ de ces actions et des topologies $p$-adiques.

On désigne par $\cS^{(r)}_{p^\infty}$ la sous-$\hRpi$-algèbre de $\cS_{p^\infty}$ définie par
\begin{equation}\label{higgs1-cohbis10bb}
\cS^{(r)}_{p^\infty}=\rS_{\hRpi}(p^r\xi^{-1}\tOmega^1_{R/\co_K}\otimes_R\hRpi),
\end{equation}
et par $\hcS^{(r)}_{p^\infty}$ son séparé complété $p$-adique.
L'algèbre $\cS^{(r)}_{p^\infty}$ vérifie des propriétés analogues à celles vérifiées par $\cS^{(r)}_{\infty}$. 
En particulier, $\cS^{(r)}_{p^\infty}$ est stable par l'action de $\Delta_{p^\infty}$ sur 
$\cS_{p^\infty}=\cS^{(0)}_{p^\infty}$, et les actions induites de $\Delta_{p^\infty}$
sur $\cS^{(r)}_{p^\infty}$ et $\hcS^{(r)}_{p^\infty}$ sont continues pour les topologies $p$-adiques.  
Sauf mention expresse du contraire, on munit $\cS^{(r)}_{p^\infty}$ et $\hcS^{(r)}_{p^\infty}$ de ces actions
et des topologies $p$-adiques.

\begin{lem}\label{higgs1-cohbis11}
Pour tout nombre rationnel $r>0$ et tout entier $i\geq 0$, le morphisme canonique 
\begin{equation}\label{higgs1-cohbis11aa}
\rH^i_\cont(\Delta_{p^\infty},\hcS^{(r)}_{p^\infty})\rightarrow \rH^i_\cont(\Delta_\infty,\hcS^{(r)}_\infty)
\end{equation}
est un isomorphisme, et le morphisme canonique 
\begin{equation}\label{higgs1-cohbis11a}
\rH^i_\cont(\Delta_\infty,\hcS^{(r)}_\infty)\rightarrow \rH^i_\cont(\Delta,\hcS^{(r)})
\end{equation}
est un presque-isomorphisme.
\end{lem} 

Soit $n$ un entier $\geq 0$. L'homomorphisme canonique
\begin{equation}\label{higgs1-cohbis11bb}
\cS^{(r)}_{p^\infty}/p^n\cS^{(r)}_{p^\infty}\rightarrow (\cS^{(r)}_\infty/p^n\cS^{(r)}_\infty)^{\Sigma_0}
\end{equation}
est un isomorphisme en vertu de \ref{higgs1-cg07}. Il en résulte par \eqref{higgs1-gal3a} que le morphisme
canonique 
\begin{equation}\label{higgs1-cohbis11bc}
\rH^i(\Delta_{p^\infty},\cS^{(r)}_{p^\infty}/p^n\cS^{(r)}_{p^\infty})\rightarrow \rH^i(\Delta_\infty,\cS^{(r)}_\infty/p^n\cS^{(r)}_\infty)
\end{equation}
est un isomorphisme. On en déduit par \eqref{higgs1-limproj2c} et \eqref{higgs1-limproj2d} que le morphisme \eqref{higgs1-cohbis11aa}
est un isomorphisme. 

D'autre part, l'homomorphisme canonique
\begin{equation}\label{higgs1-cohbis11b}
\cS^{(r)}_\infty/p^n\cS^{(r)}_\infty\rightarrow (\cS^{(r)}/p^n\cS^{(r)})^\Sigma
\end{equation}
est un presque-isomorphisme en vertu de \ref{higgs1-pur12}. 
On en déduit par \ref{higgs1-pur6} que le morphisme canonique 
\begin{equation}\label{higgs1-cohbis11d}
\psi_n\colon \rH^i(\Delta_\infty,\cS^{(r)}_\infty/p^n\cS^{(r)}_\infty)\rightarrow \rH^i(\Delta,\cS^{(r)}/p^n\cS^{(r)})
\end{equation}
est un presque-isomorphisme. Notons $A_n$ (resp. $C_n$) le noyau (resp. conoyau) de $\psi_n$. 
Alors les $\co_\oK$-modules
\[
\underset{\underset{n\geq 0}{\longleftarrow}}{\lim}\ A_n, \ \ \ \underset{\underset{n\geq 0}{\longleftarrow}}{\lim}\ C_n, 
\ \ \ \rR^1\underset{\underset{n\geq 0}{\longleftarrow}}{\lim}\ A_n, \ \ \ 
\rR^1\underset{\underset{n\geq 0}{\longleftarrow}}{\lim}\ C_n
\]
sont presque nuls en vertu de (\cite{gr} 2.4.2(ii)). Par suite, les morphismes 
\[
\underset{\underset{n\geq 0}{\longleftarrow}}{\lim}\ \psi_n
\ \ \ {\rm et} \ \ \ \rR^1\underset{\underset{n\geq 0}{\longleftarrow}}{\lim}\ \psi_n
\]
sont des presque-isomorphismes. On en déduit par \eqref{higgs1-limproj2c} et \eqref{higgs1-limproj2d} que le morphisme 
\eqref{higgs1-cohbis11a} est un isomorphisme. 

\subsection{}\label{higgs1-cohbis12}
Soient $t_1,\dots,t_d\in P^\gp$ tels que leurs images dans $(P^\gp/\mZ\lambda)\otimes_\mZ\mZ_p$ 
forment une $\mZ_p$-base, 
$(\chi_{t_i})_{1\leq i\leq d}$ leurs images dans $\Hom(\Delta_{p^\infty},\mZ_p(1))$ \eqref{higgs1-ext-log14e},
$\zeta$ une $\mZ_p$-base de $\mZ_p(1)$. 
Les $(\chi_{t_i})_{1\leq i\leq d}$ forment une $\mZ_p$-base de $\Hom(\Delta_{p^\infty},\mZ_p(1))$ \eqref{higgs1-ext-log14gh}.
Il existe donc une unique $\mZ_p$-base $(\gamma_i)_{1\leq i\leq d}$  de $\Delta_{p^\infty}$ 
telle que  $\chi_{t_i}(\gamma_j)=\delta_{ij}\zeta$ pour tous $1\leq i,j\leq d$  \eqref{higgs1-gal2}. 
Pour tout entier $0\leq i\leq d$, on désigne par ${_i\Xi}_{p^\infty}$ le sous-groupe de \eqref{higgs1-cg3km}
\begin{equation}
{\Xi}_{p^\infty}=\Hom(\Delta_{p^\infty},\mu_{p^\infty}(\co_\oK))
\end{equation}
formé des homomorphismes $\nu\colon \Delta_{p^\infty}\rightarrow \mu_{p^\infty}(\co_\oK)$ tels que 
$\nu(\gamma_j)=1$ pour tout $1\leq j\leq i$.

Les $(d\log(t_i))_{1\leq i\leq d}$ forment une $R$-base de $\tOmega^1_{R/\co_K}$ \eqref{higgs1-ext-log13e}. 
Pour tout $1\leq i\leq d$ et tout $\un=(n_1,\dots,n_d)\in \mN^d$, 
posons $y_i=\xi^{-1}d\log(t_i)\in \xi^{-1}\tOmega^1_{R/\co_K}\subset \cS_{p^\infty}$,
$|\un|=\sum_{i=1}^d n_i$ et $\uy^{\un}=\prod_{i=1}^dy_i^{n_i}\in \cS_{p^\infty}$.
On observera que $R_{p^\infty}$ est séparé pour la topologie $p$-adique 
et s'identifie donc à un sous-anneau de $\hRpi$~; cela résulte par exemple de \eqref{higgs1-cg16ae}, \eqref{higgs1-cg16i}
et du fait que $\Spec(R_1)$ est un ouvert de $\Spec(C_1)$.
Pour tout nombre rationnel $r>0$, tout entier $0\leq i\leq d$ et tout $\nu\in {_i\Xi}_{p^\infty}$, 
tenant compte de \ref{higgs1-cg16} et avec les mêmes notations, 
on désigne par ${_i\cS}^{(r)}_{p^\infty}(\nu)$ et ${_i\cS}^{(r)}_{p^\infty}$
les sous-$R_1$-modules de $\cS^{(r)}_{p^\infty}$ définis par 
\begin{eqnarray}
{_i\cS}^{(r)}_{p^\infty}(\nu)&=&\bigoplus_{\un\in J_i}p^{r|\un|} R_{p^\infty}^{(\nu)} \uy^\un,\\
{_i\cS}^{(r)}_{p^\infty}&=&\bigoplus_{\nu'\in {_i\Xi}_{p^\infty}}{_i\cS}^{(r)}_{p^\infty}(\nu'),
\end{eqnarray}
où $J_{i}$ est le sous-ensemble de $\mN^d$ formé des éléments $\un=(n_1,\dots,n_d)$
tels que $n_1=\dots=n_{i}=0$.  
On note ${_i\hcS}^{(r)}_{p^\infty}$ le séparé complété $p$-adique de ${_i\cS}^{(r)}_{p^\infty}$. 
La topologie $p$-adique de $R_{p^\infty}$ étant induite par la topologie $p$-adique de $\hRpi$, 
on déduit facilement de \eqref{higgs1-cg16a} 
que la topologie $p$-adique de ${_i\cS}^{(r)}_{p^\infty}$ est induite par la topologie $p$-adique de $\cS^{(r)}_{p^\infty}$.
Par suite, ${_i\hcS}^{(r)}_{p^\infty}$ est l'adhérence de ${_i\cS}^{(r)}_{p^\infty}$ dans $\hcS^{(r)}_{p^\infty}$. 
Il résulte de \ref{higgs1-cg16} et \eqref{higgs1-cohbis3a} que pour tout $1\leq j\leq d$ et tout $\nu\in {_i\Xi}_{p^\infty}$, 
$\gamma_{j}$ préserve ${_{i}\cS}^{(r)}_{p^\infty}(\nu)$ et donc aussi 
${_i\cS}^{(r)}_{p^\infty}$. Si $1\leq j\leq i$, $\gamma_{j}$ fixe ${_{i}\cS}^{(r)}_{p^\infty}(\nu)$ et  
${_i\cS}^{(r)}_{p^\infty}$.

\begin{lem}\label{higgs1-cohbis14} 
Les hypothèses étant celles de \eqref{higgs1-cohbis12}, soient, de plus, $i$ un entier tel que $1\leq i\leq d$,
$r,r'$ deux nombres rationnels tels que $r'>r>0$. Alors~:
\begin{itemize}
\item[{\rm (i)}] On a ${_0\hcS}^{(r)}_{p^\infty}=\hcS^{(r)}_{p^\infty}$ et 
${_d\hcS}^{(r)}_{p^\infty}=\hRun$.
\item[{\rm (ii)}] Le noyau du morphisme 
\begin{equation}\label{higgs1-cohbis14a}
\gamma_i-\id\colon {_{(i-1)}\hcS}^{(r)}_{p^\infty}\rightarrow {_{(i-1)}\hcS}^{(r)}_{p^\infty}
\end{equation}
est égal à ${_i\hcS}^{(r)}_{p^\infty}$.
\item[{\rm (iii)}] Il existe un entier $\alpha\geq 0$, dépendant de $r$ et $r'$ 
mais pas des données \eqref{higgs1-dlog1}, tel que l'on ait
\begin{equation}\label{higgs1-cohbis14b}
p^\alpha\cdot {_{(i-1)}\hcS}^{(r')}_{p^\infty}\subset (\gamma_i-\id)({_{(i-1)}\hcS}^{(r)}_{p^\infty}).
\end{equation}
\end{itemize}
\end{lem}

(i) Comme $R_{p^\infty}^{(1)}=R_1$ \eqref{higgs1-cg16b}, on a ${_d\hcS}^{(r)}_{p^\infty}=\hRun$. D'autre part, 
on a  
\begin{equation}\label{higgs1-cohbis14c}
{_0\cS}^{(r)}_{p^\infty}=\rS_{R_{p^\infty}}(p^r\xi^{-1}\tOmega^1_{R/\co_K}\otimes_RR_{p^\infty}),
\end{equation}
ce qui implique que ${_0\hcS}^{(r)}_{p^\infty}=\hcS^{(r)}_{p^\infty}$.

(ii) Pour tout entier $n\geq 0$, on a clairement 
\begin{equation}\label{higgs1-cohbis14d}
p^n \cdot{_{i}\cS}^{(r)}_{p^\infty}={_{i}\cS}^{(r)}_{p^\infty}\cap (p^n\cdot{_{(i-1)}\cS}^{(r)}_{p^\infty}).
\end{equation} 
Il suffit donc de montrer que la suite 
\begin{equation}\label{higgs1-cohbis14e}
\xymatrix{
0\ar[r]&{{_{i}\cS}^{(r)}_{p^\infty}}\ar[r]&{{_{(i-1)}\cS}^{(r)}_{p^\infty}}\ar[r]^{\gamma_i-\id}&{{_{(i-1)}\cS}^{(r)}_{p^\infty}}}
\end{equation}
est exacte. Soient $\nu\in {_{(i-1)}\Xi}_{p^\infty}$ et 
\begin{equation}\label{higgs1-cohbis14f}
z=\sum_{\un\in J_{i-1}}p^{r|\un|}a_\un \uy^\un\in {_{(i-1)}\cS}^{(r)}_{p^\infty}(\nu).
\end{equation}
On a alors \eqref{higgs1-cohbis3a}
\begin{equation}\label{higgs1-cohbis14g}
(\gamma_i-\id)(z)=\sum_{\un\in J_{i-1}}p^{r|\un|}b_\un \uy^\un\in {_{(i-1)}\cS}^{(r)}_{p^\infty}(\nu),
\end{equation}
où pour tout $\un=(n_1,\dots,n_d)\in J_{i-1}$, 
\begin{equation}\label{higgs1-cohbis14h}
b_\un=(\nu(\gamma_i)-1)a_\un+\sum_{\um=(m_1,\dots,m_d)
\in J_{i-1}(\un)} p^{r(m_i-n_i)}\binom{m_i}{n_i} \nu(\gamma_i)a_{\um}w^{m_i-n_i},
\end{equation}
$J_{i-1}(\un)$ désigne le sous-ensemble de $J_{i-1}$ formé des éléments $\um=(m_1,\dots,m_d)$ tels que 
$m_j=n_j$ pour $j\not=i$ et $m_i>n_i$, 
et $w=\xi^{-1}\log([\zeta])$ est un élément de valuation $\frac{1}{p-1}$ de $\co_C$ \eqref{higgs1-ext111}.

Supposons que $\gamma_i(z)=z$ et $\nu(\gamma_i)\not=1$.  Alors $v(\nu(\gamma_i)-1)\leq \frac{1}{p-1}$ et
pour tout  $\un=(n_1,\dots,n_d)\in J_{i-1}$, on a 
\[
a_\un=-(\nu(\gamma_i)-1)^{-1}\sum_{\um=(m_1,\dots,m_d)
\in J_{i-1}(\un)} p^{r(m_i-n_i)}\binom{m_i}{n_i} \nu(\gamma_i)a_{\um}w^{m_i-n_i}.
\]
On en déduit que pour tout $\alpha\in \mN$ et tout $\un\in J_{i-1}$,
on a $a_\un\in p^{r\alpha} R_{p^\infty}^{(\nu)}$ (on le prouve par récurrence sur $\alpha$);
donc $z=0$ car $R_{p^\infty}^{(\nu)}$ est séparé pour la topologie $p$-adique. 

Supposons que $\gamma_i(z)=z$ et $\nu(\gamma_i)=1$. Alors pour tout $\un=(n_1,\dots,n_d)\in J_{i-1}$, si on pose 
$\un'=(n'_1,\dots,n'_d)\in J_{i-1}(\un)$ avec $n'_i=n_i+1$, on a 
\[
(n_i+1)!a_{\un'}=-\sum_{\um=(m_1,\dots,m_d)
\in J_{i-1}(\un')} p^{r(m_i-n_i-1)}m_i! a_{\um}\frac{w^{m_i-n_i-1}}{(m_i-n_i)!}.
\]
On a $w^{m-1}/m!\in \co_C$ pour tout entier $m\geq 1$. On en déduit que 
pour tout $\alpha\in \mN$ et tout $\un=(n_1,\dots,n_d)\in J_{i-1}$ tel que $n_i\geq 1$, 
on a $n_i!a_{\un}\in p^{r\alpha} R_{p^\infty}^{(\nu)}$ (on le prouve par récurrence sur $\alpha$); donc $a_\un=0$. 
Par suite $z \in {_{i}\cS}^{(r)}_{p^\infty}(\nu)$, ce qui achève la preuve de l'exactitude de la suite \eqref{higgs1-cohbis14e}.

(iii) Il suffit de montrer qu'il existe un entier $\alpha\geq 0$ ne dépendant que de $r$ et $r'$
tel que pour tout $\nu\in {_{(i-1)}\Xi}_{p^\infty}$, on ait
\begin{equation}\label{higgs1-cohbis14aa}
p^\alpha({_{(i-1)}\cS}^{(r')}_{p^\infty}(\nu))\subset (\gamma_i-\id)({_{(i-1)}\hcS}^{(r)}_{p^\infty}).
\end{equation}
En effet, si l'on pose $M=(\gamma_i-\id)({_{(i-1)}\hcS}^{(r)}_{p^\infty})$, on en déduira 
par complétion $p$-adique un diagramme commutatif 
\begin{equation}\label{higgs1-cohbis14ab}
\xymatrix{
{p^\alpha({_{(i-1)}\hcS}^{(r')}_{p^\infty}})\ar[r]&\hM\ar[r]&{{_{(i-1)}\hcS}^{(r)}_{p^\infty}}\\
&{{_{(i-1)}\hcS}^{(r)}_{p^\infty}}\ar[ru]_{\gamma_i-\id}\ar@{->>}[u]&}
\end{equation}
où la flèche verticale sera surjective en vertu (\cite{egr1} 1.8.5). 

Supposons $\nu(\gamma_i)\not=1$. D'après \eqref{higgs1-cohbis14g}, on a 
\begin{equation}\label{higgs1-cohbis14ac}
(\nu(\gamma_i)-1)({_{(i-1)}\cS}^{(r)}_{p^\infty}(\nu))\subset  (\gamma_i-\id)({_{(i-1)}\cS}^{(r)}_{p^\infty}(\nu))
+(\nu(\gamma_i)-1)p^r ({_{(i-1)}\cS}^{(r)}_{p^\infty}(\nu)).
\end{equation}
On en déduit que 
\begin{equation}\label{higgs1-cohbis14ad}
(\nu(\gamma_i)-1)({_{(i-1)}\cS}^{(r)}_{p^\infty}(\nu))\subset  (\gamma_i-\id)({_{(i-1)}\hcS}^{(r)}_{p^\infty}).
\end{equation}
On peut donc prendre $\alpha=1$ car $v(\nu(\gamma_i)-1)\leq \frac{1}{p-1}$. 

Supposons que $\nu(\gamma_i)=1$, de sorte que $\nu\in {_{i}\Xi}_{p^\infty}$.
Notons $(R_{p^\infty}^{(\nu)})^\wedge$ et $( {_{(i-1)}\cS}^{(r)}_{p^\infty}(\nu))^\wedge$
les séparés complétés $p$-adiques de $R_{p^\infty}^{(\nu)}$ et ${_{(i-1)}\cS}^{(r)}_{p^\infty}(\nu)$, respectivement. 
On observera que $( {_{(i-1)}\cS}^{(r)}_{p^\infty}(\nu))^\wedge$ s'identifie à un sous-$\hRun$-module de 
${_{(i-1)}\hcS}^{(r)}_{p^\infty}$. 
Tout élément $z$ de $( {_{(i-1)}\cS}^{(r)}_{p^\infty}(\nu))^\wedge$ peut s'écrire comme somme d'une série 
\begin{equation}
z=\sum_{\un\in J_{i-1}}p^{r|\un|}a_\un \uy^\un
\end{equation}
où $a_\un \in (R_{p^\infty}^{(\nu)})^\wedge$ et $a_\un$ tend vers $0$ quand $|\un|$ tend vers l'infini. 
Par suite, $(\gamma_i-\id)(z)$ est encore donné par la formule \eqref{higgs1-cohbis14g}. 
Il suffit donc de montrer qu'il existe un entier $\alpha\geq 0$ ne dépendant que de $r$ et $r'$ tel que pour tout
\begin{equation}\label{higgs1-cohbis14ae}
\sum_{\un\in J_{i-1}}p^{r'|\un|}b_\un \uy^\un\in {_{(i-1)}\cS}^{(r')}_{p^\infty}(\nu),
\end{equation}
le système d'équations linéaires définies, pour $\un=(n_1,\dots,n_d)\in J_{i-1}$, par 
\begin{equation}\label{higgs1-cohbis14af}
p^\alpha p^{(r'-r)|\un|}n_i!b_\un=\sum_{\um=(m_1,\dots,m_d)
\in J_{i-1}(\un)} p^{r(m_i-n_i)}m_i! a_{\um}\frac{w^{m_i-n_i}}{(m_i-n_i)!},
\end{equation}
admette une solution $a_\um\in (R_{p^\infty}^{(\nu)})^\wedge$ pour $\um\in J_{i-1}$ telle que 
$a_\um$ tend vers $0$ quand $|\um|$ tend vers l'infini. Pour $\un=(n_1,\dots,n_d) \in J_{i-1}$, posons
\begin{eqnarray}
a'_\un&=&n_i!p^{-\frac{n_i}{p-1}}a_\un,\label{higgs1-cohbis14ba}\\
b'_\un&=&p^\alpha p^{(r'-r)|\un|} n_i!p^{-\frac{n_i}{p-1}}b_\un,\label{higgs1-cohbis14bb}
\end{eqnarray}
de sorte que l'équation \eqref{higgs1-cohbis14af} devient 
\begin{equation}\label{higgs1-cohbis14ag}
b'_\un=p^{r+\frac{1}{p-1}}w\sum_{\um=(m_1,\dots,m_d)
\in J_{i-1}(\un)} p^{(r+\frac{1}{p-1})(m_i-n_i-1)}a'_{\um}\frac{w^{(m_i-n_i-1)}}{(m_i-n_i)!}. 
\end{equation}
Considérons l'endomorphisme $\hRun$-linéaire $\Phi$ de $(\oplus_{\un\in J_{i-1}}R_{p^\infty}^{(\nu)})^\wedge$ défini,  
pour une suite $(x_\un)_{\un\in J_{i-1}}$ d'éléments de $(R_{p^\infty}^{(\nu)})^\wedge$ tendant vers $0$ quand
$|\un|$ tend vers l'infini, par
\begin{equation}
\Phi(\sum_{\un\in J_{i-1}}x_\un)= \sum_{\un\in J_{i-1}}z_\un,
\end{equation}
où pour $\un=(n_1,\dots,n_d)\in J_{i-1}$, 
\begin{equation}
z_\un=\sum_{\um=(m_1,\dots,m_d)
\in \{\un\}\cup J_{i-1}(\un)} p^{(r+\frac{1}{p-1})(m_i-n_i)}x_{\um}\frac{w^{(m_i-n_i)}}{(m_i-n_i+1)!}.
\end{equation}
Comme $\Phi$ est congru à l'identité modulo $p^{r+\frac{1}{p-1}}$, il est surjectif en vertu (\cite{egr1} 1.8.5). 
Par suite, pour toute suite $b'_\un\in p^{r+\frac{1}{p-1}}w (R_{p^\infty}^{(\nu)})^\wedge$ pour $\un \in J_{i-1}$ 
tendant vers $0$ quand $|\un|$ tend vers l'infini, 
l'équation \eqref{higgs1-cohbis14ag} admet une solution $a'_\um\in (R_{p^\infty}^{(\nu)})^\wedge$ pour $\um\in J_{i-1}$ 
tendant vers $0$ quand $|\um|$ tend vers l'infini. 
D'autre part, $v(w)=\frac{1}{p-1}$ et il existe un entier $\alpha\geq 0$ tel que pour tout $n\in \mN$, on ait 
\begin{equation}
(r'-r)n+v(n!)-\frac{n}{p-1}+\alpha\geq r+\frac{2}{p-1}. 
\end{equation}
L'assertion recherchée s'ensuit en prenant pour $b'_\un$ pour
$\un\in J_{i-1}$ les éléments définis par \eqref{higgs1-cohbis14bb} (qui sont en fait nuls sauf un nombre fini).

\subsection{}\label{higgs1-cohbis141}
Conservons les hypothèses de \ref{higgs1-cohbis12}. Pour tout nombre rationnel $r>0$, 
on définit par récurrence, pour tout entier $0\leq i\leq d$, un complexe $\mK_i^{(r),\bullet}$ de $\hRpi$-représentations
continues de $\Delta_{p^\infty}$ en posant $\mK_0^{(r),\bullet}=\hcS^{(r)}_{p^\infty}[0]$ 
et pour tout $1\leq i\leq n$, $\mK_i^{(r),\bullet}$ est la fibre du morphisme 
\begin{equation}
\gamma_i-\id\colon \mK_{i-1}^{(r),\bullet}\rightarrow \mK_{i-1}^{(r),\bullet}.
\end{equation} 
Il résulte de \ref{higgs1-cg01} et \eqref{higgs1-koszul2f} 
qu'on a un isomorphisme canonique dans $\bD^+(\bMod(\hRun))$
\begin{equation}\label{higgs1-cohbis141a}
\rC_{\cont}^\bullet(\Delta_{p^\infty},\hcS^{(r)}_{p^\infty})\stackrel{\sim}{\rightarrow}\mK_d^{(r),\bullet}.
\end{equation}
Pour tous nombres rationnels $r'>r>0$, l'homomorphisme canonique $\hcS^{(r')}_{p^\infty}\rightarrow \hcS^{(r)}_{p^\infty}$
induit pour tout entier $0\leq i\leq d$ un morphisme $\mK_i^{(r'),\bullet}\rightarrow \mK_i^{(r),\bullet}$
de complexes de $\hRpi$-représentations continues de $\Delta_{p^\infty}$.

\begin{prop}\label{higgs1-cohbis142}
Les hypothèses étant celles de \eqref{higgs1-cohbis12}, soient, de plus, $i$ un entier tel que $0\leq i\leq d$,
$r, r'$ deux nombres rationnels tels que $r'>r>0$. Alors~:
\begin{itemize}
\item[{\rm (i)}] On a un isomorphisme canonique $\hRun$-linéaire et $\Delta_{p^\infty}$-équivariant 
\begin{equation}\label{higgs1-cohbis142a}
{_i\hcS}^{(r)}_{p^\infty}\stackrel{\sim}{\rightarrow}\rH^0(\mK_i^{(r),\bullet}).
\end{equation}
\item[{\rm (ii)}] Il existe un entier $\alpha_i\geq 0$, dépendant de $r$, $r'$ et $i$,
mais pas des données \eqref{higgs1-dlog1}, tel que 
pour tout entier $j\geq 1$, le morphisme canonique 
\begin{equation}\label{higgs1-cohbis142b}
\rH^j(\mK_i^{(r'),\bullet}) \rightarrow\rH^j(\mK_i^{(r),\bullet})
\end{equation}
soit annulé par $p^{\alpha_i}$. 
\end{itemize}
\end{prop}

Procédons par récurrence sur $i$. La proposition est immédiate pour $i=0$ d'après \ref{higgs1-cohbis14}(i). 
Supposons $i\geq 1$ et la proposition démontrée pour $i-1$. Le triangle distingué 
\begin{equation}\label{higgs1-cohbis142c}
\xymatrix{
{\mK_i^{(r),\bullet}}\ar[r]&{\mK_{i-1}^{(r),\bullet}}\ar[r]^-(0.4){\gamma_i-\id}& 
{\mK_{i-1}^{(r),\bullet}}\ar[r]^-(0.5){+1}&}
\end{equation}
et l'hypothèse de récurrence induisent une suite exacte
\begin{equation}\label{higgs1-cohbis142d}
\xymatrix{
0\ar[r]&{\rH^0(\mK_i^{(r),\bullet})}\ar[r]&{{_{i-1}\hcS}^{(r)}_{p^\infty}}\ar[r]^-(0.4){\gamma_i-\id}& 
{{_{i-1}\hcS}^{(r)}_{p^\infty}}},
\end{equation}
qui implique proposition (i) en vertu de \ref{higgs1-cohbis14}(ii).

Pour tout entier $j\geq 1$, le triangle distingué \eqref{higgs1-cohbis142c} induit une suite exacte de $\hRun$-modules
\begin{equation}\label{higgs1-cohbis142e}
0\rightarrow C_j^{(r)}\rightarrow \rH^j(\mK_i^{(r),\bullet}) \rightarrow D_j^{(r)} \rightarrow 0,
\end{equation}
où $C_j^{(r)}$ est un quotient de $\rH^{j-1}(\mK_{i-1}^{(r),\bullet})$ et 
$D_j^{(r)}$ est un sous-module de $\rH^j(\mK_{i-1}^{(r),\bullet})$. De plus, d'après l'hypothèse de récurrence, 
on a un isomorphisme canonique 
\begin{equation}\label{higgs1-cohbis142f}
C_1^{(r)}\stackrel{\sim}{\rightarrow}{_{i-1}\hcS}^{(r)}_{p^\infty}/(\gamma_i-1)({_{i-1}\hcS}^{(r)}_{p^\infty}).
\end{equation}
Les morphismes canoniques $\mK_{i-1}^{(r'),\bullet}\rightarrow \mK_{i-1}^{(r),\bullet}$
et $\mK_i^{(r'),\bullet}\rightarrow \mK_i^{(r),\bullet}$ induisent des morphismes $C_j^{(r')}\rightarrow C_j^{(r)}$ et 
$D_j^{(r')}\rightarrow D_j^{(r)}$ qui s'insèrent dans un diagramme commutatif 
\begin{equation}\label{higgs1-cohbis142g}
\xymatrix{
0\ar[r]&{C_j^{(r')}}\ar[r]\ar[d]&{\rH^j(\mK_i^{(r'),\bullet})}\ar[r]\ar[d]&{D_j^{(r')}}\ar[r]\ar[d]&0\\
0\ar[r]&{C_j^{(r)}}\ar[r]&{\rH^j(\mK_i^{(r),\bullet})}\ar[r]&{D_j^{(r)}}\ar[r]&0}
\end{equation}
où la flèche verticale au centre est le morphisme canonique. 

Posons $r''=(r+r')/2$. D'après l'hypothèse de récurrence, il existe un entier 
$\alpha_{i-1}\geq 0$, dépendant seulement de $r$, $r'$ et $i-1$, 
tel que pour tout entier $j\geq 1$, le morphisme $D_j^{(r')}\rightarrow D_j^{(r'')}$ soit annulé par $p^{\alpha_{i-1}}$. 
D'autre part, compte tenu de l'hypothèse de récurrence et en vertu de \eqref{higgs1-cohbis142f} et \ref{higgs1-cohbis14}(ii), il existe un 
entier $\alpha'_{i-1}\geq 0$, dépendant seulement de $r$, $r'$ et $i-1$, tel que pour tout entier $j\geq 1$, 
le morphisme $C_j^{(r'')}\rightarrow C_j^{(r)}$ soit annulé par $p^{\alpha'_{i-1}}$. La proposition (ii) s'ensuit en prenant 
$\alpha_i=\alpha_{i-1}+\alpha'_{i-1}$. 

\begin{cor}\label{higgs1-cohbis143}
Soient $r,r'$ deux nombres rationnels tels que $r'>r>0$. Alors~:
\begin{itemize}
\item[{\rm (i)}] L'homomorphisme canonique 
\begin{equation}\label{higgs1-cohbis143c}
\hRun \rightarrow (\hcS^{(r)}_{p^\infty})^{\Delta_{p^\infty}}
\end{equation}
est un isomorphisme.
\item[{\rm (ii)}] Il existe un entier $\alpha\geq 0$, dépendant de $r$, $r'$ et $d$, 
mais pas des données \eqref{higgs1-dlog1}, tel que 
pour tout entier $j\geq 1$, le morphisme canonique 
\begin{equation}\label{higgs1-cohbis143d}
\rH^j_\cont(\Delta_{p^\infty},\hcS^{(r')}_{p^\infty})\rightarrow \rH^j_\cont(\Delta_{p^\infty},\hcS^{(r)}_{p^\infty})
\end{equation}
soit annulé par $p^\alpha$. 
\end{itemize}
\end{cor}

Cela résulte de \eqref{higgs1-cohbis141a} et \ref{higgs1-cohbis142}.

\begin{cor}\label{higgs1-cohbis144}
Soient $r,r'$ deux nombres rationnels tels que $r'>r>0$. Alors~:
\begin{itemize}
\item[{\rm (i)}] L'homomorphisme canonique 
\begin{equation}\label{higgs1-cohbis144a}
\hRun \rightarrow (\hcS^{(r)})^\Delta
\end{equation}
est un presque-isomorphisme.
\item[{\rm (ii)}] Il existe un entier $\alpha\geq 0$, dépendant de $r$, $r'$ et $d$, 
mais pas des données \eqref{higgs1-dlog1}, tel que 
pour tout entier $j\geq 1$, le morphisme canonique 
\begin{equation}\label{higgs1-cohbis144b}
\rH^j_\cont(\Delta,\hcS^{(r')})\rightarrow \rH^j_\cont(\Delta,\hcS^{(r)})
\end{equation}
soit annulé par $p^\alpha$. 
\end{itemize}
\end{cor}

Cela résulte de \ref{higgs1-cohbis11} et \ref{higgs1-cohbis143}.

\subsection{}\label{higgs1-cohbis15}
Nous pouvons maintenant démontrer le théorème \ref{higgs1-cohbis9}. 
Comme $\Delta$ est compact, pour tout nombre rationnel $r>0$, le morphisme canonique \eqref{higgs1-not78}
\begin{equation}\label{higgs1-cohbis15a}
\rC^\bullet_\cont(\Delta,\hcS^{(r)})\otimes_{\mZ_p}\mQ_p\rightarrow 
\rC^\bullet_\cont(\Delta,\hcS^{(r)}\otimes_{\mZ_p}\mQ_p)
\end{equation}
est un isomorphisme. Le théorème \ref{higgs1-cohbis9} résulte donc de \ref{higgs1-cohbis144}.

\begin{prop}\label{higgs1-cohbis145}
Soient $r,r'$ deux nombres rationnels tels que $r'>r>0$. Alors~:
\begin{itemize}
\item[{\rm (i)}] Pour tout entier $n\geq 1$, le morphisme canonique 
\begin{equation}\label{higgs1-cohbis145a}
R_1/p^nR_1\rightarrow (\cS^{(r)}/p^n\cS^{(r)})^{\Delta}
\end{equation}
est presque-injectif. Notons $\cH^{(r)}_n$ son conoyau. 
\item[{\rm (ii)}] Il existe un entier $\alpha\geq 0$, dépendant de $r$, $r'$ et $d$, 
mais pas des données \eqref{higgs1-dlog1}, tel que pour tout entier $n\geq 1$, 
le morphisme canonique $\cH^{(r')}_n\rightarrow \cH^{(r)}_n$ soit annulé par $p^\alpha$. 
\item[{\rm (iii)}] Il existe un entier $\beta\geq 0$, dépendant de $r$, $r'$ et $d$, 
mais pas des données \eqref{higgs1-dlog1}, tel que pour tous entiers $n,q\geq 1$, 
le morphisme canonique
\begin{equation}\label{higgs1-cohbis145b}
\rH^q(\Delta,\cS^{(r')}/p^n\cS^{(r')})\rightarrow \rH^q(\Delta,\cS^{(r)}/p^n\cS^{(r)})
\end{equation}
soit annulé par $p^\beta$. 
\end{itemize}
\end{prop}

(i) Cela résulte de \ref{higgs1-cohbis144}(i) et de 
la suite exacte longue de cohomologie associée à la suite exacte courte de $\hoR$-représentations de $\Delta$ 
\begin{equation}\label{higgs1-cohbis145c}
0\longrightarrow \hcS^{(r)}\stackrel{\cdot p^n}{\longrightarrow} \hcS^{(r)}\longrightarrow
\hcS^{(r)}/p^n\hcS^{(r)}\longrightarrow 0.
\end{equation}
On en déduit aussi un morphisme $\hRun$-linéaire presque-injectif
\begin{equation}\label{higgs1-cohbis145d}
\cH^{(r)}_n\rightarrow \rH^1_{\cont}(\Delta,\hcS^{(r)}).
\end{equation}

(ii) Cela résulte de \eqref{higgs1-cohbis145d} et \ref{higgs1-cohbis144}(ii).

(iii) Pour tous entiers $n,q\geq 1$, la suite exacte longue de cohomologie déduite de \eqref{higgs1-cohbis145c} 
fournit une suite exacte de $\hRun$-modules
\begin{equation}
0\rightarrow \rH^q_{\cont}(\Delta,\hcS^{(r)})/p^n\rH^q_{\cont}(\Delta,\hcS^{(r)})
\rightarrow \rH^q(\Delta,\cS^{(r)}/p^n\cS^{(r)}) \rightarrow T_n^{(r),q}\rightarrow 0,
\end{equation}
où $T_n^{(r),q}$ est un sous-module de $p^n$-torsion de $\rH^{q+1}_{\cont}(\Delta,\hcS^{(r)})$.
Posons $r''=(r+r')/2$. D'après \ref{higgs1-cohbis144}(ii), il existe un entier $\beta'>0$, dépendant 
seulement de $r$, $r'$ et $d$, tel que
pour tout entier $q\geq 1$, les morphismes canoniques 
\[
\rH^q_{\cont}(\Delta,\hcS^{(r')})\rightarrow \rH^q_{\cont}(\Delta,\hcS^{(r'')}) \ \ \ {\rm et} \ \ \ 
\rH^q_{\cont}(\Delta,\hcS^{(r'')})\rightarrow \rH^q_{\cont}(\Delta,\hcS^{(r)})
\] 
soient annulés par $p^{\beta'}$. La proposition s'ensuit en prenant $\beta=2\beta'$.

\section{Représentations de Dolbeault}\label{higgs1-dolbeault}

\subsection{}\label{higgs1-dolbeault1}\index{101002@$(\tX,\cM_\tX)$}\index{101202@$\cF^{(r)}$, $\cC^{(r)}$, $\cC^\dagger$}
\index{101203@$\alpha^{r,r'}$, $\halpha^{r,r'}$}
\index{101204@$d_{\cC^{(r)}}$, $d_{\hcC^{(r)}}$, $d_{\cC^\dagger}$}\index{101205@$\tmK^\bullet(\hcC^{(r)},p^rd_{\hcC^{(r)}})$}
\index{Higgs-Tate!4@Algèbre de ---  d'épaisseur $r$ ($\cC^{(r)}$)}
On fixe dans la suite de cet article une $(\cA_2(\oS),\cM_{\cA_2(\oS)})$-déformation lisse 
$(\tX,\cM_\tX)$ de $(\coX,\cM_{\coX})$ \eqref{higgs1-deflog1}
et on note $\cF$ la $\hoR$-extension de Higgs-Tate et
$\cC$ la $\hoR$-algèbre de Higgs-Tate associés \eqref{higgs1-tor201}.
On désigne par $\hcC$ le séparé complété $p$-adique de $\cC$. 
Pour tout nombre rationnel $r\geq 0$, on note $\cF^{(r)}$ la $\hoR$-représentation de $\Delta$ 
déduite de $\cF$ par image inverse par la multiplication par $p^r$ sur 
$\xi^{-1}\tOmega^1_{R/\co_K}\otimes_R\hoR$, de sorte qu'on a une suite exacte localement scindée de $\hoR$-modules
\begin{equation}\label{higgs1-dolbeault1a}
0\longrightarrow \hoR\longrightarrow \cF^{(r)}\stackrel{u^{(r)}}{\longrightarrow} 
\xi^{-1}\tOmega^1_{R/\co_K}\otimes_R\hoR\longrightarrow 0.
\end{equation}
D'après (\cite{illusie1} I 4.3.1.7), cette suite induit pour tout entier $n\geq 1$, une suite exacte  \eqref{higgs1-not621}
\begin{equation}\label{higgs1-dolbeault1b}
0\rightarrow \rS^{n-1}_{\hoR}(\cF^{(r)})\rightarrow \rS^{n}_{\hoR}(\cF^{(r)})\rightarrow 
\rS^n_{\hoR}(\xi^{-1}\tOmega^1_{R/\co_K}
\otimes_R\hoR)\rightarrow 0.
\end{equation}
Les $\hoR$-modules $(\rS^{n}_{\hoR}(\cF^{(r)}))_{n\in \mN}$ forment donc un système inductif filtrant, 
dont la limite inductive 
\begin{equation}\label{higgs1-dolbeault1c}
\cC^{(r)}=\underset{\underset{n\geq 0}{\longrightarrow}}\lim\ \rS^n_{\hoR}(\cF^{(r)})
\end{equation}
est naturellement munie d'une structure de $\hoR$-algèbre. 
Il existe un et un unique homomorphisme de $\hoR$-algèbres 
\begin{equation}\label{higgs1-dolbeault1d}
\mu^{(r)} \colon \cC^{(r)}\rightarrow \cC^{(r)}\otimes_\hoR\cS^{(r)},
\end{equation}
où $\cS^{(r)}$ est la $\hoR$-algèbre définie dans  \eqref{higgs1-chb1a}, 
tel que pour tout $x\in \cF^{(r)}$, on ait 
\begin{equation}\label{higgs1-dolbeault1e}
\mu^{(r)}(x)=x\otimes 1+1\otimes (p^r \cdot u^{(r)}(x)).
\end{equation} 
Celui-ci fait de $\Spec(\cC^{(r)})$ un $\Spec(\cS^{(r)})$-fibré principal homogène sur $\hY$ 
(cf. la preuve de \ref{higgs1-eph12}). 

L'action de $\Delta$ sur $\cF^{(r)}$ induit une action sur $\cC^{(r)}$ 
par des automorphismes d'anneaux, compatible avec son action sur $\hoR$, que l'on appelle {\em action canonique}.
La $\hoR$-algèbre $\cC^{(r)}$ munie de cette action
est appelée {\em l'algèbre de Higgs-Tate} d'épaisseur $r$ associée à $(\tX,\cM_\tX)$. 
On note $\hcC^{(r)}$ le séparé complété $p$-adique de $\cC^{(r)}$ que l'on suppose toujours muni de 
la topologie $p$-adique. On munit $\hcC^{(r)}\otimes_{\mZ_p}\mQ_p$ de la topologie $p$-adique \eqref{higgs1-not54}.
Compte tenu de \ref{higgs1-pur8} et de sa preuve, $\cC^{(r)}$ et $\hcC^{(r)}$ sont $\co_C$-plats.
Pour tous nombres rationnels $r'\geq r\geq 0$, on a un $\hoR$-homomorphisme canonique injectif et $\Delta$-équivariant
$\alpha^{r,r'}\colon \cC^{(r')}\rightarrow \cC^{(r)}$. On vérifie aussitôt que l'homomorphisme induit
$\halpha^{r,r'}\colon\hcC^{(r')}\rightarrow \hcC^{(r)}$ est injectif. On pose 
\begin{equation}\label{higgs1-dolbeault1f}
\cC^\dagger=\underset{\underset{r\in \mQ_{>0}}{\longrightarrow}}{\lim} \hcC^{(r)},
\end{equation}
que l'on identifie à une sous-$\hoR$-algèbre de $\hcC=\hcC^{(0)}$ par la limite inductive des homomorphismes 
$(\halpha^{0,r})_{r\in \mQ_{>0}}$. Les actions de $\Delta$ sur les anneaux $(\hcC^{(r)})_{r\in \mQ_{>0}}$ 
induisent une action sur $\cC^\dagger$ par des automorphismes d'anneaux, 
compatible avec ses actions sur $\hoR$ et sur $\hcC$.

On désigne par  
\begin{equation}\label{higgs1-dolbeault1g}
d_{\cC^{(r)}}\colon \cC^{(r)}\rightarrow \xi^{-1}\tOmega^1_{R/\co_K}\otimes_R\cC^{(r)}
\end{equation}
la $\hoR$-dérivation universelle de $\cC^{(r)}$ et par 
\begin{equation}\label{higgs1-dolbeault1h}
d_{\hcC^{(r)}}\colon \hcC^{(r)}\rightarrow \xi^{-1}\tOmega^1_{R/\co_K}\otimes_R\hcC^{(r)}
\end{equation}
son prolongement aux complétés (on notera que le $R$-module $\tOmega^1_{R/\co_K}$ est libre de type fini). 
On voit aussitôt, comme pour $d_{\cC}$ \eqref{higgs1-tor23d}, que les dérivations 
$d_{\cC^{(r)}}$ et $d_{\hcC^{(r)}}$ sont $\Delta$-équivariantes. 
Par ailleurs, $d_{\cC^{(r)}}$ et $d_{\hcC^{(r)}}$ sont également des $\hoR$-champs 
de Higgs à coefficients dans $\xi^{-1}\tOmega^1_{R/\co_K}$ 
puisque $\xi^{-1}\tOmega^1_{R/\co_K}\otimes_R \hoR=d_{\cC^{(r)}}(\cF^{(r)})\subset d_{\cC^{(r)}}(\cC^{(r)})$ 
(cf. \ref{higgs1-not16} et \ref{higgs1-not22}). On désigne par $\mK^\bullet(\hcC^{(r)},p^rd_{\hcC^{(r)}})$
le complexe de Dolbeault de $(\hcC^{(r)},p^rd_{\hcC^{(r)}})$  \eqref{higgs1-not61}
et par $\tmK^\bullet(\hcC^{(r)},p^rd_{\hcC^{(r)}})$ 
le complexe de Dolbeault augmenté 
\begin{equation}\label{higgs1-dolbeault1i}
\hoR\rightarrow \mK^0(\hcC^{(r)},p^rd_{\hcC^{(r)}})\rightarrow \mK^1(\hcC^{(r)},p^rd_{\hcC^{(r)}})\rightarrow \dots
\rightarrow \mK^n(\hcC^{(r)},p^rd_{\hcC^{(r)}})\rightarrow \dots,
\end{equation}
où $\hoR$ est placé en degré $-1$ et la différentielle $\hoR\rightarrow\hcC^{(r)}$ est l'homomorphisme canonique. 

Pour tous nombres rationnels $r'\geq r\geq 0$, on a 
\begin{equation}\label{higgs1-dolbeault1j}
p^{r'}(\id \times \alpha^{r,r'}) \circ d_{\cC^{(r')}}=p^rd_{\cC^{(r)}}\circ \alpha^{r,r'}.
\end{equation}
Par suite, $\halpha^{r,r'}$ induit un morphisme de complexes 
\begin{equation}\label{higgs1-dolbeault1k}
\nu^{r,r'}\colon \tmK^\bullet(\hcC^{(r')},p^{r'}d_{\hcC^{(r')}})\rightarrow \tmK^\bullet(\hcC^{(r)},p^rd_{\hcC^{(r)}}).
\end{equation}

D'après \eqref{higgs1-dolbeault1j}, les dérivations $p^rd_{\hcC^{(r)}}$ induisent une $\hoR$-dérivation 
\begin{equation}\label{higgs1-dolbeault1l}
d_{\cC^\dagger}\colon \cC^\dagger\rightarrow \xi^{-1}\tOmega^1_{R/\co_K} \otimes_R\cC^\dagger,
\end{equation}
qui n'est autre que la restriction de $d_{\hcC}$ à $\cC^\dagger$. 
C'est également un $\hoR$-champ de Higgs à coefficients dans $\xi^{-1}\tOmega^1_{R/\co_K}$. 
On note $\mK^\bullet(\cC^\dagger,d_{\cC^\dagger})$ 
le complexe de Dolbeault de $(\cC^\dagger,d_{\cC^\dagger})$. 
Comme $\hoR$ est $\co_C$-plat \eqref{higgs1-pur8}, pour tout nombre rationnel $r\geq 0$, on a 
\begin{equation}\label{higgs1-dolbeault1m}
\ker(d_{\cC^\dagger})=\ker(d_{\hcC^{(r)}})=\hoR.
\end{equation}

\subsection{}\label{higgs1-dolb43}\index{101210@$\cC_0$, $\cC_0^{(r)}$}
Lorsque $(\tX,\cM_\tX)$ est la déformation $(\tX_0,\cM_{\tX_0})$
de $(\coX,\cM_\coX)$ définie par la carte $(P,\gamma)$ \eqref{higgs1-ext24ac}, 
on équipe d'un indice $_0$ les objets définis dans \eqref{higgs1-dolbeault1}~: $\cC_0$, $\cC_0^{(r)}$... 
La section $\psi_0\in \cL_0(\hY)$ définie par la carte $(P,\gamma)$ \eqref{higgs1-ext24},
induit alors un isomorphisme de $\hoR$-algèbres 
\begin{equation}\label{higgs1-dolb43a}
\cS\stackrel{\sim}{\rightarrow}\cC_0,
\end{equation}
où $\cS$ est la $\hoR$-algèbre définie dans \eqref{higgs1-tor2d}.
Celui-ci est $\Delta$-équivariant lorsque l'on munit $\cS$ de l'action $\varphi_0=\varphi_{\psi_0}$
de $\Delta$ induite par $\psi_0$ \eqref{higgs1-tor202g}. Il est par ailleurs 
compatible aux dérivations universelles d'après \eqref{higgs1-tor23c}. 
Pour tout nombre rationnel $r\geq 0$, $\psi_0$ induit un $\hoR$-homomorphisme  
$\cC_0^{(r)}\rightarrow \hoR$ et par suite un isomorphisme de $\hoR$-algèbres 
\begin{equation}\label{higgs1-dolb43b}
\cS^{(r)}\stackrel{\sim}{\rightarrow}\cC_0^{(r)},
\end{equation}
où $\cS^{(r)}$ est la $\hoR$-algèbre définie dans  \eqref{higgs1-chb1a}.
Pour tous nombres rationnels $r'\geq r\geq 0$, le diagramme 
\begin{equation}\label{higgs1-dolb43c}
\xymatrix{
{\cS^{(r')}}\ar[d]\ar[r]^\sim&{\cC_0^{(r')}}\ar[d]\\
{\cS^{(r)}}\ar[r]^\sim&{\cC_0^{(r)}}}
\end{equation}
où les flèches verticales sont les homomorphismes canoniques,
est commutatif. On en déduit que l'isomorphisme \eqref{higgs1-dolb43b} est $\Delta$-équivariant  
lorsque l'on munit $\cS^{(r)}$ de l'action de $\Delta$ induite par $\varphi_0$ \eqref{higgs1-cohbis3}.
Par ailleurs, on vérifie aussitôt que le diagramme \eqref{higgs1-chb1d}
\begin{equation}\label{higgs1-dolb43d}
\xymatrix{
{\cS^{(r)}}\ar[d]_{d_{\cS^{(r)}}}\ar[r]^-(0.5)\sim&{\cC_0^{(r)}}\ar[d]^{d_{\cC_0^{(r)}}}\\
{\xi^{-1}\tOmega^1_{R/\co_K}\otimes_R\cS^{(r)}}
&{\xi^{-1}\tOmega^1_{R/\co_K}\otimes_R\cC_0^{(r)}}\ar@{=}[l]}
\end{equation}
est commutatif. 

\vspace{2mm}

Pour des références ultérieures, nous réécrivons dans les cinq propositions qui suivent les principaux résultats du § \ref{higgs1-CGII} 
relatifs à l'algèbre $\cS$ en termes de l'algèbre $\cC$, pour une déformation générale 
$(\tX,\cM_{\tX})$ de $(\coX,\cM_\coX)$, en tenant compte de \ref{higgs1-dolb43} et \ref{higgs1-tor25}.  

\begin{prop}\label{higgs1-dolb44}
\begin{itemize}
\item[{\rm (i)}] Pour tous nombres rationnels $r'>r>0$, il existe un nombre rationnel $\alpha\geq 0$ dépendant de $r$ et $r'$, 
mais pas des données \eqref{higgs1-dlog1}, tel que
\begin{equation}\label{higgs1-dolb44cc}
p^\alpha\nu^{r,r'}\colon \tmK^\bullet(\hcC^{(r')},p^{r'}d_{\hcC^{(r')}})\rightarrow 
\tmK^\bullet(\hcC^{(r)},p^rd_{\hcC^{(r)}}),
\end{equation}
où $\nu^{r,r'}$ est le morphisme  \eqref{higgs1-dolbeault1k}, 
soit homotope à $0$ par une homotopie $\hoR$-linéaire. 

\item[{\rm (ii)}] Pour tous nombres rationnels $r'>r>0$, le morphisme canonique
\begin{equation}\label{higgs1-dolb44c}
\nu^{r,r'}\otimes_{\mZ_p}\mQ_p\colon \tmK^\bullet(\hcC^{(r')},p^{r'}d_{\hcC^{(r')}})\otimes_{\mZ_p}\mQ_p\rightarrow 
\tmK^\bullet(\hcC^{(r)},p^rd_{\hcC^{(r)}})\otimes_{\mZ_p}\mQ_p
\end{equation}
est homotope à $0$ par une homotopie continue. 

\item[{\rm (iii)}] Le complexe $\mK^\bullet(\cC^\dagger,d_{\cC^\dagger})\otimes_{\mZ_p}\mQ_p$ est une résolution 
de $\hoR[\frac 1 p]$.
\end{itemize}
\end{prop}
Cela résulte de \ref{higgs1-cohbis70}, \ref{higgs1-cohbis7} et \ref{higgs1-cohbis8}.

\begin{prop}\label{higgs1-dolb45}
Pour tout nombre rationnel $r\geq 0$, les actions de $\Delta$ sur $\cC^{(r)}$ et $\hcC^{(r)}$ sont 
continues pour les topologies $p$-adiques. 
\end{prop}

Cela résulte de \ref{higgs1-cohbis3}. 

\begin{prop}\label{higgs1-dolb46}
Soit $r$ un nombre rationnel $>0$. Alors~:
\begin{itemize}
\item[{\rm (i)}] Le morphisme canonique 
\begin{equation}\label{higgs1-dolb46a}
\hRun\otimes_{\mZ_p}\mQ_p\rightarrow (\hcC^{(r)}\otimes_{\mZ_p}\mQ_p)^\Delta
\end{equation}
est un isomorphisme. 
\item[{\rm (ii)}] Pour tout entier $i\geq 1$, on a 
\begin{equation}\label{higgs1-dolb46b}
\underset{\underset{r\in \mQ_{>0}}{\longrightarrow}}{\lim}\ 
\rH^i_\cont(\Delta,\hcC^{(r)}\otimes_{\mZ_p}\mQ_p)=0.
\end{equation}
\end{itemize}
\end{prop}

Cela résulte de \ref{higgs1-cohbis9}.

\begin{cor}\label{higgs1-dolb47}
Pour tout nombre rationnel $r> 0$, on a $(\cC^\dagger)^\Delta=(\hcC^{(r)})^\Delta=\hRun$.
\end{cor}

Cela résulte de \ref{higgs1-cohbis99} (ou de \ref{higgs1-dolb46}).

\begin{prop}\label{higgs1-dolb48}
Soient $r,r'$ deux nombres rationnels tels que $r'>r>0$. Alors~:
\begin{itemize}
\item[{\rm (i)}] Pour tout entier $n\geq 1$, l'homomorphisme canonique 
\begin{equation}\label{higgs1-dolb48a}
R_1/p^nR_1\rightarrow (\cC^{(r)}/p^n\cC^{(r)})^{\Delta}
\end{equation}
est presque-injectif. Notons $\cH^{(r)}_n$ son conoyau. 
\item[{\rm (ii)}] Il existe un entier $\alpha\geq 0$, dépendant de $r$, $r'$ et $d$, 
mais pas des données \eqref{higgs1-dlog1}, tel que pour tout entier $n\geq 1$, 
le morphisme canonique $\cH^{(r')}_n\rightarrow \cH^{(r)}_n$ soit annulé par $p^\alpha$. 
\item[{\rm (iii)}] Il existe un entier $\beta\geq 0$, dépendant de $r$, $r'$ et $d$, 
mais pas des données \eqref{higgs1-dlog1}, tel que pour tous entiers $n,q\geq 1$, 
le morphisme canonique
\begin{equation}\label{higgs1-dolb48b}
\rH^q(\Delta,\cC^{(r')}/p^n\cC^{(r')})\rightarrow \rH^q(\Delta,\cC^{(r)}/p^n\cC^{(r)})
\end{equation}
soit annulé par $p^\beta$. 
\end{itemize}
\end{prop}

Cela résulte de \ref{higgs1-cohbis145}.

\subsection{}\label{higgs1-dolb2}\index{101220@$\mH(M)$}
Pour toute $\hoR$-représentation $M$ de $\Delta$, on note $\mH(M)$ le $\hRun$-module défini par 
\begin{equation}\label{higgs1-dolb2a}
\mH(M)=(M\otimes_{\hoR}\cC^\dagger)^\Delta.
\end{equation}
On le munit du $\hRun$-champ de Higgs 
à coefficients dans $\xi^{-1}\tOmega^1_{R/\co_K}$ induit par $d_{\cC^\dagger}$ \eqref{higgs1-dolbeault1l} (cf. \ref{higgs1-eip75}).   
On définit ainsi un foncteur 
\begin{equation}\label{higgs1-dolb2c}
\mH\colon \bRep_{\hoR}(\Delta) \rightarrow \bMH(\hRun,\xi^{-1}\tOmega^1_{R/\co_K}).
\end{equation}

\subsection{}\label{higgs1-dolb3}\index{101222@$\mV(N)$}
Pour tout $\hRun$-module de Higgs $(N,\theta)$ à coefficients dans $\xi^{-1}\tOmega^1_{R/\co_K}$ \eqref{higgs1-eip75}, 
on note $\mV(N)$ le $\hoR$-module défini par 
\begin{equation}\label{higgs1-dolb3a}
\mV(N)=(N\otimes_{\hRun}\cC^\dagger)^{\theta_\tot=0},
\end{equation}
où $\theta_\tot=\theta\otimes \id+\id\otimes d_{\cC^\dagger}$ est le $\hRun$-champ de Higgs total sur 
$N\otimes_{\hRun}\cC^\dagger$ \eqref{higgs1-not6ab}.
On le munit de l'action $\hoR$-semi-linéaire de $\Delta$ induite
par son action naturelle sur $\cC^\dagger$. 
On définit ainsi un foncteur 
\begin{equation}\label{higgs1-dolb3c}
\mV\colon \bMH(\hRun,\xi^{-1}\tOmega^1_{R/\co_K})\rightarrow \bRep_{\hoR}(\Delta).
\end{equation}

\begin{remas}\label{higgs1-dolb32}
{\rm (i)}\ Il résulte de \ref{higgs1-tor25} que les foncteurs $\mH$ et $\mV$ ne dépendent pas du choix 
de la déformation $(\tX,\cM_\tX)$, à isomorphisme non-canonique près.

(ii)\ Pour toute $\hoR$-représentation $M$ de $\Delta$, le morphisme canonique 
\begin{equation}\label{higgs1-dolb23a}
\mH(M)\otimes_{\hRun}\hRun[\frac 1 p]\rightarrow \mH(M\otimes_{\hoR}\hoR[\frac 1 p])
\end{equation}
est un isomorphisme.   

{\rm (iii)}\ Pour tout $\hRun$-module de Higgs $(N,\theta)$ à coefficients dans $\xi^{-1}\tOmega^1_{R/\co_K}$, 
le morphisme canonique 
\begin{equation}\label{higgs1-dolb32a}
\mV(N)\otimes_{\hoR}\hoR[\frac 1 p]\rightarrow \mV(N\otimes_{\hRun}\hRun[\frac 1 p])
\end{equation}
est un isomorphisme.   
\end{remas}

\begin{defi}\label{higgs1-dolb21}\index{Representation de Dolbeault a@$\hoR$-représentation de Dolbeault}
On dit qu'une $\hoR$-représentation continue $M$ de $\Delta$ 
est {\em de Dolbeault} si les conditions suivantes sont remplies~:
\begin{itemize}
\item[{\rm (i)}] $M$ est un $\hoR$-module projectif de type fini, muni de la topologie $p$-adique~;
\item[{\rm (ii)}] $\mH(M)$ est un $\hRun$-module projectif de type fini~; 
\item[{\rm (iii)}] le morphisme $\cC^\dagger$-linéaire canonique 
\begin{equation}\label{higgs1-dolb21b}
\mH(M) \otimes_{\hRun}\cC^\dagger\rightarrow  M\otimes_{\hoR}\cC^\dagger
\end{equation}
est un isomorphisme.
\end{itemize}
\end{defi}

Cette notion ne dépend pas du choix de $(\tX,\cM_\tX)$ (\ref{higgs1-tor25} et \ref{higgs1-dolb32}(i)).

\begin{defi}\label{higgs1-dolb31}\index{Module de Higgs soluble a@$\hRun$-module de Higgs soluble} 
On dit qu'un $\hRun$-module de Higgs $(N,\theta)$ à coefficients dans $\xi^{-1}\tOmega^1_{R/\co_K}$
est {\em soluble} si les conditions suivantes sont remplies~:
\begin{itemize}
\item[{\rm (i)}] $N$ est un $\hRun$-module projectif de type fini~;
\item[{\rm (ii)}] $\mV(N)$ est un $\hoR$-module projectif de type fini~;
\item[{\rm (iii)}] le morphisme $\cC^\dagger$-linéaire canonique 
\begin{equation}\label{higgs1-dolb31c}
\mV(N) \otimes_{\hoR}\cC^\dagger\rightarrow  N\otimes_{\hRun}\cC^\dagger
\end{equation}
est un isomorphisme.
\end{itemize}
\end{defi}

Cette notion ne dépend pas du choix de $(\tX,\cM_\tX)$ (\ref{higgs1-tor25} et \ref{higgs1-dolb32}(i)).

\begin{lem}\label{higgs1-dolb6}
Soit $M$ une $\hoR$-représentation de Dolbeault de $\Delta$. Alors le $\hRun$-module de Higgs 
$\mH(M)$ est soluble, et on a un $\hoR$-isomorphisme canonique fonctoriel et $\Delta$-équivariant
\begin{equation}\label{higgs1-dolb6a}
\mV(\mH(M))\stackrel{\sim}{\rightarrow} M.
\end{equation}
\end{lem}

En effet, le morphisme $\cC^\dagger$-linéaire canonique 
\begin{equation}\label{higgs1-dolb6b}
\mH(M) \otimes_{\hRun}\cC^\dagger\rightarrow  M\otimes_{\hoR}\cC^\dagger
\end{equation}
est un isomorphisme $\Delta$-équivariant 
de $\hoR$-modules de Higgs à coefficients dans $\xi^{-1}\tOmega^1_{R/\co_K}$,
où $\cC^\dagger$ est muni du champ de Higgs $d_{\cC^\dagger}$ \eqref{higgs1-dolbeault1l}, 
$\mH(M)$ est muni de l'action triviale de $\Delta$ et $M$ est muni du champ de Higgs nul (cf. \ref{higgs1-not17}). 
Comme $M$ est $\hoR$-plat et que $\ker(d_{\cC^\dagger})=\hoR$, on en déduit 
un $\hoR$-isomorphisme $\Delta$-équivariant $\mV(\mH(M))\stackrel{\sim}{\rightarrow} M$.
Le morphisme $\cC^\dagger$-linéaire canonique 
\begin{equation}
\mV(\mH(M))\otimes_{\hoR}\cC^\dagger\rightarrow \mH(M)\otimes_{\hoR}\cC^\dagger
\end{equation}
s'identifie alors à l'inverse de l'isomorphisme \eqref{higgs1-dolb6b}, ce qui montre que $\mH(M)$ est soluble.

\begin{lem}\label{higgs1-dolb5}
Soit $(N,\theta)$ un $\hRun$-module de Higgs soluble à coefficients dans $\xi^{-1}\tOmega^1_{R/\co_K}$.
Alors la $\hoR$-représentation $\mV(N)$ de $\Delta$ est de Dolbeault,  
et on a un isomorphisme canonique fonctoriel de $\hRun$-modules de Higgs
\begin{equation}\label{higgs1-dolb5a}
\mH(\mV(N))\stackrel{\sim}{\rightarrow} N.
\end{equation}
\end{lem}

Pour tous nombres rationnels $r'\geq r \geq 0$, le morphisme canonique 
$N\otimes_{\hRun}\hcC^{(r')}\rightarrow N\otimes_{\hRun}\hcC^{(r)}$ est injectif car $N$ est $\hRun$-plat. 
Comme $\mV(N)$ est un $\hoR$-module de type fini, il existe un nombre rationnel $r>0$ tel que l'on ait
\begin{equation}
\mV(N)=(N\otimes_{\hRun}\hcC^{(r)})^{\theta^{(r)}_\tot=0},
\end{equation}
où $\theta^{(r)}_\tot=\theta\otimes \id+p^r\id\otimes d_{\hcC^{(r)}}$ est le $\hRun$-champ de Higgs total sur 
$N\otimes_{\hRun}\hcC^{(r)}$. D'autre part, comme $\hcC^{(r)}$ est $\co_C$-plat et que $N$ est $\hRun$-plat,
$\xi^{-1}\tOmega^1_{R/\co_K}\otimes_RN\otimes_{\hRun}\hcC^{(r)}$ est $\co_C$-plat. Par suite, pour tout $n\geq 0$
et tout $x\in N\otimes_{\hRun}\hcC^{(r)}$ tels que $\theta^{(r)}_\tot(p^nx)=0$, on a $\theta^{(r)}_\tot(x)=0$.  
On en déduit que  
\begin{equation}
(p^nN\otimes_{\hRun}\hcC^{(r)})\cap \mV(N)=p^n\mV(N).
\end{equation}
Donc la topologie $p$-adique sur $\mV(N)$ est induite par celle de $N\otimes_{\hRun}\hcC^{(r)}$. 
Il résulte alors de \ref{higgs1-dolb45} que l'action de $\Delta$ sur $\mV(N)$ est continue pour la topologie $p$-adique. 

Le morphisme $\cC^\dagger$-linéaire canonique 
\begin{equation}\label{higgs1-dolb5b}
\mV(N) \otimes_{\hoR}\cC^\dagger\rightarrow  N\otimes_{\hRun}\cC^\dagger
\end{equation}
est un isomorphisme $\Delta$-équivariant 
de $\hoR$-modules de Higgs à coefficients dans $\xi^{-1}\tOmega^1_{R/\co_K}$, 
où $\cC^\dagger$ est muni du champ de Higgs $d_{\cC^\dagger}$,
$\mV(N)$ est muni du champ de Higgs nul et $N$ est muni de l'action triviale de $\Delta$ (cf. \ref{higgs1-not17}).
Comme $N$ est un facteur direct d'un $\hRun$-module libre de type fini, on a 
$(N\otimes_{\hRun}\cC^\dagger)^\Delta=N$ \eqref{higgs1-dolb47}. 
On en déduit un isomorphisme
de $\hRun$-modules de Higgs $\mH(\mV(N))\stackrel{\sim}{\rightarrow} N$.
Le morphisme $\cC^\dagger$-linéaire canonique 
\begin{equation}
\mH(\mV(M))\otimes_{\hRun}\cC^\dagger\rightarrow \mV(M)\otimes_{\hRun}\cC^\dagger
\end{equation}
s'identifie alors à l'inverse de l'isomorphisme \eqref{higgs1-dolb5b}, ce qui montre que $\mV(N)$ est de Dolbeault.

\begin{prop}\label{higgs1-dolb7}
Les foncteurs $\mV$ et $\mH$ induisent des équivalences de catégories quasi-inverses l'une de l'autre,
entre la catégorie des $\hoR$-représentations de Dolbeault de $\Delta$ et celle des $\hRun$-modules de Higgs
solubles à coefficients dans $\xi^{-1}\tOmega^1_{R/\co_K}$.
\end{prop}
Cela résulte de \ref{higgs1-dolb6} et \ref{higgs1-dolb5}.

\begin{defi}\label{higgs1-drt2}\index{Representation de Dolbeault b@$\hoR[\frac 1 p]$-représentation de Dolbeault}
On dit qu'une $\hoR[\frac 1 p]$-représentation continue $M$ de $\Delta$ est {\em de Dolbeault} 
si les conditions suivantes sont remplies~:
\begin{itemize}
\item[(i)] $M$ est un $\hoR[\frac 1 p]$-module projectif de type fini, muni de la topologie $p$-adique \eqref{higgs1-not54};
\item[(ii)] $\mH(M)$ est un $\hRun[\frac 1 p]$-module projectif de type fini \eqref{higgs1-dolb2a};
\item[(iii)] le morphisme canonique 
\begin{equation}\label{higgs1-drt2a}
\mH(M) \otimes_{\hRun}\cC^\dagger\rightarrow  M\otimes_{\hoR}\cC^\dagger
\end{equation}
est un isomorphisme.
\end{itemize}
\end{defi}

Cette notion ne dépend pas du choix de $(\tX,\cM_{\tX})$ (\ref{higgs1-tor25} et \ref{higgs1-dolb32}(i)).

\begin{rema}\label{higgs1-drt200}
Pour qu'une $\hoR[\frac 1 p]$-représentation continue $M$ de $\Delta$ 
soit de Dolbeault, il faut et il suffit qu'elle vérifie les conditions (i) et (ii) de \ref{higgs1-drt2} ainsi que la condition
suivante~:
\begin{itemize}
\item[{\rm (iii')}] il existe un nombre rationnel $r>0$ tel que $\mH(M)$ soit contenu dans 
$M\otimes_{\hoR}\hcC^{(r)}$ \eqref{higgs1-dolbeault1c}, et que le morphisme canonique 
\begin{equation}\label{higgs1-drt200a}
\mH(M) \otimes_{\hRun}\hcC^{(r)}\rightarrow  M\otimes_{\hoR}\hcC^{(r)}
\end{equation}
soit un isomorphisme.
\end{itemize}
En effet, pour tout nombre rationnel $r> 0$, le morphisme canonique 
$M\otimes_\hoR\hcC^{(r)}\rightarrow M\otimes_\hoR\cC^\dagger$ est injectif car $M$ est $\hoR$-plat.
La condition (iii') implique clairement la condition \ref{higgs1-drt2}(iii).
Inversement, supposons la condition \ref{higgs1-drt2}(iii) remplie.  
Comme $\mH(M)$ est un $\hRun[\frac 1 p]$-module de type fini,
il existe un nombre rationnel $r>0$ tel que $\mH(M)$ soit contenu dans 
$M\otimes_{\hoR}\hcC^{(r)}$. Comme $M$ est de type fini sur $\hoR[\frac 1 p]$, 
quitte à diminuer $r$, on peut supposer le morphisme \eqref{higgs1-drt200a} surjectif. 
Par ailleurs, $\mH(M)$ étant $\hRun$-plat, pour tout nombre rationnel $r>0$, 
le morphisme canonique 
\begin{equation}
\mH(M)\otimes_\hRun\hcC^{(r)}\rightarrow \mH(M)\otimes_\hRun\cC^\dagger
\end{equation} 
est injectif. Le morphisme \eqref{higgs1-drt200a} est donc injectif en vertu de \ref{higgs1-drt2}(iii). 
\end{rema}

\begin{defi}\label{higgs1-drt3}\index{Module de Higgs soluble b@$\hRun[\frac 1 p]$-module de Higgs soluble} 
On dit qu'un $\hRun[\frac 1 p]$-module de Higgs $(N,\theta)$ à coefficients dans $\xi^{-1}\tOmega^1_{R/\co_K}$
est {\em soluble} si les conditions suivantes sont remplies~:
\begin{itemize}
\item[(i)] $N$ est un $\hRun[\frac 1 p]$-module projectif de type fini~;
\item[(ii)] $\mV(N)$ est un $\hoR[\frac 1 p]$-module projectif de type fini~;
\item[(iii)] le morphisme canonique 
\begin{equation}\label{higgs1-drt3a}
\mV(N) \otimes_{\hoR}\cC^\dagger\rightarrow  N\otimes_{\hRun}\cC^\dagger
\end{equation}
est un isomorphisme.
\end{itemize}
\end{defi}

Cette notion ne dépend pas du choix de $(\tX,\cM_\tX)$ (\ref{higgs1-tor25} et \ref{higgs1-dolb32}(i)).

\begin{rema}\label{higgs1-drt300}
Pour qu'un $\hRun[\frac 1 p]$-module de Higgs $(N,\theta)$ à coefficients dans $\xi^{-1}\tOmega^1_{R/\co_K}$
soit soluble, il faut et il suffit qu'il vérifie les conditions (i) et (ii) de \ref{higgs1-drt3} ainsi que la condition
suivante~:
\begin{itemize}
\item[{\rm (iii')}] il existe un nombre rationnel $r>0$ tel que $\mV(N)$ soit contenu dans 
$N\otimes_{\hRun}\hcC^{(r)}$, et que le morphisme canonique 
\begin{equation}\label{higgs1-drt300a}
\mV(N) \otimes_{\hoR}\hcC^{(r)}\rightarrow  N\otimes_{\hRun}\hcC^{(r)}
\end{equation}
soit un isomorphisme.
\end{itemize}
La preuve, similaire à celle de \ref{higgs1-drt200}, est laissée au lecteur. 
\end{rema}

\begin{lem}\label{higgs1-drt4}
Soit $M$ une $\hoR[\frac 1 p]$-représentation de Dolbeault de $\Delta$. 
Alors le $\hRun[\frac 1p]$-module de Higgs $\mH(M)$ est soluble et 
on a un $\hoR[\frac 1 p]$-isomorphisme canonique fonctoriel et $\Delta$-équivariant
\begin{equation}\label{higgs1-drt4a}
\mV(\mH(M))\stackrel{\sim}{\rightarrow} M. 
\end{equation}
\end{lem}

En effet, le morphisme $\cC^\dagger$-linéaire canonique 
\begin{equation}\label{higgs1-drt4b}
\mH(M) \otimes_{\hRun}\cC^\dagger\rightarrow  M\otimes_{\hoR}\cC^\dagger
\end{equation}
est un $\cC^\dagger$-isomorphisme $\Delta$-équivariant 
de $\hoR$-modules de Higgs à coefficients dans $\xi^{-1}\tOmega^1_{R/\co_K}$,
où $\cC^\dagger$ est muni du champ de Higgs $d_{\cC^\dagger}$ \eqref{higgs1-dolbeault1l}, 
$\mH(M)$ est muni de l'action triviale de $\Delta$ et $M$ est muni du champ de Higgs nul (cf. \ref{higgs1-not17}). 
Comme $M$ est $\hoR$-plat et que $\ker(d_{\cC^\dagger})=\hoR$, on en déduit 
un $\hoR[\frac 1 p]$-isomorphisme $\Delta$-équivariant $\mV(\mH(M))\stackrel{\sim}{\rightarrow} M$. 
Le morphisme $\cC^\dagger$-linéaire canonique 
\begin{equation}\label{higgs1-drt4c}
\mV(\mH(M))\otimes_{\hoR}\cC^\dagger\rightarrow \mH(M)\otimes_{\hoR}\cC^\dagger
\end{equation}
s'identifie alors à l'inverse de l'isomorphisme \eqref{higgs1-drt4b}, ce qui montre que $\mH(M)$ est soluble.

\begin{lem}\label{higgs1-drt150}
Soient $V$ un $\hoR[\frac 1 p]$-module projectif de type fini,
$T$ un sous-$\hoR$-module de type fini de $V$ tel que $V=T\otimes_{\mZ_p}\mQ_p$, $r$ un nombre rationnel $\geq 0$. 
Posons $\cM=V\otimes_{\hoR}\hcC^{(r)}$ et notons $\cM^\circ$ 
l'image canonique de $T\otimes_{\hoR}\hcC^{(r)}$ dans $\cM$. 
Alors, le morphisme canonique $V\rightarrow \cM$ est injectif, et il existe un entier $m\geq 0$ tel que $T\subset V\cap \cM^\circ
\subset p^{-m}T$. 
\end{lem}

Choisissons un $\hoR[\frac 1 p]$-module projectif de type fini $V'$ et un $\hoR[\frac 1 p]$-isomorphisme  
$\varphi\colon V\oplus V'\stackrel{\sim}{\rightarrow} (\hoR[\frac 1 p])^n$, où $n$ est un entier $\geq 1$. 
Soit $T'$ un sous-$\hoR$-module de type fini de $V'$ tel que $V'=T'\otimes_{\mZ_p}\mQ_p$. 
Posons $\cM'=V'\otimes_{\hoR}\hcC^{(r)}$ et notons $\cM'^\circ$ 
l'image canonique de $T'\otimes_{\hoR}\hcC^{(r)}$ dans $\cM'$. L'isomorphisme $\varphi$
induit un $(\hcC^{(r)}\otimes_{\mZ_p}\mQ_p)$-isomorphisme $\phi\colon \cM\oplus\cM'\stackrel{\sim}{\rightarrow} (\hcC^{(r)}\otimes_{\mZ_p}\mQ_p)^n$.
Comme l'homomorphisme structural $\hoR\rightarrow \hcC^{(r)}$ est injectif \eqref{higgs1-dolbeault1m}, 
le morphisme canonique $V\rightarrow \cM$ est injectif. Par ailleurs, il existe un entier $j\geq 0$ tel que l'on ait
$p^j\hoR^n\subset \varphi(T\oplus T')\subset p^{-j}\hoR^n$ et donc
$p^j(\hcC^{(r)})^n\subset \phi(\cM^\circ\oplus \cM'^\circ)\subset p^{-j}(\hcC^{(r)})^n$. 
Il résulte aussitôt de \ref{higgs1-tor25} et \eqref{higgs1-dolb43b} que $\hoR[\frac 1 p]\cap \hcC^{(r)}=\hoR$. On en déduit  que 
\begin{equation}
\varphi(T\oplus T')\subset \varphi((V\cap \cM^\circ)\oplus (V'\cap \cM'^\circ))\subset p^{-j}\hoR^n\subset p^{-2j}\varphi(T\oplus T');
\end{equation} 
d'où la proposition. 

\begin{lem}\label{higgs1-drt15}
Soit $(N,\theta)$ un $\hRun[\frac 1 p]$-module de Higgs soluble à coefficients dans $\xi^{-1}\tOmega^1_{R/\co_K}$. 
Alors $\mV(N)$ est une $\hoR[\frac 1 p]$-représentation de Dolbeault de $\Delta$, et on a un isomorphisme
canonique fonctoriel de $\hRun[\frac 1 p]$-modules de Higgs
\begin{equation}\label{higgs1-drt15a}
\mH(\mV(N))\stackrel{\sim}{\rightarrow} N.
\end{equation}
\end{lem}

Montrons d'abord que $\mV(N)$ est une représentation continue de $\Delta$ pour la topologie $p$-adique \eqref{higgs1-not54}. 
D'après \ref{higgs1-drt300}, il existe un nombre rationnel $r>0$ tel que le morphisme $\hcC^{(r)}$-linéaire canonique 
\begin{equation}\label{higgs1-drt15b}
\phi\colon \mV(N) \otimes_{\hoR}\hcC^{(r)}\rightarrow  N\otimes_{\hRun}\hcC^{(r)}
\end{equation}
soit un isomorphisme $\Delta$-équivariant.
Soit $T$ un sous-$\hoR$-module de type fini de $\mV(N)$ tel que $\mV(N)=T\otimes_{\mZ_p}\mQ_p$. 
Posons $\cM=\mV(N)\otimes_{\hoR}\hcC^{(r)}$ et notons $\cM^\circ$ 
l'image canonique de $T\otimes_{\hoR}\hcC^{(r)}$ dans $\cM$. 
Soit $N^\circ$ un sous-$\hRun$-module de type fini de $N$ tel que $N=N^\circ\otimes_{\mZ_p}\mQ_p$.
Posons $\cN=N\otimes_{\hRun}\hcC^{(r)}$ et notons $\cN^\circ$ 
l'image canonique de $N^\circ\otimes_{\hoR}\hcC^{(r)}$ dans  $\cN$. 
Comme $\cM^\circ$ et $\cN^\circ$ sont des $\hoR$-modules de type fini, il existe un entier $n\geq 0$ 
tel que $p^n\cN^\circ \subset \phi(\cM^\circ)\subset p^{-n}\cN^\circ$. 
Par ailleurs, d'après \ref{higgs1-drt150}, le morphisme canonique $\mV(N)\rightarrow \cM$ est injectif et il existe un entier $m\geq 0$ tel que 
$T\subset \mV(N)\cap \cM^\circ \subset p^{-m}T$. On en déduit que 
\begin{equation}\label{higgs1-drt15d}
p^n T\subset \mV(N)\cap \phi^{-1}(\cN^\circ) \subset p^{-m-n}T. 
\end{equation}

Soient $x\in T$ et $\nu$ un entier $\geq 0$. 
D'après \eqref{higgs1-drt15d}, il existe un entier $b\geq 1$, $y_1,\dots,y_b\in N^\circ$ et $\alpha_1,\dots,\alpha_b\in \hcC^{(r)}$ tels que 
$p^n\phi(x)=\sum_{1\leq i\leq b} \alpha_i y_i$. En vertu de \ref{higgs1-dolb45}, il existe un sous-groupe ouvert $\Delta_{x,\nu}$ de $\Delta$
tel que pour tout $g\in \Delta_{x,\nu}$ et tout $1\leq i\leq b$, 
on ait $g(\alpha_i)-\alpha_i\in p^{m+2n+\nu}\hcC^{(r)}$. On en déduit que  
\begin{equation}\label{higgs1-drt15e}
g(x)-x\in  \mV(N)\cap (p^{m+n+\nu}\phi^{-1}(\cN^\circ))\subset p^{\nu} T.
\end{equation}
Comme $T$ est de type fini sur $\hoR$, on en déduit qu'il existe un sous-groupe ouvert $\Delta'$ de $\Delta$
tel que pour tout $g\in \Delta'$, on ait $g(T)\subset T$. Par suite, $T'=\sum_{g\in \Delta}g(T)$ est un sous-$\hoR$-module 
de type fini de $\mV(N)$, stable par $\Delta$. 
Remplaçant $T$ par $T'$, on se réduit au cas où $T$ est stable par l'action de $\Delta$. Il résulte alors de \eqref{higgs1-drt15e}
que l'action de $\Delta$ sur $T$ est continue pour la topologie $p$-adique, et il en est donc de même de l'action de $\Delta$
sur $\mV(N)$.

Le morphisme $\cC^\dagger$-linéaire canonique
\begin{equation}\label{higgs1-drt15c}
\mV(N) \otimes_{\hoR}\cC^\dagger\rightarrow  N\otimes_{\hRun}\cC^\dagger
\end{equation}
est un $\cC^\dagger$-isomorphisme $\Delta$-équivariant 
de $\hoR$-modules de Higgs à coefficients dans $\xi^{-1}\tOmega^1_{R/\co_K}$,
où $\cC^\dagger$ est muni du $\hoR$-champ de Higgs $d_{\cC^\dagger}$,
$N$ est muni de l'action triviale de $\Delta$ et $\mV(N)$ est muni du $\hoR$-champ de Higgs nul (cf. \ref{higgs1-not17}).
Comme $N$ est un facteur direct d'un $\hRun[\frac 1 p]$-module libre de type fini, on a 
$(N\otimes_{\hRun}\cC^\dagger)^\Delta=N$ \eqref{higgs1-dolb47}. 
On en déduit alors un isomorphisme
de $\hRun[\frac 1 p]$-modules de Higgs $\mH(\mV(N))\stackrel{\sim}{\rightarrow} N$.
Le morphisme $\cC^\dagger$-linéaire canonique 
\begin{equation}
\mH(\mV(N))\otimes_{\hRun}\cC^\dagger\rightarrow \mV(N)\otimes_{\hRun}\cC^\dagger
\end{equation}
s'identifie alors à l'inverse de l'isomorphisme \eqref{higgs1-drt15c}, ce qui montre 
que $\mV(N)$ est une $\hoR[\frac 1 p]$-représentation de Dolbeault de $\Delta$.

\begin{rema}
La preuve de \ref{higgs1-drt15} donnée ci-dessus est due à Tsuji. Elle est plus simple et plus élégante 
que la preuve que nous avons donnée dans une première version de cet article.
\end{rema}

\begin{prop}\label{higgs1-drt16}
Les foncteurs $\mH$  et $\mV$ induisent des équivalences de catégories quasi-inverses l'une de l'autre,
entre la catégorie des $\hoR[\frac 1 p]$-représentations de Dolbeault de $\Delta$ et celle des 
$\hRun[\frac 1 p]$-modules de Higgs solubles à coefficients dans $\xi^{-1}\tOmega^1_{R/\co_K}$.
\end{prop}
Cela résulte de \ref{higgs1-drt4} et \ref{higgs1-drt15}. 

\subsection{}\label{higgs1-dolb8}
Soient $M$ une $\hoR[\frac 1 p]$-représentation de Dolbeault de $\Delta$,  
$(\mH(M),\theta)$ le $\hRun[\frac 1 p]$-module de Higgs à coefficients 
dans $\xi^{-1}\tOmega^1_{R/\co_K}$ associé \eqref{higgs1-dolb2c}, 
$\theta_\tot=\theta\otimes \id+\id\otimes d_{\cC^\dagger}$
le $\hRun$-champ de Higgs total sur $\mH(M)\otimes_{\hRun}\cC^\dagger$.
Il résulte aussitôt de \eqref{higgs1-drt2a} qu'on a  un isomorphisme canonique fonctoriel 
de complexes de $\hoR$-représentations
\begin{equation}\label{higgs1-dolb8a}
\mK^\bullet(\mH(M)\otimes_{\hRun}\cC^\dagger,\theta_\tot)\stackrel{\sim}{\rightarrow} 
M\otimes_{\hoR}\mK^\bullet(\cC^\dagger,d_{\cC^\dagger}), 
\end{equation}
où $\mK^\bullet(-,-)$ désigne le complexe de Dolbeault \eqref{higgs1-not6b}. 
Pour tout $i\geq 0$, $\mH(M)\otimes_R\tOmega^i_{R/\co_K}$ 
étant un facteur direct d'un $\hRun[\frac 1 p]$-module libre de type fini, on a, d'après \ref{higgs1-dolb47},
\begin{equation}\label{higgs1-dolb8b}
(\xi^{-i}\mH(M)\otimes_{\hRun}\cC^\dagger\otimes_R\tOmega^i_{R/\co_K})^\Delta=
\xi^{-i}\mH(M)\otimes_R\tOmega^i_{R/\co_K}.
\end{equation} 
On en déduit un isomorphisme canonique fonctoriel de complexes de $\hRun[\frac 1 p]$-modules
\begin{equation}\label{higgs1-dolb8c}
\mK^\bullet(\mH(M),\theta)\stackrel{\sim}{\rightarrow} 
(M\otimes_{\hoR}\mK^\bullet(\cC^\dagger,d_{\cC^\dagger}))^\Delta, 
\end{equation}
où le foncteur $(-)^\Delta$ à droite est défini composante par composante. Ce résultat peut se raffiner 
comme suit.

\begin{prop}[\cite{faltings3} §3, \cite{tsuji3} 5.3.2]\label{higgs1-dolb9}
Soient $M$ une $\hoR[\frac 1 p]$-représentation de Dolbeault de $\Delta$,  
$(\mH(M),\theta)$ le $\hRun[\frac 1 p]$-module de Higgs à coefficients 
dans $\xi^{-1}\tOmega^1_{R/\co_K}$ associé \eqref{higgs1-dolb2c}. 
On a alors un isomorphisme canonique fonctoriel dans $\bD^+(\bMod(\hRun[\frac 1 p]))$ 
\begin{equation}\label{higgs1-dolb9a}
\rC_\cont^\bullet(\Delta, M)\stackrel{\sim}{\rightarrow} \mK^\bullet(\mH(M),\theta),
\end{equation}
où $\rC_\cont^\bullet(\Delta, M)$ est le complexe de cochaînes continues de $\Delta$ à valeurs dans $M$ \eqref{higgs1-not78}
et $\mK^\bullet(\mH(M),\theta)$ est le complexe de Dolbeault \eqref{higgs1-not6b}. 
\end{prop}

En effet, pour tout entier $i\geq 0$ et tous nombres rationnels $r'>r>0$, le morphisme canonique 
\begin{equation}\label{higgs1-dolb9b}
\rC^i_\cont(\Delta, M\otimes_{\hoR} \tmK^\bullet(\hcC^{(r')},p^{r'}d_{\hcC^{(r')}}))\rightarrow 
\rC^i_\cont(\Delta, M\otimes_{\hoR} \tmK^\bullet(\hcC^{(r)},p^rd_{\hcC^{(r)}}))
\end{equation}
est homotope à $0$, en vertu de \ref{higgs1-dolb44}(ii). Par suite, le complexe 
\begin{equation}\label{higgs1-dolb9c}
\underset{\underset{r\in \mQ_{>0}}{\longrightarrow}}{\lim}\  
\rC^i_\cont(\Delta, M\otimes_{\hoR} \tmK^\bullet(\hcC^{(r)},p^rd_{\hcC^{(r)}}))
\end{equation}
est acyclique. La première suite spectrale du bicomplexe (\cite{ega3} § 0 (11.3.2.2))
\begin{equation}\label{higgs1-dolb9d}
\underset{\underset{r\in \mQ_{>0}}{\longrightarrow}}{\lim}\  
\rC^\bullet_\cont(\Delta, M\otimes_{\hoR} \tmK^\bullet(\hcC^{(r)},p^rd_{\hcC^{(r)}}))
\end{equation}
implique alors que le complexe simple associé \eqref{higgs1-not580b}
\begin{equation}\label{higgs1-dolb9e}
\underset{\underset{r\in \mQ_{>0}}{\longrightarrow}}{\lim}\  
\int \rC^\bullet_\cont(\Delta, M\otimes_{\hoR} \tmK^\bullet(\hcC^{(r)},p^rd_{\hcC^{(r)}}))
\end{equation}
est acyclique. Les morphismes canoniques $\hoR[0]\rightarrow \mK^\bullet(\hcC^{(r)},p^rd_{\hcC^{(r)}})$ 
induisent donc un quasi-isomor\-phisme 
\begin{equation}\label{higgs1-dolb9f}
\rC^\bullet_\cont(\Delta, M)\rightarrow \underset{\underset{r\in \mQ_{>0}}{\longrightarrow}}{\lim}\ 
\int \rC^\bullet_\cont(\Delta, M\otimes_{\hoR} \mK^\bullet(\hcC^{(r)},p^rd_{\hcC^{(r)}})).
\end{equation}

D'après \ref{higgs1-drt200}, il existe un nombre rationnel $r_0>0$ tel que $\mH(M)$ soit contenu dans 
$M\otimes_{\hoR}\hcC^{(r_0)}$ et que pour tout nombre rationnel $0<r\leq r_0$, le morphisme canonique
\begin{equation}\label{higgs1-dolb9g}
\mH(M)\otimes_{\hRun}\hcC^{(r)}\rightarrow M\otimes_{\hoR}\hcC^{(r)}
\end{equation}
soit bijectif.
Notant $\theta^{(r)}_\tot=\theta\otimes \id+p^r\id\otimes d_{\hcC^{(r)}}$ le $\hRun$-champ  
de Higgs total sur $\mH(M)\otimes_{\hRun}\hcC^{(r)}$, on en déduit un isomorphisme 
\begin{equation}\label{higgs1-dolb9h}
\mK^\bullet (\mH(M)\otimes_{\hRun}\hcC^{(r)},\theta^{(r)}_\tot)
\stackrel{\sim}{\rightarrow} M\otimes_{\hoR} \mK^\bullet(\hcC^{(r)},p^rd_{\hcC^{(r)}}).
\end{equation}

Pour tout entier $i\geq 0$, le morphisme canonique 
\begin{equation}\label{higgs1-dolb9i}
\mH(M)\otimes_{\hRun}\rC^i_\cont(\Delta, \hcC^{(r)})\rightarrow 
\rC^i_\cont(\Delta,\mH(M)\otimes_{\hRun}\hcC^{(r)})
\end{equation}
est un isomorphisme. En effet, on peut se réduire au cas où 
$\mH(M)$ est un $\hRun[\frac 1 p]$-module libre de type fini, auquel cas l'assertion résulte de la compacité de $\Delta$. 
Par ailleurs, la dérivation $d_{\hcC^{(r)}}$ induit un $\hRun$-champ de Higgs $\delta^{i,(r)}$ sur 
$\rC^i_\cont(\Delta,\hcC^{(r)})$ à coefficients dans $\xi^{-1}\tOmega^1_{R/\co_K}$. 
Les différentielles du complexe $\rC^\bullet_\cont(\Delta,\hcC^{(r)})$ sont des morphismes de modules de Higgs. 
Notant $\vartheta^{i,(r)}_\tot=\theta\otimes \id+p^r\id\otimes \delta^{i,(r)}$
le $\hRun$-champ  de Higgs total sur $\mH(M)\otimes_{\hRun}\rC^i_\cont(\Delta,\hcC^{(r)})$, 
le morphisme canonique 
\begin{equation}\label{higgs1-dolb9j}
\mK^\bullet (\mH(M)\otimes_{\hRun}\rC^i_\cont(\Delta, \hcC^{(r)}),\vartheta^{i,(r)}_\tot)\rightarrow 
\rC^i_\cont(\Delta,\mK^\bullet ( \mH(M)\otimes_{\hRun}\hcC^{(r)},\theta^{(r)}_\tot))
\end{equation}
est un isomorphisme \eqref{higgs1-dolb9i}. 
Il résulte de \ref{higgs1-dolb46}, en utilisant la première suite spectrale des bicomplexes, 
que le morphisme canonique 
\begin{equation}\label{higgs1-dolb9k}
\mK^\bullet(\mH(M),\theta)\rightarrow \underset{\underset{r\in \mQ_{>0}}{\longrightarrow}}{\lim}\ 
\int\mK^\bullet(\mH(M)\otimes_{\hRun} \rC^\bullet_\cont(\Delta, \hcC^{(r)}),\vartheta^{\bullet,(r)}_\tot)
\end{equation}
est un quasi-isomorphisme. La proposition s'ensuit.

\section{Petites représentations}\label{higgs1-RP}

Les hypothèses et notations du §\ref{higgs1-dolbeault} sont en vigueur dans cette section. 

\begin{defi}\label{higgs1-higgs2}\index{Representation petite a@$A$-représentation $\alpha$-quasi-petite, $\alpha$-petite, quasi-petite, petite 
($A$ algèbre $p$-adique)}\index{101301@$\bRep^{\alpha\trqpp}_{A}(G)$, $\bRep^{\qpp}_{A}(G)$, $\bRep^{\alpha\trp}_A(G)$, $\bRep^{\p}_A(G)$}
Soient $G$ un groupe topologique, $A$ une $\co_C$-algèbre complète et séparée pour la topologie $p$-adique,
munie d'une action continue de $G$ (par des homomorphismes de $\co_C$-algèbres), 
$\alpha$ un nombre rationnel $>0$, 
$M$ une $A$-représentation continue de $G$, munie de la topologie $p$-adique.  
\begin{itemize}
\item[{\rm (i)}] On dit que $M$ est {\em $\alpha$-quasi-petite} si le $A$-module $M$ est complet et séparé pour la topologie 
$p$-adique, et est engendré par un nombre fini d'éléments $G$-invariants modulo $p^{\alpha}M$. 
\item[{\rm (ii)}] On dit que $M$ est {\em $\alpha$-petite} si $M$ est un $A$-module libre de type fini
ayant une base sur $A$ formée d'éléments $G$-invariants modulo $p^{\alpha}M$. 
\item[{\rm (iii)}] On dit que $M$ est {\em quasi-petite} (resp. {\em petite}) si elle est $\alpha'$-quasi-petite
(resp. $\alpha'$-petite) pour un nombre rationnel $\alpha'>\frac{2}{p-1}$. 
\end{itemize}

On désigne par $\bRep^{\alpha\trqpp}_{A}(G)$ (resp. $\bRep^{\qpp}_{A}(G)$) la sous-catégorie pleine de 
$\bRep^\cont_A(G)$ formée des $A$-représentations $\alpha$-quasi-petites (resp. quasi-petites) de $G$ 
dont le $A$-module sous-jacent est $\co_C$-plat, et par $\bRep^{\alpha\trp}_A(G)$ (resp. $\bRep^{\p}_A(G)$)
la sous-catégorie pleine de $\bRep^\cont_A(G)$ formée des $A$-représentations $\alpha$-petites 
(resp. petites) de $G$. 
\end{defi}

Si l'action de $G$ sur $A$ est triviale, 
pour qu'une $A$-représentation $M$ de $G$ soit $\alpha$-petite, il faut et il suffit
qu'elle soit $\alpha$-quasi-petite et que $M$ soit un $A$-module libre de type fini.

\begin{defi}\label{higgs1-drt6}\index{Representation petite b@$A[\frac 1 p]$-représentation petite ($A$ algèbre $p$-adique)}
\index{101303@$\bRep^{\p}_{A[\frac 1 p]}(G)$}
Soient $G$ un groupe topologique, $A$ une $\co_C$-algèbre complète et séparée pour la topologie $p$-adique,
munie d'une action continue de $G$ (par des homomorphismes de $\co_C$-algèbres).
On munit $A[\frac 1 p]$ de la topologie $p$-adique \eqref{higgs1-not54}.
On dit qu'une $A[\frac 1 p]$-représentation continue $M$ de $G$ est {\em petite} si
les conditions suivantes sont remplies~: 
\begin{itemize}
\item[(i)] $M$ est un $A[\frac 1 p]$-module projectif de type fini, muni de la topologie $p$-adique \eqref{higgs1-not54}~;
\item[(ii)] il existe un nombre rationnel $\alpha>\frac{2}{p-1}$
et un sous-$A$-module de type fini $M^\circ$ de $M$, stable par $G$, 
engendré par un nombre fini d'éléments $G$-invariants modulo $p^{\alpha}M^\circ$,
et qui engendre $M$ sur $A[\frac 1 p]$. 
\end{itemize}
On désigne par $\bRep^{\p}_{A[\frac 1 p]}(G)$
la sous-catégorie pleine de $\bRep^\cont_{A[\frac 1 p]}(G)$ formée des $A$-représentations petites de $G$. 
\end{defi}

On notera que contrairement au cas entier \eqref{higgs1-higgs2}, on ne demande pas que  
$M$ soit libre sur $A[\frac 1 p]$.

\begin{remas}\label{higgs1-drt65}
Soient $G$ un groupe topologique, $A$ une $\co_C$-algèbre complète et séparée pour la topologie $p$-adique,
munie d'une action continue de $G$ (par des homomorphismes de $\co_C$-algèbres).
\begin{itemize}
\item[{\rm (i)}] Soient $M$ un $A[\frac 1 p]$-module projectif de type fini, $M^\circ$ un sous-$A$-module de type fini
de $M$. Alors $M^\circ$ est complet et séparé pour la topologie $p$-adique. En effet, $M^\circ$ est complet en vertu de 
(\cite{ac} chap. III §2.12 cor.~1 de prop.~16). D'autre part, quitte à ajouter à $M$ un facteur direct, 
on peut le supposer libre de type fini sur $A[\frac 1 p]$. Par suite, il existe un entier $m\geq 0$ tel que $p^m M^\circ$
soit contenu dans un $A$-module libre de type fini $N$. Donc $\cap_{n\geq 0}p^nM^\circ
\subset \cap_{n\geq 0}p^nN=0$. 
\item[{\rm (ii)}] Soient $M$ une $A[\frac 1 p]$-représentation petite de $G$, 
$M^\circ$ un sous-$A$-module de type fini de $M$ vérifiant la condition \ref{higgs1-drt6}(ii). 
Il résulte alors de (i) que $M^\circ$ est une $A$-représentation quasi-petite de $G$. 
\end{itemize}
\end{remas}

\begin{defi}\label{higgs1-higgs4}\index{Module de Higgs petit a@$\hRun$-module de Higgs $\beta$-quasi-petit, $\beta$-petit, quasi-petit, petit}
\index{101305@$\bMH^{\beta\trp}(\hRun,\xi^{-1}\tOmega^1_{R/\co_K})$}
Soient $\beta$ un nombre rationnel $>0$, $(N,\theta)$ un $\hRun$-module de Higgs  
à coefficients dans $\xi^{-1}\tOmega^1_{R/\co_K}$ \eqref{higgs1-eip75}.
\begin{itemize} 
\item[{\rm (i)}] On dit que $(N,\theta)$ est {\em $\beta$-quasi-petit} si $N$ est de type fini sur $\hRun$
et si $\theta$ est un multiple de  $p^{\beta}$ dans $\xi^{-1}\End_{\hRun}(N)\otimes_R\tOmega^1_{R/\co_K}$
\eqref{higgs1-not65}. On dit alors aussi que le $\hRun$-champ de Higgs $\theta$
est {\em $\beta$-quasi-petit}. 
\item[{\rm (ii)}] On dit que $(N,\theta)$ est {\em $\beta$-petit} s'il est $\beta$-quasi-petit et 
si $N$ est libre de type fini sur $\hRun$. On dit alors aussi que le $\hRun$-champ de Higgs $\theta$
est {\em $\beta$-petit}. 
\item[{\rm (iii)}] On dit que $(N,\theta)$ est {\em quasi-petit} (resp. {\em petit}) s'il est $\beta'$-quasi-petit
(resp. $\beta'$-petit) pour un nombre rationnel $\beta'>\frac{1}{p-1}$. 
On dit alors aussi que le $\hRun$-champ de Higgs $\theta$ est {\em quasi-petit} (resp. {\em petit}). 
\end{itemize}

On désigne par $\bMH^{\beta\trqpp}(\hRun,\xi^{-1}\tOmega^1_{R/\co_K})$ 
la sous-catégorie pleine de $\bMH(\hRun,\xi^{-1}\tOmega^1_{R/\co_K})$  \eqref{higgs1-eip75}
formée des $\hRun$-modules de Higgs $\beta$-quasi-petits dont le $\hRun$-module sous-jacent est $\co_C$-plat
et par $\bMH^{\beta\trp}(\hRun,\xi^{-1}\tOmega^1_{R/\co_K})$ la sous-catégorie pleine de 
$\bMH(\hRun,\xi^{-1}\tOmega^1_{R/\co_K})$ 
formée des $\hRun$-modules de Higgs $\beta$-petits.
\end{defi}

\begin{defi}\label{higgs1-drt7}\index{Module de Higgs petit b@$\hRun[\frac 1 p]$-module de Higgs petit}
\index{101307@$\bMH^{\p}(\hRun[\frac 1 p],\xi^{-1}\tOmega^1_{R/\co_K})$}
On dit qu'un $\hRun[\frac 1 p]$-module de Higgs $(N,\theta)$ à coefficients dans $\xi^{-1}\tOmega^1_{R/\co_K}$ 
est {\em petit} si les conditions suivantes sont remplies~: 
\begin{itemize}
\item[(i)] $N$ est un $\hRun[\frac 1 p]$-module projectif de type fini~; 
\item[(ii)] il existe un nombre rationnel $\beta>\frac{1}{p-1}$ et 
un sous-$\hRun$-module de type fini $N^\circ$ de $N$, qui l'engendre sur $\hRun[\frac 1 p]$, tels que l'on ait 
\begin{equation}\label{higgs1-drt7a}
\theta(N^\circ)\subset p^{\beta}\xi^{-1}N^\circ \otimes_R\tOmega^1_{R/\co_K}.
\end{equation}
\end{itemize}
On désigne par $\bMH^{\p}(\hRun[\frac 1 p],\xi^{-1}\tOmega^1_{R/\co_K})$ la sous-catégorie pleine de 
$\bMH(\hRun[\frac 1 p],\xi^{-1}\tOmega^1_{R/\co_K})$ 
formée des $\hRun[\frac 1 p]$-modules de Higgs petits (cf. \ref{higgs1-eip75}). 
\end{defi}

On notera que contrairement au cas entier \eqref{higgs1-higgs4}, on ne demande pas que $N$ soit libre sur $\hRun[\frac 1 p]$.

\begin{lem}\label{higgs1-drt75}
Soient $A$ un anneau intègre de corps des fractions $L$, $B$ un sous-anneau de $L$ contenant $A$,
$N$ un $B$-module plat de type fini, $u_1,\dots,u_\ell$ des endomorphismes $B$-linéaires de $N$
qui commutent deux à deux et tels que pour tout $1\leq i\leq \ell$, le polynôme caractéristique de
l'endomorphisme $u_i\otimes \id$ de $N\otimes_BL$ soit à coefficients dans $A$. 
Alors $N$ est engendré sur $B$ par un sous-$A$-module de type fini $M$ tel que $u_i(M)\subset M$
pour tout $1\leq i\leq \ell$. 
\end{lem}

On procède par récurrence sur le nombre d'endomorphismes. 
Supposons l'assertion établie pour $\ell-1$ et montrons-là pour $\ell$. 
Soit $M'$ un sous-$A$-module de type fini de $N$ qui l'engendre sur $B$ et
tel que $u_i(M')\subset M'$ pour tout $1\leq i\leq \ell-1$. Soit 
$P(X)=X^n+a_1X^{n-1}+\dots+a_n\in A[X]$ le polynôme caractéristique de l'endomorphisme
$u_\ell\otimes \id$ de $N\otimes_BL$. Comme $N$ est $B$-plat,
on peut l'identifier à un sous-$B$-module de $N\otimes_BL$. 
Par suite, l'endomorphisme $P(u_\ell)$ de $N$ est nul, 
et le sous-$A$-module $M=\sum_{i=0}^{n-1}u^i_\ell(M')$ de $N$ répond à la question. 

\begin{lem}\label{higgs1-drt76}
Soit  $(N,\theta)$ un $\hRun[\frac 1 p]$-module de Higgs à coefficients dans $\xi^{-1}\tOmega^1_{R/\co_K}$ 
tel que les conditions suivantes soient remplies~: 
\begin{itemize}
\item[{\rm (i)}] $N$ est un $\hRun[\frac 1 p]$-module projectif de type fini~; 
\item[{\rm (ii)}] il existe un nombre rationnel $\beta>\frac{1}{p-1}$ tel que pour tout $i\geq 1$, 
le $i$-ième invariant caractéristique de $\theta$ appartienne à 
$p^{i\beta}\xi^{-i}\rS_R^i(\tOmega^1_{R/\co_K})\otimes_R\hRun$ \eqref{higgs1-not62}. 
\end{itemize}
Alors $(N,\theta)$ est petit. 
\end{lem}

Comme $\Spec(\hR_1)$ est $\co_C$-plat, 
ses points génériques sont les points génériques de $\Spec(\hR_1[\frac 1 p])$.
Comme ce dernier schéma est noethérien, $\Spec(\hR_1)$ n'a qu'un nombre fini de points génériques.
D'autre part, $\hR_1$ est normal en vertu de \ref{higgs1-pur81}. Il est donc isomorphe à un produit fini 
$\prod_{j=1}^n A_j$ d'anneaux normaux et intègres. 
D'après \ref{higgs1-drt75}, pour tout $1\leq j\leq n$, il existe un 
sous-$A_j$-module de type fini $N_j^\circ$ de $N\otimes_{\hR_1}A_j$ qui l'engendre sur $A_j[\frac 1 p]$
tel que si $\theta_j$ désigne le $A_j[\frac 1 p]$-champ de Higgs sur $N\otimes_{\hR_1}A_j$ 
à coefficients dans $\xi^{-1}\tOmega^1_{R/\co_K}$ induit par $\theta$, on ait 
\begin{equation}
\theta_j(N_j^\circ)\subset p^{\beta}\xi^{-1}N^\circ_j\otimes_R\tOmega^1_{R/\co_K}.
\end{equation}
La proposition s'ensuit.

\subsection{}\label{higgs1-higgs5}
Rappelons que sur $C$, la fonction logarithme $\log(x)$ converge lorsque $x\in 1+\fm_C$, et que la fonction 
exponentielle $\exp(x)$ converge lorsque $v(x)>\frac{1}{p-1}$. Pour tout $x\in C$ tel que $v(x)>\frac{1}{p-1}$, on a 
\begin{eqnarray*}
\exp(x)&\equiv& 1+x \mod(x\fm_C),\\
\log(1+x)&\equiv& x \mod(x\fm_C),
\end{eqnarray*}
$\exp(\log(1+x))=1+x$ et $\log(\exp(x))=x$.

\subsection{}\label{higgs1-higgs55}
Soit $M$ un $\hRun$-module de type fini. 
D'après (\cite{egr1} 1.10.2), $M$ est complet et séparé pour la topologie $p$-adique, 
et il en est alors de même de $\End_{\hRun}(M)$. Supposons, de plus, $M$ $\co_C$-plat.
Alors $\End_{\hRun}(M)$ est $\co_C$-plat et pour tout nombre rationnel $\alpha\geq 0$, 
l'homomorphisme canonique $p^{\alpha}\End_{\hRun}(M)\rightarrow \Hom_{\hRun}(M,p^{\alpha}M)$ est un isomorphisme.
Soient $u$ un endomorphisme $\hRun$-linéaire de $M$, $\alpha$ un nombre rationnel $>0$. 
Si $u$ induit un automorphisme de $M/p^\alpha M$,  
alors $u$ est un automorphisme de $M$. En effet, $u$ est surjectif par le lemme de Nakayama.
Si $x\in M$ est tel que $u(x)=0$, alors il existe $y\in M$ tel que $x=p^{\alpha}y$; comme $M$ est $\co_C$-plat,
on a $u(y)=0$; on en déduit que $x\in \cap_{n\geq 0}p^{n\alpha}M=0$, d'où l'assertion. 
On peut donc identifier $\id+p^\alpha \End_{\hRun}(M)$ à un sous-groupe de $\Aut_{\hRun}(M)$.
Supposons $\alpha>\frac{1}{p-1}$ et  $v=u-\id\in p^\alpha \End_{\hRun}(M)$. Pour tout $x\in M$, les séries 
\begin{eqnarray}
\exp(v)(x)&=&\sum_{n\geq 0} \frac{1}{n!}v^n(x),\\
\log(u)(x)&=&\sum_{n\geq 0} \frac{(-1)^n}{n+1}v^{n+1}(x),
\end{eqnarray}
convergent alors dans $M$, et définissent deux $\hRun$-endomorphismes $\exp(v)$ et $\log(u)$ de $M$. 
De plus, on a $\exp(v)\in \id + p^{\alpha}\End_{\hRun}(M)$, $\log(u) \in p^{\alpha}\End_{\hRun}(M)$ et $\exp(\log(u))=u$
et $\log(\exp(v))=v$.

\subsection{}\label{higgs1-higgs6}
Soient $M$ un $\hRun$-module de type fini et $\co_C$-plat, 
$\alpha$ un nombre rationnel $>\frac{1}{p-1}$, $\beta=\alpha-\frac{1}{p-1}$. 
On désigne par $\Psi_M$ l'isomorphisme composé
\begin{equation}\label{higgs1-higgs6d}
\xymatrix{
{\Hom_{\mZ}(\Delta_\infty, p^{\alpha}\End_{\hRun}(M))}\ar[r]^-(0.5)\sim\ar[dr]_{\Psi_M}&
{p^{\beta}\xi^{-1}\End_{\hRun}(M)\otimes_{\hRun}\Hom_{\mZ}(\Delta_\infty,\hRun(1))}
\ar[d]^{\id \otimes\tdelta}\\
&{p^{\beta}\xi^{-1}\End_{\hRun}(M)\otimes_{R}\tOmega^1_{R/\co_K}}}
\end{equation}
où $\tdelta$ est l'isomorphisme \eqref{higgs1-ext-log14h} et 
l'isomorphisme horizontal provient de \eqref{higgs1-higgs551b} et de l'isomorphisme canonique
$\co_C(1)\stackrel{\sim}{\rightarrow} p^{\frac{1}{p-1}}\xi\co_C$ \eqref{higgs1-ext11c}.
On notera que le $\hRun$-module $\End_{\hRun}(M)$ est complet et séparé pour la topologie $p$-adique 
et $\co_C$-plat \eqref{higgs1-higgs55}.
On peut décrire explicitement $\Psi_M$ comme suit. Soient $t_1,\dots,t_d\in P^\gp$ tels que leurs images  
dans $(P^\gp/\mZ\lambda)\otimes_\mZ\mZ_p$ forment une $\mZ_p$-base,
$(\chi_{t_i})_{1\leq i\leq d}$ leurs images dans $\Hom_{\mZ}(\Delta_\infty,\mZ_p(1))$ \eqref{higgs1-ext-log14e},
$(d\log(t_i))_{1\leq i\leq d}$ leurs images dans $\tOmega^1_{R/\co_K}$ \eqref{higgs1-ext-log13e}. 
Les  $(d\log(t_i))_{1\leq i\leq d}$ forment alors une $R$-base de $\tOmega^1_{R/\co_K}$ \eqref{higgs1-log-ext45b}
et les $(\chi_{t_i})_{1\leq i\leq d}$ forment une $\mZ_p$-base de 
$\Hom_\mZ(\Delta_\infty,\mZ_p(1))$ \eqref{higgs1-ext-log14gh}. 
Pour tout $1\leq i\leq d$, notons $\chi_i$ l'homomorphisme composé 
\begin{equation}\label{higgs1-higgs6gg}
\xymatrix{
{\Delta_\infty}\ar[r]^-(0.5){\chi_{t_i}}& 
{\mZ_p(1)}\ar[r]^-(0.5){\log([\ ])}&{p^{\frac{1}{p-1}}\xi\co_C}},
\end{equation}
où la seconde flèche est induite par l'isomorphisme \eqref{higgs1-ext11c}. 
Compte tenu de \eqref{higgs1-higgs551b}, pour tout homomorphisme
$\phi\colon \Delta_\infty \rightarrow p^{\alpha}\End_{\hRun}(M)$, 
il existe $\phi_i\in p^{\beta}\End_{\hRun}(M)$ $(1\leq i\leq d)$ tels que 
\begin{equation}\label{higgs1-higgs6bb}
\phi=\sum_{i=1}^d\xi^{-1}\phi_i\otimes \chi_i.
\end{equation}
Concrètement, si $\zeta$ est une $\mZ_p$-base de $\mZ_p(1)$, 
il existe des éléments $\gamma_1,\dots,\gamma_d\in \Delta_\infty$ tels que pour tous $1\leq i,j\leq d$, on ait 
$\chi_{t_i}(\gamma_j)=\delta_{ij}\zeta$. Pour tout $1\leq i\leq d$, on a donc
\begin{equation}\label{higgs1-higgs6dd}
\xi^{-1}\log([\zeta])\phi_i=\phi(\gamma_i).
\end{equation} 
On rappelle que $\xi^{-1}\log([\zeta])$ est un élément de valuation $\frac{1}{p-1}$ de $\co_C$ \eqref{higgs1-ext111}. 
Comme $\tdelta(\chi_{t_i})=d\log(t_i)$ pour tout $1\leq i\leq d$, on a  
\begin{equation}\label{higgs1-higgs6b}
\Psi_M(\phi)=\sum_{i=1}^d\xi^{-1}\phi_i \otimes d\log(t_i).
\end{equation}

Soit $\varphi$ une $\hRun$-représentation $\alpha$-quasi-petite de $\Delta_\infty$ sur $M$.
Comme $\Delta_\infty$ agit trivialement sur $\hRun$, $\varphi$ est un homomorphisme 
\begin{equation}\label{higgs1-higgs6aa}
\varphi\colon \Delta_\infty\rightarrow \Aut_{\hRun}(M)
\end{equation}
d'image contenue dans le sous-groupe $\id+p^{\alpha}\End_{\hRun}(M)$ de $\Aut_{\hRun}(M)$.
Comme $\Delta_\infty$ est abélien, on peut définir l'homomorphisme \eqref{higgs1-higgs55}
\begin{equation}\label{higgs1-higgs6a}
\log(\varphi)\colon \Delta_\infty\rightarrow p^{\alpha}\End_{\hRun}(M).
\end{equation}
De plus, il résulte de \eqref{higgs1-higgs6dd} et \eqref{higgs1-higgs6b} que $\Psi_M(\log(\varphi))\wedge \Psi_M(\log(\varphi))=0$, 
autrement dit, $\Psi_M(\log(\varphi))$ est un $\hRun$-champ de Higgs 
$\beta$-quasi-petit sur $M$ à coefficients dans $\xi^{-1}\tOmega^1_{R/\co_K}$. 
On obtient ainsi un foncteur 
\begin{equation}\label{higgs1-higgs6c}
\bRep_{\hRun}^{\alpha\trqpp}(\Delta_\infty)\rightarrow \bMH^{\beta\trqpp}(\hRun,\xi^{-1}\tOmega^1_{R/\co_K}),\ \ \ 
(M,\varphi)\mapsto (M,\Psi_M(\log(\varphi))).
\end{equation}

Soit $\theta$ un $\hRun$-champ de Higgs $\beta$-quasi-petit sur $M$ à coefficients dans 
$\xi^{-1}\tOmega^1_{R/\co_K}$. Compte tenu de \eqref{higgs1-higgs6bb} et \eqref{higgs1-higgs6dd},
comme $\theta\wedge \theta=0$,
l'image de l'homomorphisme $\Psi^{-1}_M(\theta)\colon \Delta_\infty \rightarrow p^\alpha\End_{\hRun}(M)$ 
est formée d'endomorphismes de $M$ qui commutent deux à deux. 
On peut donc définir l'homomorphisme \eqref{higgs1-higgs55}
\begin{equation}\label{higgs1-higgs7a}
\exp(\Psi^{-1}_M(\theta))\colon \Delta_\infty \rightarrow \Aut_{\hRun}(M),
\end{equation}
qui est clairement une $\hRun$-représentation $\alpha$-quasi-petite de $\Delta_\infty$ sur $M$. 
On définit ainsi un foncteur 
\begin{equation}\label{higgs1-higgs6cc}
\bMH^{\beta\trqpp}(\hRun,\xi^{-1}\tOmega^1_{R/\co_K})\rightarrow \bRep_{\hRun}^{\alpha\trqpp}(\Delta_\infty),\ \ \ 
(M,\theta)\mapsto (M,\exp(\Psi^{-1}_M(\theta))). 
\end{equation}

\begin{prop}[Faltings, \cite{faltings3}]\label{higgs1-higgs7}
Pour tous nombres rationnels $\alpha>\frac{1}{p-1}$ et $\beta=\alpha-\frac{1}{p-1}$,
les foncteurs \eqref{higgs1-higgs6c} et \eqref{higgs1-higgs6cc} 
sont des équivalences de catégories quasi-inverses l'une de l'autre. Ils induisent des équivalences de catégories 
quasi-inverses l'une de l'autre
\begin{eqnarray}
\bRep_{\hRun}^{\alpha\trp}(\Delta_\infty)\stackrel{\sim}{\rightarrow}\bMH^{\beta\trp}(\hRun,\xi^{-1}\tOmega^1_{R/\co_K}),\label{higgs1-higgs7b}\\
\bMH^{\beta\trp}(\hRun,\xi^{-1}\tOmega^1_{R/\co_K})\stackrel{\sim}{\rightarrow}\bRep_{\hRun}^{\alpha\trp}(\Delta_\infty).
\label{higgs1-higgs7bb}
\end{eqnarray}
\end{prop}

La première assertion résulte aussitôt de la définition des foncteurs \eqref{higgs1-higgs6c} et \eqref{higgs1-higgs6cc}.
La seconde assertion s'ensuit puisque $\hRun$ est $\co_C$-plat \eqref{higgs1-pur8}. 

\begin{rema}\label{higgs1-higgs75}\index{101320@$\varphi=\exp\left(\sum_{i=1}^d\xi^{-1}\theta_i\otimes \chi_{i}\right)$}
\index{101321@$\theta=\sum_{i=1}^d\xi^{-1}\theta_i\otimes d\log(t_i)$}
Sous les hypothèses de \ref{higgs1-higgs6}, 
si $(M,\varphi)$ est un objet de $\bRep_{\hRun}^{\alpha\trqpp}(\Delta_\infty)$ et $(M,\theta)$
est un objet de $\bMH^{\beta\trqpp}(\hRun,\xi^{-1}\tOmega^1_{R/\co_K})$ qui se correspondent par les foncteurs \eqref{higgs1-higgs6c}
et \eqref{higgs1-higgs6cc}, alors $\varphi$ et $\theta$ sont reliés par les formules suivantes~:
\begin{eqnarray}
\varphi&=&\exp\left(\sum_{i=1}^d\xi^{-1}\theta_i\otimes \chi_{i}\right),\label{higgs1-higgs75a}\\
\theta&=&\sum_{i=1}^d\xi^{-1}\theta_i\otimes d\log(t_i).\label{higgs1-higgs75b}
\end{eqnarray}
Comme les $\chi_i$ se factorisent à travers $\Delta_{p^\infty}$ \eqref{higgs1-ext-log14},
$\varphi$ se factorise à travers $\Delta_{p^\infty}$.
\end{rema}

\subsection{}\label{higgs1-higgs99}\index{101325@$\fS^{(r)}$}
Pour tout nombre rationnel $r\geq 0$, on désigne par $\fS^{(r)}$ la sous-$\hRun$-algèbre de $\cS^{(r)}$ \eqref{higgs1-chb1a}
définie par
\begin{equation}\label{higgs1-higgs99a}
\fS^{(r)}=\rS_{\hRun}(p^r\xi^{-1}\tOmega^1_{R/\co_K}\otimes_R\hRun),
\end{equation}
et par $\hfS^{(r)}$ son séparé complété $p$-adique, que l'on munit de la topologie $p$-adique. 
On notera que $\hfS^{(r)}$ est $\hRun$-plat en vertu de (\cite{egr1} 1.12.4) et donc $\co_C$-plat \eqref{higgs1-pur8}. 
On pose $\fS=\fS^{(0)}$ et $\hfS=\hfS^{(0)}$. 
Pour tous nombres rationnels $r'\geq r\geq 0$, on a un homomorphisme injectif canonique 
$\tta^{r,r'}\colon \fS^{(r')}\rightarrow \fS^{(r)}$. On vérifie aussitôt que l'homomorphisme 
induit $\htta^{r,r'}\colon\hfS^{(r')}\rightarrow \hfS^{(r)}$ est injectif. On a un $\fS^{(r)}$-isomorphisme canonique 
\begin{equation}
\Omega^1_{\fS^{(r)}/\hRun}\stackrel{\sim}{\rightarrow}\xi^{-1}\tOmega^1_{R/\co_K}\otimes_R\fS^{(r)}.
\end{equation} 
On note
\begin{equation}\label{higgs1-higgs99b}
d_{\fS^{(r)}}\colon \fS^{(r)}\rightarrow \xi^{-1}\tOmega^1_{R/\co_K}\otimes_R\fS^{(r)}
\end{equation}
la $\hRun$-dérivation universelle de $\fS^{(r)}$, qui est également un $\hRun$-champ de Higgs à coefficients dans 
$\xi^{-1}\tOmega^1_{R/\co_K}$ puisque $\xi^{-1}\tOmega^1_{R/\co_K}\otimes_R \hRun\subset d_{\fS^{(r)}}(\fS^{(r)})$.
On désigne par 
\begin{equation}\label{higgs1-higgs99d}
d_{\hfS^{(r)}}\colon \hfS^{(r)}\rightarrow \xi^{-1}\tOmega^1_{R/\co_K}\otimes_R\hfS^{(r)}
\end{equation}
son prolongement aux complétés. Pour tous nombres rationnels $r'\geq r\geq 0$, on a 
\begin{equation}\label{higgs1-higgs99e}
p^{r'-r}(\id \times \tta^{r,r'}) \circ d_{\fS^{(r')}}=d_{\fS^{(r)}}\circ \tta^{r,r'}.
\end{equation}

On désigne par $(\tX_0,\cM_{\tX_0})$ la $(\cA_2(\oS),\cM_{\cA_2(\oS)})$-déformation lisse de 
$(\coX,\cM_\coX)$ définie par la carte $(P,\gamma)$ \eqref{higgs1-ext24ac}, 
par $\cL_0$ le torseur de Higgs-Tate associé \eqref{higgs1-tor2}, 
par $\psi_0\in \cL_0(\hY)$ la section définie par la même carte \eqref{higgs1-ext24b}
et par $\varphi_0=\varphi_{\psi_0}$ l'action de $\Delta$ sur $\cS$ induite par $\psi_0$ \eqref{higgs1-tor202g}. 
D'après \ref{higgs1-ext26}, $\varphi_0$ préserve $\fS$ et l'action induite de $\Delta$ 
sur $\fS$ se factorise à travers $\Delta_{p^\infty}$. On désigne encore par $\varphi_0$ l'action 
de $\Delta_{p^\infty}$ sur $\fS$ ainsi définie. 

Soient $t_1,\dots,t_d\in P^\gp$ tels que leurs images dans $(P^\gp/\mZ\lambda)\otimes_\mZ\mZ_p$ 
forment une $\mZ_p$-base, 
$(\chi_{t_i})_{1\leq i\leq d}$ leurs images dans $\Hom_\mZ(\Delta_{p^\infty},\mZ_p(1))$ \eqref{higgs1-ext-log14e},
$(d\log(t_i))_{1\leq i\leq d}$ leurs images dans $\tOmega^1_{R/\co_K}$ \eqref{higgs1-ext-log13e}. 
Pour tout $1\leq i\leq d$, on pose $y_i=\xi^{-1}d\log(t_i)\in \xi^{-1}\tOmega^1_{R/\co_K}\subset \fS$,
et on note $\chi_i$ l'homomorphisme composé 
\begin{equation}\label{higgs1-higgs99f}
\xymatrix{
{\Delta_{p^\infty}}\ar[r]^-(0.5){\chi_{t_i}}& 
{\mZ_p(1)}\ar[r]^-(0.5){\log([\ ])}&{p^{\frac{1}{p-1}}\xi\co_C}},
\end{equation}
où la seconde flèche est induite par l'isomorphisme \eqref{higgs1-ext11c}. 
D'après \eqref{higgs1-ext26b}, pour tout $g\in \Delta_{p^\infty}$, on a 
\begin{equation}\label{higgs1-higgs99g}
\varphi_0(g)=\exp(-\sum_{i=1}^d\xi^{-1}\frac{\partial}{\partial y_i}\otimes \chi_i(g)).
\end{equation}
Il s'ensuit que $\varphi_0$ est continue pour la topologie $p$-adique sur $\fS$ (cf. \ref{higgs1-cohbis10}).
Pour tout nombre rationnel $r\geq 0$, $\varphi_0$ préserve $\fS^{(r)}$, 
et les actions induites de $\Delta_{p^\infty}$ sur $\fS^{(r)}$ et $\hfS^{(r)}$ 
sont continues pour les topologies $p$-adiques (cf. \ref{higgs1-cohbis3}). 
Sauf mention explicite du contraire, on munit $\fS^{(r)}$ et $\hfS^{(r)}$ de ces actions. 
Il résulte aussitôt de \eqref{higgs1-higgs99g} que $d_{\fS^{(r)}}$ et $d_{\hfS^{(r)}}$ sont $\Delta_{p^\infty}$-équivariants.

\subsection{}\label{higgs1-higgs10}\index{101330@$\exp_r(\theta)$}
Conservons les notations de \ref{higgs1-higgs99}, 
soient, de plus, $\beta$, $r$ deux nombres rationnels tels que $\beta>r+\frac{1}{p-1}$ et $r\geq 0$, 
$(N,\theta)$ un $\hRun$-module de Higgs $\beta$-quasi-petit à coefficients dans $\xi^{-1}\tOmega^1_{R/\co_K}$ \eqref{higgs1-higgs4}
tel que $N$ soit $\co_C$-plat. On peut écrire de manière unique 
\begin{equation}\label{higgs1-higgs10a}
\theta=\sum_{i=1}^d\theta_i \otimes y_i,
\end{equation}
où les $\theta_i$ sont des endomorphismes de $N$ appartenant à $p^{\beta}\End_{\hRun}(N)$  
et commutant deux à deux. Pour tout $\un=(n_1,\dots,n_d)\in \mN^d$, 
posons $|\un|=\sum_{i=1}^d n_i$, $\un!=\prod_{i=1}^d n_i!$, $\utheta^{\un}=\prod_{i=1}^d\theta_i^{n_i} \in \End_{\hRun}(N)$
et $\uy^{\un}=\prod_{i=1}^dy_i^{n_i}\in \fS$. 
On notera que $N\otimes_{\hRun}\hfS^{(r)}$ est complet
et séparé pour la topologie $p$-adique (\cite{egr1} 1.10.2), et qu'il est $\co_C$-plat
puisque $\hfS^{(r)}$ est $\hRun$-plat (\cite{egr1} 1.12.4). 
Par suite, pour tout $z\in N\otimes_{\hRun}\hfS^{(r)}$, la série 
\begin{equation}\label{higgs1-higgs10b}
\sum_{\un\in \mN^d} \frac{1}{\un!}(\utheta^{\un}\otimes\uy^{\un})(z)
\end{equation}
converge dans $N\otimes_{\hRun}\hfS^{(r)}$, et définit un endomorphisme $\hfS^{(r)}$-linéaire 
de $N\otimes_{\hRun} \hfS^{(r)}$, que l'on note 
\begin{equation}\label{higgs1-higgs10c}
\exp_r(\theta)\colon N\otimes_{\hRun} \hfS^{(r)}\rightarrow N\otimes_{\hRun} \hfS^{(r)}.
\end{equation}
Pour tout nombre rationnel $r'$ tel que $0\leq r'\leq r$, le diagramme 
\begin{equation}
\xymatrix{
{N\otimes_{\hRun} \hfS^{(r)}}\ar[rr]^-(0.5){\exp_r(\theta)}\ar[d]_-(0.4){\id \otimes \tta^{r',r}}&&
{N\otimes_{\hRun} \hfS^{(r)}}\ar[d]^-(0.4){\id \otimes \tta^{r',r}}\\
{N\otimes_{\hRun} \hfS^{(r')}}\ar[rr]^-(0.5){\exp_{r'}(\theta)}&&{N\otimes_{\hRun} \hfS^{(r')}}}
\end{equation}
est commutatif. On peut donc se permettre d'omettre l'indice $r$ de la notation $\exp_r(\theta)$ sans risque d'ambiguïté. 

\begin{prop}\label{higgs1-higgs11}
Soient $\beta$, $r$ deux nombres rationnels tels que $\beta>r+\frac{1}{p-1}$ et $r\geq 0$, 
$N$ un $\hRun$-module de type fini et $\co_C$-plat, 
$\theta$ un $\hRun$-champ de Higgs $\beta$-quasi-petit sur $N$ à coefficients dans $\xi^{-1}\tOmega^1_{R/\co_K}$,
$\varphi$ la $\hRun$-représentation quasi-petite de $\Delta_\infty$ sur $N$ associée à $\theta$ par le foncteur \eqref{higgs1-higgs6cc}. 
Alors, le $\hfS^{(r)}$-module $N\otimes_{\hRun} \hfS^{(r)}$ est complet et séparé pour la topologie $p$-adique,
et l'endomorphisme \eqref{higgs1-higgs10c}
\begin{equation}\label{higgs1-higgs11a}
\exp_r(\theta)\colon N\otimes_{\hRun} \hfS^{(r)}\rightarrow N\otimes_{\hRun} \hfS^{(r)}
\end{equation}
est un isomorphisme $\Delta_\infty$-équivariant de modules à $p^r$-connexions $p$-adiques intégrables relativement à l'extension 
$\hfS^{(r)}/\hRun$ \eqref{higgs1-not20}, où $\hfS^{(r)}$ est muni de l'action de $\Delta_{\infty}$ induite par $\varphi_0$ \eqref{higgs1-higgs99g}
et de la $p^r$-connexion $p$-adique $p^rd_{\hfS^{(r)}}$ \eqref{higgs1-higgs99d}, 
le module $N$ de la source est muni de l'action triviale de $\Delta_{\infty}$ 
et du $\hRun$-champ de Higgs $\theta$, et le module $N$ du but est muni de l'action $\varphi$
de $\Delta_{\infty}$ et du $\hRun$-champ de Higgs nul (cf. \ref{higgs1-not23}). 
En particulier, $\exp_r(\theta)$ est un isomorphisme de $\hRun$-modules de Higgs à coefficients dans
$\xi^{-1}\tOmega^1_{R/\co_K}$ \eqref{higgs1-not22}. 
\end{prop}

On observera d'abord que $N\otimes_{\hRun} \hfS^{(r)}$ est complet et séparé pour la topologie $p$-adique en vertu de 
(\cite{egr1} 1.10.2(ii)). Soient $x_1,\dots,x_n$ des $\hRun$-générateurs de $N$, et pour tout $1\leq \ell\leq d$, 
soit $A_\ell=(m_{ij}^\ell)_{1\leq i,j\leq n}$ une matrice $n\times n$ à coefficients dans $\hRun$ telle que 
pour tout $1\leq j\leq n$, on ait
\begin{equation}\label{higgs1-higgs11aa}
\theta_\ell(x_j)=p^{\beta}\sum_{i=1}^nm_{ij}^\ell x_i.
\end{equation}
On notera que les matrices $A_\ell$ ne commutent pas deux à deux en général
(mais elles commutent si $N$ est libre sur $\hRun$ de base $x_1,\dots,x_n$). 
Pour tout $\un=(n_1,\dots,n_d)\in \mN^d$, on pose 
\begin{equation}
\uA^{\un}=A_1^{n_1}\cdot A_2^{n_2}\cdots A_d^{n_d} \in \Mat_{n}(\hRun). 
\end{equation}
La série 
\begin{equation}\label{higgs1-higgs11ab}
E=\sum_{\un\in \mN^d}\frac{p^{\beta|\un|}}{\un!}\uA^{\un}\otimes\uy^{\un}
\end{equation}
définit alors une matrice $n\times n$ à coefficients dans $\hfS^{(r)}$. 
Pour tous $(a_1,\dots,a_n)\in (\hfS^{(r)})^n$, on a 
\[
\exp_r(\theta)(\sum_{i=1}^nx_i\otimes a_i)=
\sum_{i=1}^nx_i\otimes b_i\in N\otimes_{\hRun}\hfS^{(r)},
\] 
où  $(b_1,\dots,b_n)=E\cdot (a_1,\dots,a_n)$. 

Comme la matrice $E-\id$ est à coefficients dans $p^{\beta-r}\hfS^{(r)}$,
le déterminant de $E$ est inversible dans $\hfS^{(r)}$. Si $E^{-1}=(f_{ij})_{1\leq i,j\leq n}$ 
est la matrice inverse de $E$ dans $\GL_{n}(\hfS^{(r)})$, alors pour tout $1\leq i\leq d$, on a 
\begin{equation}
\exp_r(\theta)(\sum_{j=1}^dx_j\otimes f_{ji})=x_i\otimes 1.
\end{equation}
Donc $\exp_r(\theta)$ est surjectif. 

Soit $x\in N\otimes_{\hRun}\hfS^{(r)}$ tel que $\exp_r(\theta)(x)=0$. Comme $E\equiv \id \mod(p^{\beta-r})$,
il existe $y\in N\otimes_{\hRun}\hfS^{(r)}$
tel que $x=p^{\beta-r}y$. Comme $N\otimes_{\hRun}\hfS^{(r)}$ est $\co_C$-plat \eqref{higgs1-higgs10},
on a $\exp_r(\theta)(y)=0$. On en déduit que $x\in \cap_{n\geq 0}p^{n(\beta-r)}(N\otimes_{\hRun}\hfS^{(r)})$
et par suite que $x=0$ puisque $N\otimes_{\hRun}\hfS^{(r)}$ est séparé pour la topologie $p$-adique (\cite{egr1} 1.10.2). 
Donc $\exp_r(\theta)$ est bijectif. 

Il résulte aussitôt de la définition  \eqref{higgs1-higgs10b} que le diagramme 
\begin{equation}
\xymatrix{
{N \otimes_{\hRun} \hfS^{(r)}}\ar[rr]^-(0.5){\exp_r(\theta)}
\ar[d]_{\nabla^{(r)}}&&{N\otimes_{\hRun} \hfS^{(r)}}\ar[d]^{p^r\id\otimes d_{\hfS^{(r)}}}\\
{\xi^{-1}\tOmega^1_{R/\co_K}\otimes_RN \otimes_{\hRun} \hfS^{(r)}}\ar[rr]^-(0.5){\exp_r(\theta)\otimes \id}&&
{\xi^{-1}\tOmega^1_{R/\co_K}\otimes_RN\otimes_{\hRun} \hfS^{(r)}}}
\end{equation}
où  $\nabla^{(r)}=\theta\otimes \id+p^r \id\otimes d_{\hfS^{(r)}}$ \eqref{higgs1-not23b}, est commutatif modulo $p^n$ pour tout $n\geq 1$. Il est donc commutatif, autrement dit, $\exp_r(\theta)$ est un morphisme de modules 
à $p^r$-connexions $p$-adiques relativement à l'extension $\hfS^{(r)}/\hRun$.

D'après \ref{higgs1-higgs75} et \eqref{higgs1-higgs99g},  pour tout $g\in \Delta_\infty$, on a 
\begin{eqnarray}
\varphi(g)&=&\exp(\sum_{i=1}^d\xi^{-1}\theta_i\otimes \chi_{i}(g)),\label{higgs1-higgs11c}\\
\varphi_0(g)&=&\exp(-\sum_{i=1}^d\xi^{-1}\frac{\partial}{\partial y_i}\otimes \chi_{i}(g)).\label{higgs1-higgs11d}
\end{eqnarray}
D'autre part, considérant $\exp(-\xi^{-1}\frac{\partial}{\partial y_i}\otimes \chi_{i}(g))$ comme un $\hRun$-automorphisme  
de $\hfS^{(r)}$, on a 
\begin{eqnarray}
\ \ \ \lefteqn{(\id_N\otimes \exp(-\xi^{-1}\frac{\partial}{\partial y_i}\otimes \chi_{i}(g)))\circ \exp_r(\theta)}\label{higgs1-higgs11e}\\
&=&(\exp(-\xi^{-1}\theta_i\otimes \chi_{i}(g))\otimes\id_{\hfS^{(r)}})\circ \exp_r(\theta)\circ 
(\id_N\otimes \exp(-\xi^{-1}\frac{\partial}{\partial y_i}\otimes \chi_{i}(g))).\nonumber
\end{eqnarray}
En effet, il suffit de vérifier cette équation modulo $p^n$ pour tout $n\geq 1$,
ce qui résulte formellement de la définition \eqref{higgs1-higgs10b} (on utilise en premier lieu que 
$\exp(-\xi^{-1}\frac{\partial}{\partial y_i}\otimes \chi_{i}(g))$ est un homomorphisme de la $\hRun$-algèbre $\fS^{(r)}$).
Comme les endomorphismes $\theta_i$ commutent deux à deux, on en déduit que 
\begin{eqnarray}
\lefteqn{(\exp(\sum_{i=1}^d\xi^{-1}\theta_i\otimes \chi_{i}(g))
\otimes \exp(-\sum_{i=1}^d\xi^{-1}\frac{\partial}{\partial y_i}\otimes \chi_{i}(g)))\circ \exp_r(\theta)}\label{higgs1-higgs11f}\\
&=& \exp_r(\theta) \circ (\id_N\otimes \exp(-\sum_{i=1}^d\xi^{-1}\frac{\partial}{\partial y_i}\otimes \chi_{i}(g))).\nonumber
\end{eqnarray}
Par suite, le morphisme $\exp_r(\theta)$ est $\Delta_\infty$-équivariant.

\begin{cor}\label{higgs1-higgs12}
Sous les hypothèses de \eqref{higgs1-higgs11}, on a
un $\hcC^{(r)}$-isomorphisme fonctoriel et $\Delta$-équivariant de $\hoR$-modules de Higgs 
à coefficients dans $\xi^{-1}\tOmega^1_{R/\co_K}$,
\begin{equation}\label{higgs1-higgs12a}
N\otimes_{\hRun} \hcC^{(r)}\stackrel{\sim}{\rightarrow} N\otimes_{\hRun} \hcC^{(r)},
\end{equation}
où $\hcC^{(r)}$ est muni de l'action canonique de $\Delta$ et du $\hoR$-champ de Higgs $p^rd_{\hcC^{(r)}}$ \eqref{higgs1-dolbeault1h},
le module $N$ de la source est muni de l'action triviale de $\Delta_{\infty}$ 
et du $\hRun$-champ de Higgs $\theta$, et le module $N$ du but est muni de l'action $\varphi$
de $\Delta_{\infty}$ et du $\hRun$-champ de Higgs nul. 
Si, de plus, la déformation $(\tX,\cM_\tX)$  est définie par la carte $(P,\gamma)$ \eqref{higgs1-ext24ac}, l'isomorphisme est canonique. 
\end{cor}

En effet, d'après \ref{higgs1-higgs11}, $\exp_r(\theta)$ induit 
un $\hcS^{(r)}$-isomorphisme fonctoriel et $\Delta$-équivariant de $\hoR$-modules de Higgs 
à coefficients dans $\xi^{-1}\tOmega^1_{R/\co_K}$ \eqref{higgs1-chb1},
\begin{equation}\label{higgs1-higgs12b}
N\otimes_{\hRun} \hcS^{(r)}\stackrel{\sim}{\rightarrow} N\otimes_{\hRun} \hcS^{(r)},
\end{equation}
où $\hcS^{(r)}$ est muni de l'action de $\Delta$ induite par $\varphi_0$ \eqref{higgs1-cohbis5} et 
du $\hoR$-champ de Higgs $p^rd_{\hcS^{(r)}}$ \eqref{higgs1-chb1e},
le module $N$ de la source est muni de l'action triviale de $\Delta_{\infty}$ 
et du $\hRun$-champ de Higgs $\theta$, et le module $N$ du but est muni de l'action $\varphi$
de $\Delta_{\infty}$ et du $\hRun$-champ de Higgs nul. 
La proposition s'ensuit compte tenu de \ref{higgs1-tor25} et \ref{higgs1-dolb43}.

\begin{rema}\label{higgs1-higgs120}
Sous les hypothèses de \ref{higgs1-higgs11}, en vertu de \ref{higgs1-not21}, $\exp_r(\theta)$ induit un 
$\hcC^{(r)}$-isomorphisme fonctoriel et $\Delta$-équivariant de modules à $p^r$-connexions $p$-adiques intégrables relativement 
à l'extension $\hcC^{(r)}/\hoR$, 
\begin{equation}\label{higgs1-higgs120a}
N\hotimes_{\hRun} \hcC^{(r)}\stackrel{\sim}{\rightarrow} N\hotimes_{\hRun} \hcC^{(r)},
\end{equation}
où $\hcC^{(r)}$ est muni de l'action canonique de $\Delta$ et de la $p^r$-connexion $p$-adique 
$p^rd_{\hcC^{(r)}}$ \eqref{higgs1-dolbeault1h}, le module $N$ de la source est muni de l'action triviale de $\Delta_{\infty}$ 
et du $\hRun$-champ de Higgs $\theta$, et le module $N$ du but est muni de l'action $\varphi$
de $\Delta_{\infty}$ et du $\hRun$-champ de Higgs nul (cf. \ref{higgs1-not23}).  
\end{rema}

\begin{cor}\label{higgs1-higgs121}
Sous les hypothèses de \eqref{higgs1-higgs11}, on a 
un $\cC^\dagger$-isomorphisme fonctoriel et $\Delta$-équivariant de $\hoR$-modules de Higgs 
à coefficients dans $\xi^{-1}\tOmega^1_{R/\co_K}$
\begin{equation}\label{higgs1-higgs121a}
N\otimes_{\hRun} \cC^\dagger\stackrel{\sim}{\rightarrow} N\otimes_{\hRun} \cC^\dagger,
\end{equation}
où $\cC^\dagger$ est muni de l'action canonique de $\Delta$ et du $\hoR$-champ de Higgs $d_{\cC^\dagger}$ 
\eqref{higgs1-dolbeault1l},
le module $N$ de la source est muni de l'action triviale de $\Delta_{\infty}$ 
et du $\hRun$-champ de Higgs $\theta$, et le module $N$ du but est muni de l'action $\varphi$
de $\Delta_{\infty}$ et du $\hRun$-champ de Higgs nul \eqref{higgs1-not63}. 
Si, de plus, la déformation $(\tX,\cM_\tX)$  est définie par la carte $(P,\gamma)$ \eqref{higgs1-ext24ac}, l'isomorphisme est canonique. 
\end{cor}

Cela résulte de \ref{higgs1-higgs12}.

\begin{cor}\label{higgs1-higgs13}
Soient $N$ un $\hRun$-module projectif de type fini,
$\theta$ un $\hRun$-champ de Higgs quasi-petit sur $N$ à coefficients dans $\xi^{-1}\tOmega^1_{R/\co_K}$,
$\varphi$ la $\hRun$-représentation quasi-petite de $\Delta_\infty$ sur $N$ associée à $\theta$ par le foncteur \eqref{higgs1-higgs6cc}. 
On a alors un isomorphisme fonctoriel de $\hRun$-modules de Higgs à coefficients dans $\xi^{-1}\tOmega^1_{R/\co_K}$
\begin{equation}\label{higgs1-higgs13a}
\mH((N,\varphi)\otimes_{\hRun}\hoR)\stackrel{\sim}{\rightarrow}(N,\theta),
\end{equation}
où $\mH$ est le foncteur \eqref{higgs1-dolb2c}, et un isomorphisme fonctoriel de $\hoR$-représentations de $\Delta$
\begin{equation}\label{higgs1-higgs13b}
\mV(N,\theta)\stackrel{\sim}{\rightarrow} (N,\varphi)\otimes_{\hRun}\hoR,
\end{equation}
où $\mV$ est le foncteur \eqref{higgs1-dolb3c}. 
Si, de plus, la déformation $(\tX,\cM_\tX)$  est définie par la carte $(P,\gamma)$ \eqref{higgs1-ext24ac}, les
isomorphismes sont canoniques. 
\end{cor}

Cela résulte de \ref{higgs1-higgs121} et du fait que
$\ker(d_{\cC^\dagger})=\hoR$ et $(\cC^\dagger)^{\Delta}=\hRun$ \eqref{higgs1-dolb47}.

\begin{cor}\label{higgs1-higgs131}
Tout petit $\hRun$-module de Higgs à coefficients dans $\xi^{-1}\tOmega^1_{R/\co_K}$ 
est soluble, et son image par $\mV$ est une petite $\hoR$-représentation de $\Delta$. 
\end{cor}

Cela résulte de \ref{higgs1-higgs121} et \ref{higgs1-higgs13}.

\vspace{2mm}

Cet énoncé sera renforcé dans \ref{higgs1-desc51}.

\subsection{}\label{higgs1-drt81}
Soient $(N,\theta)$ un petit $\hRun[\frac 1 p]$-module de Higgs à coefficients dans $\xi^{-1}\tOmega^1_{R/\co_K}$ \eqref{higgs1-drt7},
$\beta$ un nombre rationnel $>\frac{1}{p-1}$,  $N^\circ$  un sous-$\hRun$-module de type fini 
de $N$, qui l'engendre sur $\hRun[\frac 1 p]$, tel que l'on ait 
\begin{equation}\label{higgs1-drt81a}
\theta(N^\circ)\subset p^{\beta}\xi^{-1}N^\circ \otimes_R\tOmega^1_{R/\co_K},
\end{equation}
$\theta^\circ$ le $\hRun$-champ de Higgs sur $N^\circ$ à coefficients dans $\xi^{-1}\tOmega^1_{R/\co_K}$
induit par $\theta$. Alors $(N^\circ,\theta^\circ)$ est un $\hRun$-module de Higgs quasi-petit \eqref{higgs1-higgs4}. 
On désigne par $\varphi^\circ$ la $\hRun$-représentation quasi-petite de $\Delta_\infty$ sur $N^\circ$ 
associée à  $\theta^\circ$ par le foncteur \eqref{higgs1-higgs6cc} et par $\varphi$ la petite $\hRun[\frac 1 p]$-représentation 
de $\Delta_\infty$ sur $N$ déduite de $\varphi^\circ$. 
Montrons que $\varphi$ ne dépend pas du choix de $N^\circ$, et que 
la correspondance $(N,\theta)\mapsto (N,\varphi)$ définit un foncteur \eqref{higgs1-drt6}
\begin{equation}\label{higgs1-drt81b}
\bMH^\p(\hRun[\frac 1 p],\xi^{-1}\tOmega^1_{R/\co_K})\rightarrow \bRep_{\hRun[\frac 1 p]}^\p(\Delta_\infty). 
\end{equation}

En effet, soient $u\colon (N,\theta)\rightarrow (N_1,\theta_1)$ un morphisme de 
petits $\hRun[\frac 1 p]$-modules de Higgs à coefficients dans $\xi^{-1}\tOmega^1_{R/\co_K}$, 
$N^\circ_1$ un sous-$\hRun$-module de type fini de $N_1$, qui l'engendre sur $\hRun[\frac 1 p]$, tel que l'on ait 
$u(N^\circ)\subset N^\circ_1$ et  
\begin{equation}\label{higgs1-drt81c}
\theta_1(N^\circ_1)\subset p^{\beta}\xi^{-1}N^\circ_1 \otimes_R\tOmega^1_{R/\co_K}.
\end{equation}
On note $\theta^\circ_1$ le $\hRun$-champ de Higgs sur $N^\circ_1$ à coefficients dans $\xi^{-1}\tOmega^1_{R/\co_K}$
induit par $\theta_1$, $\varphi^\circ_1$ la $\hRun$-représentation de $\Delta_\infty$ sur $N^\circ_1$ 
associée à  $\theta^\circ_1$ par le foncteur \eqref{higgs1-higgs6cc} et $\varphi_1$ la $\hRun[\frac 1 p]$-représentation 
de $\Delta_\infty$ sur $N_1$ induite par $\varphi^\circ_1$.  Le morphisme  
$(N^\circ,\varphi^\circ)\rightarrow (N^\circ_1,\varphi^\circ_1)$ induit par $u$ étant clairement $\Delta_\infty$-équivariant,
il en est alors de même de $u\colon (N,\varphi)\rightarrow (N_1,\varphi_1)$.

Montrons que $\varphi$ ne dépend pas du choix de $N^\circ$.
Soient $\gamma$ un nombre rationnel $>\frac{1}{p-1}$,  $N^\star$ un sous-$\hRun$-module de type fini 
de $N$, qui l'engendre sur $\hRun[\frac 1 p]$, tels que l'on ait 
\begin{equation}\label{higgs1-drt81d}
\theta(N^\star)\subset p^{\gamma}\xi^{-1}N^\star \otimes_R\tOmega^1_{R/\co_K}.
\end{equation}
Remplaçant $\beta$ et $\gamma$ par le plus petit d'entre eux, et $N^\star$ par $N^\circ+N^\star$, 
on peut supposer $\beta=\gamma$ et $N^\circ \subset N^\star$.
Appliquant alors ce qui précède à l'identité de $N$, on déduit que $\varphi$ ne dépend pas du choix de $N^\circ$. 

Montrons que pour tout morphisme $v\colon (N,\theta)\rightarrow (N',\theta')$
de petits $\hRun[\frac 1 p]$-modules de Higgs à coefficients dans $\xi^{-1}\tOmega^1_{R/\co_K}$, 
si $\varphi'$ est la petite $\hRun[\frac 1 p]$-représentation de $\Delta_\infty$ sur $N'$ associée à $\theta'$, 
alors $v\colon (N,\varphi)\rightarrow (N',\varphi')$ est $\Delta_\infty$-équivariant. 
Soient $\beta'$ un nombre rationnel $>\frac{1}{p-1}$,  $N'^\circ$ un sous-$\hRun$-module de type fini 
de $N'$, qui l'engendre sur $\hRun[\frac 1 p]$, tel que l'on ait 
\begin{equation}\label{higgs1-drt81e}
\theta'(N'^\circ)\subset p^{\beta'}\xi^{-1}N'^\circ \otimes_R\tOmega^1_{R/\co_K}.
\end{equation}
Remplaçant $\beta$ et $\beta'$ par le plus petit d'entre eux, 
et $N'^\circ$ par $v(N^\circ)+N'^\circ$, on peut supposer $\beta=\beta'$ et $v(N^\circ) \subset N'^\circ$.
On conclut alors comme plus haut que  $v\colon (N,\varphi)\rightarrow (N',\varphi')$ est $\Delta_\infty$-équivariant. 

\subsection{}\label{higgs1-drt82}
Soient $(N,\varphi)$ une petite $\hRun[\frac 1 p]$-représentation de $\Delta_\infty$ \eqref{higgs1-drt6}, 
$\alpha$ un nombre rationnel $>\frac{2}{p-1}$,  $N^\circ$ un sous-$\hRun$-module de type fini de $N$, 
stable par $\Delta_\infty$, engendré par un nombre fini d'éléments $\Delta_\infty$-invariants modulo 
$p^{\alpha}N^\circ$, et qui engendre $N$ sur $\hRun[\frac 1 p]$. D'après \ref{higgs1-drt65}(ii), la $\hRun$-représentation 
$\varphi^\circ$ de $\Delta_\infty$ sur $N^\circ$, induite par $\varphi$, est quasi-petite. Notons $\theta^\circ$
le $\hRun$-champ de Higgs quasi-petit sur $N^\circ$ à coefficients dans $\xi^{-1}\tOmega^1_{R/\co_K}$ 
associé à $\varphi^\circ$ par le foncteur \eqref{higgs1-higgs6c}, et $\theta$ le petit 
$\hRun[\frac 1 p]$-champ de Higgs sur $N$ induit par $\theta^\circ$. 
Procédant comme dans \ref{higgs1-drt81}, on montre que $(N,\theta)$ ne dépend pas du choix $N^\circ$ et que 
la correspondance $(N,\varphi)\mapsto (N,\theta)$ définit un foncteur \eqref{higgs1-drt7}
\begin{equation}\label{higgs1-drt82b}
\bRep_{\hRun[\frac 1 p]}^\p(\Delta_\infty) \rightarrow \bMH^\p(\hRun[\frac 1 p],\xi^{-1}\tOmega^1_{R/\co_K}). 
\end{equation}
Il résulte aussitôt de \ref{higgs1-higgs6} que les foncteurs \eqref{higgs1-drt81b} et \eqref{higgs1-drt82b} 
sont quasi-inverses l'un de l'autre.

\begin{prop}\label{higgs1-drt8}
Soient $(N,\theta)$ un petit $\hRun[\frac 1 p]$-module de Higgs à coefficients dans $\xi^{-1}\tOmega^1_{R/\co_K}$ \eqref{higgs1-drt7}, 
$\varphi$ la petite $\hRun[\frac 1 p]$-représentation de $\Delta_\infty$ sur $N$ associée à $\theta$ par le foncteur \eqref{higgs1-drt81b}. 
Alors~:
\begin{itemize}
\item[{\rm (i)}] On a un $\cC^\dagger$-isomorphisme fonctoriel 
$\Delta$-équivariant de $\hoR$-modules de Higgs à coefficients dans 
$\xi^{-1}\tOmega^1_{R/\co_K}$
\begin{equation}\label{higgs1-drt8aa}
N \otimes_{\hRun} \cC^\dagger\stackrel{\sim}{\rightarrow} N \otimes_{\hRun} \cC^\dagger,
\end{equation}
où $\cC^\dagger$ est muni de l'action canonique de $\Delta$ et du $\hoR$-champ de Higgs $d_{\cC^\dagger}$ 
\eqref{higgs1-dolbeault1l}, le module $N$ de la source est muni de l'action triviale de $\Delta_{\infty}$ 
et du $\hRun$-champ de Higgs $\theta$, et le module $N$ du but est muni de l'action $\varphi$
de $\Delta_{\infty}$ et du $\hRun$-champ de Higgs nul. 
Si, de plus, la déformation $(\tX,\cM_\tX)$  est définie par la carte $(P,\gamma)$ \eqref{higgs1-ext24ac}, 
l'isomorphisme est canonique.
\item[{\rm (ii)}] Le $\hRun[\frac 1 p]$-module de Higgs $(N,\theta)$ est soluble,
et on a un $\hoR[\frac 1 p]$-isomorphisme $\Delta$-équivariant fonctoriel 
\begin{equation}\label{higgs1-drt8b}
\mV(N)\stackrel{\sim}{\rightarrow} (N,\varphi)\otimes_{\hRun}\hoR,
\end{equation}
où $\mV$ est le foncteur \eqref{higgs1-dolb3c}. 
Si, de plus, la déformation $(\tX,\cM_\tX)$  est définie par la carte $(P,\gamma)$ \eqref{higgs1-ext24ac}, 
l'isomorphisme est canonique.

\item[{\rm (iii)}]  La $\hoR[\frac 1 p]$-représentation $\mV(N)$ de $\Delta$ est petite et de Dolbeault,
et on a un isomorphisme fonctoriel de $\hRun[\frac 1 p]$-modules de Higgs
\begin{equation}\label{higgs1-drt8a}
\mH(\mV(N))\stackrel{\sim}{\rightarrow}(N,\theta),
\end{equation}
où $\mH$ est le foncteur \eqref{higgs1-dolb2c}. 
Si, de plus, la déformation $(\tX,\cM_\tX)$  est définie par la carte $(P,\gamma)$ \eqref{higgs1-ext24ac}, 
l'isomorphisme est canonique.
\end{itemize}
\end{prop}

L'isomorphisme \eqref{higgs1-drt8aa} résulte de \ref{higgs1-higgs121}; le caractère fonctoriel se démontre 
comme dans \ref{higgs1-drt81}. Les autres assertions s'en déduisent compte tenu du fait que
$\ker(d_{\cC^\dagger})=\hoR$ et $(\cC^\dagger)^{\Delta}=\hRun$ \eqref{higgs1-dolb47}.

\begin{prop}[\cite{tsuji3} 5.3.10]\label{higgs1-drt19}
Soient $N$ un $\hRun[\frac 1 p]$-module projectif de type fini, 
$\theta$ un $\hRun[\frac 1 p]$-champ de Higgs sur $N$ à coefficients dans $\xi^{-1}\tOmega^1_{R/\co_K}$,
$M$ un $\hoR[\frac 1 p]$-module, $r$ un nombre rationnel $>0$, 
\begin{equation}\label{higgs1-drt19a}
N\otimes_{\hRun}\hcC^{(r)}\stackrel{\sim}{\rightarrow}M\otimes_{\hoR}\hcC^{(r)}
\end{equation}
un isomorphisme $\hcC^{(r)}$-linéaire de $\hoR$-modules de Higgs à coefficients dans $\xi^{-1}\tOmega^1_{R/\co_K}$,
où $\hcC^{(r)}$ est muni du champ de Higgs $p^rd_{\hcC^{(r)}}$ \eqref{higgs1-dolbeault1h} et $M$ est muni du champ de Higgs nul.
Alors, $(N,\theta)$ est un petit $\hRun[\frac 1 p]$-module de Higgs à coefficients 
dans $\xi^{-1}\tOmega^1_{R/\co_K}$ \eqref{higgs1-drt7}.
\end{prop}

En effet, compte tenu de \ref{higgs1-tor25} et \ref{higgs1-dolb43}, l'isomorphisme \eqref{higgs1-drt19a} 
induit un isomorphisme $\hcS^{(r)}$-linéaire de $\hoR$-modules de Higgs à coefficients dans $\xi^{-1}\tOmega^1_{R/\co_K}$,
\begin{equation}\label{higgs1-drt19b}
N\otimes_{\hRun}\hcS^{(r)}\stackrel{\sim}{\rightarrow}M\otimes_{\hoR}\hcS^{(r)},
\end{equation}
où $\hcS^{(r)}$ est muni du champ de Higgs $p^rd_{\hcS^{(r)}}$ \eqref{higgs1-chb1e} et $M$ est muni du champ de Higgs nul.
Celui-ci induit un isomorphisme $\hoR$-linéaire 
\begin{equation}\label{higgs1-drt19c}
N\otimes_{\hRun}\hoR\stackrel{\sim}{\rightarrow}M.
\end{equation}

Soient $t_1,\dots,t_d\in P^\gp$ tels que leurs images  
dans $(P^\gp/\mZ\lambda)\otimes_\mZ\mZ_p$ forment une $\mZ_p$-base, de sorte que 
$(d\log(t_i))_{1\leq i\leq d}$ est une $R$-base de $\tOmega^1_{R/\co_K}$ \eqref{higgs1-log-ext45b}.
Pour tout $1\leq i\leq d$, on pose $y_i=\xi^{-1}d\log(t_i)\in \xi^{-1}\tOmega^1_{R/\co_K}\subset \cS$. 
Pour tout $\un=(n_1,\dots,n_d)\in \mN^d$, 
on pose $|\un|=\sum_{i=1}^d n_i$ et $\uy^{\un}=\prod_{i=1}^dy_i^{n_i}\in \cS$. 
Munissant $\hoR[\frac 1 p]$ de la topologie $p$-adique \eqref{higgs1-not54}, 
la $\hoR[\frac 1 p]$-algèbre $\hcS^{(r)}[\frac 1 p]$ s'identifie canoniquement à
\begin{equation}\label{higgs1-drt19d}
\{\sum_{\un\in \mN^d} a_\un\uy^\un\in \hoR[\frac 1 p][[y_1,\dots,y_d]] \ | \ \underset{|\un|\rightarrow +\infty}{\lim}\ 
p^{-r|\un|}a_\un=0\}.
\end{equation}
On notera qu'avec les conventions de \eqref{higgs1-chb1e}, on a
$d_{\hcS^{(r)}}(p^ry_i)=1\otimes y_i$ pour tout $1\leq i\leq d$.

Le $\hoR[\frac 1 p][[y_1,\dots,y_d]]$-module
$N\otimes_{\hRun[\frac 1 p]}\hoR[\frac 1 p][[y_1,\dots,y_d]]$ est complet et séparé pour la topologie $(y_1,\dots,y_d)$-adique. 
En effet, il est complet en vertu de (\cite{ac} chap. III §2.12 cor.~1 de prop.~16) et il est séparé en tant que facteur direct 
d'un module libre de type fini. Par suite, toute série formelle 
$\sum_{\un\in \mN^d} x_\un\otimes \uy^\un$, où $x_\un\in N\otimes_{\hRun}\hoR$, converge dans 
$N\otimes_{\hRun[\frac 1 p]}\hoR[\frac 1 p][[y_1,\dots,y_d]]$. Inversement, tout élément $x$
de $N\otimes_{\hRun[\frac 1 p]}\hoR[\frac 1 p][[y_1,\dots,y_d]]$ peut s'écrire uniquement sous la forme 
\begin{equation}\label{higgs1-drt19f}
x=\sum_{\un\in \mN^d} x_\un\otimes \uy^\un,
\end{equation}
où $x_\un\in N\otimes_{\hRun}\hoR$. 
Munissant le $\hoR[\frac 1 p]$-module $N\otimes_{\hRun}\hoR$ de la topologie $p$-adique \eqref{higgs1-not54},
le $\hcS^{(r)}[\frac 1 p]$-module $N\otimes_{\hRun}\hcS^{(r)}$ s'identifie au $\hoR[\frac 1 p]$-module 
\begin{equation}\label{higgs1-drt19g}
\{\sum_{\un\in \mN^d} x_\un\otimes \uy^\un \ | \ 
x_\un\in N\otimes_{\hRun}\hoR, \ \underset{|\un|\rightarrow +\infty}{\lim}\ p^{-r|\un|} x_\un=0\}.
\end{equation}
En effet, ce dernier contient clairement $N\otimes_{\hRun}\hcS^{(r)}$. Pour prouver l'égalité, on peut se réduire
au cas où $N$ est libre de type fini sur $\hRun[\frac 1 p]$, auquel cas l'assertion est évidente \eqref{higgs1-drt19d}.

Soient $x\in N$, $z\in M$ l'image de $x\otimes 1$ par l'isomorphisme \eqref{higgs1-drt19c},
\begin{equation}\label{higgs1-drt19h}
\sum_{\un\in \mN^d} x_\un\otimes \uy^\un \in N\otimes_{\hRun}\hcS^{(r)} 
\end{equation}
l'image de $z\otimes 1\in M\otimes_\hoR\hcS^{(r)}$ par l'isomorphisme inverse de \eqref{higgs1-drt19b},
où $x_\un\in N\otimes_{\hRun}\hoR$. On a clairement
$x_{\underline{0}}=x$. \'Ecrivons 
\begin{equation}\label{higgs1-drt19i}
\theta=\sum_{i=1}^d \theta_i\otimes y_i,
\end{equation}
où les $\theta_i$ sont des $\hRun[\frac 1 p]$-endomorphismes de $N$ qui commutent deux à deux. 
Comme la section \eqref{higgs1-drt19h} est annulée par $\theta\otimes 1+1\otimes p^rd_{\hcS^{(r)}}$, 
pour tout $\un=(n_1,\dots,n_d)\in \mN^d$ et tout $1\leq i\leq d$, on a 
\begin{equation}
\theta_i(x_\un)+(n_i+1)x_{\un+\underline{1}_i}=0,
\end{equation}
où $\underline{1}_i$ est l'élément de $\mN^d$ dont toutes 
les composantes sont nulles sauf la $i$-ième qui vaut $1$. On en déduit que 
\begin{equation}
x_\un=(-1)^{|\un|}(\prod_{1\leq i\leq d}\frac{1}{n_i!}\theta_i^{n_i})(x).
\end{equation}
Par conséquent, $p^{-r|\un|}(\prod_{1\leq i\leq d}\frac{1}{n_i!}\theta_i^{n_i})(x)$ tend vers $0$ 
dans $N\otimes_\hRun\hoR$ quand $|\un|$ tend l'infini \eqref{higgs1-drt19g}. 

La topologie $p$-adique de $N$ est induite par la topologie $p$-adique de $N\otimes_{\hRun}\hoR$. 
En effet, comme $N$ est projectif de type fini sur $\hRun[\frac 1 p]$, on peut se réduire au cas où $N$ est libre de type fini,
et même au cas où $N=\hRun[\frac 1 p]$.  
Pour tout $n\geq 1$, comme l'homomorphisme canonique $\hRun/p^n\hRun\rightarrow \hoR/p^n\hoR$ est injectif 
(cf. la preuve de \ref{higgs1-pur8}), on a $p^n\hoR\cap \hRun[\frac 1 p]=p^n\hRun$, d'où l'assertion. 

On déduit de ce qui précède que pour tout $x\in N$, $p^{-r|\un|}(\prod_{1\leq i\leq d}\frac{1}{n_i!}\theta_i^{n_i})(x)$ tend vers $0$ 
dans $N$ quand $|\un|$ tend l'infini. 

Soient $N_0$ un sous-$\hRun$-module de type fini de $N$ qui l'engendre sur $\hRun[\frac 1 p]$, 
$\varepsilon$ un nombre rationnel tel que $\frac{1}{p-1}<\varepsilon<r+\frac{1}{p-1}$. 
Comme la suite $p^{(r-\varepsilon)n}n!$ tend vers $0$ dans $\co_C$ quand $n$ tend vers l'infini, pour tout $x\in N$, 
$p^{-\varepsilon |\un|}(\prod_{1\leq i\leq d}\theta_i^{n_i})(x)$ tend vers $0$ dans $N$ quand $|\un|$ tend l'infini.
On peut donc considérer le sous-$\hRun$-module 
\begin{equation}
N^\circ=\sum_{\un\in \mN^d}p^{-\varepsilon |\un|}(\prod_{1\leq i\leq d}\theta_i^{n_i})(N_0)
\end{equation}
de $N$. Il est de type fini sur $\hRun$ et il engendre $N$ sur $\hRun[\frac 1 p]$. Comme on a 
\begin{equation}
\theta(N^\circ)\subset p^\varepsilon \xi^{-1} N^\circ\otimes_R\tOmega^1_{R/\co_K},
\end{equation} 
$(N,\theta)$ est un petit $\hRun[\frac 1 p]$-module de Higgs à coefficients dans $\xi^{-1}\tOmega^1_{R/\co_K}$.

\begin{cor}\label{higgs1-drt20}
Pour qu'un $\hRun[\frac 1 p]$-module de Higgs à coefficients dans $\xi^{-1}\tOmega^1_{R/\co_K}$ 
soit soluble \eqref{higgs1-drt3}, il faut et il suffit qu'il soit petit \eqref{higgs1-drt7}.  
\end{cor}

Cela résulte de \ref{higgs1-drt8}(ii) et \ref{higgs1-drt19}.

\begin{cor}\label{higgs1-drt21}
Toute $\hoR[\frac 1 p]$-représentation de Dolbeault de $\Delta$ \eqref{higgs1-drt2} est petite \eqref{higgs1-drt6}. 
\end{cor}

Cela résulte de \ref{higgs1-drt16}, \ref{higgs1-drt20} et \ref{higgs1-drt8}(iii)

\vspace{2mm}

Cet énoncé sera conditionnellement renforcé dans \ref{higgs1-dsct4}.

\section{Descente des petites représentations et applications}\label{higgs1-descente}

Les hypothèses et notations du §\ref{higgs1-dolbeault} sont en vigueur dans cette section.

\begin{prop}\label{higgs1-desc1}
Soient $a$ un élément non nul de $\co_\oK$, $\alpha$ un nombre rationnel $>\frac{1}{p-1}$,   
$M_1$, $M_2$ deux $(R_1/aR_1)$-représentations $\alpha$-petites de $\Delta_{p^\infty}$ \eqref{higgs1-higgs2},  
\begin{equation}\label{higgs1-desc1a}
\ou\colon M_1\otimes_{R_1}\oR\rightarrow M_2\otimes_{R_1}\oR
\end{equation}
un morphisme  $\oR$-linéaire et $\Delta$-équivariant.  
Supposons que $v(a)>\frac{1}{p-1}+\alpha$ et posons $b=ap^{-\alpha-\frac{1}{p-1}}$.
Alors il existe un et un unique morphisme  $R_1$-linéaire et $\Delta_{p^\infty}$-équivariant 
\begin{equation}\label{higgs1-desc1b}
u\colon M_1/bM_1\rightarrow M_2/bM_2
\end{equation}
tel que $u\otimes_{R_1}\oR\equiv \ou \mod bM_2\otimes_{R_1}\oR$.
\end{prop}

Notons $M$ le $R_1$-$\Delta_{p^\infty}$-module discret $\Hom_{R_1}(M_1,M_2)$,
qui est donc une $(R_1/aR_1)$-représentation $\alpha$-petite de $\Delta_{p^\infty}$. 
Comme $R_1$ est normal d'après \ref{higgs1-gal1}(ii), on a $p^\alpha b R_1=(p^\alpha b \oR)\cap R_1$. 
Par suite, le morphisme canonique 
\begin{equation}\label{higgs1-desc1d}
\rH^0(\Delta_{p^\infty},M/p^\alpha bM)\rightarrow \rH^0(\Delta,(M/p^\alpha bM)\otimes_{R_1}\oR)
\end{equation}
est injectif et son conoyau est annulé par $p^{\frac{1}{p-1}}\fm_\oK$ en vertu de \ref{higgs1-cg15}. 
Il existe alors un morphisme  $R_1$-linéaire et $\Delta_{p^\infty}$-équivariant 
\begin{equation}\label{higgs1-desc1e}
v\colon M_1/p^\alpha bM_1\rightarrow M_2/p^\alpha bM_2
\end{equation}
tel que $v\otimes_{R_1}\oR\equiv p^\alpha \ou \mod p^\alpha bM_2\otimes_{R_1}\oR$.
Comme $p^\alpha R_1=(p^\alpha \oR)\cap R_1$ et que $R_1$ est $\co_\oK$-plat, 
représentant $v$ par une matrice à coefficients dans $R_1/p^\alpha b R_1$, 
on voit qu'il existe un unique morphisme  $R_1$-linéaire et $\Delta_{p^\infty}$-équivariant 
\begin{equation}\label{higgs1-desc1f}
u\colon M_1/bM_1\rightarrow M_2/ bM_2
\end{equation}
tel que $v=p^\alpha u$. Utilisant de nouveau que $\oR$ est $\co_\oK$-plat, on en déduit que 
$u\otimes_{R_1}\oR\equiv \ou \mod bM_2\otimes_{R_1}\oR$.
L'unicité de $u$ résulte du fait que l'homomorphisme canonique $R_1/bR_1\rightarrow \oR/b\oR$
est injectif (cf. la preuve de \ref{higgs1-pur8}).

\begin{prop}\label{higgs1-desc2}
Soient $\alpha$ un nombre rationnel $>\frac{1}{p-1}$,
$a$ un élément non nul de $\co_\oK$ tel que $v(a)>\alpha$, 
$M$ un $(\oR/a\oR)$-module libre de rang $r\geq 1$, muni de la topologie discrète 
et d'une action $\oR$-semi-linéaire et continue de $\Delta$ telle que $M$ admette une base $e_1,\dots,e_r$
formée d'éléments $\Delta$-invariants modulo $p^{2\alpha} M$. 
Alors il existe un $R_1$-$\Delta_{p^\infty}$-module discret $N$ dont le module sous-jacent est
libre de rang $r$ sur $R_1/ap^{-\alpha}R_1$, 
ayant une base $f_1,\dots,f_r$ formée d'éléments $\Delta_{p^\infty}$-invariants modulo $p^{\alpha} N$, 
et un isomorphisme $\oR$-linéaire et $\Delta$-équivariant
\begin{equation}\label{higgs1-desc2a}
N\otimes_{R_1}\oR\stackrel{\sim}{\rightarrow}M/ap^{-\alpha}M
\end{equation}
qui transforme $f_j\otimes 1 \mod p^{\alpha}N\otimes_{R_1}\oR$ en $e_j\mod p^{\alpha}M$ 
pour tout $1\leq j\leq r$.  
\end{prop}

La proposition est évidente si $v(a)\leq 3\alpha$, auquel cas on prendra $N=(R_1/ap^{-\alpha}R_1)^r$ 
muni de la représentation triviale de $\Delta_{p^\infty}$ et de la base canonique. 
Supposons donc $v(a)> 3\alpha$. Soit $n$ un entier $\geq 1$ tel que 
\begin{equation}\label{higgs1-desc2b}
\varepsilon=\frac{v(a)-3\alpha}{n}<\frac{1}{3}(\alpha-\frac{1}{p-1}).
\end{equation}
Montrons par une récurrence finie que pour tout $0\leq i\leq n$, 
la proposition vaut pour la $\oR$-représentation $M/p^{3\alpha+i\varepsilon}M$, autrement dit, 
qu'il existe un $R_1$-$\Delta_{p^\infty}$-module discret  $N_i$ dont le module sous-jacent 
est libre de rang $r$ sur $R_1/p^{2\alpha+i\varepsilon}R_1$, ayant une base $f_1^{(i)},\dots,f_r^{(i)}$
formée d'éléments $\Delta_{p^\infty}$-invariants modulo $p^{\alpha} N_i$ et 
un isomorphisme $\oR$-linéaire et $\Delta$-équivariant
\begin{equation}\label{higgs1-desc2c}
N_i\otimes_{R_1}\oR\stackrel{\sim}{\rightarrow}M/p^{2\alpha+i\varepsilon}M
\end{equation}
qui transforme $f_j^{(i)}\otimes 1 \mod  p^{\alpha} N_i\otimes_{R_1}\oR$ 
en $e_j\mod  p^{\alpha}M$ pour tout $1\leq j\leq r$. 
La représentation $N=N_n$ répondra alors à la question puisque $2\alpha+n\varepsilon=v(a)-\alpha$. 

On prend $N_0=(R_1/p^{2\alpha} R_1)^r$ muni de la représentation triviale de $\Delta_{p^\infty}$ et de la base canonique.
Supposons $N_i$ construit avec $0\leq i<n$ et construisons $N_{i+1}$. 
D'après \ref{higgs1-cnab5}, l'obstruction à relever $N_i$ en un $(R_1/p^{3\alpha+i\varepsilon}R_1)$-$\Delta_{p^\infty}$-module 
discret dont le module sous-jacent est libre de type fini sur $R_1/p^{3\alpha+i\varepsilon}R_1$ est un élément $\fo$ de 
$\rH^2(\Delta_{p^\infty},\Mat_r(R_1/p^{\alpha} R_1))$. D'autre part, le morphisme 
\begin{equation}\label{higgs1-desc2d}
\rH^2(\Delta_{p^\infty},\Mat_r(R_1/p^{\alpha} R_1))\rightarrow \rH^2(\Delta,\Mat_r(\oR/p^{\alpha} \oR))
\end{equation}
est presque-injectif en vertu de \ref{higgs1-cg17}.
L'image de $\fo$ dans $\rH^2(\Delta,\Mat_r(\oR/p^{\alpha} \oR))$ est nulle puisque 
la représentation $M/p^{2\alpha+i\varepsilon}M$
se relève en $M/p^{3\alpha+i\varepsilon}M$; donc $p^\varepsilon \fo=0$. 
Par suite, en vertu de \eqref{higgs1-cnab5d}, $N_i/p^{2\alpha+(i-1)\varepsilon}N_i$
se relève en un $(R_1/p^{3\alpha+(i-1)\varepsilon}R_1)$-$\Delta_{p^\infty}$-module discret $N'_{i+1}$ 
dont le module sous-jacent est libre de type fini sur $R_1/p^{3\alpha+(i-1)\varepsilon}R_1$. 
D'après \ref{higgs1-cnab4}, le relèvement $N'_{i+1}\otimes_{R_1}\oR$ de $M/p^{2\alpha+(i-1)\varepsilon} M$
se déduit de $M/p^{3\alpha+(i-1)\varepsilon} M$ par torsion par un élément 
$\ofc$ de $\rH^1(\Delta,\Mat_r(\oR/p^{\alpha} \oR))$. 
En vertu de \ref{higgs1-cg17}, le conoyau du morphisme canonique
\begin{equation}\label{higgs1-desc2e}
\rH^1(\Delta_{p^\infty},\Mat_r(R_1/p^{\alpha} R_1))\rightarrow \rH^1(\Delta,\Mat_r(\oR/p^{\alpha} \oR))
\end{equation}
est annulé par $p^{\frac{1}{p-1}+\varepsilon}$. 
Comme $\alpha>\frac{1}{p-1}+3\varepsilon$, $p^{\alpha-2\varepsilon}\ofc$ est l'image d'un élément 
\[
\fc'\in \rH^1(\Delta_{p^\infty},\Mat_r(R_1/p^{\alpha} R_1)).
\] 
Modifiant $N'_{i+1}/p^{2\alpha+(i+1)\varepsilon}N'_{i+1}$ par torsion par $-\fc'$, 
on obtient un $(R_1/p^{2\alpha+(i+1)\varepsilon}R_1)$-$\Delta_{p^\infty}$-module discret $N_{i+1}$ 
dont le module sous-jacent est libre de type fini sur $R_1/p^{2\alpha+(i+1)\varepsilon}R_1$, qui relève 
$N_i/p^{\alpha+(i+1)\varepsilon}N_i$. On notera que $2\alpha+(i-1)\varepsilon>\alpha+(i+1)\varepsilon$. 
Alors $N_{i+1}$ répond à la question d'après la dernière assertion du \ref{higgs1-cnab4}. 

\begin{prop}\label{higgs1-desc3}
Soient $a$ un élément non nul de $\co_\oK$, $\alpha,\beta$ deux nombres rationnels tels que 
$v(a)>\alpha>\beta>\frac{1}{p-1}$, $M$ un $(\oR/a\oR)$-module libre de rang $r\geq 1$, muni de la topologie discrète 
et d'une action $\oR$-semi-linéaire et continue de $\Delta$ telle que $M$ admette une base 
formée d'éléments $\Delta$-invariants modulo $p^{2\alpha} M$. 
Alors il existe un $R_1$-$\Delta_{p^\infty}$-module discret $N$ dont le module sous-jacent est
libre de rang $r$ sur $R_1/ap^{-\alpha}R_1$, 
ayant une base formée d'éléments $\Delta_{p^\infty}$-invariants modulo $p^{2\beta} N$, 
et un isomorphisme $\oR$-linéaire et $\Delta$-équivariant
\begin{equation}\label{higgs1-desc3a}
N\otimes_{R_1}\oR\stackrel{\sim}{\rightarrow}M/ap^{-\alpha}M.
\end{equation}
\end{prop}

La proposition est évidente si $v(a)\leq 3\alpha$, auquel cas on prendra $N=(R_1/ap^{-\alpha}R_1)^r$ 
muni de la représentation triviale de $\Delta_{p^\infty}$. 
Supposons donc $v(a)> 3\alpha$. Soit $n$ un entier $\geq 1$ tel que 
\begin{equation}\label{higgs1-desc3b}
\varepsilon=\frac{v(a)-3\alpha}{n}<\inf(\frac{1}{3}(\alpha-\frac{1}{p-1}),\alpha-\beta).
\end{equation}
Montrons par une récurrence finie que pour tout $0\leq i\leq n$, 
la proposition vaut pour la $\oR$-représentation $M/p^{3\alpha+i\varepsilon}M$, autrement dit, 
qu'il existe un $R_1$-$\Delta_{p^\infty}$-module discret  $N_i$ dont le module sous-jacent 
est libre de rang $r$ sur $R_1/p^{2\alpha+i\varepsilon}R_1$, ayant une base
formée d'éléments $\Delta_{p^\infty}$-invariants modulo $p^{2(\alpha-\varepsilon)} N_i$ et 
un isomorphisme $\oR$-linéaire et $\Delta$-équivariant
\begin{equation}\label{higgs1-desc3c}
N_i\otimes_{R_1}\oR\stackrel{\sim}{\rightarrow}M/p^{2\alpha+i\varepsilon}M.
\end{equation}
La représentation $N=N_n$ répondra alors à la question puisque $2\alpha+n\varepsilon=v(a)-\alpha$
et $\alpha-\varepsilon>\beta$. 

On prend $N_0=(R_1/p^{2\alpha} R_1)^r$ muni de la représentation triviale de $\Delta_{p^\infty}$ et de la base canonique.
Supposons $N_i$ construit avec $0\leq i<n$ et construisons $N_{i+1}$. 
D'après \ref{higgs1-cnab5} (appliqué avec $a^q=p^{3\alpha+i\varepsilon}$, $a^n=p^{2\alpha+i\varepsilon}$
et $a^m=p^{2(\alpha-\varepsilon)}$),
l'obstruction à relever $N_i$ en un $(R_1/p^{3\alpha+i\varepsilon}R_1)$-$\Delta_{p^\infty}$-module discret 
dont le module sous-jacent est libre de type fini sur $R_1/p^{3\alpha+i\varepsilon}R_1$ est un élément $\fo$ de 
$\rH^2(\Delta_{p^\infty},\Mat_r(R_1/p^{\alpha} R_1))$. 
D'autre part, le morphisme 
\begin{equation}\label{higgs1-desc3d}
\rH^2(\Delta_{p^\infty},\Mat_r(R_1/p^{\alpha} R_1))\rightarrow \rH^2(\Delta,\Mat_r(\oR/p^{\alpha} \oR))
\end{equation}
est presque-injectif en vertu de \ref{higgs1-cg17}.
L'image de $\fo$ dans $\rH^2(\Delta,\Mat_r(\oR/p^{\alpha} \oR))$ est nulle puisque 
la représentation $M/p^{2\alpha+i\varepsilon}M$
se relève en $M/p^{3\alpha+i\varepsilon}M$; donc $p^\varepsilon \fo=0$. 
Par suite, d'après \eqref{higgs1-cnab5d}, $N_i/p^{2\alpha+(i-1)\varepsilon}N_i$
se relève en un $(R_1/p^{3\alpha+(i-1)\varepsilon}R_1)$-$\Delta_{p^\infty}$-module discret $N'_{i+1}$ 
dont le module sous-jacent est libre de type fini sur $R_1/p^{3\alpha+(i-1)\varepsilon}R_1$. 
D'après \ref{higgs1-cnab4}, le relèvement $N'_{i+1}\otimes_{R_1}\oR$ de $M/p^{2\alpha+(i-1)\varepsilon} M$ 
se déduit du relèvement $M/p^{3\alpha+(i-1)\varepsilon} M$ par torsion par un élément $\ofc$ de 
$\rH^1(\Delta,\Mat_r(\oR/p^{\alpha} \oR))$. 
En vertu de \ref{higgs1-pur121}, le conoyau du morphisme canonique
\begin{equation}\label{higgs1-desc3e}
\rH^1(\Delta_{p^\infty},\Mat_r(R_{p^\infty}/p^{\alpha} R_{p^\infty}))\rightarrow \rH^1(\Delta,\Mat_r(\oR/p^{\alpha} \oR))
\end{equation}
est annulé par $p^\varepsilon$. Donc $p^{\varepsilon}\ofc$ est l'image d'un élément 
\[
\fc'\in \rH^1(\Delta_{p^\infty},\Mat_r(R_{p^\infty}/p^{\alpha} R_{p^\infty})).
\] 
D'après la preuve de \ref{higgs1-cg3}(i), on a $\fc'=\fc'_1+\fc'_2$, où $p^{\frac{1}{p-1}}\fc'_2=0$ et $\fc'_1$ est l'image d'un élément
\[
\fc_1\in \rH^1(\Delta_{p^\infty},\Mat_r(R_1/p^{\alpha} R_1)).
\] 
Modifiant $N'_{i+1}/p^{3\alpha+(i-2)\varepsilon} N'_{i+1}$ par torsion par $-\fc_1$, 
on obtient un $(R_1/p^{3\alpha+(i-2)\varepsilon}R_1)$-$\Delta_{p^\infty}$-module discret $N''_{i+1}$ 
dont le module sous-jacent est libre de type fini sur $R_1/p^{3\alpha+(i-2)\varepsilon}R_1$, qui relève 
$N_i/p^{2\alpha+(i-2)\varepsilon}N_i$. D'après la dernière assertion de \ref{higgs1-cnab4}, le relèvement 
$N''_{i+1}\otimes_{R_1}\oR$ de $M/p^{2\alpha+(i-2)\varepsilon} M$ se déduit du relèvement 
$M/p^{3\alpha+(i-2)\varepsilon} M$ par torsion par l'image $\ofc_2$ de $\fc'_2$ 
dans $\rH^1(\Delta,\Mat_r(\oR/p^{\alpha} \oR))$.  
Comme $p^{\frac{1}{p-1}}\ofc_2=0$, on en déduit encore par la dernière assertion de \ref{higgs1-cnab4} 
qu'il existe un isomorphisme $\oR$-linéaire et $\Delta$-équivariant
\begin{equation}\label{higgs1-desc3f}
(N''_{i+1}/p^{3\alpha-\frac{1}{p-1}+(i-2)\varepsilon} N''_{i+1})\otimes_{R_1}\oR\stackrel{\sim}{\rightarrow}M/
p^{3\alpha-\frac{1}{p-1}+(i-2)\varepsilon}M.
\end{equation}
Comme $\alpha>\frac{1}{p-1}+3\varepsilon$, la représentation 
$N_{i+1}=N''_{i+1}/p^{2\alpha+(i+1)\varepsilon} N''_{i+1}$ répond alors à la question.

\begin{prop}[Faltings, \cite{faltings3}]\label{higgs1-desc4} 
Le foncteur \eqref{higgs1-higgs2}
\begin{equation}\label{higgs1-desc4a}
\bRep_{\hRun}^\p(\Delta_{\infty})\rightarrow \bRep_{\hoR}^\p(\Delta), \ \ \
M\mapsto M\otimes_{\hRun}\hoR
\end{equation}
est une équivalence de catégories.  
\end{prop}

En effet, comme le foncteur canonique
\begin{equation}
\bRep_{\hRun}^{\p}(\Delta_{p^\infty})\rightarrow \bRep_{\hRun}^{\p}(\Delta_{\infty})\label{higgs1-higgs8a}\\
\end{equation}
est une équivalence de catégories \eqref{higgs1-higgs75}, il suffit de montrer que le foncteur 
\begin{equation}\label{higgs1-desc4e}
\bRep_{\hRun}^\p(\Delta_{p^\infty})\rightarrow \bRep_{\hoR}^\p(\Delta), \ \ \
M\mapsto M\otimes_{\hRun}\hoR
\end{equation}
est une équivalence de catégories. Ce foncteur est pleinement fidèle en vertu de \ref{higgs1-desc1}. 
Montrons qu'il est essentiellement surjectif. Soient $\alpha,\beta$ deux nombres rationnels 
tels que $\alpha>\beta>\frac{1}{p-1}$, $M$ une $\hoR$-représentation $(2\alpha)$-petite de rang $r\geq 1$. 
D'après \ref{higgs1-desc3}, pour tout entier $n>\alpha$, 
il existe un $R_1$-$\Delta_{p^\infty}$-module discret $N_n$ dont le module sous-jacent est
libre de rang $r$ sur $R_1/p^{n-\alpha}R_1$, 
ayant une base formée d'éléments $\Delta_{p^\infty}$-invariants modulo $p^{2\beta} N$, 
et un isomorphisme $\oR$-linéaire et $\Delta_{p^\infty}$-équivariant
\begin{equation}\label{higgs1-desc4b}
N_n\otimes_{R_1}\oR\stackrel{\sim}{\rightarrow}M/p^{n-\alpha}M.
\end{equation}
En vertu de \ref{higgs1-desc1}, pour tous entiers $n\geq m>\alpha$, il existe un unique isomorphisme 
$R_1$-linéaire et $\Delta_{p^\infty}$-équivariant
\begin{equation}\label{higgs1-desc4c}
N_n/p^{m-\alpha}N_n\stackrel{\sim}{\rightarrow}N_m
\end{equation}
compatible avec les isomorphismes \eqref{higgs1-desc4b}. Par suite, les $R_1$-modules $(N_n)_{n>\alpha}$
forment un système projectif, et si $N$ est sa limite projective, on a un isomorphisme 
$\hoR$-linéaire et $\Delta$-équivariant
\begin{equation}\label{higgs1-desc4d}
N\otimes_{\hRun}\hoR\stackrel{\sim}{\rightarrow}M.
\end{equation}

\begin{remas}
(i)\  Bien que parallèles, les preuves des propositions \ref{higgs1-desc2} et \ref{higgs1-desc3} diffèrent significativement 
dans leurs dernières étapes qui consistent à définir $N_{i+1}$ à partir de $N'_{i+1}$. 
On notera qu'en ce qui concerne la descente de la base, la preuve de \ref{higgs1-desc3} 
ne donne pas mieux que celle de \ref{higgs1-desc2}, c'est pourquoi nous avons omis cette partie de l'énoncé \ref{higgs1-desc3}. 
La proposition \ref{higgs1-desc2} est due à Faltings (\cite{faltings3} Lemme 1). 

(ii)\ L'énoncé de descente établi par Faltings dans \cite{faltings3} est en fait légèrement plus faible que \ref{higgs1-desc4}.
\end{remas}

\begin{cor}\label{higgs1-desc5}
Toute petite $\hoR$-représentation de $\Delta$ est de Dolbeault, et son image par 
$\mH$ \eqref{higgs1-dolb2c} est un petit $\hRun$-module de Higgs à coefficients dans $\xi^{-1}\tOmega^1_{R/\co_K}$.  
\end{cor}

Cela résulte de \ref{higgs1-desc4}, \ref{higgs1-higgs7}, \ref{higgs1-higgs121} et \ref{higgs1-higgs13}.

\begin{cor}\label{higgs1-desc51}
Les foncteurs $\mH$ \eqref{higgs1-dolb2c} et $\mV$ \eqref{higgs1-dolb3c} induisent des équivalences de catégories quasi-inverses l'une de l'autre,
entre la catégorie des petites $\hoR$-représentations de $\Delta$ et celle des petits $\hRun$-modules de Higgs  
à coefficients dans $\xi^{-1}\tOmega^1_{R/\co_K}$.
\end{cor}

Cela résulte de \ref{higgs1-dolb7} et \ref{higgs1-higgs131}, \ref{higgs1-desc5}.

\begin{prop}\label{higgs1-dsct4}
Les propositions suivantes sont équivalentes~: 
\begin{itemize}
\item[{\rm (i)}] Pour toute petite $\hoR[\frac 1 p]$-représentation $M$ de $\Delta$ \eqref{higgs1-drt6}, 
il existe  un petit $\hRun[\frac 1 p]$-module de Higgs $N$ à coefficients dans 
$\xi^{-1}\tOmega^1_{R/\co_K}$ \eqref{higgs1-drt7} et un $\hoR[\frac 1 p]$-isomorphisme $\Delta$-équivariant
\begin{equation}\label{higgs1-dsct4a}
M\stackrel{\sim}{\rightarrow} \mV(N).
\end{equation}
\item[{\rm (ii)}] Toute petite $\hoR[\frac 1 p]$-représentation de $\Delta$  est de Dolbeault \eqref{higgs1-drt2}.
\item[{\rm (iii)}] Pour qu'une $\hoR[\frac 1 p]$-représentation de $\Delta$ soit de  de Dolbeault,  
il faut et il suffit qu'elle soit petite.
\item[{\rm (iv)}] Le foncteur \eqref{higgs1-drt6}
\begin{equation}\label{higgs1-dsct4b}
\bRep_{\hRun[\frac 1 p]}^\p(\Delta_{\infty})\rightarrow \bRep_{\hoR[\frac 1 p]}^\p(\Delta), \ \ \
M\mapsto M\otimes_{\hRun}\hoR
\end{equation}
est une équivalence de catégories. 
\end{itemize}
\end{prop}

En effet, l'implication (i)$\Rightarrow$(ii) résulte de \ref{higgs1-drt8}(iii), l'implication (ii)$\Rightarrow$(iii)
est une conséquence de \ref{higgs1-drt21}, et l'implication (iii)$\Rightarrow$(i) résulte de \ref{higgs1-dolb7} et \ref{higgs1-drt20}. 

Montrons ensuite (i)$\Rightarrow$(iv). Pour toute petite $\hoR[\frac 1 p]$-représentation $M$
de $\Delta$, le $\hRun[\frac 1 p]$-module de Higgs $\mH(M)$ est petit d'après (i) et \ref{higgs1-drt8}(iii). Il lui est donc associé 
par le foncteur \eqref{higgs1-drt81b} une petite $\hRun[\frac 1 p]$-représentation $\varphi$ de $\Delta_\infty$ sur $\mH(M)$.  
La correspondance $M\mapsto (\mH(M),\varphi)$ définit un foncteur 
\begin{equation}\label{higgs1-dsct4c}
\bRep_{\hoR[\frac 1 p]}^\p(\Delta) \rightarrow \bRep_{\hRun[\frac 1 p]}^\p(\Delta_{\infty}).
\end{equation}
Montrons que les foncteurs \eqref{higgs1-dsct4b} et \eqref{higgs1-dsct4c} sont quasi-inverses l'un de l'autre. 
Pour toute petite $\hoR[\frac 1 p]$-représentation $M$ de $\Delta$, on a des $\hoR$-isomorphismes 
$\Delta$-équivariants fonctoriels \eqref{higgs1-drt8b}
\begin{equation}\label{higgs1-dsct4d}
\mH(M)\otimes_{\hRun}\hoR\stackrel{\sim}{\rightarrow}\mV(\mH(M))
\stackrel{\sim}{\rightarrow} M.
\end{equation}
Par ailleurs, soient $(N,\varphi)$ une petite $\hRun[\frac 1 p]$-représentation de $\Delta_\infty$, $\theta$ le 
petit $\hRun[\frac 1 p]$-champ de Higgs sur $N$ à coefficients dans $\xi^{-1}\tOmega^1_{R/\co_K}$ associé à $\varphi$
par le foncteur \eqref{higgs1-drt82b}. 
D'après \ref{higgs1-drt8}, on a un isomorphisme fonctoriel de $\hRun$-modules de Higgs 
\begin{equation}\label{higgs1-dsct4e}
\mH((N,\varphi)\otimes_{\hRun}\hoR)\stackrel{\sim}{\rightarrow}(N,\theta).
\end{equation}
Appliquant le foncteur \eqref{higgs1-drt81b}, quasi-inverse du foncteur \eqref{higgs1-drt82b}, on obtient un $\hRun$-isomor\-phisme 
$\Delta_\infty$-équivariant et fonctoriel 
\begin{equation}\label{higgs1-dsct4g}
\mH((N,\varphi)\otimes_{\hRun}\hoR)\stackrel{\sim}{\rightarrow} (N,\varphi).
\end{equation}
Les foncteurs \eqref{higgs1-dsct4b} et \eqref{higgs1-dsct4c} sont donc quasi-inverses l'un de l'autre. 

Montrons enfin (iv)$\Rightarrow$(i). Soit $M$ une petite $\hoR[\frac 1 p]$-représentation de $\Delta$.
D'après (iv), il existe une petite $\hRun[\frac 1 p]$-représentation $(N,\varphi)$ de $\Delta_\infty$ et un 
$\hoR$-isomorphisme $\Delta$-équivariant 
\begin{equation}
M\stackrel{\sim}{\rightarrow}(N,\varphi)\otimes_{\hRun}\hoR.
\end{equation} 
Soit $\theta$ le petit $\hRun[\frac 1 p]$-champ de Higgs sur $N$ à coefficients dans $\xi^{-1}\tOmega^1_{R/\co_K}$ 
associé à $\varphi$ par le foncteur \eqref{higgs1-drt82b}. En vertu de \ref{higgs1-drt8}(ii), comme les foncteurs \eqref{higgs1-drt81b} et 
\eqref{higgs1-drt82b} sont quasi-inverses l'un de l'autre, on a un $\hoR[\frac 1 p]$-isomorphisme  $\Delta$-équivariant 
$M\stackrel{\sim}{\rightarrow} \mV(N,\theta)$, d'où la proposition. 

\begin{cor}\label{higgs1-dsct7}
Supposons les propositions équivalentes de \eqref{higgs1-dsct4} satisfaites.  
Alors pour toute petite $\hoR[\frac 1 p]$-représentation $M$ de $\Delta$, 
la $\hoR[\frac 1 p]$-représentation $\Hom_{\hoR[\frac 1 p]}(M,\hoR[\frac 1 p])$ de $\Delta$ est petite.
\end{cor}

En effet, il existe une petite $\hRun[\frac 1 p]$-représentation $N$ de $\Delta_\infty$ et un $\hoR$-isomorphisme 
$\Delta$-équivariant $M\stackrel{\sim}{\rightarrow}N\otimes_{\hRun}\hoR$. On en déduit un 
$\hoR$-isomorphisme $\Delta$-équivariant 
\begin{equation}
\Hom_{\hoR[\frac 1 p]}(M,\hoR[\frac 1 p])\stackrel{\sim}{\rightarrow}
\Hom_{\hRun[\frac 1 p]}(N,\hRun[\frac 1 p])\otimes_{\hRun}\hoR.
\end{equation}
Soient $\alpha$ un nombre rationnel $>\frac{2}{p-1}$, 
$N^\circ$ un sous-$\hRun$-module de type fini de $N$, stable par $\Delta_\infty$, 
engendré par un nombre fini d'éléments $\Delta_\infty$-invariants modulo 
$p^{\alpha}N^\circ$, et qui engendre $N$ sur $\hRun[\frac 1 p]$. 
En vertu de (\cite{egr1} 10.10.2(iii)), $N^\circ$ est un $\hRun$-module cohérent. 
Il en est donc de même de  $\Hom_{\hRun}(N^\circ,\hRun)$, et le morphisme canonique 
\begin{equation}
\Hom_{\hRun}(N^\circ,\hRun)\otimes_{\hRun}\hRun[\frac 1 p] \rightarrow \Hom_{\hRun[\frac 1 p]}(N,\hRun[\frac 1 p])
\end{equation}
est un isomorphisme. D'autre part, comme $\hRun$ est $\co_\oK$-plat, pour tout nombre rationnel $\beta>0$, 
le morphisme canonique 
\begin{equation}
\Hom_{\hRun}(N^\circ,\hRun)\otimes_{\co_\oK}\co_\oK/p^\beta\co_\oK \rightarrow 
\Hom_{R_1}(N^\circ/p^\beta N^\circ,R_1/p^\beta R_1)
\end{equation}
est injectif. Il s'ensuit que la représentation de $\Delta_\infty$ sur $\Hom_{\hRun}(N^\circ,\hRun)$ est continue pour la topologie $p$-adique, et est quasi-petite, d'où la proposition.

\begin{cor}\label{higgs1-dsct6}
Supposons les propositions équivalentes de \eqref{higgs1-dsct4} satisfaites.  
Alors pour toutes petites $\hoR[\frac 1 p]$-représentations $M$ et $M'$ de $\Delta$, 
et tout morphisme $\hoR[\frac 1 p]$-linéaire, $\Delta$-équivariant et surjectif $u\colon M'\rightarrow M$, 
la $\hoR[\frac 1 p]$-représentation de $\Delta$ sur le noyau de $u$ est petite. 
\end{cor}

Cela résulte de \ref{higgs1-dsct7}.

\begin{lem}\label{higgs1-desc6}
Soient $\alpha,\varepsilon$ deux nombres rationnels tels que $0<\varepsilon <\alpha$, 
$M$ une $\hoR$-représentation $\alpha$-quasi-petite de $\Delta$. Alors~: 
\begin{itemize}
\item[{\rm (i)}]  Le $R_\infty$-module $(M/p^{\alpha} M)^\Sigma$ est presque de type fini. 
\item[{\rm (ii)}] Le $\hoR$-module $p^\varepsilon M$ est engendré par un nombre fini d'éléments qui sont d'une part 
$\Delta$-invariants modulo $p^{\alpha}M$ et d'autre part $\Sigma$-invariants.
\item[{\rm (iii)}] Le morphisme canonique $M^\Sigma\otimes_{\hRi}\hoR\rightarrow M$ est presque surjectif. 
\end{itemize}
\end{lem}

Soient $x_1,\dots,x_d$ des générateurs de $M$ sur $\hoR$ qui sont $\Delta$-invariants modulo 
$p^{\alpha}M$. 

(i) On désigne par $u\colon (\oR/p^{\alpha}\oR)^d\rightarrow M/p^{\alpha}M$ le morphisme $\oR$-linéaire, 
$\Delta$-équivariant et surjectif défini par les $(x_i)_{1\leq i\leq d}$ et par $C$ son noyau. 
Alors $\rH^1(\Sigma,C)$ est presque nul en vertu de \ref{higgs1-pur5}.
D'autre part, le morphisme canonique 
\[
R_\infty/p^{\alpha}R_\infty \rightarrow (\oR/p^{\alpha}\oR)^\Sigma
\]
est un presque-isomorphisme \eqref{higgs1-pur12}. Par suite, le morphisme $R_\infty$-linéaire
\begin{equation}
(R_\infty/p^{\alpha} R_\infty)^d\rightarrow (M/p^{\alpha}M)^\Sigma
\end{equation}
défini par les $(x_i)_{1\leq i\leq d}$ est presque surjectif. 

(ii) Le $\hoR$-module $p^{\alpha}M$ est complet et séparé pour la topologie $p$-adique~; il  
est complet en vertu de (\cite{ac} chap. III §2.12 cor.~1 de prop.~16), et est séparé en tant que sous-module de $M$. 
Donc $\rH^1_\cont(\Sigma,p^{\alpha}M)$ est presque nul en vertu de \ref{higgs1-pur55}. Par suite, 
le morphisme canonique 
\begin{equation}
M^\Sigma\rightarrow (M/p^{\alpha}M)^\Sigma
\end{equation}
est presque surjectif. Les classes des éléments $p^{\varepsilon} x_1,\dots,p^{\varepsilon} x_n$
dans $(M/p^{\alpha}M)^\Sigma$ se relèvent en des éléments $x'_1,\dots,x'_d\in M^\Sigma$. 
Pour tout $1\leq i\leq d$, on a 
\begin{equation}
x'_i\in M^\Sigma\cap (p^\varepsilon M)=(p^\varepsilon M)^\Sigma.
\end{equation}
D'autre part, par  le lemme de Nakayama, les $(x'_i)_{1\leq i\leq d}$ engendrent $p^{\varepsilon}M$ sur $\hoR$. 

(iii) Cela résulte de (ii). 

\begin{lem}\label{higgs1-desc66}
Soient $M$ un $\hoR$-module complet et séparé pour la topologie $p$-adique, muni d'une action 
$\hoR$-semi-linéaire et continue de $\Sigma$, $x_1,\dots,x_d$ des éléments de $M^\Sigma$ 
qui engendrent $M$ sur $\hoR$. Alors le morphisme $\hRi$-linéaire $\hRi^d\rightarrow M^\Sigma$
défini par les $(x_i)_{1\leq i\leq d}$, est presque surjectif.
\end{lem} 

Notons $C$ le noyau du morphisme $\hoR$-linéaire, $\Sigma$-équivariant et surjectif $\hoR^d\rightarrow M$
défini par les $(x_i)_{1\leq i\leq d}$. 
Comme $M$ est séparé, $C$ est fermé dans $\hoR^d$ pour la topologie $p$-adique. 
Il est donc complet et séparé pour 
la topologie induite par la topologie $p$-adique de $\hoR^d$; autrement dit, $C$ est isomorphe à la limite
projective du système projectif des $\hoR$-modules $(C/(C\cap p^n\hoR^d))_{n\in \mN}$. 
Par suite, $\rH^1_\cont(\Sigma,C)$ est presque nul
en vertu de \ref{higgs1-pur55}. Par ailleurs, l'homomorphisme canonique $\hRi\rightarrow \hoR^\Sigma$
est un presque-isomorphisme d'après \ref{higgs1-pur122}, d'où la proposition. 

\begin{lem}\label{higgs1-dsct0}
Soient $\alpha$ un nombre rationnel $>0$, 
$M$ une $\hRi$-représentation $\alpha$-quasi-petite de $\Delta_\infty$ telle que $M$ soit $\mZ_p$-plat.
Alors, le $\hRpi$-module $M^{\Sigma_0}$, muni de l'action induite de $\Delta_{p^\infty}$,
est une $\hRpi$-représentation $\alpha$-quasi-petite de $\Delta_{p^\infty}$, et le morphisme canonique 
\begin{equation}\label{higgs1-dsct0a}
M^{\Sigma_0}\otimes_{\hRpi}\hRi\rightarrow M
\end{equation}
est surjectif.
\end{lem}

D'après \eqref{higgs1-limproj2c} et \eqref{higgs1-limproj2d}, on a  
\begin{equation}\label{higgs1-dsct0b}
0\rightarrow \rR^1\underset{\underset{n}{\longleftarrow}}{\lim}\ (M/p^nM)^{\Sigma_0}\rightarrow 
\rH^1_\cont(\Sigma_0,M)\rightarrow \underset{\underset{n}{\longleftarrow}}{\lim}\ 
\rH^1(\Sigma_0,M/p^nM)\rightarrow 0.
\end{equation} 
D'autre part, comme $\Sigma_0$ est un groupe profini d'ordre premier à $p$,
$\rH^1(\Sigma_0,M/p^nM)=0$ pour tout $n\geq 0$,  et le système projectif
$((M/p^nM)^{\Sigma_0})_{n\geq 0}$ vérifie la condition de Mittag-Leffler. 
On en déduit que $\rH^1_\cont(\Sigma_0,M)=0$. 
De même, on a $\rH^1_\cont(\Sigma_0,p^\alpha M)=0$ puisque $p^\alpha M$ est complet et séparé 
pour la topologie $p$-adique. Par suite le morphisme canonique 
\begin{equation}\label{higgs1-dsct0c}
M^{\Sigma_0}\rightarrow (M/p^\alpha M)^{\Sigma_0}
\end{equation}
est surjectif. 

Soient $x_1,\dots,x_d$ des générateurs de $M$ sur $R_\infty$ qui sont $\Delta_\infty$-invariants modulo $p^\alpha M$.
Leurs classes dans $(M/p^\alpha M)^{\Sigma_0}$ se relèvent en des éléments $x'_1,\dots,x'_d$ de $M^{\Sigma_0}$, 
qui engendrent $M$ sur $\hRi$. 
Donc le morphisme canonique $M^{\Sigma_0}\otimes_{\hRpi}\hRi\rightarrow M$ est surjectif.
Calquant la preuve de \ref{higgs1-desc66}, on montre que le morphisme 
\begin{equation}\label{higgs1-dsct0d}
\hRpi^d\rightarrow M^{\Sigma_0}
\end{equation}
défini par les $(x'_i)_{1\leq i\leq d}$ est surjectif. Par suite, $M^{\Sigma_0}$ est un $\hRpi$-module de type fini. 
Il est donc complet et séparé pour la topologie $p$-adique~; il  
est complet en vertu de (\cite{ac} chap. III §2.12 cor.~1 de prop.~16), et est séparé en tant que sous-module de $M$. 
D'autre part, comme $M$ est $\mZ_p$-plat, 
la topologie $p$-adique de $M^{\Sigma_0}$ est clairement induite par la topologie $p$-adique de $M$.
Donc la représentation de $\Delta_{p^\infty}$ sur $M^{\Sigma_0}$, induite par celle de $\Delta$,  
est continue pour la topologie $p$-adique.
Pour tout $1\leq i\leq d$ et tout $g\in \Delta_{p^\infty}$, on a $g(x'_i)-x'_i\in M^{\Sigma_0}\cap p^\alpha M=
p^\alpha M^{\Sigma_0}$. On en déduit que $M^{\Sigma_0}$ est une $\hRpi$-représentation $\alpha$-quasi-petite de $\Delta_{p^\infty}$,
d'où la proposition.

\begin{prop}\label{higgs1-dsct1}
Soient $\alpha,\varepsilon$ deux nombres rationnels tels que $0<2\varepsilon <\alpha$, 
$M$ une $\hoR$-représentation $\alpha$-quasi-petite de $\Delta$ telle que $M$ soit $\mZ_p$-plat. 
Alors, il existe une $\hRpi$-représentation $(\alpha-2\varepsilon)$-quasi-petite $M'$ de $\Delta_{p^\infty}$ 
et un morphisme $\hoR$-linéaire, $\Delta$-équivariant et surjectif
\begin{equation}\label{higgs1-dsct1a}
M'\otimes_{\hRpi}\hoR\rightarrow p^\varepsilon M.
\end{equation}
\end{prop}

D'après \ref{higgs1-desc6}(ii), il existe des générateurs $x_1,\dots,x_d$ de $p^\varepsilon M$ sur $\hoR$ 
qui sont d'une part $\Delta$-invariants modulo 
$p^{\alpha}M$ et d'autre part $\Sigma$-invariants. Pour tout $1\leq i\leq d$ et tout $g\in \Delta$, 
$g(x_i)-x_i\in  p^{\alpha}M \cap M^\Sigma=p^\alpha M^\Sigma$. 
Donc en vertu de \ref{higgs1-desc66} (appliqué à la $\hoR$-représentation 
$(\alpha-\varepsilon)$-quasi-petite $p^\varepsilon M$ de $\Sigma$), 
il existe $(a_{ij})_{1\leq j\leq d}\in \hRi^d$ tels que 
\begin{equation}\label{higgs1-dsct1b}
g(x_i)-x_i=\sum_{j=1}^dp^{\alpha-2\varepsilon}a_{ij}x_j
\end{equation}
Notons $M_1$ le sous-$\hRi$-module de $p^{\varepsilon} M$ engendré par $x_1,\dots, x_d$. 
D'après \eqref{higgs1-dsct1b}, $M_1$ est stable par l'action de $\Delta$.
Comme $M_1\subset M^\Sigma$, l'action induite de $\Delta$ sur $M_1$ se factorise à travers $\Delta_\infty$.

Le $\hRi$-module $M_1$ est complet et séparé pour la topologie $p$-adique~; il  
est complet en vertu de (\cite{ac} chap. III §2.12 cor.~1 de prop.~16), et est séparé en tant que sous-module de $M$.
Pour tout entier $n\geq 2\varepsilon$, on a 
\begin{equation}\label{higgs1-dsct1c}
p^n M_1\subset M_1\cap p^nM\subset p^nM^\Sigma \subset p^{n-2\varepsilon}M_1,
\end{equation}
la dernière inclusion étant une conséquence de \ref{higgs1-desc66}. 
Par suite, la topologie $p$-adique de $M_1$ est induite par la topologie $p$-adique de $M$. 
Compte tenu de \eqref{higgs1-dsct1b}, on en déduit que $M_1$ est une $\hRi$-représentation 
$(\alpha-2\varepsilon)$-quasi-petite de $\Delta_\infty$. Comme le morphisme canonique 
\begin{equation}\label{higgs1-dsct1d}
M_1\otimes_{\hRi}\hoR\rightarrow p^\varepsilon M
\end{equation}
est clairement surjectif, la proposition résulte de \ref{higgs1-dsct0} appliqué à $M_1$.

\begin{prop}\label{higgs1-dsct2}
Supposons $d=1$, et soit $M$ une petite $\hoR[\frac 1 p]$-représentation de $\Delta$. 
Alors il existe une petite $\hRun$-représentation $M'$ de $\Delta_\infty$, 
et un morphisme $\hoR[\frac 1 p]$-linéaire, $\Delta$-équivariant et surjectif 
\begin{equation}\label{higgs1-dsct2a}
M'\otimes_{\hRun}\hoR[\frac 1 p]\rightarrow M.
\end{equation}
\end{prop}
En effet, d'après \ref{higgs1-drt65}, il existe une $\hoR$-représentation quasi-petite $M^\circ$ de $\Delta$ telle que 
$M^\circ$ soit $\mZ_p$-plat, et un isomorphisme $\hoR[\frac 1 p]$-linéaire et $\Delta$-équivariant 
\begin{equation}\label{higgs1-dsct2b}
M^\circ\otimes_{\hoR}\hoR[\frac 1 p] \stackrel{\sim}{\rightarrow} M. 
\end{equation}
En vertu de \ref{higgs1-dsct1}, quitte à remplacer 
$M^\circ$ par $p^\varepsilon M^\circ$, pour un nombre rationnel $\varepsilon >0$, on peut supposer qu'il existe 
une $\hRpi$-représentation quasi-petite $N$ de $\Delta_{p^\infty}$, et un morphisme $\hoR$-linéaire, 
$\Delta$-équivariant et surjectif 
\begin{equation}\label{higgs1-dsct2c}
N\otimes_{\hRpi}\hoR\rightarrow M^\circ.
\end{equation}
Il existe un nombre rationnel $\alpha>\frac{1}{p-1}$ et des générateurs $x_1,\dots,x_d$ de $N$ sur $\hRpi$
qui sont $\Delta_{p^\infty}$-invariants modulo $p^{2\alpha}N$. On désigne par $N'$ le $\hRpi$-module libre
de base $e_1,\dots,e_d$ et par $\sigma\colon N'\rightarrow N$ le morphisme $\hRpi$-linéaire qui envoie
$e_i$ sur $x_i$ pour tout $1\leq i\leq d$. Notons 
\begin{equation}\label{higgs1-dsct2d}
\varphi_0\colon \Delta_{p^\infty}\rightarrow \Aut_{\hRun}(N')
\end{equation}
la $\hRpi$-représentation triviale de $\Delta_{p^\infty}$ relativement à la base $e_1,\dots,e_d$.  
Choisissons une $\mZ_p$-base $\gamma$ de $\Delta_{p^\infty}$. 
Il existe un automorphisme $\hRpi$-linéaire $u$ de $N'$ tel que pour tout $1\leq i\leq d$, on ait 
$u(e_i)-e_i\in p^{2\alpha} N'$ et 
\begin{equation}\label{higgs1-dsct2e}
\sigma(u(e_i))=\gamma(x_i).
\end{equation}  
Pour tout $g\in \Delta_{p^\infty}$, on note ${^g u}$  l'automorphisme $\hRpi$-linéaire 
$\varphi_0(g)\circ u\circ \varphi_0(g^{-1})$ de $N'$. 
Si $A\in \GL_d(\hRpi)$ est la matrice de $u$ relativement à la base $e_1,\dots, e_d$ de $N'$, 
alors $g(A)$ est la matrice de ${^g u}$ relativement à la même base.  
Soient $r$ un entier $\geq 0$, $A_r$ la classe de $A$ dans $\GL_d(R_{p^\infty}/p^rR_{p^\infty})$. 
Il existe un entier $m$ tel que $A_r$ appartiennent à $\GL_d(R_{p^m}/p^rR_{p^m})$, de sorte que 
$\gamma^{p^m}(A_r)=A_r$. Donc pour tout $n\geq 0$, on a 
\begin{equation}\label{higgs1-dsct2f}
A_r\gamma(A_r)\gamma^2(A_r)\dots \gamma^{p^{m+n}-1}(A_r)= 
(A_r\gamma(A_r)\gamma^2(A_r)\dots \gamma^{p^{m}-1}(A_r))^{p^n}.
\end{equation}
Comme $A\equiv \id \mod(p^{2\alpha} \hRpi)$, pour $n$ suffisamment grand, le produit \eqref{higgs1-dsct2f}
est égal à l'identité de $\GL_d(R_{p^\infty}/p^rR_{p^\infty})$. 
Par suite, pour tout $y\in N'$, la suite d'éléments de $N'$
\begin{equation}
n\mapsto (u\circ \varphi_0(\gamma))^{p^n}(y)=u\circ ({^\gamma u})\circ \dots\circ 
({^{\gamma^{p^n-1}} u})\circ\varphi_0(\gamma^{p^n})(y)
\end{equation}
converge vers $y$, pour la topologie $p$-adique. On en déduit que l'homomorphisme 
\begin{equation}
\mZ \rightarrow \Aut_{\hRun}(N'),\ \ \ n\mapsto (u\circ \varphi_0(\gamma))^n
\end{equation}
se prolonge en une $\hRpi$-représentation $\varphi$ de $\Delta_{p^\infty}$ sur $N'$, 
où l'on identifie $\mZ$ à un sous-groupe de $\Delta_{p^\infty}$ par l'injection $n\mapsto \gamma^n$. 
Il est clair que $\varphi$ est une représentation continue pour la topologie $p$-adique de $N'$, et est même petite. 
De plus, comme $\varphi(\gamma)=u\circ \varphi_0(\gamma)$, le morphisme $\sigma$ est $\Delta_{p^\infty}$-équivariant
\eqref{higgs1-dsct2e}. La proposition résulte alors de \ref{higgs1-desc4}.

\begin{prop}\label{higgs1-dsct3}
Supposons $d=1$ et $\hoR[\frac 1 p]$ fidèlement plat sur $\hRun[\frac 1 p]$ {\rm (cf. \ref{higgs1-drt18}(i))}. 
Alors, les propositions de \eqref{higgs1-dsct4} sont équivalentes à la proposition suivante~: 
\begin{itemize}
\item[$(\star)$] Pour toutes petites $\hoR[\frac 1 p]$-représentations $M$ et $M'$ de $\Delta$, 
et tout morphisme $\hoR[\frac 1 p]$-linéaire, $\Delta$-équivariant et surjectif $u\colon M'\rightarrow M$, 
la $\hoR[\frac 1 p]$-représentation de $\Delta$ sur le noyau de $u$ est petite. 
\end{itemize}
\end{prop}

D'après \ref{higgs1-dsct4} et \ref{higgs1-dsct6}, il suffit de montrer que $(\star)$ implique  \ref{higgs1-dsct4}(i). 
Soit $M$ une petite $\hoR[\frac 1 p]$-représentation de $\Delta$. 
D'après \ref{higgs1-dsct2}, il existe une petite $\hoR$-représentation $M'$ de $\Delta$ et 
un morphisme $\hoR[\frac 1 p]$-linéaire, $\Delta$-équivariant et surjectif
\begin{equation}\label{higgs1-dsct3b}
u\colon M'\otimes_{\hoR}\hoR[\frac 1 p]\rightarrow M.
\end{equation}
Par hypothèse, la $\hoR[\frac 1 p]$-représentation de $\Delta$ sur le noyau de $u$ est petite.
Appliquant de nouveau \ref{higgs1-dsct2}, on obtient une petite $\hoR$-représentation $M''$ de $\Delta$ 
et une suite exacte de morphismes $\hoR[\frac 1 p]$-linéaires et $\Delta$-équivariants
\begin{equation}\label{higgs1-dsct3c}
M''\otimes_{\hoR}\hoR[\frac 1 p]\stackrel{v}{\longrightarrow} M'\otimes_{\hoR}\hoR[\frac 1 p]
\stackrel{u}{\longrightarrow} M\longrightarrow 0.
\end{equation}
Remplaçant $M''$ par $p^n M''$ pour un entier $n\geq 0$, on peut supposer que $v(M'')\subset M'$. 
On désigne par $(N,\theta)$ le $\hRun[\frac 1 p]$-module de Higgs à coefficients dans $\xi^{-1}\tOmega^1_{R/\co_K}$,
conoyau du morphisme
\begin{equation}\label{higgs1-dsct3d}
\mH(v)\colon \mH(M''\otimes_{\hoR}\hoR[\frac 1 p])\rightarrow \mH(M'\otimes_{\hoR}\hoR[\frac 1 p]).
\end{equation}
En vertu de \ref{higgs1-desc5}, la $\hoR$-représentation $M'$ de $\Delta$ est de Dolbeault, le $\hRun$-module de Higgs 
$\mH(M')$ est petit et soluble, et on a un $\hoR$-isomorphisme canonique fonctoriel et $\Delta$-équivariant
\eqref{higgs1-dolb6a}
\begin{equation}\label{higgs1-dsct3e}
\mV(\mH(M'))\stackrel{\sim}{\rightarrow}M'. 
\end{equation}
Par suite, d'après \eqref{higgs1-higgs13b}, on a un $\hoR$-isomorphisme fonctoriel
\begin{equation}\label{higgs1-dsct3f}
\mH(M')\otimes_{\hRun}\hoR \stackrel{\sim}{\rightarrow}M'.
\end{equation}
On en déduit un isomorphisme 
\begin{equation}\label{higgs1-dsct3g}
N\otimes_{\hRun}\hoR \stackrel{\sim}{\rightarrow}M.
\end{equation}
Comme $\hoR[\frac 1 p]$ est fidèlement plat sur $\hRun[\frac 1 p]$, 
$N$ est projectif de type fini sur $\hRun[\frac 1p]$. 
Par suite, en tant que quotient de $\mH(M'\otimes_{\hoR}\hoR[\frac 1 p])$, le $\hRun[\frac 1 p]$-module de Higgs
$(N,\theta)$ est petit. 
Donc en vertu de \ref{higgs1-drt8}(ii) et \eqref{higgs1-dsct3e}, on a un $\hoR[\frac 1 p]$-isomorphisme  $\Delta$-équivariant
$M \stackrel{\sim}{\rightarrow} \mV(N)$, d'où la proposition \ref{higgs1-dsct4}(i).

\begin{remas}\label{higgs1-drt18} 
{\rm (i)}\ On s'attend à ce que $\hoR[\frac 1 p]$ soit toujours fidèlement plat sur $\hRun[\frac 1 p]$. 
Cet énoncé a été établi par Tsuji si $\cM_X$ est défini par un diviseur à croisements normaux strict sur $X$. 

{\rm (ii)}\ Les énoncés \ref{higgs1-dsct2} et \ref{higgs1-dsct3} sont directement inspirés de
l'approche de Faltings  (\cite{faltings3} theo.~3, page 852). 

{\rm (iii)}\ On s'attend à ce que les propositions équivalentes de \ref{higgs1-dsct4} soient vérifiées pour tout $d$. 
\end{remas}

\section{Représentations de Hodge-Tate}\label{higgs1-RHT}

\subsection{}\label{higgs1-hyodo1}\index{Hyodo (Anneau de ---)}\index{101501@$\cC_\HT$}
Considérons la suite exacte \eqref{higgs1-log-ext17b}
\begin{equation}\label{higgs1-hyodo1a}
0\rightarrow (\pi\rho)^{-1}  \hoR\rightarrow \cE \rightarrow
\tOmega^1_{R/\co_K}\otimes_R \hoR(-1)\rightarrow 0
\end{equation}
et notons aussi $1\in \cE$ l'image de $1\in (\pi\rho)^{-1}  \hoR$.
Pour tout $n\geq 0$, posons \eqref{higgs1-not621}
\begin{equation}\label{higgs1-hyodo1b}
\iota_n\colon \rS^n(\cE)\rightarrow \rS^{n+1}(\cE)
\end{equation}
le morphisme $\hoR$-linéaire défini de la façon suivante.
Pour $n=0$, $\iota_0$ est le composé des injections canoniques
$\hoR\rightarrow  (\pi\rho)^{-1}  \hoR\rightarrow \cE$. Pour $n\geq 1$ et pour tous $x_1,\dots,x_n\in \cE$, on a
\begin{equation}\label{higgs1-hyodo1c}
\iota_n([x_1\otimes\dots\otimes x_n])= [1\otimes x_1\otimes\dots\otimes x_n].
\end{equation}
Les $\cH_n$ forment alors un système inductif. On pose 
\begin{equation}\label{higgs1-hyodo1d}
\cC_\HT=\underset{\underset{n\geq 0}{\longrightarrow}}{\lim}\ \rS^n(\cE).
\end{equation}
Les morphismes de multiplication $\rS^m(\cE)\otimes_{\hoR}\rS^n(\cE)\rightarrow \rS^{m+n}(\cE)$ 
($m,n\geq 0$) induisent une structure de $\hoR$-algèbre sur $\cC_\HT$. 
On notera que $\cC_\HT$ est naturellement une $\hoR$-représentation de $\Gamma$.
La représentation ainsi définie est appelée l'{\em anneau de Hyodo} (\cite{hyodo1} §1).

\begin{rema}\label{higgs1-hyodo2bis}
On notera que $p$ est inversible dans $\cC_\HT$, de sorte que $\cC_\HT\simeq \cC_\HT[\frac 1 p]$.
En effet, on a $1=\pi\rho\cdot (\pi\rho)^{-1}$, où $(\pi\rho)^{-1}\in (\pi\rho)^{-1}\hoR\subset \cE\subset \cC_\HT$.  
\end{rema}

\subsection{} \label{higgs1-hyodo5}
Pour tous entiers $i\geq 0$ et $n\geq 0$, on définit un morphisme $\hoR$-linéaire et $\Gamma$-équivariant
\begin{equation}\label{higgs1-hyodo5a}
\varkappa_{i,n}\colon \rS^n(\cE)\otimes_R\tOmega^i_{R/\co_K}(-i)\rightarrow \rS^{n-1}(\cE)\otimes_R\tOmega^{i+1}_{R/\co_K}(-i-1)
\end{equation}
par $\varkappa_{i,0}=0$, et si $n\geq 1$ par
\begin{equation}\label{higgs1-hyodo5b}
\varkappa_{i,n}([x_1\otimes\dots\otimes x_n]\otimes \omega)= \sum_{1\leq i\leq n} 
[x_1\otimes\dots\otimes x_{i-1}\otimes x_{i+1}\otimes\dots\otimes x_n]\otimes (u(x_i)\wedge \omega),
\end{equation}
où $x_1,\dots,x_n\in \cE$, $\omega\in \tOmega^n_{R/\co_K}(-n)$ et $u\colon  \cE \rightarrow
\tOmega^1_{R/\co_K}\otimes_R \hoR(-1)$ le morphisme canonique 
\eqref{higgs1-hyodo1a}. On vérifie aussitôt que l'on a $\varkappa_{i+1,n-1}\circ \varkappa_{i,n}=0$ et que le diagramme 
\[
\xymatrix{
{\rS^n(\cE)\otimes_R\tOmega^i_{R/\co_K}(-i)}\ar[r]^-(0.5){\varkappa_{i,n}}\ar[d]_{\iota_n\otimes \id}
&{\rS^{n-1}(\cE)\otimes_R\tOmega^{i+1}_{R/\co_K}(-i-1)}\ar[d]^{\iota_{n-1}\otimes \id}\\
{\rS^{n+1}(\cE)\otimes_R\tOmega^i_{R/\co_K}(-i)}\ar[r]^-(0.5){\varkappa_{i,n+1}}&
{\rS^n(\cE)\otimes_R\tOmega^{i+1}_{R/\co_K}(-i-1)}}
\]
est commutatif. On obtient alors par passage à la limite inductive un morphisme $\hoR[\frac 1 p]$-linéaire 
et $\Gamma$-équivariant
\begin{equation}\label{higgs1-hyodo5c}
\varkappa_{i}\colon \cC_\HT\otimes_R\tOmega^i_{R/\co_K}(-i)\rightarrow \cC_\HT\otimes_R\tOmega^{i+1}_{R/\co_K}(-i-1).
\end{equation}
Ces morphismes définissent un complexe $\cK^\bullet$ de $\hoR[\frac 1 p]$-modules en posant 
\begin{equation}\label{higgs1-hyodo5d}
\cK^i=\left\{
\begin{array}{clcr}
\cC_\HT\otimes_R\tOmega^i_{R/\co_K}(-i) &{\rm si}\ i\geq 0,\\
0&{\rm si}\ i< 0.
\end{array}
\right.
\end{equation}

\subsection{}\label{higgs1-hyodo7}\index{101507@$\rD^i(V)$}
Soit $V$ une $\mQ_p$-représentation de $\Gamma$. Pour tout entier $i$, on note $\rD^i(V)$ 
le $\hR[\frac 1 p]$-module défini par  
\begin{equation}\label{higgs1-hyodo7a}
\rD^i(V)=(V\otimes_{\mQ_p}\cC_\HT(i))^\Gamma.
\end{equation}
En prenant les $\Gamma$-invariants du complexe $V\otimes_{\mQ_p}\cK^\bullet(i)$ \eqref{higgs1-hyodo5d}, on obtient un complexe 
$\hR[\frac 1 p]$-linéaire, noté $\mD^i(V)$,  
\begin{equation}\label{higgs1-hyodo7b}
\rD^i(V)\rightarrow \rD^{i-1}(V)\otimes_R\tOmega^1_{R/\co_K}\rightarrow \rD^{i-2}(V)\otimes_R\tOmega^2_{R/\co_K}
\rightarrow \dots,
\end{equation}
où $\rD^i(V)$ est placé en degré $0$. On munit le $\hR[\frac 1 p]$-module $\oplus_{i\in \mZ}\rD^i(V)(-i)$ du  
champ de Higgs à coefficients dans $\tOmega^1_{R/\co_K}(-1)$ induit par $\varkappa_0$ \eqref{higgs1-hyodo5c}.  
Le complexe de Dolbeault du module de Higgs ainsi défini est $\oplus_{i\in \mZ}\mD^i(V)(-i)$ \eqref{higgs1-not6b}.

\begin{defi}[\cite{hyodo1} 2.1]\label{higgs1-hyodo8}\index{representation de Hodge-Tate@$\mQ_p$-représentation de Hodge-Tate}
On dit qu'une $\mQ_p$-représentation continue $V$ de $\Gamma$ est de {\em Hodge-Tate} 
si les conditions suivantes sont remplies~:  
\begin{itemize}
\item[{\rm (i)}] $V$ est un $\mQ_p$-espace vectoriel de dimension finie, muni de la topologie $p$-adique \eqref{higgs1-not54}. 
\item[{\rm (ii)}] Le morphisme canonique 
\begin{equation}\label{higgs1-hyodo8a}
\oplus_{i\in \mZ}\rD^i(V)\otimes_{\hR[\frac 1 p]}\cC_\HT(-i)\rightarrow V\otimes_{\mQ_p}\cC_\HT
\end{equation}
est un isomorphisme.
\end{itemize} 
\end{defi}

On peut faire les remarques suivantes~:

\subsubsection{}\label{higgs1-hyodo8b} Hyodo montre dans {\em loc. cit.} que pour toute 
$\mQ_p$-représentation continue de dimension finie $V$ de $\Gamma$, le morphisme \eqref{higgs1-hyodo8a} est injectif.  

\subsubsection{}\label{higgs1-hyodo8c}  Pour toute $\mQ_p$-représentation de Hodge-Tate $V$ de $\Gamma$,  
les $\hR[\frac 1 p]$-modules $\rD^i(V)$ $(i\in \mZ)$ sont localement libres de type fini 
et sont nuls sauf pour un nombre fini d'entiers $i$, appelés les {\em poids de Hodge-Tate} de $V$ 
(cf. \cite{brinon2} 4.2.7).

\begin{prop}\label{higgs1-hyodo9}
Soient $(\tX_0,\cM_{\tX_0})$ la $(\cA_2(\oS),\cM_{\cA_2(\oS)})$-déformation lisse de 
$(\coX,\cM_\coX)$ 
définie par la carte $(P,\gamma)$ \eqref{higgs1-ext24ac}, 
$(P^\gp/\mZ\lambda)_{\lib}$ le quotient de $P^\gp/\mZ\lambda$ par son sous-module de torsion, 
\begin{equation}
w\colon (P^\gp/\mZ\lambda)_{\lib}\rightarrow P^\gp
\end{equation}
un inverse à droite du morphisme canonique $P^\gp\rightarrow (P^\gp/\mZ\lambda)_\lib$.  
On note $\cC_0$ la $\hoR$-algèbre de Higgs-Tate \eqref{higgs1-tor201}
et $\cF_0$ la $\hoR$-extension de Higgs-Tate \eqref{higgs1-tor201} associées à $(\tX_0,\cM_{\tX_0})$.
Il existe alors un isomorphisme canonique de $\hoR[\frac 1 p]$-algèbres 
\begin{equation}\label{higgs1-hyodo9b}
\tbeta_w\colon \cC_0[\frac 1 p]\rightarrow \cC_\HT
\end{equation}
tel que pour tout $x\in p^{-\frac{1}{p-1}}\cF_0\subset \cC_0[\frac 1 p]$, on ait
\begin{equation}\label{higgs1-hyodo9c}
\tbeta_w(x)=\beta_w(x)\in \cE\subset \cC_\HT,
\end{equation}
où $\beta_w$ est le morphisme \eqref{higgs1-ext246b}. 
De plus, $\tbeta_w$ est $\Delta$-équivariant et le diagramme 
\begin{equation}\label{higgs1-hyodo9d}
\xymatrix{
{\cC_0}\ar[rr]^{\tbeta_w}\ar[d]_{d_{\cC_0}}&&{\cC_\HT}\ar[d]^{\varkappa_0}\\
{\xi^{-1}\cC_0\otimes_R\tOmega^1_{R/\co_K}}\ar[r]^-(0.5)v_-(0.5)\sim&
{p^{-\frac{1}{p-1}}\cC_0\otimes_R\tOmega^1_{R/\co_K}(-1)}\ar[r]^-(0.5){\tbeta_w\otimes \id}&
{\cC_\HT\otimes_R\tOmega^1_{R/\co_K}(-1)}}
\end{equation}
où $v$ est l'isomorphisme induit par \eqref{higgs1-ext11c}  est commutatif.
\end{prop}

Cela résulte de \eqref{higgs1-tor2g} et \ref{higgs1-ext246}.

\begin{prop}\label{higgs1-hyodo10}
Conservons les hypothèses de \eqref{higgs1-hyodo9}, soient de plus  
$V$ une $\mQ_p$-représentation de Hodge-Tate de $\Gamma$,
$(\mH_0(V\otimes_{\mZ_p}\hoR),\theta)$ le module de Higgs à coefficients dans 
$\xi^{-1}\tOmega^1_{R/\co_K}$ associé à la $\hoR[\frac 1 p]$-représentation 
$V\otimes_{\mZ_p}\hoR$ de $\Delta$ par le foncteur  \eqref{higgs1-dolb2c} relatif à la déformation $(\tX_0,\cM_{\tX_0})$, 
$\theta'$ le champ de Higgs sur $\mH_0(V\otimes_{\mZ_p}\hoR)$ 
à coefficients dans $\tOmega^1_{R/\co_K}(-1)$ déduit de $\theta$ et de l'isomorphisme
$\hoR(1)\stackrel{\sim}{\rightarrow} p^{\frac{1}{p-1}}\xi \hoR$ \eqref{higgs1-ext11c}.  
Alors $V\otimes_{\mZ_p}\hoR$ est une $\hoR[\frac 1 p]$-représentation de Dolbeault de $\Delta$, et 
l'on a un $\hRun[\frac 1 p]$-isomorphisme canonique fonctoriel de $\hRun[\frac 1 p]$-modules de Higgs à
coefficients dans $\tOmega^1_{R/\co_K}(-1)$ 
\begin{equation}\label{higgs1-hyodo10a}
\oplus_{i\in \mZ}\rD^i(V)\otimes_{\hR}\hRun(-i)\stackrel{\sim}{\rightarrow}\mH_0(V\otimes_{\mZ_p}\hoR),
\end{equation}
où le membre de gauche est muni du champ de Higgs induit par $\varkappa_0$ \eqref{higgs1-hyodo5c} 
et le membre de droite est muni de $\theta'$.  
\end{prop} 

On notera d'abord que $V\otimes_{\mZ_p}\hoR$ est une $\hoR[\frac 1 p]$-représentation continue de $\Delta$. 
Si l'on munit $\rD^i(V)$ de l'action triviale de $\Gamma$, $V\otimes_{\mZ_p}\hoR$ du champ de Higgs trivial
et $\oplus_{i\in \mZ}\rD^i(V)\otimes_{\hR}\hoR(-i)$ et $\cC_\HT$ des champs de Higgs induits par $\varkappa_0$, 
le morphisme canonique 
\begin{equation}\label{higgs1-hyodo10b}
\oplus_{i\in \mZ}\rD^i(V)\otimes_{\hR[\frac 1 p]}\cC_\HT(-i) \rightarrow V\otimes_{\mQ_p}\cC_\HT
\end{equation}
est un isomorphisme $\cC_\HT$-linéaire et $\Gamma$-équivariant de $\hoR[\frac 1 p]$-modules de Higgs à coefficients dans 
$\tOmega^1_{R/\co_K}(-1)$. 
Compte tenu de \ref{higgs1-hyodo9}, 
on en déduit un isomorphisme $\cC_0^\dagger$-linéaire et $\Delta$-équivariant de 
$\hoR[\frac 1 p]$-modules de Higgs à coefficients dans $\tOmega^1_{R/\co_K}(-1)$
\begin{equation}\label{higgs1-hyodo10c}
\oplus_{i\in \mZ}\rD^i(V)\otimes_{\hR}\cC_0^\dagger(-i) \stackrel{\sim}{\rightarrow} 
V\otimes_{\mZ_p}\cC_0^\dagger,
\end{equation}
où $\cC_0^\dagger[\frac 1 p]$ est muni du champ de Higgs induit par $d_{\cC_0^\dagger}$ et par l'isomorphisme
$\hoR(1)\stackrel{\sim}{\rightarrow} p^{\frac{1}{p-1}}\xi \hoR$ \eqref{higgs1-ext11c}.  
Comme $\rD^i(V)$ est un facteur direct d'un $\hR[\frac 1 p]$-module libre de type fini \eqref{higgs1-hyodo8c}, 
on a $(\rD^i(V)\otimes_{\hR}\cC_0^\dagger)^\Delta=\rD^i(V)\otimes_{\hR}\hRun$ en vertu de \ref{higgs1-dolb47}. 
Par suite, en prenant dans \eqref{higgs1-hyodo10c} les invariants sous $\Delta$, on obtient un isomorphisme 
de $\hRun[\frac 1 p]$-modules de Higgs à coefficients dans $\tOmega^1_{R/\co_K}(-1)$
\begin{equation}\label{higgs1-hyodo10d}
\oplus_{i\in \mZ}\rD^i(V)\otimes_{\hR}\hRun(-i) \stackrel{\sim}{\rightarrow} 
(V\otimes_{\mZ_p}\cC_0^\dagger)^\Delta=\mH_0(V\otimes_{\mZ_p}\hoR).
\end{equation}
De plus, le morphisme canonique  
\begin{equation}\label{higgs1-hyodo10e}
\mH_0(V\otimes_{\mZ_p}\hoR)\otimes_{\hRun}\cC_0^\dagger\rightarrow V\otimes_{\mZ_p}\cC_0^\dagger
\end{equation}
s'identifie à l'isomorphisme \eqref{higgs1-hyodo10c}. Donc $V\otimes_{\mZ_p}\hoR$
est une $\hoR[\frac 1 p]$-représentation de Dolbeault de~$\Delta$.

\printindex

\end{document}